%% file: main.tex
\documentclass[
]{chairXthesis}
\pdfoutput=1 
\usepackage[utf8]{inputenc}    
\usepackage[T1]{fontenc}       
\usepackage{longtable}         
\usepackage{calc}
\usepackage{exscale}           
\usepackage{array}             
\usepackage{subcaption} 
\usepackage{wrapfig} 
\usepackage{gitinfo2} 
\usepackage{xspace}            
\usepackage{stackengine} 	
\setstackEOL{\\}
\usepackage{ifdraft}           
\usepackage{xparse} 
\usepackage{imakeidx} 
\usepackage{appendix}
\usepackage[symbols, record]{glossaries-extra}
\usepackage[expansion=false    
           ]{microtype}        
\usepackage[backend = biber,
backref = true,
style = alphabetic,
maxcitenames = 99,
maxbibnames = 99,
mincitenames = 5]{biblatex}
\usepackage{subfiles} 		
\usepackage[
draft,		
final,                    
plainpages=false,              
pdfpagelabels=true,            
pdfencoding=auto,              
unicode=true,                  
hypertexnames=true,            
naturalnames=true,              
colorlinks,
linkcolor={blue!70!black},
citecolor={blue!50!black},
urlcolor={blue!80!black}
]{hyperref}                    

\usetikzlibrary[patterns,patterns.meta]

\colorlet{c1}{blue}
\colorlet{c2}{green}
\colorlet{c3}{red!30!white}

\ifdraft{\synctex=1}{}

\input{names-phd-thesis}

\input{macros-phd-thesis}

\DeclareNameAlias{default}{last-first}

\makeindex[intoc,options= {-s phd-index.ist}]

\glsxtrsetgrouptitle{symbols}{Symbols}
\glsxtrsetgrouptitle{categories}{Categories}

\GlsXtrLoadResources[
src = {symbols},
sort = {none},
category = {symbols},
group = {symbols},
]
\GlsXtrLoadResources[
src = {categories},
sort = {none},
category = {symbols},
group = {categories},
]

\begin{document}

\selectlanguage{english}

%
%
\begin{titlingpage}
	\maketitle
	
	Updated version of the dissertation for the attainment of the doctoral degree at the University of Würzburg.
	Original version can be found at\newline
	\url{https://opus.bibliothek.uni-wuerzburg.de/frontdoor/index/index/docId/30167}.
	
\end{titlingpage}

%
%
\begin{titlingpage}
	\thesistitlepage
\end{titlingpage}

\thispagestyle{empty}
\noindent
\begin{tabular}{ l l }
	\textbf{Betreuer:}  & Prof. Dr. Stefan Waldmann (Julius-Maximilians-Universität Würzburg)  \\ 
	\textbf{1. Gutachter:} & Prof. Dr. Martin Bordemann  (Université de Haute-Alsace, Mulhouse)\\  
	\textbf{2. Gutachter:} & Prof. Dr. Luca Vitagliano (Università degli Studi di Salerno)
\end{tabular}
 
\clearpage

%
%
\begin{titlingpage}
	\innertitlepage
\end{titlingpage}

%
%
\renewcommand {\epigraphflush} {center}
\setlength{\epigraphwidth}{8cm}

\begin{vplace}[0.7]
	\epigraph{Everyone generalizes from one example.\newline At least, I do.}{ Vlad Taltos (\textit{Issola} by Steven Brust)}
\end{vplace}
\thispagestyle{empty}
\clearpage

\newpage\null\thispagestyle{empty}\newpage

%
%
\frontmatter

%
%
\tableofcontents*
\thispagestyle{plain}
\clearpage

%
%

\mainmatter
%
%
\chapter*{Introduction}
\addcontentsline{toc}{chapter}{Introduction}
\markboth{INTRODUCTION}{INTRODUCTION} 
\label{chap:Introduction}
\counterwithout{equation}{chapter} 
\counterwithout{figure}{chapter} 

\setcounter{page}{1}

\subfile{introduction}

\counterwithin{equation}{section}
\counterwithin{figure}{section}

%
%
\chapter{Constraint Algebraic Structures}
\label{chap:ConstraintAlgebraicStuctures}

\subfile{algebraic-structures}

%
%
\chapter{Constraint Geometric Structures}
\label{chap:ConstraintGeometricStructures}
\subfile{geometric-structures}

%
%
\chapter{Deformation Theory of Constraint Algebras}
\label{chap:DeformationTheory}

\subfile{deformation-theory}
\phantomsection
\chapter*{Outlook}
\label{chap:Outlook}
\addcontentsline{toc}{chapter}{Outlook}

\subfile{outlook}

%
%
\begin{appendices}

\chapter{Categorical Tools}
\label{chap:CategoryTheory}

\subfile{appendix-monoids}

\chapter{Poisson Geometry}
\label{chap:PoissonGeometry}
\renewcommand{\theequation}{\Alph{chapter}.\arabic{equation}}
\renewcommand{\thedefinition}{\Alph{chapter}.\arabic{definition}}
\renewcommand{\theproposition}{\Alph{chapter}.\arabic{proposition}}
\renewcommand{\thetheorem}{\Alph{chapter}.\arabic{theorem}}

\subfile{poisson-geometry}

\end{appendices}

%
%
\backmatter

%
%

{
  \footnotesize
  \printbibliography[heading=bibintoc]
}

%
%
\apptoglossarypreamble{\glsfindwidesttoplevelname[\currentglossary]}
\printunsrtglossary[style=alttreegroup]

%
%
\printindex

%
%

\ifdraft{\clearpage}
\ifdraft{\phantomsection}
\ifdraft{\addcontentsline{toc}{chapter}{List of Corrections}}
\ifdraft{\listoffixmes}

%
%

\end{document}


%% file: names-phd-thesis.tex
\title{Constraint Reduction in Algebra,\\Geometry and Deformation Theory}
\date{January 2023}
\author{Marvin Dippell}
\university{Julius-Maximilians-Universität Würzburg}
\institute{Institut für Mathematik}

%% file: macros-phd-thesis.tex
\makeatletter

\newcommand{\coloniff}{\mathrel{\vcentcolon\Leftrightarrow}}
\newcommand{\Eq}{\operatorname{eq}}
\newcommand{\Coeq}{\operatorname{coeq}}
\newcommand{\Ker}{\operatorname{Ker}}
\newcommand{\Coker}{\operatorname{Coker}}
\newcommand{\coev}{\operatorname{coev}}
\newcommand{\regimage}{\operator{regim}}
\newcommand{\vanishing}{\spacename{I}}
\newcommand{\Id}{\mathrm{Id}}
\newcommand{\lidentity}{\operatorname{left}}
\newcommand{\ridentity}{\operatorname{right}}
\newcommand{\commutative}{\mathrm{com}}
\newcommand{\normalizer}{\operatorname{N}}
\newcommand{\Shuffle}{\mathrm{Sh}}
\newcommand{\len}{\operatorname{len}}
\newcommand{\deform}[1]{\boldsymbol{#1}}
\newcommand{\discrete}{{\scriptstyle{\mathrm{dis}}}}
\newcommand{\standard}{\scriptstyle{\mathrm{std}}}

\newcommand{\MCset}{\functor{MC}}
\newcommand{\Def}{\functor{Def}}
\newcommand{\Hochschild}{\functor{HH}}
\newcommand{\diff}{\mathrm{diff}}
\newcommand{\qalgebra}[1]{{\pmb{\algebra{#1}}}}
\newcommand{\qvanishing}{{\pmb{\vanishing}}}
\newcommand{\qPnormalizer}{{\pmb{\Pnormalizer}}}

\newcommand{\Bott}{{\scriptstyle\mathrm{Bott}}}
\newcommand{\Ann}{\operatorname{Ann}}
\newcommand{\Orb}{\functor{Orb}}
\newcommand{\VecFields}{\mathfrak{X}}
\newcommand{\SymCoDer}{\mathrm{D}}
\newcommand{\deRham}{{\scriptstyle\mathrm{dR}}}
\newcommand{\Pnormalizer}{\mathcal{B}}
\renewcommand{\Schouten}[2][]{\ch@irxbbracket[#1]{#2}}

\newcommand{\Cat}{\categoryname{Cat}}
\newcommand{\Equiv}{\categoryname{Equiv}}
\newcommand{\LModules}[1]{\deco{}{#1}{\Modules}{}{}}
\renewcommand{\LeftModules}[1]{{#1}\textsf{-}\categoryname{Mod}}
\renewcommand{\RightModules}[2][]{\categoryname{Mod}_{#1}\textsf{-}#2}
\newcommand{\Proj}{\categoryname{Proj}}
\newcommand{\Monoids}{\mathsf{Mon}}

\newcommand{\CoChains}{\categoryname{Ch}}
\newcommand{\LieAlg}{\categoryname{LieAlg}}
\newcommand{\DGLieAlg}{\categoryname{DGLA}}
\newcommand{\GroupAct}{\Groups\categoryname{Act}}

\newcommand{\Naturals}{\mathbb{N}}
\newcommand{\Integers}{\mathbb{Z}}
\newcommand{\Reals}{\mathbb{R}}
\newcommand{\ComplexNum}{\mathbb{C}}

\newcommand{\Total}{{\scriptscriptstyle\script{T}}}
\newcommand{\Wobs}{{\scriptscriptstyle\script{N}}}
\newcommand{\Null}{{\scriptscriptstyle\script{0}}}
\newcommand{\TOTAL}{\script{T}}
\newcommand{\WOBS}{\script{N}}
\newcommand{\NULL}{\script{0}}
\newcommand{\str}{{\scriptstyle\mathrm{str}}}
\newcommand{\inj}{{\scriptstyle\mathrm{emb}}}

\newcommand{\sym}{{\scriptstyle\mathrm{sym}}}
\newcommand{\ext}{{\scriptstyle\mathrm{ext}}}

\newcommand{\Con}{\categoryname{C}}
\newcommand{\injCon}{\Con^\inj}
\newcommand{\strCon}{\Con_\str}
\newcommand{\injstrCon}{\Con^\inj_\str}

\newcommand{\ConSet}{{\Con\Sets}}
\newcommand{\injConSet}{{\injCon\Sets}}
\newcommand{\strConSet}{{\strCon\Sets}}
\newcommand{\injstrConSet}{{\injstrCon\Sets}}
\newcommand{\ConIndSet}{{\Con_{\scriptstyle\mathrm{ind}}\Sets}}
\newcommand{\injConIndSet}{{\Con^\inj_{\scriptstyle\mathrm{ind}}\Sets}}
\newcommand{\ConAlg}{{\Con\Algebras}}
\newcommand{\injConAlg}{{\injCon\Algebras}}
\newcommand{\strConAlg}{{\strCon\Algebras}}
\newcommand{\injstrConAlg}{{\injstrCon\Algebras}}

\newcommand{\ConGroups}{{\Con\Groups}}

\newcommand{\strConGroups}{{\strCon\Groups}}
\newcommand{\ConGroupAct}{{\Con\GroupAct}}
\newcommand{\ConMod}{\Con\Modules}
\newcommand{\injConMod}{{\injCon\Modules}}
\newcommand{\strConMod}{{\strCon\Modules}}
\newcommand{\injstrConMod}{{\injstrCon\Modules}}

\newcommand{\ConLMod}[1]{\deco{}{#1}{\Con\Modules}{}{}}

\newcommand{\strConLMod}[1]{\deco{}{#1}{\strCon\Modules}{}{}}

\newcommand{\ConRMod}[1]{\deco{}{}{\Con\Modules}{}{#1}}
\newcommand{\injConRMod}[1]{\deco{}{}{\injCon\Modules}{}{#1}}
\newcommand{\strConRMod}[1]{\deco{}{}{\strCon\Modules}{}{#1}}
\newcommand{\injstrConRMod}[1]{\deco{}{}{\injstrCon\Modules}{}{#1}}

\newcommand{\ConBimod}{{\Con\Bimodules}}
\newcommand{\injConBimod}{{\injCon\Bimodules}}
\newcommand{\strConBimod}{{\strCon\Bimodules}}
\newcommand{\injstrConBimod}{{\injstrCon\Bimodules}}

\newcommand{\ConProj}{\Con\Proj}
\newcommand{\strConProj}{\strCon\Proj}
\newcommand{\ConMfld}{{\Con\Manifolds}}

\newcommand{\strConMfld}{{\strCon\Manifolds}}
\newcommand{\ConVect}{{\Con\Vect}}

\newcommand{\ConLieAlg}{\Con\LieAlg}
\newcommand{\ConDGLieAlg}{\Con\DGLieAlg}

\newcommand{\ConMap}{{\Con\Map}}
\newcommand{\strConMap}{\strCon\Map}
\newcommand{\ConHom}{{\Con\!\Hom}}
\newcommand{\strConHom}{{\strCon\!\Hom}}
\newcommand{\ConCinfty}{\Con\Cinfty}
\newcommand{\strConCinfty}{\strCon\Cinfty}
\newcommand{\ConSecinfty}{\Con\Secinfty}
\newcommand{\strtensor}[1][{}]{\mathbin{\boxtimes_{\scriptscriptstyle{#1}}}}
\newcommand{\injstrtensor}[1][{}]{\mathbin{\boxtimes^\inj_{\scriptscriptstyle{#1}}}}
\newcommand{\injtensor}[1][{}]{\mathbin{\otimes^\inj_{\scriptscriptstyle{#1}}}}
\newcommand{\strcirctensor}[1][{}]{\mathbin{\otimes^\str_{\scriptscriptstyle{#1}}}}
\newcommand{\ConEnd}{{\Con\!\End}}

\newcommand{\ConAut}{{\Con\!\Aut}}
\newcommand{\ConDer}{{\Con\!\Der}}

\newcommand{\ConInnDer}{{\Con\!\InnDer}}
\newcommand{\ConOrb}{{\Con\Orb}}
\newcommand{\ConVecFields}{\Con\VecFields}
\newcommand{\ConForms}{\Con\Forms}
\newcommand{\ConDiffop}{\Con\!\Diffop}

\newlength{\grid}
\NewDocumentCommand{\ConGrid}{ O{0} O{0} O{0} O{0} O{0} O{0} O{0} 
O{0} O{0} }{
\mathbin{
\begin{tikzpicture}[x=9/ 3,y=9/3,rotate=45,baseline=3pt]
	\ifthenelse{#1=1}{\fill[black] (0,0) rectangle (1,1);}{}
	\ifthenelse{#1=2}{\fill[black!50] (0,0) rectangle (1,1);}{}
	\ifthenelse{#2=1}{\fill[black] (1,0) rectangle (2,1);}{}
	\ifthenelse{#2=2}{\fill[black!50] (1,0) rectangle (2,1);}{}
	\ifthenelse{#3=1}{\fill[black] (2,0) rectangle (3,1);}{}
	\ifthenelse{#3=2}{\fill[black!50] (2,0) rectangle (3,1);}{}
	\ifthenelse{#4=1}{\fill[black] (0,1) rectangle (1,2);}{}
	\ifthenelse{#4=2}{\fill[black!50] (0,1) rectangle (1,2);}{}
	\ifthenelse{#5=1}{\fill[black] (1,1) rectangle (2,2);}{}
	\ifthenelse{#5=2}{\fill[black!50] (1,1) rectangle (2,2);}{}
	\ifthenelse{#6=1}{\fill[black] (2,1) rectangle (3,2);}{}
	\ifthenelse{#6=2}{\fill[black!50] (2,1) rectangle (3,2);}{}
	\ifthenelse{#7=1}{\fill[black] (0,2) rectangle (1,3);}{}
	\ifthenelse{#7=2}{\fill[black!50] (0,2) rectangle (1,3);}{}
	\ifthenelse{#8=1}{\fill[black] (1,2) rectangle (2,3);}{}
	\ifthenelse{#8=2}{\fill[black!50] (1,2) rectangle (2,3);}{}
	\ifthenelse{#9=1}{\fill[black] (2,2) rectangle (3,3);}{}
	\ifthenelse{#9=2}{\fill[black!50] (2,2) rectangle (3,3);}{}
	
	\draw[step=1] (0,0) grid (3,3);
\end{tikzpicture}
}}

\makeatother


%% file: introduction.tex
Since the advent of quantum mechanics in the early twentieth century, physicists struggled to find a general scheme to
construct a  quantum mechanical analogue to a given classical system.
Such a \emph{quantization} procedure becomes necessary since,
even though nature seems to be inherently governed by quantum physics, our inability to directly perceive quantum mechanical
features forces us to rely on classical physics as a guideline to construct and interpret quantum mechanical systems.
In mathematical terms a classical mechanical system is often described by a manifold $M$, to be understood as the phase space of the
system, with points in $M$ representing individual states, together with a symplectic structure $\omega$ (or more generally a Poisson structure $\pi$) governing the dynamics of the
system \cite{abraham.marsden:1985a}.
The commutative algebra \glsadd{Cinfty}$\Cinfty(M)$
of real or complex functions on $M$  together with its Poisson bracket
\glsadd{PoissonBracket}$\{\argument, \argument\}$
is then interpreted as the algebra of observables.
On the other hand, the states of a quantum mechanical system are given by unit vectors in some Hilbert space $\mathcal{H}$
and its algebra of observables is given by the \emph{non-commutative} algebra of operators $\mathfrak{B}(\mathcal{H})$
on $\mathcal{H}$ with the induced commutator \glsadd{LieBracket}$[\argument, \argument]$.
A quantization is then generally supposed to yield, for a given classical system $(M,\omega)$,  a Hilbert space
$\mathcal{H}$ and a linear quantization map $Q \colon \Cinfty(M) \to \mathfrak{B}(\mathcal{H})$
fulfilling the following property \cite{ali.englis:2005a}:
\glsadd{hbar}
\begin{equation} \label{eq:QuantumClassicalBracket}
	Q(\{f,g\}) = \frac{1}{\I \hbar} [Q(f),Q(g)].
\end{equation}
The hope to find such a perfect quantization is destroyed by various no-go theorems, such as  the Groenewold-van Hove Theorem \cite{groenewold:1946a},
which forces us to weaken some of our assumptions.

\paragraph{Deformation Quantization}
Over the years many quantization schemes have been proposed, such as geometric quantization \cite{woodhouse:1997a}, $C^*$-algebraic deformation quantization \cite{rieffel:1994a}, strict deformation quantization \cite{landsman:1998a, rieffel:1989a} and convergent deformation quantization \cite{waldmann:2019a}.
In this thesis we will focus on formal deformation quantization.
In formal deformation quantization, as introduced in \cite{bayen.flato.fronsdal.lichnerowicz:1978a}, one assumes that \eqref{eq:QuantumClassicalBracket} only holds
asymptotically.
More precisely, given a Poisson manifold $M$ with Poisson bracket $\{\argument, \argument\}$ on $\Cinfty(M)$,
the algebra of complex-valued functions on $M$,
a \emph{star product} is an associative multiplication
\glsadd{starProduct}
\begin{equation} \label{eq:StarProduct}
	\star = \sum_{r=0}^{\infty} \colon \Cinfty(M)\formal{\lambda} \tensor[\mathbb{C}] \Cinfty(M)\formal{\lambda} \to \Cinfty(M)\formal{\lambda}
\end{equation}
of the form
\begin{equation}
	f \star g = \sum_{r=0}^{\infty} \lambda^r C_r(f,g)
\end{equation}
with bidifferential operators $C_r \colon \Cinfty(M) \tensor[\mathbb{C}] \Cinfty(M) \to \Cinfty(M)$
such that
\begin{propertieslist}
	\item $f \star g = fg + \sum_{r=1}^\infty \lambda^r C_r(f,g)$,
	\item $\frac{1}{\I \lambda}[f,g]_\star = \{f,g\} + \lambda(\dots)$,
	\item $1 \star f = f = f \star 1$.
\end{propertieslist}
The in general non-commutative algebra $(\Cinfty(M)\formal{\lambda},\star)$ is then interpreted as the observable algebra
of the quantized system.
In order to get in contact with the standard formulation of quantum mechanics we interpret $\lambda$ as a formal
replacement of $\hbar$, but then we still need to establish a suitable notion of convergence and find a representation on a (pre-)Hilbert space.
This leads to strict deformation quantization, which we will not discuss here, see
\cite{waldmann:2019a} for an overview.
Such star products are nothing but associative deformations of the commutative algebra $\Cinfty(M)$ in the sense of Gerstenhaber \cite{gerstenhaber:1964a}, and this deformation problem is governed by the Hochschild cohomology $\Hochschild(\Cinfty(M))$ \cite{hochschild:1945a}.
The existence and classification of star products was proved over the years for many situations,
see e.g. \cite{cahen.gutt:1982a, dewilde.lecomte:1983a} for the existence of star product on cotangent bundles.
In \cite{dewilde.lecomte:1983b} and \cite{fedosov:1994a} the existence on general symplectic manifolds was proven.
This development culminated in Kontsevich's Formality Theorem establishing the existence and classification of formal star products on general Poisson manifolds \cite{kontsevich:2003a}.

\paragraph{Geometric Reduction}
In classical mechanics symmetry reduction plays an important role.
Mathematically, this is usually phrased in terms of Marsden-Weinstein reduction \cite{marsden.weinstein:1974a} on a symplectic manifold $(M,\omega)$.
For this assume that a connected Lie group $\group{G}$ acts on $M$ in a Hamiltonian fashion, i.e.
there exists a momentum map $J \colon M \to \liealg{g}^*$, with $\liealg{g}$ denoting the Lie algebra of $\group{G}$,
such that
\begin{equation}
	\phi(\xi) = X_{J_\xi}
\end{equation}
for $\xi \in \liealg{g}$, and $\phi$ denoting the infinitesimal action of $\liealg{g}$.
If $0 \in \liealg{g}$ is a value and regular value of $J$, then $C \coloneqq J^{-1}(\{0\})$
is a closed submanifold of $M$.
Moreover, suppose that $\group{G}$ acts freely and properly on $C$, then
\begin{equation}
	M_\red \coloneqq C / \group{G}
\end{equation}
is a symplectic manifold with symplectic form $\omega_\red$ fulfilling $\pi^*\omega_\red = \iota^*\omega$,
with $\iota \colon C \to M$ the inclusion and $\pi \colon C \to M_\red$ the canonical projection.
It turns out that $C$ is a coisotropic submanifold of $M$ and that the above reduction procedure can actually be done for any coisotropic submanifold of a Poisson manifold.
Such coisotropic submanifolds and their reduction were introduced by Weinstein in \cite{weinstein:1988a} based on ideas
of Poisson reduction from \cite{marsden.ratiu:1086a}, see also \cite{stasheff:1997a}.
For this consider a Poisson manifold $(M,\pi)$, with $\argument^\sharp \colon T^*M \to TM$ denoting the corresponding musical homomorphism.
Then a submanifold $C$ of $M$ is coisotropic if and only if 
\begin{equation}
	\Ann(T_pC)^\# \subseteq T_pC
\end{equation}
for all $p \in C$.
Every such coisotropic submanifold carries a so-called characteristic distribution $D \subseteq TC$ spanned by the 
Hamiltonian vector fields $X_f$ for all functions $f$ vanishing on $C$.
\begin{figure}[t]
	\centering
	\includegraphics{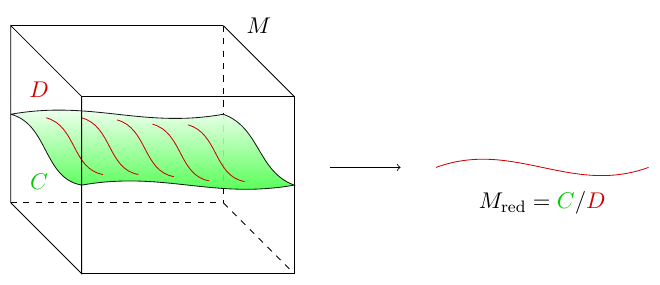}
	\caption{Reduction of coisotropic submanifold $C \subseteq M$ with characteristic distribution $D$.}
	\label{fig:ConstraintManifold}
\end{figure}
In the above case of the Marsden-Weinstein reduction of a symplectic manifold the leaves of this distribution
are given by the orbits of the group action on $C$.
If the characteristic distribution is nice enough, we can construct a reduced manifold
\begin{equation}
	M_\red \coloneqq C / D,
\end{equation}
see \autoref{fig:ConstraintManifold}, which carries a Poisson structure $\pi_\red$ induced by the Poisson structure $\pi$ on $M$.

\paragraph{Quantization vs. Reduction}
The question arises how quantization relates to (symmetry) reduction.
In particular, does quantization commute with reduction?
Can we quantize the reduction data on the classical side to allow for some kind of reduction procedure on the quantized system such that it does not matter whether we first quantize and then reduce or first reduce and then construct the corresponding quantum system?
In other words: does the following diagram commute?
\begin{equation}
\begin{tikzcd}
	\Centerstack[l]{classical physics\\with symmetries}
		\arrow[r,"Q"]
		\arrow[d,"\red"{swap}]
	&\Centerstack[l]{quantum physics\\with symmetries}
		\arrow[d,"\red"]\\	
	\text{classical physics}
		\arrow[r,"Q"]
	&\text{quantum physics}
\end{tikzcd}
\end{equation}
This question has been asked, and sometimes answered, for many different notions of reduction and quantization, see e.g.
\cite{meinrenken:1996a, guillemin.sternberg:1982c} for the case of geometric quantization,
\cite{bordemann.herbig.waldmann:2000a}
for a BRST-type reduction in deformation quantization
and \cite{bordemann:2004a:pre, bordemann:2005a}
for the symplectic case in deformation quantization.
In this thesis we want to focus on coisotropic reduction in the Poisson setting and its relation to formal deformation quantization.
More precisely, we ask under which conditions a given star product $\star$ on a Poisson manifold $(M,\pi)$ equipped with a coisotropic submanifold $C \subseteq M$ induces a star product $\star_\red$ on the reduced manifold $M_\red$.
Moreover, we want to clarify if such compatible star products exist and how equivalence of such star products may be investigated.

In \cite{cattaneo.felder:2007a} a similar situation is considered.
There a resolution of $\Cinfty(C)$ by means of the conormal bundle of $C$ is constructed.
This resolution carries a $P_\infty$-structure which, under certain conditions, can be shown to induce a deformation
of $\Cinfty(C)$.
Note however, that this approach only uses an infinitesimal neighbourhood of $C$, while we are interested in honest global star products allowing for a reduction.

\paragraph{Algebraic Reduction}
The general strategy is now to reformulate the geometric situation of a coisotropic submanifold equipped with its characteristic distributions in algebraic terms, similar to the way the algebra of observables $\Cinfty(M)$ is used to algebraically describe the manifold $M$.
Any closed submanifold $C \subseteq M$ can be described in terms of functions by its vanishing ideal
\glsadd{vanishingIdeal}
\begin{equation}
	\vanishing_C = \left\{ f \in \Cinfty(M) \bigm| f\at{C} = 0 \right\} \subseteq \Cinfty(M).
\end{equation}
Similarly, the foliation induced by any distribution $D \subseteq TM$ can be encoded by the subalgebra
\glsadd{invariantFunctions}
\begin{equation}
	\Cinfty(M)^D = \{ f \in \Cinfty(M) \mid \Lie_X f = 0 \text{ for all } X \in \Secinfty(D) \} \subseteq \Cinfty(M).
\end{equation}
Now for a coisotropic submanifold $C \subseteq (M,\pi)$ the characteristic distribution $D$ is only defined on $C$.
This leads us to consider the subalgebra
\glsadd{locallyInvariantFunctions}
\begin{equation}
	\Cinfty_D(M) = \{ f \in \Cinfty(M) \mid \Lie_X f\at{C} = 0 \text{ for all } X \in \Secinfty(D) \} \subseteq \Cinfty(M)
\end{equation}
instead.
Note that the vanishing ideal $\vanishing_C$ is contained in $\Cinfty_D(M)$.
Thus we have established a correspondence
\begin{equation} \label{eq:GeometryAlgebraCorrespondenceConstraint}
	(M,C,D)  \leftrightsquigarrow (\Cinfty(M), \Cinfty_D(M), \vanishing_C)
\end{equation}
between a manifold $M$ equipped with a closed submanifold $C$ on one side and a distribution $D$ on $C$ and
its algebra of functions $\Cinfty(M)$ equipped with the subalgebra $\Cinfty_D(M)$
of functions which are invariant on $C$ and the vanishing ideal $\vanishing_C$ on the other side.
Motivated by coisotropic reduction there is a reduction procedure for both sides of this correspondence.
Namely, on the geometric side, under the assumption of a simple distribution, 
we can construct the reduced manifold $M_\red \coloneqq C / D$ as before, 
while on the algebraic side we can always construct the reduced algebra
\begin{equation}
	\Cinfty(M)_\red = \Cinfty_D(M) / \vanishing_C.
\end{equation}
It is then easy to see that $\Cinfty(M)_\red \simeq \Cinfty(M_\red)$
is just the algebra of functions on the reduced manifold.
If we consider again the setting of a coisotropic submanifold $C \subseteq M$
then one can show that $\vanishing_C$ is a Poisson subalgebra of $(\Cinfty(M), \{\argument, \argument\})$
and that $\Cinfty_D(M)$ coincides with the Poisson normalizer
\begin{equation}
	\Pnormalizer_C = \{ f \in \Cinfty(M) \mid \{f,g\} \in \vanishing_C \text{ for all } g \in \vanishing_C\}	
\end{equation}
of $\vanishing_C$.
In particular, $\vanishing_C$ becomes a Poisson ideal in the Poisson subalgebra $\Pnormalizer_C$, and
therefore $\Cinfty(M)_\red$ carries itself a Poisson bracket, which turns
$M_\red$ into a Poisson manifold.

\paragraph{Constraint Algebras and their Deformations}
We now seek to carry over the basic ideas of deformation quantization to this more structured situation.
This means we want to treat the triple $(\Cinfty(M), \Cinfty_D(M),\vanishing_C)$
as a single algebraic entity and study deformations of it.
Thus we use $(\Cinfty(M), \Cinfty_D(M),\vanishing_C)$ as the motivating example to define an (embedded)
\emph{constraint algebra} $\algebra{A}$ as consisting of a unital associative algebra $\algebra{A}_\Total$
together with a unital subalgebra $\algebra{A}_\Wobs$ and an ideal $\algebra{A}_\Null \subseteq \algebra{A}_\Wobs$.
The subscript $\WOBS$ is supposed to remind the reader of the coisotropic situation, where $\algebra{A}_\Wobs$ is given as 
the Poisson normalizer of the Poisson subalgebra $\algebra{A}_\Null$.

In a next step we can try to define formal deformations of constraint algebras by taking the classical definition of a formal deformation and formally replace algebras by constraint algebras.
In particular, replacing $\Cinfty(M)$ by the constraint algebra $(\Cinfty(M), \Cinfty_D(M), \vanishing_C)$
in \eqref{eq:StarProduct} should yield the definition of a constraint star product.
To make sense of this we first need to clarify some notions:
\begin{cptitem}
	\item What are modules over constraint algebras and their tensor products?
	\item What are constraint multidifferential operators?
	\item Is there a cohomology theory governing the deformation problem of constraint algebras?
\end{cptitem}
We will answer these questions by taking a categorical point of view:
constraint algebras can be realized as monoid objects internal to a certain monoidal category 
$\ConMod_\field{k}$ equipped with a tensor product $\tensor[\field{k}]$, whose objects will be called
constraint $\field{k}$-modules.
By abstract categorical considerations, the definition and some first properties
of constraint modules over constraint algebras, as well as their tensor products, are then fixed.
In contrast to classical categories of modules, the categories of constraint modules will not form abelian categories.
This will lead to effects not present in classical module theory, and forces us to thoroughly examine even the most basic constructions of constraint modules.

This categorical approach will immediately allow us to find constraint analogues of many other classical algebraic concepts, such as derivations, groups, vector spaces, Lie algebras etc.
All these constraint notions will consist of a classical object as $\TOTAL$-component, together with a subobject as
$\WOBS$-component and an equivalence relation or ideal as $\NULL$-component.
Then, by construction, there is always a reduction procedure, defined by taking the quotient of the $\WOBS$-component by the $\NULL$-component, which by definition always yields a classical object.
It should be noted that the motivating example $(\Cinfty(M), \Cinfty_D(M),\vanishing_C)$
has additional properties not accounted for in the definition of constraint algebras, namely that
$\vanishing_C$ is an ideal not only in $\Cinfty_D(M)$ but in all of $\Cinfty(M)$.
Such constraint algebras will be called strong, and their modules will allow for two different canonical tensor products
$\tensor$ and $\strtensor$, whose interplay will be an important piece of study.
Note however that, since we are interested in non-commutative deformations of constraint algebras we should not expect
$\algebra{A}_\Null$ to stay a two-sided ideal in $\algebra{A}_\Total$ after deformation, see Lu's coisotropic creed \cite{lu:1993a}.
Even in the classical geometric situation we will encounter examples of honest non-strong constraint algebras.

Having found suitable notions of modules over constraint algebras we can introduce constraint differential operators
using Grothendieck's algebraic definition and thus arrive at the definition of constraint star products in analogy to
\eqref{eq:StarProduct}, which is nothing but a formal deformation of the constraint algebra $(\Cinfty(M), \Cinfty_D(M),\vanishing_C)$
by constraint differential operators.

The classical deformation theory of $\Cinfty(M)$ is governed by the differentiable Hochschild cohomology
$\Hochschild_\diff(\Cinfty(M))$, which can be computed by the Hochschild-Kostant-Rosenberg Theorem \cite{hochschild.kostant.rosenberg:1962a,gutt.rawnsley:1999a}, proving the existence of an isomorphism
\glsadd{sections}
\begin{equation} \label{eq:HKR}
	\Hochschild_\diff^\bullet(\Cinfty(M)) \simeq \Secinfty(\Anti^\bullet TM)
\end{equation}
of Gerstenhaber algebras, identifying the Hochschild cohomology with multivector fields on $M$.
When we want to find a constraint analogue of the classical HKR Theorem, we have to make sense of both sides
of \eqref{eq:HKR} in the constraint setting.

\paragraph{Constraint Manifolds and Vector Bundles}
Consider again the correspondence \eqref{eq:GeometryAlgebraCorrespondenceConstraint}.
Here, on the algebraic side, we have a subobject together with an equivalence relation on the subobject which is compatible with the structure of the subobject in a suitable sense.
In our example we have a subalgebra and an ideal inside this subalgebra.
On the geometric side, the triple $(M,C,D)$ carries the same underlying structure:
A subobject $C \subseteq M$, i.e. a submanifold, together with an equivalence relation on $C$, which in our case comes from a distribution $D$ on $C$.
We will understand in the course of this thesis that both the geometric and the algebraic side of \eqref{eq:GeometryAlgebraCorrespondenceConstraint} can be derived from the notion of constraint sets.
In particular, $(M,C,D)$ can be understood as a constraint set equipped with geometric structure, while the constraint algebra
$(\Cinfty(M), \Cinfty_D(M), \vanishing_C)$ can be seen as a constraint set equipped with algebraic structure.
Therefore we will call $\mathcal{M} = (M,C,D)$ a constraint manifold.
From this point of view we can reformulate the correspondence \eqref{eq:GeometryAlgebraCorrespondenceConstraint}
as a functor
\begin{equation} \label{eq:ConstraintFunctionsFunctor}
\begin{split}
	\ConCinfty &\colon \ConMfld \to \ConAlg, \\
	\ConCinfty(\mathcal{M}) &\coloneqq (\Cinfty(M), \Cinfty_D(M),\vanishing_C),
\end{split}
\end{equation}
from the category of constraint manifolds to the category of constraint algebras.

The notion of constraint manifolds encompasses two extreme, but important cases, namely that of a submanifold 
$C \subseteq M$ without a distribution, described by $(M,C,0)$, and that of a distribution $D$ on $M$ without an additional submanifold, described by $(M,M,D)$.
When applying the functor $\ConCinfty$ we obtain
\begin{align}
	\ConCinfty(M,C,0) &= (\Cinfty(M), \Cinfty(M), \vanishing_C)\\
	\shortintertext{and}
	\ConCinfty(M,M,D) &= (\Cinfty(M), \Cinfty(M)^D, 0).
\end{align}
Thus we see that the information of the $\WOBS$-component on the geometric side, i.e. the submanifold $C$, is encoded in the $\NULL$-component on the algebraic side.
Conversely, the geometric $\NULL$-component, i.e. the distribution $D$, is described by the $\WOBS$-component of the algebra of functions.
Therefore, if we are searching for a common framework including both the geometric and algebraic information, even if we are only interested in submanifolds \emph{or} in distributions, we have to consider the full constraint triple.
In particular, we can \emph{not} expect to be able to separate the reduction problem into two independent problems taking care of restriction and quotients separately.

The notion of constraint manifolds suggests to also introduce constraint versions of other geometric concepts,
such as constraint vector bundles and, in particular, constraint tangent and cotangent bundles.
A constraint vector bundle $E$ over a constraint manifold $\mathcal{M} = (M,C,D)$ will consist
of a vector bundle $E_\Total \to M$ with a subbundle
$E_\Wobs \to C$ of the restricted vector bundle $\iota^\#E_\Total \to C$,
a subbundle $E_\Null \subseteq E_\Wobs$ and a holonomy-free partial $D$-connection $\nabla$ on
$E_\Wobs / E_\Null$.
One should think of $E_\Null$ and $\nabla$ to define an equivalence relation on $E_\Wobs$ such that the 
quotient is a vector bundle.
We thus get the reduced vector bundle
\begin{equation}
	E_\red = (E_\Wobs / E_\Null ) /\nabla
\end{equation}
by identifying fibres of $E_\Wobs / E_\Null$ along the leaves using the parallel transport of $\nabla$.
Examples of constraint vector bundles have been considered before under various names,
e.g. quotient data \cite{cabrera.ortiz:2022a} and infinitesimal ideal systems in
\cite{jotzlean.ortiz:2014a}.
See also \cite{marrero.padron.rodriguez-almos:2012a} for related structures in the study of Marsden-Weinstein reduction for symplectic-like Lie algebroids.
Similar to \eqref{eq:ConstraintFunctionsFunctor} we will obtain a constraint sections functor
\begin{equation}
	\ConSecinfty \colon \ConVect(\mathcal{M}) \to \ConMod(\ConCinfty(\mathcal{M})),
\end{equation}
which yields for any constraint vector bundle a constraint $\ConCinfty(\mathcal{M})$-module of sections.
These constraint modules of sections will allow for clear geometric interpretations.
In particular, the sections of the constraint tangent bundle $T\mathcal{M}$ will be given by
\begin{equation}
	\begin{split}
		\ConSecinfty(T\mathcal{M})_\Total 
		&= \Secinfty(TM),\\
		\ConSecinfty(T\mathcal{M})_\Wobs
		&= \left\{ X \in \Secinfty(TM) \mid X\at{C} \in \Secinfty(TC), [X,Y] \in \Secinfty(D) \text{ for all } Y \in \Secinfty(D) \right\},\\
		\ConSecinfty(T\mathcal{M})_\Null
		&= \left\{ X \in \Secinfty(TM) \mid X\at{C} \in \Secinfty(D) \right\}.
	\end{split}
\end{equation}
Here the partial $D$-connection is given by the Bott connection, which is holonomy-free if the leaf space is smooth.
Motivated by the classical Serre-Swan Theorem \cite{swan:1962a, nestruev:2020a} we will identify sections of constraint vector bundles as a certain class
of projective constraint modules. This will lead us to the first main theorem (see \autoref{thm:strConSerreSwan}):

\begin{maintheorem}[Constraint Serre-Swan Theorem]
The monoidal category of constraint vector bundles over a constraint manifold $\mathcal{M}$
is equivalent to the monoidal category of projective strong constraint modules over the constraint algebra
$\ConCinfty(\mathcal{M})$.
\end{maintheorem}

We obtained a similar result for projective non-strong constraint modules in \cite{dippell.menke.waldmann:2022a}, where the equivalence
to a category of certain systems of vector bundles was shown.
In our present terms these could be understood as strong constraint vector bundles over strong constraint manifolds, but these objects will not be studied in this thesis.

With the constraint Serre-Swan Theorem we can, at least roughly, make sense of the right-hand side of \eqref{eq:HKR}.
Moreover, since all constraint notions are by definition equipped with a reduction functor and all constraint analogues of classical constructions, such as taking sections, are designed to be compatible with reduction, we will be able to show that
taking sections of (constraint) vector bundles commutes with reduction.

\paragraph{Constraint Differential Operators and Symbol Calculus}
In order to understand the left-hand side of \eqref{eq:HKR} in the constraint world we want to study more deeply the
constraint differential operators $\ConDiffop^\bullet(\mathcal{M})$ of $\ConCinfty(\mathcal{M})$.
It will turn out that 
$\ConDiffop^\bullet(\mathcal{M})$  can be understood in geometric terms by using a constraint covariant derivative,
which will lead us to the second major result (see \autoref{thm:ConMultSymbolCalculus}):

\begin{maintheorem}[Constraint symbol calculus]
	Given a constraint covariant derivative $\nabla$ on a constraint manifold $\mathcal{M} = (M,C,D)$
	there is an isomorphism
	\begin{equation}
		\Op \colon \ConSecinfty(\Sym^\bullet_{\tensor}T\mathcal{M} \strtensor \cdots \strtensor \Sym^\bullet_{\tensor}T\mathcal{M})
		\to \ConDiffop^\bullet(\mathcal{M}).
	\end{equation}
\end{maintheorem}

Every symmetric tensor power $\Sym^\bullet_{\tensor} T\mathcal{M}$ corresponds to a differential operator with a single input.
To obtain general multidifferential operators we need to combine these using the tensor product $\strtensor$,
which will be defined for constraint vector bundles in a similar manner as for constraint modules.
Thus, for understanding the constraint symbol calculus we need to study both tensor products $\tensor$ and $\strtensor$ and their relationship.

\paragraph{Constraint Hochschild Cohomology}
Motivated by classical deformation theory we consider the constraint Hochschild complex
of the constraint algebra $\algebra{A}$ given by
\begin{equation}
	C^\bullet(\algebra{A}) \coloneqq \ConHom_\field{k}(\algebra{A}^{\tensor n},\algebra{A})
\end{equation}
and we will show that this actually carries a compatible Hochschild differential $\delta$, allowing to define the constraint Hochschild cohomology by
\begin{equation}
	\Hochschild^\bullet(\algebra{A}) \coloneqq \frac{\ker \delta}{\image \delta}.
\end{equation}
This constraint Hochschild cohomology will be shown to govern the deformation theory of $\algebra{A}$
in familiar ways. 
To make this more precise, note that $\Hochschild^\bullet(\algebra{A})$
is constructed out of the constraint algebra $\algebra{A}$
and thus carries itself the structure of a graded constraint module,
meaning that $\Hochschild^\bullet(\algebra{A})$
consists of a $\TOTAL$-, $\WOBS$- and $\NULL$-component.
The $\TOTAL$-component is just given by the classical Hochschild cohomology of $\algebra{A}$,
and therefore contains information about the deformation theory of 
$\algebra{A}_\Total$ without taking into account the additional reduction information.
This additional structure is now incorporated into the $\WOBS$-component.
In particular, $\Hochschild^2(\algebra{A})_\Wobs$
can be identified with equivalence classes of infinitesimal deformations of $\algebra{A}_\Total$ which 
preserve the reduction data, and thus can be reduced to infinitesimal deformations of
$\algebra{A}_\red$.
Similarly, $\Hochschild^3(\algebra{A})_\Wobs$ gives obstructions
to extending deformations which are compatible with reduction in such a way that they stay compatible.
Finally, $\Hochschild^2(\algebra{A})_\Null$
and $\Hochschild^3(\algebra{A})_\Null$
give those infinitesimal deformations and obstructions that vanish after reduction.
We will be able to identify the zeroth Hochschild cohomology with the constraint version of centre and
the first Hochschild cohomology with the constraint derivations.

Additionally, the constraint Hochschild complex $C^\bullet(\algebra{A})$ will admit a Gerstenhaber bracket
allowing us to interpret formal deformations of $\algebra{A}$ as Maurer-Cartan elements in an associated constraint differential graded Lie algebra.
The equivalence of formal deformations can then be reformulated using a suitable gauge action on the constraint set of Maurer-Cartan elements.

It should be noted that there exists a well-established deformation theory for diagrams of algebras, 
see \cite{fregier.markl.yau:2009,gerstenhaber.schack:1983a}.
Interpreting a constraint algebra $\algebra{A}$ as a span
${\algebra{A}_\red \twoheadleftarrow \algebra{A}_\Wobs \to \algebra{A}_\Total}$
one might consider deformations of this diagram as a deformation of constraint algebras.
However, the category of modules over a diagram, which is the main ingredient used in \cite{gerstenhaber.schack:1983a},
is always abelian and hence cannot agree with our notion of constraint modules.
It remains to be seen if elements of this theory can help to compute constraint Hochschild cohomology.

We will apply the above results on constraint Hochschild cohomology
to the constraint algebra $\ConCinfty(\mathcal{M})$
but restricting ourselves to the differentiable Hochschild complex
\begin{equation}
	C_\diff^\bullet(\ConCinfty(\mathcal{M})) \coloneqq \ConDiffop^\bullet(\mathcal{M})
\end{equation}
in order to obtain constraint star products.

Even though the interpretation of constraint Hochschild cohomology fits well into our general constraint scheme,
here also something unexpected happens:
despite all constraint objects so far had their $\WOBS$-component embedded as a subobject into the $\TOTAL$-component, this
will in general not be true for the constraint Hochschild cohomology.
There will still exist a map
\begin{equation}
	\iota_\Hochschild \colon \Hochschild_\diff^\bullet(\ConCinfty(\mathcal{M}))_\Wobs \to \Hochschild_\diff^\bullet(\ConCinfty(\mathcal{M}))_\Total,
\end{equation}
but it might not be injective.
This immediately leads to problems when searching for a constraint analogue of the HKR Theorem, because while
the $\WOBS$-component of the left hand side of \eqref{eq:HKR} seems not to be injected into the $\TOTAL$-component in general,
the obvious constraint generalizations of the right-hand side will.

Even though we will not be able to fully solve the problem of finding a constraint analogue of the HKR Theorem in this thesis,
we can get deeper insights into the problem by considering the situation of flat space.
Thus we want to study the constraint Hochschild cohomology for 
\begin{equation}
	\mathcal{M} = \Reals^n \coloneqq (\Reals^{n_\Total}, \Reals^{n_\Wobs} ,\Reals^{n_\Null})
	\quad\text{ with } n_\Total \geq n_\Wobs \geq n_\Null.
\end{equation}
We will be able to compute the constraint Hochschild cohomology up to degree two in this case,
which gives the final main result of this thesis (see \autoref{thm:LocalSecondConHochschild}):

\begin{maintheorem}
The second constraint Hochschild cohomology for $\Reals^n = (\Reals^{n_\Total}, \Reals^{n_\Wobs}, \Reals^{n_\Null})$
is given by
\begin{equation}
\begin{split}
	\Hochschild^2_\diff(\ConCinfty(\Reals^n))_{\Wobs}
	&\simeq \big( \Anti^2\ConSecinfty(T\Reals^n)_\Wobs 
	+ \ConSecinfty(T\Reals^{n_\Total}) \wedge \ConSecinfty(T\Reals^n)_\Null \big) \\
	&\quad\oplus \left( \bigoplus_{k=1}^\infty \Sym^k \Secinfty(T\Reals^{n_\Null}\at{\Reals^{n_\Wobs}}) \vee \Secinfty(T\Reals^{n_\Total - n_\Wobs}\at{\Reals^{n_\Wobs}})\right).
\end{split}
\end{equation}
\end{maintheorem}

The term $\Anti^2\ConSecinfty(T\Reals^n)_\Wobs$ should be interpreted as bivector fields on $\Reals^{n_\Total}$
for which both legs are separately compatible with reduction, and hence these contributions will yield
bivector fields on the reduced manifold.
In contrast $\ConSecinfty(T\Reals^{n_\Total}) \wedge \ConSecinfty(T\Reals^n)_\Null$
describes bivector fields for which at least one leg vanishes after reduction, 
meaning that these contributions will reduce to zero.
The third summand is symmetric and therefore it is not a bivector field, but should rather be interpreted as a higher order differential operator.
This shows that $\Hochschild_\diff^2(\ConCinfty(\Reals^n))_\Wobs$
can not sit injectively inside $\Hochschild_\diff^2(\ConCinfty(\Reals^n))_\Total = \Hochschild_\diff^2(\Cinfty(\Reals^{n_\Total}))$.
Moreover, these terms, when interpreted as bi\-differential operators, can in principal have arbitrary degrees of differentiation while the classical HKR Theorem tells us that only multidifferential operators of order one in each slot appear in cohomology.

Beside their applications in the study of reduction of star products, many of the introduced concepts 
lend themselves for the study of reduction in other areas.
For example, the introduction of constraint bimodules and their tensor product naturally leads to the question
of Morita theory of constraint algebras, which itself could be an important part of the study of representations of algebras compatible with reduction.
Some first result can be found in \cite{dippell:2018a, dippell.esposito.waldmann:2019a}.
Moreover, the notion of constraint projective module can be used to introduce and study $K_0$-theory compatible with reduction.
On the geometric side, constraint vector bundles and constraint Lie algebras could be used to study the reduction of Lie algebroids and related geometric objects.

\subsubsection*{Structure of the Thesis}
This thesis is structured into three chapters:
\begin{cptitem}
	\item \autoref{chap:ConstraintAlgebraicStuctures}:
	Starting from the notion of constraint sets we develop in a mostly categorical fashion
	various constraint versions of well-known classical notions, such as groups, $\field{k}$-modules, algebras, modules over algebras etc.
	Required basics from category theory are recalled in \hyperref[chap:CategoryTheory]{Appendix~\ref{chap:CategoryTheory}}.
	From the beginning we will introduce slight variations for every constraint notion, namely that of strong constraint and embedded (strong) constraint objects,
	where strong constraint objects are constraint objects with the additional property that the $\NULL$-component defines also an equivalence relation on the $\TOTAL$-component, and embedded (strong) constraint objects are (strong) constraint objects
	with $\WOBS$-component embedded into the $\TOTAL$-component.
	The necessity to study also non-embedded constraint objects comes from the constraint Hochschild cohomology as introduced above.
	Having defined these basic constraint notions we will study free and projective modules over
	(embedded strong) constraint algebras.
	This will lead to a characterization of projective modules by constraint versions of the dual basis theorem.
	\item \autoref{chap:ConstraintGeometricStructures}:
	In this second chapter we introduce and study constraint manifolds and vector bundles as geometric counterparts of the algebraic constraint objects in the first chapter.
	A constraint version of the Serre-Swan Theorem will make this duality between algebra and geometry precise.
	Building on this, we will introduce constraint differential forms and \mbox{(multi-)}vector fields, establishing a Cartan calculus on constraint manifolds.
	We will then use constraint covariant derivatives to establish a symbol calculus for constraint multidifferential operators on constraint manifolds.
	
	Readers mostly interested in the geometric side of the story can directly begin with this second chapter.
	However, some definitions, like that of constraint algebras and modules, will be needed to follow the exposition.
	Basics on coisotropic reduction for Poisson manifolds can be found in
	\hyperref[chap:PoissonGeometry]{Appendix~\ref{chap:PoissonGeometry}}
	\item \autoref{chap:DeformationTheory}:
	We will bring together the geometric and algebraic objects introduced in the first two chapters to study star products compatible with reduction.
	For this we will introduce constraint versions of Hochschild cohomology and study deformations of constraint algebras using techniques from the theory of differential graded Lie algebras.
	Finally, we will compute the lowest constraint Hochschild cohomologies in the flat case.
\end{cptitem}

Afterwards we will give an \hyperref[chap:Outlook]{outlook} on related open questions and possible paths for further studies.

\subsubsection*{Bibliographical Notes}

This thesis is based on three publications \cite{dippell.esposito.waldmann:2019a}, \cite{dippell.menke.waldmann:2022a}
and \cite{dippell.esposito.waldmann:2022a}.
Since the basic notions used there have somewhat changed over time, let us comment a bit on their relation to the current thesis.

A first version of the notion of constraint algebra was introduced in \cite{dippell:2018a} and \cite{dippell.esposito.waldmann:2019a}
under the name of \emph{coisotropic triple of algebras} as a tool to study the behaviour of Morita equivalence under reduction.
These coisotropic triples would now be called embedded constraint algebras, with the additional property of
$\algebra{A}_\Null$ being a left ideal in $\algebra{A}_\Total$.
The notion of bimodules over coisotropic triples of algebras as used in
\cite{dippell.esposito.waldmann:2019a} already coincides with the notion of constraint bimodules over constraint algebras,
and reduction functors for constraint algebras and modules were already defined.
Thus the bicategory of bimodules over coisotropic triples of algebras as constructed in
\cite{dippell.esposito.waldmann:2019a} can be understood as a subbicategory of the bicategory of constraint bimodules over constraint algebras.
The proofs can easily be carried over to the more general situation of constraint algebras.

In \cite{dippell.menke.waldmann:2022a} the notions of coisotropic triples of algebras and modules were replaced by
\emph{coisotropic algebras} and \emph{coisotropic modules}, which agree with what we call constraint algebras and constraint modules in this thesis.
Here also coisotropic index sets, now called constraint index sets, were introduced in order to study
free and projective coisotropic modules.
These results can be found in \autoref{sec:ConIndSets}, \autoref{sec:FreeConMod} and
\autoref{sec:ProjectiveConModules}.
The goal of \cite{dippell.menke.waldmann:2022a} was then to find a suitable notion of vector bundles
over what we would now call a constraint manifold, such that sections of these vector bundles correspond to 
constraint modules by some sort of Serre-Swan Theorem.
These vector bundles are similar to the constraint vector bundles we introduce in
\autoref{sec:ConVectorBundles}, but they do still differ in important aspects.
In particular, there seems to be no good notion of tangent bundles or dual bundles.
We will see in the course of this thesis that these deficiencies come from the fact that the algebraic analogue of constraint manifolds
is given by \emph{strong} constraint algebras, and hence when searching for the correct notion of vector bundles one should also consider \emph{strong} constraint modules.
Although the vector bundles introduced in \cite{dippell.menke.waldmann:2022a} do not agree with our objects of study many ideas and smaller results used in \autoref{chap:ConstraintGeometricStructures} are based on \cite{dippell.menke.waldmann:2022a}.

Finally, in \cite{dippell.esposito.waldmann:2022a} the formal deformation theory of what was still called coisotropic algebras was studied.
The introduction of constraint Hochschild cohomology and the deformation functor based on constraint DGLAs is based on
\cite{dippell.esposito.waldmann:2022a}.

Even though this thesis is based on these three publications, a considerable amount does appear here for the first time.
In particular the notions of \emph{strong} constraint algebras and related objects, together with the strong tensor product, as well as the notion of constraint vector bundles have not been studied before.
Also the three main results as introduced above have not appeared elsewhere.

\subsubsection*{Notation and Conventions}

We adopt the following conventions:
\begin{cptitem}
	\item If not specified otherwise,
	\glsadd{ring}$\field{k}$ denotes a commutative unital ring, and
	\glsadd{field}$\field{K}$ denotes an arbitrary field.
	\item We will often use the term \emph{classical} to denote standard, non-constraint objects.
	For example, a constraint algebra consists of three classical algebras.
	\item Constraint analogues of classical categories or functors will be denoted by the classical symbol with a preceding $\Con$.
	For example $\Algebras$ denotes the category of classical algebras while $\ConAlg$ denotes the category of constraint algebras.
	Similarly, $\Cinfty(M)$ is the classical algebra of functions on a manifold $M$, while
	$\ConCinfty(\mathcal{M})$ denotes the constraint algebra of functions on the constraint manifold $\mathcal{M}$.
	\item Forgetful functors are often denoted by 
	\glsadd{forgetfulFunctor}$\functor{U}$ and their left adjoint free functors by
	\glsadd{freeFunctor}$\functor{F}$.
	Exceptions occur when these functors need to be referenced at a later stage.
	\item All constraint constructions will admit a reduction functor.
	Every such reduction functor is denoted by $\red$, and we will specify its domain only if necessary.
	\item Manifolds are considered to be connected, smooth, and in particular 
	Hausdorff and second countable.
\end{cptitem}

\subsubsection*{Acknowledgements}
First of all I would like to thank my supervisor Stefan Waldmann for his guidance and support throughout the years,
for keeping the doors to his (virtual) office open, and for the many hours spent in front of the blackboard discussing about mathematics.

Great thanks go to all current and former members of ChairX, especially to Felix Menke and Matthias Schötz for sharing the office with me, to Knut Hüper for remembering his sugary debts and to Francesco Cattafi, Martina Flammer, Madeleine Jotz Lean, Gregor Schaumann, Markus Schlarb and Christiaan van de Ven for making ChairX my second home. 
Special thanks to Chiara Esposito who despite the distance was always close when advice was needed.

I also want to thank Martina, David Kern, Felix, Markus and Matthias for our online seminar that kept my mathematical spirits high during the years of the pandemic.
I am also grateful for the weekly meetings with Andreas Kraft and Jonas Schnitzer.
Let me also thank Anthony Giaquinto and Severin Barmeier for helpful remarks.

I would like to express my gratitude to Simone Gutt and Martin Bordemann for the opportunity to visit Brussels and Mulhouse,
as well as Paolo Aschieri for my time in Alessandria.
Particularly, I would like to mention Thomas Weber who introduced me to the beautiful cities of Naples and Turin.
A big thanks again to Chiara and the whole group in Salerno who hosted me for a few weeks and introduced me to the Italian way of life.

I also have to thank all my friends who kept me sane during the years of the pandemic and beyond, especially the ``outlaws'', Marc Herrmann, Veronika Karl and Frank Pörner, who even dared to take a walk with me during the lockdown.

I am grateful for all the helpful remarks by the referees Martin Bordemann and Luca Vitagliano.

Finally, I want to thank Sarah Gaß for designing the title page and all the busy proofreaders who, in an incredibly short time, found way more errors than I ever expected to be able to sneak into this text, especially Michael Lopin for trying to teach me English comma placement and Thorben Römer, to whom I owe the fact that I used the word ``constraint'' exactly $2136$ times in this thesis; well, now it is $2137$ times.

%% file: algebraic-structures.tex
We introduce numerous algebraic objects for which a reduction 
procedure can be defined.
These algebraic structures naturally appear in our study of 
deformation quantization of coisotropic submanifolds.
Even more, constraint algebras and their deformation theory will be 
our main objects of interest.

We start by introducing constraint sets in 
\autoref{sec:ConstraintSets}, which should be thought of
as traditional sets equipped with the structures required to allow
for a notion of reduction.
More explicitly, these will consist of a set $M_\Total$ together with 
a map $\iota_M \colon M_\Wobs \to M_\Total$ from another set 
$M_\Wobs$ which is itself equipped with an equivalence relation
$\sim_M$.
The reduction of constraint sets is then defined as
\begin{equation}
	M_\red \coloneqq M_\Wobs / \mathord{\sim}_M.
\end{equation}
A constraint set with injective $\iota_M$ will be called 
\emph{embedded}, while an additional equivalence relation on 
$M_\Total$ leads to the notion of \emph{strong} constraint sets.
Most examples from geometry will lead to embedded strong constraint 
sets, nevertheless, honest constraint sets, even non-embedded ones, 
will naturally appear.
Canonical constructions, such as limits and colimits as well as 
mono-, epimorphisms and notions of image and subobject will be 
studied for the various flavours of constraint sets.
It will be apparent that even though $\ConSet$ shares a lot of 
features with $\Sets$, it differs at important points, giving a first 
hint that introducing classical mathematical objects internal to 
$\ConSet$ might produce some unfamiliar effects.
 
We do not investigate constraint sets for their own sake, but as the 
foundation for all following notions appearing in this thesis.
Beginning with \autoref{sec:ConstraintkModules} we introduce
additional algebraic structure on constraint sets.
The idea is to follow the classical hierarchy of algebraic notions
but implement their categorical definitions in the category $\ConSet$ 
instead of the classical category $\Sets$.
Therefore, we start with constructing constraint (abelian) groups,
followed by constraint $\field{k}$-modules and their strong 
constraint cousins.
Unsurprisingly, it will turn out that these derived constraint notions
share a structural similarity with constraint sets.
For example a constraint $\field{k}$-module will be given by a 
$\field{k}$-module
$\module{E}_\Total$ together with a module morphism
$\iota_\module{E} \colon \module{E}_\Wobs \to \module{E}_\Total$
and an equivalence relation on $\module{E}_\Wobs$
compatible with the $\field{k}$-module structure.
Since in most algebraic categories equivalence relations compatible
with the algebraic structure can be understood as subobjects of a 
certain type, e.g. normal subgroups, submodules, ideals, etc.,
the equivalence relation on the $\WOBS$-component will mostly be 
replaced by such a subobject.
For constraint $\field{k}$-modules this means we consider 
a submodule $\module{E}_\Null$ of $\module{E}_\Wobs$.
The trinity of $\TOTAL$-, $\WOBS$- and $\NULL$-component will be 
prevalent in the rest of this work.
At this point we already encounter the two different tensor products
$\tensor$ and $\strtensor$ for constraint modules.
Their interplay and their mismatch alike
will have tremendous impact on the later chapters.

Before continuing to introduce constraint algebras and their modules 
we pause to take a closer look at constraint vector spaces and their 
bases in \autoref{sec:ConstraintLinAlg}.
For this it will be useful to first study constraint index sets in
\autoref{sec:ConIndSets}.
Then in \autoref{sec:ConVectorSpaces} the relation of $\tensor$ and 
$\strtensor$ for constraint vector spaces will become apparent.
These results will serve as a guideline for the later study of free 
and projective constraint modules.

In \autoref{sec:ConstraintAlgebrasModules} we proceed to define
(strong) constraint algebras as monoid objects internal to the 
category $\ConMod_{\field{k}}$ of constraint $\field{k}$-modules and 
introduce modules over such constraint algebras.

Even though the main example for constraint modules over constraint 
algebras, namely that of constraint manifolds and their vector 
bundles, will not be introduced until 
\autoref{chap:ConstraintGeometricStructures} the reader acquainted 
with classical differential geometry will anticipate the relevance of 
free and projective constraint modules.
Therefore, we will study these notions in both the strong and 
non-strong case and will find characterizations of projective 
constraint modules analogous to the classical situation, using a 
lifting property, as summands of free modules and as allowing for a 
sort of dual basis.

In the last section of this chapter we collect additional 
constraint notions which will be useful later on but which are either 
special cases of objects we studied before or whose definition and 
properties follow in a more or less straightforward way from what has 
been done before.
In particular, \autoref{sec:ConHomologicalAlgebra} contains the 
basics of graded constraint modules and foundational results for 
homological algebra of those.
At last, in \autoref{sec:ConDGLAs}, constraint (differential graded) 
Lie algebras and related structures are introduced.

\section{Constraint Sets}
\label{sec:ConstraintSets}

\subfile{constraint-sets}

\section{Constraint $\field{k}$-Modules}
\label{sec:ConstraintkModules}

\subfile{constraint-modules}

\section{Interlude: Constraint Linear Algebra}
\label{sec:ConstraintLinAlg}

\subfile{constraint-linalg}

\section{Constraint Algebras and Modules}
\label{sec:ConstraintAlgebrasModules}

\subfile{constraint-algebra-modules}

\section{Regular Projective Modules}
\label{sec:RegularProjectiveModules}

\subfile{projective-modules}

\section{More Constraint Structures}
\label{sec:MoreConstraintStructures}

\subfile{constraint-misc}

%% file: constraint-sets.tex
Consider the motivating example of a coisotropic submanifold $C$ of a 
Poisson manifold $M$.
Forgetting all the geometric structure and just remembering the bare 
set-theoretic minimum needed for reduction leaves us with the set 
$M$, a subset $C$ and an equivalence relation defined on $C$ given by 
the characteristic distribution.
This motivates the following definition.

\begin{definition}[Constraint set]\
	\label{def:ConstraintSet}
	\begin{definitionlist}
		\item A \emph{constraint set}
		\index{constraint!set}
		$M$ consists of a map
		$\iota_M \colon M_\Wobs \to M_\Total$ of sets, together with 
		an
		equivalence relation $\sim_M$ on $M_\Wobs$.
		\item A morphism $f \colon M \to N$ of constraint sets
		(or \emph{constraint morphism})
		consists of maps
		$f_\Total \colon M_\Total \to N_\Total$
		and
		$f_\Wobs \colon M_\Wobs \to N_\Wobs$
		such that
		$f_\Total \circ \iota_M = \iota_N \circ f_\Wobs$
		and
		$f_\Wobs$ preserves the equivalence relation, i.e.
		$f_\Wobs(x) \sim_N f_\Wobs(y)$ for all $x \sim_M y$.
		The set of constraint morphisms from $M$ to $N$ is denoted by
		\glsadd{maps}$\Map(M,N)$.
		\item The category of constraint sets and their morphisms
		is denoted by
		\glsadd{ConSet}$\ConSet$.
	\end{definitionlist}
\end{definition}

We will often suppress the map included in the definition of 
constraint sets and just write $M = (M_\Total,M_\Wobs,\sim_M)$.
Following our motivation it would be natural to include injectivity
of $\iota_M$ in the definition of constraint sets.
In fact, most examples of coisotropic sets will be of this form and
thus they will get their own name later on.
But injectivity of $\iota_M$ is not preserved under some important 
categorical constructions which are compatible with reduction.
Hence we excluded this property from the definition of constraint 
sets.

Let us collect some important properties of the category $\ConSet$.
For this we need the notion of pushforward and pullback of 
equivalence relations:
Let $f \colon M \to N$ be a map between sets, and let $\sim_M$ and 
$\sim_N$ be equivalence relations on $M$ and $N$, respectively.
We denote by $f^*(\sim_N) = \mathbin{\sim_{f^*}}$ the pullback equivalence 
relation on $M$ 
defined by
\begin{equation}
	x \sim_{f^*} x' \coloniff f(x) \sim_N f(x').
\end{equation}
In general we say that a map is compatible with the equivalence 
relations if $\mathbin{\sim_M} \subseteq f^*(\sim_N)$.
Moreover, by $f_*(\sim_M) = \mathbin{\sim_{f_*}}$ we denote the pushforward 
equivalence relation 
on $N$ given as the equivalence relation generated by
\begin{equation}
	f(x) \sim_{f_*} f(x')  \text{ for all } x \sim_M x'.
\end{equation}
Note that this implies that $f_*(\sim_M)$ is discrete outside of 
$\image(f)$.
The discrete equivalence relation will always be denoted by $\sim_\discrete$.
With this we can give a description of useful co/limits in
$\ConSet$,
see \autoref{ex:limitsandcolimits} for the general definitions.

\begin{proposition}[Co/limits in $\ConSet$]
	\label{prop:CoLimitsCSet}
Let $M, N, P$ be constraint sets, and let
$f,g \colon M \to N$ as well as
$h \colon P \to N$ be constraint morphisms.
\begin{propositionlist}
\item\label{prop:CoLimitsCSet_initial}
\index{initial object!constraint sets}
The initial object in $\ConSet$ is given by $(\emptyset,\emptyset, 
\sim)$, with $\sim$ the unique equivalence relation on $\emptyset$.
\item\label{prop:CoLimitsCSet_final}
\index{final object!constraint sets}
The final object in $\ConSet$ is given by 
$1 \coloneqq (\{\pt\},\{\pt\},\sim)$, with $\{\pt\}$ any one-element set and
$\sim$ the unique 
equivalence relation 
on $\{\pt\}$.
\item\label{prop:CoLimitsCSet_product}
\index{product!constraint sets}
The product is given by
\glsadd{product}
\begin{equation}
\begin{split}
	(M \times N)_\Total &= M_\Total \times N_\Total, \\
	(M \times N)_\Wobs &= M_\Wobs \times N_\Wobs,
\end{split}
\end{equation}
with the product map
$\iota_{M \times N}  = \iota_M \times \iota_N \colon M_\Wobs 
\times N_\Wobs
\longrightarrow M_\Total \times N_\Total$
and the product relation 
$\sim_{M \times N}$ 
given by 
\begin{equation}
	(x_1, y_1) \sim_{M \times N} (x_2,y_2)
	\coloniff
	x_1 \sim_M x_2 \text{ and } y_1 \sim_N y_2.
\end{equation}
\item\label{prop:CoLimitsCSet_coproduct}
\index{coproduct!constraint sets}
The coproduct is given by
\glsadd{coproduct}
\begin{equation}
\begin{split}
	(M \sqcup N)_\Total &= M_\Total \sqcup N_\Total, \\
	(M \sqcup N)_\Wobs &= M_\Wobs \sqcup N_\Wobs,
\end{split}
\end{equation}
with the coproduct map
$\iota_M \sqcup \iota_N \colon M_\Wobs \sqcup N_\Wobs 
\longrightarrow M_\Total \sqcup N_\Total$
and the relation 
\begin{equation}
	x \sim_{M \sqcup N} y
	\coloniff x \sim_M y \text{ or } x \sim_N y.
\end{equation}
Here $\sqcup$ denotes the disjoint union of sets.
\item\label{prop:CoLimitsCSet_pullback}
\index{pullback!constraint sets}
The pullback of $f$ and $h$ is given by the constraint set
\begin{equation}
\begin{split}
	(M \decorate*[_{f}]{\times}{_{h}} P)_\Total
	&= M_\Total \decorate*[_{f_\Total}]{\times}{_{h_\Total}} P_\Total
	= \left\{ (x,y) \in M_\Total \times P_\Total \mid f_\Total(x) = h_\Total(y) \right\},\\
	(M \decorate*[_f]{\times}{_h} P)_\Wobs
	&= M_\Wobs \decorate*[_{f_\Wobs}]{\times}{_{h_\Wobs}} P_\Wobs
	= \left\{ (x,y) \in M_\Wobs \times P_\Wobs \mid f_\Wobs(x) = h_\Wobs(y) \right\},
\end{split}
\end{equation}
with the relation $\sim_{\decorate*[_f]{\times}{_h}}$ given by
\begin{equation} \label{eq:EquivalenceRelPullbackConSet}
	(x_1,y_1) \sim_{\decorate*[_f]{\times}{_h}} (x_2, y_2)
	\quad\Leftrightarrow\quad
	x_1 \sim_M x_2 \text{ and } y_1 \sim_P y_2
\end{equation}
and projection maps
\begin{align}
	(\pr_\Total^M,\pr_\Wobs^M) 
	&\colon (M \decorate*[_{f}]{\times}{_{h}} P) \longrightarrow 
	M,\\
	(\pr_\Total^P,\pr_\Wobs^P) 
	&\colon (M \decorate*[_{f}]{\times}{_{h}} P) \longrightarrow 
	P.
\end{align}
\item\label{prop:CoLimitsCSet_equalizer}
\index{equalizer!constraint sets}
The equalizer of $f$ and $g$ is given by the constraint 
set
\glsadd{equalizer}
\begin{equation}
\begin{split}
	\Eq(f,g)_\Total
	&= \Eq(f_\Total, g_\Total) = \left\{ x \in M_\Total \mid f_\Total(x) = g_\Total(x)  \right\}, \\
	\Eq(f,g)_\Wobs
	&= \Eq(f_\Wobs, g_\Wobs)= \left\{ x \in M_\Wobs \mid f_\Wobs(x) = g_\Wobs(x)  \right\},
\end{split}
\end{equation}
with the equivalence relation given by the restriction of $\sim_M$
and the morphism\linebreak
$i =(i_\Total, i_\Wobs) \colon \Eq(f,g) \to M$
given by the inclusions $i_\Total$ and $i_\Wobs$
of $\Eq(f_\Total,g_\Total)$ 
and $\Eq(f_\Wobs,g_\Wobs)$
into $M_\Total$ and $M_\Wobs$, respectively.
\item\label{prop:CoLimitsCSet_coequalizer}
\index{coequalizer!constraint sets}
The coequalizer of $f$ and $g$ is given by the constraint 
set
\glsadd{coequalizer}
\begin{equation}
\begin{split}
	\Coeq(f,g)_\Total
	&= \Coeq(f_\Total, g_\Total), \\
	\Coeq(f,g)_\Wobs
	&= \Coeq(f_\Wobs, g_\Wobs),
\end{split}
\end{equation}
with the equivalence relation given by $(q_\Wobs)_\ast(\sim_N)$
with
$q_\Wobs \colon N_\Wobs \to	\Coeq(f,g)_\Wobs$
and the morphism
$q = (q_\Total,q_\Wobs) \colon N \to \Coeq(f,g)$
of constraint sets.
Here $q_\Wobs \colon N_\Wobs \to \Coeq(f_\Wobs,g_\Wobs)$
and $q_\Total \colon N_\Total \to \Coeq(f_\Total,g_\Total)$
denote the coequalizer in $\Sets$ of $f_\Wobs$, $g_\Wobs$ and
$f_\Total$, $g_\Total$, respectively.
More explicitly, $\Coeq(f_\Total,g_\Total)$ is given by
$N/\mathord{\sim}$ with $\sim $ the equivalence relation generated by
$y_1 \sim y_2$ if and only if 
there exist $x \in M$ such that $f(x) = y_1$ and $g(x) = y_2$.
\item\label{prop:CoLimitsCSet_coLimits}
The category $\ConSet$ has all finite limits and colimits.
\end{propositionlist}
\end{proposition}

\begin{proof}
Since the strategy to prove these statements is always the same,
we will not perform everything at great length.
Let us instead prove \ref{prop:CoLimitsCSet_initial} and 
\ref{prop:CoLimitsCSet_pullback} in detail, then the rest should be 
clear.

Since $\emptyset$ is the initial object in $\Sets$ we know that there 
exist unique maps $\emptyset \to M_\Total$ and
$\emptyset \to M_\Wobs$.
It also follows by the uniqueness that 
\begin{equation*}
\begin{tikzcd}
	\emptyset
		\arrow[r]
	& M_\Total \\
	\emptyset
		\arrow[u,"\id"]
		\arrow[r]
	& M_\Wobs
		\arrow[u,"\iota_M"{swap}]
\end{tikzcd}
\end{equation*}
commutes.
Moreover, $\emptyset \to M_\Wobs$ is clearly compatible with the equivalence relations.
Thus we obtain a constraint morphism $(\emptyset, \emptyset,\sim_\emptyset) \to M$
which is unique, since its components are.

For \ref{prop:CoLimitsCSet_pullback} consider another constraint set $X$ with constraint morphisms
$\phi \colon X \to M$ and $\psi \colon X \to P$ such that
$f \circ \phi = h \circ \psi$.
This means we have the diagram:

\begin{equation*}
\begin{tikzcd}[row sep = 1.7em]
	X_\Total
		\arrow[rrrrdd, bend left = 10pt,"\psi_\Total"{near start}]
		\arrow[ddrr,dashed]
	&
	&
	&
	&\\
	{}
	&
	&
	&
	&\\
	X_\Wobs
		\arrow[ddrr,dashed]
		\arrow[rrrrdd, bend left = 10pt,"\psi_\Wobs"{near start}]
		\arrow[uu]
		\arrow[dddr,bend right =10pt,"\phi_\Wobs"{swap}]
	&
	& (M \decorate*[_{f}]{\times}{_{h}} P)_\Total
		\arrow[rr]
		\arrow[from = dd, crossing over]
		\arrow[dl,crossing over]
	&
	& P_\Total
		\arrow[dl,"h_\Total"{swap}]\\
	{}
	&M_\Total
		\arrow[from = uuul,bend right =10pt,crossing over,"\phi_\Total"]
		\arrow[rr,crossing over,"f_\Total"{near start}]
	&
	& N_\Total
	&\\
	{}
	&
	& (M \decorate*[_{f}]{\times}{_{h}} P)_\Wobs
		\arrow[rr]
		\arrow[dl]
	&
	& P_\Wobs 
		\arrow[uu]
		\arrow[dl,"h_\Wobs"{swap}]\\
	{}
	&M_\Wobs
		\arrow[uu]
		\arrow[rr,"f_\Wobs"{near start}]
	&
	& N_\Wobs
		\arrow[uu,crossing over]
	&
\end{tikzcd}
\end{equation*}
Since $(M \decorate*[_{f}]{\times}{_{h}} P)_\Total$ and
$(M \decorate*[_{f}]{\times}{_{h}} P)_\Wobs$
are pullbacks of sets, there exist unique 
$\kappa_\Total \colon X_\Total \to (M\decorate*[_{f}]{\times}{_{h}} P)_\Total$
and 
$\kappa_\Wobs \colon X_\Wobs \to (M\decorate*[_{f}]{\times}{_{h}} P)_\Wobs$
making the $\TOTAL$- and $\WOBS$-planes in the above diagram commute.
Again, by the universal property of $(M \decorate*[_{f}]{\times}{_{h}} P)_\Total$, we see that
$\iota_{\decorate*[_{f}]{\times}{_{h}}} \circ \kappa_\Wobs = \kappa_\Total \circ \iota_X$.
It remains to show that $\kappa_\Wobs$ is compatible with the equivalence relations.
For this consider $x_1 \sim_X x_2$.
Since $\phi$ and $\psi$ are constraint maps, we have
$\phi_\Wobs(x_1) \sim_M \phi_\Wobs(x_2)$
and
$\psi_\Wobs(x_1) \sim_P \psi_\Wobs(x_2)$.
Then by \eqref{eq:EquivalenceRelPullbackConSet} we get
$\kappa_\Wobs(x_1) \sim_{\decorate*[_{f}]{\times}{_{h}}} \kappa_\Wobs(x_2)$.
Thus $\kappa$ is a constraint morphism.

Parts \ref{prop:CoLimitsCSet_final}, \ref{prop:CoLimitsCSet_product}, \ref{prop:CoLimitsCSet_coproduct}
and \ref{prop:CoLimitsCSet_coequalizer} follow analogously.
The category $\ConSet$ has all finite limits and colimits for general reasons, since it has pullbacks and a terminal object as well as coequalizers, coproducts and an initial object, see \cite[Prop. 2.8.2]{borceux:1994a}.
\end{proof}

\begin{remark}
In \cite{kinoshita.power:2014a} the closely related category $\Equiv$ 
of sets equipped with an equivalence relation is examined.
Many of our results can be derived by understanding $\ConSet$ as a 
comma category of $\Sets$ and $\Equiv$.
\end{remark}

Up to now constraint sets seem to be very well behaved.
Indeed all the above constructions can be understood as combining 
easy constructions of sets and equivalence classes.
Therefore one might expect $\ConSet$ to resemble the category 
$\Sets$, but this is only partially true as the following 
characterization of 
(regular) monomorphisms and epimorphisms shows,
see
\hyperref[chap:CategoryTheory]{Appendix~\ref{chap:CategoryTheory}}
or the abstract definitions.
\begin{proposition}[Mono- and epimorphisms in $\ConSet$]
	\label{prop:MonoEpiCSet}
Let $f \colon M \to N$ be a morphism of constraint sets.
\begin{propositionlist}
	\item \label{prop:MonoEpiCSet_1}
	\index{monomorphism!constraint sets}
	The morphism $f$ is a monomorphism if and only if
	$f_\Total$ and $f_\Wobs$ are injective.
	\item \label{prop:MonoEpiCSet_2}
	\index{epimorphism!constraint sets}
	The morphism $f$ is an epimorphism if and only if
	$f_\Total$ and $f_\Wobs$ are surjective.
	\item \label{prop:MonoEpiCSet_3}
	\index{monomorphism!regular!constraint sets}
	The morphism $f$ is a regular monomorphism if and only if
	$f_\Total$ and $f_\Wobs$ are injective and 
	$(f_\Wobs)^*\left(\sim_N\right) = \mathbin{\sim_M}$.
	\item \label{prop:MonoEpiCSet_4}
	\index{epimorphism!regular!constraint sets}
	The morphism $f$ is a regular epimorphism if and only 
	if $f_\Total$ and $f_\Wobs$ are surjective and\linebreak
	$\mathbin{\sim_N} = (f_\Wobs)_\ast(\sim_M)$.
\end{propositionlist}
\end{proposition}

\begin{proof}
We only show \ref{prop:MonoEpiCSet_1} and
\ref{prop:MonoEpiCSet_3},
the statements for epimorphisms follow analogously.

Let $g_1,g_2 \colon X \to M$ with $f \circ g_1 = f \circ g_2$ be given
and assume that $f_\Total$ and $f_\Wobs$ are injective.
Then it follows from $f_\Total \circ (g_1)_\Total = f_\Total \circ (g_2)_\Total$
and $f_\Wobs \circ (g_2)_\Wobs = f_\Wobs \circ (g_2)_\Wobs$
that $(g_1)_\Total = (g_2)_\Total$
and $(g_1)_\Wobs = (g_2)_\Wobs$
hold, and thus $g_1 = g_2$ follows.
For the other implication suppose that $f$ is a monomorphism.
Let now $g_1, g_2 \colon X' \to M_\Total$ be given with $f_\Total \circ g_1 = f_\Total \circ g_2$.
Define
\begin{equation*}
	U \coloneqq \left\{ (m_1,m_2,x) \in M_\Wobs \times M_\Wobs \times X' \mid
	g_1(x) = \iota_M(m_1), g_2(x) = \iota_M(m_2) \text{ and } f_\Wobs(m_1) = f_\Wobs(m_2) \right\}.
\end{equation*}
Then $X = (X', U,\sim_\discrete)$ with $\iota_X = \pr_3$ 
is a constraint set.
Moreover, $(g_1, \pr_1) \colon X \to M$ and $(g_2,\pr_2) \colon X \to M$ are 
constraint morphisms with $f \circ (g_1,\pr_1) = f \circ (g_2,\pr_2)$.
Since $f$ is a monomorphism by assumption, it follows $g_1 = g_2$, and 
thus $f_\Total$ is injective.
To show that $f_\Wobs$ is injective let
$g_1, g_2 \colon X' \to M_\Wobs$ with $f_\Wobs \circ g_1 = f_\Wobs \circ g_2$ be given.
Then $X = (X',X',\sim_\discrete)$ with $\iota_X = \id_{X'}$
is a constraint set.
Moreover, $(\iota_M \circ g_1,g_1) \colon X \to M$
and $(\iota_M \circ g_2,g_2) \colon X \to M$ are constraint morphisms
with $f \circ (\iota_M \circ g_1,g_1) = f \circ (\iota_M \circ g_2,g_2)$.
Since $f$ is a monomorphism it follows that $g_1 = g_2$ and hence
$f_\Wobs$ is injective.
This shows the first part.

For the second part, recall that a regular monomorphism
is the equalizer of some pair of parallel morphisms.
Suppose $f$ is a regular monomorphism, then it is a monomorphism
by a general result from category theory, see \cite[Prop. 2.4.3]{borceux:1994a}.
Moreover, there exist
$h_1,h_2 \colon N \to Y$ such that
$M = \Eq(h_1,h_2)$ and $f = i$,
with $i$ as in \autoref{prop:CoLimitsCSet} \ref{prop:CoLimitsCSet_equalizer}.
Then $\sim_M$ is just the restriction of $\sim_N$.
In other words, $f_\Wobs^*(\sim_N) = \mathbin{\sim_M}$.
For the reverse implication assume that $f_\Total$ and $f_\Wobs$
are injective and $f_\Wobs^*(\sim_N) = \mathbin{\sim_M}$.
In $\Sets$ every injective function can be written as an equalizer of 
its characteristic function and the function that is constant $1$.
In our situation this means
$f_\Total = \Eq(1,\xi_{M_\Total})$
and
$f_\Wobs = \Eq(1, \chi_{M_\Wobs})$
with $\chi_\Total \colon N_\Total \to \{0,1\}$
the characteristic function of the subset $M_\Total \subseteq N_\Total$
and  $\chi_\Wobs \colon N_\Wobs \to \{0,1\}$
the characteristic function of $M_\Wobs \subseteq N_\Wobs$.
Then the combined characteristic functions
$(\chi_\Total,\chi_\Wobs) , (1,1) \colon N \to (\{0,1\},\{0,1\},\sim_\discrete)$
form constraint morphisms.
By \autoref{prop:CoLimitsCSet} \ref{prop:CoLimitsCSet_equalizer},
we have $f = \Eq((1,1),(\chi_\Total,\chi_\Wobs))$.	
\end{proof}

\begin{remark}
In a general category there exist many variations of monos (epis),
e.g. extremal, strong, strict, effective.
In $\Sets$ all these notions agree, while in $\ConSet$ they form two 
classes.
Since it can be shown that all other notions of monos (epis) are 
equivalent to either regular or plain monos (epis), we only need to 
consider these two.
\end{remark}

The fact that not every monomorphism of constraint sets is regular 
will have far reaching consequences for all further investigations.
A first noteworthy consequence is that a morphism which is mono and 
epi need not be an isomorphism:

\begin{example}
	Consider the constraint sets $M$ and $N$ with
	$M_\Total = M_\Wobs = N_\Total = N_\Wobs = \{1,2\}$
	and $\sim_M$ the discrete and $\sim_N$ the trivial equivalence 
	relation.
	Then $f = (\id, \id) \colon M \to N$ is monomorphism and 
	epimorphism of constraint sets.
	But it is not an isomorphism, as
	$\id \colon N_\Wobs \to M_\Wobs$ does not preserve the 
	equivalence relation.
\end{example}

This example directly shows that $\ConSet$ is not balanced,
meaning that a morphism which is mono and epi is not necessarily an isomorphism,
and thus, 
in contrast to $\Sets$, cannot be a topos,
see \cite{johnstone:2014a} for details on topoi.
Nevertheless, constraint isomorphisms can be characterized using
regular mono- and epimorphisms:

\begin{lemma}
	\label{lem:ConIsos}
Let $f \colon M \to N$ be a morphism of constraint sets.
The following statements are equivalent:
\begin{lemmalist}
	\item\label{lem:ConIsos_1}
	The constraint morphism $f$ is an isomorphism.
	\item\label{lem:ConIsos_2}
	The constraint morphism $f$ is a monomorphism and a regular epimorphism.
	\item\label{lem:ConIsos_3}
	The constraint morphism $f$ is a regular monomorphism and an epimorphism.
\end{lemmalist}
\end{lemma}

\begin{proof}
Suppose $f$ is an isomorphism, then there exists an inverse constraint morphism\linebreak
$f^{-1} \colon N \to M$.
Thus $f_\Total$ and $f_\Wobs$ are invertible and hence surjective and injective.
Now suppose $(x,x') \in f_\Wobs^*(\sim_N)$.
Then by definition $f_\Wobs(x) \sim_N f_\Wobs(x')$
and applying $f^{-1}_\Wobs$ yields $x \sim_M x'$.
Hence $f_\Wobs^*(\sim_N) = \mathbin{\sim_M}$, and thus $f$ is a regular monomorphism.
Suppose $(y,y') \in \mathbin{\sim_M}$, then $f^{-1}_\Wobs(y) \sim_N f^{-1}_\Wobs(y')$.
Applying $f_\Wobs$ shows $(y,y') \in (f_\Wobs)_*(\sim_M)$
and thus $(f_\Wobs)_*(\sim_M) = \mathbin{\sim_N}$.
Hence $f$ is also a regular epimorphism.
This shows \ref{lem:ConIsos_1}$\implies$\ref{lem:ConIsos_2}
and \ref{lem:ConIsos_1}$\implies$\ref{lem:ConIsos_3}.

Suppose \ref{lem:ConIsos_2}.
By definition $f_\Total$ and $f_\Wobs$ are isomorphisms.
It only remains to show that $f^{-1}_\Wobs$ is compatible with the equivalence relations.
For this let $y,y' \in N_\Wobs$ with $y \sim_N y'$ be given.
Since $f_\Wobs$ is an isomorphism there exist unique $x,x' \in M_\Wobs$
such that $f_\Wobs(x) = y$ and $f_\Wobs(x') = y'$.
Moreover, since $f$ is a regular monomorphism we know that $f_\Wobs^*(\sim_N) = \mathbin{\sim_M}$, meaning that $x \sim_M x'$.
Hence, $f^{-1}$ is a constraint morphism, and therefore $f$ is a constraint isomorphism.

The implication \ref{lem:ConIsos_3}$\implies$\ref{lem:ConIsos_1} follows analogously.
\end{proof}

Related to this mismatch of regular and plain monomorphisms is the definition of a subset of a constraint 
set.
We could either define a subset as an equivalence class of 
monomorphisms or of regular monomorphisms.
It is common to choose regular monomorphisms in such a situation and 
we will follow this strategy.

\begin{definition}[Constraint subset]
	\label{def:ConstraintSubset}
	\index{subset!constraint set}
	A \emph{constraint subset} of a constraint set $M$ consists of 
	subsets
	$U_\Total \subseteq M_\Total$ and $U_\Wobs \subseteq M_\Wobs$
	such that $\iota_M(U_\Wobs) \subseteq U_\Total$.
\end{definition}

Every constraint subset $U$ of $M$ defines a constraint set
$U = (U_\Total, U_\Wobs, \sim_M\!\!\at[\small]{U_\Wobs})$.
The obvious inclusion $i \colon U \to M$ is a regular monomorphism.
With this we can now define the image and preimage of a constraint 
morphism.

\begin{definition}[Image and preimage]
	\label{def:ImagePreimageCSet}
	Let $f \colon M \to N$ be a morphism of constraint sets.
	\begin{definitionlist}
		\item \index{preimage of constraint subset}
		Let $U \subset N$ be a constraint subset with inclusion
		$i \colon U \to N$.
		The \emph{preimage} of $U$ along $f$ is defined by
		\begin{equation}
			f^{-1}(U) \coloneqq M \decorate*[_f]{\times}{_i} U.
		\end{equation}
		More explicitly, we have
		\begin{equation}
			f^{-1}(U) 
			= \Big(f_\Total^{-1}(U_\Total),\,\,
			f_\Wobs^{-1}(U_\Wobs),\,\,
			\sim_M\at{f_\Wobs^{-1}(U_\Wobs)} \Big)
		\end{equation}
		\item \index{image!constraint sets}
		The \emph{image} of $f$ is defined by
		\glsadd{image}
		\begin{equation}
			\image(f) \coloneqq \big(\image(f_\Total),\,\, \image(f_\Wobs),\,\, \sim_{\image}\big)
		\end{equation}
		with $\sim_{\image} = (f_\Wobs)_\ast(\sim_M)$.
		\item \index{image!regular!constraint sets}
		The \emph{regular image} of $f$ is defined by
		\glsadd{regimage}
		\begin{equation}
			\regimage(f) \coloneqq \big(\image(f_\Total),\,\, \image(f_\Wobs),\,\, \sim_{\regimage}\big)
		\end{equation}
		with $f(x_1) \sim_{\regimage} f(x_2)$ if and only if
		$f(x_1) \sim_N f(x_2)$.
	\end{definitionlist}
\end{definition}

\begin{example}
	\label{ex:RegularImageVSImage}
Image and regular image of a constraint morphism do not agree in general.
To see this let
$M = (\{1,2\}, \{1,2\}, \sim_\discrete)$ and
$N = (\{1,2\},\{1,2\},\sim_N)$ with $1 \sim_N 2$
be given
and consider the constraint morphism
$f = (\id_{\{1,2\}}, \id_{\{1,2\}}) \colon M \to N$.
Then $\image(f) = M$ while $\regimage(f) = N$.
\end{example}

Using the image we can factorize every constraint morphism as a 
regular epimorphism followed by a monomorphism, while the regular 
image yields a factorization as an epimorphism followed by a regular 
monomorphism.
We will mainly use the regular image, since using our definition of 
constraint subset it is in fact a constraint subset of the codomain, 
while the image is not.

Let us now turn our attention to the set of all constraint morphisms 
between constraint sets.
This set can actually be upgraded to a constraint set itself.

\begin{proposition}[Closed monoidal structure on $\ConSet$]
	\label{def:ClosedMonoidalStructureCSet}
Let $M$ and $N$ be constraint sets.
\begin{propositionlist}
	\item Setting
	\glsadd{ConMap}
	\begin{equation}
	\begin{split}
		\ConMap(M,N)_\Total &\coloneqq \Map(M_\Total,N_\Wobs),\\
		\ConMap(M,N)_\Wobs &\coloneqq \Map(M,N),
	\end{split}
	\end{equation}
	together with the inclusion
	$\iota \colon \Map(M,N) \to \Map(M_\Total,N_\Total)$
	given by $\iota((f_\Total,f_\Wobs)) = f_\Total$
	and the equivalence relation on $\ConMap(M,N)$ given by
	\begin{equation}
		f \sim g \coloniff \forall x \in M_\Wobs : f(x) \sim_N g(x),
	\end{equation}
	defines a constraint set $\ConMap(M,N)$.
	\item The functor
	$\ConMap(M, \argument) \colon \ConSet \to \ConSet$
	is right adjoint to the functor
	$\argument \times M \colon \ConSet \to \ConSet$, i.e.
	$\ConSet$ is a cartesian closed category.
\end{propositionlist}
\end{proposition}

\begin{proof}
The first part is a simple check.
For the second part recall the definition of adjoint functors from
\autoref{def:AdjointFunctors}.
Fix $X \in \ConSet$
and define functors
$\functor{F} = \argument \times X$
and
$\functor{G} = \ConMap(X,\argument)$.
Then
$\ev(M) \colon \ConMap(X,M) \times X \to M$ given by
\index{evaluation}
\glsadd{evaluation}
\begin{align*}
	\ev(M)_\Total \colon \Map(X_\Total,M_\Wobs) \times X_\Total \ni (f,x) &\mapsto f(x) \in M_\Total, \\
	\ev(M)_\Wobs \colon \Map(X,M) \times X_\Wobs \ni (f,x) &\mapsto f_\Wobs(x) \in M_\Wobs
\end{align*}
is a constraint map.
Similarly,
$\coev(M) \colon M \to \ConMap(X,M \times X)$
defined by
\index{coevaluation}
\glsadd{coevaluation}
\begin{align*}
	\coev(M)_\Total \colon M_\Total \ni m &\mapsto (x \mapsto (m,x)) \in \Map(X_\Total, M_\Total \times X_\Total), \\
	\coev(M)_\Wobs \colon M_\Wobs \ni m &\mapsto \big(x \mapsto (x,\iota_M(m)), x \mapsto (m,x)\big) \in \Map(X, M \times X)
\end{align*}
is a constraint map.
In particular, these are compatible with the equivalence relations.
They define natural transformations, since for a morphism $f \colon M \to N$
the diagrams 
\begin{equation*}
\begin{tikzcd}
	\ConMap(X,M) \times X
		\arrow[r,"\ev(M)"]
		\arrow[d,"{\ConMap(X,f) \times X}"{swap}]
	& M
		\arrow[d,"f"] \\
	\ConMap(X,N) \times X
		\arrow[r,"\ev(N)"]
	&N
\end{tikzcd}
\quad\text{and}\quad
\begin{tikzcd}
	M
		\arrow[r,"\coev(M)"]
		\arrow[d,"f"{swap}]
	& \ConMap(X,M \times X)
		\arrow[d,"{\ConMap(X,f \times X)}"] \\
	N
		\arrow[r,"\coev(N)"]
	& \ConMap(X,N\times X)
\end{tikzcd}
\end{equation*}
commute.
It remains to check that
\begin{equation*}
	\id_{M\times X} = \ev(M \times X) \circ \functor{F}(\coev(M))
	\quad\text{and}\quad 
	\id_{\ConMap(X,M)} = \functor{G}(\ev(M)) \circ \coev(\ConMap(X,M))
\end{equation*}
hold.
We need to check this separately on the $\TOTAL$- and $\WOBS$-component.
Thus let $(m,x) \in M_\Total \times X_\Total$ and $f \colon X_\Total \to M_\Total$ be given.
Then
\begin{align*}
	\ev(M\times X)_\Total\big( \functor{F}(\coev(M)_\Total)\big)(m,x)
	= \ev(M\times X)_\Total\big(\coev(M)_\Total(m),x\big)
	= (m,x) \\
	\shortintertext{and}
	\Big(\functor{G}\big(\ev(M)_\Total\big)\big(\coev\big(\ConMap(X,M)\big)_\Total\big)(f)\Big)(x)
	= \big(x' \mapsto \ev(M)_\Total(f,x')\big)(x)
	= f(x).
\end{align*}
The same computations hold for the $\WOBS$-component, which finally shows that 
$\functor{G}$ is indeed right adjoint to $\functor{F}$
and we obtain a cartesian closed category.
\end{proof}

Before we turn our attention to a more special class of constraint sets let us
investigate more closely the relationship of constraint sets and classical sets.
We have obvious forgetful functors
\begin{align} \label{eq:ForgetfulToTotalConSets}
	\functor{U}_\Total &\colon \ConSet \to \Sets,\qquad M \mapsto M_\Total \\
	\shortintertext{ and }
	\functor{U}_\Wobs &\colon \ConSet \to \Sets,\qquad M \mapsto M_\Wobs
\end{align}
forgetting everything but the indicated components.
When looking at \autoref{prop:CoLimitsCSet} it becomes clear that the $\TOTAL$-components of definitions and constructions
internal to $\ConSet$ will just be the classical definitions and constructions for the $\TOTAL$-components.
We summarize this:

\begin{lemma}
	The forgetful functor $\functor{U}_\Total \colon \ConSet \to \Sets$
	is cartesian closed and preserves finite limits and colimits.
\end{lemma}

It will be a recurring theme for all our constraint definitions, constructions and theorems that
their $\TOTAL$-components will recover their classical analogues.

The forgetful functor $\functor{U}_\Total$ has an obvious left adjoint given by
$\functor{F}_\Total(M) \coloneqq (M,M,\sim_\discrete)$.
Thus we can also understand $\Sets$ as the full subcategory of $\ConSet$
consisting of constraint sets $M$ with $M_\Total = M_\Wobs$
and $\mathbin{\sim_M} = \mathbin{\sim_\discrete}$.

\renewcommand{\thesubsubsection}{\thesection.\arabic{subsubsection}}
\subsubsection{Embedded Constraint Sets}

Most examples of constraint sets as they appear in \autoref{chap:ConstraintGeometricStructures} will 
exhibit $M_\Wobs$ as a subset of $M_\Total$.

\begin{definition}[Embedded constraint set]\
	\label{def:InjConstraintSet}
\begin{definitionlist}
	\item \index{embedded constraint!set}
	\index{constraint!embedded|see {embedded constraint}}
	A constraint set $M$ with injective $\iota_M$ is called 
	an \emph{embedded constraint set}.
	\item The full subcategory of $\ConSet$ consisting of
	embedded constraint sets is denoted by
	\glsadd{injConSet}$\injConSet$.
\end{definitionlist}
\end{definition}

Note that for a morphism $f = (f_\Total,f_\Wobs)$ of embedded 
constraint sets the map $f_\Wobs$ is completely determined by 
$f_\Total$.
Hence we will often identify $f$ with $f_\Total$.
Then $f_\Wobs$ is just the restriction of $f$
to the $\WOBS$-component.

\begin{proposition}[The category $\injConSet$]\
	\label{prop:CategoryInjCSet}
	\begin{propositionlist}
		\item The subcategory $\injConSet$ of $\ConSet$ is closed 
		under 
		finite limits and has all finite colimits.
		\item The subcategory $\injConSet$ is an exponential ideal in $\ConSet$, this means
		for all $X \in \ConSet$ and $M \in \injConSet$ we have
		$\ConMap(X,M) \in \injConSet$.
		\item The category $\injConSet$ is cartesian closed.
	\end{propositionlist}
\end{proposition}

\begin{proof}
We show that $\injConSet$ is a reflective subcategory of 
$\ConSet$.
Denote by $\functor{I} \colon \injConSet \to \ConSet$
the inclusion.
Mapping a constraint set
$M = (M_\Total,M_\Wobs,\sim_M)$
to
\begin{equation*}
	M^\inj \coloneqq {(M_\Total, \iota_M(M_\Wobs), (\iota_M)_*\sim_M)}
\end{equation*}
defines a functor
$\argument^\inj \colon \ConSet \to \injConSet$,
with $(\iota_M)_*\sim_M$
the induced equivalence relation on the image
of $\iota_M$.
The functor $\argument^\inj$ is left adjoint to $\functor{I}$,
thus $\injConSet$ is a reflective subcategory of $\ConSet$ and 
hence is closed under finite limits and has all finite colimits,
see \cite[Sec. 3.5]{borceux:1994a}.
For the second part note that by \cite[Prop. 4.3.1]{johnstone:2002a}
it would be enough to show that $\argument^\inj$ preserves finite 
products.
But let us show this more directly:
Let $f, g \in \ConMap(X,M)$ be given with $f_\Total = g_\Total$.
Diagrammatically we have:
\begin{equation*}
\begin{tikzcd}
	X_\Total 
		\arrow[r,"f_\Total",shift left=3pt]
		\arrow[r,"g_\Total"{swap},shift right=3pt]
	&M_\Total \\
	X_\Wobs
		\arrow[u,"\iota_X"]
		\arrow[r,"f_\Wobs",shift left=3pt]
		\arrow[r,"g_\Wobs"{swap},shift right=3pt]
	& M_\Wobs
		\arrow[u,"\iota_M"{swap},hookrightarrow]
\end{tikzcd}
\end{equation*}
Since $f_\Total = g_\Total$ we have $\iota_M \circ f_\Wobs = \iota_M g_\Wobs$
and thus by the injectivity of $\iota_M$ we obtain
$f_\Wobs = g_\Wobs$.
Thus $\iota \colon \Map(X,M) \to \Map(X_\Total,M_\Total)$
as defined in \autoref{def:ClosedMonoidalStructureCSet} is injective.
This shows the second part.
The third part is now a direct consequence of the second.
\end{proof}

At this point it seems that we could restrict ourselves to the 
category $\injConSet$ since all categorical constructions exist in 
this category.
However, note that even though colimits exist in $\injConSet$ they do not 
necessarily agree with the respective colimits in the surrounding 
category $\ConSet$, as the next example illustrates:

\begin{example} \label{ex:injConSetIsNotClosedUnderColimits}
Consider two embedded constraint sets $M$ and $N$ given by
$M \coloneqq (\{\pt\}, \emptyset, \sim)$
and $N \coloneqq (\{0,1\},\{0,1\},\sim_\discrete)$
together with the constraint maps
$f \equiv 0$ and $g \equiv 1$
from $M$ to $N$.
Their coequalizer is then given by
$\Coeq(f,g) = (\{0\},\{0,1\},\sim_\discrete)$,
which is obviously not embedded.
\end{example}

This will have consequences for the reduction of (embedded) constraint sets
as we will shortly discuss.

The subcategory $\injConSet$ of embedded constraint sets can also
be characterized by general categorical terms as the subcategory
of regular projective objects:

\begin{proposition}
	\label{prop:ProjectiveConSets}%
Let $P \in \ConSet$ be a constraint set.
Then the following statements are equivalent:
\begin{propositionlist}
	\item\label{prop:ProjectiveConSets_1}
	Every regular epimorphism $M \to P$ splits.
	\item\label{prop:ProjectiveConSets_2}
	\index{projective!constraint set}
	$P$ is a regular projective object in $\ConSet$, i.e.
	for every regular epimorphism $f \colon M \to N$ and every morphism
	$g \colon P \to N$ there exists a morphism
	$h \colon P \to M$ such that $f \circ h = g$.
	\item\label{prop:ProjectiveConSets_3}
	We have $P \in \injConSet$.
\end{propositionlist}
\end{proposition}

\begin{proof}
Assume \ref{prop:ProjectiveConSets_1}.
Given $f$ and $g$ as in \ref{prop:ProjectiveConSets_2}
consider the pullback
$P \deco{}{g}{\times}{}{f} M$.
It is easy to see that
$\pr_1 \colon P \deco{}{g}{\times}{}{f} M \to P$
is a regular epimorphism.
By assumption $\pr_1$ splits, i.e. there exists
$i \colon P \to P \deco{}{g}{\times}{}{f} M$
such that
$\mathord{\pr_1} \circ \mathop{i} = \id_P$.
Then $\chi = \pr_2 \circ \mathop{i}$ gives the desired morphism.
Conversely,  choosing
$g = \id_P$ in \ref{prop:ProjectiveConSets_2} directly yields \ref{prop:ProjectiveConSets_1}.
Now assume again \ref{prop:ProjectiveConSets_2}.
We want to show that
$\iota_P \colon P_\Wobs \to P_\Total$
is injective.
For this consider
$M_\Total \coloneqq P_\Total \times P_\Wobs$,
$M_\Wobs \coloneqq P_\Wobs$
and $\iota_M \coloneqq \iota_P \times \id_{P_\Wobs}$.
Then
$f = (\pr_1, \id_{p_\Wobs}) \colon M \to P$
is a regular epimorphism and hence splits by assumption.
Therefore, there exists
$h \colon P \to M$ with $f \circ h = \id_P$.
Thus $h_\Total \circ \iota_P = \iota_M \circ h_\Wobs$
is injective since $\iota_M$ and $h_\Wobs$ are injective.
Then $\iota_P$ must also be injective.
Finally, assume \ref{prop:ProjectiveConSets_3} and let
$f \colon M \to P$ be a regular epimorphism.
It follows that $f\at{M_\Null} \colon M_\Null \to P_\Null$ is surjective and thus
there exists a splitting $h \colon P_\Null \to M_\Null$.
Now we can extend $h$ successively to $P_\Wobs$ and
$P_\Total$, obtaining a splitting of $f$.
\end{proof}

\subsubsection{Reduction of Constraint Sets}
Constraint sets were introduced in order to formalize the set 
theoretic information underlying geometric reduction principles.
Thus they are defined in such a way to allow for a reduction 
procedure already on this set theoretic level.

\begin{definition}[Reduction functor]
	\label{def:ReductionFunctorCSet}
	\index{reduction!sets}
The functor $\red \colon \ConSet \to \Sets$ given by mapping a 
constraint set $M$ to
$M_\red \coloneqq M_\Wobs / \mathord{\sim}_M$ and a constraint
morphism $f \colon M \to N$ to the induced morphism
$f_\red \colon M_\red \to N_\red$
is called \emph{reduction functor}.
\end{definition}

This reduction procedure can now be shown to be compatible with the 
various constructions from \autoref{prop:CoLimitsCSet}:

\begin{proposition}[Properties of reduction]\
	\label{prop:ReductionPropsCSet}
\begin{propositionlist}
\item The functor $\red \colon \ConSet \to \Sets$ preserves all finite
limits and colimits.
\item The functor $\red \colon \ConSet \to \Sets$ is cartesian closed.
\end{propositionlist}
\end{proposition}

\begin{proof}
A straightforward computation shows that $\red$ preserves the 
final object and pullbacks, and thus preserves all finite limits.
Moreover, it preserves coproducts as well as coequalizer, and hence 
preserves all finite colimits.
For the second part note that since $\red$ preserves products it 
is a cartesian functor.
Moreover, since the final object $1$ is the unit of the monoidal 
structure, the first part shows that $\red$ preserves this unit.
Finally, we have a canonical injection
\begin{equation*}
	\ConMap(M,N)_\red \hookrightarrow \Map(M_\red,N_\red),
\end{equation*}
which is also surjective, since using the axiom of choice any 
morphism in $\Map(M_\red,N_\red)$ can be lifted to a morphism in
$\Map(M_\Wobs,N_\Wobs)$ compatible with the equivalence 
relations, and then be extended to $M_\Total$.
\end{proof}

This result shows that $\ConSet$ is the correct category for studying 
constructions compatible with reduction.

\begin{remark}\
	\label{rem:ReductionConSets}
\begin{remarklist}
	\item The reduction functor can be understood as first forgetting the $\TOTAL$-component
	and then computing a coequalizer in the resulting category, whose objects have been called
	coisotropic pairs in \cite{dippell.esposito.waldmann:2019a}.
	Since taking colimits commutes with colimits, whenever the forgetful functor commutes with
	colimits, so does the whole reduction.
	This point of view could lead to a more general theory for
	the relation of reduction to co/limits.
	\item \label{rem:RedcutionOfEquations}
	Let $f,g \colon M \to N$ be maps between sets.
	Then their pullback is given by the subset $\{x \in M \mid f(x) = g(x)\} \subseteq M$.
	Thus pullbacks can be understood as describing subsets of elements fulfilling a given equation.
	Since the reduction of constraint sets commutes with limits, pullbacks reduce to pullbacks.
	In other words, elements of a constraint set fulfilling the equation $f(x) = g(x)$ will reduce
	to elements satisfying the reduced equation $f_\red([x]) = g_\red([x])$.
	However it is important to note that the functor $\red$ does not \emph{reflect} limits, meaning that
	even if the reduced equation $f_\red([x]) = g_\red([x])$ is fulfilled we can \emph{not}
	infer that also $f(x) = g(x)$ must hold.
	\item \label{rem:NonEmbeddedConSetsAreImportant}
	By contrast, the reduction of embedded constraint sets, which is given by the
	composition $\red \circ \functor{I}$ of the inclusion
	$\functor{I} \colon \injConSet \to \ConSet$ with the above reduction functor,
	may not preserve colimits, since $\functor{I}$ does not, as shown in \autoref{ex:injConSetIsNotClosedUnderColimits}.
	Thus even if we are mainly interested in examples which yield embedded constraint sets,
	the moment we construct colimits we are forced to work in the bigger category
	$\ConSet$ if we want our construction to stay compatible with reduction.
\end{remarklist}
\end{remark}

\subsubsection{Strong Constraint Sets}

Another special type of constraint sets appears in the setting of 
Hamiltonian actions of Lie groups $\group{G}$ on a symplectic or 
Poisson manifold $M$.
 In this case the coisotropic submanifold is given by the zero level set
 $C$ of the momentum map, but the equivalence relation on
 $C$ can be viewed as the restriction of the orbit relation on $M$.
 In this situation the underlying constraint set carries an additional
 equivalence relation on the $\TOTAL$-component.

\begin{definition}[Strong constraint set]\
	\label{def:StrConstraintSet}
	\begin{definitionlist}
		\item \index{strong constraint!set}
		\index{constraint!strong|see{strong constraint}}
		A constraint set $M$ together with an 
		equivalence relation $\sim_M^\Total$ on $M_\Total$
		\index{saturated subset}
		such that $\image(\iota_M)$ is saturated,
		i.e. from $\iota_M(x) \sim_M^\Total y$ follows $y \in 
		\image(\iota_M)$ for all $x \in M_\Wobs$, $y \in M_\Total$,
		and $\sim_M^\Total$ restricts to $(\iota_M)_\ast(\sim_M)$ on 
		$\image(\iota_M)$,
		is called a \emph{strong constraint set}.
		\item A morphism $f \colon M \to N$ of strong constraint sets
		(or \emph{constraint morphism})
		is a morphism of constraint sets with $f_\Total$ preserving 
		the equivalence relation, i.e.
		$\sim_M^\Total \subseteq f_\Total^*(\sim_N^\Total)$.
		\item The category of strong constraint sets and their 
		morphisms is denoted by
		\glsadd{strConSet}$\strConSet$.
		The category of embedded strong constraint sets, i.e. those with injective
		$\iota_M \colon M_\Wobs \to M_\Total$,
		is denoted by
		\glsadd{strConSet}$\injstrConSet$.
	\end{definitionlist}
\end{definition}

Observe that embedded strong constraint sets are just given by a subset
$M_\Wobs \subseteq M_\Total$ together with an equivalence relation
$\sim_\Total$ on $M_\Total$ such that $M_\Wobs$ is saturated with respect
to $\sim_\Total$.
Moreover, morphisms of strong constraint sets are again completely determined
by their $\TOTAL$-components.

Even though strong constraint sets will appear as the structure underlying 
many objects of interest (in particular the functions on constraint manifolds, see \autoref{prop:FunctionsOnConManifolds}),
we will not investigate them in full detail.
This is justified by the fact that in geometric situations we will be confronted only with embedded
strong constraint sets, and in algebraic situations it is easier to work with subobjects instead of equivalence relations.

The next proposition clarifies the relation between constraint and 
strong constraint sets.

\begin{proposition}[The category $\strConSet$]\
	\label{prop:CategoryStrCSet}
\begin{propositionlist}
	\item \label{prop:CategoryStrCSet_1}
	Forgetting the equivalence relation on the
	$\TOTAL$-component yields a functor
	$\functor{U} \colon \strConSet \to \ConSet$.
	\item \label{prop:CategoryStrCSet_2}
	\index{strong hull!constraint set}
	The functor
	$\functor{U} \colon \strConSet \to \ConSet$
	has a left adjoint
	$\argument^\str \colon \ConSet \to \strConSet$
	given on objects by
	$M^\str = M$,
	together with the equivalence relation $(\iota_M)_*(\sim_M)$ 
	on $M_\Total$.
	\item \label{prop:CategoryStrCSet_3}
	The category $\strConSet$ is $\ConSet$-enriched with
	\glsadd{strConMap}
	\begin{equation}
	\begin{split}
		\strConMap(M,N)_\Total
		&\coloneqq \Map(M_\Total, N_\Total), \\
		\strConMap(M,N)_\Wobs
		&\coloneqq \left\{ f \in \ConMap(\functor{U}(M),\functor{U}(N))_\Wobs \mid
		f_\Total(x) \sim_N^\Total f_\Total(y) \text{ for all } x \sim_M^\Total y \right\},
	\end{split}
	\end{equation}
	with the obvious inclusion
	$\iota \colon \strConMap(M,N)_\Wobs \to 
	\strConMap(M,N)_\Total$
	and the equivalence relation on $\strConMap(M,N)_\Wobs$ 
	given by
	\begin{equation}
		f \sim g \coloniff \forall x \in M_\Total :
		f_\Total(x) \sim_N^\Total g_\Total(x)
		\text{ and }
		\forall x \in M_\Wobs : f_\Wobs(x) \sim_N g_\Wobs(x)
	\end{equation}
	for $M,N \in \strConSet$.
	\item \label{prop:CategoryStrCSet_4}
	The functor $\functor{U}$ is $\ConSet$-enriched.
\end{propositionlist}
\end{proposition}

\begin{proof}
The first part is clear.
For the second part, choose the constraint morphisms
\begin{equation*}
	\varepsilon_M \colon \functor{U}(M)^\str \to M
	\qquad\text{ and }\qquad
	\eta_M \colon M \to \functor{U}(M^\str)
\end{equation*}
to be the identity on both
$M_\Total$ and $M_\Wobs$.
Hence they
clearly define the evaluation and coevaluation
of the adjunction.
The third part is an easy check, using the usual composition of maps
as composition in the enriched category.
The last part is then just the fact that
we have a canonical morphism
$\strConMap(M,N) \to \ConMap(\functor{U}(M), \functor{U}(N))$.
\end{proof}

It is important to note that the $\ConSet$-enrichment of 
$\strConSet$ does not agree with its internal hom with respect to the 
cartesian monoidal structure, which we have not spelled out.
The reason we do not consider the closed structure is that the 
forgetful functor $\functor{U} \colon \strConSet \to \ConSet$ is 
not closed, and hence the internal hom will not be compatible with 
reduction.
Whereas, considering the $\ConSet$-enrichment we can define a functor 
of reduction on $\strConSet$ by simply forgetting to $\ConSet$ first,
and thus obtain a co/limit-preserving reduction
\index{reduction!strong constraint sets}
\begin{equation}
	\red \colon \strConSet \to \Sets.
\end{equation}

There is another important relation between strong constraint and 
constraint sets: The $\ConSet$-enriched category $\strConSet$ is 
powered and copowered, cf. \cite[Chap. 6.5]{borceux:1994b},
meaning that morphisms into and products with a strong constraint set 
can be equipped with the structure of a strong constraint set.

\begin{proposition}[Co/Power in $\strConSet$]
	\label{prop:CoPowerCSet}
	Let $M \in \ConSet$ and $N \in \strConSet$ be given.
	\begin{propositionlist}
		\item We have $\ConMap(M,\functor{U}(N)) \in \strConSet$ with
		\begin{equation}
			f \sim^\Total g \coloniff \forall x \in M_\Total : f(x) 
			\sim_N^\Total 
			g(x)
		\end{equation}
		for $f,g \in \Map(M_\Total,N_\Total)$.
		\item We have $M \times \functor{U}(N) \in \strConSet$ with
		\begin{equation}
			(x,y) \sim^\Total (x',y') \coloniff
			\begin{cases}
				x = x' \text{ and } y \sim_N^\Total y' &\text{ if } x 
				\notin 
				\image(\iota_M) \text{ or } x' \notin \image(\iota_M) 
				\\
				x \sim_{(\iota_M)_*} x' \text{ and } y \sim_N^\Total 
				y' &\text{ else}.
			\end{cases}
		\end{equation}
	\end{propositionlist}
\end{proposition}

We will often suppress the forgetful functor $\functor{U}$ in our 
notation.
With the two previous propositions we get for a fixed strong 
constraint set $N$ functors
\begin{equation}
	\ConMap(\argument, N) \colon \ConSet^\opp \to \strConSet
\end{equation}
and
\begin{equation}
	\strConMap(\argument, N) \colon \strConSet^\opp \to \ConSet.
\end{equation}

%% file: constraint-modules.tex
After having defined the category $\ConSet$ as a replacement for 
$\Sets$ which admits a reduction procedure, we can now start to 
implement virtually all the classical mathematical objects internal 
to this category.
In this chapter we will concentrate on algebraic notions.
Thus we could proceed as follows: Since $\ConSet$ is (cartesian) 
monoidal we can construct the category of monoids internal to 
$\ConSet$, giving a notion of constraint monoids.
Requiring invertibility leads us to constraint groups.
Now considering monoids internal to the category of constraint 
abelian groups 
yields constraint rings, and additionally constraint modules over 
such.
Continuing, we obtain categories of constraint algebras, constraint 
modules over algebras etc.
Constructing all these algebraic notions in this way has the 
advantage that all resulting structures will automatically come 
equipped with a functor of reduction.

Since not all intermediate steps will be needed in this thesis we 
will only spell out those constructions important for our discussion.
In \autoref{sec:ConGroups} we introduce constraint groups and their 
actions. 
On one hand these will be the basis to define constraint 
$\field{k}$-modules in \autoref{sec:ConKModules}, on the other hand constraint groups
will feature prominently as the gauge group acting on Maurer-Cartan 
elements of constraint differential graded Lie algebras, see 
\autoref{sec:ConDeformationFunctor}.

\subsection{Constraint Groups}
\label{sec:ConGroups}

If we consider groups internal to the category $\ConSet$ of 
constraint sets we would obtain a group homomorphism
$\iota_\group{G} \colon \group{G}_\Wobs \to \group{G}_\Total$
together with an equivalence relation on $\group{G}_\Wobs$ compatible 
with the group structure.
Such equivalence relations can equivalently be given by normal 
subgroups, leading us to the following definition.

\begin{definition}[Constraint group]\
	\label{def:ConGroup}%
	\begin{definitionlist}
		\item \index{constraint!group}
		A \emph{constraint group} is given by a
		triple of groups
		$\group{G} = (\group{G}_\Total, \group{G}_\Wobs, 
		\group{G}_\Null)$,
		with $\group{G}_\Null \subseteq \group{G}_\Wobs$ 
		a normal subgroup, together with a group homomorphism
		$\iota_\group{G} \colon \group{G}_\Wobs \to \group{G}_\Total$.	
		\item A \emph{morphism
			$\Phi \colon \group{G} \to \group{H}$
			of constraint groups} $\group{G}$ and $\group{H}$ is given by a 
		pair of group homomorphisms
		$\Phi_\Total \colon \group{G}_\Total \to \group{H}_\Total$ and
		$\Phi_\Wobs \colon \group{G}_\Wobs \to \group{H}_\Wobs$ such
		that
		$\Phi_\Total \circ \iota_\group{G} = \iota_\group{H} \circ 
		\Phi_\Wobs$
		and
		$\Phi_\Wobs(\group{G}_\Null) \subseteq \group{H}_\Null$.
		\item The category of constraint groups is denoted by
		\glsadd{ConGroups}$\ConGroups$.
	\end{definitionlist}
\end{definition}

\begin{example}\
	\label{ex:ConGroups}
\begin{examplelist}
	\item \label{ex:ConGroups_1}
	\index{constraint!automorphisms}
	Let $M \in \ConSet$ be a constraint set.
	The invertible constraint endomorphisms of $M$
	define a constraint subset
	\glsadd{ConAut}$\ConAut(M) \subseteq \ConMap(M,M)$.
	They form a constraint group by considering the equivalence 
	relation on $\ConAut(M)_N$ as the normal subgroup
	\begin{equation}
		\ConAut(M)_\Null = \left\{ f \in \ConAut(M)_\Wobs \mid 
		\forall x 
		\in M_\Wobs : f(x) \sim_M x \right\}.
	\end{equation}
	\item \label{ex:ConGroups_2}
	Let
	$M \in \injConSet$ be an embedded constraint set.
	Let	furthermore $\group{G}$ be a group acting on $M_\Total$ via
	$\Phi \colon \group{G} \times M_\Total \to M_\Total$.
	Then $(\group{G}, \group{G}_{M_\Wobs}, \group{G}_\sim)$,
	with $\group{G}_{M_\Wobs}$ the stabilizer subgroup of the subset
	$M_\Wobs$ and $\group{G}_\sim$ the normal subgroup of
	$\group{G}_{M_\Wobs}$ consisting of all
	$g \in \group{G}_{M_\Wobs}$ such that $\Phi_g(x) \sim x$ for
	all $x \in M_\Wobs$, is a constraint group.
\end{examplelist}
\end{example}

Using the constraint automorphism group we could define an action of 
a group $\group{G}$ on a constraint set $M$ to be a constraint 
group morphism $\Phi \colon \group{G} \to \ConAut(M)$.
To phrase this in more elementary terms note that the equivalence 
relation on the product of two constraint groups $\group{G}$
and $\group{H}$ is given by the normal subgroup
\begin{equation}
	(\group{G} \times \group{H})_\Null 
	= \group{G}_\Null \times \group{H}_\Null.
\end{equation}

\begin{definition}[Action of constraint group]\
	\label{def:ConGroupAction}%
	\begin{definitionlist}
		\item \index{constraint!group!action}
		Let $\group{G}$ be a constraint group and $M$ a 
		constraint 
		set.
		An \emph{action} of $\group{G}$ on $M$ is given by an action
		$\Phi_\Total \colon \group{G}_\Total \times M_\Total \to 
		M_\Total$
		of $\group{G}_\Total$ on $M_\Total$
		and an action
		$\Phi_\Wobs \colon \group{G}_\Wobs \times M_\Wobs \to M_\Wobs$
		of
		$\group{G}_\Wobs$ on $M_\Wobs$
		such that
		$\iota_M \circ \Phi_\Wobs 
		= \Phi_\Total \circ (\iota_\group{G}\times \iota_M)$
		and
		$(\Phi_\Wobs)_g(x) \sim_M x$
		for all
		$g \in \group{G}_\Null$ and $x \in M_\Wobs$.
		\item Let $\group{G}$ and $\group{H}$ be constraint groups 
		acting
		on constraint sets $M$ and $N$, respectively.
		A \emph{morphism of constraint group actions}
		is given by a pair $(\phi, f)$ consisting of a constraint 
		group 
		morphism $\phi \colon \group{G} \to \group{H}$ and a morphism $f \colon M \to 
		N$ 
		of constraint sets, such that
		\begin{equation}
			f_\Total\big((\Phi^\group{G}_\Total)_g(x)\big)
			= (\Phi^\group{H}_\Total)_{\phi(g)}\big(f_\Total(x)\big)
		\end{equation}
		for all $g \in \group{G}_\Total$, $x \in M_\Total$
		and
		\begin{equation}
			f_\Wobs\big((\Phi^\group{G}_\Wobs)_g(x)\big)
			= (\Phi^\group{H}_\Wobs)_{\phi(g)}\big(f_\Wobs(x)\big)
		\end{equation}
		for all $g \in \group{G}_\Wobs$, $x \in M_\Wobs$
		holds.
		Such a map $f$ will also be called \emph{equivariant along 
		$\phi$}.
		\item The category of actions of constraint groups on 
		constraint 
		sets together with the above defined morphisms is denoted by
		\glsadd{ConGroupAct}$\ConGroupAct$.
	\end{definitionlist}
\end{definition}

Nevertheless, it is sometimes useful to think of a group action in 
terms of a morphism $\Phi \colon \group{G} \to \ConAut(M)$,
or, equivalently, as a morphism
$\Phi \colon \group{G} \times M \to M$ of constraint sets
fulfilling the usual properties of group actions in every component.
As is commonly done, we will often use $\acts$ for a generic group 
action, and sometimes even omit writing out the action entirely.

\begin{example}\
Let
$(\group{G}, \group{G}_{M_\Wobs}, \group{G}_\sim)$
be the constraint group constructed from a group action
of $\group{G}$ on $M_\Total$ as in
\autoref{ex:ConGroups} \ref{ex:ConGroups_2}.
Then $(\Phi,\Phi\at{\group{G}_{M_\Wobs}})$ clearly
gives a constraint action on $M$.
\end{example}

Next we want to consider constraint orbit spaces 
of constraint group actions.

\begin{lemma}[Constraint orbit space]
	\label{lem:ConOrbitSpacec}%
Let $M \in \ConSet$ together with an action of a constraint group 
$\group{G}$ on $M$ be given.
Then $M / \group{G}$ defined by
\begin{equation}
\begin{split}
	(M / \group{G})_\Total 
	\coloneqq M_\Total / \group{G}_\Total,\\
	(M / \group{G})_\Wobs
	\coloneqq M_\Wobs / \group{G}_\Wobs,
\end{split}
\end{equation}
together with
\begin{equation}
	\iota_{M/\group{G}} \colon (M / \group{G})_\Wobs 
	\to (M / \group{G})_\Total,
	\qquad
	\iota_{M/\group{G}}(\group{G}_\Wobs x) \coloneqq 
	\group{G}_\Total\iota_M(x)
\end{equation}
and equivalence relation on $\group{G}_\Wobs$ given by
\begin{equation}
	\group{G}_\Wobs x \sim_{M/\group{G}} \group{G}_\Wobs y  
	\coloniff \exists g,g' \in \group{G}_\Wobs : 
	(g \acts x) \sim_M (g' \acts y)
\end{equation}
for all $x, y \in M_\Wobs / \group{G}_\Wobs$,
is a constraint set.
\end{lemma}

\begin{proof}
	The map $\iota_{M/\group{G}}$ is well-defined since
	$\iota_M \circ \Phi_\Wobs = \Phi_\Total \circ 
	(\iota_\group{G}\times \iota_M)$
	holds by the definition of constraint group action.
	Moreover, it is easy to check that
	$\sim_{M/\group{G}}$ defines an equivalence 
	relation on $M_\Wobs / \group{G}_\Wobs$.
\end{proof}

We call $M/\group{G}$ the \emph{constraint orbit space}
of the action of $\group{G}$ on $M$.

\begin{remark}
	For a given constraint group action of $\group{G}$ on $M$ we 
	could also construct an equivalence relation internal to 
	$\ConSet$, i.e. a constraint subset $R_\group{G} \subseteq M 
	\times M$ with the usual properties.
	Then $M/\group{G}$ as defined above is indeed the coequalizer of 
	this internal equivalence relation, and $\sim_{M / \group{G}}$
	is just the pushforward relation of $\sim_M$ along the quotient 
	map.
\end{remark}

Constructing the constraint orbit space from a constraint group 
action is actually functorial.

\begin{proposition}[Orbit space functor]
	\label{prop:OrbitSpaceFunctor}
Mapping every constraint group action
$\Phi$ of $\group{G}$ on $M$ to its orbit space
$M / \group{G}$  defines a functor
$\ConOrb \colon \ConGroupAct \to \ConSet$.
\end{proposition}

\begin{proof}
	Consider an equivariant map $f \colon M \to N$ along a morphism
	$\phi \colon G \to H$ of constraint groups.
	By the classical theory we know that $f_\Total$ and $f_\Wobs$ 
	induce maps
	$\check{f}_\Total \colon M_\Total / \group{G}_\Total \to N_\Total
	/ \group{H}_\Total$
	and
	$\check{f}_\Wobs \colon M_\Wobs / \group{G}_\Wobs \to N_\Wobs
	/ \group{H}_\Wobs$
	which are compatible with $\iota_{M / \group{G}}$ and
	$\iota_{N / \group{H}}$.
	It remains to show that $\check{f}_\Wobs$ is compatible with the 
	equivalence relations.
	For this let $\group{G}_\Wobs x, \group{G}_\Wobs y \in M_\Wobs / 
	\group{G}_\Wobs$ be given with
	$\group{G}_\Wobs x \sim_{M / \group{G}} \group{G}_\Wobs y$.
	Hence there exist
	$g,g' \in \group{G}_\Wobs$ such that
	$g \acts x \sim_M g' \acts y$.
	Then, since $f_\Wobs$ is compatible with the equivalence 
	relations, we get
	\begin{align*}
		\phi(g) \acts f_\Wobs(x)
		= f_\Wobs(g \acts x)
		\sim_N f_\Wobs(g' \acts y)
		= \phi(g') \acts f_\Wobs(y),
	\end{align*}
	showing that 
	$\check{f}_\Wobs(\group{G}_\Wobs x)
	= \group{H}_\Wobs f_\Wobs(x) \sim_{N / \group{H}} 
	\group{H}_\Wobs f_\Wobs(y)
	= \check{f}_\Wobs(\group{G}_\Wobs y)$.
	Thus $\check{f}_\Wobs$ is a morphism of constraint sets.
\end{proof}

\subsubsection{Reduction of Constraint Groups}

As in the case of constraint sets we have a reduction functor
\glsadd{Groups}$\red \colon \ConGroups \to \Groups$ given by
\index{reduction!groups}
\begin{equation}
	\group{G}_\red = \group{G}_\Wobs / \group{G}_\Null.
\end{equation}
Note that $\group{G}_\Null$ is exactly the kernel of the projection 
map $\pi \colon \group{G}_\Wobs \to \group{G}_\red$.
Thus we immediately get
\index{reduction!automorphisms}
\begin{equation}
	\label{eq:ReducedAutGroup}
	\ConAut(M)_\red \subseteq \Aut(M_\red).
\end{equation}
The next example shows that, in general, we cannot 
expect more.

\begin{example}
Let $M_\Total = M_\Wobs = \left\{1,2,3\right\}$ with 
equivalence 
relation $\sim$ given by the only non trivial relation $2 
\sim 3$.
Then $M_\red = \{[1],[2]\}$.
The map $f([1]) = [2]$, $f([2]) = [1]$ is obviously 
invertible on 
$M_\red$, but there cannot exist an automorphism $g$ of $M$ 
with 
$g_\red = f$, since from this it would follow that $g(2) = 1 = g(3)$.
\end{example}

\begin{remark}\
	\label{rem:ReductionCommutesWithConstruction}
\begin{remarklist}
	\item It will be a recurring theme that an (often functorial) 
	construction on certain objects which we can also define for 
	their constraint analogues will commute with reduction.
	To be a bit more precise, consider the following picture:
	Assume we have a functorial construction
	$\functor{F} \colon \category{C} \to \category{D}$
	on a category $\category{C}$ with values in the category 
	$\category{D}$ and its constraint analogue 
	$\Con\functor{F} \colon \Con\category{C} \to \Con\category{D}$
	on the category of constraint objects ``internal'' to 
	$\category{C}$.
	Then the diagram
	\begin{equation} \label{diag:ConstructionVSReduction}
		\begin{tikzcd}
			\Con\category{C}
			\arrow[r,"\Con\functor{F}"]
			\arrow[d,"\red",swap]
			&\Con\category{D} 
			\arrow[d,"\red"]
			\arrow[dl,Rightarrow,"\eta"{swap}]\\
			\category{C}
			\arrow[r,"\functor{F}"]
			&\category{D}
		\end{tikzcd}
	\end{equation}
	will often commute up to a natural isomorphism.
	Nevertheless, there will also occur situations in which \eqref{diag:ConstructionVSReduction}
	only commutes up to an \emph{injective} natural transformation $\eta$.
	This typically happens when $\functor{F}$ and $\Con\functor{F}$
	construct certain limits in $\category{C}$ and $\Con\category{C}$, respectively.
	Since $\red$ does not necessarily commute with taking limits this leads 
	to $\eta$ not being an isomorphism.
	For us this will be of interest when $\functor{F}$, and hence $\Con\functor{F}$
	map sets to their subsets fulfilling a given equation.
	\item To make the above assignment of a constraint category $\Con\category{C}$ to a given category $\category{C}$
	precise we would need to restrict ourselves to categories allowing for a well-behaved notion of equivalence relation.
	Then we expect $\Con$ to be functorial, and every $\Con\category{C}$ would automatically admit a reduction functor
	$\red \colon \Con\category{C} \to \category{C}$.
	
\end{remarklist}
\end{remark}

Next we want to investigate how constraint group actions behave with 
respect to reduction.

\begin{lemma}[Reduction of group actions]
	\index{reduction!group action}
Let $\group{G}$ be a constraint group acting via\linebreak
$\Phi \colon \group{G} \to \ConAut(M)$
on a constraint set $M$.
Then $\Phi_\red$ defines an action of $\group{G}_\red$
on $M_\red$.
\end{lemma}

\begin{proof}
	Since reduction is functorial on the category of constraint 
	groups, and with the help of
	\eqref{eq:ReducedAutGroup} we see immediately that
	$\Phi$ reduces to
	\glsadd{Aut}$\Phi_\red \colon \group{G}_\red \to \Aut(M_\red)$,
	giving a group action of $\group{G}_\red$ on
	$M_\red$.
\end{proof}

Again, this is functorial.
To state this, we denote by
\glsadd{GroupAct}$\GroupAct$ the category of classical 
group actions and equivariant maps along group morphisms between them.

\begin{proposition}
Reducing constraint group actions defines a functor
\begin{equation}
	\red \colon \ConGroupAct \to \GroupAct.
\end{equation}
\end{proposition}

\begin{proof}
	Let $\group{G}$ and $\group{H}$ be constraint groups acting on
	constraint sets $M$ and $N$, respectively.
	Moreover, let $f \colon M \to N$ be an equivariant constraint map 
	along a constraint group morphism
	$\phi \colon \group{G} \to \group{H}$.
	These reduce to a map $f_\red \colon M_\red \to N_\red$
	and a group morphism $\phi_\red \colon \group{G}_\red \to 
	\group{H}_\red$, with
	\begin{align*}
		f_\red([g] \acts [x])
		= [f(g \acts x)]
		= [\phi(g) \acts f(x)]
		= \phi_\red([g]) \acts f_\red([x]),
	\end{align*}
	showing that $f_\red$ is equivariant along $\phi_\red$.
\end{proof}

This raises directly the question if constructing orbit spaces is 
compatible with reduction.
For this denote by $\Orb \colon \GroupAct \to \Sets$
the classical construction of the orbit space.
We obtain the following result, cf. 
\autoref{rem:ReductionCommutesWithConstruction}.

\begin{proposition}[Orbit spaces vs. reduction]
	\label{prop:OrbitSpaceVSReduction}
There exists a natural isomorphism $\eta$ making the following diagram 
commute:
\begin{equation}
\begin{tikzcd}
	\ConGroupAct
		\arrow[r,"\ConOrb"]
		\arrow[d,"\red"{swap}]
	& \ConSet 
		\arrow[d,"\red"]
		\arrow[dl,Rightarrow,"\eta"{swap}]\\
	\GroupAct
		\arrow[r,"\Orb"]
	& \Sets
\end{tikzcd}
\end{equation}
\end{proposition}

\begin{proof}
Define
$\eta \colon \ConOrb \circ \red \Longrightarrow \red \circ \Orb$
for every constraint action $\Phi \colon \group{G} \to \ConAut(M)$
by
\begin{equation*}
	\eta_\Phi \colon (M / \group{G})_\red \to M_\red / \group{G}_\red,
	\qquad
	\eta_\Phi([\group{G}_\Wobs x]) \coloneqq
	\group{G}_\red[x].
\end{equation*}
To see that this map is well-defined consider
$\group{G}_\Wobs y \sim_{M / \group{G}} \group{G}_\Wobs x$.
Then there exist $g,g' \in \group{G}_\Wobs$ such that
$g \acts x \sim_M g' \acts y$.
Hence
\begin{align*}
	\group{G}_\red[x]
	= \group{G}_\red[g^{-1}g' \acts y]
	= \group{G}_\red([g^{-1}g'] \acts [y])
	= \group{G}_\red[y]
\end{align*}
holds, showing that $\eta_\Phi$ does not depend on the choice of a 
representative.
Now $\eta_\Phi$ is obviously invertible with inverse given by
$\eta_\Phi^{-1}(\group{G}_\red[x]) = [\group{G}_\Wobs x]$.
To show that $\eta$ is a natural transformation consider another 
constraint action
$\Psi \colon \group{H} \to \ConAut(N)$
and an equivariant constraint map $f \colon M \to N$
along a constraint group morphism
$\phi \colon \group{G} \to \group{H}$.
Then for all $[\group{G}_\Wobs x] \in (M / \group{G})_\red$
we have
\begin{align*}
	\eta_\Psi\left(\check{f}_\red([\group{G}_\Wobs x]) \right)
	&= \eta_\Psi\left( [\group{H}_\Wobs f_\Wobs(x)] \right)
	= \group{H}_\red\left( [f_\Wobs(x)] \right)
	= \group{H}_\red \left( f_\red([x]) \right)\\
	&= \check{f}_\red\left( \group{G}_\red [x] \right)
	= \check{f}_\red\left( \eta_\Phi([\group{G}_\Wobs x]) \right).
\end{align*}
\end{proof}

\subsubsection{Strong constraint groups}
For completeness let us also remark on strong constraint groups.
As a group internal to $\strConSet$ a strong constraint group
consists of a group morphism $\iota_\group{G} \colon \group{G}_\Wobs 
\to \group{G}_\Total$ and normal subgroups
$\group{G}_\Null \subseteq \group{G}_\Wobs$
and $\group{G}^\Total_\Null \subseteq \group{G}_\Total$
such that
$\iota_\group{G}(\group{G}_\Null) \subseteq \group{G}^\Total_\Null$.
Moreover,
$\iota_\group{G}(\group{G}_\Wobs) \subseteq \group{G}_\Total$
needs to be saturated.
It is easy to see that this enforces $\group{G}^\Total_\Null = 
\iota_\group{G}(\group{G}_\Null)$.
Thus strong constraint groups can be defined as follows:

\begin{definition}[Strong constraint group]\
	\label{def:StrConGroup}
	\index{strong constraint!group}
\begin{definitionlist}
	\item A constraint group $G$ such that 
	$\iota_\group{G}(\group{G}_\Null) \subseteq \group{G}_\Total$
	is a normal subgroup is called \emph{strong constraint group}.
	\item A morphism of strong constraint groups is just a morphism 
	of constraint groups.
	\item The category of strong constraint groups will be denoted by
	$\strConGroups$.
\end{definitionlist}
\end{definition}

Note that in contrast to strong constraint sets the morphisms between 
strong constraint groups are just the morphisms of their underlying 
constraint groups.
Hence, we will write $\ConHom(\group{G},\group{H})$ instead of
$\strConHom(\group{G},\group{H})$
for the constraint set of constraint group homomorphisms,
cf. \autoref{prop:CategoryStrCSet} \ref{prop:CategoryStrCSet_3}.

From the definition it is clear that there exists a forgetful functor
$\functor{U} \colon \strConGroups \to \ConGroups$.
The reduction of strong constraint groups is then given by first 
forgetting to the category of constraint groups:
\begin{equation}
	\red = \red \circ \functor{U} \colon \strConGroups \to \ConGroups.
\end{equation}

\subsection{Constraint $\field{k}$-Modules}
	\label{sec:ConKModules}

We could continue by defining constraint rings as monoids internal to 
the category of abelian constraint groups.
Since we will not need these objects during this thesis we instead
move on to constraint modules over a certain class of rings.
For the rest of this chapter let $\field{k}$ be a commutative unital ring.

\begin{definition}[Constraint $\field{k}$-modules]\
	\label{def:ConkModules}%
	\glsadd{conModulesFont}
\begin{definitionlist}
	\item \index{constraint!$\field{k}$-module}A \emph{constraint $\field{k}$-module}
	is given by a triple
	$\module{E} =
	(\module{E}_\Total, \module{E}_\Wobs,\module{E}_\Null)$
	of $\field{k}$-modules together with a
	module homomorphism
	$\iota_\module{E} \colon \module{E}_\Wobs \to \module{E}_\Total$ 
	such that
	$\module{E}_\Null \subseteq \module{E}_\Wobs$ is a submodule.
	\item A \emph{morphism 
	$\Phi\colon \module{E} \to \module{F}$
	of constraint $\field{k}$-modules}
	is a pair
	$(\Phi_\Total, \Phi_\Wobs)$
	of module homomorphisms
	$\Phi_\Total \colon \module{E}_\Total \to \module{F}_\Total$
	and
	$\Phi_\Wobs\colon \module{E}_\Wobs \to \module{F}_\Wobs$
	such that
	$\Phi_\Total \circ \iota_\module{E}
	= \iota_{\module{F}} \circ \Phi_\Wobs$
	and
	$\Phi_\Wobs(\module{E}_\Null) \subseteq \module{F}_\Null$.
	\item The category of constraint $\field{k}$-modules is denoted by
	\glsadd{ConModk}$\ConMod_\field{k}$
	and the set of morphisms
	between constraint $\field{k}$-modules $\module{E}$ and
	$\module{F}$ is denoted by
	$\Hom_\field{k}(\module{E},\module{F})$.
\end{definitionlist}
\end{definition} 

There is an obvious forgetful functor
$\functor{U} \colon \ConMod_\field{k} \to \ConSet$,
forgetting all algebraic structures.
The equivalence relation on $\functor{U}(\module{E})_\Wobs$
is induced by the submodule $\module{E}_\Null \subseteq 
\module{E}_\Wobs$.
It can be shown that $(\ConMod_\field{k},\functor{U})$
is an algebraic category in the sense of
\cite{adamek.herrlich.strecker:1990a}
and hence behaves in many respects as we would expect from a category 
of objects equipped with algebraic structure.
In particular, as we will see in the next proposition,
many categorical constructions in $\ConMod_\field{k}$ are given by
the corresponding constructions in $\ConSet$ equipped with the 
structure of a constraint $\field{k}$-module.

\begin{proposition}[Co/limits in $\ConMod_\field{k}$]
	\label{prop:CoLimitsConModk}
Let $\module{E}$, $\module{F}$ and $\module{G}$ be constraint 
$\field{k}$-modules and let
$\Phi,\Psi \colon \module{E} \to \module{F}$
as well as
$\Theta \colon \module{G} \to \module{F}$
be constraint morphisms.
\begin{propositionlist}
	\item\label{prop:CoLimitsConModk_InitalFinal}
	\index{initial object!constraint $\field{k}$-modules}
	\index{final object!constraint $\field{k}$-modules}
	The initial and final object in $\ConMod_\field{k}$ agree 
	and are given by $0 \coloneqq (0,0,0)$.
	\item\label{prop:CoLimitsConModk_biproduct}
	\index{product!constraint $\field{k}$-modules}
	\index{coproduct!constraint $\field{k}$-modules}
	\index{direct sum!constraint $\field{k}$-modules}
	The binary product and binary coproduct in $\ConMod_\field{k}$
	agree and are given by
	\glsadd{directSum}
	\begin{equation}
	\begin{split}
		(\module{E} \oplus \module{F})_\Total 
		= \module{E}_\Total \oplus \module{F}_\Total, \\
		(\module{E} \oplus \module{F})_\Wobs 
		= \module{E}_\Wobs \oplus \module{F}_\Wobs, \\
		(\module{E} \oplus \module{F})_\Null 
		= \module{E}_\Null \oplus \module{F}_\Null,
	\end{split}
	\end{equation}
	with $\iota_\oplus = \iota_\module{E} + \iota_\module{F}
	\colon \module{E}_\Wobs \oplus \module{F}_\Wobs
	\to \module{E}_\Total \oplus \module{F}_\Total$.
	\item\label{prop:CoLimitsConModk_pullback}
	\index{pullback!constraint $\field{k}$-modules}
	The pullback of $\Phi$ and $\Theta$ is given by the 
	constraint $\field{k}$-module
	\begin{equation}
	\begin{split}
		(\module{E} \decorate*[_{\Phi}]{\times}{_{\Theta}} 
		\module{G})_\Total
		&= \module{E}_\Total 
		\decorate*[_{\Phi\Total}]{\times}{_{\Theta_\Total}} 
		\module{G}_\Total, \\
		(\module{E} \decorate*[_\Phi]{\times}{_\Theta} 
		\module{G})_\Wobs
		&= \module{E}_\Wobs 
		\decorate*[_{\Phi\Wobs}]{\times}{_{\Theta_\Wobs}} 
		\module{G}_\Wobs, \\
		(\module{E} \decorate*[_\Phi]{\times}{_\Theta} 
		\module{G})_\Null
		&= \module{E}_\Null 
		\decorate*[_{\Phi\Wobs}]{\times}{_{\Theta_\Wobs}} 
		\module{G}_\Null,
	\end{split}
	\end{equation}
	with projection maps
	\begin{align}
		(\pr_\Total^\module{E},\pr_\Wobs^\module{E}) 
		&\colon (\module{E} \decorate*[_{\Phi}]{\times}{_{\Theta}} 
		\module{G}) \longrightarrow \module{E}, \\
		(\pr_\Total^\module{G},\pr_\Wobs^\module{G}) 
		&\colon (\module{E} \decorate*[_{\Phi}]{\times}{_{\Theta}} 
		\module{G}) \longrightarrow \module{G}.
	\end{align}
	\item\label{prop:CoLimitsConModk_kernel}
	\index{kernel!constraint $\field{k}$-modules}
	The kernel of $\Phi$ is given by the constraint 
	$\field{k}$-module
	\glsadd{kernel}
	\begin{equation}
	\begin{split}
		\ker(\Phi)_\Total &= \ker(\Phi_\Total), \\ 
		\ker(\Phi)_\Wobs &= \ker(\Phi_\Wobs), \\ 
		\ker(\Phi)_\Null &= \ker(\Phi_\Wobs) \cap \module{E}_\Null,
	\end{split}
	\end{equation}	
	with $\iota_{\ker} \colon \ker(\Phi_\Wobs) \to \ker(\Phi_\Total)$ 
	the morphism induced by $\iota_\module{E}$.
	\item\label{prop:CoLimitsConModk_cokernel}
	\index{cokernel!constraint $\field{k}$-modules}
	The cokernel of $\Phi$ is given by the constraint 
	$\field{k}$-module
	\glsadd{cokernel}
	\begin{equation}
	\begin{split}
		\coker(\Phi)_\Total &= \module{F}_\Total / 
		\image(\Phi_\Total), \\
		\coker(\Phi)_\Wobs &= \module{F}_\Wobs / \image(\Phi_\Wobs), 
		\\
		\coker(\Phi)_\Null &= \module{F}_\Null / \image(\Phi_\Wobs),
	\end{split}
	\end{equation}
	with $\iota_{\coker} \colon \module{F}_\Wobs / \image(\Phi_\Wobs) 
	\to \module{F}_\Total / \image(\Phi_\Total)$
	the morphism induced by $\iota_\module{F}$.
	\item\label{prop:CoLimitsConModk_coequalizer}
	\index{coequalizer!constraint $\field{k}$-modules}
	The coequalizer of $\Phi$ and $\Psi$ is given by the 
	constraint $\field{k}$-module
	\glsadd{coequalizer}
	\begin{equation}
	\begin{split}
		\Coeq(\Phi,\Psi)_\Total
		&= \Coeq(\Phi_\Total, \Psi_\Total), \\
		\Coeq(\Phi,\Psi)_\Wobs
		&= \Coeq(\Phi_\Wobs, \Psi_\Wobs), \\
		\Coeq(\Phi,\Psi)_\Null
		&= q_\Wobs(\module{F}_\Null),
	\end{split}
	\end{equation}
	with $q = (q_\Total,q_\Wobs) \colon \module{F} \to 
	\Coeq(\Phi,\Psi)$.
	Here $q_\Total$ and $q_\Wobs$ denote the coequalizer morphisms
	of $\Phi_\Total$, $\Psi_\Total$ and $\Phi_\Wobs$, $\Psi_\Wobs$,
	respectively.
	\item\label{prop:CoLimitsConModk_hasCoLimits}
	The category $\ConMod_\field{k}$ has all finite limits and 
	colimits.
\end{propositionlist}
\end{proposition}

\begin{proof}
The form of the $\TOTAL$- and $\WOBS$-components follow
directly from \autoref{prop:CoLimitsCSet} and the classical 
characterization
of co/limits of $\field{k}$-modules.
It only remains to show that in each case the equivalence relation 
as given in \autoref{prop:CoLimitsCSet}
translates to the correct $\NULL$-components.
We show this for \ref{prop:CoLimitsConModk_kernel}, the rest can be 
done analogously.
For this note that $\ker(\Phi)$ is the equalizer of
$0 \colon \module{E} \to \module{F}$ and $\Phi$.
The equivalence relation on $\Eq(\Phi_\Wobs,0)$
is given by the restriction of the equivalence relation 
on $\module{E}$, thus $\ker(\Phi)_\Null = \ker(\Phi)_\Wobs \cap 
\module{E}_\Null$.
\end{proof}

In categories of sets equipped with algebraic structure, like groups, 
module algebras etc., we are used to the fact that a morphism 
respecting the algebraic structure is mono (epi) if and only if its 
underlying map of sets is mono (epi).
The same holds for $\ConMod_\field{k}$, forcing us to distinguish
regular from plain monos and epis.

\begin{proposition}[Mono- and epimorphisms in $\ConMod_\field{k}$]
	\label{prop:MonoEpisConModk}
Let $\Phi \colon \module{E} \to \module{F}$ be a morphism of 
constraint $\field{k}$-modules.
\begin{propositionlist}
	\item \index{monomorphism!constraint $\field{k}$-modules}
	$\Phi$ is a monomorphism if and only if $\Phi_\Total$ and 
	$\Phi_\Wobs$ are injective module homomorphisms.
	\item \index{epimorphism!constraint $\field{k}$-modules}
	$\Phi$ is an epimorphism if and only if $\Phi_\Total$ and 
	$\Phi_\Wobs$ are surjective module homomorphisms.
	\item \index{monomorphism!regular!constraint $\field{k}$-modules}
	$\Phi$ is a regular monomorphism if and only if it is a 
	monomorphism with
	$\Phi_\Wobs^{-1}(\module{F}_\Null) = \module{E}_\Null$.
	\item \index{epimorphism!regular!constraint $\field{k}$-modules}
	$\Phi$ is a regular epimorphism if and only if it is an 
	epimorphism with
	$\Phi_\Wobs(\module{E}_\Null) = \module{F}_\Null$.
\end{propositionlist}
\end{proposition}

\begin{proof}
This is just a repetition of the arguments used in the proof of
\autoref{prop:MonoEpiCSet}.
The conditions $(\Phi_\Wobs)^*(\sim_\module{F}) = \mathbin{\sim_\module{E}}$
and $(\Phi_\Wobs)_*(\sim_\module{E}) = \mathbin{\sim_\module{F}}$
for regular mono- and epimorphisms translate to
$\Phi_\Wobs^{-1}(\module{F}_\Null) = \module{E}_\Null$
and
$\Phi_\Wobs(\module{E}_\Null) = \module{F}_\Null$, respectively.
\end{proof}

\begin{example}\
\begin{examplelist}
	\item By the explicit formulas of \autoref{prop:CoLimitsConModk}
	\ref{prop:CoLimitsConModk_kernel} we see that the canonical
	inclusion $i \colon \ker(\Phi) \to \module{E}$
	is a regular monomorphism.
	\item By the explicit formulas of \autoref{prop:CoLimitsConModk}
	\ref{prop:CoLimitsConModk_coequalizer} we see that
	$q \colon \module{F} \to \Coeq(\Phi,\Psi)$
	is a regular epimorphism.
\end{examplelist}
\end{example}

\begin{remark}
Big parts of classical homological algebra solely rely on the fact 
that the usual categories of modules form abelian categories.
Since in $\ConMod_\field{k}$ regular and plain monos (or epis)
do not agree in general, it follows directly that $\ConMod_\field{k}$
is \emph{not} abelian.
This is the reason why in the theory of constraint algebraic objects
many effects appear which are unfamiliar if viewed from the point of 
view of classical algebra.
Another consequence is that we cannot rely on general techniques 
from abelian categories and hence we need to thoroughly examine even 
the most basic constructions in our categories of constraint algebraic 
objects.
\end{remark}

In any abelian category there is a canonical epi-mono factorization as
\begin{equation}
	\coker(\ker(\Phi)) \simeq \ker(\coker(\Phi))
\end{equation}
for every morphism $\Phi$.
\index{image!factorization}
In the non-abelian category $\ConMod_\field{k}$ there is no such 
canonical isomorphism, leading to two different factorizations:
We can either use $\coker(\ker(\Phi))$ to obtain an epi-regular mono 
factorization or we can use $\ker(\coker(\Phi))$ and get a 
regular epi-mono factorization.
These factorizations correspond to the image and regular image, 
respectively, see \autoref{def:ImagePreimageCSet}.
\filbreak
\begin{proposition}[Image and regular image]
	\label{def:ImageCModk}
Let $\Phi \colon \module{E} \to \module{F}$
be a morphism of constraint $\field{k}$-modules.
\begin{propositionlist}
	\item \index{image!constraint $\field{k}$-modules}
	The image of $\Phi$, as a morphism of constraint sets,
	is a constraint $\field{k}$-module given by\linebreak
	$\coker(\ker(\Phi))$.
	More explicitly:
	\glsadd{image}
	\begin{equation}
		\image(\Phi) \simeq 
		\big( \image(\Phi_\Total), \,\,
		\image(\Phi_\Wobs), \,\,
		\image(\Phi_\Wobs\at{\module{E}_\Null}) \big).
	\end{equation}
	\item \index{image!regular!constraint $\field{k}$-modules}
	The regular image of $\Phi$, as a morphism of constraint 
	sets, is a constraint $\field{k}$-module given by
	$\ker(\coker(\Phi))$.
	More explicitly:
	\glsadd{regimage}
	\begin{equation}
		\regimage(\Phi) \simeq 
		\big( \image(\Phi_\Total), \,\,
		\image(\Phi_\Wobs), \,\,
		\image(\Phi_\Wobs) \cap \module{F}_\Null \big).
	\end{equation}
\end{propositionlist}
\end{proposition}

\begin{proof}
Again the $\TOTAL$- and $\WOBS$-components are clear, since
$\ker(\coker)$ and $\coker(\ker)$ agree for classical
module morphisms.
For the $\NULL$-component we have by \autoref{prop:CoLimitsConModk}
\begin{equation*}
	\image(\Phi)_\Null = \coker(\ker \Phi)_\Null
	= \module{E}_\Null / \ker(\Phi_\Wobs)
	\simeq \image(\Phi_\Wobs\at{\module{E}_\Null})
\end{equation*}
and
\begin{equation*}
	\regimage(\Phi)_\Null
	= \ker(\coker \Phi)_\Null
	= \ker(\coker(\Phi)_\Wobs) \cap \module{F}_\Null
	\simeq \image(\Phi_\Wobs) \cap \module{F}_\Null.
\end{equation*}
\end{proof}

In analogy to constraint sets we can now define constraint submodules 
as follows.

\begin{definition}[Constraint submodule]
	\index{submodule of constraint module}
Let $\module{E}$ be a constraint $\field{k}$-module.
A \emph{constraint submodule} of $\module{E}$ 
consists of submodules
$\module{F}_\Total \subseteq \module{E}_\Total$
and $\module{F}_\Wobs \subseteq \module{E}_\Wobs$
such that $\iota_\module{E}(\module{F}_\Wobs) \subseteq 
\module{F}_\Total$.
\end{definition}

Every submodule can be understood as a regular monomorphism
$i \colon \module{F} \to \module{E}$
defined on the constraint module
$\module{F} = (\module{F}_\Total, \, \module{F}_\Wobs, \, i_\Wobs^{-1}(\module{E}_\Null))$.
Observe that the regular image of a morphism of constraint 
$\field{k}$-modules is a constraint submodule, while the image is not.

The existence of zero morphisms and coequalizers allows us to 
introduce quotients of constraint modules.

\begin{definition}[Quotient module]
	\label{def:ConQuotientkModule}
	\index{quotient of constraint modules}
Let $\module{F} \subseteq \module{E}$ be a constraint submodule.
The \emph{quotient $\module{E}/\module{F}$} is defined as the 
coequalizer of the inclusion $i \colon \module{F} \to \module{E}$
and the zero morphism $0 \colon \module{F} \to \module{E}$.
More explicitly:
\begin{align}
	\module{E} / \module{F}
	= \left( \module{E}_\Total / \module{F}_\Total,\,\,
	\module{E}_\Wobs / \module{F}_\Wobs,\,\,
	\module{E}_\Null /\module{F}_\Wobs \right).
\end{align}
\end{definition}

Here $\module{E}_\Null / \module{F}_\Wobs$ denotes the submodules of
$\module{E}_\Wobs / \module{F}_\Wobs$ generated by
equivalence classes $[x]$ of $x \in \module{E}_\Null$.
We could also define a quotient module with respect to more general 
submodules, i.e. non regular monomorphisms, but since
the coequalizer does not depend on the $\NULL$-component of 
$\module{F}$ this will not make a difference.
The independence of the quotient on the $\NULL$-component of the 
divisor will be important when we define constraint cohomology, see
\autoref{sec:ConHomologicalAlgebra}.

Let us now equip the category $\ConMod_\field{k}$ with the additional 
structure of a monoidal category, see \hyperref[sec:MonoidsModules]{Appendix~\ref{sec:MonoidsModules}} for the definition of 
monoidal categories.

\begin{proposition}[Monoidal structure on $\ConMod_\field{k}$]\
	\label{prop:MonoidalStructureCModk}
	\index{tensor product!constraint $\field{k}$-modules}
	\glsadd{tensork}
\begin{propositionlist}
	\item \label{prop:MonoidalStructureCModk_1}
	Let $\module{E}, \module{F} \in \ConMod_\field{k}$.
	Then
	\begin{equation}
	\begin{split}
		(\module{E} \tensor[\field{k}] \module{F})_\Total
		&\coloneqq \module{E}_\Total \tensor[\field{k}] 
		\module{F}_\Total, \\
		(\module{E} \tensor[\field{k}] \module{F})_\Wobs
		&\coloneqq \module{E}_\Wobs \tensor[\field{k}] 
		\module{F}_\Wobs, \\
		(\module{E} \tensor[\field{k}] \module{F})_\Null
		&\coloneqq \module{E}_\Null \tensor[\field{k}] 
		\module{F}_\Wobs
		+ \module{E}_\Wobs \tensor[\field{k}] \module{F}_\Null,
	\end{split}
	\end{equation}
	with $\iota_{\tensor} = \iota_\module{E} \tensor \iota_\module{F} 
	\colon \module{E}_\Wobs \tensor[\field{k}] \module{F}_\Wobs 
	\to \module{E}_\Total \tensor[\field{k}] \module{F}_\Total$,
	is a constraint $\field{k}$-module.
	\item \label{prop:MonoidalStructureCModk_2}
	The category $\ConMod_\field{k}$ equipped with the tensor 
	product
	$\tensor[\field{k}]$
	and unit $(\field{k}, \field{k},0)$
	is a symmetric monoidal category.
\end{propositionlist}
\end{proposition}

\begin{proof}
For the first part note that
$\module{E}_\Null \tensor[\field{k}] \module{F}_\Wobs$
denotes the submodule of $\module{E}_\Wobs \tensor[\field{k}] 
\module{F}_\Wobs$
generated by elements of the form $x \tensor y \in \module{E}_\Null 
\tensor[\field{k}] \module{F}_\Wobs$.
Then \ref{prop:MonoidalStructureCModk_1} is clear.
The constraint module $(\field{k},\field{k},0)$
is obviously the unit for $\tensor$.
It remains to show that there is an associativity isomorphism for
$\tensor$.
This is given by the usual associativity isomorphism on the
$\TOTAL$- and $\WOBS$-component, and it preserves the 
$\NULL$-component:
\begin{align*}
	((\module{E} \tensor[\field{k}] \module{F}) \tensor[\field{k}] 
	\module{G})_\Null
	= \module{E}_\Null \tensor[\field{k}] \module{F}_\Wobs 
	\tensor[\field{k}] \module{G}_\Wobs
	+ \module{E}_\Wobs \tensor[\field{k}] \module{F}_\Null 
	\tensor[\field{k}] \module{G}_\Wobs
	+ \module{E}_\Wobs \tensor[\field{k}] \module{F}_\Wobs 
	\tensor[\field{k}] \module{G}_\Null
	= (\module{E} \tensor[\field{k}] (\module{F} \tensor[\field{k}] 
	\module{G}))_\Null.
\end{align*}
{\ }
\end{proof}

It is easy to see that the set
$\Hom_\field{k}(\module{E},\module{F})$ 
of constraint morphisms between constraint $\field{k}$-modules 
carries the structure of a $\field{k}$-module, leading to a 
$\Modules_\field{k}$ enrichment on $\ConMod_\field{k}$.
The module $\Hom_\field{k}(\module{E},\module{F})$ can be enhanced to 
a constraint $\field{k}$-module and this internal hom turns out to be 
compatible with the tensor product of constraint modules.

\begin{proposition}[Internal hom in $\ConMod_\field{k}$]\
	\label{prop:InterlHomCModk}
	\index{internal hom!constraint $\field{k}$-modules}
\begin{propositionlist}
	\item \label{prop:InterlHomCModk_1}
	Let $\module{E},\module{F} \in \ConMod_\field{k}$.
	Then
	\glsadd{ConHomk}
	\begin{equation}\label{eq:InterlHomCModk}
	\begin{split}
		\ConHom_\field{k}(\module{E},\module{F})_\Total
		&\coloneqq \Hom_\field{k}(\module{E}_\Total, 
		\module{F}_\Total),\\
		\ConHom_\field{k}(\module{E},\module{F})_\Wobs
		&\coloneqq \Hom_\field{k}(\module{E},\module{F}),\\
		\ConHom_\field{k}(\module{E},\module{F})_\Null
		&\coloneqq \left\{ \Phi \in 
		\Hom_\field{k}(\module{E},\module{F}) \mid 
		\Phi_\Wobs(\module{E}_\Wobs) \subseteq \module{F}_\Null 
		\right\},
	\end{split}
	\end{equation}
	with
	$\iota_{\Hom} \colon \Hom_\field{k}(\module{E},\module{F})
	\ni (\Phi_\Total, \Phi_\Wobs) \mapsto \Phi_\Total \in 
	\Hom_\field{k}(\module{E}_\Total,\module{F}_\Total)$,
	is a constraint $\field{k}$-module.
	\item \label{prop:InterlHomCModk_2}
	For fixed $\module{E} \in \ConMod_\field{k}$ the functor
	$(\argument \tensor[\field{k}] \module{E}) \colon \ConMod_\field{k} 
	\to \ConMod_\field{k}$
	is left adjoint to $\ConHom_\field{k}(\module{E},\argument)$,
	i.e. $\ConMod_\field{k}$ is closed monoidal.
\end{propositionlist}
\end{proposition}

\begin{proof}
The proof is completely analogous to that of 
\autoref{def:ClosedMonoidalStructureCSet}.
Alternatively, observe that
\eqref{eq:InterlHomCModk}
is a constraint subset of $\ConMap(\module{E},\module{F})$
which is compatible with composition.
\end{proof}

The fact that $\ConMod_\field{k}$ is closed monoidal implies that 
there is a natural isomorphism
\begin{equation} \label{eq:TensorHomCModk}
	\Hom_\field{k}(\module{E} \tensor[\field{k}] 
	\module{F},\module{G}) \simeq 
	\Hom_\field{k}(\module{E},\ConHom_\field{k}(\module{F},\module{G}))
\end{equation}
for all $\module{E},\module{F},\module{G} \in \ConMod_\field{k}$.
A straightforward computation shows that \eqref{eq:TensorHomCModk}
can be enhanced to an isomorphism
\begin{equation} \label{eq:InternalTensorHomCModk}
	\ConHom_\field{k}(\module{E} \tensor[\field{k}] 
	\module{F},\module{G}) \simeq 
	\ConHom_\field{k}(\module{E},\ConHom_\field{k}(\module{F},\module{G}))
\end{equation}
of constraint $\field{k}$-modules.
Here the $\TOTAL$-component is just the usual tensor-hom adjunction of
$\field{k}$-modules and the $\WOBS$-component is exactly 
\eqref{eq:TensorHomCModk}.
It is also worth noting that as part of the adjunction we obtain
the evaluation map
\index{evaluation}
\glsadd{evaluation}
\begin{equation} \label{eq:EvaluationMap}
	\ev \colon \ConHom_\field{k}(\module{E},\module{F}) 
	\tensor[\field{k}] \module{E} \to \module{F},\qquad
	\ev_{\Total/\Wobs}(\Phi \tensor x) = \Phi(x)
\end{equation}
and the coevaluation map
\index{coevaluation}
\glsadd{coevaluation}
\begin{equation} \label{eq:CoevaluationMap}
	\coev \colon \module{F} \to
	\ConHom_\field{k}(\module{E},\module{E} \tensor[\field{k}] 
	\module{F}),\qquad
	\coev_{\Total/\Wobs}(y)(x) = x \tensor y.
\end{equation}
After investigating properties of the category $\ConMod_\field{k}$ 
itself, let us next look at how we can relate it to other known 
categories.
By the way we defined constraint $\field{k}$-modules it is clear that 
forgetting all algebraic structure and using the equivalence relation 
induced by the $\NULL$-component yields a functor
\begin{equation}
	\functor{U} \colon \ConMod_\field{k} \to \ConSet.
\end{equation}
It is then easy to see that 
it preserves finite limits and is lax closed,
since $\ConHom_\field{k}(\module{E},\module{F})
\subseteq \ConMap(\functor{U}(\module{E}), 
\functor{U}(\module{F}))$.
There is also a forgetful functor to the $\TOTAL$-component
$\functor{U}_\Total \colon \ConMod_{\field{k}} \to 
\Modules_\field{k}$, similar to \eqref{eq:ForgetfulToTotalConSets}.
Moreover, we can identify the category $\Modules_\field{k}$
of classical $\field{k}$-modules with the subcategory of 
$\ConMod_{\field{k}}$ consisting of constraint modules of the form
$(\module{E},\module{E},0)$.
We will often use this identification implicitly.
In particular we will write
$\field{k} = (\field{k},\field{k},0)$.

\subsubsection{Embedded Constraint Modules}

Similar to the case of embedded constraint sets we can also 
consider constraint $\field{k}$-modules $\module{E}$ with injective 
module morphism $\iota_\module{E}$.
We will denote the subcategory of $\ConMod_{\field{k}}$ consisting of 
such embedded constraint $\field{k}$-modules by
\glsadd{injConModk}$\injConMod_{\field{k}}$.

\begin{proposition}[The category $\injConMod_\field{k}$]\
	\label{prop:CategoryInjCMod}
	\index{embedded constraint!$\field{k}$-module}
\begin{propositionlist}
	\item $\injConMod_{\field{k}}$ is a reflective subcategory of
	$\ConMod_\field{k}$
	with reflector $\cdot^\inj \colon \ConMod_{\field{k}} \to 
	\injConMod_{\field{k}}$ given by
	\begin{equation}
		\module{E}^\inj \coloneqq (\module{E}_\Total, 
		\iota_\module{E}(\module{E}_\Wobs),
		\iota_\module{E}(\module{E}_\Null) ).
	\end{equation}
	\item The subcategory $\injConMod_{\field{k}}$
	of $\ConMod_\field{k}$ is closed under 
	finite limits.
	\item $\injConMod_{\field{k}}$ is closed symmetric monoidal with 
	respect to $\injtensor[\field{k}]$ defined by
	\glsadd{injTensork}
	\begin{equation}
		\module{E} \injtensor[\field{k}] \module{F} \coloneqq 
		(\module{E} 
		\tensor[\field{k}] \module{F})^\inj.
	\end{equation}
	\item The functor $\argument^\inj \colon 
	(\ConMod_{\field{k}},\tensor[\field{k}]) 
	\to (\injConMod_{\field{k}}, \injtensor[\field{k}])$
	is monoidal.
\end{propositionlist}
\end{proposition}

\begin{proof}
By definition $\injConMod_{\field{k}}$ is a full subcategory of
$\ConMod_{\field{k}}$.
To show that $\argument^\inj$ is left adjoint to the embedding
$\functor{U} \colon \injConMod_{\field{k}} \to 
\ConMod_{\field{k}}$
consider the natural transformations
\begin{equation*}
	\varepsilon \colon (\argument^\inj) \circ \functor{U} \Longrightarrow \id_{\injConMod_{\field{k}}}
	\quad\text{ and }\quad
	\eta \colon \id_{\ConMod_{\field{k}}} \Longrightarrow \functor{U} \circ (\argument^\inj)
\end{equation*}
given by $\varepsilon_\module{E} \coloneqq \id_{\module{E}}$
and
$\eta_\module{E} \coloneqq (\id_{\module{E}_\Total}, \iota_\module{E})$.
The triangle identities are then easily checked, and we 
immediately see that $\varepsilon$ is a natural isomorphism.
This yields the first part.
Since every reflective subcategory is closed under limits, the 
second part follows directly from the first.
For the last two parts we use a simple version of Day's reflection theorem, see
\autoref{thm:ReflectionTheorem}.
To see that $\argument^\inj$ is monoidal, 
see \autoref{def:LaxMonoidalFunctor},
it remains to show that 
\begin{equation*}
	(\eta_\module{E} \tensor \eta_\module{F})^\inj
	\colon (\module{E} \tensor[\field{k}] \module{F})^\inj 
	\to (\module{E}^\inj \tensor[\field{k}] \module{F}^\inj)^\inj
\end{equation*}
is an isomorphism for all $\module{E}, \module{F} \in 
\ConMod_{\field{k}}$.
This is clear, since
$(\eta_\module{E} \tensor \eta_\module{F})^\inj_\Wobs
	(\iota_\module{E}(x) \tensor \iota_\module{F}(y)) = 
	\iota_\module{E}(x) \tensor \iota_\module{F}(y)$.
\end{proof}

The tensor product of two injective module maps is in general not injective.
Thus $\module{E} \tensor[\field{k}] \module{F}$ might not be embedded, even if
$\module{E}$ and $\module{F}$ are.
The definition of $\injtensor[\field{k}]$ cures this defect.
However, this results in
$\functor{U} \colon \injConMod_{\field{k}} \to \ConMod_{\field{k}}$
not being a monoidal functor.
Moreover, $\injConMod_{\field{k}}$ is not closed under colimits as the next example shows, cf. \autoref{ex:injConSetIsNotClosedUnderColimits}.

\begin{example}
	\label{ex:QuotientEmbeddedConModk}
Consider the embedded constraint $\Reals$-module
$\module{E} = (\Reals^2, \Reals^2,0)$
and its embedded constraint submodule
$\module{F} = (\Reals, 0, 0)$.
Then its quotient
$\module{E} / \module{F} = (\Reals,\Reals^2,0)$
is not embedded.
\end{example}

\subsubsection{Reduction on $\ConMod_{\field{k}}$}

Since the $\NULL$-component encodes the equivalence relation on the 
$\WOBS$-components, reduction is simply given by their quotient.
We collect properties of the reduction functor: 

\begin{proposition}[Reduction on $\ConMod_\field{k}$]\
	\label{prop:ReductionOnCModk}
	\index{reduction!$\field{k}$-modules}
\begin{propositionlist}
	\item Mapping constraint $\field{k}$-modules $\module{E}$
	to $\module{E}_\red \coloneqq \module{E}_\Wobs / \module{E}_\Null$
	and constraint morphisms to the induced morphisms on the quotient 
	we obtain a functor
	\glsadd{Modk}
	\begin{equation}
		\red \colon \ConMod_\field{k} \to \Modules_\field{k}.
	\end{equation}
	\item The functor
	$\red \colon \ConMod_\field{k} \to \Modules_\field{k}$
	is monoidal.
	\item The functor
	$\red \colon \ConMod_\field{k} \to \Modules_\field{k}$
	is lax closed with injective natural
	transformation
	$\red \circ \ConHom_\field{k} \Rightarrow \Hom_\field{k} \circ (\red \times \red)$.
	\item The functor
	$\red \colon \ConMod_\field{k} \to \Modules_\field{k}$
	preserves finite limits and colimits.
\end{propositionlist}
\end{proposition}

\begin{proof}
The first part is obvious.
To show that $\red$ is monoidal, observe that $\field{k}_\red = \field{k} / 0 \simeq \field{k}$.
Moreover, $[x \tensor y] \mapsto [x] \tensor[] [y]$ gives an isomorphism
\begin{align*}
	(\module{E} \tensor[\field{k}] \module{F})_\red
	= (\module{E}_\Wobs \tensor[\field{k}] \module{F}_\Wobs) / (\module{E}_\Null \tensor[\field{k}] \module{F}_\Wobs
	+ \module{E}_\Wobs \tensor[\field{k}] \module{F}_\Null)
	\simeq \module{E}_\Wobs / \module{E}_\Null \tensor[\field{k}] \module{F}_\Wobs / \module{F}_\Null
	= \module{E}_\red \tensor[\field{k}] \module{F}_\red.
\end{align*}
These isomorphisms are clearly natural.
Since morphisms of constraint modules preserve the $\NULL$-component we obtain a morphism
$\eta_{\module{E},\module{F}} \colon \ConHom_\field{k}(\module{E},\module{F})_\red \to \Hom_\field{k}(\module{E}_\red,\module{F}_\red)$,
which is injective since $\ConHom_\field{k}(\module{E},\module{F})$ contains exactly those
morphisms which vanish after reduction.
This shows the lax closedness of $\red$.
It is easy to see that reduction preserves the co/limits listed in 
\autoref{prop:CoLimitsConModk}.
From this it follows directly that $\red$ preserves all finite limits and colimits.
\end{proof}

By contrast, reduction on the monoidal category $(\injConMod_\field{k},\injtensor[\field{k}])$
is in general not monoidal, since $\functor{U} \colon \injConMod_{\field{k}} \to \ConMod_{\field{k}}$ is not.

\subsection{Strong Constraint $\field{k}$-Modules}
\label{sec:strConModk}

Considering $\field{k}$-modules constructed internal to the category
$\strConSet$ of strong constraint sets we would obtain an abelian
strong constraint group
$\module{E} = (\module{E}_\Total,\module{E}_\Wobs,\module{E}_\Null)$
together with $\field{k}$-multiplications.
From \autoref{def:StrConGroup} it is clear that abelian strong constraint groups
do not differ from abelian constraint groups.
Thus, as objects, strong constraint $\field{k}$-modules coincide 
with constraint $\field{k}$-modules.
However, thinking of the $\NULL$-component as defining an equivalence 
relation on the $\TOTAL$-component leads to a different kind of 
tensor product.

To motivate the definition 
we anticipate the introduction of (strong) constraint algebras in \autoref{sec:ConAlgebras}:
A constraint algebra $\algebra{A}$ will be defined as a constraint $\field{k}$-module together with a multiplication
$\mu \colon \algebra{A} \tensor[\field{k}] \algebra{A} \to  \algebra{A}$.
By definition of $\tensor[\field{k}]$ this will implement $\algebra{A}_\Null$ as a two-sided ideal in $\algebra{A}_\Wobs$.
Now for a strong constraint algebra $\algebra{A}$ we expect $\algebra{A}_\Null$ to behave like a two-sided ideal in 
$\algebra{A}_\Total$.
To implement this idea, at least for embedded modules, we need to modify our tensor product to
\begin{equation} \label{eq:MotivationStrTensor_1}
	(\module{E} \strtensor[\field{k}] \module{F})_\Null
	= \module{E}_\Null \tensor[\field{k}] \module{F}_\Total + 
	\module{E}_\Total \tensor[\field{k}] \module{F}_\Null.
\end{equation}
In order to turn $\module{E} \strtensor[\field{k}] \module{F}$ into a 
constraint module we have to enlarge the $\WOBS$-component to
\begin{equation} \label{eq:MotivationStrTensor_2}
	(\module{E} \strtensor[\field{k}] \module{F})_\Wobs = 
	\module{E}_\Wobs \tensor[\field{k}] 
	\module{F}_\Wobs
	+ \module{E}_\Null \tensor[\field{k}] \module{F}_\Total 
	+ \module{E}_\Total \tensor[\field{k}] \module{F}_\Null.
\end{equation}
If we want to implement this tensor product also for non-embedded modules, we have to replace
the internal sum of submodules in \eqref{eq:MotivationStrTensor_1} and \eqref{eq:MotivationStrTensor_2} by an external direct sum.
To prevent counting elements in $\module{E}_\Null \tensor[\field{k}] \module{F}_\Wobs$
and elements in $\module{E}_\Wobs \tensor[\field{k}] \module{F}_\Null$ twice we have to quotient by an appropriate ideal.
This leads to the following definition:

\begin{proposition}[Strong tensor product]\
	\label{prop:MonoidalStructureCstrModk}
	\index{strong tensor product!constraint $\field{k}$-modules}
	\index{tensor product!strong|see{strong tensor product}}
\begin{propositionlist}
	\item Let $\module{E},\module{F} \in \ConMod_{\field{k}}$.
	Then
	\glsadd{strtensork}
	\begin{equation}
	\begin{split}
		(\module{E} \strtensor[\field{k}] \module{F})_\Total
		&\coloneqq \module{E}_\Total \tensor[\field{k}] 
		\module{F}_\Total, \\
		(\module{E} \strtensor[\field{k}] \module{F})_\Wobs
		&\coloneqq
		\frac{(\module{E}_\Wobs \tensor[\field{k}] 
		\module{F}_\Wobs)
		\oplus (\module{E}_\Null \tensor[\field{k}] \module{F}_\Total )
		\oplus (\module{E}_\Total \tensor[\field{k}] 
		\module{F}_\Null)}{\module{I}^\module{E}_\module{F}}, 
		\\
		(\module{E} \strtensor[\field{k}] \module{F})_\Null
		&\coloneqq
		\frac{\left( \module{E}_\Null \tensor[\field{k}] 
		\module{F}_\Wobs + \module{E}_\Wobs \tensor[\field{k}] 
		\module{F}_\Null \right)
		\oplus
		(\module{E}_\Null \tensor[\field{k}] \module{F}_\Total)
		\oplus (\module{E}_\Total \tensor[\field{k}] 
		\module{F}_\Null)}{\module{I}^\module{E}_\module{F}},
	\end{split}
	\end{equation}
	with
	\begin{equation}
	\begin{split}
		\module{I}^\module{E}_\module{F}
		\coloneqq &\Span_\field{k}\left\{
		(x_0 \tensor y,0,0) - (0,x_0 \tensor \iota_\module{F}(y),0)
		\mid x_0 \in \module{E}_\Null, y \in 
		\module{F}_\Wobs
		\right\} \\
		&+ \Span_\field{k}\left\{
		(x \tensor y_0,0,0) - (0,0,\iota_\module{E}(x) \tensor y_0)
		\mid x \in \module{E}_\Wobs, y_0 \in \module{F}_\Null
		\right\}
	\end{split}
	\end{equation}
	and $\iota_{\strtensor} = \iota_\module{E} \tensor 
	\iota_\module{F} + \iota_\module{E} \tensor \id_{\module{F}_\Total} + 
	\id_{\module{E}_\Total} \tensor \iota_\module{F}$,
	is a constraint $\field{k}$-module.
	\item The category $\ConMod_\field{k}$
	equipped with the tensor product
	$\strtensor[\field{k}]$
	is a symmetric monoidal category with unit
	$(\field{k}, \field{k},0)$.
\end{propositionlist}
\end{proposition}

\begin{proof}
	The first part is clear.
	For the second part consider the constraint $\field{k}$-module
	$\module{E} \strtensor[\field{k}] \module{F} \strtensor[\field{k}] \module{G}$
	defined by
		\begin{align*}
		(\module{E} \strtensor[\field{k}] \module{F} \strtensor[\field{k}] \module{G})_\Total
		&\coloneqq \module{E}_\Total \tensor[\field{k}] 
		\module{F}_\Total \tensor[\field{k}] \module{G}_\Total, \\
		(\module{E} \strtensor[\field{k}] \module{F} \strtensor[\field{k}] \module{G})_\Wobs
		&\coloneqq
		\frac{(\module{E}_\Wobs \tensor[\field{k}] 
			\module{F}_\Wobs \tensor[\field{k}] \module{G}_\Wobs)
			\oplus (\module{E}_\Null \tensor[\field{k}] \module{F}_\Total \tensor[\field{k}] \module{G}_\Total)
			\oplus (\module{E}_\Total \tensor[\field{k}] 
			\module{F}_\Null \tensor[\field{k}] \module{G}_\Total)
			\oplus (\module{E}_\Total \tensor[\field{k}] 
			\module{F}_\Total \tensor[\field{k}] \module{G}_\Null)}{\module{J}}, 
		\\
		(\module{E} \strtensor[\field{k}] \module{F} \strtensor[\field{k}] \module{G})_\Null
		&\coloneqq
		\frac{\left( \module{E} \tensor[\field{k}] 
			\module{F} \tensor[\field{k}] \module{G} \right)_\Null
			\oplus
			(\module{E}_\Null \tensor[\field{k}] \module{F}_\Total \tensor[\field{k}] \module{G}_\Total)
			\oplus (\module{E}_\Total \tensor[\field{k}] 
			\module{F}_\Null \tensor[\field{k}] \module{G}_\Total)
			\oplus (\module{E}_\Total \tensor[\field{k}] 
			\module{F}_\Total \tensor[\field{k}] \module{G}_\Null)}{\module{J}},
	\end{align*}
	with
	\begin{align*}
		\module{J} \coloneqq 
		&\Span_\field{k}\left\{ ((x_0 \tensor y \tensor z),0,0,0) - (0,(x_0 \tensor \iota_\module{F}(y) \tensor \iota_\module{G}(z)),0,0) \mid
		x_0 \in \module{E}_\Null, y \in \module{F}_\Wobs, z \in \module{G}_\Wobs \right\} \\
		&+\Span_\field{k}\left\{ ((x \tensor y_0 \tensor z),0,0,0) - (0,0,(\iota_\module{E}(x) \tensor y_0 \tensor \iota_\module{G}(z)),0) \mid
		x \in \module{E}_\Wobs, y \in \module{F}_\Null, z \in \module{G}_\Wobs \right\} \\
		&+\Span_\field{k}\left\{ ((x \tensor y \tensor z_0),0,0,0) - (0,0,0,(\iota_\module{E}(x) \tensor \iota_\module{F}(y) \tensor z)) \mid
		x \in \module{E}_\Wobs, y \in \module{F}_\Wobs, z_0 \in \module{G}_\Null \right\}.
	\end{align*}
	Note that we implicitly use the associativity of the classical tensor product.
	It is now easy to write down canonical isomorphisms between 
	$\module{E} \strtensor[\field{k}] \module{F} \strtensor[\field{k}] \module{G}$
	and
	$(\module{E} \strtensor[\field{k}] \module{F}) \strtensor[\field{k}] \module{G}$
	as well as
	$\module{E} \strtensor[\field{k}] (\module{F} \strtensor[\field{k}] \module{G})$
	by specifying it on every direct summand separately
	and checking that it is well-defined on the quotient by $\module{J}$.
	It is then a straightforward but incredibly tedious task to check all properties of
	a monoidal category, see
	\autoref{def: MonoidalCategory}
\end{proof}

\begin{definition}[The category $\strConMod_{\field{k}}$]
	\label{def:strConMod}
	\index{strong constraint!$\field{k}$-module}
	\glsadd{strConModk}
We call $\strtensor[\field{k}]$ the \emph{strong tensor product}
of constraint $\field{k}$-modules and denote the monoidal category 
$(\ConMod_{\field{k}},\strtensor[\field{k}])$
by $\strConMod_{\field{k}}$. 
\end{definition}

Note that as categories $\ConMod_{\field{k}}$ and $\strConMod_{\field{k}}$ are the same, they only differ by their monoidal structure.
Even though there is no difference between modules from $\ConMod_{\field{k}}$
and $\strConMod_{\field{k}}$ we will write $\module{E} \in 
\strConMod_{\field{k}}$
and call $\module{E}$ a strong constraint $\field{k}$-module if we want to stress that the tensor product to be used is $\strtensor[\field{k}]$.
Later on, when we consider modules over constraint algebras, we will need to distinguish
strong constraint from constraint modules more carefully.

Note that we can easily reformulate the $\NULL$-component as
\begin{equation}
	(\module{E} \strtensor[\field{k}] \module{F})_\Null
	\simeq \frac{(\module{E}_\Null \tensor[\field{k}] \module{F}_\Total)
	\oplus (\module{E}_\Total \tensor[\field{k}] \module{F}_\Null)}
	{\Span_\field{k} \left\{ (x_0 \tensor \iota_\module{F}(y_0),0) - 
		(0,(\iota_\module{E}(x_0) \tensor y_0)) \right\}}.
\end{equation}

To unwind the definition of $\strtensor[\field{k}]$ observe that 
$(\module{E} \strtensor[\field{k}] \module{F})_\Wobs$ is a colimit of 
$\field{k}$-modules:
\begin{equation}
	\label{diag:strTensorColimit}
\begin{tikzcd}
	{}
	&\module{E}_\Null \tensor[\field{k}] \module{F}_\Wobs
		\arrow[r,"\id\tensor \iota_\module{F}"]
		\arrow[d]
	& \module{E}_\Null \tensor[\field{k}] \module{F}_\Total
		\arrow[dd,dashed, bend left=15pt]\\
	\module{E}_\Wobs \tensor[\field{k}] \module{F}_\Null
		\arrow[r]
		\arrow[d,"\iota_\module{E} \tensor \id"{swap}]
	&\module{E}_\Wobs \tensor[\field{k}] \module{F}_\Wobs
		\arrow[dr,dashed]
	&{} \\
	\module{E}_\Total \tensor[\field{k}] \module{F}_\Null
		\arrow[rr,dashed,bend right=15pt]
	&{}
	&(\module{E} \strtensor[\field{k}] \module{F})_\Wobs
\end{tikzcd}
\end{equation}
From this the following characterization of morphisms on 
strong tensor products follows directly.

\begin{lemma}
	\label{lem:MorphismsOnStrConModules}
	Let $\module{E}, \module{F}, \module{G} \in 
	\strConMod_{\field{k}}$ be given.
	A constraint morphism 
	$\Phi \colon \module{E} \strtensor[\field{k}] \module{F}
	\to \module{G}$
	is equivalently given by module morphisms
	\begin{align}
		\Phi_\Total &\colon \module{E}_\Total \tensor[\field{k}] 
		\module{F}_\Total \to \module{G}_\Total,\\
		\Phi_\Wobs^{\Wobs\Wobs} &\colon \module{E}_\Wobs 
		\tensor[\field{k}] \module{F}_\Wobs \to \module{G}_\Wobs,\\
		\Phi_\Wobs^{\Total\Null} &\colon \module{E}_\Total 
		\tensor[\field{k}] \module{F}_\Null \to \module{G}_\Null,\\
		\Phi_\Wobs^{\Null\Total} &\colon \module{E}_\Null 
		\tensor[\field{k}] \module{F}_\Total \to \module{G}_\Null,
	\end{align}
such that
\begin{align}
	\iota_\module{G} \circ \Phi_\Wobs^{\Wobs\Wobs}
	&= \Phi_\Total \circ (\iota_\module{E} \tensor \iota_\module{F}),\\
	\iota_\module{G} \circ \Phi_\Wobs^{\Total\Null}
	&= \Phi_\Total \circ (\id_{\module{E}_\Total} \tensor 
	\iota_\module{F}),\\
	\iota_\module{G} \circ \Phi_\Wobs^{\Null\Total}
	&= \Phi_\Total \circ (\iota_\module{E} \tensor 
	\id_{\module{F}_\Total})
\end{align}
hold.
\end{lemma}

Since $\strConMod_\field{k}$ and $\ConMod_\field{k}$
are the same as categories, we see that $\strConMod$ obtains actually 
two monoidal structures $\strtensor[\field{k}]$
and $\tensor[\field{k}]$.
These are obviously not independent.

\begin{proposition}
	\label{prop:CompatibilityOfStrAndTensorCModk}
The identity functor $\strConMod_{\field{k}} \to 
\ConMod_{\field{k}}$ is lax monoidal with the morphism
$\module{E} \tensor[\field{k}] \module{F} \to \module{E} 
\strtensor[\field{k}] \module{F}$
given by the identity on the $\TOTAL$-component and
the inclusion in the first summand in the $\WOBS$-component
for all $\module{E},\module{F} \in \strConMod_{\field{k}}$.
\end{proposition}

\begin{proof}
Since the units for $\tensor[\field{k}]$ and $\strtensor[\field{k}]$ 
agree, the identity functor clearly preserves them.
For all $\module{E},\module{F} \in \strConMod_{\field{k}}$
the map $\mu_{\module{E},\module{F}} \colon \module{E} 
\tensor[\field{k}] \module{F} \to \module{E} 
\strtensor[\field{k}] \module{F}$ is defined by 
$(\mu_{\module{E},\module{F}})_\Total = \id$ and
$(\mu_{\module{E},\module{F}})_\Wobs = \pr \circ i_1$, with
\begin{equation*}
	i_1 \colon \module{E}_\Wobs \tensor[\field{k}] \module{F}_\Wobs \to
	(\module{E}_\Wobs \tensor[\field{k}] \module{F}_\Wobs)
	\oplus (\module{E}_\Null \tensor[\field{k}] \module{F}_\Total )
	\oplus (\module{E}_\Total \tensor[\field{k}] \module{F}_\Null)
\end{equation*}
the inclusion into the first component
and $\pr$ the projection on the quotient
as a constraint $\field{k}$-module morphism.
It is now a straightforward check that 
$\mu_{\module{E},\module{F}}$ fulfils the properties of a lax monoidal
functor, see \cite[Def. 7.5.1]{borceux:1994a}.
\end{proof}

Moreover, $\ConMod_{\field{k}}$, and therefore $\strConMod_{\field{k}}$ too, are $\ConMod_{\field{k}}$-enriched categories.
However, $\strConMod_{\field{k}}$ is not closed monoidal with respect to $\strtensor[\field{k}]$ but only with respect to
$\tensor[\field{k}]$.
Thus we will not repeat the structure and compatibilities for the internal hom.

\subsubsection{Embedded Strong Constraint Modules}

Recall that the subcategory of constraint modules $\module{E}$
with injective $\iota_\module{E}$ is denoted by
$\injConMod_\field{k}$.
In general we know that $\module{E} \strtensor[\field{k}] \module{F}$
might not be embedded, even if $\module{E}$ and $\module{F}$ are.
However, by \autoref{prop:CategoryInjCMod}  we see that
$\injConMod_{\field{k}}$ is a reflective subcategory of
$\ConMod_{\field{k}}$, and we can use this to define a new tensor product
on $\injConMod_\field{k}$.

\begin{proposition}\
	\label{prop:strTensorInjectiveModules}
	\index{embedded strong constraint!tensor product!$\field{k}$-modules}
	\index{tensor product!embedded strong|see{embedded strong constraint}}
\begin{propositionlist}
	\item $\injConMod_{\field{k}}$ is symmetric monoidal with respect to
	$\injstrtensor[\field{k}]$ defined by
	\glsadd{injstrtensork}
	\begin{equation}
		\module{E} \injstrtensor[\field{k}] \module{F}
		\coloneqq (\module{E} \strtensor[\field{k}] \module{F})^\inj.
	\end{equation}
	\item  The functor 
	$\argument^\inj \colon (\ConMod_{\field{k}}, \strtensor[\field{k}]) \to (\injConMod_{\field{k}},\injstrtensor[\field{k}])$
	is monoidal.
\end{propositionlist}
\end{proposition}

\begin{proof}
We use again Day's reflection theorem, see
\autoref{thm:ReflectionTheorem}.
The only thing left to show then is
$(\eta \tensor \eta)^\inj \colon (\module{E} \strtensor[\field{k}] \module{F})^\inj
\to (\module{E}^\inj \strtensor[\field{k}] \module{F}^\inj)^\inj$,
with $\eta_\module{E} \colon \module{E} \to \module{E}^\inj$
given by $(\eta_\module{E})_\Total = \id_{\module{E}_\Total}$
and $(\eta_\module{F})_\Wobs(x) = \iota_\module{E}(x)$,
is an isomorphism for all $\module{E}, \module{F} \in \ConMod_{\field{k}}$.
This is clear since
\begin{equation*}
	(\eta \tensor \eta)^\inj(\iota_\module{E}(x) \tensor \iota_\module{F}(y))
	= \iota_\module{E}(x) \tensor \iota_\module{F}(y).
\end{equation*}
\end{proof}

More explicitly, we have
\begin{equation}
	\label{eq:InjStrTensorCModk}
\begin{split}
	(\module{E} \injstrtensor[\field{k}] \module{F})_\Total
	&= \module{E}_\Total \tensor[\field{k}] \module{F}_\Total, \\
	(\module{E} \injstrtensor[\field{k}] \module{F})_\Wobs
	&= \module{E}_\Wobs \tensor[\field{k}] \module{F}_\Wobs 
	+\module{E}_\Total \tensor[\field{k}] \module{F}_\Null
	+\module{E}_\Null \tensor[\field{k}] \module{F}_\Total,\\
	(\module{E} \injstrtensor[\field{k}] \module{F})_\Total
	&=\module{E}_\Null \tensor[\field{k}] \module{F}_\Null
	+\module{E}_\Null \tensor[\field{k}] \module{F}_\Total,
\end{split}
\end{equation}
where we consider the $\WOBS$- and $\NULL$-components to be the 
submodules generated by elements of the given form.
This is exactly what we expected when motivating the definition of
$\strtensor[\field{k}]$.

\begin{definition}[The category $\injstrConMod_\field{k}$]
	\label{def:injstrConModk}
	\index{embedded strong constraint!$\field{k}$-module}
The monoidal category $(\injConMod_{\field{k}}, \injstrtensor[\field{k}])$
is denoted by
\glsadd{injstrConModk}$\injstrConMod_\field{k}$.
\end{definition}

In analogy to constraint modules it should be noted that the forgetful functor
$\functor{U} \colon \injstrConMod_{\field{k}} \to  \strConMod_{\field{k}}$
is not monoidal.
Summarizing, we showed that the two monoidal structures
$\tensor[\field{k}]$ and $\strtensor[\field{k}]$ induce
monoidal structures
$\injtensor[\field{k}]$ and $\injstrtensor[\field{k}]$
on $\injstrConMod_\field{k}$.
The compatibility from 
\autoref{prop:CompatibilityOfStrAndTensorCModk} carries over to the 
embedded
modules.

\begin{proposition}\
		\label{prop:QuotientOfTensorsConModk}
\begin{propositionlist}
	\item The identity functor
	$\injstrConMod_{\field{k}} \to \injConMod_{\field{k}}$
	is lax monoidal with
	$\module{E} \injtensor[\field{k}] \module{F} \to \module{E} 
	\injstrtensor[\field{k}] \module{F}$
	for all $\module{E},\module{F} \in \injstrConMod_{\field{k}}$
	given by the identity on the $\TOTAL$-component and
	the inclusion in the $\WOBS$-component.
	\item Let $\module{E},\module{F} \in \injstrConMod_\field{k}$
	be given.
	Then there is an isomorphism of constraint 
	$\field{k}$-modules 
	such that
	\begin{equation}
		\begin{split}
			\left( \frac{\module{E} \injstrtensor[k] 
			\module{F}}{\module{E} 
				\injtensor[\field{k}] \module{F}} \right)_\Total
			&\simeq 0, \\
			\left( \frac{\module{E} \injstrtensor[k] 
			\module{F}}{\module{E} 
			\injtensor[\field{k}] \module{F}} \right)_\Wobs
			&\simeq \left(\module{E}_\Null \injtensor[\field{k}] 
			\frac{\module{F}_\Total}{\module{F}_\Wobs}\right)
			\oplus
			\left(
			\frac{\module{E}_\Total}{\module{E}_\Wobs}
			\injtensor[\field{k}] \module{E}_\Null \right),	\\
			\left( \frac{\module{E} \injstrtensor[k] 
			\module{F}}{\module{E} 
			\injtensor[\field{k}] \module{F}} \right)_\Null
			&\simeq \left( \frac{\module{E} \injstrtensor[k] 
			\module{F}}{\module{E} 
			\injtensor[\field{k}] \module{F}} \right)_\Wobs.	
		\end{split}
	\end{equation}
\end{propositionlist}
\end{proposition}

\begin{proof}
The proof for the first part is completely analogous to that of
\autoref{prop:CompatibilityOfStrAndTensorCModk}.
The second part follows from the explicit description of
$\injstrtensor[\field{k}]$ in
\eqref{eq:InjStrTensorCModk}
and the definition of
$\injtensor[\field{k}]$ in  
\autoref{prop:MonoidalStructureCModk}.
\end{proof}

\subsubsection{Reduction}

The categories $\strConMod_{\field{k}}$ and $\ConMod_{\field{k}}$ only differ by their monoidal structure.
Hence the results from \autoref{prop:ReductionOnCModk} apply also for the reduction of strong constraint modules.
The only new structure introduced, namely $\strtensor[\field{k}]$,
coincides with $\tensor[\field{k}]$ after reduction as already indicated in
\autoref{prop:QuotientOfTensorsConModk}:

\begin{proposition}[Reduction on $\strConMod_{\field{k}}$]
	\label{prop:ReductionOnStrCModk}
	\index{reduction!tensor product!$\field{k}$-modules}
	\index{reduction!strong tensor product!$\field{k}$-modules}
The functor reduction functor
\begin{equation}
	\red \colon (\strConMod_{\field{k}},\strtensor[\field{k}]) \to 
	(\Modules_\field{k},\tensor[\field{k}])
\end{equation}
is monoidal.
In particular we have
\begin{equation}
	(\module{E} \strtensor[\field{k}] \module{F})_\red 
	\simeq (\module{E} \tensor[\field{k}] \module{F})_\red
\end{equation}
for all $\module{E}, \module{F} \in \strConMod_{\field{k}}$.
\end{proposition}

\begin{proof}
We directly have
\begin{align*}
	(\module{E} \strtensor[\field{k}] \module{F})_\red 
	= \frac{(\module{E} \strtensor[\field{k}] 
	\module{F})_\Wobs}{(\module{E} \strtensor[\field{k}] 
	\module{F})_\Null}
	\simeq \frac{\module{E}_\Wobs \tensor[\field{k}] 
	\module{F}_\Wobs}{\module{E}_\Null \tensor[\field{k}] 
	\module{F}_\Wobs + \module{E}_\Wobs \tensor[\field{k}] 
	\module{F}_\Null}
	= (\module{E} \tensor[\field{k}] \module{F})_\red
\end{align*}
for $\module{E},\module{F} \in \strConMod_{\field{k}}$.
The properties of a monoidal functor can by checked directly by writing out the above isomorphism on elements, see 
\autoref{def:LaxMonoidalFunctor} for the definition of monoidal functor.
\end{proof}




%% file: constraint-linalg.tex
Before we continue to construct (strong) constraint algebras in \autoref{sec:ConstraintAlgebrasModules},
let us take a step back and examine the structures present on the categories of (strong) constraint
$\field{k}$-modules in the special case of $\field{k}$ being a field.
In other words we will consider constraint vector spaces.

One of the main features distinguishing vector spaces from general modules is the existence of a basis.
For a constraint vector space $V$ a constraint basis should be given by a constraint subset $B \subseteq V$.
Though we could introduce bases this way, it is more convenient to use a slightly different notion of sets in the constraint setting, namely that of constraint index sets.
These will also play an important role in our study of free and projective constraint modules in
\autoref{sec:RegularProjectiveModules}.
We will introduce and study constraint index sets in \autoref{sec:ConIndSets} before we come back to constraint vector spaces 
in \autoref{sec:ConVectorSpaces}.
	
\subsection{Constraint Index Sets}
\label{sec:ConIndSets}

Recall that our algebraic notions, like constraint groups and modules, have underlying constraint sets.
This can be understood as a consequence of constructing these objects 
as algebraic objects internal to the category $\ConSet$ and its 
relatives.
Nevertheless, in our definitions of constraint groups and modules
we rephrased the equivalence relation on the 
$\WOBS$-component in terms of normal subgroups and submodules.
Thus instead of forgetting to the underlying constraint set, where we 
recover the equivalence relation from the $\NULL$-component, we could 
also forget all algebraic structure but keep the $\NULL$-component as 
a subset of the $\WOBS$-component.
This leads to a different notion of underlying set for constraint objects:
\filbreak
\begin{definition}[Constraint index sets]\
	\label{def:ConstraintIndexSets}
	\index{constraint!index sets}
\begin{definitionlist}
	\item A \emph{constraint index set} consists of a map
	$\iota_M \colon M_\Wobs \to M_\Total$ of sets together with a 
	subset\linebreak
	$M_\Null \subseteq M_\Wobs$.
	\item A morphism $f \colon M \to N$ of constraint sets $M$ and $N$
	(or \emph{constraint morphism}) consists of maps
	$f_\Total \colon M_\Total \to N_\Total$
	and 
	$f_\Wobs \colon M_\Wobs \to N_\Wobs$
	such that $f_\Total \circ \iota_M = \iota_N \circ f_\Wobs$
	and $f_\Wobs(M_\Null) \subseteq N_\Null$.
	\item The category of constraint index sets and their morphisms 
	is denoted by
	\glsadd{ConIndSet}$\ConIndSet$.
\end{definitionlist}
\end{definition}

\begin{example}
	There are obvious forgetful functors from the categories
	$\ConGroups$ and $\ConMod_{\field{k}}$ to
	$\ConIndSet$ by forgetting all algebraic structure.
\end{example}

\begin{proposition}[Co/limits in $\ConIndSet$]
	\label{prop:CoLimitsCIndSet}
Let $M$, $N$ and $P$ be constraint index sets and let
$f,g \colon M \to N$ as well as
$h \colon P \to N$ be constraint morphisms.
\begin{propositionlist}
	\item\label{prop:CoLimitsCIndSet_initial}
	\index{initial object!constraint index sets}
	The initial object in $\ConIndSet$ is given by
	$(\emptyset,\emptyset, \emptyset)$.
	\item\label{prop:CoLimitsCIndSet_final}
	\index{final object!constraint index sets}
	The final object in $\ConIndSet$ is given by 
	$(\{\pt\},\{\pt\},\{\pt\})$.
	\item\label{prop:CoLimitsCIndSet_product}
	\index{product!constraint index sets}
	The product is given by
	\glsadd{product}
	\begin{equation}
	\begin{split}
		(M \times N)_\Total &= M_\Total \times N_\Total, \\
		(M \times N)_\Wobs &= M_\Wobs \times N_\Wobs, \\
		(M \times N)_\Null &= M_\Null \times N_\Null,
	\end{split}
	\end{equation}
	with the product map
	$\iota_{M \times N}  = \iota_M \times \iota_N \colon M_\Wobs 
	\times N_\Wobs
	\longrightarrow M_\Total \times N_\Total$.
	\item\label{prop:CoLimitsCIndSet_coproduct}
	\index{coproduct!constraint index sets}
	The coproduct is given by
	\glsadd{coproduct}
	\begin{equation}
	\begin{split}
		(M \sqcup N)_\Total &= M_\Total \sqcup N_\Total, \\
		(M \sqcup N)_\Wobs &= M_\Wobs \sqcup N_\Wobs, \\
		(M \sqcup N)_\Null &= M_\Null \sqcup N_\Null,
	\end{split}
	\end{equation}
	with the coproduct map
	$\iota_M \sqcup \iota_N \colon M_\Wobs \sqcup N_\Wobs 
	\longrightarrow M_\Total \sqcup N_\Total$.
	\item\label{prop:CoLimitsCIndSet_pullback}
	\index{pullback!constraint index sets}
	The pullback of $f$ and $h$ is given by
	\begin{equation}
	\begin{split}
		(M \decorate*[_{f}]{\times}{_{h}} P)_\Total
		&= M_\Total \decorate*[_{f_\Total}]{\times}{_{h_\Total}} P_\Total,\\
		(M \decorate*[_f]{\times}{_h} P)_\Wobs
		&= M_\Wobs \decorate*[_{f_\Wobs}]{\times}{_{h_\Wobs}} P_\Wobs, \\
		(M \decorate*[_f]{\times}{_h} P)_\Null
		&= (f_\Wobs)^{-1}(N_\Null) \decorate*[_{f_\Wobs}]{\times}{_{h_\Wobs}} (h_\Wobs)^{-1}(N_\Null),
	\end{split}
	\end{equation}
	with projection maps
	\begin{align}
		(\pr_\Total^M,\pr_\Wobs^M) 
		&\colon (M \decorate*[_{f}]{\times}{_{h}} P) \longrightarrow 
		M,\\
		(\pr_\Total^P,\pr_\Wobs^P) 
		&\colon (M \decorate*[_{f}]{\times}{_{h}} P) \longrightarrow 
		N.
	\end{align}
	\item\label{prop:CoLimitsCIndSet_equalizer}
	\index{equalizer!constraint index sets}
	The equalizer of $f$ and $g$ is given by
	\glsadd{equalizer}
	\begin{equation}
	\begin{split}
		\Eq(f,g)_\Total
		&= \Eq(f_\Total, g_\Total) = \left\{ x \in M_\Total \mid f_\Total(x) = g_\Total(x)  \right\}, \\
		\Eq(f,g)_\Wobs
		&= \Eq(f_\Wobs, g_\Wobs)= \left\{ x \in M_\Wobs \mid f_\Wobs(x) = g_\Wobs(x)  \right\}, \\
		\Eq(f,g)_\Null
		&= i_\Wobs^{-1}(M_\Null)= \left\{ x \in M_\Null \mid f_\Wobs(x) = g_\Wobs(x)  \right\},
	\end{split}
	\end{equation}
	with
	$i =(i_\Total, i_\Wobs) \colon \Eq(f,g) \to M$
	given by the inclusions $i_\Total$ and $i_\Wobs$
	of $\Eq(f_\Total,g_\Total)$ 
	and $\Eq(f_\Wobs,g_\Wobs)$
	into $M_\Total$ and $M_\Wobs$, respectively.
	\item\label{prop:CoLimitsCIndSet_coequalizer}
	\index{coequalizer!constraint index sets}
	The coequalizer of $f$ and $g$ is given by
	\glsadd{coequalizer}
	\begin{equation}
	\begin{split}
		\Coeq(f,g)_\Total
		&= \Coeq(f_\Total, g_\Total), \\
		\Coeq(f,g)_\Wobs
		&= \Coeq(f_\Wobs, g_\Wobs),\\
		\Coeq(f,g)_\Null
		&= q_\Wobs(N_\Null)
	\end{split}
	\end{equation}
	with the morphism
	$q = (q_\Total,q_\Wobs) \colon N \to \Coeq(f,g)$
	of constraint index sets.
	Here the maps $q_\Wobs \colon N_\Wobs \to \Coeq(f_\Wobs,g_\Wobs)$
	and $q_\Total \colon N_\Total \to \Coeq(f_\Total,g_\Total)$
	denote the coequalizer in $\Sets$ of $f_\Wobs$, $g_\Wobs$ and
	$f_\Total$, $g_\Total$, respectively.
	\item\label{prop:CoLimitsCIndSet_coLimits}
	The category $\ConIndSet$ has all finite limits and colimits.
\end{propositionlist}
\end{proposition}

\begin{proof}
The proof follows from the same arguments as the proof of
\autoref{prop:CoLimitsCSet}.
In particular, the $\TOTAL$- and $\WOBS$-components are given
by the classical statements in $\Sets$.
The $\NULL$-component is then always given by the smallest subset
of the $\WOBS$-component such that the involved morphisms become 
constraint.
\end{proof}

As for constraint sets and constraint modules we have to distinguish
between monos (epis) and regular monos (epis).

\begin{proposition}[Mono- and epimorphisms in $\ConIndSet$]
	\label{prop:MonoEpiConIndSet}
Let $f \colon M \to N$ be a constraint morphism between 
constraint index sets.
\begin{propositionlist}
	\item \label{prop:MonoEpiConIndSet_1}
	\index{monomorphism!constraint index sets}
	$f$ is a monomorphism if and only if $f_\Total$ and 
	$f_\Wobs$ are injective maps.
	\item \label{prop:MonoEpiConIndSet_2}
	\index{epimorphism!constraint index sets}
	$f$ is an epimorphism if and only if $f_\Total$ and 
	$f_\Wobs$ are surjective maps.
	\item \label{prop:MonoEpiConIndSet_3}
	\index{monomorphism!regular!constraint index sets}
	$f$ is a regular monomorphism if and only if it is a 
	monomorphism with
	$f_\Wobs^{-1}(N_\Null) = M_\Null$.
	\item \label{prop:MonoEpiConIndSet_4}
	\index{epimorphism!regular!constraint index sets}
	$f$ is a regular epimorphism if and only if it is an 
	epimorphism with
	$f_\Wobs(M_\Null) = N_\Null$.
\end{propositionlist}
\end{proposition}

\begin{proof}
	Statements \ref{prop:MonoEpiConIndSet_1} and \ref{prop:MonoEpiConIndSet_2}
	follow by the same arguments used in \autoref{prop:MonoEpiCSet}.
	Then \ref{prop:MonoEpiConIndSet_3} and \ref{prop:MonoEpiConIndSet_4}
	follow by the characterization of equalizer and coequalizer in
	\autoref{prop:CoLimitsCIndSet}.
\end{proof}

Similarly to the case of constraint sets, it is not enough for a constraint morphism between constraint index sets to be an epimorphism and monomorphism in order to be invertible, cf. \autoref{lem:ConIsos}:

\begin{lemma}
	\label{lem:indConIsos}
Let $f \colon M \to N$ be a constraint morphism between constraint index sets.
The following statements are equivalent:
\begin{lemmalist}
	\item The constraint morphism $f$ is an isomorphism.
	\item The constraint morphism $f$ is a regular monomorphism and an epimorphism.
	\item The constraint morphism $f$ is a monomorphism and a regular epimorphism.
\end{lemmalist}
\end{lemma}

\begin{proof}
A constraint morphism is an isomorphism if and only if it is a bijection
on the $\TOTAL$-, $\WOBS$- and $\NULL$-components.
Being a bijection on $\TOTAL$- and $\WOBS$-components amounts to
$f_\Total$ and $f_\Wobs$ being bijective.
Moreover, $f_\Wobs$ restricts to a bijection on the $\NULL$-component
if and only if $f_\Wobs(M_\Null) = N_\Null$ or equivalently
$f_\Wobs^{-1}(N_\Null) = M_\Null$.
\end{proof}

We define subsets of constraint index sets as images of regular monomorphisms.

\begin{definition}[Constraint index subsets]
	\index{subset!constraint index set}
A \emph{constraint subset} of a constraint index set $M$ consists of subsets
$U_\Total \subseteq M_\Total$ and $U_\Wobs \subseteq M_\Wobs$
such that $\iota_M(U_\Wobs) \subseteq U_\Total$.
\end{definition}

We can view a constraint subset $(U_\Total,U_\Wobs)$ of a constraint index set $M$ itself
as a constraint index set $U = (U_\Total, U_\Wobs, U_\Wobs \cap M_\Null)$
with a regular monomorphism $i \colon U \to M$ as embedding.
For constraint subsets of constraint index set the following definitions will be useful:

\begin{definition}[Union and Intersection of constraint subsets]
Let $M \in \ConIndSet$ and constraint subsets $U,V \subseteq M$ be given.
\begin{definitionlist}
	\item \index{union of constraint index subsets}
	The \emph{intersection} of $U$ and $V$ is defined by
	\glsadd{intersection}
	\begin{equation}
		U \cap V \coloneqq (U_\Total \cap V_\Total, U_\Wobs \cap 
		V_\Wobs, U_\Null \cap V_\Null),
	\end{equation}
	with $\iota_{U\cap V} = \iota_M\at{U_\Wobs \cap V_\Wobs}$.
	\item \index{intersection of constraint index subsets}
	The \emph{union} of $U$ and $V$ is defined by
	\glsadd{union}
	\begin{equation}
		U \cup V \coloneqq (U_\Total \cup V_\Total, U_\Wobs \cup 
		V_\Wobs, U_\Null \cup V_\Null),
	\end{equation}
	with $\iota_{U \cup V} = \iota_M\at{U_\Wobs \cup V_\Wobs}$.
\end{definitionlist}
\end{definition}

Note that $U\cup V$ and $U \cap V$ form again subsets of $M$ since
$U_\Null \cap V_\Null = (U_\Wobs \cap V_\Wobs) \cap M_\Null$
and
$U_\Null \cup V_\Null = (U_\Wobs \cup V_\Wobs) \cup M_\Null$.

Let us from now on focus on \emph{embedded} constraint index sets, i.e. those constraint index sets $M$
with injective
$\iota_M \colon M_\Wobs \to M_\Total$.
We will denote their category by
\glsadd{injConIndSet}$\injConIndSet$.
Even though most of what follows can be also considered inside the bigger category
$\ConIndSet$ this would only complicate the exposition, and it will not be needed in the rest of the thesis.

For (embedded) strong constraint $\field{k}$-modules we have constructed two different kinds of tensor products.
There are now similar constructions available for embedded constraint index sets, which are not present in the classical category of sets.

\begin{definition}[Tensor products and dual]
	\label{def:TensorProdDualConIndSet}
Let $M,N \in \injConIndSet$.
\begin{definitionlist}
	\item \index{tensor product!constraint index sets}
	The \emph{tensor product} of $M$ and $N$ is defined by
	\glsadd{tensor}
	\begin{equation}
		\begin{split}
			(M \tensor N)_\Total 
			&\coloneqq M_\Total \times N_\Total, \\
			(M \tensor N)_\Wobs
			&\coloneqq M_\Wobs \times N_\Wobs, \\
			(M \tensor N)_\Null 
			&\coloneqq (M_\Wobs \times N_\Null) \cup (M_\Null \times N_\Wobs).
		\end{split}
	\end{equation}
	\item \index{strong tensor product!constraint index sets}
	The \emph{strong tensor product} of $M$ and $N$ is defined by
	\glsadd{strtensor}
	\begin{equation}
		\begin{split}
			(M \strtensor N)_\Total  
			&\coloneqq M_\Total \times N_\Total,\\
			(M \strtensor N)_\Wobs 
			&\coloneqq (M_\Wobs \times N_\Wobs) \cup (M_\Total \times N_\Null) \cup (M_\Null \times N_\Total), \\
			(M \strtensor N)_\Null
			&\coloneqq (M_\Total \times N_\Null) \cup (M_\Null \times N_\Total).
		\end{split}
	\end{equation}
	\item \index{dual!constraint index set}
	The \emph{dual} of $M$ is defined by
	\begin{equation}
		\begin{split}
			(M^*)_\Total &\coloneqq M_\Total,\\
			(M^*)_\Wobs &\coloneqq M_\Total \setminus M_\Null,\\
			(M^*)_\Null &\coloneqq M_\Total \setminus M_\Wobs.
		\end{split}
	\end{equation}
	\item \index{reduction!constraint index sets}
	The \emph{reduction} of $M$ is defined by
	\begin{equation}
		M_\red \coloneqq M_\Wobs \setminus M_\Null.
	\end{equation}
\end{definitionlist}
\end{definition}

Since there is no dual for the sets $M_\Wobs$ and $N_\Null$
we will often write $M^*_\Wobs$ and $M^*_\Null$
instead of $(M^*)_\Wobs$ and $(M^*)_\Null$, respectively.
All of the above constructions can be shown to be functorial.
Moreover, it is easy to see that $\times$, $\tensor$ and $\strtensor$ yield
monoidal structures on $\injConIndSet$.
Using the dual we can decompose the (strong) tensor product as follows.

\begin{lemma}
	\label{lem:DecomposeStrTensorConIndSets}
Let $M,N \in \injConIndSet$.
\begin{lemmalist}
	\item It holds that
	\begin{equation}
		\begin{split}
			(M \tensor N)_\Total
			&= (M \tensor N)_\Wobs \sqcup (M^*_\Null \times N_\Total) \sqcup (M_\Total \times N^*_\Null),\\
			(M \tensor N)_\Wobs
			&= (M \tensor N)_\Null \sqcup \left(M_\red \times N_\red\right).
		\end{split}
	\end{equation}
	\item It holds that
	\begin{equation}
		\begin{split}
			(M \strtensor N)_\Total
			&= (M \strtensor N)_\Wobs \sqcup (M^*_\Null \times N^*_\Wobs)
			\sqcup (M^*_\Wobs \times N^*_\Null),\\
			(M \strtensor N)_\Wobs
			&= (M \strtensor N)_\Null \sqcup \left(M_\red \times N_\red\right)\\
			&= (M \tensor N)_\Wobs \sqcup (M^*_\Null \times N_\Null) \sqcup (M_\Null \times N^*_\Null),\\
			(M \strtensor N)_\Null
			&= (M \tensor N)_\Null \sqcup (M^*_\Null \times N_\Null) \sqcup (M_\Null \times N^*_\Null).
		\end{split}
	\end{equation}
\end{lemmalist}
\end{lemma}

It can be useful to picture the components of $\tensor$ and $\strtensor$ as subsets
of the cartesian product as follows:

\begin{figure}[h]
	\centering
	\begin{subfigure}{0.4\linewidth}
		\begin{tikzpicture}[scale=1.5]
			\draw (1.5,3) node [anchor=south] {$M \tensor N$};
			
			\fill[pattern={Lines[angle=45,distance={20pt},line width={10pt}]},pattern color=c1] (0,0)
			-- (0,2) -- (1,2) -- (1,1) -- (2,1) -- (2,0) -- (0,0);
			\fill[pattern={Lines[angle=45,distance={20pt},line width={10pt},yshift={14pt}]},pattern color=c2] (0,0) rectangle (2,2);
			
			\draw[step=1cm] (0,0) grid (3,3);
			
			\draw (0,0.5) node [anchor=east] {$N_\Null$};
			\draw (0,1.5) node [anchor=east] {$N_\Wobs \setminus N_\Null$};
			\draw (0,2.5) node [anchor=east] {$N_\Total \setminus N_\Null$};
			
			\draw (0.5,0) node [anchor=north] {$M_\Null$};
			\draw (1.5,0) node [anchor=north] {$M_\Wobs \setminus M_\Null$};
			\draw (2.5,0) node [anchor=north] {$M_\Total \setminus M_\Null$};
		\end{tikzpicture}
	\end{subfigure}
	\begin{subfigure}{0.4\linewidth}
		\begin{tikzpicture}[scale=1.5]
			\draw (1.5,3) node [anchor=south] {$M \strtensor N$};
			
			\fill[pattern={Lines[angle=45,distance={20pt},line width={10pt}]},pattern color=c1] (0,0)
			-- (0,3) -- (1,3) -- (1,1) -- (3,1) -- (3,0) -- (0,0);
			\fill[pattern={Lines[angle=45,distance={20pt},line width={10pt},yshift={14pt}]},pattern color=c2] (0,0) -- (0,3) -- (1,3) -- (1,2) -- (2,2) -- (2,1) -- (3,1) -- (3,0) -- (0,0);
			
			\draw[step=1cm] (0,0) grid (3,3);
			
			\draw (0,0.5) node [anchor=east] {$N_\Null$};
			\draw (0,1.5) node [anchor=east] {$N_\Wobs \setminus N_\Null$};
			\draw (0,2.5) node [anchor=east] {$N_\Total \setminus N_\Null$};
			
			\draw (0.5,0) node [anchor=north] {$M_\Null$};
			\draw (1.5,0) node [anchor=north] {$M_\Wobs \setminus M_\Null$};
			\draw (2.5,0) node [anchor=north] {$M_\Total \setminus M_\Null$};
		\end{tikzpicture}
	\end{subfigure}
	\caption{$M \tensor N$ and $M \strtensor N$ as subsets of $M \times N$. $\WOBS$-components in green, $\NULL$-components in blue.}
	\label{fig:ProductsConIndSets}
\end{figure}

\begin{notation}
	\label{not:ConIndSetProducts}
We will use scaled down versions of 
\autoref{fig:ProductsConIndSets}.
These will be rotated by $45^\circ$ counter-clockwise, such that
$M_\Null \times N_\Null$ is represented by the 
bottom diamond.
For example, for constraint index sets $M$ and $N$  we write
\glsadd{condiamond}
\begin{align*}
	(M \tensor N)_\Wobs &= M \ConGrid[2][2][0][2][2] N &&
	&(M\strtensor N)_\Wobs &= M \ConGrid[2][2][2][2][2][0][2] N \\
	(M \tensor N)_\Null &= N \ConGrid[2][2][0][2] N &&
	&(M\strtensor N)_\Null &= M \ConGrid[2][2][2][2][0][0][2] N
\end{align*}
as subsets of $M_\Total \times N_\Total$.

We can also combine the whole constraint index set into one picture by using an overlay of $\WOBS$-and $\NULL$-component:
\begin{align*}
	M \tensor N &= M \ConGrid[1][1][0][1][2] N &&
	&M \strtensor N &= M \ConGrid[1][1][1][1][2][0][1] N.
\end{align*}
Observe that in this notation the dual is given by inverting the colours, i.e.
white becomes black, black becomes white and grey stays grey.
We can also replace $M$ and $N$ by their duals if we also reflect the diamond along its horizontal axis, e.g.
$M \ConGrid[1][1][0][2][2] N = M^* \ConGrid[0][0][0][0][2][1][0][2][1] N^*$.
The reduction of a diamond is given by its grey parts, e.g.
$(M \ConGrid[0][2][1][2][0][0][1] N )_\red = M \ConGrid[0][2][0][2] N$.
\end{notation}

With this notation the following compatibilities are easy to prove.

\begin{proposition}
	\label{prop:DualsConIndSet}
Let $M,N \in \injConIndSet$.
\begin{propositionlist}
	\item \label{prop:DualsConIndSet_1}
	We have
	\begin{equation}
		(M \tensor N)_\red = M_\red \times N_\red.
	\end{equation}
	\item \label{prop:DualsConIndSet_2}
	We have
	\begin{equation}
		(M \strtensor N)_\red = M_\red \times N_\red.
	\end{equation}
	\item \label{prop:DualsConIndSet_3}
	We have 
	\begin{equation}
		(M^*)_\red = M_\red.
	\end{equation}
	\item \label{prop:DualsConIndSet_4}
	We have
	\begin{equation}
		(M\tensor N)^* = M^* \strtensor N^*.
	\end{equation}
	\item \label{prop:DualsConIndSet_5}
	We have
	\begin{equation}
		(M \strtensor N)^* = M^* \tensor N^*.
	\end{equation}
	\item \label{prop:DualsConIndSet_6}
	We have
	\begin{equation}
		(M^*)^* = M.
	\end{equation}
\end{propositionlist}
\end{proposition}

\begin{proof}
We compute \ref{prop:DualsConIndSet_3} and \ref{prop:DualsConIndSet_6} explicitly:
We have
\begin{align*}
	(M^*)_\red 
	&= (M^*)_\Wobs \setminus (M^*)_\Null
	= (M_\Total \setminus M_\Null) \setminus (M_\Total \setminus M_\Wobs)
	= M_\Wobs \setminus M_\Null
	= M_\red
	\shortintertext{and}
	\big((M^*)^*\big)_\Wobs 
	&= (M^*)_\Total \setminus (M^*)_\Null
	= M_\Total \setminus (M_\Total \setminus M_\Wobs)
	= M_\Wobs, \\
	\big((M^*)^*\big)_\Null
	&= (M^*)_\Total \setminus (M^*)_\Wobs
	= M_\Total \setminus (M_\Total \setminus M_\Null)
	= M_\Null.
\end{align*}
The rest is a straightforward application of the notation introduced in 
\autoref{not:ConIndSetProducts}:
\begin{align*}
	(M \tensor N)_\red 
	&= (M \ConGrid[1][1][0][1][2] N)_\red
	= M \ConGrid[0][0][0][0][2] N
	= M_\red \times N_\red, \\
	(M \strtensor N)_\red 
	&= (M \ConGrid[1][1][1][1][2][0][1] N)_\red
	= M \ConGrid[0][0][0][0][2] N
	= M_\red \times N_\red,\\
	\shortintertext{and}
	(M \tensor N)^*
	&= (M \ConGrid[1][1][0][1][2] N)^*
	= M \ConGrid[0][0][1][0][2][1][1][1][1] N
	= M^* \ConGrid[1][1][0][1][2] N^*
	= M^* \strtensor N^*,\\
	(M \strtensor N)^*
	&= (M \ConGrid[1][1][1][1][2][0][1] N)^*
	= M \ConGrid[0][0][0][0][2][1][0][1][1] N
	= M^* \ConGrid[1][1][0][1][2] N^*
	= M^* \tensor N^*.
\end{align*}
\end{proof}

For a finite constraint index set $M = (M_\Total,M_\Wobs,M_\Null)$ we can define its \emph{cardinality} as
\index{cardinality of constraint index set}
\begin{equation}
	\label{eq:ConIndSetCardinality}
	\abs{M} \coloneqq (\abs{M_\Total}, \abs{M_\Wobs}, \abs{M_\Null}).
\end{equation}
Thus every finite constraint index set $M$ has an associated cardinality consisting
of three natural numbers $\abs{M}_\Total \coloneqq \abs{M_\Total}$,
$\abs{M}_\Wobs \coloneqq \abs{M_\Wobs}$
and $\abs{M}_\Null \coloneqq \abs{M_\Null}$
with $\abs{M}_\Null \leq \abs{M}_\Wobs$.
If $M$ is embedded we additionally have $\abs{M}_\Wobs \leq \abs{M}_\Total$.

\begin{corollary}
	\label{cor:CardinalityOfConstructionsConIndSet}
Let $M = (M_\Total, M_\Wobs, M_\Null)$ and
$N = (N_\Total, N_\Wobs, N_\Null)$
be finite embedded constraint index sets.
\begin{corollarylist}
	\item The cardinality of the product of $M$ and $N$ is 
	given by
	\begin{equation}
		\begin{split}
			\abs{M \times N}_\Total 
			&= \abs{M}_\Total \cdot \abs{N}_\Total, \\
			\abs{M \times N}_\Wobs
			&= \abs{M}_\Wobs \cdot \abs{N}_\Wobs, \\
			\abs{M \times N}_\Null 
			&= \abs{M}_\Null \cdot \abs{N}_\Null.
		\end{split}
	\end{equation}
	\item The cardinality of the disjoint union of $M$ and $N$ is 
	given by
	\begin{equation}
		\begin{split}
			\abs{M \sqcup N}_\Total 
			&= \abs{M}_\Total + \abs{N}_\Total, \\
			\abs{M \sqcup N}_\Wobs
			&= \abs{M}_\Wobs + \abs{N}_\Wobs, \\
			\abs{M \sqcup N}_\Null 
			&= \abs{M}_\Null + \abs{N}_\Null.
		\end{split}
	\end{equation}
	\item The cardinality of the tensor product of $M$ and $N$ is 
	given by
	\begin{equation}
		\begin{split}
			\abs{M \tensor N}_\Total 
			&= \abs{M}_\Total \cdot \abs{N}_\Total, \\
			\abs{M \tensor N}_\Wobs
			&= \abs{M}_\Wobs \cdot \abs{N}_\Wobs, \\
			\abs{M \tensor N}_\Null 
			&= \abs{M}_\Wobs \cdot \abs{N}_\Null + \abs{M}_\Null \cdot 
			\abs{N}_\Wobs
			- \abs{M}_\Null \cdot \abs{N}_\Null.
		\end{split}
	\end{equation}
	\item The cardinality of the strong tensor product of $M$ and $N$ 
	is given by
	\begin{equation}
		\begin{split}
			\abs{M \strtensor N}_\Total  
			&= \abs{M}_\Total \cdot \abs{N}_\Total,\\
			\abs{M \strtensor N}_\Wobs 
			&= \abs{M}_\Wobs \cdot \abs{N}_\Wobs + (\abs{M}_\Total - 
			\abs{M}_\Wobs) \cdot \abs{N}_\Null 
			+ \abs{M}_\Null \cdot ( \abs{N}_\Total - \abs{N}_\Wobs), \\
			\abs{M \strtensor N}_\Null
			&= \abs{M}_\Total \cdot \abs{N}_\Null + \abs{M}_\Null \cdot 
			\abs{N}_\Total - \abs{M}_\Null \cdot \abs{N}_\Null.
		\end{split}
	\end{equation}
	\item The cardinality of the dual of $M$ is given by
	\begin{equation}
		\begin{split}
			\abs{M^*}_\Total &= \abs{M}_\Total,\\
			\abs{M^*}_\Wobs &= \abs{M}_\Total - \abs{M}_\Null,\\
			\abs{M^*}_\Null &= \abs{M}_\Total - \abs{M}_\Wobs.
		\end{split}
	\end{equation}
	\item The cardinality of the reduction of $M$ is given by
	\begin{equation}
		\abs{M_\red} = \abs{M}_\Wobs - \abs{M}_\Null.
	\end{equation}
\end{corollarylist}
\end{corollary}

For finite embedded constraint index sets we will more suggestively write
$M + N$ for the disjoint union and $M \cdot N$ for their product.

\begin{remark}
The cardinality $\abs{\argument}$ yields a map from finite embedded constraint index sets to
$\Con\Naturals_0^3 \coloneqq \{(n_\Total, n_\Wobs, n_\Null) \in \Naturals_0^3 \mid n_\Null \leq n_\Wobs \leq n_\Total\}$,
and isomorphic constraint index sets have the same cardinality.
Conversely, to every $n \in \Con\Naturals_0^3$ we can associate the finite embedded constraint index set
$(\{1, \dotsc, n_\Null\}, \{1,\dotsc, n_\Wobs\},\{1, \dotsc, n_\Total\})$.
We will often use this identification implicitly and for example write
$k \in n_\Wobs$ instead of $k \in \{1, \dotsc, n_\Wobs\}$.
In particular, when we apply the above constructions to triples of natural numbers
this means we apply them to their associated finite embedded constraint index sets,
e.g. $n_\red = n_\Wobs \setminus n_\Null = \{n_\Null+1, \dotsc, n_\Wobs\}$.
\end{remark}

In contrast to constraint sets, the reduction of constraint index sets 
does \emph{not} commute with forgetting algebraic 
structure.
This is not surprising since forgetting to constraint index sets also 
forgets the equivalence relation needed for reduction.
Nevertheless, as we will see soon, when considered as bases of constraint vector spaces
the reduction of constraint index sets is compatible and yields the correct basis of the reduced 
space.

\begin{remark}
	\label{rem:ConSetVSConIndSet}
	Note that given a constraint set $M = (M_\Total, M_\Wobs, \sim_M)$
	we can construct a constraint index set out of it:
	Choose (using the axiom of choice) a splitting
	$s_M \colon M_\red \to M_\Wobs$ of the projection
	$\pr_M \colon M_\Wobs \to M_\red$,
	then $M' \coloneqq (M_\Total,M_\Wobs, M_\Wobs \setminus \image s_M)$
	is a constraint index set with $M'_\red \simeq M_\red$.
	This procedure is in general not functorial, since this would require 
	a coherent choice of splitting for all constraint sets.
	Moreover, there is in general no way to reconstruct the equivalence 
	relation $\sim_M$ from $M'$.
	Thus $\ConSet$ and $\ConIndSet$ should not be considered equivalent.
\end{remark}

In the category $\Sets$ of sets the axiom of choice is equivalent to 
the statement that every epimorphism splits.
Even though we assume the axiom of choice to hold in $\Sets$,
this does not imply an equivalent statement about the splitting
of regular or plain epimorphisms in $\ConIndSet$.

\begin{example}\
	\begin{examplelist}
		\item 	Consider the constraint embedded index sets
		\begin{equation}
			M = \big(\{1,2\},\, \{1,2\},\, \{1\}\big)
			\qquad\text{ and }\qquad
			N = \big(\{1,2\},\,\{1,2\},\, \{1,2\}\big).
		\end{equation}
		Then $f = (\id, \id) \colon M \to N$ is an epimorphism,
		but it does not split since
		$\id \colon N_\Wobs \to M_\Wobs$ does not preserve the $\NULL$-component.
		\item Consider the constraint index sets 
		\begin{equation}
			M = \big(\{1,2\},\,\{1,2\},\, \{1,2\}\big)
			\qquad\text{ and }\qquad
			N = \big(\{1\},\, \{1,2\},\, \{1,2\} \big).
		\end{equation}
		Then $f = (1, \id) \colon M \to N$ is a regular epimorphism
		but it does not split, since there exists no constraint morphism
		extending
		$\id \colon N_\Wobs \to M_\Wobs$.
	\end{examplelist}
\end{example}

It turns out that constraint index sets for which every regular epimorphism
into them splits are exactly the embedded constraint index sets, cf. \autoref{prop:ProjectiveConSets}.

\begin{proposition}
	\label{prop:ProjectiveIndexSets}
	\index{projective!constraint index set}
Let $P \in \ConIndSet$ be a constraint index set.  Then the
following statements are equivalent:
\begin{propositionlist}
	\item\label{prop:ProjectiveIndexSets_1} Every regular epimorphism
	$M \to P$ splits.
	\item\label{prop:ProjectiveIndexSets_2} For every regular
	epimorphism $\Phi \colon M \to N$ and every morphism
	$\Psi \colon P \to N$ there exists a morphism
	$\chi \colon P \to M$ such that $\Phi \circ \chi = \Psi$.
	\item\label{prop:ProjectiveIndexSets_3} We have $P \in \injConIndSet$.
\end{propositionlist}
\end{proposition}

\begin{proof}
The proof is completely analogous to the one of \autoref{prop:ProjectiveConSets}.
\end{proof}

\subsection{Constraint Vector Spaces}
\label{sec:ConVectorSpaces}

Let $\field{K}$ be a field.
We want to study (embedded) constraint $\field{K}$-vector spaces in this section.
On the one hand, these will give us a first impression about what kind of effects we can expect for
free and projective constraint modules over more general constraint rings or algebras.
On the other hand, these constraint vector spaces will describe the pointwise structure given by
constraint vector bundles, which will be introduced in \autoref{sec:ConVectorBundles}.
For simplicity, we define constraint vector spaces to be embedded from the start.

\begin{definition}[Constraint vector space]
	\label{def:ConVectSpace}
	\index{constraint!vector space}
Let $\field{K}$ be a field.
\begin{definitionlist}
	\item An embedded constraint $\field{K}$-module is called 
	\emph{constraint vector space} over $\field{K}$.
	\item The category of constraint $\field{K}$-vector spaces
	is denoted by
	\glsadd{ConVectK}$\ConVect_\field{K}$.	
\end{definitionlist}	
\end{definition}

Thus a constraint vector space $V$ simply consists
of a $\field{K}$-vector space $V_\Total$ together with subspaces
$V_\Null \subseteq V_\Wobs \subseteq V_\Total$.
It is now easy to see that every constraint vector space is free in the following sense:

\begin{proposition}[Constraint vector spaces are free]
	\label{prop:ConVectSpacesAreFree}
Every constraint $\field{K}$-vector space is free, i.e. there exists a constraint subset
$i \colon B \hookrightarrow V$ such that for every constraint map
$\phi \colon B \to W$ there exists a unique 
linear constraint map
$\Phi \colon V \to W$ such that
$\Phi \circ i = \phi$.
\end{proposition}

\begin{proof}
	Choose a basis for $V_\Null$ and extend it successively to $V_\Wobs$ and $V_\Total$.
\end{proof}

We will call such a constraint subset $i \colon B \hookrightarrow V$ a
\index{basis of constraint vector space}
\emph{constraint basis}
of $V$.
Since for vector spaces the cardinality of all bases agree,
the same is true for constraint vector spaces, allowing us to define 
the dimension of a constraint vector space by
\index{dimension!constraint vector space}
\glsadd{dim}
\begin{equation}
	\dim(V) \coloneqq \big(\dim(V_\Total),\, \dim(V_\Wobs),\, \dim(V_\Null)\big).
\end{equation}
As usual we call $V$ \emph{finite dimensional} if $\dim(V)$ is a 
finite constraint index set.

\begin{example}
	\label{ex:FreeConVectSpace}
For $n_\Null \leq n_\Wobs \leq n_\Total \in \Naturals_\Null$
there is a constraint vector space
$\Reals^n \coloneqq (\Reals^{n_\Total}, \Reals^{n_\Wobs}, \Reals^{n_\Null})$.
By \autoref{prop:ConVectSpacesAreFree} every finite dimensional constraint vector
space is of this form.
\end{example}

Let us quickly recall some constructions for constraint vector spaces,
known already from constraint modules.
Since we only consider $\field{K}$-vector spaces, we drop the index for the tensor products.

\begin{proposition}
	\label{prop:ConstuctionsCVectSpaces}
Let $V,W \in \ConVect_\field{K}$ be finite dimensional constraint vector spaces
and let $B_V$ and $B_W$ be constraint bases of $V$ and $W$, respectively.
\begin{propositionlist}
	\item \label{prop:ConstuctionsCVectSpaces_1}
	\index{direct sum!constraint vector spaces}
	$B_V \sqcup B_W$ is a constraint basis for 
	\glsadd{directSum}
	\begin{equation}
		V \oplus W = (V_\Total \oplus W_\Total, \,\, V_\Wobs \oplus W_\Wobs ,\,\, V_\Null \oplus W_\Null)
	\end{equation}
	and we have $\dim(V \oplus W) = \dim(V) + \dim(W)$.
	\item \label{prop:ConstuctionsCVectSpaces_2}
	\index{tensor product!constraint vector spaces}
	$B_V \tensor B_W$ is a constraint basis for
	\glsadd{tensor}
	\begin{equation}
		V \tensor W = (V_\Total \tensor W_\Total ,\,\, V_\Wobs \tensor W_\Wobs,\,\, V_\Wobs \tensor W_\Null + V_\Null \tensor W_\Wobs)
	\end{equation}
	and we have $\dim(V \tensor W) = \dim(V) \tensor \dim(W)$.
	\item \label{prop:ConstuctionsCVectSpaces_3}
	\index{strong tensor product!constraint vector spaces}
	$B_V \strtensor B_W$ is a constraint basis for
	\glsadd{strtensor}
	\begin{equation}
		V \strtensor W = \big(V_\Total \tensor W_\Total,\,\, 
		(V_\Wobs \tensor W_\Wobs) + (V_\Total \tensor W_\Null) + (V_\Null \tensor W_\Total),\,\,
		(V_\Total \tensor W_\Null) + (V_\Null \tensor W_\Total)\big)
	\end{equation}
	and we have $\dim(V \strtensor W) = \dim(V) \strtensor \dim(W)$.
	\item  \label{prop:ConstuctionsCVectSpaces_4}
	\index{dual!constraint vector space}
	$(B_V)^*$, i.e. the constraint dual set of $B_V$, is a constraint basis for
	\begin{equation}
		V^* = \ConHom_\field{K}(V,\field{K})
		= \big((V_\Total)^*,\, \Ann_{V_\Total}(V_\Null),\, \Ann_{V_\Total}(V_\Wobs)\big),
	\end{equation}
	where
	\glsadd{annihilator}$\Ann_{V_\Total}(V_\Null)$ and $\Ann_{V_\Total}(V_\Wobs)$
	denote the annihilators of $V_\Null$ and $V_\Wobs$ considered as subspaces of $(V_\Total)^*$
	and we have $\dim(V^*) = \dim(V)^*$.
\end{propositionlist}
\end{proposition}

\begin{proof}
These are all simple checks.
For \ref{prop:ConstuctionsCVectSpaces_4} recall that by the definition of constraint internal hom we have
\begin{equation*}
	(V^*)_\Wobs = \{ \alpha \in (V_\Total)^* \mid \alpha(V_\Wobs) \subseteq \field{K},\, \alpha(V_\Null) \subseteq 0\}
	= \Ann_{V_\Total}(V_\Null)
\end{equation*}
since $\field{K} = (\field{K}, \field{K}, 0)$ as a constraint vector space.
Similarly,
\begin{equation*}
	(V^*)_\Null = \{ \alpha \in (V_\Total)^* \mid \alpha(V_\Wobs) \subseteq 0\}
	= \Ann_{V_\Total}(V_\Wobs).
\end{equation*}
Note that we use the identification of
$\dim(V)$ and $\dim(W)$ with finite embedded constraint index sets
as well as their compositions from \autoref{def:TensorProdDualConIndSet}.
\end{proof}

In the following we check some of the well-known compatibilities of dualizing with the different
notions of tensor products.
\filbreak
\begin{proposition}
	\label{prop:CompatibilityDualTensorConVecSpaces}
Let $V,W \in \ConVect_\field{K}$ be finite dimensional constraint vector spaces.
\begin{propositionlist}
	\item \label{prop:CompatibilityDualTensorConVecSpaces_1}
	We have canonically
	\begin{equation}
		(V \oplus W)^* \simeq V^* \oplus W^*
	\end{equation}
	and therefore $\dim\!\big((V \oplus W)^*\big) = \dim(V)^* + \dim(W)^*$.
	\item \label{prop:CompatibilityDualTensorConVecSpaces_2}
	We have canonically
	\begin{equation} \label{eq:CompatibilityDualTensorConVecSpaces_DualOfTensor}
		(V \tensor W)^* \simeq V^* \strtensor W^*
	\end{equation}
	and therefore $\dim\!\big((V \tensor W)^*\big) = \dim(V)^* \strtensor 
	\dim(W)^*$.
	\item \label{prop:CompatibilityDualTensorConVecSpaces_3}
	We have canonically
	\begin{equation}
		(V \strtensor W)^* \simeq V^* \tensor W^*
	\end{equation}
	and therefore $\dim\!\big((V \strtensor W)^*\big) = \dim(V)^* \tensor 
	\dim(W)^*$.
	\item \label{prop:CompatibilityDualTensorConVecSpaces_4}
	We have canonically
	\begin{equation}
		\ConHom(V,W) \simeq W \strtensor V^*
	\end{equation}
	and therefore $\dim\!\big(\ConHom(V,W)\big) = \dim(W) \strtensor \dim(V)^*$.
\end{propositionlist}
\end{proposition}

\begin{proof}
Except for \ref{prop:CompatibilityDualTensorConVecSpaces_4}
these can be shown by choosing constraint dual bases of $V$ and $W$
and the use of \autoref{prop:ConstuctionsCVectSpaces}.
Then the dimensions follow from \autoref{cor:CardinalityOfConstructionsConIndSet}.
For the last part note that for $w \tensor \alpha \in W_\Total \tensor V^*_\Null + W_\Null \tensor V^*_\Total$
we have
$(w \tensor \alpha)(v) = w \cdot \alpha(v) \in W_\Null$ for all $\alpha \in V_\Wobs$.
With this it is easy to see that $B_W \strtensor (B_V)^*$ is a basis for
$\ConHom(V,W)$.
\end{proof}

This result shows that the two tensor products $\tensor$ and $\strtensor$ are intimately linked.
In particular, \eqref{eq:CompatibilityDualTensorConVecSpaces_DualOfTensor}
shows that we could have deduced $\strtensor$ from $\tensor$ by
\emph{defining} $V \strtensor W \coloneqq (V^* \tensor W^*)^*$,
at least in the finite-dimensional case.
Moreover, by definition, these tensor products are related by an injective morphism
\begin{equation}
	V \tensor W \hookrightarrow V \strtensor W,
\end{equation}
and they distribute in the sense that there exists a morphism
\begin{equation} \label{eq:RebracketingConTensors}
	U \tensor (V \strtensor W) \to (U \tensor V) \strtensor W.
\end{equation}
Both morphisms are not isomorphisms in general, as the next example shows:

\begin{example}
Consider the constraint vector space $\Reals^{n} = (\Reals^3,\Reals^2,\Reals^1)$ from
\autoref{ex:FreeConVectSpace} with $n = (3,2,1)$.
\begin{examplelist}
	\item 
	Then it holds that
	\begin{align}
		\dim(\Reals^n \tensor \Reals^n) &= n \tensor n = (9,4,3)
		\shortintertext{and}
		\dim(\Reals^n \strtensor \Reals^n) &= n \strtensor n = (9,6,5)
	\end{align}
	by \autoref{prop:ConstuctionsCVectSpaces} and \autoref{cor:CardinalityOfConstructionsConIndSet}.
	Thus $\Reals^n \tensor \Reals^n$ and $\Reals^n \strtensor \Reals^n$ cannot be isomorphic.
	\item We have
	\begin{align}
		\dim(\Reals^n \tensor (\Reals^n \strtensor \Reals^n))
		&= n \tensor (n \strtensor n)
		= (27, 12,11)
		\shortintertext{and}
		\dim\!\big((\Reals^n \tensor \Reals^n) \strtensor \Reals^n\big)
		&= (n \tensor n) \strtensor n= (27, 16,15).
	\end{align}
	This shows that $\Reals^n \tensor (\Reals^n \strtensor \Reals^n)$
	and $(\Reals^n \tensor \Reals^n) \strtensor \Reals^n$ are not isomorphic.
\end{examplelist}
\end{example}

\begin{remark}
The relations between $\tensor$ and $\strtensor$ can be derived from the fact that
$\ConVect_\field{K}$ together with $\tensor$ and $\argument^*$ forms a $*$-autonomous category, see
\cite{barr:1979a}.
\end{remark}

%% file: constraint-algebra-modules.tex
In this section we will define our main objects of interest: 
constraint algebras and their modules.
Following our philosophy to construct new constraint notions as objects internal to 
certain constraint categories, we will define constraint algebras as monoids internal 
to the categories of constraint $\field{k}$-modules introduced in \autoref{sec:ConAlgebras}.
Since $\ConMod_{\field{k}}$ carries two different monoidal structures, this will lead to the definitions
of constraint algebras in \autoref{sec:ConAlgebras} and strong constraint algebras in \autoref{sec:StrConAlgebras}.
In both cases we can also consider the subcategories of embedded constraint $\field{k}$-modules
which lead to embedded (strong) constraint algebras.
In these sections we will also introduce the corresponding notion of (strong) constraint modules over
(strong) constraint algebras.

\subsection{Constraint Algebras and their Modules}
\label{sec:ConAlgebras}

The following definition is just a reformulation of monoids internal to
$(\ConMod_{\field{k}},\tensor[\field{k}])$,
cf. \hyperref[sec:MonoidsModules]{Appendix~\ref{sec:MonoidsModules}} for the definition of monoids internal to a monoidal category.

\begin{definition}[Constraint algebra]\
	\label{def:ConstraintAlgebras}%
	\glsadd{conAlgebrasFont}
\begin{definitionlist}
	\item \index{constraint!algebra}
	A \emph{constraint $\field{k}$-algebra}
	is a triple
	$\algebra{A} = (\algebra{A}_\Total, \algebra{A}_\Wobs,
	\algebra{A}_\Null)$ consisting of unital associative 
	$\field{k}$-algebras
	$\algebra{A}_\Total$ and $\algebra{A}_\Wobs$ together with a two-sided
	ideal $\algebra{A}_\Null \subseteq \algebra{A}_\Wobs$ and
	a unital algebra homomorphism
	$\iota \colon \algebra{A}_\Wobs \to \algebra{A}_\Total$.
	\item A \emph{morphism
	$\phi \colon \algebra{A} \to \algebra{B}$
	of constraint $\field{k}$-algebras}
	is given by a pair of unital algebra homomorphisms
	$\phi_\Total \colon \algebra{A}_\Total \to \algebra{B}_\Total$
	and
	$\phi_\Wobs \colon \algebra{A}_\Wobs \to \algebra{B}_\Wobs$
	such that
	$\iota_\algebra{B} \circ \phi_\Wobs = \phi_\Total \circ
	\iota_\algebra{A}$ and
	$\phi_\Wobs(\algebra{A}_\Null) \subseteq \algebra{B}_\Null$.
	\item The category of constraint $\field{k}$-algebras is denoted by
	\glsadd{ConAlgk}$\ConAlg_\field{k}$.
\end{definitionlist}
\end{definition}

When the underlying ring $\field{k}$ is clear from context, we will write simply
$\ConAlg$ for the category of constraint algebras.

\begin{example}
	\label{ex:ConAlgebras}
Let $M \in \ConSet$ be a constraint set.
\begin{examplelist}
	\item \label{ex:FunctionsOnCSetI}
	\index{constraint!functions on a set}
	Consider the ring $\field{k}$
	as a constraint set $(\field{k},\field{k},\sim_\discrete)$.
	Then $\ConMap(M,\field{k})$ is a constraint algebra given by
	\begin{equation}
		\begin{split}
			\ConMap(M,\field{k})_\Total
			&= \Map(M,\field{k}), \\
			\ConMap(M,\field{k})_\Wobs
			&= \left\{ f \in \Map(M_\Total,\field{k}) \mid 
			f(\iota_M(x)) = f(\iota_M(y)) \text{ for all } x \sim_M y \right\},\\
			\ConMap(M,\field{k})_\Null
			&= \big\{ f \in \Map(M_\Total,\field{k}) \bigm| 
			f\at{\image \iota_M} = 0 \big\}.
		\end{split}
	\end{equation}
	\item Every ring $\field{k}$ can be seen as a constraint algebra
	$\field{k} = (\field{k}, \field{k},0)$.
	\item \index{constraint!endomorphisms}
	Let $V$ be a constraint vector space over a field $\field{K}$.
	Then
	\glsadd{ConEnd}$\ConEnd(V)$ is a constraint algebra with respect to composition of constraint morphisms.
\end{examplelist}
\end{example}

Since $\ConMod_{\field{k}}$ is a symmetric monoidal category, we can define commutative constraint algebras.
In this case we can define a constraint version of the center.

\begin{proposition}[Constraint center]
	\label{prop:ConCenter}
	\index{center}
Let $\field{k}$ be a commutative ring and let $\algebra{A} \in \ConAlg$ be a constraint algebra.
Then
\glsadd{ConCenter}$\Center(\algebra{A})$, defined by
\begin{equation}
\begin{split}
	\Center(\algebra{A})_\Total
	&\coloneqq \Center(\algebra{A}_\Total),\\
	\Center(\algebra{A})_\Wobs
	&\coloneqq \Center(\algebra{A}_\Wobs),\\
	\Center(\algebra{A})_\Null
	&\coloneqq \Center(\algebra{A}_\Wobs) \cap \algebra{A}_\Null,
\end{split}
\end{equation}
is a commutative constraint algebra, the \emph{center} of $\algebra{A}$.
\end{proposition}

From an algebraic point of view it is natural to study modules 
over constraint algebras.
Since \autoref{ex:ConAlgebras} \ref{ex:FunctionsOnCSetI} shows
that constraint algebras encode 
the algebraic structure of functions on a space allowing for 
reduction and vector bundles should correspond to certain modules, 
the next definition is also interesting from a geometric standpoint.

\begin{definition}[Modules over constraint algebras]
	\label{def:ConstraintModules}%
	Let $\algebra{A}, \algebra{B} \in \ConAlg$
	be constraint algebras.
	\begin{definitionlist}
		\item \index{constraint!right module}
		A \emph{constraint right $\algebra{A}$-module}
		is a constraint $\field{k}$-module 
		$\module{E} = (\module{E}_\Total, \module{E}_\Wobs,
		\module{E}_\Null)$ with a right $\algebra{A}_\Total$-module 
		structure on $\module{E}_\Total$
		and a right $\algebra{A}_\Wobs$-module structure on 
		$\module{E}_\Wobs$ such that
		$\module{E}_\Null \subseteq \module{E}_\Wobs$
		is an $\algebra{A}_\Wobs$-submodule,
		$\iota_\module{E} \colon \module{E}_\Wobs \to \module{E}_\Total$
		is an $\algebra{A}_\Wobs$-module morphism and
		$\module{E}_\Wobs \cdot \algebra{A}_\Null \subseteq 
		\module{E}_\Null$.
		\item \index{constraint!left module}
		A \emph{constraint left $\algebra{B}$-module}
		is a constraint $\field{k}$-module 
		$\module{E} = (\module{E}_\Total, \module{E}_\Wobs, 
		\module{E}_\Null)$
		with a left $\algebra{B}_\Total$-module 
		structure on $\module{E}_\Total$
		and a left $\algebra{B}_\Wobs$-module structure on 
		$\module{E}_\Wobs$ such that
		$\module{E}_\Null \subseteq \module{E}_\Wobs$
		is an $\algebra{B}_\Wobs$-submodule,
		$\iota_\module{E} \colon \module{E}_\Wobs \to \module{E}_\Total$
		is an $\algebra{B}_\Wobs$-module morphism and
		$\module{B}_\Null \cdot \algebra{E}_\Wobs \subseteq 
		\module{E}_\Null$.
		\item \index{constraint!bimodule}
		A \emph{constraint $(\algebra{B},\algebra{A})$-bimodule}
		is a constraint $\field{k}$-module 
		$\module{E} = (\module{E}_\Total, \module{E}_\Wobs,
		\module{E}_\Null)$ with commuting constraint left $\algebra{B}$- 
		and right $\algebra{A}$-module structures. 
		\item A \emph{morphism 
			$\Phi\colon \module{E} \to \module{F}$
			between constraint left-/right-/bi-modules} is a pair
		$(\Phi_\Total, \Phi_\Wobs)$
		of left-/right-/bi-module morphisms
		$\Phi_\Total \colon \module{E}_\Total \to \module{F}_\Total$
		and
		$\Phi\colon \module{E}_\Wobs \to \module{F}_\Wobs$
		such that
		$\Phi_\Total \circ \iota_\module{E} = \iota_{\module{F}} \circ 
		\Phi_\Wobs$
		and
		$\Phi_\Wobs(\module{E}_\Null) \subseteq \module{F}_\Null$.
		\item The categories of	constraint right $\algebra{A}$-modules,
		left $\algebra{B}$-modules and
		$(\algebra{B}, \algebra{A})$-bimodules are denoted by
		\glsadd{ConRModA}$\ConRMod{\algebra{A}}$,
		\glsadd{ConLModB}$\ConLMod{\algebra{B}}$ 
		and
		\glsadd{ConBimodBA}$\ConBimod(\algebra{B},\algebra{A})$, respectively.
	\end{definitionlist}
\end{definition}

Again, this definition can also be understood as modules internal to 
the monoidal category $(\ConMod_{\field{k}}, \tensor[\field{k}])$.
As we would expect, right $\algebra{A}$-modules can be understood as 
$(\field{k},\algebra{A})$-bimodules, writing again $\field{k} = 
(\field{k},\field{k},0)$, and similarly for left modules.
Moreover, constraint $\field{k}$-modules as defined in 
\autoref{sec:ConstraintkModules} are nothing but constraint 
$(\field{k},\field{k})$-bimodules.

We will not go into details of the construction of limits and colimits for 
modules over constraint algebras here.
Suffice to say that the underlying constraint $\field{k}$-modules are 
given by the corresponding construction from 
\autoref{sec:ConstraintkModules} and the $\algebra{A}$-module 
structures on the respective components are the obvious ones.
Since the tensor product of constraint modules over a constraint 
algebra will be important later on, we spell it out in detail.

\begin{proposition}[Tensor product of constraint modules]
	\label{lemma:TensorProductBimodules}%
	\index{tensor product!constraint bimodules}
Let  $\algebra{A}, \algebra{B}, \algebra{C} \in \ConAlg$
be given and let
$\module{F} \in \ConBimod(\algebra{C},\algebra{B})$
as well as
$\module{E} \in \ConBimod(\algebra{B},\algebra{A})$
be constraint bimodules.
Then the constraint $(\algebra{C},\algebra{A})$-bimodule
$\module{F} \tensor[\algebra{B}] \module{E}$
is given by
\glsadd{tensorA}
\begin{equation}
	\begin{split}
		\left(
		\module{F}
		\tensor[\algebra{B}]
		\module{E}
		\right)_{\Total}
		&=
		\module{F}_{\Total}
		\tensor[\algebra{B}_{\Total}]
		\module{E}_{\Total}, \\
		\left(
		\module{F}
		\tensor[\algebra{B}]
		\module{E}
		\right)_{\Wobs}
		&=
		\module{F}_{\Wobs}
		\tensor[\algebra{B}_{\Wobs}]
		\module{E}_{\Wobs}, \\
		\left(
		\module{F}
		\tensor[\algebra{B}]
		\module{E}
		\right)_{\Null}
		&=
		\module{F}_{\Wobs}
		\tensor[\algebra{B}_\Wobs]
		\module{E}_{\Null}
		+
		\module{F}_{\Null}
		\tensor[\algebra{B}_\Wobs]
		\module{E}_\Wobs,
	\end{split}
\end{equation}
with $\iota_{\tensor} = \iota_\module{F} \tensor \iota_\module{E}$.
\end{proposition}

\begin{proof}
	Denote by
	$\lambda \colon \algebra{B} \tensor[\field{k}] \module{E} \to 
	\module{E}$
	the left $\algebra{B}$-multiplication on $\module{E}$ and by
	$\rho \colon \module{F} \tensor[\field{k}] \algebra{B} \to 
	\module{F}$
	the right $\algebra{B}$-multiplication on $\module{F}$.
	Then $\module{F} \tensor[\algebra{B}] \module{E}$
	is defined as the coequalizer of
	$\id_\module{F} \tensor \lambda$ and
	$\rho \tensor \id_\module{E}$
	as constraint morphisms from
	$\module{F} \tensor[\tensor_\field{k}] \algebra{B} 
	\tensor[\field{k}] \module{E}$
	to 
	$\module{F} \tensor[\field{k}] \module{E}$.
	Applying \autoref{prop:MonoidalStructureCModk}
	and \autoref{prop:CoLimitsConModk} 
	\ref{prop:CoLimitsConModk_coequalizer}
	gives the desired result.
\end{proof}

With this tensor product we can construct a bicategory of constraint 
modules analogous to the classical bicategory of bimodules, see
\cite{johnson.yau:2021a}
for a modern treatment of bicategories.

\begin{proposition}[The bicategory $\ConBimod{}{}$]
	Using constraint algebras as objects, constraint bimodules as 
	$1$-morphisms, morphisms of constraint bimodules as $2$-morphisms 
	and the tensor product of constraint modules as composition 
	defines a bicategory
	\glsadd{ConBimod}$\ConBimod$.
\end{proposition}

\begin{remark}
	\label{rem:MoritaEquivalence}
	\index{Morita equivalence}
In classical algebra two algebras $\algebra{A}$ and $\algebra{B}$ are considered to be Morita equivalent if their respective categories of representations, i.e. their categories of (right-)modules, are equivalent.
This can then be reformulated to the fact that $\algebra{A}$ and $\algebra{B}$ are equivalent internal to $\Bimodules$, meaning that there exists an invertible $1$-morphism between $\algebra{A}$ and $\algebra{B}$.
It turns out that these invertible $1$-morphisms are given by finitely generated projective full
$(\algebra{B},\algebra{A})$-bimodules.
The bicategory $\ConBimod$ now opens up a way to study Morita theory of constraint algebras
by defining constraint algebras to be Morita equivalent if they are equivalent internal to
$\ConBimod$.
In order to characterize constraint Morita equivalence bimodules it seems reasonable to study 
finitely generated projective constraint modules first.
Even though we will not be concerned with Morita theory this can be seen as a motivation for
\autoref{sec:RegularProjectiveModules}.
The Morita theory of a subcategory of constraint algebras has been studied in 
\cite{dippell:2018a,dippell.esposito.waldmann:2019a} under the name of Morita equivalence for coisotropic algebras.
There you can also find a more detailed construction of 
$\ConBimod$.
\end{remark}

The internal hom of constraint $\field{k}$-modules carries over to a 
constraint module structure on the homomorphisms of constraint 
modules over constraint algebras.

\begin{proposition}[Module structure on module morphisms]
	\label{prop:ModuleStructureOnModuleHom}
	\index{internal hom!constraint modules}
Let $\algebra{A}$, $\algebra{B}$ and $\algebra{C}$
be constraint algebras and let
$\module{E} \in \ConBimod(\algebra{B},\algebra{A})$
as well as $\module{F} \in \ConBimod(\algebra{C},\algebra{A})$.
Then the right $\algebra{A}$-module morphisms form a constraint 
$(\algebra{C},\algebra{B})$-bimodule given by
\glsadd{ConHomA}
\begin{equation}
	\begin{split}
		\ConHom_{\algebra{A}}(\module{E},\module{F})_\Total
		&\coloneqq \Hom_{\algebra{A}_\Total}(\module{E}_\Total, 
		\module{F}_\Total),\\
		\ConHom_{\algebra{A}}(\module{E},\module{F})_\Wobs
		&\coloneqq \big\{ (\Phi_\Total, \Phi_\Wobs) \in 
		\Hom_{\algebra{A}_\Total}(\module{E}_\Total, 
		\module{F}_\Total) \times 
		\Hom_{\algebra{A}_\Wobs}(\module{E}_\Wobs, \module{F}_\Wobs)\mid\\
		&\qquad\qquad\qquad\qquad\Phi_\Total \circ \iota_\module{E} = \iota_\module{F} 
		\circ \Phi_\Wobs \text{ and } \Phi_\Wobs(\module{E}_\Null) 
		\subseteq \module{F}_\Null \big\},\\
		\ConHom_{\algebra{A}}(\module{E},\module{F})_\Null
		&\coloneqq \left\{ \Phi \in 
		\ConHom_{\algebra{A}}(\module{E},\module{F})_\Wobs \mid 
		\Phi_\Wobs(\module{E}_\Wobs) \subseteq \module{F}_\Null 
		\right\}.
	\end{split}
\end{equation}
\end{proposition}

With this proposition it is clear that the categories 
$\ConLMod{\algebra{B}}$ and $\ConRMod{\algebra{A}}$
of constraint left and right modules are enriched over 
$\ConMod_{\field{k}}$.
Moreover, we can define duals for constraint modules.

\begin{definition}[Dual module]
	\label{def:DualModule}
	\index{dual module}
Let $\algebra{A} \in \ConAlg$ and $\module{E} \in 
\ConRMod{\algebra{A}}$.
We call the constraint left $\algebra{A}$-module
\glsadd{DualModule}$\module{E}^* \coloneqq \ConHom_{\algebra{A}}(\module{E},\algebra{A})$
the \emph{dual module of $\module{E}$}.
\end{definition}

To give a first example of a constraint module over a constraint 
algebra we introduce the notion of derivations of constraint algebras.

\begin{definition}[Derivation]
	\label{def:ConDerivation}%
	\index{derivations}
Let $\algebra{A} \in \ConAlg$ be a constraint algebra and let 
$\module{M} \in \ConBimod(\algebra{A},\algebra{A})$ be
an $\algebra{A}$-bimodule. 
A \emph{derivation with values in} $\module{M}$ is a morphism
$D \colon \algebra{A} \longrightarrow \module{M}$ of constraint
$\field{k}$-modules such that
\begin{equation}
	\label{eq:Derivation}
	D \circ \mu
	=
	\ell \circ (\id \tensor D) + r \circ (D \tensor \id)
\end{equation}
holds, where $r$ and $\ell$ denote the right and left
$\algebra{A}$-multiplications of $\module{M}$, respectively, 
and $\mu$ is the multiplication of $\algebra{A}$.
The set of derivations will be denoted by
\glsadd{DerA}$\Der(\algebra{A},\module{M})$.
If $\module{M} = \algebra{A}$ we write $\Der(\algebra{A})$.
\end{definition}

\begin{lemma}
Let $\algebra{A} \in \ConAlg$ be a constraint algebra and let 
$\module{M} \in \ConBimod(\algebra{A},\algebra{A})$ be
an $\algebra{A}$-bimodule. 
A derivation $D = (D_\Total,D_\Wobs)$ with values in $\module{M}$ is a morphism of 
constraint $\field{k}$-modules such that
\begin{equation}
	D_\Total(ab) = aD_\Total(b) + D_\Total(a)b
\end{equation}
holds for all $a,b \in \algebra{A}_\Total$ and 
\begin{equation}
	D_\Wobs(ab) = aD_\Wobs(b) + D_\Wobs(a)b
\end{equation}
holds for all $a,b \in \algebra{A}_\Wobs$.
\end{lemma}

\begin{proof}
This is exactly the componentwise evaluation of \eqref{eq:Derivation} on elements.
\end{proof}

We can arrange the constraint derivations as a constraint 
submodule of the internal homomorphism
$\ConHom_\field{k}(\algebra{A},\module{M})$ as follows.

\begin{proposition}[$\field{k}$-module of derivations]
	\label{prop:ConDer}%
	\index{derivations}
Let $\algebra{A} \in \ConAlg$ be a constraint algebra and let 
$\module{M} \in \ConBimod(\algebra{A},\algebra{A})$ be
a constraint $\algebra{A}$-bimodule.  Then
\glsadd{ConDerA}
\begin{equation}
	\begin{split}
		\ConDer(\algebra{A},\module{M})_\Total
		&:= \Der(\algebra{A}_\Total, \module{M}_\Total),\\
		\ConDer(\algebra{A},\module{M})_\Wobs
		&:=	\big\{(D_\Total,D_\Wobs) \in 
		\Hom_\field{k}(\algebra{A},\module{M})
		\; \big| \;
		D_\Total \in \Der(\algebra{A}_\Total,\module{M}_\Total),\\
		&\qquad\qquad\qquad\qquad\qquad\qquad\qquad
		D_\Wobs \in \Der(\algebra{A}_\Wobs, \module{M}_\Wobs) \big\},	\\
		\ConDer(\algebra{A},\module{M})_\Null
		&:=	\big\{(D_\Total, D_\Wobs ) \in 
		\Der(\algebra{A},\module{M})_\Wobs
		\; \big| \;
		D_\Wobs(\algebra{A}_\Wobs) \subseteq \module{M}_\Null
		\big\}
	\end{split}
\end{equation}
defines a constraint $\field{k}$-module
$\ConDer(\algebra{A},\module{M})$.
\end{proposition}

One needs to be careful with the notation, since
$\Der(\algebra{A})$ has different meanings depending on whether
$\algebra{A}$ is a constraint or a classical algebra.
Also note that
$\ConDer(\algebra{A},\module{M})_\Wobs = \Der(\algebra{A},\module{M})$
is just the set of derivations of a	constraint algebra $\algebra{A}$ 
with values in the constraint
module $\module{M}$ as given in \autoref{def:ConDerivation}.
The constraint $\field{k}$-module of derivations on $\algebra{A}$
with values in $\algebra{A}$ is denoted by 
$\ConDer(\algebra{A})$.

As for classical algebras the derivations turn out to be a bimodule 
if the algebra is commutative:

\begin{corollary}[$\algebra{A}$-module of derivations]
	\label{cor:DerivationsAsAModule}%
	Let $\algebra{A} \in \ConAlg$ be a commutative	constraint algebra.
	Then $\ConDer(\algebra{A})$ is a constraint $\algebra{A}$-bimodule.
\end{corollary}

\subsubsection{Embedded Constraint Algebras and their Modules}

The subcategory of constraint algebras $\algebra{A}$ with injective
$\iota_\algebra{A} \colon \algebra{A}_\Wobs \to \algebra{A}_\Total$
will be denoted by  \index{embedded constraint!algebra}
\glsadd{injConAlg}$\injConAlg$.

\begin{corollary}
Let $\algebra{A} \in \injConMod_{\field{k}}$ be a constraint module.
Then a monoid structure on $\algebra{A}$
internal to $(\injConMod_{\field{k}},\injtensor[\field{k}])$
is equivalently given by an algebra structure on $\algebra{A}_\Total$
such that $\algebra{A}_\Wobs \subseteq \algebra{A}_\Total$
is a subalgebra and $\algebra{A}_\Null \subseteq \algebra{A}_\Wobs$
is a two-sided ideal.
\end{corollary}

\begin{proof}
This is clear by the definition of $\injtensor[\field{k}]$ in \autoref{prop:CategoryInjCMod}.
\end{proof}

In other words, $\injConAlg$ is exactly the category 
of monoids internal to $\injConMod_{\field{k}}$ with tensor product
$\injtensor[\field{k}]$.
Again by \autoref{prop:CategoryInjCMod} it is easy to see that $\injConAlg$ is a reflective subcategory of
$\ConAlg$.
Unsurprisingly, we will call such algebras \emph{embedded}.

\begin{example}\
	\label{ex:injConAlg}
	\index{constraint!functions on a set}
\begin{examplelist}
	\item \label{ex:injConAlg_1}
	Let $M \in \ConSet$ be a constraint set.
	Then the constraint algebra $\Map(M,\field{k})$, as already considered in \autoref{ex:ConAlgebras}
	\ref{ex:FunctionsOnCSetI}, is embedded.
	This can be understood as a consequence of $\field{k} = (\field{k}, \field{k},0)$ being an embedded constraint algebra,
	see also \autoref{prop:CoPowerCSet}.
	\item \label{ex:injConAlg_2}
	Let $M \in \strConSet$ be a strong constraint set.
	Then $\strConMap(M,\field{k})$ is an embedded constraint algebra 
	given by
	\begin{equation}
		\begin{split}
			\strConMap(M,\field{k})_\Total
			&= \Map(M_\Total, \field{k}),\\
			\strConMap(M,\field{k})_\Wobs
			&= \left\{ f \in \Map(M_\Total, \field{k}) \mid 
			f(x) = f(y) \text{ for all } x \sim_M^\Total y \right\},\\
			\strConMap(M,\field{k})_\Null
			&= \left\{ f \in \strConMap(M,\field{k})_\Wobs \bigm|
			f\at{\image \iota_M} = 0 \right\}.
		\end{split}
	\end{equation} 
	\item \label{ex:injConAlg_3}
	Let $M \in \injstrConSet$ be an embedded strong constraint set with 
	inclusion $M_\Wobs \subseteq M_\Total$.
	Then $\strConMap(M,\field{k})_\Wobs$ is the subalgebra of 
	functions constant along the equivalence classes of 
	$\sim_M^\Total$
	and $\strConMap(M,\field{k})_\Null$ is the intersection of this 
	subalgebra with the vanishing ideal of $M_\Wobs$.
\end{examplelist}
\end{example}

\begin{remark}
In \cite{dippell.esposito.waldmann:2019a} so called coisotropic triples of algebras were considered.
These are embedded constraint algebras with the additional property
of $\algebra{A}_\Null$ being a left ideal in $\algebra{A}_\Total$.
Note, that $\strConMap(M,\field{k})$ from \autoref{ex:injConAlg} \ref{ex:injConAlg_3} is not of this form,
since $\strConMap(M,\field{k})_\Null$ is not an ideal in $\algebra{A}_\Total$.
\end{remark}

Considering embedded constraint bimodules leads to the bicategory
\glsadd{injConBimod}$\injConBimod$.
Since we will not need the full bicategory, let us stick to the 
following situation:

\begin{proposition}[The category $\injConBimod(\algebra{A})$]\
	\label{prop:CategoryInjCBimod}
	\index{embedded constraint!bimodule}
Let $\algebra{A} \in \injConAlg$ be given.
\begin{propositionlist}
	\item The category $\injConBimod(\algebra{A})$ is a reflective subcategory 
	of
	$\ConBimod(\algebra{A})$
	with reflector\linebreak
	$\cdot^\inj \colon \ConBimod(\algebra{A})
	\to \injConBimod(\algebra{A})$ given by
	\begin{equation}
		\module{E}^\inj \coloneqq (\module{E}_\Total, 
		\iota_\module{E}(\module{E}_\Wobs),
		\iota_\module{E}(\module{E}_\Null) ).
	\end{equation}
	\item The subcategory $\injConBimod(\algebra{A})$
	is closed under finite limits.
	\item $\injConBimod(\algebra{A})$ is closed monoidal with 
	respect to $\injtensor[\algebra{A}]$ defined by
	\glsadd{injTensorA}
	\begin{equation}
		\module{E} \injtensor[\algebra{A}] \module{F} \coloneqq 
		(\module{E} 
		\tensor[\algebra{A}] \module{F})^\inj.
	\end{equation}
	\item The functor $\argument^\inj \colon 
	(\ConBimod(\algebra{A}),\tensor[\algebra{A}])
	\to (\injConBimod(\algebra{A}), \injtensor[\algebra{A}])$
	is monoidal, and the functor
	$\functor{U} \colon 
	(\injConBimod(\algebra{A}),\injtensor[\algebra{A}]) 
	\to (\ConBimod(\algebra{A}), \tensor[\algebra{A}])$
	is lax monoidal.
\end{propositionlist}
\end{proposition}

\begin{proof}
	The proof is completely analogous to the one of 
	\autoref{prop:CategoryInjCMod}.
	More conceptually, one could carry over the monoidal adjunction 
	from \autoref{prop:CategoryInjCMod} to realize
	$\injConBimod$ as a reflective sub-bicategory of $\ConBimod$.
	Then $\injConBimod(\algebra{A})$ is automatically 
	a reflective subcategory of $\ConBimod(\algebra{A})$.
\end{proof}

\subsubsection{Reduction}

By the definition of constraint algebras
internal to the monoidal category 
$(\ConMod_{\field{k}},\tensor[\field{k}])$ together with the fact that
$\red \colon \ConMod_{\field{k}} \to \Modules_\field{k}$
is monoidal it induces a reduction 
functor
\index{reduction!algebras}
\glsadd{Alg}
\begin{equation}
	\red \colon \ConAlg \to \Algebras
\end{equation}
given by $\algebra{A}_\red = \algebra{A}_\Wobs / \algebra{A}_\Null$.

Similarly, we obtain an induced reduction functor on constraint bimodules.
For the sake of exposition let us spell this out.
\filbreak
\begin{proposition}[Reduction of constraint bimodules]\
	\label{prop:ReductionConBimod}
	\index{reduction!bimodules}
\begin{propositionlist}
	\item Let $\algebra{A}, \algebra{B} \in \ConAlg$ be constraint 
	algebras and $\module{E} \in \ConBimod(\algebra{B},\algebra{A})$.
	Then $\module{E}_\red = \module{E}_\Wobs / \module{E}_\Null$ is a 
	$(\algebra{B}_\red,\algebra{A}_\red)$-bimodule.
	\item Reduction defines a functor of bicategories
	\glsadd{Bimod}$\red \colon \ConBimod \to \Bimodules$
	to the bicategory of algebras and bimodules.	
	\item Let $\algebra{A} \in \ConAlg$ be a constraint algebra.
	The functor
	$\red \colon \ConMod_\algebra{A} \to \Modules_{\algebra{A}_\red}$
	is lax closed with injective natural 
	transformation
	$\red \circ \ConHom_\algebra{A} \Rightarrow \Hom_{\algebra{A}_\red} \circ (\red \times \red)$.
\end{propositionlist}
\end{proposition}

\begin{proof}
Since $\algebra{B}_\Null \cdot \module{E}_\Wobs \subseteq 
\module{E}_\Null$ and $\module{E}_\Wobs \cdot \algebra{A}_\Null 
\subseteq \module{E}_\Null$ hold by definition of a constraint 
bimodule, we get a well-defined 
$(\algebra{B}_\red,\algebra{A}_\red)$-bimodule structure on 
$\module{E}_\red$.
The proof of the second part can be found,
for the special case of embedded constraint algebras with
$\algebra{A}_\Null \subseteq \algebra{A}_\Total$
a left ideal,
in detail in \cite{dippell:2018a,dippell.esposito.waldmann:2019a}.
This proof directly carries over to our situation.
For the last part it is easy to see that 
there is a morphism
$\ConHom_\algebra{A}(\module{E},\module{F})_\Wobs \to \ConHom_{\algebra{A}_\red}(\module{E}_\red, \module{F}_\red)$
whose kernel is exactly given by $\ConHom_\algebra{A}(\module{E},\module{F})_\Null$,
cf. \autoref{prop:ReductionOnCModk}.
\end{proof}

The reduction of constraint left or right modules is then to be 
understood as a special case of reduction of bimodules.
In particular we get from \autoref{prop:ReductionConBimod}
also the existence of reduction functors
\begin{equation}
	\red \colon \ConRMod{\algebra{A}} \to \Modules_{\algebra{A}_\red}
\end{equation}
and
\begin{equation}
	\red \colon \ConLMod{\algebra{B}} \to \LModules{\algebra{B}_\red}.
\end{equation}

\begin{example}
	\label{ex:ReductionConDer}
	\index{reduction!derivations}
Let $\algebra{A} \in \ConAlg$ be given.
Every $(D_\Total, D_\Wobs) \in \Der(\algebra{A})_\Wobs$
defines a derivation on
$\algebra{A}_\red = \algebra{A}_\Wobs / \algebra{A}_\Null$
since the condition
$D_\Wobs(\algebra{A}_\Null) \subseteq \algebra{A}_\Null$ is
automatically satisfied.
Hence we have a $\field{k}$-linear map
$\ConDer(\algebra{A})_\Wobs \to \Der(\algebra{A}_\red)$.
The kernel of this linear map is exactly given by 
$\ConDer(\algebra{A})_\Null$,
thus there exists an injective module homomorphism
\begin{equation}
	\label{eq:ReducedDerivation}
	\ConDer(\algebra{A})_\red \hookrightarrow \Der(\algebra{A}_\red).
\end{equation}
This is simply the restriction of the canonical morphism
$\ConHom_\field{k}(\algebra{A},\algebra{A})_\red \hookrightarrow
\Hom_\field{k}(\algebra{A}_\red,\algebra{A}_\red)$
from \autoref{prop:ReductionOnCModk}
to the submodule $\ConDer(\algebra{A})_\red$.
\end{example}

\begin{example}
	Our notion of a constraint algebra generalizes and unifies
	previous notions used in non-commutative geometry referring to
	features of the derivations:
	\begin{examplelist}
		\item \label{item:SubmanifoldAlgebra}
		\index{submanifold algebra}
		A \emph{submanifold algebra} in the sense of \cite{masson:1996a} 
		and \cite{dandrea:2020a} can equivalently be described as a 
		constraint algebra $\algebra{A}$ with
		$\algebra{A}_\Total = \algebra{A}_\Wobs$
		such that the canonical module morphism
		\eqref{eq:ReducedDerivation} is an isomorphism.
		\item \label{item:QuotientManifoldAlgebra}
		\index{quotient manifold algebra}
		A \emph{quotient manifold algebra} in the sense of 
		\cite{masson:1996a} can	equivalently be described as a 
		constraint algebra $\algebra{A}$ with 
		$\algebra{A}_\Wobs \subseteq \algebra{A}_\Total$
		a subalgebra and $\algebra{A}_\Null = 0$ such that
		$\Center(\algebra{A}_\red) \simeq \Center(\algebra{A})_\red$,
		$\Der(\algebra{A}_\red) \simeq \ConDer(\algebra{A})_\red$
		via \eqref{eq:ReducedDerivation} and
		\begin{equation}
			\algebra{A}_\Wobs = \{ a \in \algebra{A}_\Total \mid 
			D_\Total(a) = 0 \textrm{ for all } (D_\Total,D_\Wobs) \in
			\ConDer(\algebra{A})_\Null\}
		\end{equation}
		holds.
		Here $\Center(\algebra{A})$ denotes the constraint center of the 
		constraint algebra $\algebra{A}$,
		see \autoref{prop:ConCenter}.
	\end{examplelist}
\end{example}

\subsection{Strong Constraint Algebras and their Modules}
\label{sec:StrConAlgebras}

We replace now the tensor product $\tensor[\field{k}]$ on $\ConMod_{\field{k}}$ by the strong tensor product
$\strtensor[\field{k}]$.
Even though in later chapters we will only need the embedded situation, let us, for conecptual reasons, quickly introduce non-embedded strong constraint algebras and modules.

\begin{definition}[Strong constraint algebra]\
	\label{def:strConAlg}
	\index{strong constraint!algebra}
\begin{definitionlist}
	\item A \emph{strong constraint algebra} is a monoid object 
	internal to the category 
	$\strConMod_{\field{k}}$
	equipped with the strong tensor product $\strtensor[\field{k}]$.
	\item The category of strong constraint algebras is denoted by
	\glsadd{strConAlg}$\strConAlg$.
\end{definitionlist}
\end{definition}

Despite being conceptually clear, we need to unwrap the definition in 
order to be able to actually work with it.
We expect a strong constraint algebra to resemble a constraint algebra $\algebra{A}$
with the additional property that $\algebra{A}_\Null$ behaves like a two-sided ideal in
$\algebra{A}_\Total$.
For embedded algebras this will be true, but if $\iota_\algebra{A} \colon \algebra{A}_\Wobs \to \algebra{A}_\Total$
is not injective, it will turn out to be more complicated.

\begin{proposition}
	Let $\algebra{A}, \algebra{B} \in \strConMod_\field{k}$ and
	$f \colon \algebra{A} \to \algebra{B}$ a morphism of constraint
	$\field{k}$-modules.
	\begin{propositionlist}
		\item The structure of a strong constraint algebra on
		$\algebra{A}$ is equivalently given by
		the following data:
		\begin{propertieslist}
			\item an algebra structure $(\mu_\Total, 1_\Total)$
			on $\algebra{A}_\Total$,
			\item an algebra structure 
			$(\mu_\Wobs^{\Wobs\Wobs}, 1_\Wobs^{\Wobs\Wobs})$
			on $\algebra{A}_\Wobs$,
			\item an $(\algebra{A}_\Total, \algebra{A}_\Total)$-bimodule
			structure $(\mu_\Wobs^{\Total\Null}, \mu_\Wobs^{\Null\Total})$
			on $\algebra{A}_\Null$,
		\end{propertieslist}
		such that
		\begin{propertieslist}[resume]
			\item $\iota_\algebra{A} \colon \algebra{A}_\Wobs \to 
			\algebra{A}_\Total$
			is an algebra homomorphism,
			\item $\iota_\algebra{A}\at{\algebra{A}_\Null} \colon \algebra{A}_\Null \to \algebra{A}_\Total$
			is a morphism of $(\algebra{A}_\Total, \algebra{A}_\Total)$-bimodules,
			\item the $(\algebra{A}_\Wobs,\algebra{A}_\Wobs)$-bimodule
			structure on $\algebra{A}_\Null$ defined by
			$(\mu_\Wobs^{\Total\Null} \circ (\iota_\algebra{A} \tensor 
			\id), \mu_\Wobs^{\Null\Total} \circ (\id \tensor 
			\iota_\algebra{A}))$
			agrees with the one defined by the restriction of
			$\mu_\Wobs^{\Wobs\Wobs}$.
		\end{propertieslist}
	\item The morphism $f$ being a morphism of strong constraint algebras is equivalent to
	the following properties:
	\begin{propertieslist}
		\item $f_\Total \colon \algebra{A}_\Total \to \algebra{B}_\Total$
		is an algebra homomorphism.
		\item $f_\Wobs \colon \algebra{A}_\Wobs \to \algebra{B}_\Wobs$
		is an algebra homomorphism.
		\item $f_\Wobs\at{\algebra{A}_\Null} \colon \algebra{A}_\Null \to \algebra{B}_\Null$
		is a morphism of $(\algebra{A}_\Total,\algebra{A}_\Total)$-bimodules.		
	\end{propertieslist}
	\end{propositionlist}
\end{proposition}

\begin{proof}
	Consider constraint $\field{k}$-module maps
	$\mu = (\mu_\Total, \mu_\Wobs) \colon \algebra{A} 
	\strtensor[\field{k}] \algebra{A} \to \algebra{A}$
	and 
	$1 \colon (\field{k},\field{k},0) \to \algebra{A}$.
	Then by \autoref{lem:MorphismsOnStrConModules}
	we know that $\mu$ is given by
	\begin{align*}
		\mu_\Total &\colon \module{A}_\Total
		\tensor[\field{k}] 
		\module{A}_\Total \to \module{A}_\Total, \quad &
		\mu_\Wobs^{\Wobs\Wobs} &\colon \module{A}_\Wobs 
		\tensor[\field{k}] 
		\module{A}_\Wobs \to \module{A}_\Wobs, \\
		\mu_\Wobs^{\Total\Null} &\colon \module{A}_\Total 
		\tensor[\field{k}] 
		\module{A}_\Null \to \module{A}_\Null, \quad &
		\mu_\Wobs^{\Null\Total} &\colon \module{A}_\Null 
		\tensor[\field{k}] 
		\module{A}_\Total \to \module{A}_\Null.
	\end{align*}
	Writing out the associativity and unit diagrams of 
	\autoref{def:Monoid} in terms of these maps we obtain
	a $\field{k}$-algebra structure on $\algebra{A}_\Total$ by considering
	$\mu_\Total$ and $1_\Total$.
	Similarly, $\mu_\Wobs^{\Wobs\Wobs}$ and $1_\Wobs^{\Wobs\Wobs}$
	yield the algebra structure on $\algebra{A}_\Wobs$
	and $\mu_\Wobs^{\Total\Null}$, $\mu_\Wobs^{\Null\Total}$
	give the right- and left module structure on $\algebra{A}_\Null$, respectively.
	The compatibilities are then required to turn everything into morphisms of constraint modules.
	The second part follows directly by spelling out \autoref{def:MorphismOfMonoids}
	in terms of the different components.
\end{proof}

Strong constraint algebras were defined as internal monoids with respect to
$\strtensor_\field{k}$.
Continuing, we obtain strong constraint modules internal to 
$(\strConMod_{\field{k}}, \strtensor_\field{k})$.
\filbreak
\begin{definition}[Modules over strong constraint algebras]
	\label{def:strConstraintModules}%
	\index{strong constraint!bimodule}
Let $\algebra{A}, \algebra{B} \in \strConAlg$
be strong constraint algebras.
\begin{definitionlist}
	\item A \emph{strong constraint $(\algebra{B},\algebra{A})$-bimodule}
	is a $(\algebra{B},\algebra{A})$-bimodule internal to the 
	monoidal category $(\strConMod_{\field{k}}, 
	\strtensor_\field{k})$.
	\item A \emph{strong constraint left $\algebra{B}$-module}
	is a strong constraint $(\algebra{B},\field{k})$-module.
	\item A \emph{strong constraint right $\algebra{A}$-module}
	is a strong constraint $(\field{k},\algebra{A})$-module.
	\item The categories of	strong constraint right 
	$\algebra{A}$-modules,
	left $\algebra{B}$-modules and
	$(\algebra{B}, \algebra{A})$-bimodules are denoted by
	\glsadd{strConRModA}$\strConRMod{\algebra{A}}$,
	\glsadd{strConLModB}$\strConLMod{\algebra{B}}$ 
	and
	\glsadd{strConBimodBA}$\strConBimod(\algebra{B},\algebra{A})$, respectively.
\end{definitionlist}
\end{definition}

We will denote the set of constraint morphisms between strong constraint 
right $\algebra{A}$-modules $\module{E}$ and $\module{F}$ by
$\Hom_{\algebra{A}}(\module{E},\module{F})$.

Let us take a closer look at strong constraint right 
$\algebra{A}$-modules for a strong constraint algebra $\algebra{A}$.
The structure for left- and bimodules then follows analogously.

\begin{proposition}
	\label{prop:StrConModulesStructure}
Let $\algebra{A} \in \strConAlg$ and $\module{E} \in 
\strConMod_\field{k}$.
Then the structure of a strong constraint right 
$\algebra{A}$-module
on $\module{E}$ is equivalently given by the following data:
\begin{propositionlist}
	\item \label{prop:StrConModulesStructure_1}
	an $\algebra{A}_\Total$-module structure
	$\rho_\Total \colon \module{E}_\Total \tensor[\field{k}] \algebra{A}_ \Total \to \module{E}_\Total$ on $\module{E}_\Total$,
	\item \label{prop:StrConModulesStructure_2}
	an $\algebra{A}_\Wobs$-module structure 
	$\rho_\Wobs^{\Wobs\Wobs} \colon \module{E}_\Wobs \tensor[\field{k}] \algebra{A}_ \Wobs \to \module{E}_\Wobs$ on $\module{E}_\Wobs$,
	\item an $\algebra{A}_\Total$-module structure 
	$\rho_\Wobs^{\Null\Total} \colon \module{E}_\Null \tensor[\field{k}] \algebra{A}_ \Total \to \module{E}_\Null$ on $\module{E}_\Null$,
	\item \label{prop:StrConModulesStructure_3}
	a morphism $\rho_\Wobs^{\Total\Null} \colon 
	\module{E}_\Total \tensor[\field{k}] \algebra{A}_\Null \to 
	\module{E}_\Null$ of right $\algebra{A}_\Total$- and 
	$\algebra{A}_\Wobs$-modules,
\end{propositionlist}
such that
\begin{propositionlist}[resume]
	\item \label{prop:StrConModulesStructure_4}
	$\iota_\module{E} \colon \module{E}_\Wobs \to 
	\module{E}_\Total$ is a morphism of right 
	$\algebra{A}_\Wobs$-modules,
	\item \label{prop:StrConModulesStructure_5}
	$\module{E}_\Null \subseteq \module{E}_\Wobs$ is an
	$\algebra{A}_\Wobs$-submodule,
	\item \label{prop:StrConModulesStructure_6}
	$\iota_\module{E}\at{\module{E}_\Null} \colon 
	\module{E}_\Null \to \module{E}_\Total$
	is a morphism of $\algebra{A}_\Total$-modules.
\end{propositionlist}
\end{proposition}

\begin{proof}
A strong constraint right $\algebra{A}$-module $\module{E}$ is 
given 
by a constraint $\field{k}$-module $\module{E}$ together with a 
constraint map $\rho \colon \module{E} \strtensor[\field{k}] 
\algebra{A} \to \algebra{A}$ fulfilling the usual axioms for a 
right action, see \autoref{def:RightModuleOverMonoid}.
By \autoref{prop:MonoidalStructureCstrModk} and the fact that the 
strong tensor product is given by a colimit, see 
\eqref{diag:strTensorColimit}, the map $\rho_\Wobs$ is 
equivalently 
described by $\field{k}$-module morphisms
\begin{align*}
	\rho_\Wobs^{\Wobs\Wobs} \colon \module{E}_\Wobs 
	\tensor[\field{k}] \algebra{A}_\Wobs \to \module{E}_\Wobs,
	\qquad
	\rho_\Wobs^{\Total\Null} \colon \module{E}_\Total 
	\tensor[\field{k}] \algebra{A}_\Null \to \module{E}_\Null,
	\qquad\text{and}\qquad
	\rho_\Wobs^{\Null\Total} \colon \module{E}_\Null 
	\tensor[\field{k}] \algebra{A}_\Total \to \module{E}_\Null,
\end{align*}
fulfilling
\begin{align*}
	\rho_\Wobs^{\Total\Null}(\iota_\module{E}(x), a)
	&= \rho_\Wobs^{\Wobs\Wobs}(x,a)
	&&\text{for all } x \in \module{E}_\Wobs, a \in 
	\algebra{A}_\Null,\\
	\rho_\Wobs^{\Null\Total}(x, \iota_\algebra{A}(a))
	&= \rho_\Wobs^{\Wobs\Wobs}(x,a)
	&&\text{for all } x \in \module{E}_\Null, a \in 
	\algebra{A}_\Wobs.
\end{align*}
From the fact that $\rho$ defines a module structure on 
$\module{E}$ 
it follows directly that $\rho_\Total$ and 
$\rho_\Wobs^{\Total\Null}$ 
define right $\algebra{A}_\Total$-module structures on 
$\module{E}_\Total$ and $\module{E}_\Null$, respectively.
Moreover, $\module{E}_\Wobs$ becomes a right
$\algebra{A}_\Wobs$-module via $\rho_\Wobs^{\Wobs\Wobs}$.
\end{proof}

The strong tensor product of strong constraint $\field{k}$-modules 
now carries over to strong constraint $\algebra{A}$-modules.
The following is a reformulation of the tensor product of internal modules,
see \autoref{prop: tensor product of bimodules}, internal to $(\ConMod_{\field{k}}, \strtensor[\field{k}])$ and
spelled out in components.
\filbreak
\begin{proposition}[Strong tensor product of strong constraint modules]
	\label{lemma:strTensorProductBimodules}%
	\index{strong tensor product!strong constraint bimodules}
Let $\algebra{A}, \algebra{B}$ and $\algebra{C}$
be strong constraint algebras and let
$\module{F} \in \strConBimod(\algebra{C},\algebra{B})$
and
$\module{E} \in \strConBimod(\algebra{B},\algebra{A})$
be strong constraint bimodules.
Then the strong constraint $(\algebra{C},\algebra{A})$-bimodule
$\module{F} \strtensor[\algebra{B}] \module{E}$
is given by
\begin{equation}
	\begin{split}
		(\module{F} \strtensor[\algebra{B}] \module{E})_\Total
		&\coloneqq \module{F}_\Total \tensor[\algebra{B}_\Total] 
		\module{E}_\Total, \\
		(\module{F} \strtensor[\algebra{B}] \module{E})_\Wobs
		&\coloneqq
		\frac{(\module{F}_\Wobs \tensor[\algebra{B}_\Wobs] \module{E}_\Wobs)
			\oplus (\module{F}_\Null \tensor[\algebra{B}_\Total] \module{E}_\Total )
			\oplus (\module{F}_\Total \tensor[\algebra{B}_\Total] \module{E}_\Null)}
		{\module{I}^\algebra{B}_{\module{F},\module{E}}}, 
		\\
		(\module{F} \strtensor[\algebra{B}] \module{E})_\Null
		&\coloneqq
		\frac{\left( \module{F}_\Null \tensor[\algebra{B}_\Wobs] 
			\module{E}_\Wobs + \module{F}_\Wobs \tensor[\algebra{B}_\Wobs] 
			\module{E}_\Null \right)
			\oplus
			(\module{F}_\Null \tensor[\algebra{B}_\Total] \module{E}_\Total)
			\oplus (\module{F}_\Total \tensor[\algebra{B}_\Total] \module{E}_\Null)}
		{\module{I}^\algebra{B}_{\module{F},\module{E}}},
	\end{split}
\end{equation}
with
\begin{equation}
	\begin{split}
		\module{I}^\algebra{B}_{\module{F},\module{E}}
		\coloneqq &\Span_\field{k}\left\{
		(x_0 \tensor y,0,0) - (0,x_0 \tensor \iota_\module{E}(y),0)
		\mid x_0 \in \module{F}_\Null, y \in 
		\module{E}_\Wobs
		\right\} \\
		&+ \Span_\field{k}\left\{
		(x \tensor y_0,0,0) - (0,0,\iota_\module{F}(x) \tensor y_0)
		\mid x \in \module{F}_\Wobs, y_0 \in \module{E}_\Null
		\right\}.
	\end{split}
\end{equation}
\end{proposition}

\subsubsection{Embedded Strong Constraint Algebras and their Modules}

Let $\algebra{A} \in \strConAlg$ be a strong constraint algebra with multiplication
$\mu \colon \algebra{A} \strtensor[\field{k}] \algebra{A} \to \algebra{A}$.
If $\iota_\algebra{A} \colon \algebra{A}_\Wobs \to \algebra{A}_\Total$
is injective, then $\mu_\Wobs^{\Wobs\Wobs}$, $\mu_\Wobs^{\Null\Total}$
and $\mu_\Wobs^{\Total\Null}$ are completely determined by
$\mu_\Total$.
Hence in this case the notion of strong constraint algebras simplifies 
drastically.

\begin{corollary}
	\index{embedded strong constraint!algebra}
Let $\algebra{A} \in \injConMod_\field{k}$.
Then a strong constraint algebra structure on $\algebra{A}$ is 
equivalently given by an algebra structure on 
$\algebra{A}_\Total$ such that $\algebra{A}_\Wobs \subseteq 
\algebra{A}_\Total$ is a subalgebra and $\algebra{A}_\Null 
\subseteq \algebra{A}_\Total$ is a two-sided ideal with
$\algebra{A}_\Null \subseteq \algebra{A}_\Wobs$.
\end{corollary}

Note that non-embedded strong constraint algebras carry additional structure with respect to constraint algebras,
while embedded strong constraint algebras do not.
They just fulfil the additional property of $\algebra{A}_\Null$ being a two-sided ideal in $\algebra{A}_\Total$.

\begin{example}\
	\label{ex:injstrConAlg}
	\index{constraint!functions on a set}
\begin{examplelist}
	\item \label{ex:injstrConAlg_1}
	Let $M \in \injConSet$ be an embedded constraint set.
	Then $\ConMap(M,\field{k})$ is an embedded strong constraint algebra,
	cf. \autoref{ex:injConAlg} \ref{ex:injConAlg_1}.
	\item Let $M \in \injstrConSet$ be an embedded strong constraint set.
	Then $\strConMap(M,\field{k})$
	is an embedded constraint algebra which is in general \emph{not}
	strong constraint, since $\strConMap(M,\field{k})_\Null$
	is not a two-sided ideal in $\strConMap(M,\field{k})_\Total$ in general,
	cf. \autoref{ex:injConAlg} \ref{ex:injConAlg_3}.
\end{examplelist}
\end{example}

\begin{remark}
Non-commutative examples of constraint algebras will rarely be strong,
see e.g. the coisotropic creed in \cite{lu:1993a}.
We will come back to this in \autoref{chap:DeformationTheory}.
\end{remark}

Let us now turn to modules.
We call a strong constraint bimodule $\module{E}$ \emph{embedded} if $\iota_\module{E}$
is injective,
and we denote the category of embedded strong constraint $(\algebra{B},\algebra{A})$-bimodules
by
\glsadd{injstrConBimodBA}$\injstrConBimod(\algebra{B,\algebra{A}})$, etc.
In this case the various left and right multiplications in \autoref{prop:StrConModulesStructure} are determined by their $\TOTAL$-components.
This gives the following characterization:

\begin{lemma}
	\label{lem:EnmbStrConModules}
	\index{embedded strong constraint!bimodule}
	\glsadd{injstrConAlg}
Let $\algebra{A}, \algebra{B} \in \injstrConAlg$ and
$\module{E} \in \injstrConMod_\field{k}$.
Then the structure of a strong constraint
$(\algebra{B},\algebra{A})$-module
is equivalently given by a 
$(\algebra{B}_\Total,\algebra{A}_\Total)$-bimodule structure
on $\module{E}_\Total$ such that
\begin{lemmalist}
	\item $\module{E}_\Wobs \subseteq \module{E}_\Total$
	is a $(\algebra{B}_\Wobs,\algebra{A}_\Wobs)$-submodule,
	\item $\module{E}_\Null \subseteq \module{E}_\Total$
	is a $(\algebra{B}_\Total,\algebra{A}_\Total)$-submodule,
	\item $\module{E}_\Null \subseteq \module{E}_\Wobs$
	is a $(\algebra{B}_\Wobs,\algebra{A}_\Wobs)$-submodule.
\end{lemmalist}
\end{lemma}

Similarly to embedded strong constraint algebras also embedded strong constraint modules
are just constraint modules with an additional property instead of additional structure as in the non-embedded situation.

\begin{lemma}
	\label{lem:EmbStrConModulesMorphisms}
	Let $\algebra{A},\algebra{B} \in \injstrConAlg$ and let
	\glsadd{injstrConBimodBA}$\module{E}, \module{F} \in 
	\injstrConBimod(\algebra{B},\algebra{A})$ be 
	strong constraint bimodules.
	Then a bimodule morphism
	$\Phi \colon \module{E} \to \module{F}$
	is equivalently given by a
	$(\algebra{B}_\Total, \algebra{A}_\Total)$-bimodule morphism
	$\Phi_\Total \colon \module{E}_\Total \to \module{F}_\Total$
	such that $\Phi_\Total(\module{E}_\Wobs) \subseteq \module{F}_\Wobs$
	and $\Phi_\Total(\module{E}_\Null) \subseteq \module{F}_\Null$.
\end{lemma}

Since constraint morphisms of embedded constraint modules
are determined by their behaviour on the $\TOTAL$-components,
it is clear that also $\ConHom_\algebra{A}(\module{E},\module{F})$
is embedded if $\module{E}$ and $\module{F}$ are embedded.

\begin{proposition}
	\label{prop:strInternalHom}
	\index{internal hom!embedded strong constraint bimodules}
Let $\algebra{A}, \algebra{B}, \algebra{C} \in \injstrConAlg$ be embedded strong 
constraint algebras and let
$\module{E} \in \injstrConBimod(\algebra{B},\algebra{A})$
as well as 
$\module{F} \in \injstrConBimod(\algebra{C},\algebra{A})$
be embedded strong constraint bimodules.
Then the right $\algebra{A}$-module morphisms 
$\ConHom_\algebra{A}(\module{E},\module{F})$
form an embedded strong constraint $(\algebra{C}, \algebra{B})$-bimodule.
\end{proposition}

\begin{proof}
It is clear that $\ConHom_\algebra{A}(\module{E},\module{F})$
is an embedded constraint $(\algebra{C},\algebra{B})$-bimodule.
To see that it is a strong bimodule, consider
$\Phi \in \ConHom_\algebra{A}(\module{E},\module{F})_\Null$
and $c \in \algebra{C}_\Total$.
Then for all $x \in \module{E}_\Wobs$ we have
$(c \cdot \Phi) (x) = c \cdot \Phi(x) \in \algebra{C}_\Total \cdot \module{F}_\Null \subseteq \module{F}_\Null$
and thus $c\cdot \Phi \in \ConHom_\algebra{A}(\module{E},\module{F})_\Null$.
Analogously, we obtain $\Phi \cdot b \in \ConHom_\algebra{A}(\module{E},\module{F})$ for all
$b \in \algebra{B}_\Total$.
\end{proof}

Even though for a strong constraint module
$\module{E} \in \strConBimod(\algebra{B}, \algebra{A})$
the constraint endomorphisms
$\ConEnd_\algebra{A}(\module{E})
= \ConHom_\algebra{A}(\module{E},\module{E})$
form a strong constraint module, they will in general not be a strong
constraint algebra, since the composition of
$\Phi \in \ConEnd_\algebra{A}(\module{E})_\Total$
with $\Psi \in \ConEnd_\algebra{A}(\module{E})_\Null$
might not end up in the $\NULL$-component.
Nevertheless, $\ConEnd_\algebra{A}(\module{E})$
is still a constraint algebra with respect to composition.

By \autoref{prop:strInternalHom} the dual module
$\module{E}^* = \strConHom_\algebra{A}(\module{E},\algebra{A})$
in particular
is an embedded strong constraint module.
It is also easy to see that the direct sum $\module{E} \oplus \module{F}$ of two embedded strong constraint algebras is again embedded strong constraint.

\begin{proposition}
	Let $\algebra{A}, \algebra{B} \in \injstrConAlg$ and $\module{E},\module{F} \in \injstrConBimod(\algebra{B}, \algebra{A})$.
	Then there exists a canonical isomorphism
	\begin{equation}
		(\module{E} \oplus \module{F})^*
		\simeq \module{E}^* \oplus \module{F}^*
	\end{equation}
	of embedded strong constraint $(\algebra{B},\algebra{A})$-bimodules.
\end{proposition}

\begin{proof}
	The bimodule morphism
	$\Phi \colon \module{E}^*_\Total \oplus \module{F}^*_\Total \to (\module{E} \oplus \module{F})^*_\Total$
	given by
	$\Phi(\alpha + \beta)(v + w) \coloneqq \alpha(v) + \beta(w)$
	for all $\alpha \in \module{E}^*$, $\beta \in \module{F}^*$, $v \in \module{E}$ and
	$w \in \module{F}$
	is invertible, with inverse given by
	$\Phi^{-1}(\eta) = (\eta \circ i_\module{E}, \eta \circ i_\module{F})$.
	Here $i_\module{E}$ and $i_\module{F}$ denote the canonical inclusions of $\module{E}$
	and $\module{F}$ in $\module{E} \oplus \module{F}$.
	It is now a straightforward proof to show that both $\Phi$ and $\Phi^{-1}$ are constraint morphisms.
\end{proof}

The strong tensor product of two embedded strong constraint modules will in general not be embedded.
In other words, $\injstrConBimod(\algebra{A})$ is not a monoidal subcategory
of $(\strConBimod(\algebra{A}), \strtensor[\algebra{A}])$,
nevertheless, it carries enough structure to transfer $\strtensor[\algebra{A}]$
to a monoidal structure on $\injstrConBimod(\algebra{A})$,
similar to \autoref{prop:strTensorInjectiveModules}:

\begin{proposition}[The category $\injstrConBimod(\algebra{A})$]\
	\label{prop:CategoryInjstrCBimod}
	Let $\algebra{A} \in \injstrConAlg$ be given.
	\begin{propositionlist}
		\item \label{prop:CategoryInjstrCBimod_1}
		The category $\injstrConBimod(\algebra{A})$ is a reflective subcategory 
		of
		$\strConBimod(\algebra{A})$
		with reflector\linebreak
		$\cdot^\inj \colon \strConBimod(\algebra{A})
		\to \injstrConBimod(\algebra{A})$ given by
		\begin{equation}
			\module{E}^\inj \coloneqq (\module{E}_\Total, 
			\iota_\module{E}(\module{E}_\Wobs),
			\iota_\module{E}(\module{E}_\Null) ).
		\end{equation}
		\item \label{prop:CategoryInjstrCBimod_2}
		The subcategory $\injstrConBimod(\algebra{A})$
		is closed under finite limits.
		\item \label{prop:CategoryInjstrCBimod_3}
		The category $\injstrConBimod(\algebra{A})$ is monoidal with 
		respect to $\injstrtensor[\algebra{A}]$ defined by
		\index{strong tensor product!embedded}
		\glsadd{injStrTensorA}
		\begin{equation}
			\module{E} \injstrtensor[\algebra{A}] \module{F} \coloneqq 
			(\module{E} 
			\strtensor[\algebra{A}] \module{F})^\inj.
		\end{equation}
		\item \label{prop:CategoryInjstrCBimod_4}
		The functor $\argument^\inj \colon 
		(\strConBimod(\algebra{A}),\strtensor[\algebra{A}])
		\to (\injstrConBimod(\algebra{A}), \injstrtensor[\algebra{A}])$
		is monoidal, and the functor
		$\functor{U} \colon 
		(\injstrConBimod(\algebra{A}),\injstrtensor[\algebra{A}]) 
		\to (\strConBimod(\algebra{A}), \strtensor[\algebra{A}])$
		is lax monoidal.
	\end{propositionlist}
\end{proposition}

\begin{proof}
	The properties \ref{prop:StrConModulesStructure_4} to \ref{prop:StrConModulesStructure_6}
	from \autoref{prop:StrConModulesStructure}
	ensure that $\module{E}^\inj$ is an embedded strong constraint
	module as in \autoref{lem:EnmbStrConModules}.
	With this \ref{prop:CategoryInjstrCBimod_1}
	is clear, and \ref{prop:CategoryInjstrCBimod_2}
	follows directly.
	Part \ref{prop:CategoryInjstrCBimod_3} and \ref{prop:CategoryInjstrCBimod_4}
	follow again from Day's reflection theorem, see \autoref{thm:ReflectionTheorem}.
	For this we only need to see that
	$(\eta_\module{E} \tensor \eta_\module{F})^\inj 
	\colon (\module{E} \strtensor[\algebra{A}] \module{F})^\inj
	\to (\module{E}^\inj \strtensor[\algebra{A}] \module{F}^\inj)^\inj$,
	with $\eta_\module{E} \colon \module{E} \to \module{E}^\inj$
	given by $(\eta_\module{E})_\Total = \id_{\module{E}_\Total}$
	and $(\eta_\module{F})_\Wobs(x) = \iota_\module{E}(x)$,
	is an isomorphism for all $\module{E}, \module{F} \in \strConBimod(\algebra{A})$.
	This is clear since
	\begin{equation*}
		(\eta_\module{E} \tensor \eta_\module{F})^\inj\big(\iota_\module{E}(x) \tensor \iota_\module{F}(y)\big)
		= \iota_\module{E}(x) \tensor \iota_\module{F}(y).
	\end{equation*}
\end{proof}

This embedded strong constraint tensor product resembles the motivating formulas for
the strong tensor product, see \eqref{eq:MotivationStrTensor_1} and \eqref{eq:MotivationStrTensor_2}:

\begin{corollary}
	\label{cor:injStrTensorExplicit}
Let $\algebra{A} \in \injstrConAlg$,
and let $\module{E}, \module{F} \in \injstrConBimod(\algebra{A})$
be $\algebra{A}$-bimodules.
Then we have
\begin{equation}
\begin{split}
	(\module{E} \injstrtensor[\algebra{A}] \module{F})_\Total
	&= \module{E}_\Total \tensor[\algebra{A}_\Total] 
	\module{F}_\Total,\\
	(\module{E} \injstrtensor[\algebra{A}] \module{F})_\Wobs
	&= \module{E}_\Wobs \tensor[\algebra{A}_\Total] \module{F}_\Wobs 
	+ \module{E}_\Null \tensor[\algebra{A}_\Total] \module{F}_\Total 
	+ \module{E}_\Total \tensor[\algebra{A}_\Total] \module{F}_\Null,
	\\
	(\module{E} \injstrtensor[\algebra{A}] \module{F})_\Null
	&= \module{E}_\Null \tensor[\algebra{A}_\Total] \module{F}_\Total 
	+ \module{E}_\Total \tensor[\algebra{A}_\Total] \module{F}_\Null.
\end{split}
\end{equation}
\end{corollary}

Having the strong tensor product and duals at hand, we obtain a canonical morphism
resembling \autoref{prop:CompatibilityDualTensorConVecSpaces} \ref{prop:CompatibilityDualTensorConVecSpaces_4}
from constraint vector spaces.

\begin{proposition}
	\label{prop:ConHomAndTensorDual}
	Let $\algebra{A} \in \injstrConAlg$
	and let $\module{E}, \module{F} \in \injstrConBimod(\algebra{A})$
	be $\algebra{A}$-bimodules.
	Then 
	\begin{equation}
		\label{eq:ConHomAndTensorDual}
		\module{F}_\Total \tensor[\algebra{A}_\Total] 
		\module{E}^*_\Total
		\ni y \tensor \alpha 
		\mapsto (x \mapsto y \cdot \alpha(x)) 
		\in 
		\Hom_{\algebra{A}_\Total}(\module{E}_\Total,\module{F}_\Total)
	\end{equation}
	defines a constraint morphism
	$\module{F} \injstrtensor[\algebra{A}] \module{E}^*
	\to \ConHom_\algebra{A}(\module{E},\module{F})_\Total$.
\end{proposition}

\begin{proof}
	The map \eqref{eq:ConHomAndTensorDual} is the canonical 
	$\algebra{A}_\Total$-module morphism from classical algebra.
	To show that it is a constraint $\algebra{A}$-module morphism
	consider first
	$y \tensor \alpha \in (\module{F} \injstrtensor 
	\module{E}^*)_\Null
	= \module{F} \ConGrid[2][2][2][2][0][0][2] \module{E}^*$.
	Here we use the notation introduced in 
	\autoref{not:ConIndSetProducts}.
	If $x \in \module{E}_\Wobs$, then
	$y \cdot \alpha(x) \in \module{F}_\Null \cdot 
	\algebra{A}_\Total + \module{F}_\Total \cdot \algebra{A}_\Null
	\subseteq \module{F}_\Null$.
	Hence \eqref{eq:ConHomAndTensorDual} maps $\NULL$-component to
	$\NULL$-component.
	Now let $y \tensor \alpha \in 
	\module{F} \ConGrid[0][0][0][0][2] \module{E}^*
	\subseteq (\module{F} \injstrtensor[\algebra{A}] 
	\module{E}^*)_\Wobs$.
	If $x \in \module{E}_\Null$, then
	$y \cdot \alpha(x) \in \module{F}_\Wobs \cdot \algebra{A}_\Null 
	\subseteq \module{F}_\Null$.
	If $x \in \module{E}_\Wobs$, then
	$y \cdot \alpha(x) \in \module{F}_\Wobs \cdot \algebra{A}_\Wobs
	\subseteq \module{F}_\Wobs$.
	Thus \eqref{eq:ConHomAndTensorDual}
	is a constraint morphism.
\end{proof}


\subsubsection{Strong Hull}

We obtain a forgetful functor
$\functor{U} \colon \strConAlg \to \ConAlg$
by mapping a strong constraint algebra $(\algebra{A},\mu)$
to the constraint algebra $\algebra{A}$ obtained by
dismissing the $\algebra{A}_\Total$-bimodule structure on 
$\algebra{A}_\Null$.
This functor obviously restricts to
$\functor{U} \colon \injstrConAlg \to \injConAlg$.
In this case we can easily describe its corresponding free 
construction.
\filbreak
\begin{proposition}\
	\label{prop:StrongificationAlgebras}
	\index{strong hull!constraint algebra}
\begin{propositionlist}
	\item Let $\algebra{A} \in \injConAlg$, then
	\glsadd{strHull}
	\begin{equation}
		\begin{split}
			\algebra{A}^\str_\Total
			&= \algebra{A}_\Total,\\
			\algebra{A}^\str_\Wobs
			&= \algebra{A}_\Wobs
			+ \algebra{A}_\Total \cdot \algebra{A}_\Null \cdot 
			\algebra{A}_\Total,\\
			\algebra{A}^\str_\Null
			&= \algebra{A}_\Total \cdot \algebra{A}_\Null \cdot 
			\algebra{A}_\Total
		\end{split}
	\end{equation}
	is a strong constraint algebra $\algebra{A}^\str \in 
	\injstrConAlg$.
	\item Mapping $\algebra{A} \in \injConAlg$ to
	$\algebra{A}^\str \in \injstrConAlg$ and morphisms
	$\phi \colon \algebra{A} \to \algebra{B}$ to
	$\phi^\str \colon \algebra{A}^\str \to \algebra{B}^\str$ given 
	by $\phi^\str = \phi$
	defines a functor
	$\argument^\str \colon \injConAlg \to \injstrConAlg$.
	\item The functor $\argument^\str$ is left adjoint to the 
	forgetful functor
	$\functor{U} \colon \injstrConAlg \to \injConAlg$.
	\item $\injstrConAlg$ is a reflective subcategory of
	$\injConAlg$.
\end{propositionlist}
\end{proposition}

\begin{proof}
	The first and second part are straightforward checks.
	For the third part let
	$\eta_\algebra{A} \colon \algebra{A} \to \functor{U}(\algebra{A}^\str)$
	be the obvious inclusion for every $\algebra{A} \in \ConAlg$,
	and let 
	$\varepsilon_\algebra{B} \colon (\functor{U}\algebra{B})^\str \to 
	\algebra{B}$
	be the identity for every $\algebra{B} \in \injstrConAlg$.
	It is now easy to check that the maps $\eta$ and $\varepsilon$
	are the unit and counit of the required adjunction.
	The last part is clear, since the counit $\varepsilon$
	is just the identity.
\end{proof}

\begin{figure}[t]
	\centering
	\includegraphics{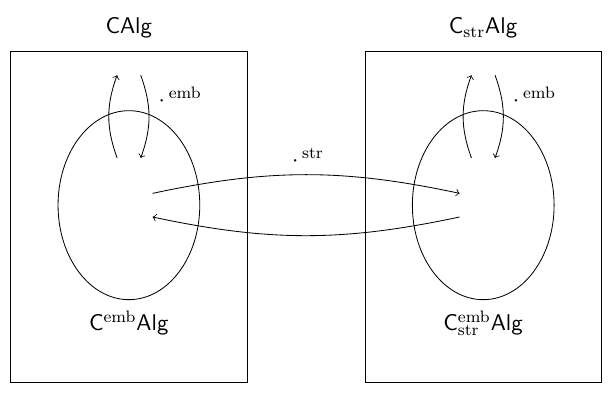}
	\caption{Overview of the different categories of constraint algebras. Unnamed arrows denote forgetful functors.}
	\label{fig:ConAlgebras}
\end{figure}

We will call $\algebra{A}^\str$ the \emph{strong hull of}
$\algebra{A}$.
See \autoref{fig:ConAlgebras} for an overview of the various categories of constraint algebras and their relationship.
For functions on embedded strong constraint sets the construction of 
the strong hull can be viewed as the algebraic analogue of forgetting 
the equivalence relation outside of the subset.

\begin{proposition}
	\label{prop:ForgettingVSstrongHull}
Let $\mathbb{K}$ be a field.
The diagram
\begin{equation}
	\begin{tikzcd}[column sep = huge]
		\injstrConSet
		\arrow[r,"{\strConMap(\argument,\field{K})}"]
		\arrow[d,"\functor{U}"{swap}]
		&\injConAlg
		\arrow[d,"\argument^\str"]\\
		\injConSet
		\arrow[r,"{\ConMap(\argument,\field{K})}"]
		&\injstrConAlg
	\end{tikzcd}
\end{equation}
commutes up to a natural isomorphism.
Here $\functor{U} \colon \injstrConSet \to \injConSet$ denotes
the functor forgetting the equivalence relation outside of the 
$\WOBS$-component, see \autoref{prop:CategoryStrCSet}.	
\end{proposition}

\begin{proof}
Since $\field{K}$ is a field, \autoref{prop:StrongificationAlgebras}
applies.
Now on the $\TOTAL$-component the diagram commutes strictly.
On the one hand we know by \autoref{ex:ConAlgebras} \ref{ex:FunctionsOnCSetI} that for every
embedded strong constraint set $M$ the ideal
$\ConMap(\functor{U}(M),\field{k})_\Null$
is just the vanishing ideal $\vanishing_{M_\Wobs}$ of $M_\Wobs$.
On the other hand \autoref{prop:CategoryStrCSet} \ref{prop:CategoryStrCSet_3} characterizes 
$\strConMap(M,\field{k})_\Null$ as those functions vanishing on 
$M_\Wobs$ which are constant along the equivalence classes on 
$M_\Total$.
Using the characteristic function $\chi_{M_\Wobs} \colon M_\Total 
\to \{0,1\}$ we can write every $f \in \vanishing_{M_\Wobs}$
as $f = \chi_{M_\Wobs} \cdot f$, with $\chi \in 
\strConMap(M,\field{k})_\Null$ and $f \in \algebra{A}_\Total$.
Hence the $\NULL$-components agree.
Similarly, every $g \in \strConMap(M,\field{k})_\Wobs$
is constant along equivalence classes on $M_\Wobs$ and can be 
written as
$g = (1 - \chi_{M_\Wobs}) \cdot g + \chi_{M_\Wobs} \cdot g$,
with $(1 - \chi_{M_\Wobs}) \cdot g \in \ConMap(M,\field{k})_\Wobs$
and $\chi_{M_\Wobs} \cdot g \in \ConMap(M,\field{k})^\str_\Null$.
\end{proof}

There is again the obvious forgetful functor
$\functor{U} \colon \strConRMod{\algebra{A}} \to 
\ConRMod{\algebra{A}}$
by forgetting the module structure $\rho = (\rho_\Total,\rho_\Wobs)$ 
to $(\rho_\Total,\rho_\Wobs^{\Wobs\Wobs})$.
Analogously to the algebra case
there is also a way to construct strong constraint modules out of 
non-strong ones if we assume the algebra and the module to be 
embedded.

\begin{proposition}[Strong hull]\
	\label{prop:StrongHullModules}
	\index{strong hull!constraint bimodule}
Let $\algebra{A},\algebra{B} \in \injstrConAlg$.
\begin{propositionlist}
\item \label{prop:StrongHullModules_1}
 Let $\module{E} \in 
\injConBimod(\algebra{B},\algebra{A})$.
Then
\begin{equation}
\begin{split}
	\module{E}^\str_\Total
	&\coloneqq \module{E}_\Total, \\
	\module{E}^\str_\Wobs
	&\coloneqq \module{E}_\Wobs + \algebra{B}_\Total \cdot 
	\module{E}_\Null \cdot \algebra{A}_\Total + \algebra{B}_\Null \cdot \module{E}_\Total + \module{E}_\Total \cdot \algebra{A}_\Null, \\
	\module{E}^\str_\Null
	&\coloneqq  \algebra{B}_\Total \cdot \module{E}_\Null \cdot 
	\algebra{A}_\Total + \algebra{B}_\Null \cdot \module{E}_\Total + \module{E}_\Total \cdot \algebra{A}_\Null 
\end{split}
\end{equation}
is a strong constraint $(\algebra{B},\algebra{A})$-bimodule.
\item \label{prop:StrongHullModules_2}
Mapping $\module{E} \in \injConBimod(\algebra{B},\algebra{A})$ 
to 
$\module{E}^\str \in \injstrConBimod(\algebra{B},\algebra{A})$ and 
morphisms
$\Phi \colon \module{E} \to \module{F}$ to 
$\Phi^\str \colon \module{E}^\str \to \module{F}^\str$ given by
$\Phi^\str \coloneqq \Phi$
defines a functor
\begin{equation}
	\argument^\str \colon \injConBimod(\algebra{B},\algebra{A}) \to 
	\injstrConBimod(\algebra{B},\algebra{A}).
\end{equation}
\item \label{prop:StrongHullModules_3}
$\injstrConBimod(\algebra{B},\algebra{A})$ is a reflective 
subcategory
of $\injConBimod(\algebra{B},\algebra{A})$ with reflector 
$\argument^\str$.
\end{propositionlist}
\end{proposition}

\begin{proof}
	The first and second part are clear.
	For the third part, note that by 
	\autoref{lem:EnmbStrConModules} and
	\autoref{lem:EmbStrConModulesMorphisms}
	the category $\injstrConBimod(\algebra{B},\algebra{A})$
	is a full subcategory of $\injConBimod(\algebra{B},\algebra{A})$.
	To show that $\argument^\str$ is left adjoint to the embedding
	$\functor{U} \colon \injstrConBimod(\algebra{B},\algebra{A})
	\to \injConBimod(\algebra{B},\algebra{A})$
	consider the counit $\varepsilon$, defined for every 
	$\module{E} \in \injstrConBimod(\algebra{B},\algebra{A})$
	by
	$\varepsilon_\module{E} \coloneqq \id_{\module{E}}$,
	and the unit $\eta$,
	defined for every
	$\module{F} \in \injConBimod(\algebra{B},\algebra{A})$
	by the obvious inclusion
	$\eta_\module{F} \colon \module{F} \to \module{F}^\str$.
	The triangle identities are then easily verified, and since $\varepsilon$
	is an isomorphism, we see that
	$\injstrConBimod(\algebra{B},\algebra{A})$ is
	a reflective subcategory.
\end{proof}

For any $\algebra{A} \in \injstrConAlg$ we know that
$\injConBimod(\algebra{A})$ is a monoidal category
with tensor product $\injtensor[\algebra{A}]$ as
given in \autoref{prop:CategoryInjCBimod}.
We would now like to transport this monoidal structure to the
reflective subcategory $\injstrConBimod(\algebra{A})$,
but in this generality Day's reflection theorem does not apply.
Nevertheless, when we restrict ourselves to symmetric bimodules
over commutative strong constraint algebras
this can be achieved.
We denote the category of symmetric embedded strong constraint bimodules by
\glsadd{syminjstrConBimodA}$\injstrConBimod(\algebra{A})_\sym$.

\begin{proposition}
	Let $\algebra{A} \in \injstrConAlg$ be a
	commutative embedded strong constraint algebra.
	\begin{propositionlist}
		\item \label{prop:StrongHullModules_4}
		$\injstrConBimod(\algebra{A})_\sym$ is a monoidal category with 
		respect to $\strcirctensor[\algebra{A}]$ defined by
		\begin{equation}
			\module{E} \strcirctensor[\algebra{A}] \module{F} \coloneqq 
			(\module{E} 
			\injtensor[\algebra{A}] \module{F})^\str.
		\end{equation}
		\item \label{prop:StrongHullModules_5}
		The functor $\argument^\str \colon 
		(\injConBimod(\algebra{A})_\sym,\injtensor[\algebra{A}])
		\to (\injstrConBimod(\algebra{A})_\sym, \strcirctensor[\algebra{A}])$
		is monoidal.
	\end{propositionlist}
\end{proposition}

\begin{proof}
	We use again Day's reflection theorem,
	see \autoref{thm:ReflectionTheorem}.
	For this note that \autoref{prop:StrongHullModules}
	restricts to the subcategories of symmetric bimodules.
	Hence $\injConBimod(\algebra{A})_\sym$
	is a reflective subcategory of $\injstrConBimod(\algebra{A})$.
	Furthermore, we have canonically
	\begin{align*}
		(\module{E}^\str \injtensor[\algebra{A}] \module{F}^\str)^\str_\Wobs
		&= (\module{E}^\str \injtensor[\algebra{A}] \module{F}^\str)_\Wobs
		+\algebra{A}_\Total \cdot (\module{E}^\str \injtensor[\algebra{A}] \module{F}^\str)_\Null \cdot \algebra{A}_\Total \\
		&\quad + \algebra{A}_\Null \cdot (\module{E}^\str \injtensor[\algebra{A}] \module{E}^\str)_\Total
		+ (\module{E}^\str \injtensor[\algebra{A}] \module{E}^\str)_\Total \cdot \algebra{A}_\Null\\
		&= \module{E}^\str_\Wobs \tensor[\algebra{A}_\Total] \module{F}^\str_\Wobs
		+ \algebra{A}_\Total \cdot (\module{E}^\str_\Null \tensor[\algebra{A}_\Total] \module{F}^\str_\Wobs) \cdot \algebra{A}_\Total
		+ \algebra{A}_\Total \cdot (\module{E}^\str_\Wobs \tensor[\algebra{A}_\Total] \module{F}^\str_\Null) \cdot \algebra{A}_\Total \\
		&\quad + \algebra{A}_\Null \cdot (\module{E}_\Total \tensor[\algebra{A}_\Total] \module{E}_\Total)
		+ (\module{E}_\Total \tensor[\algebra{A}_\Total] \module{E}_\Total) \cdot \algebra{A}_\Null\\
		&= \module{E}_\Wobs \tensor[\algebra{A}_\Total] \module{F}_\Wobs
		+ \algebra{A}_\Total \cdot (\module{E}_\Null \tensor[\algebra{A}_\Total] \module{F}_\Wobs ) \cdot \algebra{A}_\Total
		+ \algebra{A}_\Total \cdot (\module{E}_\Wobs \tensor[\algebra{A}_\Total] \module{F}_\Null ) \cdot \algebra{A}_\Total \\
		&\quad + \algebra{A}_\Null \cdot (\module{E}_\Total \tensor[\algebra{A}_\Total] \module{E}_\Total)
		+ (\module{E}_\Total \tensor[\algebra{A}_\Total] \module{E}_\Total) \cdot \algebra{A}_\Null\\
		&= (\module{E} \injtensor[\algebra{A}] \module{F})^\str 
	\end{align*}
	for all symmetric bimodules $\module{E}$ and $\module{F}$.
	Thus by \autoref{thm:ReflectionTheorem}
	we see that $\strcirctensor[\algebra{A}]$ defines a monoidal structure
	on $\injstrConBimod(\algebra{A})_\sym$ such that
	$\argument^\str$ becomes a monoidal functor.
\end{proof}

Using the definition of the strong hull and $\argument^\inj$ as defined in
\autoref{prop:CategoryInjCMod}
directly yields the following explicit description of $\strcirctensor$:

\begin{corollary}
	\label{cor:strCircTensorExplicit}
Let $\algebra{A} \in \injstrConAlg$ be commutative,
and let $\module{E}, \module{F} \in \injstrConBimod(\algebra{A})_\sym$
be symmetric $\algebra{A}$-bimodules.
Then we have
\begin{equation}
\begin{split}
	(\module{E} \strcirctensor[\algebra{A}] \module{F})_\Total
	&= \module{E}_\Total \tensor[\algebra{A}_\Total] 
	\module{F}_\Total, \\
	(\module{E} \strcirctensor[\algebra{A}] \module{F})_\Wobs
	&= \module{E}_\Wobs \tensor[\algebra{A}_\Total] \module{F}_\Wobs,\\
	(\module{E} \strcirctensor[\algebra{A}] \module{F})_\Null
	&= \module{E}_\Null \tensor[\algebra{A}_\Total] \module{F}_\Wobs 
	+ \module{E}_\Wobs \tensor[\algebra{A}_\Total] \module{F}_\Null.
\end{split}
\end{equation}
\end{corollary}

\begin{proposition}
	\label{prop:injstrInternalHom}
	Let $\algebra{A} \in \injstrConAlg$ be commutative.
	The category $\injstrConBimod(\algebra{A})_\sym$ is closed monoidal with respect to
	$\strcirctensor[\algebra{A}]$.
	The internal hom is given by $\ConHom_{\algebra{A}}(\module{E},\module{F})$.
\end{proposition}

\begin{proof}
In \autoref{prop:strInternalHom} we showed that
$\ConHom_{\algebra{A}}(\module{E},\module{F})$
is again an embedded strong constraint bimodule.
The symmetry is clear.
For the $\TOTAL$-component we have the classical
evaluation $\ev_{\module{F}_\Total} \colon \Hom_{\algebra{A}_\Total}(\module{E}_\Total,\module{F}_\Total) \tensor[\algebra{A}_\Total] \module{E}_\Total \to \module{E}_\Total$
and coevaluation 
$\coev_{\module{F}_\Total} \colon \module{F}_\Total \to \Hom_{\algebra{A}_\Total}(\module{E}_\Total, \module{F}_\Total \tensor[\algebra{A}_\Total] \module{E}_\Total)$.
These are easily seen to be constraint morphisms.
\end{proof}

Clearly, the two monoidal structures on 
$\injstrConBimod(\algebra{A})_\sym$ are not unrelated.
Looking at \autoref{cor:injStrTensorExplicit} and
\autoref{cor:strCircTensorExplicit} we see that for
$\module{E}, \module{F} \in \injstrConBimod(\algebra{A})_\sym$
we have a canonical inclusion
\begin{equation*}
	\module{E} \strcirctensor[\algebra{A}] \module{F}
	\hookrightarrow \module{E} \injstrtensor[\algebra{A}] \module{F}.
\end{equation*}
Moreover, dualizing turns one tensor product into another.

\begin{proposition}
	\label{prop:DualOfTensor}
	Let $\algebra{A} \in \injstrConAlg$ be commutative, and let
	$\module{E}, \module{F},\module{G} \in \injstrConBimod(\algebra{A})_\sym$
	be symmetric $\algebra{A}$-bimodules.
	\begin{propositionlist}
		\item \label{prop:DualOfTensor_1}There is a canonical morphism
		$\strConHom_\algebra{A}(\module{E} \injstrtensor[\algebra{A}] \module{F},\module{G})
		\to \strConHom_\algebra{A}(\module{E} \strcirctensor[\algebra{A}] \module{F},\module{G})$
		of constraint $\algebra{A}$-bimodules
		induced by $\module{E} \strcirctensor[\algebra{A}] \module{F} \to \module{E} \injstrtensor[\algebra{A}] \module{F}$.
		\item There is a canonical morphism 
		$\module{E}^* \strcirctensor[\algebra{A}] \module{F}^* 
		\to (\module{E} \injstrtensor[\algebra{A}] \module{F})^*$
		of constraint $\algebra{A}$-bimodules given by 
		\begin{equation}
			\label{eq:DualOfStrTensor}
			\module{E}^*_\Total \tensor[\algebra{A}_\Total] \module{F}^*_\Total
			\ni \alpha \tensor \beta \mapsto
			(x \tensor y \mapsto \alpha(x)\cdot \beta(y))
			\in (\module{E}_\Total \tensor[\algebra{A}_\Total] \module{F}_\Total)^*.
		\end{equation}		
		\item There is a canonical morphism 
		$\module{E}^* \injstrtensor[\algebra{A}] \module{F}^* 
		\to (\module{E} \strcirctensor[\algebra{A}] \module{F})^*$
		of constraint $\algebra{A}$-bimodules given by 
		\begin{equation}
			\label{eq:DualOfTensor}
			\module{E}^*_\Total \tensor[\algebra{A}_\Total] \module{F}^*_\Total
			\ni \alpha \tensor \beta \mapsto
			(x \tensor y \mapsto \alpha(x)\cdot \beta(y))
			\in (\module{E}_\Total \tensor[\algebra{A}_\Total] \module{F}_\Total)^*.
		\end{equation}
	\end{propositionlist}
\end{proposition}

\begin{proof}
	On the $\TOTAL$-components both maps are defined by the canonical map from classical algebra, which is clearly an $\algebra{A}_\Total$-bimodule morphism.
	It remains to show that the maps preserve the substructures.
	For the first part let
	$\alpha \tensor \beta \in (\module{E}^* \strcirctensor[\algebra{A}] \module{F}^*)_\Null
	= \module{E}^* \ConGrid[2][2][0][2] \module{F}^*$.
	\begin{cptitem}
		\item For $x \tensor y \in \module{E} \ConGrid[2][2][2][2][0][0][2] \module{F}
		= (\module{E} \injstrtensor[\algebra{A}] \module{F})_\Null$
		we have
		$\alpha(x)\cdot\beta(y) 
		\in \algebra{A}_\Null \cdot \algebra{A}_\Total + \algebra{A}_\Total \cdot \algebra{A}_\Null
		= \algebra{A}_\Null$.
		\item For $x \tensor y \in \module{E} \ConGrid[0][0][0][0][2] \module{F}
		\subseteq (\module{E} \injstrtensor[\algebra{A}] \module{F})_\Wobs$
		we have
		$\alpha(x)\cdot\beta(y) 
		\in \algebra{A}_\Null \cdot \algebra{A}_\Wobs + \algebra{A}_\Wobs \cdot \algebra{A}_\Null
		= \algebra{A}_\Null$.
	\end{cptitem}
	Thus \eqref{eq:DualOfStrTensor} maps $\NULL$-component to $\NULL$-component.
	Next suppose
	$\alpha \tensor \beta \in 
	\module{E}^* \ConGrid[0][0][0][0][2] \module{F}^*
	\subseteq  (\module{E}^* \strcirctensor[\algebra{A}] \module{F}^*)_\Wobs$.
	\begin{cptitem}
		\item For $x \tensor y \in \module{E} \ConGrid[2][2][2][2][0][0][2] \module{F}
		= (\module{E} \injstrtensor[\algebra{A}] \module{F})_\Null$
		we have
		$\alpha(x)\cdot\beta(y) 
		\in \algebra{A}_\Null \cdot \algebra{A}_\Total + \algebra{A}_\Total \cdot \algebra{A}_\Null
		= \algebra{A}_\Null$.
		\item For $x \tensor y \in \module{E} \ConGrid[0][0][0][0][2] \module{F}
		\subseteq (\module{E} \injstrtensor[\algebra{A}] \module{F})_\Wobs$
		we have
		$\alpha(x)\cdot\beta(y) 
		\in \algebra{A}_\Wobs \cdot \algebra{A}_\Wobs = \algebra{A}_\Wobs$.
	\end{cptitem}
	This shows that \eqref{eq:DualOfStrTensor} also maps $\WOBS$-component to
	$\WOBS$-component and therefore is a constraint morphism.
	For the second part consider at first
	$\alpha \tensor \beta \in (\module{E}^* \injstrtensor[\algebra{A}] \module{F}^*)_\Null
	= \module{E}^* \ConGrid[2][2][2][2][0][0][2] \module{F}^*$.
	\begin{cptitem}
		\item For $x \tensor y \in \module{E} \ConGrid[2][2][0][2][2] \module{F}
		= (\module{E} \strcirctensor[\algebra{A}] \module{F})_\Wobs$
		we have
		$\alpha(x)\cdot\beta(y) 
		\in \algebra{A}_\Null \cdot \algebra{A}_\Total + \algebra{A}_\Total \cdot \algebra{A}_\Null
		= \algebra{A}_\Null$,
	\end{cptitem}
	showing that \eqref{eq:DualOfTensor} maps $\NULL$-component to
	$\NULL$-component.
	Next choose $\alpha \tensor \beta \in \module{E}^* \ConGrid[0][0][0][0][2] \module{F}^*
	\subseteq (\module{E} \strcirctensor[\algebra{A}] \module{F})_\Wobs$.
	\begin{cptitem}
		\item For $x \tensor y \in \module{E} \ConGrid[2][2][0][2] \module{F}
		= (\module{E} \strcirctensor[\algebra{A}] \module{F})_\Null$
		we have
		$\alpha(x)\cdot\beta(y) 
		\in \algebra{A}_\Null \cdot \algebra{A}_\Wobs + \algebra{A}_\Wobs \cdot \algebra{A}_\Null
		= \algebra{A}_\Null$.
		\item For $x \tensor y \in \module{E} \ConGrid[0][0][0][0][2] \module{F}
		\subseteq (\module{E} \strcirctensor[\algebra{A}] \module{F})_\Wobs$
		we have
		$\alpha(x)\cdot\beta(y) 
		\in \algebra{A}_\Wobs \cdot \algebra{A}_\Wobs = \algebra{A}_\Wobs$.
	\end{cptitem}
	This shows that \eqref{eq:DualOfTensor} also preserves the $\WOBS$-components and
	hence is a constraint morphism.
\end{proof}

We cannot expect \eqref{eq:DualOfStrTensor} and \eqref{eq:DualOfTensor}
to be isomorphisms, since they do not even need to be isomorphisms on the $\TOTAL$-component.
However, in classical algebra we know that for finitely generated projective modules these indeed become isomorphisms,
see also \autoref{prop:CompatibilityDualTensorConVecSpaces} for the case of finite dimensional constraint vector spaces.
We will see in \autoref{sec:RegularProjectiveModules}
that they also become constraint isomorphisms if the involved modules
are finitely generated projective as constraint modules.

Since $\injstrConBimod(\algebra{A})_\sym$ is closed monoidal,
we can define an insertion morphism as the composition
\index{insertion}
\begin{equation} \label{eq:Insertion}
	\begin{split}
		\ins \colon \strConHom_\algebra{A}(\module{E} \tensor[\algebra{A}] \module{F}, \module{G}) \tensor[\algebra{A}] \module{E}
		&\longrightarrow \strConHom(\module{E}, \strConHom_\algebra{A}(\module{F},\module{G})) \tensor[\algebra{A}] \module{E}\\
		&\overset{\ev}{\longrightarrow} \strConHom(\module{F, \module{G}}),
	\end{split}
\end{equation}
and we will write $\ins_X(\Phi) \colon \module{F} \to \module{G}$
for $X \in \module{E}$ and $\Phi \colon \module{E} \tensor[\algebra{A}] \module{F} \to \module{G}$.
Using \autoref{prop:DualOfTensor} \ref{prop:DualOfTensor_1} we can define a constraint insertion morphism
as
\begin{equation}
\begin{split}
	\ins \colon \strConHom_\algebra{A}(\module{E} \strtensor[\algebra{A}] \module{F}, \module{G}) \tensor[\algebra{A}] \module{E}	
	&\longrightarrow \strConHom_\algebra{A}(\module{E} \tensor[\algebra{A}] \module{F}, \module{G}) \tensor[\algebra{A}] \module{E} \\
	&\longrightarrow \strConHom(\module{E}, \strConHom_\algebra{A}(\module{F},\module{G})) \tensor[\algebra{A}] \module{E}\\
	&\overset{\ev}{\longrightarrow} \strConHom(\module{F, \module{G}}),
\end{split}
\end{equation}
similar to \eqref{eq:Insertion}.

%

\subsubsection{Reduction}

The reduction of strong constraint modules is again given by first applying the forgetful functor
$\functor{U} \colon \strConBimod \to \ConBimod$
and then using the reduction functor on $\ConBimod$, see \autoref{prop:ReductionConBimod}.
Similar to \autoref{prop:ReductionOnStrCModk} we can show that the tensor product and strong tensor product do not differ after reduction:

\begin{proposition}[Reduction on $\strConBimod$]\
	\label{prop:ReductionOnStrConBimod}
	\index{reduction!strong constraint bimodules}
Let $\algebra{A}, \algebra{B}, \algebra{C} \in \strConAlg$ be given,
and let $\module{E} \in \strConBimod(\algebra{C},\algebra{B})$
and $\module{F} \in \strConBimod(\algebra{B}, \algebra{A})$.
Then there is a canonical isomorphism
\begin{equation}
	(\module{E} \strtensor[\algebra{B}] \module{F})_\red 
	\simeq \module{E}_\red \tensor[\algebra{B}_\red] \module{F}_\red.
\end{equation}
\end{proposition}

\begin{proof}
Recall the definition of $\strtensor[\algebra{B}]$ from \autoref{lemma:strTensorProductBimodules}:
\begin{equation*}
	\begin{split}
		(\module{E} \strtensor[\algebra{B}] \module{F})_\Wobs
		&=
		\frac{(\module{E}_\Wobs \tensor[\algebra{B}_\Wobs] \module{F}_\Wobs)
			\oplus (\module{E}_\Null \tensor[\algebra{B}_\Total] \module{F}_\Total )
			\oplus (\module{E}_\Total \tensor[\algebra{B}_\Total] \module{F}_\Null)}
		{\module{I}^\algebra{B}_{\module{E},\module{F}}}, 
		\\
		(\module{E} \strtensor[\algebra{B}] \module{F})_\Null
		&=
		\frac{\left( \module{E}_\Null \tensor[\algebra{B}_\Wobs] 
			\module{F}_\Wobs + \module{E}_\Wobs \tensor[\algebra{B}_\Wobs] 
			\module{F}_\Null \right)
			\oplus
			(\module{E}_\Null \tensor[\algebra{B}_\Total] \module{F}_\Total)
			\oplus (\module{E}_\Total \tensor[\algebra{B}_\Total] \module{F}_\Null)}
		{\module{I}^\algebra{B}_{\module{E},\module{F}}},
	\end{split}
\end{equation*}
with
\begin{equation*}
	\begin{split}
		\module{I}^\algebra{B}_{\module{E},\module{F}}
		= &\Span_\field{k}\left\{
		(x_0 \tensor y,0,0) - (0,x_0 \tensor \iota_\module{F}(y),0)
		\mid x_0 \in \module{E}_\Null, y \in 
		\module{F}_\Wobs
		\right\} \\
		&+ \Span_\field{k}\left\{
		(x \tensor y_0,0,0) - (0,0,\iota_\module{E}(x) \tensor y_0)
		\mid x \in \module{E}_\Wobs, y_0 \in \module{F}_\Null
		\right\}.
	\end{split}
\end{equation*}
Note that the second and third term in $(\module{F} \strtensor[\algebra{B}] \module{F})_\Wobs$
directly vanish after reduction.
Then the obvious map $\module{E}_\Wobs \tensor[\algebra{B}_\Wobs] \module{F}_\Wobs \to (\module{E} \strtensor[\algebra{B}] \module{F})_\red$, obtained ny mapping into the first component, yields the desired isomorphism.
\end{proof}

For embedded strong constraint algebras and modules note again that reduction will not be compatible with many constructions, since the embedding of $\injstrConAlg$ into $\strConAlg$ or of $\injstrConBimod$ into $\strConBimod$
will, in general, not be monoidal and not preserve colimits.

%% file: projective-modules.tex
In classical geometry projective modules over the algebra of functions play an important role since they can be identified with vector bundles over smooth manifolds and thus serve as the algebraic description of vector bundles.
We will examine the constraint analogue of this relationship in \autoref{sec:SectionsOfConVect}.
From an algebraic point of view projective modules can be understood as a slight generalization of the concept of free modules.
Therefore we will investigate free (strong) constraint modules in \autoref{sec:FreeConMod} and  \autoref{sec:FreeStrConMod} before we focus on projective (strong) constraint modules in \autoref{sec:ProjectiveConModules} and \autoref{sec:ProjStrConMod}.
It will turn out that projective (strong) constraint modules can be characterized in several different ways similarly to the classical situation: by a lifting property, as direct sums of free modules, or by a dual basis lemma.

\subsection{Free Constraint Modules}
	\label{sec:FreeConMod}
	
As a first important family of constraint modules we will introduce free modules in this section.
Morally, these should be constraint modules with a constraint basis.
For this we need to specify a category of objects of potential bases.
We start with the obvious choice of constraint sets and the forgetful functor
$\functor{U} \colon \ConRMod{\algebra{A}} \to \ConSet$.
Then we search for a left adjoint to this.
In the following we use brackets in the exponent of $\algebra{A}^{(M)}$ to indicate the use of direct sums instead of products indexed by $M$.

\begin{proposition}[$\ConSet$-free constraint module]
	\label{prop:ConSetFreeConMod}
	\index{free!constraint module}
Let $\algebra{A} \in \ConAlg$ be a constraint algebra.
\begin{propositionlist}
	\item For every constraint set $M \in \ConSet$ setting
	\glsadd{freeModule}
	\begin{equation}
	\begin{split}
		(\algebra{A}^{(M)})_\Total 
		&\coloneqq \algebra{A}_\Total^{(M_\Total)},\\
		(\algebra{A}^{(M)})_\Wobs 
		&\coloneqq \algebra{A}_\Wobs^{(M_\Wobs)},\\
		(\algebra{A}^{(M)})_\Null
		&\coloneqq \Big\lbrace x \in \algebra{A}_\Wobs^{(M_\Wobs)} \Bigm| \forall m \in 
		M_\Wobs \colon \sum_{n \sim_M m} x^n \in \algebra{A}_\Null \Big\rbrace,
	\end{split}
	\end{equation}
	together with the map
	$\iota_{\algebra{A}^{(M)}} \colon 
	(\algebra{A}^{(M)})_\Wobs \to 
	(\algebra{A}^{(M)})_\Total$
	given by
	\begin{equation}\label{prop:ConSetFreeConMod_Iota}
		\iota_{\algebra{A}^{(M)}}\Big( \sum_{m \in M_\Wobs}  
		b_m^\Wobs x^m\Big)
		\coloneqq \sum_{m \in M_\Wobs} b_{\iota_M(m)}^\Total x^m,
	\end{equation}
	defines a constraint right $\algebra{A}$-module.
	Here $b^\Total_m$ and $b^\Wobs_m$ denote the basis 
	elements of 
	the free modules $\algebra{A}_\Total^{(M_\Total)}$
	and $\algebra{A}_\Wobs^{(M_\Wobs)}$, respectively, and by
	$x^m$ we denote the corresponding coefficients.
	\item For every constraint set $M \in \ConSet$
	the constraint right $\algebra{A}$-module 
	$\algebra{A}^{(M)}$
	satisfies the following universal property:
	For every $\module{E} \in \ConRMod{\algebra{A}}$
	and $f \colon M \to \module{E}$
	there exists a unique morphism
	$\Phi \colon \algebra{A}^{(M)} \to \module{E}$
	of constraint right $\algebra{A}$-modules such that
	\begin{equation}
		\begin{tikzcd}
			\algebra{A}^{(M)}
			\arrow[r,"\Phi"]
			&\module{E} \\
			M
			\arrow[u,"i"]
			\arrow[ur,"f"{swap}]
			&
		\end{tikzcd}
	\end{equation}
	commutes, where $i \colon M \to \algebra{A}^{(M)}$
	is given by
	$i_{\Total/\Wobs}(m) \coloneqq b^{\Total/\Wobs}_m$.
	\item The functor $\functor{F} \colon \ConSet \to 
	\ConRMod{\algebra{A}}$
	given by
	\begin{equation}
		\functor{F}(M) \coloneqq \algebra{A}^{(M)}
	\end{equation}
	on objects and
	\begin{equation}
		\functor{F}(f) \colon \algebra{A}^{(M)} \to 
		\algebra{A}^{(N)},
		\qquad
		\functor{F}(f)(b^{\Total/\Wobs}_m) \coloneqq 
		b^{\Total/\Wobs}_{f(m)}
	\end{equation}
	on morphisms is left adjoint to the forgetful functor
	$\functor{U} \colon \ConRMod{\algebra{A}} \to \ConSet$.
\end{propositionlist}
\end{proposition}

\begin{proof}
The first part is a simple check of 
\autoref{def:ConstraintModules}.
For the second part note that $i$ is indeed a map of 
constraint 
sets, since for $x \sim_M y$ and $m \in M_\Wobs$ arbitrary we 
have
\begin{equation*}
	\sum_{n \sim_M m} (i_\Wobs(x) - i_\Wobs(y))_n
	= \sum_{n \sim_M m} ( b^\Wobs_x- b^\Wobs_y)_n
	= \sum_{n \sim_M m} \delta_{xn} - \delta_{yn}
	= \begin{cases}
		\delta_{xx} - \delta_{yy}=0 &\text{ if } m \sim_M x,\\
		0&\text{ else}.
	\end{cases}
\end{equation*}
Since $\algebra{A}_\Wobs^{(M_\Wobs)}$
and $\algebra{A}_\Total^{(M_\Total)}$
are free modules we get by the classical universal properties 
module morphisms
$\Phi_\Wobs \colon \algebra{A}_\Wobs^{(M_\Wobs)} \to 
\module{E}_\Wobs$
and
$\Phi_\Total \colon \algebra{A}_\Total^{(M_\Total)} \to 
\module{E}_\Total$.
Moreover, we have
$\iota_\module{E} \circ f_\Wobs 
= \iota_\module{E} \circ \Phi_\Wobs \circ i_\Wobs$
and
$\iota_\module{E} \circ f_\Wobs 
= f_\Total \circ \iota_M
= \Phi_\Total \circ i_\Total \circ \iota_M
= \Phi_\Total \circ \iota_{\algebra{A}^{(M)}} \circ i_\Wobs$,
and hence the universal property of 
$\algebra{A}_\Wobs^{(M_\Wobs)}$
together with the injectivity of $i_\Wobs$
ensures
$\iota_\module{E} \circ \Phi_\Wobs = \Phi_\Total \circ 
\iota_{\algebra{A}^{(M)}}$.
To show that $\Phi_\Wobs$ preserves the $\NULL$-component let
$x \in (\algebra{A}^{(M)})_\Null$ be given.
Then
\begin{align*}
	\Phi_\Wobs(x)
	&= \sum_{m \in M_\Wobs} \Phi_\Wobs(b^\Wobs_m) x^m \\
	&= \sum_{[m] \in M_\Wobs / \sim_M} \sum_{n \sim_M m} 
	f_\Wobs(n)x^n \\
	&= \sum_{[m] \in M_\Wobs / \sim_M} \sum_{n \sim_M m} 
	\left(f_\Wobs(n) - f_\Wobs(m)\right)x^n + f_\Wobs(m)x^n\\
	&= \sum_{[m] \in M_\Wobs / \sim_M} \sum_{n \sim_M m} 
	\underbrace{\left(f_\Wobs(n) - f_\Wobs(m)\right)}_{\in 
		\module{E}_\Null}x^n + 
	f_\Wobs(m) \cdot \sum_{[m] \in M_\Wobs / \sim_M} 
	\underbrace{\sum_{n \sim_M m}x^n}_{\in \algebra{A}_\Null}.
\end{align*}
Thus $\Phi \coloneqq (\Phi_\Total,\Phi_\Wobs)$ is a constraint 
morphism.
Finally, the uniqueness is clear since the uniqueness of 
$\Phi_\Total$ and $\Phi_\Wobs$ is guaranteed by the classical 
universal property.
The third part is just the usual reformulation of universal 
properties via adjoint functors.
\end{proof}

By \eqref{prop:ConSetFreeConMod_Iota} it is clear that an injective 
constraint set $M$ yields an injective constraint algebra 
$\algebra{A}^{(M)}$.
	
\begin{definition}[$\ConSet$-free constraint module]
	\label{def:ConSetFreeConMod}
	Let $\algebra{A} \in \ConAlg$ be a constraint algebra.
	A constraint $\algebra{A}$-module
	$\module{E} \in \ConRMod{\algebra{A}}$ is called 
	\emph{$\ConSet$-free} if there 
	exists a constraint set $M \in \ConSet$ such that
	$\module{E} \simeq \algebra{A}^{(M)}$.
\end{definition}

Though this yields a conceptually clear notion of free constraint 
modules, it is sort of clumsy to work with, since the 
$\NULL$-component is defined using an equivalence relation on 
$M_\Wobs$.
To remedy this deficiency we can use constraint index sets instead:

\begin{proposition}[$\ConIndSet$-free constraint module]
	\label{prop:ConIndSetFreeModule}
	\index{free!constraint module}
Let $\algebra{A} \in \ConAlg$ be a constraint algebra.
\begin{propositionlist}
	\item For every constraint index set $M \in \ConIndSet$ setting
	\begin{equation}
	\begin{split}
		(\algebra{A}^{(M)})_\Total 
		&\coloneqq \algebra{A}_\Total^{(M_\Total)},\\
		(\algebra{A}^{(M)})_\Wobs
		&\coloneqq \algebra{A}_\Wobs^{(M_\Wobs)}, \\
		(\algebra{A}^{(M)})_\Null 
		&\coloneqq 	\algebra{A}_\Null^{(M_\Wobs\setminus M_\Null)} \oplus 
		\algebra{A}_\Wobs^{(M_\Null)},
	\end{split}
	\end{equation}
	together with the map
	$\iota_{\algebra{A}^{(M)}} \colon 
	(\algebra{A}^{(M)})_\Wobs \to 
	(\algebra{A}^{(M)})_\Total$
	given by
	\begin{equation}\label{prop:ConIndSetFreeModule_Iota}
		\iota_{\algebra{A}^{(M)}}\Big( \sum_{m \in M_\Wobs}  
		b_m^\Wobs x^m\Big)
		\coloneqq \sum_{m \in M_\Wobs} b_{\iota_M(m)}^\Total x^m,
	\end{equation}
	defines a constraint right $\algebra{A}$-module.
	Here $b^\Total_m$ and $b^\Wobs_m$ denote the basis 
	elements of 
	the free modules $\algebra{A}_\Total^{(M_\Total)}$
	and $\algebra{A}_\Wobs^{(M_\Wobs)}$, respectively.
	\item For every constraint index set $M \in \ConIndSet$
	the constraint right $\algebra{A}$-module 
	$\algebra{A}^{(M)}$
	satisfies the following universal property:
	For every $\module{E} \in \ConRMod{\algebra{A}}$
	and $f \colon M \to \module{E}$
	there exists a unique morphism
	$\Phi \colon \algebra{A}^{(M)} \to \module{E}$
	of constraint right $\algebra{A}$-modules such that
	\begin{equation}
		\begin{tikzcd}
			\algebra{A}^{(M)}
			\arrow[r,"\Phi"]
			&\module{E} \\
			M
			\arrow[u,"i"]
			\arrow[ur,"f"{swap}]
			&
		\end{tikzcd}
	\end{equation}
	commutes, where $i \colon M \to \algebra{A}^{(M)}$
	is given by
	$i_{\Total/\Wobs}(m) \coloneqq b^{\Total/\Wobs}_m$.
	\item The functor $\functor{F} \colon \ConIndSet \to 
	\ConRMod{\algebra{A}}$
	given by
	\begin{equation}
		\functor{F}(M) \coloneqq \algebra{A}^{(M)}
	\end{equation}
	on objects and
	\begin{equation}
		\functor{F}(f) \colon \algebra{A}^{(M)} \to 
		\algebra{A}^{(N)},
		\qquad
		\functor{F}(f)(b^{\Total/\Wobs}_m) \coloneqq 
		b^{\Total/\Wobs}_{f(m)}
	\end{equation}
	on morphisms is left adjoint to the forgetful functor
	$\functor{U} \colon \ConRMod{\algebra{A}} \to \ConIndSet$.
\end{propositionlist}
\end{proposition}

\begin{proof}
The proof works similar to that of \autoref{prop:ConSetFreeConMod}:
The first part is a simple check of the definition of constraint right $\algebra{A}$-modules.
For the second part note that $i$ is indeed a map of constraint index sets.
Then $\Phi_\Total$ and $\Phi_\Wobs$ are given by the unique morphisms that exist by the universal properties
of $\algebra{A}_\Total^{(M_\Total)}$
and $\algebra{A}_\Wobs^{(M_\Wobs)}$,
and $\Phi_\Wobs$ preserves the $\NULL$-component, since
\begin{align*}
	\Phi_\Wobs\Big(\sum_{m \in M_\Wobs\setminus M_\Null} b^\Wobs_m x^m + \sum_{m \in M_\Null} b^\Wobs_m x^m\Big)
	&= \sum_{m \in M_\Wobs\setminus M_\Null} \Phi_\Wobs(b^\Wobs_m) x^m + \sum_{m \in M_\Null} \Phi(b^\Wobs_m) x^m\\
	&= \sum_{m \in M_\Wobs\setminus M_\Null} f_\Wobs(m) \underbrace{x^m}_{\in \algebra{A}_\Null}
	+ \sum_{m \in M_\Null} \underbrace{f_\Wobs(m)}_{\in \module{E}_\Null} x^m
\end{align*}
if $x^m \in \algebra{A}_\Wobs$ for all $m \in M_\Wobs$ and $x^m \in \algebra{A}_\Null$ for all $m \in M_\Null$.
The third part follows again by abstract nonsense.
\end{proof}

As for $\ConSet$-free modules, we also get by 
\eqref{prop:ConIndSetFreeModule_Iota} that an embedded constraint 
index set $M$ yields an embedded constraint module 
$\algebra{A}^{(M)}$.

\begin{definition}[$\ConIndSet$-Free constraint module]
	\label{def:ConIndFreeConMod}
	Let $\algebra{A} \in \ConAlg$ be a constraint algebra.
	A constraint $\algebra{A}$-module
	$\module{E} \in \ConRMod{\algebra{A}}$ is called 
	\emph{$\ConIndSet$-free} if there 
	exists a constraint index set $M \in \ConIndSet$ such that
	$\module{E} \simeq \algebra{A}^{(M)}$.
	Every such $M$ is called a \emph{constraint basis} of 
	$\module{E}$, and if $M$ is finite we call $\algebra{A}^{(M)}$
	\emph{finitely generated free}.
\end{definition}

\begin{example}
Every constraint $\field{K}$-vector space is a free strong constraint $\field{K}$-module by
\autoref{prop:ConVectSpacesAreFree}, and the notions of bases agree.
\end{example}

While the categories $\ConSet$ and $\ConIndSet$ are not equivalent, 
cf. \autoref{rem:ConSetVSConIndSet}, the respective free modules are 
closely related, as the next results show.

\begin{lemma}[From $\ConSet$-free to $\ConIndSet$-free modules]
Let $\algebra{A} \in \ConAlg$ be a constraint algebra and
$M \in \ConSet$.
\begin{lemmalist}
	\item There exist $\hat{M} \in \injConIndSet$ and a regular 
	epimorphism
	$\Phi \colon \algebra{A}^{(\hat{M})} \to \algebra{A}^{(M)}$.
	\item If $M \in \injConSet$, then $\Phi$ can be chosen to be an 
	isomorphism.
\end{lemmalist}
\end{lemma}

\begin{proof}
Choose a splitting $s \colon M_\red \to M_\Wobs$
of the projection
$\pr_M \colon M_\Wobs \to M_\red$.
Then we define $\hat{M} \in \injConIndSet$ by
\begin{align*}
	\hat{M}_\Total &\coloneqq M_\Wobs \sqcup
	( M_\Total \setminus \iota_M(M_\Wobs)),\\
	\hat{M}_\Wobs &\coloneqq M_\Wobs,\\
	\hat{M}_\Null &\coloneqq M_\Wobs \setminus \image s
\end{align*}
and denote $q \coloneqq s \circ \pr_M \colon M_\Wobs \to \image s$.
Now define $\Phi \colon \algebra{A}^{(\hat{M})} \to 
\algebra{A}^{(M)}$ by
\begin{align*}
	\Phi_\Total(x) &\coloneqq
	\sum_{i \in \image s} b_{\iota_M(i)} x^i
	+ \sum_{i \in M_\Wobs \setminus \image s}
	\big(b_{\iota_M(i)} - b_{\iota_M(q(i))} \big) x^i
	+ \sum_{i \in M_\Total \setminus \iota_M(M_\Wobs)}
	b_i x^i, \\
	\Phi_\Wobs(x) &\coloneqq
	\sum_{i \in \image s} b_i x^i 
	+ \sum_{i \in M_\Wobs \setminus \image s}
	\big( b_i - b_{q(i)} \big) x^i.
\end{align*}
To see that $\Phi$ is indeed a constraint morphism we compute
\begin{equation*}
	\Phi_\Total (\iota_{\algebra{A}^{(\hat{M})}}(b_i))
	= \Phi_\Total (b_{\iota_M(i)})
	= \begin{cases}
		b_{\iota_M(i)} & \text{ if } i \in \image s \\
		b_{\iota_M(i)} - b_{\iota_M(q(i))} & \text{ if } i \in 
		M_\Wobs \setminus \image s
	\end{cases}
\end{equation*}
and
\begin{equation*}
	\iota_{\algebra{A}^{(M)}}(\Phi_\Wobs(b_i))
	= \begin{cases}
		\iota_{\algebra{A}^{(M)}}(b_i) & \text{ if } i \in 
		\image s \\
		\iota_{\algebra{A}^{(M)}}(b_i - b_{q(i)}) & \text{ 
		if } i \in M_\Wobs \setminus \image s
	\end{cases}
	= \begin{cases}
		b_{\iota_M(i)} & \text{ if } i \in \image s \\
		b_{\iota_M(i)} - b_{\iota_M(q(i))} & \text{ if } i \in 
		M_\Wobs \setminus \image s.
	\end{cases}
\end{equation*}
Moreover, for 
$x \in \algebra{A}^{(\hat{M})}_\Null$
we know $x^i \in \algebra{A}_\Null$ if $i \in \image s$ 
and hence for fixed $j \in M_\Wobs$ we get
\begin{align*}
	\sum_{i \sim_M j} \big( \Phi_\Wobs(x)\big)^i
	&= \sum_{i \sim_M j \cap (M_\Wobs \setminus \image s)} 
	x^i
	+ \big(\Phi_\Wobs(x)\big)^{q(j)} \\
	&= \sum_{i \in [j] \cap (M_\Wobs \setminus \image s)} 
	x^i
	+ x^{q(j)}
	- \sum_{i \in [j] \cap (M_\Wobs \setminus \image s)} 
	x^i \\
	&= x^{q(j)} \in \algebra{A}_\Null,
\end{align*}
where $[j]$ denotes the equivalence class of $j$.
By the definition of $\Phi$ it is clear that $\Phi_\Total$ and 
$\Phi_\Wobs$ are surjective.
Additionally, we have
$\Phi_\Wobs(\algebra{A}^{(\hat{M})}_\Null) = 
\algebra{A}^{(M)}_\Null$, because for
$x \in \algebra{A}^{(M)}_\Null$
we can define
$y_\Wobs \in \algebra{A}^{(\hat{M})}_\Null$
by $y^i \coloneqq \sum_{j \in [i]} x^j$
for $i \in \image s$ and
$y^i \coloneqq x^i$ for
$i \in M_\Wobs \setminus \image s$.
Then by construction $y \in \algebra{A}^{(\hat{M})}_\Null$
and $\Phi_\Wobs (y) = x$.
Thus $\Phi$ is a regular epimorphism.
Finally, $\Phi_\Wobs$ is always injective.
In case that $M \in \injConSet$, i.e. $\iota_M$ is injective, we also 
get that $\Phi_\Total$ is injective, and therefore $\Phi$ is an 
isomorphism of constraint modules.
\end{proof}

We see that at least every $\injConSet$-free module is also 
$\injConIndSet$-free.
However, as the next lemma shows even in the embedded situation this 
correspondence is not perfect.

\begin{lemma}[From $\ConIndSet$-free to $\ConSet$-free modules]
Let $\algebra{A} \in \ConAlg$ be a constraint algebra and
$M \in \ConIndSet$.
\begin{lemmalist}
	\item There exist $\check{M} \in \injConSet$ and a regular 
	epimorphism
	$\Phi \colon \algebra{A}^{(\check{M})} \to \algebra{A}^{(M)}$.
	\item If $M \in \injConIndSet$ and $M_\Wobs \neq M_\Null$,
	then $\check{M}$ can be chosen in such a way that $\Phi$ is an isomorphism.
	\item If $M \in \injConIndSet$ and $M_\Wobs = M_\Null$,
	then there exists an isomorphism between
	$\algebra{A}^{(\check{M})}$
	and
	$\algebra{A}^{(M \sqcup \{\pt\})}$.
\end{lemmalist}
\end{lemma}

\begin{proof}
First, assume that $M_\Wobs \setminus M_\Null$ is non-empty.
Then choose an element $k \in M_\Wobs \setminus M_\Null$
and define
$\check{M}_\Total \coloneqq M_\Wobs \sqcup ( M_\Total 
\setminus \iota_M(M_\Wobs))$,
$\check{M}_\Wobs \coloneqq M_\Wobs$
and $\sim_{\check{M}}$
as the equivalence relation generated by
$i \sim k$ for all $i \in M_\Null$.
Now define
$\Phi \colon \algebra{A}^{(\check{M})} \to \algebra{A}^{(M)}$
by
\begin{align*}
	\Phi_\Total(x) &\coloneqq
	\sum_{i \in M_\Null} \big(b_{\iota_M(i)} + b_{\iota_M(k)} 
	\big) x^i
	+ \sum_{i \in M_\Wobs \setminus M_\Null}
	b_{\iota_M(i)}	x^i
	+ \sum_{i \in M_\Total \setminus \iota_M(M_\Wobs)}
	b_i x^i,\\
	\Phi_\Wobs(x) &\coloneqq
	\sum_{i \in M_\Null} \big( b_i + b_{k} \big) x^i
	+ \sum_{i \in M_\Wobs \setminus M_\Null} b_i x^i.
\end{align*}
To see that $\Phi$ is indeed a constraint morphism we compute
\begin{equation*}
	\Phi_\Total (\iota_{\algebra{A}^{(\check{M})}}(b_i))
	= \Phi_\Total (b_i)
	= \begin{cases}
		b_{\iota_M(i)} & \text{ if } i \in M_\Wobs \setminus 
		M_\Null \\
		b_{\iota_M(i)} + b_{\iota_M(k)} & \text{ if } i \in 
		M_\Null
	\end{cases}
\end{equation*}
and
\begin{equation*}
	\iota_{\algebra{A}^{(M)}}(\Phi_\Wobs(b_i))
	= \begin{cases}
		\iota_{\algebra{A}^{(M)}}(b_i) & \text{ if } i 
		\in M_\Wobs \setminus M_\Null \\
		\iota_{\algebra{A}^{(M)}}(b_i + b_{k}) & 
		\text{ if } i \in M_\Null
	\end{cases}
	= \begin{cases}
		b_{\iota_M(i)} & \text{ if } i \in M_\Wobs \setminus 
		M_\Null \\
		b_{\iota_M(i)} + c_{\iota_M(k)} & \text{ if } i \in 
		M_\Null.
	\end{cases}
\end{equation*}
Moreover, for $x \in \algebra{A}^{(\check{M})}_\Null$
and $i \in M_\Wobs \setminus M_\Null$ we have
$\big(\Phi_\Wobs(x)\big)^i = x^i \in \algebra{A}_\Null$
if $i \neq k$, and if $i=k$ we have
$\big(\Phi_\Wobs(x)\big)^i = x^{k}
+ \sum_{j \in M_\Null} x^j \in \algebra{A}_\Null$.
Surjectivity of $\Phi_\Total$ and $\Phi_\Wobs$
follow directly.
Finally, for $x \in \algebra{A}^{(M)}_\Null$
we can define $y \in \algebra{A}^{(\check{M})}_\Null$
by $y^{k} \coloneqq x^{k} - \sum_{j \in M_\Null} x^j$
and $y^i \coloneqq x^i$ for $i \neq i_0$.
Then by construction $y \in \algebra{A}^{(\check{M})}_\Null$.
Moreover, $\Phi_\Wobs (y) = x$.
This shows that $\Phi$ is a regular epimorphism.
If $\iota_M$ is injective then one can check that $\Phi_\Total$
and $\Phi_\Wobs$ are injective too and therefore $\Phi$ is an 
isomorphism of constraint modules.

If $M_\Wobs = M_\Null$ define
$M^\prime \coloneqq
(M_\Total \sqcup \{\pt\},\,
M_\Wobs \sqcup \{ \pt \},\,
M_\Null)$.
Then
\begin{equation*}
	\Psi_\Total(x) \coloneqq
	\sum_{i \in M_\Total}
	b_i	x^i, \qquad
	\Psi_\Wobs(x) \coloneqq
	\sum_{i \in M_\Wobs}
	b_i x^i
\end{equation*}
defines a regular epimorphism $\Psi \colon \algebra{A}^{(M^\prime)} \to 
\algebra{A}^{(M)}$.
Applying the first part to $M^\prime$ and composing 
the regular epimorphisms we obtain a suitable $\Phi$ in this case as well.
\end{proof}

Though the notions of $\injConIndSet$- and $\injConSet$-free modules 
are not equivalent, in most cases this difference will not be crucial.
Since the interpretation of the generating set is, especially in 
geometric situations, more intuitive we will in the following mostly 
consider $\injConIndSet$-free modules.

Let us now investigate how free constraint modules behave with 
respect to some important constructions we have introduced for 
general modules before.

\begin{proposition}
	\label{lem:SumsOfFreeModules}
Let $\algebra{A} \in \ConAlg$ and $M,N \in \ConIndSet$ be given.
\begin{propositionlist}
	\item We have
	\begin{equation}
		\algebra{A}^{(M)} \oplus \algebra{A}^{(N)} \simeq \algebra{A}^{(M \sqcup N)}.
	\end{equation}
	\item We have
	\begin{equation}
		\algebra{A}^{(M)} \tensor \algebra{A}^{(N)} \simeq \algebra{A}^{(M \tensor N)}.
	\end{equation}
\end{propositionlist}
\end{proposition}

\begin{proof}
The first part follows directly from the fact that $\sqcup$ is the coproduct in 
$\ConIndSet$ and the free functor 
$\functor{F} \colon \ConIndSet \to \ConRMod{\algebra{A}}$
is left adjoint and hence preserves colimits.
For the second part the $\TOTAL$- and $\WOBS$-components are clear.
For the $\NULL$-component note that
\begin{equation*}
	(\algebra{A}_\Null^{(M_\Wobs\setminus M_\Null)} \oplus \algebra{A}_\Wobs^{(M_\Null)}) \tensor \algebra{A}_\Wobs^{(N_\Wobs)}
	\simeq \algebra{A}_\Null^{(M_\Wobs \times N_\Wobs\setminus N_\Null)}
	\oplus \algebra{A}_\Wobs^{(M_\Wobs \times N_\Null)}
\end{equation*}
and
\begin{equation*}
	\algebra{A}_\Wobs^{(M_\Wobs)} \tensor	(\algebra{A}_\Null^{(N_\Wobs\setminus N_\Null)} \oplus \algebra{A}_\Wobs^{(N_\Null)}) 
	\simeq	\algebra{A}_\Wobs^{(N_\Wobs \times M_\Null)} 
	\oplus \algebra{A}_\Null^{(N_\Wobs \times M_\Wobs\setminus M_\Null)}.
\end{equation*}
This leads to
\begin{equation*}
	(\algebra{A}^{(M)} \tensor \algebra{A}^{(N)})_\Null
	= \algebra{A}_\Null^{(M_\Wobs\setminus M_\Null) \times (N_\Wobs\setminus N_\Null)}
	\oplus \algebra{A}_\Wobs^{(M_\Wobs\times N_\Null) \cup (M_\Null\times M_\Wobs)}
\end{equation*}
as expected.
\end{proof}

In classical algebra we know that duals of finitely generated free 
modules are again free.
This is no longer true for free constraint modules.

\begin{proposition}[Duals of free constraint modules]
	\label{prop:DualsOfConFreeModules}
	\index{dual module}
Let $\algebra{A} \in \injConAlg$ be given.
For finite $n \in \injConIndSet$ we have
\begin{equation}
\begin{split}
	\label{eqprop:DualsOfConFreeModules}
	(\algebra{A}^{n})^*_\Total
	&\simeq \algebra{A}_\Total^{n_\Total}, \\
	(\algebra{A}^{n})^*_\Wobs
	&\simeq \algebra{A}_\Null^{n_\Null} \oplus 
	\algebra{A}_\Wobs^{n_\Wobs - n_\Null} \oplus 
	\algebra{A}_\Total^{n_\Total- n_\Wobs}, \\
	(\algebra{A}^{n})^*_\Null
	&\simeq \algebra{A}_\Null^{n_\Wobs} \oplus 
	\algebra{A}_\Total^{n_\Total - n_\Wobs}.
\end{split}
\end{equation}
\end{proposition}

\begin{proof}
From classical algebra we know that
$(\algebra{A}_\Total^{n_\Total})^*$ is free with dual basis
$(b^i)_{i = 1}^{n_\Total}$.
Let $\alpha = \sum_{i=1}^{n_\Total}\alpha_i b^i \in 
(\algebra{A}^n)^*_\Wobs$ be given.
Then from $\alpha(\algebra{A}_\Wobs^{n_\Wobs}) \subseteq \algebra{A}_\Wobs$
it follows that
$\alpha_i \in \algebra{A}_\Wobs$ for all $i = 1,\dotsc, n_\Wobs$.
Since $\alpha(\algebra{A}_\Wobs^{n_\Null}) \subseteq 
\algebra{A}_\Null$ we additionally get
$\alpha_i \in \algebra{A}_\Null$
for all $i=1,\dotsc, n_\Null$.
This shows the $\WOBS$-component.
Now let $\alpha \in (\algebra{A}^n)^*_\Null$ be given.
Then the $\NULL$-component follows from $\alpha(\algebra{A}_\Wobs^{n_\Wobs}) \subseteq 
\algebra{A}_\Null$.
\end{proof}

For non-embedded constraint algebras $(\algebra{A}^n)^*$ would look 
more complicated, since the $\WOBS$-\-com\-po\-nent then consists of pairs 
of functionals.
Modules of this particular form will again show up when we look at free strong constraint modules.

\subsubsection{Reduction}

The reduction of $\ConSet$- and $\ConIndSet$-free modules yields classical free modules
over the reduced algebra:

\begin{proposition}[Reduction of free constraint modules]
	\label{prop:ReductionFreeConMod}
	\index{reduction!free constraint modules}
Let $\algebra{A} \in \ConAlg$ be a constraint algebra.
\begin{propositionlist}
	\item There exists a natural isomorphism	making the diagram
	\begin{equation}
		\begin{tikzcd}
			\ConSet
			\arrow[r,"\functor{F}"]
			\arrow[d,"\red"{swap}]
			&\ConRMod{\algebra{A}}
			\arrow[d,"\red"]\\
			\Sets
			\arrow[r,"\functor{F}"]
			&\Modules_{\algebra{A}_\red}	
		\end{tikzcd}
	\end{equation}
	commute, where $\functor{F}$ denotes the respective free 
	construction.
	In particular we have
	\begin{equation}
		(\algebra{A}^{(M)})_\red \simeq (\algebra{A}_\red)^{(M_\red)}
	\end{equation}
	for all $M \in \ConSet$.
	\item There exists a natural isomorphism making the diagram
	\begin{equation}
		\begin{tikzcd}
			\ConIndSet
			\arrow[r,"\functor{F}"]
			\arrow[d,"\red"{swap}]
			&\ConRMod{\algebra{A}}
			\arrow[d,"\red"]\\
			\Sets
			\arrow[r,"\functor{F}"]
			&\Modules_{\algebra{A}_\red}	
		\end{tikzcd}
	\end{equation}
	commute, where $\functor{F}$ denotes the respective free 
	construction.
	In particular we have
	\begin{equation}
		(\algebra{A}^{(M)})_\red \simeq (\algebra{A}_\red)^{(M_\red)}
	\end{equation}
	for all $M \in \ConIndSet$.
\end{propositionlist}
\end{proposition}

\begin{proof}
	For $M \in \ConSet$ define
	$\eta_M \colon (\algebra{A}^{(M)})_\red \to 
	(\algebra{A}_\red)^{(M_\red)}$
	by
	\begin{equation*}
		\eta_M\Big( \Big[ \sum_{m \in M_\Wobs} b_m x^m 
		\Big] \Big)
		\coloneqq \sum_{[m] \in M_\red} b_{[m]}
		\big[ \sum_{m'\sim_M m} x^{m'} \big],
	\end{equation*}
	where $b_{[m]}$ denotes the basis element of 
	$(\algebra{A}_\red)^{(M_\red)}$ corresponding to the equivalence 
	class $[m]$.
	This map is clearly well-defined on $(\algebra{A}^{(M)})_\red$
	and injective.
	For surjectivity let
	$\sum_{[m] \in M_\red} b_{[m]} [x^{[m]}] \in 
	(\algebra{A}_\red)^{(M_\red)}$ be given.
	Using the axiom of choice, choose a splitting
	$i \colon M_\red \to M_\Wobs$ of the quotient map $M_\Wobs \to 
	M_\red$.
	Then $\sum_{m \in \image(i)} b_m x^{[m]}$
	is a suitable preimage of $\sum_{[m] \in M_\red} b_{[m]} [x^{[m]}]$.
	Naturality can now be checked by a direct computation:
	Let $f \colon M \to N$ be a morphism of constraint sets,
	then
	\begin{align*}
		\left(\functor{F}(f_\red) \circ \eta_M\right)\Big( 
		\Big[ \sum_{m \in M_\Wobs} b_m x^m \Big] \Big)
		&= \functor{F}(f_\red)\Big( \sum_{[m] \in M_\red} b_{[m]}
		\big[ \sum_{m'\sim_M m} x^{m'} \big] \Big)\\
		&= \sum_{[n] \in N_\red} b_{[n]}
		\Big(\sum_{[m] \in f_\red^{-1}([n])}\big[ \sum_{m'\sim_M m} 
		x^{m'}\big] \Big)\\
		&= \sum_{[n] \in N_\red} b_{[n]}
		\Big(\sum_{n' \sim_N n}\big[ \sum_{m \in f^{-1}(n)} 
		x^{m}\big] \Big)\\
		&= \eta_N\Big( \Big[\sum_{n \in n_\Wobs} b_n \sum_{m \in 
			f^{-1(n)}} x^m \Big] \Big)\\
		&= \left(\eta_N \circ \functor{F}(f)_\red\right)\Big( \Big[ 
		\sum_{m \in 
			M_\Wobs} b_m x^m \Big]\Big).
	\end{align*}
We can use the same isomorphism $\eta$ for a constraint index set $M$.
Alternatively, we can see it more directly:
\begin{equation*}
	(\algebra{A}^{(M)})_\red
	= \frac{\algebra{A}_\Wobs^{(M_\Wobs)}}{\algebra{A}_\Null^{(M_\Wobs\setminus M_\Null)} \oplus \algebra{A}_\Wobs^{(M_\Null)}}
	\simeq \frac{\algebra{A}_\Wobs^{(M_\Wobs\setminus M_\Null)}}{\algebra{A}_\Null^{(M_\Wobs\setminus M_\Null)}}
	\simeq \algebra{A}_\red^{(M_\red)}.
\end{equation*}
\end{proof}

\subsection{Free Strong Constraint Modules}
	\label{sec:FreeStrConMod}

For strong constraint modules we will focus on free modules generated by constraint index sets.
Thus we are searching for a left adjoint functor to the forgetful functor
$\functor{U} \colon \strConRMod{\algebra{A}} \to \ConIndSet$
for a fixed strong constraint algebra $\algebra{A}$.
Note that this functor factors through
$\ConRMod{\algebra{A}}$
by first forgetting to constraint $\algebra{A}$-modules and then to their underlying
constraint index sets.
For both of those forgetful functors we have already  found left adjoints in
\autoref{prop:StrongHullModules} and
\autoref{prop:ConIndSetFreeModule}.

\begin{lemma}[$\ConIndSet$-free strong constraint module]
	\index{free!strong constraint module}
Let $\algebra{A} \in \strConAlg$ be a strong constraint algebra.
\begin{lemmalist}
\item For every constraint index set $M \in \ConIndSet$
the strong constraint right $\algebra{A}$-module 
$\big(\functor{U}\left(\algebra{A}\right)^{(M)}\big)^\str$
satisfies the following universal property:
For every $\module{E} \in \strConRMod{\algebra{A}}$
and $f \colon M \to \module{E}$
there exists a unique morphism
$\Phi \colon \big(\functor{U}\left(\algebra{A}\right)^{(M)}\big)^\str \to \module{E}$
of strong constraint right $\algebra{A}$-modules such that
\begin{equation}
\begin{tikzcd}
	\big(\functor{U}\left(\algebra{A}\right)^{(M)}\big)^\str
	\arrow[r,"\Phi"]
	&\module{E} \\
	M
	\arrow[u,"i"]
	\arrow[ur,"f"{swap}]
	&
\end{tikzcd}
\end{equation}
commutes, where $i \colon M \to \big(\functor{U}\left(\algebra{A}\right)^{(M)}\big)^\str$
is given by
$i_{\Total/\Wobs}(m) \coloneqq b^{\Total/\Wobs}_m$.
\item The functor $\functor{F} \colon \ConIndSet \to 
\strConRMod{\algebra{A}}$
given by
\begin{equation}
\functor{F}(M) \coloneqq \Big(\functor{U}\left(\algebra{A}\right)^{(M)}\Big)^\str
\end{equation}
is left adjoint to the forgetful functor $\functor{U} \colon \strConRMod{\algebra{A}} \to \ConIndSet$.
\end{lemmalist}
\end{lemma}

\begin{proof}
	The first and second statement are equivalent by general category theory, while the second part holds
	since $\functor{F}$ is defined as the composition of the left adjoints of
	$\functor{U} \colon \ConIndSet \to \ConRMod{\algebra{A}}$
	and
	$\functor{U} \colon \ConRMod{\algebra{A}} \to \strConRMod{\algebra{A}}$.
\end{proof}

For a strong constraint algebra $\algebra{A}$ we will write
$\algebra{A}^{(M)} \coloneqq \functor{U}(\algebra{A})^{(M)}$ for any 
$M \in \ConIndSet$.
No confusion should arise, since from the context it is clear whether $\algebra{A}$ is
a strong constraint or a plain constraint algebra.

\begin{definition}[$\ConIndSet$-free strong constraint module]
	\label{def:ConIndFreeStrConMod}
	Let $\algebra{A} \in \strConAlg$ be a strong constraint algebra.
	A strong constraint $\algebra{A}$-module
	$\module{E} \in \strConRMod{\algebra{A}}$ is called 
	\emph{$\ConIndSet$-free} if there 
	exists a constraint index set $M \in \ConIndSet$ such that
	$\module{E} \simeq \algebra{A}^{(M)}$.
	Every such $M$ is called a \emph{constraint basis} of 
	$\module{E}$, and if $M$ is finite we call $\algebra{A}^{(M)}$
	\emph{finitely generated free}.
\end{definition}

\begin{remark}
Free constraint modules have been introduced in \cite{menke:2020a,dippell.menke.waldmann:2022a},
but the relation to free strong constraint modules had not been developed.
\end{remark}

In the embedded case $\ConIndSet$-free strong constraint modules take on an easy form.

\begin{lemma}
\label{lem:FreeStrConModules}
Let $\algebra{A} \in \strConAlg$ be a strong constraint algebra 
and $M \in \injConIndSet$.
\begin{lemmalist}
	\item \label{lem:FreeStrConModules_1}
	We have
	\begin{equation}
	\begin{split}
		\algebra{A}^{(M)}_\Total 
		&= \algebra{A}_\Total^{(M_\Total)}, \\
		\algebra{A}^{(M)}_\Wobs 
		&= 
		\iota_{\algebra{A}}\left(\algebra{A}_\Null\right)^{(M_\Total\setminus
		 M_\Wobs)}
		\oplus \frac{\algebra{A}_\Wobs^{(M_\Wobs\setminus M_\Null)} 
		\oplus 
		\iota_{\algebra{A}}\left(\algebra{A}_\Null\right)^{(M_\Wobs 
		\setminus M_\Null)}}{\Span_\field{k}\left\{(x,0) - 
		(0,\iota(x)) \bigm| x \in \algebra{A}_\Null^{(M_\Wobs\setminus 
		M_\Null)}\right\}}
		\oplus \algebra{A}_\Total^{(M_\Null)},\\
		\algebra{A}^{(M)}_\Null 
		&= 
		\iota_{\algebra{A}}\left(\algebra{A}_\Null\right)^{(M_\Total\setminus
		 M_\Wobs)}
		\oplus 
		\iota_{\algebra{A}}\left(\algebra{A}_\Null\right)^{(M_\Wobs\setminus
		 M_\Null)}
		\oplus \algebra{A}_\Total^{(M_\Null)}.
	\end{split}
	\end{equation}
	\item \label{lem:FreeStrConModules_2}
	If additionally $\iota_\algebra{A} \colon \algebra{A}_\Wobs \to 
	\algebra{A}_\Total$
	is injective it holds that
	\begin{equation}
	\begin{split} 
		\algebra{A}^{(M)}_\Total 
		&= \algebra{A}_\Total^{(M_\Total)},\\
		\algebra{A}^{(M)}_\Wobs 
		&= \algebra{A}_\Null^{(M_\Total\setminus M_\Wobs)}
		\oplus \algebra{A}_\Wobs^{(M_\Wobs\setminus M_\Null)} \oplus 
		\algebra{A}_\Total^{(M_\Null)},\\
		\algebra{A}^{(M)}_\Null 
		&= \algebra{A}_\Null^{(M_\Total\setminus M_\Null)}
		\oplus \algebra{A}_\Total^{(M_\Null)}.
	\end{split}
	\end{equation}
\end{lemmalist}
\end{lemma}

\begin{proof}
		The first part follows directly from the construction of $\functor{U}(\algebra{A})^{(M)}$ in 
	\autoref{prop:ConIndSetFreeModule} and the definition of the strong hull
	in \autoref{prop:StrongHullModules}.
	The second part then follows immediately.
\end{proof}

The next result shows that, at least in the finitely generated case, $\injConIndSet$-free strong constraint modules
over an embedded strong constraint algebra are closed under many operations, such as
direct sums, tensor products, strong tensor products and duals.

\begin{proposition}[Duals of free modules]
	\label{prop:DualsOfStrFree}
	\index{dual module}
Let $\algebra{A} \in \injstrConAlg$ be an embedded 
strong constraint algebra and let
$n,m \in \injConIndSet$ be finite.
\begin{propositionlist}
	\item \label{prop:DualsOfStrFree_1}
	$(\algebra{A}^{n})^*$ is free and
	$(\algebra{A}^{n})^* \simeq \algebra{A}^{(n^*)}$.
	\item \label{prop:DualsOfStrFree_2}
	$\algebra{A}^n \oplus \algebra{A}^m$
	is free and
	$\algebra{A}^n \oplus \algebra{A}^m \simeq \algebra{A}^{n \sqcup m}$.
	\item \label{prop:DualsOfStrFree_3}
	$\algebra{A}^n \injtensor[\algebra{A}] \algebra{A}^m$
	is free and
	$\algebra{A}^n \injtensor[\algebra{A}] \algebra{A}^m \simeq \algebra{A}^{n \tensor m}$.
	\item \label{prop:DualsOfStrFree_4}
	$\algebra{A}^n \injstrtensor[\algebra{A}] \algebra{A}^m$
	is free and
	$\algebra{A}^n \injstrtensor[\algebra{A}] \algebra{A}^m \simeq \algebra{A}^{n \strtensor m}$.
	\item \label{prop:DualsOfStrFree_5}
	If $m \subseteq n$ is a constraint index subset, then
	$\algebra{A}^n / \algebra{A}^m$ is free and
	$\algebra{A}^n / \algebra{A}^m \simeq \algebra{A}^{n\setminus m}$.
\end{propositionlist}
\end{proposition}

\begin{proof}
For the $\TOTAL$-component all the above identities hold by the classical theory.
Part \ref{prop:DualsOfStrFree_2} follows from the fact that the free functor
$\injConIndSet \to \injstrConMod_\algebra{A}$
is a left adjoint, and thus preserves colimits.
This also explains \ref{prop:DualsOfStrFree_5}.
For \ref{prop:DualsOfStrFree_3} and \ref{prop:DualsOfStrFree_4}
it is straightforward to check that 
$n \tensor m \ni (i,j) \mapsto b_i \tensor b_j \in \algebra{A}^n \tensor \algebra{A}^m$
and
$n \strtensor m \ni (i,j) \mapsto b_i \tensor b_j \in \algebra{A}^n \strtensor \algebra{A}^m$
fulfil the universal properties of $\algebra{A}^{n \tensor m}$
and $\algebra{A}^{n \strtensor m}$, respectively.
\end{proof}

Recall from \autoref{prop:DualsOfConFreeModules} that duals of free constraint modules are in general not
free again.
For a given strong constraint algebra $\algebra{A} \in \strConAlg$ we can consider the free module
$\functor{U}(\algebra{A})^n$ of its underlying constraint algebra.
Then it turns out that its dual will still not be free as a constraint module,
but it will be free as a \emph{strong} constraint module.
\filbreak
\begin{proposition}
Let $\algebra{A} \in \injstrConAlg$ be an embedded strong constraint algebra and let
$n \in \injConIndSet$ be finite.
\begin{propositionlist}
	\item 	The dual $(\functor{U}(\algebra{A})^{n})^*$ is a free strong constraint $\algebra{A}$-module with
	\begin{equation}
		(\functor{U}(\algebra{A})^{n})^* \simeq \algebra{A}^{(n^*)}.
	\end{equation}
	\item The strong hull $(\functor{U}(\algebra{A})^n)^\str$ is a free strong constraint $\algebra{A}$-module
	with
	\begin{equation}
		(\functor{U}(\algebra{A})^n)^\str \simeq \algebra{A}^n.
	\end{equation}
\end{propositionlist}
\end{proposition}

\begin{proof}
The first part follows directly from \autoref{prop:DualsOfConFreeModules}.
For the second part we have
\begin{align*}
	\big((\functor{U}(\algebra{A})^n)^\str\big)_\Wobs
	&= \algebra{A}_\Wobs^{n_\Wobs} + (\algebra{A}_\Null^{n_\Wobs - n_\Null} \oplus \algebra{A}_\Wobs^{n_\Null})\cdot \algebra{A}_\Total + \algebra{A}_\Total^{n_\Total} \cdot \algebra{A}_\Null \\
	&= \algebra{A}_\Wobs^{n_\Wobs} + (\algebra{A}_\Null^{n_\Wobs - n_\Null} \oplus \algebra{A}_\Total^{n_\Null})
	+ \algebra{A}_\Null^{n_\Total}\\
	&= \algebra{A}_\Wobs^{n_\Wobs - n_\Null} \oplus \algebra{A}_\Total^{n_\Null} \oplus \algebra{A}_\Null^{n_\Total -n_\Wobs}\\
	&= (\algebra{A}^n)_\Wobs.
\end{align*}
A similar computation yields the correct $\NULL$-component.
\end{proof}

\subsubsection{Reduction}

As we expect, free strong constraint modules reduce to free modules over the reduced algebra.

\begin{proposition}[Reduction of free strong constraint modules]
	\label{prop:ReductionFreeStrConMod}
	\index{reduction!free strong constraint modules}
Let $\algebra{A} \in \injstrConAlg$ be an embedded strong constraint algebra.
There exists a natural isomorphism	making the diagram
	\begin{equation}
		\begin{tikzcd}
			\ConIndSet
			\arrow[r,"\functor{F}"]
			\arrow[d,"\red"{swap}]
			&\injstrConRMod{\algebra{A}}
			\arrow[d,"\red"]\\
			\Sets
			\arrow[r,"\functor{F}"]
			&\Modules_{\algebra{A}_\red}	
		\end{tikzcd}
	\end{equation}
	commute, where $\functor{F}$ denotes the respective free 
	construction.
	In particular we have
	\begin{equation}
		(\algebra{A}^{(M)})_\red \simeq (\algebra{A}_\red)^{(M_\red)}
	\end{equation}
	for all $M \in \ConIndSet$.
\end{proposition}

\begin{proof}
We have a canonical isomorphism
\begin{equation*}
	(\algebra{A}^{(M)})_\red
	= \frac{\algebra{A}_\Null^{(M_\Total\setminus M_\Wobs)} \oplus \algebra{A}_\Wobs^{(M_\Wobs\setminus M_\Null)} \oplus \algebra{A}_\Total^{(M_\Null)}}{\algebra{A}_\Null^{(M_\Total\setminus M_\Null)} \oplus \algebra{A}_\Total^{(M_\Null)}}
	\simeq \frac{\algebra{A}_\Wobs^{(M_\Wobs\setminus M_\Null)}}{\algebra{A}_\Null^{(M_\Wobs\setminus M_\Null)}}
	\simeq \algebra{A}_\red^{(M_\red)},
\end{equation*}
for which it is straightforward to see that it forms a natural transformation, see also
\autoref{prop:ReductionFreeConMod}.
\end{proof}

\subsection{Projective Constraint Modules}
	\label{sec:ProjectiveConModules}

Classical projective modules over an algebra $\algebra{A}$ can be described in
several equivalent ways.
They can be understood as projective objects in the abelian category 
$\Modules_{\algebra{A}}$, as direct summands of free modules
or as modules allowing for a dual basis (in the sense of the dual basis lemma).
When defining projective constraint modules it seems most natural to start with
the most abstract, categorical point of view.
An object in a given category is called projective if it satisfies a certain lifting property
with respect to epimorphisms.
Considering the fact, that in the category of constraint modules we have to 
distinguish different kinds of epimorphisms, see \autoref{prop:MonoEpisConModk},
there also exist different notions of projective constraint modules.
As usual we will use the stronger notion of regular epimorphisms.
%
\begin{definition}[Projective module]
	\label{def:ProjectiveModule}
	\index{projective!constraint module}
Let $\algebra{A} \in \ConAlg$ be a constraint algebra.
A constraint $\algebra{A}$-module
$\module{P} \in \ConRMod{\algebra{A}}$ is called
\emph{projective} if for every
$\module{E}, \module{F} \in \ConRMod{\algebra{A}}$,  morphism
$\Psi \colon \module{P} \longrightarrow \module{F}$ and 
regular epimorphism
\index{epimorphism!regular}
$\Phi \colon \module{E} \longrightarrow \module{F}$ there exists a
morphism $\chi \colon \module{P} \longrightarrow \module{E}$ such
that $\Phi \circ \chi = \Psi$.  Diagrammatically:
\begin{equation}
	\begin{tikzcd}
		{ }
		& \module{E}
		\arrow[twoheadrightarrow]{d}{\Phi} \\
		\module{P}
		\arrow{r}{\Psi}
		\arrow[dashed]{ur}{\chi}
		& \module{F}
	\end{tikzcd}
\end{equation}
\end{definition}

\begin{remark}
	Requiring the lifting property only for plain epimorphisms instead of regular ones would yield a too restrictive class of objects.
	To see this, assume that $\module{P}$ has the lifting property with respect to all epimorphisms.
	Consider now $\module{E} \coloneqq (\module{P}_\Total, \module{P}_\Wobs,0)$,
	$\module{F} \coloneqq \module{P}$ and $\Phi = (\id_{\module{P}_\Total},\id_{\module{P}_\Wobs})$.
	By assumption there exists a splitting $\chi_\Wobs$ of $\id_{\module{P}_\Wobs}$, which is then given by
	$\chi_\Wobs = \id_{\module{P}_\Wobs}$, such that $\chi_\Wobs(\module{P}_\Null) \subseteq 0$.
	Hence such projective modules would only allow for trivial $\NULL$-components.
\end{remark}

%

In classical algebra every free module is projective.
The following example shows that this fails in general for constraint modules.

\begin{example}
	Consider the constraint algebra $\mathbb{R} = (\mathbb{R},\mathbb{R},0)$ and the constraint index set
	$M = (\left\{ 0\right\}, \left\{0,1\right\}, \emptyset)$.
	Note that the unique map $\iota_M \colon \left\{0,1\right\} \to \{0\}$ is surjective but not injective.
	Thus we obtain a free constraint $\mathbb{R}$-module
	$\mathbb{R}^{M} \simeq (\mathbb{R}^1,\mathbb{R}^2,0)$ with
	$\iota_{\mathbb{R}^M}(x,y) = x+y$.
	This constraint module is not projective, since for $\mathbb{R}^2 = (\mathbb{R}^2,\mathbb{R}^2,0)$
	the constraint morphism $\Phi = (\iota_{\mathbb{R}^M},\id_{\mathbb{R}^2}) \colon \mathbb{R}^2 \to \mathbb{R}^M$
	is a regular epimorphism, but there cannot exist a constraint splitting $\chi$ of $\Phi$, because such a splitting would fulfil
	$\chi_\Total \circ \iota_{\mathbb{R}^M} = \id_{\mathbb{R}^2} \circ \chi_\Wobs = \id_{\mathbb{R}^2}$, in conflict with the fact
	that $\iota_{\mathbb{R}^M}$ is not injective.
\end{example}

In the above example the projectivity of $\algebra{A}^{(M)}$ fails due to the non-injectivity of $\iota_{\algebra{A}^{(M)}}$.
For general $\algebra{A} \in \ConAlg$ and $M \in \ConIndSet$ the free module
$\algebra{A}^{(M)}$ has non-injective $\iota_{\algebra{A}^{(M)}}$.
But if both $\algebra{A}$ and $M$ are embedded,
the corresponding free modules are indeed projective:

\begin{lemma}
	\label{lem:FreeConModulesAreProjective}
Let $\algebra{A} \in \injConAlg$ be an embedded constraint algebra.
For every $M \in \injConIndSet$ the free constraint module $\algebra{A}^{(M)}$ is projective.
\end{lemma}

\begin{proof}
	Let $\algebra{A}^{(M)}$ be a free constraint module with
	$M \in \injConIndSet$.  Suppose the following morphisms are given
	\begin{equation*}
		\begin{tikzcd}
			{ }
			& \module{E}
			\arrow[twoheadrightarrow]{d}{\Phi} \\
			\algebra{A}^{(M)}
			\arrow{r}{\Psi}
			& \module{F}
		\end{tikzcd}
	\end{equation*}
	with $\Phi$ a regular epimorphism. 
	Since $\Phi$ and $\Psi$ induce	morphisms $\phi \colon \module{E} \to \module{F}$ and
	$\psi \colon M \to \module{F}$ of constraint index sets we know
	by \autoref{prop:ProjectiveIndexSets} that there exists
	$\xi \colon M \to \module{E}$ such that $\phi \circ \xi = \psi$.
	Then by the freeness of $\algebra{A}^{(M)}$ there exists
	$\Xi \colon \algebra{A}^{(M)} \to \module{E}$ such that
	$\Phi \circ \Xi$ restricted to $M$ is just $\psi$.
	 Hence	$\Phi \circ \Xi = \Psi$.
\end{proof}

In the following we will concentrate on the case $\algebra{A} \in \injConAlg$.
With this we can show that the category $\ConRMod{\algebra{A}}$ has 
enough projectives in the following sense:

\begin{proposition}
	\label{prop:ModulesAreFreeImages}%
Let $\algebra{A} \in \injConAlg$ be a constraint algebra.
For every constraint module
$\module{E} \in \ConRMod{\algebra{A}}$ there exists
$M \in \injConIndSet$ and a regular epimorphism
$\Phi \colon \algebra{A}^{(M)} \to \module{E}$.
\end{proposition}

\begin{proof}
	Consider the constraint set $M$ given by
	$M_\Wobs = \module{E}_\Wobs$,
	$M_\Total = \module{E}_\Wobs \times \module{E}_\Total$
	and $\iota_M = \id_{\module{E}_\Wobs} \times \iota_\module{E}$,
	which is injective.
	Then $\phi = (\pr_2,\id_{\module{E}_\Wobs}) \colon M \to \module{E}$
	is a regular epimorphism.
	By	the universal property of $\algebra{A}^{(M)}$ there exists
	$\Phi \colon \algebra{A}^{(M)} \longrightarrow \module{E}$ such
	that $\Phi \circ i = \phi$. 
	Then $\Phi$ is a regular epimorphism since so is $\phi$.
\end{proof}

We can now use \autoref{prop:ModulesAreFreeImages} to show that projective
modules over $\algebra{A} \in \injConAlg$ are always embedded.

\begin{lemma}
	Let $\algebra{A} \in \injConAlg$ be an embedded constraint algebra and let
	$\module{P} \in \ConRMod{\algebra{A}}$ be projective.
	Then $\iota_\module{P} \colon \module{P}_\Wobs \to \module{P}_\Total$
	is injective, i.e. $\module{P} \in \injConRMod{\algebra{A}}$.
\end{lemma}

\begin{proof}
	By \autoref{prop:ModulesAreFreeImages} there exists $\module{E} \in \injConRMod{\algebra{A}}$ and a regular epimorphism
	$\Phi \colon \module{E} \to \module{P}$.
	Since $\module{P}$ is projective there exists 
	$\chi \colon \module{P} \to \module{E}$
	such that $\Phi \circ \chi = \id_\module{P}$.
	In particular, $\chi_\Wobs$ is injective and thus from
	$\chi_\Total \circ \iota_\module{P} = \iota_\module{E} \circ \chi_\Wobs$
	it follows that $\iota_\module{P}$ is injective.
\end{proof}

Another important notion in the characterization of projective
constraint modules is that of a split exact sequence.
A sequence of morphisms of constraint modules
\begin{equation}
	\begin{tikzcd}
		0
		\arrow{r}{}
		& \module{E}
		\arrow{r}{\Phi}
		& \module{F}
		\arrow{r}{\Psi}
		& \module{G}
		\arrow{r}{}
		& 0
	\end{tikzcd}
\end{equation}
is a called \emph{short exact} if $\Phi$ is a monomorphism,
$\image(\Phi) = \ker(\Psi)$, and $\Psi$ is a regular epimorphism.
It is called \emph{split exact} if in addition there exists
$\chi \colon \module{G} \to \module{F}$
such that $\Psi \circ \chi = \id_\module{G}$.

\begin{remark}
It can be shown that $\injConMod_\algebra{A}$
is a homological category in the sense of
\cite[Lemma 4.1.6]{borceux.bourn:2004a}.
The above definition of short exact sequences is in line
with the definition of short exact sequences in general
homological categories.
\end{remark}

\begin{proposition}
	\label{prop:RegProjGiveSplit}%
Let $\algebra{A} \in \ConAlg$ be a constraint algebra and
let $\module{P} \in \ConRMod{\algebra{A}}$ be a projective module.
Then every short exact sequence of the form
\begin{equation}
\begin{tikzcd}
	0
		\arrow[r]
	&\module{E}
		\arrow[r,"\Phi",hookrightarrow]
	&\module{F}
		\arrow[r,"\Psi",twoheadrightarrow]
	&\module{P}
		\arrow[r]
	&0
\end{tikzcd}
\end{equation}
is split exact.
\end{proposition}

\begin{proof}
Since $\module{P}$ is projective and $\Psi$ is a regular
epimorphism the sequence splits by the universal property of 
$\module{P}$.
\end{proof}

Despite $\ConRMod{\algebra{A}}$ not being an abelian category, the
splitting lemma nevertheless holds for constraint modules.

\begin{proposition}[Splitting lemma in $\ConRMod{\algebra{A}}$]
	\index{splitting lemma}
Let $\algebra{A} \in \ConAlg$ be a constraint algebra.
A short exact sequence
\begin{equation}
	\label{prop:SplittingLemma_shortexseq}
\begin{tikzcd}
	0
		\arrow[r]
	&\module{E}
		\arrow[r,"\Phi",hookrightarrow]
	&\module{F}
		\arrow[r,"\Psi",twoheadrightarrow]
	&\module{G}
		\arrow[r]
	&0
\end{tikzcd}
\end{equation}
in $\ConRMod{\algebra{A}}$ splits if and only if it is
isomorphic as a sequence to
\begin{equation}
\begin{tikzcd}
	0
		\arrow[r]
	&\module{E}
	\arrow[r,"i_\module{E}",hookrightarrow]
	&\module{E} \oplus \module{G}
	\arrow[r,"\pr_\module{G}",twoheadrightarrow]
	&\module{G}
	\arrow[r]
	&0
\end{tikzcd}
\end{equation}
with the canonical inclusion $i_\module{E}$ and projection
$\pr_\module{G}$.
\end{proposition}

\begin{proof}
Suppose there exists
$\chi \colon \module{F} \longrightarrow \module{E}$ such that
$\Psi \circ \chi = \id_\module{G}$.
Then we know that
$\module{F}_\Total \simeq \module{E}_\Total \oplus \module{G}_\Total$ 
and
$\module{F}_\Wobs \simeq \module{E}_\Wobs \oplus \module{G}_\Wobs$
by the splitting lemma in the respective categories of modules.
We denote these isomorphisms by $\theta_\Total$ and
$\theta_\Wobs$, respectively.
To show that these form a constraint morphism consider that
$\theta = (\Phi \circ \pr_1) + (\chi \circ \pr_2)$
is a composition of constraint morphisms, thus so is $\theta$ itself.
Moreover,
for every $y \in \module{F}_\Null$ we have
$y = (y - (\chi \circ \Psi)(y)) + (\chi \circ \Psi)(y) \in \module{E}_\Null \oplus \module{G}_\Null$,
hence $\theta$ is an isomorphism of constraint modules.
Conversely, suppose
$\theta \colon \module{E} \oplus \module{G} \to  \module{F}$
is an isomorphism such that
$\theta \circ i_\module{E} = \Phi$ and
$\Psi \circ \theta =  \pr_\module{G}$.
Then $\theta \circ i_\module{G}$ is clearly a splitting for
\eqref{prop:SplittingLemma_shortexseq}.
\end{proof}

The following result shows that projective modules
can be described as direct summands of $\injConIndSet$-free modules.
The proof is completely analogous to the usual case, see e.g.
\cite[Prop. 3.10]{jacobson:1989a}.

\begin{theorem}[Projective modules]
	\label{thm:Projective}%
	\index{projective!constraint module}
Let $\algebra{A} \in \injConAlg$ be a constraint algebra and $\module{P} 
\in \ConRMod{\algebra{A}}$ be given.
The following statements are equivalent:
\begin{theoremlist}
\item \label{thm:Projective_1}
	The module $\module{P}$ is projective.
\item \label{thm:Projective_2}
	Every short exact sequence
	$0 \rightarrow \module{E} \rightarrow \module{F} \rightarrow
	\module{P} \rightarrow 0$ splits.
\item \label{thm:Projective_3}
	The module $\module{P}$ is a
	direct summand of a $\injConIndSet$-free module, i.e.  there
	exists $M \in \injConIndSet$ and
	$\module{E} \in \ConRMod{\algebra{A}}$ such that
	$\algebra{A}^{(M)} \simeq \module{P} \oplus \module{E}$.
\item \label{thm:Projective_4}
	There exist
	$M \in \injConIndSet$ and
	$e = (e_\Total, e_\Wobs) \in
	\ConEnd_\algebra{A}(\algebra{A}^{(M)})$ such that $e^2 = e$ and
	$\module{P} \simeq e\algebra{A}^{(M)} = \image(e)$.
\end{theoremlist}
\end{theorem}

\begin{proof}
\ref{thm:Projective_1}\,$\Rightarrow$\,\ref{thm:Projective_2}: This is exactly
\autoref{prop:RegProjGiveSplit}.

\ref{thm:Projective_2}\,$\Rightarrow$\,\ref{thm:Projective_3}:
By \autoref{prop:ModulesAreFreeImages} there exists a short exact
sequence
$0 \rightarrow \module{E} \rightarrow \algebra{A}^{(M)} \rightarrow
\module{P} \rightarrow 0$
with $M \in \injConIndSet$.
This sequence splits by assumption, and therefore by the splitting
lemma we have
$\algebra{A}^{(M)} \simeq \module{E} \oplus \module{P}$.

\ref{thm:Projective_3}\,$\Rightarrow$\,\ref{thm:Projective_1}:
We have a split exact sequence
$0 \rightarrow \module{E} \rightarrow \algebra{A}^{(M)}
\rightarrow \module{P} \rightarrow 0$
with $M \in \injConIndSet$.
Let $\Psi \colon \module{P} \longrightarrow \module{F}$
and
$\Phi \colon \module{G} \longrightarrow \module{F}$
be given with $\Phi$ a regular epimorphism.
We get the following diagram:
\begin{equation*}
\begin{tikzcd}
	0
	\arrow[r]
	&  \module{E}
	\arrow[r,"\iota"]
	& \algebra{A}^{(M)}
	\arrow[r,"\pi"]
	& \module{P}
	\arrow[r]
	\arrow[d,"\Psi"]
	\arrow[l,"\sigma",bend left = 10, shift left = 2pt]
	& 0 \\
	{ }
	& { }
	& \module{G}
	\arrow[r,"\Phi",twoheadrightarrow]
	& \module{F}
	& { }
\end{tikzcd}
\end{equation*}
Since $\algebra{A}^{(M)}$ is projective there exists a
morphism
$\eta \colon \algebra{A}^{(M)} \to \module{G}$
such that $\Phi \circ \eta = \Psi \circ \pi$.
Then
$\eta \circ \sigma \colon \module{P} \to \module{G}$
yields the desired morphism making $\module{P}$ projective.

\ref{thm:Projective_3}\,$\Leftrightarrow$\,\ref{thm:Projective_4}:
If $\algebra{A}^{(M)} \simeq \module{P} \oplus \module{E}$,
then choose for $e \in \End_\algebra{A}(\algebra{A}^{(M)})$
the projection on $\module{P}$.
If $\module{P} \simeq e\algebra{A}^{(M)}$, then
$\module{E} \coloneqq \ker(e)$ gives the correct direct summand.
\end{proof}

\begin{definition}[Finitely generated projective modules]
	Let $\algebra{A} \in \injConAlg$ and
	let $\module{P} \in \injConMod_\algebra{A}$ be a projective constraint module.
	\begin{definitionlist}
		\item An embedded constraint index set $M$ such that
		$\algebra{A}^{(M)} \simeq \module{P} \oplus \module{E}$
		is called \emph{generating set} of the projective module 
		$\module{P}$.
		\item If $M$ can be chosen finite we call $\module{P}$ \emph{finitely 
			generated projective}.
		\item The category of finitely generated projective constraint modules
		over $\algebra{A}$ is denoted by
		\glsadd{ConProjA}$\ConProj(\algebra{A})$.
	\end{definitionlist}
\end{definition}

\begin{remark}
With the help of \autoref{thm:Projective} \ref{thm:Projective_3}
it is easy to see that direct sums of projective constraint modules are again projective.
This directly opens the possibility to define constraint $K_0$-theory for constraint algebras.
\end{remark}

In addition to the above characterizations of projective
modules we can also use a constraint version of a dual basis.

\begin{proposition}[Dual basis]
	\label{prop:DualBasis}
	\index{dual basis}
Let $\algebra{A} \in \injConAlg$ be a constraint algebra and $\module{P} 
\in \injConRMod{\algebra{A}}$.
The following statements are equivalent:
\begin{propositionlist}
	\item $\module{P}$ is projective with generating set
	$M \in \injConIndSet$.
	\item There exist families
	$(e_n)_{n \in M_\Total} \subseteq \module{P}_\Total$ and
	$(e^n)_{n \in M_\Total} \subseteq (\module{P}_\Total)^* =
	\Hom_{\algebra{A}_\Total}(\module{P}_\Total , \algebra{A}_\Total)$
	such that
	\begin{equation}
		\label{prop:DualBasis_Eq1}
		x = \sum_{n \in M_\Total} e_n e^n(x)
		\quad
	\end{equation}
	for all
	$x \in \module{P}_\Total$
	where for fixed $x$ only finitely many of the
	$e^n(x)$ differ from $0$.
	Moreover, the following properties need to be satisfied:
	\begin{propertieslist}
		\item \label{prop:DualBasis_1}
		One has
		$e_{n} \in \module{P}_\Wobs$ for
		$n \in M_\Wobs$.
		\item \label{prop:DualBasis_2}
		One has $e_n \in \module{P}_\Null$
		for $n \in M_\Null$.
		\item \label{prop:DualBasis_3}
		One has $e^n(\module{P}_\Wobs) \subseteq \algebra{A}_\Wobs$
		for $n \in M_\Total$.	
		\item \label{prop:DualBasis_4}
		One has
		$e^{n} \in (\module{P}^*)_\Wobs$ for $n \in M_\Total 
		\setminus M_\Null = (M^*)_\Wobs$.	
		\item \label{prop:DualBasis_5}
		One has
		$e^n(\module{P}_\Wobs) = 0$ for
		$n \in M_\Total \setminus M_\Wobs = (M^*)_\Null$.
	\end{propertieslist}
\end{propositionlist}
\end{proposition}

\begin{proof}
	Let $\module{P} \simeq e\algebra{A}^{(M)}$ be projective
	with idempotent $e \in \ConEnd_\algebra{A}(\algebra{A}^{(M)})$ and
	generating set $M \in \injConIndSet$. 
	Denote by
	$b_n \in \algebra{A}_\Total^{(M_\Total)}$ the standard basis and by
	$b^n$ the canonical coordinate functionals. 
	Defining
	$e_n = e_\Total(b_n)$ for
	$n \in M_\Total$ as well as
	$e^n = b^n\at{e\algebra{A}^{(M)}}$ gives a usual dual basis for
	$\algebra{A}_\Total^{(M_\Total)}$.
	 Thus we get \eqref{prop:DualBasis_Eq1}. 
	Since $e$ is a constraint morphism and it holds
	$b_n \in (\algebra{A}^{(M)})_\Wobs$ for $n \in M_\Wobs$
	and $b_n \in (\algebra{A}^{M})_\Null$ for $n \in M_\Null$ we get \ref{prop:DualBasis_1} and
	\ref{prop:DualBasis_2}.
	 For
	$x \in \algebra{A}_\Wobs^{(M_\Wobs)}$ it holds that
	$b^n(x) \in \algebra{A}_\Wobs$ for all $n \in M_\Total$
	and $b^n(x) = 0$ for all $n \in M_\Total \setminus M_\Wobs$.
	Moreover, if $x \in \algebra{A}_\Null^{(M_\Wobs \setminus M_\Null)} \oplus \algebra{A}_\Total^{(M_\Null)}$
	we get $b^n(x) \in \algebra{A}_\Null$ for all $n \in M_\Total \setminus M_\Null$.
	Hence \ref{prop:DualBasis_3}, \ref{prop:DualBasis_4} and \ref{prop:DualBasis_5} follow.
	Let now such a dual basis in the above sense be given.
	 The map $M \to \module{P}$ defined by
	 $n \mapsto e_n$ is a morphism of constraint index sets because of
	\ref{prop:DualBasis_1} and \ref{prop:DualBasis_2}.
	 By the universal property of free constraint modules we thus get an
	induced morphism $q \colon \algebra{A}^{(M)} \to \module{P}$.
	We define $\ins \colon \module{P} \to \algebra{A}^{(M)}$ by
	\begin{equation*}
		\ins(x) \coloneqq \sum_{n \in M_\Total} b_n e^n(x).
	\end{equation*}
	The map $\ins$ is clearly a module morphism as the $e^n$ are, and it is a constraint 
	morphism by \ref{prop:DualBasis_3}, \ref{prop:DualBasis_4} and \ref{prop:DualBasis_5}.
	We now show
	$q \circ \ins = \id_\module{P}$: For
	$x \in \module{P}_\Total$ we have
	\begin{equation*}
		q (\ins(x))
		=
		q\Bigl(\sum_{n \in M_\Total} b_n e^n(x)\Bigr)
		=
		\sum_{n \in M_\Total} e_n e^n(x)
		=
		x
	\end{equation*}
	by assumption.
	Thus the constraint endomorphism
	$e \coloneqq \ins \circ \mathop{q} \in \ConEnd_\algebra{A}(\algebra{A}^{(M)})$
	is an idempotent and $\module{P} \simeq e\algebra{A}^{(M)}$ via
	the maps $\ins$ and $q\at{e\algebra{A}^{(M)}}$.
	Hence $\module{P}$ is projective.
\end{proof}

In \autoref{prop:DualsOfConFreeModules} we have seen that duals of free constraint
modules need not be free in general.
One might hope that duals of free modules are at least projective.
But even this fails as the next example shows.
In particular we also see that duals of projective constraint modules need not be projective.

\begin{example}
	Consider a constraint algebra $\algebra{A} \in \injConAlg$ with
	$\algebra{A}_\Wobs \neq \algebra{A}_\Total$
	and finite $n \in \injConIndSet$ with $n_\Total \neq n_\Wobs$.
	Then $\algebra{A}^{n}$ is projective by \autoref{lem:FreeConModulesAreProjective}.
	We know
	\begin{equation}
		\big((\algebra{A}^n)^*\big)_\Wobs
		= \algebra{A}_\Null^{n_\Null} \oplus \algebra{A}_\Wobs^{n_\Wobs - n_\Null} \oplus \algebra{A}_\Total^{n_\Total - n_\Wobs}
	\end{equation}
	from \autoref{prop:DualsOfConFreeModules},
	which can never be a direct summand of some $\algebra{A}_\Wobs^m$.
	Thus it follows from \autoref{thm:Projective} that $\algebra{A}^n$ cannot be projective.
\end{example}

\subsubsection{Reduction}

The notion of  projectivity is compatible with the reduction functor of constraint modules.

\begin{proposition}[Reduction of projective constraint modules]
	\label{prop:ProjectiveReduction}%
	\index{reduction!projective modules}
Let $\algebra{A} \in \injConAlg$ be a constraint algebra and let $\module{P} \in \ConRMod{\algebra{A}}$ be
projective.
Then $\module{P}_\red$ is projective, and if
$\module{P} \simeq e\algebra{A}^{(M)}$ for some $M \in \injConIndSet$, then
$\module{P}_\red \simeq e_\red\algebra{A}_\red^{(M_\red)}$.
\end{proposition}

\begin{proof}
Suppose $\module{P} \simeq \module{P} \oplus \module{E}$,
with $e \in \ConEnd_\algebra{A}(\algebra{A}^{(M)})$
the projection onto $\module{P}$.
Then $(\algebra{A}^{(M)})_\red \simeq \module{P}_\red \oplus 
\module{E}_\red$
with $e_\red \in \End_{\algebra{A}_\red}((\algebra{A}^{(M)})_\red)$
the corresponding projection.
Since \autoref{prop:ReductionFreeStrConMod}
yields $(\algebra{A}^{(M)})_\red \simeq \algebra{A}_\red^{(M_\red)}$
the claim holds.
\end{proof}

\subsection{Projective Strong Constraint Modules}
	\label{sec:ProjStrConMod}

For strong constraint modules over a strong constraint algebra $\algebra{A}$ the situation
is quite similar to that of non-strong modules.
Therefore, in this section we will omit proofs that can be carried over from \autoref{sec:ProjectiveConModules}
word by word.

\begin{definition}[Projective strong module]
	\label{def:ProjectiveStrConModule}
	\index{projective!strong constraint modules}
Let $\algebra{A} \in \strConAlg$ be a strong constraint algebra.
A strong constraint $\algebra{A}$-module
$\module{P} \in \strConRMod{\algebra{A}}$ is called
\emph{projective} if for every
$\module{E}, \module{F} \in \strConRMod{\algebra{A}}$,  morphism
$\Psi \colon \module{P} \longrightarrow \module{F}$ and 
regular epimorphism
$\Phi \colon \module{E} \longrightarrow \module{F}$ there exists a
morphism $\chi \colon \module{P} \longrightarrow \module{E}$ such
that $\Phi \circ \chi = \Psi$.  Diagrammatically:
\begin{equation}
	\begin{tikzcd}
		{ }
		& \module{E}
		\arrow[twoheadrightarrow]{d}{\Phi} \\
		\module{P}
		\arrow{r}{\Psi}
		\arrow[dashed]{ur}{\chi}
		& \module{F}
	\end{tikzcd}
\end{equation}
\end{definition}

The category of strong constraint modules over an embedded strong 
constraint algebra
$\algebra{A} \in \injstrConAlg$ has enough projectives as the next proposition shows.

\begin{proposition}
	\label{lem:FreeStrConModulesAreProjective}
	Let $\algebra{A} \in \injstrConAlg$ be an embedded strong 
	constraint algebra.
	\begin{propositionlist}
		\item \label{lem:FreeStrConModulesAreProjective_1}
		For every $M \in \injConIndSet$ the free constraint module $\algebra{A}^{(M)}$ is projective.
		\item \label{lem:FreeStrConModulesAreProjective_2}
		For every strong constraint module
		$\module{E} \in \strConRMod{\algebra{A}}$ there exists
		$M \in \injConIndSet$ and a regular epimorphism
		$\Phi \colon \algebra{A}^{(M)} \to \module{E}$.
		\item \label{lem:FreeStrConModulesAreProjective_3}
		If $\module{P} \in \strConRMod{\algebra{A}}$ is projective, then
		$\module{P} \in \injstrConRMod{\algebra{A}}$.
	\end{propositionlist}
\end{proposition}

Note that here the free module $\algebra{A}^{(M)}$ is the free strong constraint module
in the sense of \autoref{lem:FreeStrConModules}.
With this the usual characterization of projective modules in terms of summands of free modules and projections also holds in the case of strong constraint modules.

\begin{theorem}[Projective strong modules]
	\label{thm:StrProjective}%
	Let $\algebra{A} \in \injstrConAlg$ be an embedded strong 
	constraint algebra and $\module{P} 
	\in \strConRMod{\algebra{A}}$ be given.
	The following statements are equivalent:
	\begin{theoremlist}
		\item \label{thm:StrProjective_1}
		The module $\module{P}$ is projective.
		\item \label{thm:StrProjective_2}
		Every short exact sequence
		$0 \rightarrow \module{E} \rightarrow \module{F} \rightarrow
		\module{P} \rightarrow 0$ splits.
		\item \label{thm:StrProjective_3}
		The module $\module{P}$ is a
		direct summand of a $\injConIndSet$-free module, i.e.  there
		exists $M \in \injConIndSet$ and
		$\module{E} \in \strConRMod{\algebra{A}}$ such that
		$\algebra{A}^{(M)} \simeq \module{P} \oplus \module{E}$.
		\item \label{thm:StrProjective_4}
		There exist
		$M \in \injConIndSet$ and
		$e = (e_\Total, e_\Wobs) \in
		\ConEnd_\algebra{A}(\algebra{A}^{(M)})$ such that $e^2 = e$ and
		$\module{P} \simeq e\algebra{A}^{(M)} = \image(e)$.
	\end{theoremlist}
\end{theorem}

\begin{definition}[Finitely generated projective modules]
	Let $\algebra{A} \in \injstrConAlg$ and
	let $\module{P} \in \injstrConMod_\algebra{A}$ be a projective strong constraint module.
	\begin{definitionlist}
		\item An embedded constraint index set $M$ such that
		$\algebra{A}^{(M)} \simeq \module{P} \oplus \module{E}$
		is called \emph{generating set} of the projective module 
		$\module{P}$.
		\item If $M$ can be chosen finite we call $\module{P}$ \emph{finitely 
			generated projective}.
		\item The category of finitely generated projective strong constraint modules
		over $\algebra{A}$ is denoted by
		\glsadd{strConProjA}$\strConProj(\algebra{A})$.
	\end{definitionlist}
\end{definition}

\begin{remark}\
\begin{remarklist}
	\item Projective constraint modules have been introduced in \cite{menke:2020a,dippell.menke.waldmann:2022a},
	while the notion of projective strong constraint modules appears here for the first time.
	\item Similar to the situation of projective constraint modules, the direct sum of projective strong constraint modules is again projective.
	This allows for the introduction of $K_0$-theory of strong constraint algebras, which will in general differ from the $K_0$-theory of constraint algebras.
\end{remarklist}
\end{remark}

There exists again a characterization in terms of a dual basis, but it differs slightly from the
dual basis for non-strong projective modules, cf. \autoref{prop:DualBasis}.

\begin{proposition}[Dual basis]
	\label{prop:StrDualBasis}
	\index{dual basis}
Let $\algebra{A} \in \injstrConAlg$ be an embedded strong 
constraint algebra and $\module{P} 
\in \injstrConRMod{\algebra{A}}$.
Then the following statements are equivalent:
\begin{propositionlist}
	\item $\module{P}$ is projective with generating set
	$M \in \injConIndSet$.
	\item There exist families
	$(e_n)_{n \in M_\Total} \subseteq \module{P}_\Total$ and
	$(e^n)_{n \in M_\Total} \subseteq (\module{P}_\Total)^* =
	\Hom_{\algebra{A}_\Total}(\module{P}_\Total , \algebra{A}_\Total)$
	such that
	\begin{equation}
		\label{prop:StrDualBasis_Eq1}
		x = \sum_{n \in M_\Total} e_n e^n(x)
		\quad
	\end{equation}
	for all
	$x \in \module{P}_\Total$
	where for fixed $x$ only finitely many of the
	$e^n(x)$ differ from $0$.
	Moreover, the following properties need to be satisfied:
	\begin{propertieslist}
		\item \label{prop:StrDualBasis_1}
		One has
		$e_{n} \in \module{P}_\Wobs$ for
		$n \in M_\Wobs$.
		\item \label{prop:StrDualBasis_2}
		One has $e_n \in \module{P}_\Null$
		for $n \in M_\Null$.
		\item \label{prop:StrDualBasis_3}
		One has
		$e^{n} \in (\module{P}^*)_\Wobs$ for $n \in M_\Total 
		\setminus M_\Null = (M^*)_\Wobs$.	
		\item \label{prop:StrDualBasis_4}
		One has
		$e^n \in (\module{P}^*)_\Null$ for
		$n \in M_\Total \setminus M_\Wobs = (M^*)_\Null$.
	\end{propertieslist}
\end{propositionlist} 
\end{proposition}

\begin{proof}
	Let $\module{P} \simeq e\algebra{A}^{(M)}$ be projective
	with idempotent $e \in \ConEnd_\algebra{A}(\algebra{A}^{(M)})$ and
	generating set $M \in \injConIndSet$. 
	Denote by
	$b_n \in \algebra{A}_\Total^{(M_\Total)}$ the standard basis and by
	$b^n$ the canonical coordinate functionals. 
	Defining
	$e_n = e_\Total(b_n)$ for
	$n \in M_\Total$ as well as
	$e^n = b^n\at{e\algebra{A}^{(M)}}$ gives a usual dual basis for
	$\algebra{A}_\Total^{(M_\Total)}$.
	Thus we get \eqref{prop:StrDualBasis_Eq1}. 
	Since $e$ is a constraint morphism and it holds
	$b_n \in (\algebra{A}^{(M)})_\Wobs$ for $n \in M_\Wobs$
	and $b_n \in (\algebra{A}^{M})_\Null$ for $n \in M_\Null$ we get \ref{prop:StrDualBasis_1} and
	\ref{prop:StrDualBasis_2}.
	For
	$x \in \algebra{A}_\Null^{(M_\Total\setminus M_\Wobs)} \oplus \algebra{A}_\Wobs^{(M_\Wobs\setminus M_\Null)} \oplus
	\algebra{A}_\Total^{(M_\Null)}$ it holds that
	$b^n(x) \in \algebra{A}_\Wobs$ for all $n \in M_\Total\setminus M_\Null$
	and $b^n(x) \in \algebra{A}_\Null$ for all $n \in M_\Total \setminus M_\Wobs$.
	Moreover, if $x \in \algebra{A}_\Null^{(M_\Wobs \setminus M_\Null)} \oplus \algebra{A}_\Total^{(M_\Null)}$
	we get $b^n(x) \in \algebra{A}_\Null$ for all $n \in M_\Total \setminus M_\Null$.
	Hence \ref{prop:StrDualBasis_3} and \ref{prop:StrDualBasis_4}.
	Let now such a dual basis in the above sense be given.
	The map $M \to \module{P}$ defined by
	$n \mapsto e_n$ is a morphism of constraint index sets because of
	\ref{prop:StrDualBasis_1} and \ref{prop:StrDualBasis_2}.
	By the universal property of free constraint modules we thus get an
	induced morphism $q \colon \algebra{A}^{(M)} \to \module{P}$.
	We define $\ins \colon \module{P} \to \algebra{A}^{(M)}$ by
	\begin{equation*}
		\ins(x) \coloneqq \sum_{n \in M_\Total} b_n e^n(x).
	\end{equation*}
	The map $\ins$ is clearly a module morphism as the $e^n$ are, and it is a constraint 
	morphism by \ref{prop:StrDualBasis_3} and \ref{prop:StrDualBasis_4}.
	We now show
	$q \circ \ins = \id_\module{P}$: For
	$x \in \module{P}_\Total$ we have
	\begin{equation*}
		q (\ins(x))
		=
		q\Bigl(\sum_{n \in M_\Total} b_n e^n(x)\Bigr)
		=
		\sum_{n \in M_\Total} e_n e^n(x)
		=
		x
	\end{equation*}
	by assumption.
	Thus the constraint endomorphism
	$e \coloneqq \ins \circ q \in \ConEnd_\algebra{A}(\algebra{A}^{(M)})$
	is an idempotent and $\module{P} \simeq e\algebra{A}^{(M)}$ via
	the maps $\ins$ and $q\at{e\algebra{A}^{(M)}}$.
	Hence $\module{P}$ is projective.
\end{proof}

We can view such a constraint dual basis as a pair $(\{e_n\}_{n \in M}, \{e^n\}_{n \in M^*})$
of constraint subsets indexed by $M$ and $M^*$, respectively.
By a constraint indexed subset $\{x_i\}_{i \in I}$ of a constraint set $X$ indexed by a constraint index set $I$
we simply mean a constraint map
$I \to X$.

\begin{proposition}[Duals of projective modules]
	\index{dual module}
Let $\algebra{A} \in \injstrConAlg$ be a strong constraint algebra
and let $\module{P} \in \strConProj(\algebra{A})$ be finitely 
generated projective.
\begin{propositionlist}
	\item $\module{P}^*$ is finitely generated projective.
	\item If $(\{e_i\}_{i \in M},\{e^i\}_{i \in M^*} )$ is a constraint dual basis of $\module{P}$, then
	$(\{e^i\}_{i \in M^*},\{e_i\}_{i \in M})$ is a constraint dual basis for $\module{P}^*$.
\end{propositionlist}
\end{proposition}

\begin{proof}
	By \autoref{thm:StrProjective} 
	\ref{thm:StrProjective_3} we know that there exists a 
	finite $M \in \injConIndSet$ and a $\algebra{A}$-module
	$\module{E}$ such that $\algebra{A}^M \simeq \module{P} \oplus 
	\module{E}$.
	Then by \autoref{prop:DualsOfStrFree} we have
	$\algebra{A}^{M^*} \simeq (\algebra{A}^M)^*\simeq \module{P}^* 
	\oplus \module{E}^*$,
	and therefore $\module{P}^*$ is again finitely generated 
	projective.
	For the second part recall that we know from classical algebra that
	$(\{e^i\},\{e_i\})_{i \in M_\Total}$ is a dual basis for 
	$\module{P}_\Total^*$, by identifying $e_i$ with its
	insertion functional $\delta_{e_i}$.
	Then using $(M^*)^* = M$ we see that properties \ref{prop:StrDualBasis_1} and \ref{prop:StrDualBasis_2}
	of \autoref{prop:StrDualBasis} for the dual basis of $\module{P}$
	exactly give \ref{prop:StrDualBasis_3} and \ref{prop:StrDualBasis_4} of
	\autoref{prop:StrDualBasis} for $\module{P}^*$, and vice versa. 
\end{proof}

By \autoref{lem:FreeStrConModulesAreProjective} \ref{lem:FreeStrConModulesAreProjective_3}
we know that for a commutative $\algebra{A} \in \injstrConAlg$ the category $\strConProj(\algebra{A})$ is a full subcategory of
$\injstrConBimod(\algebra{A})_\sym$, using the identification of one-sided modules over a commutative algebra with symmetric bimodules.
The category $\injstrConBimod(\algebra{A})_\sym$
carries two distinct monoidal structures: $\strcirctensor$ and 
$\injstrtensor$, 
see \autoref{sec:StrConAlgebras}.
We want to understand if $\strConProj(\algebra{A})$
is closed under taking $\strcirctensor$ and $\injstrtensor$ products.

\begin{proposition}[Tensor product on $\strConProj(\algebra{A})$]
	\label{prop:TensorStrConProj}
Let $\algebra{A} \in \injstrConAlg$ be commutative and  $\module{E}, \module{F} \in \strConProj(\algebra{A})$.
\begin{propositionlist}
	\item \label{prop:TensorStrConProj_1}
	$\strConProj(\algebra{A})$ is a monoidal subcategory of
	$(\injstrConBimod(\algebra{A})_\sym, \strcirctensor[\algebra{A}])$.
	In particular, $\module{E} \strcirctensor[\algebra{A}] \module{F}$ is
	finitely generated projective.
	\item \label{prop:TensorStrConProj_2}
	If $(\{e_i\}_{i \in M}, \{e^i\}_{i \in M^*})$
	and
	$(\{f_j\}_{j \in N}, \{f^j\}_{j \in N^*})$
	are dual bases of $\module{E}$ and $\module{F}$, respectively, 
	then 
	$(\{e_i \tensor f_j\}_{(i,j) \in M \tensor N}, \{e^i \tensor f^j\}_{(i,j) \in (M \tensor N)^*})$
	is a dual basis for $\module{E} \strcirctensor[\algebra{A}] \module{F}$.
\end{propositionlist}
\end{proposition}

\begin{proof}
We first prove the second part.
From classical algebra we know that
\begin{equation*}
	\big(\{e_i \tensor f_j\}_{(i,j) \in M_\Total \times N_\Total}, \{e^i \tensor f^j\}_{(i,j) \in M_\Total \times N_\Total}\big)
\end{equation*}
is a dual basis for $\module{E}_\Total \tensor[\algebra{A}_\Total] \module{F}_\Total$.
We need to check properties \ref{prop:StrDualBasis_1} to \ref{prop:StrDualBasis_4}
from \autoref{prop:StrDualBasis}:
For this recall that with the notation of \autoref{not:ConIndSetProducts} we have
\begin{align*}
	(M \tensor N)_\Wobs &= M \ConGrid[2][2][0][2][2] N,
	&&&(M \tensor N)^*_\Wobs &= M \ConGrid[0][0][2][0][2][2][2][2][2] N,\\
	(M \tensor N)_\Null &= M \ConGrid[2][2][0][2] N,
	&&&(M \tensor N)^*_\Null &= M \ConGrid[0][0][2][0][0][2][2][2][2] N.
\end{align*}
and
\begin{equation*}
	(\module{E} \strcirctensor[\algebra{A}] \module{F})_\Wobs
	= \module{E} \ConGrid[2][2][0][2][2] \module{F} \quad\text{ and }\quad
	(\module{E} \strcirctensor[\algebra{A}] \module{F})_\Null
	= \module{E} \ConGrid[2][2][0][2] \module{F}.
\end{equation*}
With this we can go through all the different cases:
\begin{cptitem}
	\item $(i,j) \in M \ConGrid[2][2][0][2][2] N$: Then
	$e_i \tensor f_j \in \module{E} \ConGrid[2][2][0][2][2] \module{F}
	= (\module{E} \strcirctensor[\algebra{A}] \module{F})_\Wobs$ holds.
	\item $(i,j) \in M \ConGrid[2][2][0][2] N$: We clearly have
	$e_i \tensor f_j \in\module{E} \ConGrid[2][2][0][2] \module{F} 
	=  (\module{E} \strcirctensor[\algebra{A}] \module{F})_\Null$, since at least
	one of $e_i$ and $f_j$ lies in the $\NULL$-component.
	\item $(i,j) \in M \ConGrid[0][0][2][0][0][2][2][2][2] N$: Suppose
	$x \tensor y \in \module{E} \ConGrid[2][2][0][2][2] \module{F}$.
	Then 
	\begin{equation*}
		(e^i \tensor f^j)(x \tensor y )
		= e^i(x) \cdot f^j(y)
		\in \algebra{A}_\Total \cdot \algebra{A}_\Null
		+ \algebra{A}_\Null \cdot \algebra{A}_\Total
		= \algebra{A}_\Null,
	\end{equation*}
	and thus $e^i \tensor f^j \in (\module{E} \tensor[\algebra{A}] \module{F})^*_\Null$.
	\item $(i,j) \in M \ConGrid[0][0][0][0][2] N$:
	For this let first $x \tensor y \in \module{E} \ConGrid[2][2][0][2][2] \module{F}$ be given, then 
	\begin{equation*}
		(e^i \tensor f^j)(x \tensor y) = e^i(x) \cdot f^j(y)
		\in \algebra{A}_\Wobs \cdot \algebra{A}_\Wobs
		= \algebra{A}_\Wobs.
	\end{equation*}
	Moreover, for $x \tensor y \in \module{E} \ConGrid[2][2][0][2] \module{F}$
	we have
	\begin{equation*}
		(e^i \tensor f^j)(x \tensor y) = e^i(x) \cdot f^j(y)
		\in \algebra{A}_\Null \cdot \algebra{A}_\Wobs
		+ \algebra{A}_\Wobs \cdot \algebra{A}_\Null
		= \algebra{A}_\Null.
	\end{equation*}
	Thus we get $e^i \tensor f^j \in (\module{E} \strcirctensor[\algebra{A}] \module{F})^*_\Wobs$.
\end{cptitem}	
This shows \ref{prop:TensorStrConProj_2}.
Hence we have $\module{E} \strcirctensor[\algebra{A}] \module{F} \in \strConProj(\algebra{A})$, and since also $\algebra{A} \in \strConProj(\algebra{A})$ holds
by \autoref{lem:FreeStrConModulesAreProjective} \ref{lem:FreeStrConModulesAreProjective_1}, we see that $\strConProj(\algebra{A})$ is a monoidal subcategory of
$\injstrConBimod(\algebra{A})_\sym$.
\end{proof}

\begin{proposition}[Strong Tensor product on $\strConProj(\algebra{A})$]
	\label{prop:StrTensorStrConProj}
Let $\algebra{A} \in \injstrConAlg$ be commutative and  $\module{E}, \module{F} \in \strConProj(\algebra{A})$.
\begin{propositionlist}
	\item \label{prop:StrTensorStrConProj_1}
	$\strConProj(\algebra{A})$ is a monoidal subcategory of
	$(\injstrConBimod(\algebra{A})_\sym, \injstrtensor[\algebra{A}])$.
	In particular, ${\module{E} \injstrtensor[\algebra{A}] \module{F}}$ is
	finitely generated projective.
	\item \label{prop:StrTensorStrConProj_2}
	If $(\{e_i\}_{i \in M}, \{e^i\}_{i \in M^*})$
	and
	$(\{f_j\}_{j \in N}, \{f^j\}_{j \in N^*})$
	are dual bases of $\module{E}$ and $\module{F}$, respectively, 
	then 
	$(\{e_i \tensor f_j\}_{(i,j) \in M \strtensor N}, \{e^i \tensor f^j\}_{(i,j) \in (M \strtensor N)^*})$
	is a dual basis for $\module{E} \injstrtensor[\algebra{A}] \module{F}$.
\end{propositionlist}
\end{proposition}

\begin{proof}
We first prove the second part.
From classical algebra we know that
\begin{equation*}
	(\{e_i \tensor f_j\}_{(i,j) \in M_\Total \times N_\Total}, \{e^i \tensor f^j\}_{(i,j) \in M_\Total \times N_\Total})
\end{equation*}
is a dual basis for $\module{E}_\Total \tensor[\algebra{A}_\Total] \module{F}_\Total$.
We need to check properties \ref{prop:StrDualBasis_1} to \ref{prop:StrDualBasis_4}
from \autoref{prop:StrDualBasis}:
For this recall that with the notation of \autoref{not:ConIndSetProducts} we have
\begin{align*}
	(M\strtensor N)_\Wobs &= M \ConGrid[2][2][2][2][2][0][2] N,
	&&&(M \strtensor N)^*_\Wobs &= M \ConGrid[0][0][0][0][2][2][0][2][2] N,\\
	(M\strtensor N)_\Null &= M \ConGrid[2][2][2][2][0][0][2] N, 
	&&&(M \strtensor N)^*_\Null &= M \ConGrid[0][0][0][0][0][2][0][2][2] N
\end{align*}
and
\begin{equation*}
	(\module{E} \injstrtensor[\algebra{A}] \module{F})_\Wobs
	= \module{E} \ConGrid[2][2][2][2][2][0][2] \module{F} \quad\text{ and }\quad
	(\module{E} \injstrtensor[\algebra{A}] \module{F})_\Null
	= \module{E} \ConGrid[2][2][2][2][0][0][2] \module{F}.
\end{equation*}
With this we can go through all the different cases:
\begin{cptitem}
	\item $(i,j) \in M \ConGrid[2][2][2][2][0][0][2] N = (M \strtensor N)_\Null$: Then
	$e_i \tensor f_j
	\in \module{E} \ConGrid[2][2][2][2][0][0][2] \module{F}
	= (\module{E} \injstrtensor[\algebra{A}] \module{F})_\Null$ holds,
	since at least
	one of $e_i$ and $f_j$ lies in the $\NULL$-component.
	\item $(i,j) \in M \ConGrid[0][0][0][0][2] N 
	\subseteq (M \strtensor N)_\Wobs$:
	We clearly have
	$e_i \tensor f_j \in \module{E} \ConGrid[0][0][0][0][2] \module{F}
	\subseteq (\module{E} \injstrtensor[\algebra{A}_\Total] \module{F})_\Wobs$.
	\item $(i,j) \in M \ConGrid[0][0][0][0][0][2][0][2][2] N 
	= (M \strtensor N)^*_\Null$:
	Suppose
	$x \tensor y \in \module{E} \ConGrid[2][2][2][2][0][0][2] \module{F} 
	= (\module{E} \injstrtensor[\algebra{A}] \module{F})_\Null$,
	then 
	\begin{equation*}
		(e^i \tensor f^j)(x \tensor y )
		= e^i(x) \cdot f^j(y)
		\in \algebra{A}_\Total \cdot \algebra{A}_\Null
		+ \algebra{A}_\Null \cdot \algebra{A}_\Total
		= \algebra{A}_\Null,
	\end{equation*}
	since both $e^i$ and $f^j$ map $\NULL$-components to $\NULL$-components.
	Moreover, for
	$x \tensor y \in \module{E} \ConGrid[0][0][0][0][2] \module{F}
	\subseteq (\module{E} \injstrtensor[\algebra{A}] \module{F})_\Wobs$
	we have
	\begin{equation*}
		(e^i \tensor f^j)(x \tensor y )
		= e^i(x) \cdot f^j(y)
		\in \algebra{A}_\Wobs \cdot \algebra{A}_\Null
		+ \algebra{A}_\Null \cdot \algebra{A}_\Wobs
		= \algebra{A}_\Null,
	\end{equation*}
	and thus
	$e^i \tensor f^j \in (\module{E} \injstrtensor[\algebra{A}_\Total] \module{F})^*_\Null$.
	\item $(i,j) \in M \ConGrid[0][0][0][0][2] N \subseteq (M \strtensor N)^*_\Wobs$:
	For this let first
	$x \tensor y \in \module{E} \ConGrid[2][2][2][2][0][0][2] \module{F} 
	= (\module{E} \injstrtensor[\algebra{A}] \module{F})_\Null$
	be given, then 
	\begin{equation*}
		(e^i \tensor f^j)(x \tensor y )
		= e^i(x) \cdot f^j(y)
		\in \algebra{A}_\Total \cdot \algebra{A}_\Null
		+ \algebra{A}_\Null \cdot \algebra{A}_\Total
		= \algebra{A}_\Null,
	\end{equation*}
	since both $e^i$ and $f^j$ map $\NULL$-components to $\NULL$-components.
	Moreover, for
	$x \tensor y \in \module{E} \ConGrid[0][0][0][0][2] \module{F}
	\subseteq (\module{E} \injstrtensor[\algebra{A}] \module{F})_\Wobs$
	we have
	\begin{equation*}
		(e^i \tensor f^j)(x \tensor y) = e^i(x) \cdot f^j(y)
		\in \algebra{A}_\Wobs \cdot \algebra{A}_\Wobs
		= \algebra{A}_\Wobs.
	\end{equation*}
	Thus we get $e^i \tensor f^j \in (\module{E} \injstrtensor[\algebra{A}] \module{F})^*_\Wobs$.
\end{cptitem}
This shows \ref{prop:TensorStrConProj_2}.
Hence we have $\module{E} \injstrtensor[\algebra{A}] \module{F} \in \strConProj(\algebra{A})$, and since also $\algebra{A} \in \strConProj(\algebra{A})$ holds
by \autoref{lem:FreeStrConModulesAreProjective} \ref{lem:FreeStrConModulesAreProjective_1}, we see that $\strConProj(\algebra{A})$ is a monoidal subcategory of
$\injstrConBimod(\algebra{A})$.
\end{proof}

We are now in a position to show that the canonical morphisms
from \autoref{prop:ConHomAndTensorDual} 
and \autoref{prop:DualOfTensor} are in fact isomorphisms, when restricting to
finitely generated projective modules.

\begin{proposition}
	\label{prop:DualTensorHomIsos}
Let $\algebra{A} \in \injstrConAlg$ and  $\module{E}, \module{F} \in 
\strConProj(\algebra{A})$.
\begin{propositionlist}
	\item \label{prop:DualTensorHomIsos_1}
	The canonical morphism
	$\module{F} \injstrtensor[\algebra{A}] \module{E}^* \to 
	\strConHom_\algebra{A}(\module{E},\module{F})$
	given by \eqref{eq:ConHomAndTensorDual}
	is an isomorphism.
	\item \label{prop:DualTensorHomIsos_2}
	The canonical morphism 
	$\module{E}^* \strcirctensor[\algebra{A}] \module{F}^* \to 
	(\module{E} \injstrtensor[\algebra{A}] \module{F})^*$
	given by \eqref{eq:DualOfStrTensor}
	is an isomorphism.
	\item \label{prop:DualTensorHomIsos_3}
	The canonical morphism 
	$\module{E}^* \injstrtensor[\algebra{A}] \module{F}^* \to 
	(\module{E} \strcirctensor[\algebra{A}] \module{F})^*$
	given by \eqref{eq:DualOfTensor}
	is an isomorphism.
\end{propositionlist}
\end{proposition}

\begin{proof}
	In all three cases we show that the well-known inverse maps on the $\TOTAL$-components
	are in fact constraint maps, and therefore yield constraint 
	inverses.
	To do this we need to fix dual bases
	$(\{e_i\}_{i \in M}, \{e^i\}_{i \in M^*})$
	and
	$(\{f_j\}_{j \in N}, \{f^j\}_{j \in N^*})$
	of $\module{E}$ and $\module{F}$, respectively.
	For the first part consider the map 
	\begin{equation*}
		\label{eq:DualTensorHomIsos_InverseHom}
		\Hom_{\algebra{A}_\Total}(\module{E}_\Total,\module{F}_\Total)
		\ni \Phi \mapsto
		\sum_{i \in M_\Total} \Phi(e_i) \tensor e^i
		\in \module{F}_\Total \tensor \module{E}^*_\Total.
		\tag{$*$}
	\end{equation*}
	This is the inverse to \eqref{eq:ConHomAndTensorDual}
	on the $\TOTAL$-component.
	Hence we need to show that \eqref{eq:DualTensorHomIsos_InverseHom}
	is a constraint morphism.
	For this let $\Phi \in 
	\strConHom_\algebra{A}(\module{E},\module{F})_\Wobs$
	be given.
	\begin{cptitem}
		\item If $i \in M_\Total\setminus M_\Wobs = M^*_\Null$, then
		$\Phi(e_i) \tensor e^i \in \module{F}_\Total \tensor 
		\module{E}^*_\Null \subseteq
		(\module{F} \injstrtensor[\algebra{A}] \module{E}^*)_\Null$.
		\item If $i \in M_\Wobs\setminus M_\Null \subseteq 
		M^*_\Wobs$, then, since in particular $i \in M_\Wobs$ holds,
		we obtain
		$\Phi(e_i) \tensor e^i \in \module{F}_\Wobs \tensor 
		\module{E}^*_\Wobs \subseteq
		(\module{F} \injstrtensor[\algebra{A}] \module{E}^*)_\Wobs$.
		\item If $i \in M_\Null$, then
		$\Phi(e_i) \tensor e^i \in \module{F}_\Null \tensor 
		\module{E}^*_\Total \subseteq
		(\module{F} \injstrtensor[\algebra{A}] \module{E}^*)_\Null$.
	\end{cptitem}
	Hence \eqref{eq:DualTensorHomIsos_InverseHom} preserves the 
	$\WOBS$-component.
	To show that it also preserves the $\NULL$-component, we only 
	need to reconsider the second case from above.
	\begin{cptitem}
		\item If $i \in M_\Wobs\setminus M_\Null \subseteq 
		M^*_\Wobs$, then, since in particular $i \in M_\Wobs$ holds,
		we obtain
		$\Phi(e_i) \tensor e^i \in \module{F}_\Wobs \tensor 
		\module{E}^*_\Null \subseteq
		(\module{F} \injstrtensor[\algebra{A}] \module{E}^*)_\Null$.
	\end{cptitem}
	This shows that \eqref{eq:DualTensorHomIsos_InverseHom}
	is a constraint inverse to \eqref{eq:ConHomAndTensorDual}.
	
	For part \ref{prop:DualTensorHomIsos_2} we need to show that 
	\begin{equation*}
		\label{eq:DualTensorHomIsos_InverseTensor}
		(\module{E}_\Total 
		\tensor[\algebra{A}_\Total] \module{F}_\Total)^*
		\ni
		\alpha
		\mapsto \sum_{(i,j) \in M_\Total \times N_\Total}
		\alpha(e_i \tensor f_j) \cdot e^i \tensor f^j
		\in \module{E}_\Total^* \tensor[\algebra{A}_\Total] 
		\module{F}_\Total^*
		\tag{$**$}
	\end{equation*}
	defines a constraint morphism 
	$(\module{E} \injstrtensor[\algebra{A}] \module{F})^*
	\to \module{E}^* \strcirctensor[\algebra{A}] \module{F}^*$.
	Recall that the families\linebreak
	$(\{ e_i \tensor f_j \}_{(i,j) \in M \strtensor N}, 
	\{e^i \tensor f^j\}_{(i,j) \in M^* \tensor N^*})$
	form a dual basis of
	$\module{E} \injstrtensor[\algebra{A}] \module{F}$,
	while the families
	$(\{ e^i \tensor f^j \}_{(i,j) \in M^* \tensor N^*}, 
	\{e_i \tensor f_j\}_{(i,j) \in M \strtensor N})$
	are a dual basis of 
	$\module{E}^* \strcirctensor[\algebra{A}] \module{F}^*$.
	Now suppose $\alpha \in (\module{E} 
	\injstrtensor[\algebra{A}] \module{F})^*_\Wobs$.
	\begin{cptitem}
		\item If $(i,j) \in M \ConGrid[0][0][0][0][0][2][0][2][2] N
		= (M^* \tensor N^*)_\Null$, then
		$\alpha(e_i \tensor f_j) \cdot e^i \tensor f^j
		\in \algebra{A}_\Total \cdot (\module{E}^* 
		\strcirctensor[\algebra{A}] \module{F}^*)_\Null
		= (\module{E}^* 
		\strcirctensor[\algebra{A}] \module{F}^*)_\Null$.
		\item If $(i,j) \in M \ConGrid[0][0][0][0][2] N
		\subseteq (M^* \tensor N^*)_\Wobs$, then,
		since also $M \ConGrid[0][0][0][0][2] N
		\subseteq (M \strtensor N)_\Wobs$ holds, 
		we get
		\begin{equation*}
			\alpha(\underbrace{e_i \tensor f_j}_{\in (\module{E} 
			\strtensor \module{F})_\Wobs}) \cdot e^i 
			\tensor f^j
			\in \algebra{A}_\Wobs \cdot
			(\module{E}^* \strcirctensor[\algebra{A}] 
			\module{F}^*)_\Wobs
			= (\module{E}^* \strcirctensor[\algebra{A}] 
			\module{F}^*)_\Wobs.
		\end{equation*}
		\item If $(i,j) \in M \ConGrid[2][2][2][2][0][0][2] N
		= (M \strtensor N)_\Null$, then
		\begin{equation*}
			\alpha(\underbrace{e_i \tensor f_j}_{(\module{E} 
			\strtensor 
				\module{F})_\Null}) \cdot e^i \tensor f^j
			\in \algebra{A}_\Null \cdot (\module{E}^* 
			\strcirctensor[\algebra{A}] \module{F}^*)_\Total
			\subseteq (\module{E}^* \strcirctensor[\algebra{A}] 
			\module{F}^*)_\Null.
		\end{equation*}
	\end{cptitem}
	Thus \eqref{eq:DualTensorHomIsos_InverseTensor}
	preserves the $\WOBS$-component.
	Next take $\alpha \in (\module{E} \injstrtensor[\algebra{A}] 
	\module{F})^*_\Null$.
	Since $(\module{E} \injstrtensor[\algebra{A}] \module{F})^*_\Null
	\subseteq (\module{E} \injstrtensor[\algebra{A}] 
	\module{F})^*_\Wobs$ we only need to check one of the above cases.
	\begin{cptitem}
		\item If $(i,j) \in M \ConGrid[0][0][0][0][2] N
		\subseteq (M^* \tensor N^*)_\Wobs$, then,
		since also since $M \ConGrid[0][0][0][0][2] N
		\subseteq (M \strtensor N)_\Wobs$ holds, 
		we get
		\begin{equation*}
			\alpha(\underbrace{e_i \tensor f_j}_{\in (\module{E} 
				\strtensor \module{F})_\Wobs}) \cdot e^i 
			\tensor f^j
			\in \algebra{A}_\Null \cdot
			(\module{E}^* \strcirctensor[\algebra{A}] 
			\module{F}^*)_\Wobs
			= (\module{E}^* \strcirctensor[\algebra{A}] 
			\module{F}^*)_\Null.
		\end{equation*}
	\end{cptitem}
	Hence \eqref{eq:DualTensorHomIsos_InverseTensor}
	also preserves the $\NULL$-component, showing that it is 
	a constraint inverse to \eqref{eq:DualOfStrTensor}.
	
	For part \ref{prop:DualTensorHomIsos_3} we proceed similarly.
	But this time we need to show that 
	\eqref{eq:DualTensorHomIsos_InverseTensor}
	defines a constraint morphism
	$(\module{E} \strcirctensor[\algebra{A}] \module{F})^*
	\to \module{E}^* \injstrtensor[\algebra{A}] \module{F}^*$.
	Recall that
	$(\{ e_i \tensor f_j \}_{(i,j) \in M \tensor N}, 
	\{e^i \tensor f^j\}_{(i,j) \in M^* \strtensor N^*})$
	is a dual basis of
	$\module{E} \strcirctensor[\algebra{A}] \module{F}$
	and
	$(\{ e^i \tensor f^j \}_{(i,j) \in M^* \strtensor N^*}, 
	\{e_i \tensor f_j\}_{(i,j) \in M \tensor N})$
	is a dual basis of 
	$\module{E}^* \injstrtensor[\algebra{A}] \module{F}^*$.
	Now suppose
	$\alpha \in (\module{E} \strcirctensor[\algebra{A}] 
	\module{F})^*_\Wobs$.
	\begin{cptitem}
		\item If $(i,j) \in M \ConGrid[0][0][2][0][0][2][2][2][2] N
		= (M^* \strtensor N^*)_\Null$, then
		$\alpha(e_i \tensor f_j) \cdot e^i \tensor f^j
		\in \algebra{A}_\Total \cdot (\module{E}^* 
		\injstrtensor[\algebra{A}] \module{F}^*)_\Null
		= (\module{E}^* 
		\injstrtensor[\algebra{A}] \module{F}^*)_\Null$.
		\item If $(i,j) \in M \ConGrid[0][0][0][0][2] N
		\subseteq (M^* \strtensor N^*)_\Wobs$, then,
		since also $M \ConGrid[0][0][0][0][2] N
		\subseteq (M \tensor N)_\Wobs$ holds, 
		we get
		\begin{equation*}
			\alpha(\underbrace{e_i \tensor f_j}_{\in (\module{E} 
				\tensor \module{F})_\Wobs}) \cdot e^i 
			\tensor f^j
			\in \algebra{A}_\Wobs \cdot
			(\module{E}^* \injstrtensor[\algebra{A}] 
			\module{F}^*)_\Wobs
			= (\module{E}^* \injstrtensor[\algebra{A}] 
			\module{F}^*)_\Wobs.
		\end{equation*}
		\item If $(i,j) \in M \ConGrid[2][2][0][2] N
		= (M \tensor N)_\Null$, then
		\begin{equation*}
			\alpha(\underbrace{e_i \tensor f_j}_{(\module{E} 
				\tensor 
				\module{F})_\Null}) \cdot e^i \tensor f^j
			\in \algebra{A}_\Null \cdot (\module{E}^* 
			\injstrtensor[\algebra{A}] \module{F}^*)_\Total
			\subseteq (\module{E}^* \injstrtensor[\algebra{A}] 
			\module{F}^*)_\Null.
		\end{equation*}
	\end{cptitem}
	Thus \eqref{eq:DualTensorHomIsos_InverseTensor}
	preserves the $\WOBS$-component.
	Next take $\alpha \in (\module{E} \strcirctensor[\algebra{A}] 
	\module{F})^*_\Null$.
	Since $(\module{E} \strcirctensor[\algebra{A}] \module{F})^*_\Null
	\subseteq (\module{E} \strcirctensor[\algebra{A}] 
	\module{F})^*_\Wobs$ we only need to check one of the above cases.
	\begin{cptitem}
		\item If $(i,j) \in M \ConGrid[0][0][0][0][2] N
		\subseteq (M^* \strtensor N^*)_\Wobs$, then,
		since also since $M \ConGrid[0][0][0][0][2] N
		\subseteq (M \tensor N)_\Wobs$ holds, 
		we get
		\begin{equation*}
			\alpha(\underbrace{e_i \tensor f_j}_{\in (\module{E} 
				\tensor \module{F})_\Wobs}) \cdot e^i 
			\tensor f^j
			\in \algebra{A}_\Null \cdot
			(\module{E}^* \injstrtensor[\algebra{A}] 
			\module{F}^*)_\Wobs
			= (\module{E}^* \injstrtensor[\algebra{A}] 
			\module{F}^*)_\Null.
		\end{equation*}
	\end{cptitem}
	Hence \eqref{eq:DualTensorHomIsos_InverseTensor}
	also preserves the $\NULL$-component, showing that it is 
	a constraint inverse to \eqref{eq:DualOfTensor}.
\end{proof}

\autoref{prop:DualTensorHomIsos} \ref{prop:DualTensorHomIsos_1} shows 
that 
$\strConHom(\module{E},\module{F})$
is again finitely generated projective.

\begin{corollary}
	\label{cor:StrHomOfTensor}
Let $\algebra{A} \in \injstrConAlg$ and
$\module{E}, \module{F}, \module{G}, \module{H}
\in \strConProj(\algebra{A})$.
\begin{corollarylist}
	\item There exists a canonical isomorphism
	\begin{equation}
		\strConHom(\module{E} \strcirctensor[\algebra{A}] \module{F}, 
		\module{G})
		\simeq \strConHom(\module{F}, \module{G} 
		\injstrtensor[\algebra{A}] \module{E}^*).
	\end{equation}
	\item There exists a canonical isomorphism
	\begin{equation}
		\strConHom(\module{E}, \module{F})
		\injstrtensor[\algebra{A}] \strConHom(\module{G},\module{H})
		\simeq \strConHom(\module{E} \strcirctensor[\algebra{A}] 
		\module{G}, 
		\module{F} \injstrtensor[\algebra{A}] \module{H}).
	\end{equation}
\end{corollarylist}
\end{corollary}

\begin{proof}
	By \autoref{prop:DualTensorHomIsos} we have canonical isomorphisms
	\begin{align*}
		\strConHom(\module{E} \strcirctensor[\algebra{A}] \module{F}, 
		\module{G})
		&\simeq \module{G} \injstrtensor[\algebra{A}] (\module{E}
		\strcirctensor[\algebra{A}] \module{F})^*\\
		&\simeq \module{G} \injstrtensor[\algebra{A}] \module{E}^*
		\injstrtensor[\algebra{A}] \module{F}^*\\
		&\simeq \strConHom(\module{F}, \module{G} 
		\injstrtensor[\algebra{A}] \module{E}^*),
	\end{align*}
	and
	\begin{align*}
		\strConHom(\module{E}, \module{F})
		\injstrtensor[\algebra{A}] \strConHom(\module{G},\module{H})
		&\simeq \module{E}^* \injstrtensor[\algebra{A}] \module{F}
		\injstrtensor[\algebra{A}] \module{G}^*
		\injstrtensor[\algebra{A}] \module{H} \\
		&\simeq (\module{E} \strcirctensor[\algebra{A}] \module{G})^*
		\injstrtensor[\algebra{A}] \module{F} 
		\injstrtensor[\algebra{A}] \module{H} \\
		&\simeq \strConHom(\module{E} \strcirctensor[\algebra{A}] 
		\module{G}, 
		\module{F} \injstrtensor[\algebra{A}] \module{H}).
	\end{align*}
\end{proof}

\begin{remark}
	With \autoref{cor:StrHomOfTensor} it is easy to show that
	$\strConProj(\algebra{A})$ forms a $*$-autonomous category,
	see \cite{barr:1979a}.
	Moreover, $\strConProj(\algebra{A})$
	can be understood as the category of linear adjoints in the 
	linear distributive category $\injstrConBimod(\algebra{A})_\sym$,
	analogous to the classical fact that finitely generated 
	projective modules can be considered as the dualizable objects in 
	the monoidal category of modules,
	cf. \cite{egger:2010a, cockett:1999a}.
	This suggests that most of the structure on 
	$\strConProj(\algebra{A})$ can actually be derived in the more 
	abstract setting of linear distributive categories.
	But, at the moment, there seems to exist no fleshed out theory of 
	monoids, modules and their linear duals internal to linear 
	distributive categories.
	
\end{remark}

\subsubsection{Reduction}

The notion of  projectivity is compatible with the reduction functor 
of strong constraint modules.

\begin{proposition}[Reduction of projective strong constraint modules]
	\label{prop:ProjectiveStrReduction}%
	\index{reduction!projective strong constraint modules}
Let $\algebra{A}$ be an\linebreak
embedded strong constraint algebra 
and $\module{P} \in \strConMod_{\algebra{A}}$.
Then $\module{P}_\red$ is projective,
and if $\module{P} \simeq e\algebra{A}^{(M)}$
for some $M \in \injConIndSet$, then
$\module{P}_\red \simeq e_\red \algebra{A}_\red^{(M_\red)}$.
\end{proposition}

%% file: constraint-misc.tex
In this short, last section of the first part, we collect some more 
constraint algebraic notions which will be needed in 
\autoref{chap:ConstraintGeometricStructures} and \autoref{chap:DeformationTheory}.
The definitions and properties should not be surprising at this point.
Thus we will restrict ourselves to what is needed later on, instead 
of giving the full-fledged theories.
In particular, we introduce constraint cochain complexes and their 
cohomology in \autoref{sec:ConHomologicalAlgebra}, while 
\autoref{sec:ConDGLAs} is concerned with 
non-associative constraint algebraic structures, such as
(differential graded) Lie algebras and Poisson algebras.

\subsection{Constraint Cochain Complexes}
\label{sec:ConHomologicalAlgebra}

Let us start to introduce $\Integers$-graded constraint modules.
Even though we could also allow for a grading by a more general 
constraint set or constraint group, this is not necessary at this 
point.

\begin{definition}[Graded constraint module]\
	\label{def:GradedModules}%
	\index{graded constraint module}
	\begin{definitionlist}
		\item \label{item:GradedCoisoModule}
		A \emph{($\Integers$-)graded constraint $\field{k}$-module} 
		is a 
		$\Integers$-indexed family
		$\{\module{M}^i\}_{i\in \Integers}$ of constraint 
		$\field{k}$-modules
		$\module{M}^i \in \ConMod_\field{k}$.
		\item \label{item:MorphismGradedCoiso}
		A \emph{morphism $\{\module{M}^i\}_{i \in \Integers} 
			\longrightarrow	\{\module{N}^i\}_{i \in \Integers}$
			of graded constraint $\field{k}$-modules}
		is given by a $\Integers$-indexed family
		$\{\Phi^i\}_{i \in \Integers}$ of morphisms
		$\Phi^i \colon \module{M}^i \longrightarrow \module{N}^i$.
		\item \label{item:GradedCoisoCat}
		We denote the category of graded constraint 
		$\field{k}$-modules by
		\glsadd{ConModkbullet}$\ConMod_\field{k}^\bullet$.
	\end{definitionlist}
\end{definition}

We can always combine the indexed family of a graded constraint 
module into a single constraint module
$\module{M}^\bullet = \bigoplus_{i \in \Integers} \module{M}^i$.
Conversely, if a given constraint module $\module{M}$ decomposes 
into a direct sum indexed by $\Integers$ we write 
$\module{M}^\bullet$ if we want to emphasize the graded structure. 
This way, every constraint module can be viewed as a graded 
constraint module by placing it at $i=0$ with all other degrees being 
trivial.

A more flexible notion of morphism between graded constraint 
modules is given by a \emph{morphism of degree} $k$, i.e.  a family
$\Phi^i \colon \module{M}^i \longrightarrow \module{N}^{i+k}$.

We will use the usual induced tensor products
\begin{align}
	\label{eq:GradedTensorProduct}
	\module{M} \tensor[\field{k}] \module{N}
	= \bigoplus_{n \in \mathbb{Z}}
	\Big(\bigoplus_{k + \ell = n} \module{M}^k \tensor[\field{k}]
	\module{N}^\ell \Big) \\
	\shortintertext{and}
	\module{M} \strtensor[\field{k}] \module{N}
	= \bigoplus_{n \in \mathbb{Z}}
	\Big(\bigoplus_{k + \ell = n} \module{M}^k \strtensor[\field{k}]
	\module{N}^\ell \Big)
\end{align}
and the symmetry with the usual Koszul signs.
This turns $\ConMod_{\field{k}}^\bullet$ into a monoidal category,
which is symmetric when considering symmetric modules.

We can now introduce constraint complexes as graded 
constraint modules together with a constraint differential.

\begin{definition}[Constraint complex]\
	\label{def:ConComplex}%
	\index{constraint!complex}
\begin{definitionlist}
	\item \label{def:ConComplex_1}
	A \emph{constraint complex}
	is a graded constraint module $\module{M}^\bullet$
	together with a constraint degree $+1$ morphism
	$\delta^\bullet \colon \module{M}^\bullet
	\longrightarrow\module{M}^{\bullet+1}$
	such that $\delta \circ \delta = 0$.
	\item \label{def:ConComplex_2}
	A \emph{morphism of constraint complexes}
	is a morphism
	$\Phi \colon \module{M}^\bullet \longrightarrow
	\module{N}^\bullet$
	of graded constraint modules, such that
	$\Phi \circ \delta_\module{M}
	= \delta_\module{N} \circ \Phi$.
	\item \label{def:ConComplex_3}
	The category of constraint complexes is denoted by
	\glsadd{CoChainsConModk}$\CoChains(\ConMod_\field{k})$.
\end{definitionlist}
\end{definition}

Since morphisms of complexes commute with the differential
$\delta$, it is easy to see that we obtain a functor by
constructing the cohomology of the constraint complex.

\begin{proposition}[Constraint cohomology]
	\label{prop:ConCohomology}%
	\index{constraint!cohomology}
Let $\module{M}^\bullet \in \CoChains(\ConMod_\field{k})$
be a constraint cochain complex with differential $\delta$.  
The maps
\begin{align}
	\module{M}^i
	\longmapsto
	\functor{H}^i(\module{M}, \delta)
	=
	\ker \delta^i / \image \delta^{i-1}
\end{align}
for $i \in \Integers$ define a functor
\glsadd{cohomologyFunctor}$\functor{H} \colon \CoChains(\ConMod_\field{k})
\longrightarrow \ConMod_\field{k}^\bullet$.
\end{proposition}

\begin{remark}[(Regular) image]
	\label{remark:RegularImageAgain}%
	\index{image!regular}
Note that constraint cohomology is defined by using the
\emph{image} of morphisms of constraint modules and not the
\emph{regular image}, see \autoref{def:ImageCModk}.
However, choosing the regular image instead
would not make a difference since the $\NULL$-component of the
denominator is not used in the quotient of constraint modules,
see \autoref{def:ConQuotientkModule}.
Moreover, note that this means that in general we cannot decide 
whether
$\ker \delta = \image \delta$
by computing cohomology, but we can decide if
$\ker \delta = \regimage \delta$
holds.
\end{remark}

\subsubsection{Reduction}

Since graded constraint modules and constraint complexes are 
given by $\Integers$-indexed families of constraint modules it should 
be clear that applying the reduction functor in every degree yields
functors
$\red \colon \ConMod_\field{k}^\bullet \to
\Modules_\field{k}^\bullet$ and
$\red \colon \CoChains(\ConMod_\field{k}) \to
\CoChains(\Modules_\field{k})$.

\begin{proposition}[Cohomology vs. reduction]
	\label{prop:CohomologyCommutesWithReduction}%
	\index{reduction!cohomology}
There exists a natural isomorphism
such that
\begin{equation}
	\begin{tikzcd}
		\CoChains(\ConMod_\field{k})
		\arrow[r,"\functor{H}"]
		\arrow[d,"\red"{swap}]
		& \ConMod_\field{k}
		\arrow[d,"\red"] \\
		\CoChains(\Modules_\field{k})
		\arrow[r,"\functor{H}"]
		&\Modules_\field{k}
	\end{tikzcd}
\end{equation}
commutes.
\end{proposition}

\begin{proof}
Define $\eta$ for every
$\module{M} \in \CoChains(\ConMod_\field{k})$ by
\begin{align*}
	\label{prop:CohomologyCommutesWithReduction_eq}
	\eta(\module{M}) \colon \functor{H}(\module{M})_\red
	\ni \big[[x]_\functor{H} \big]_\red
	\mapsto \big[ [x]_\red\big]_\functor{H}
	\in \functor{H}(\module{M}_\red).
\end{align*}
For
$\delta_\Wobs^{i-1}y \in \image \delta_\Wobs^{i-1}$
we have
$[\delta_\Wobs^{i-1}y ]_\red = \delta_\red^{i-1}[y]_\red$
and hence
$[[\delta_\Wobs^{i-1}y]_\red]_\functor{H} = 0$.
Moreover, for
$[x_0]_\functor{H} \in \functor{H}(\module{M})_\Null$
we have
$x_0 \in \module{M}^i_\Null$
and hence
$[[x_0]_\red]_\functor{H} = 0$.
Thus $\eta$ is well-defined.
Similarly, it can be shown that the inverse
$\eta^{-1}(\module{M}) \colon \functor{H}(\module{M}_\red) 
\longrightarrow \functor{H}(\module{M})_\red$
given by
$[[x]_\red]_\functor{H} \mapsto [[x]_\functor{H}]_\red$
is well-defined.
Finally, for
$\Phi \colon \module{M}^\bullet \longrightarrow 
\module{N}^\bullet$
we have
\begin{align*}
	\Big(\eta(\module{N})\circ
	\big[[\Phi^i]_\functor{H}\big]_\red\Big)
	\big(\big[[x]_\functor{H}\big]_\red\big)
	&= \Big(\eta(\module{N}) \Big)
	\big(\big[[\Phi^i(x)]_\functor{H}\big]_\red\big)\\
	&= \big[[\Phi^i(x)]_\red\big]_\functor{H}\\
	&= \Big(\big[[\Phi^i]_\red\big]_\functor{H} 
	\circ \eta(\module{M}) \Big)
	\big(\big[[x]_\functor{H}\big]_\red\big),
\end{align*}
showing that $\eta \colon \red \circ \functor{H} \Rightarrow \functor{H} \circ \red$ is indeed a natural isomorphism.
\end{proof}

A morphism
$\Phi \colon \module{M}^\bullet \to \module{N}^\bullet$
of constraint cochain complexes is called a
\emph{quasi-isomorphism} if the induced map 
$\functor{H}(\Phi)$ is an isomorphism of constraint modules.
We remark that the reduction functor
$\red \colon \CoChains(\ConMod_\field{k}) \to 
\CoChains(\Modules_\field{k})$
maps quasi-isomorphisms of constraint complexes to 
quasi-isomorphisms of cochain complexes.

\subsection{Constraint Lie Algebras}
\label{sec:ConDGLAs}

Let us collect some constraint notions involving brackets instead of associative compositions.
These notions will be important for our notions of constraint vector fields as introduced in
\autoref{sec:ConMultiVect} and constraint deformation theory, see \autoref{sec:ConDeformationFunctor}.

\begin{definition}[Constraint Lie algebra]\
	\label{def:ConLieAlg}
	\index{constraint!Lie algebra}
A \emph{constraint Lie algebra} is a constraint $\field{k}$-module 
$\liealg{g}$ together with a bracket
\glsadd{LieBracket}
\begin{equation}
	[\argument, \argument] \colon \liealg{g} \tensor[\field{k}] 
	\liealg{g} \to \liealg{g},
\end{equation}
with $[\argument, \argument] \circ \Delta = 0$
and fulfilling the usual Jacobi identity in every component.
Here $\Delta(\xi) = \xi \tensor \xi$
denotes the usual diagonal.
\end{definition}

Equivalently, a constraint Lie algebra is given by a Lie algebra 
morphism $\iota_\liealg{g} \colon \liealg{g}_\Wobs \to 
\liealg{g}_\Total$
between two Lie algebras $\liealg{g}_\Wobs$ and $\liealg{g}_\Total$
together with a Lie ideal $\liealg{g}_\Null \subseteq 
\liealg{g}_\Wobs$.
Then a morphism of constraint Lie algebras is simply a morphism of constraint $\field{k}$-modules
such that it is a Lie algebra morphism on both $\TOTAL$- and $\WOBS$-component.
We denote the category of constraint Lie algebras by \glsadd{ConLieAlg}$\ConLieAlg$.

\begin{example}\
	\label{ex:ConLieAlg}
\begin{examplelist}
	\item \index{constraint!endomorphisms}
	Let $\module{E}$ be a constraint $\field{k}$-module.
	The internal endomorphisms
	$\ConEnd_\field{k}(\module{E})$
	are a constraint Lie algebra given by the usual commutator
	$[\argument,\argument]^{\module{E}_\Total}$ on
	$\ConEnd_\field{k}(\algebra{E})_\Total$
	and the pair
	$([\argument, \argument]^{\module{E}_\Total}, 
	[\argument,\argument]^{\module{E}_\Wobs})$
	on
	$\ConEnd_\field{k}(\algebra{E})_\Wobs$.
	\item \index{constraint!derivations}
	Let $\algebra{A} \in \ConAlg$ be a constraint algebra.
	The constraint derivations $\ConDer(\algebra{A})$
	as introduced in \autoref{prop:ConDer}
	forms a constraint Lie algebra
	which can be seen as a constraint Lie subalgebra
	of $\ConEnd_\field{k}(\algebra{A})$.
\end{examplelist}
\end{example}

Note that even for a strong constraint algebra, we only obtain
a constraint Lie algebra $\ConDer(\algebra{A})$, and not 
a strong constraint Lie algebra, i.e. a constraint Lie algebra
with bracket defined on $\liealg{g} \strtensor[\field{k}] \liealg{g}$.

We can now state the definition of a constraint Lie-Rinehart algebra,
cf. \cite{rinehart:1963a,huebschmann:2003a} for the classical notion.

\begin{definition}[Constraint Lie-Rinehart algebra]
	\label{def:ConLieRinehartAlg}
	\index{constraint!Lie-Rinehart algebra}
A \emph{constraint Lie-Rinehart algebra}
consists of the following data:
\begin{definitionlist}
	\item A commutative constraint algebra $\algebra{A}$.
	\item A constraint $\algebra{A}$-module
	$\liealg{g}$ together with a Lie algebra structure
	$[\argument, \argument]$.
	\item A constraint morphism $\rho \colon \liealg{g} \to 
	\ConDer(\algebra{A})$ of constraint Lie algebras and constraint
	$\algebra{A}$-modules.
\end{definitionlist}
such that
\begin{equation}
	[\xi, a \cdot \eta] = \rho(\xi)(a) \cdot \eta + a[\xi,\eta]
\end{equation}
holds for all $\xi, \eta \in \liealg{g}_{\Total/\Wobs}$ and $a \in 
\algebra{A}_{\Total/\Wobs}$.
\end{definition}

Let us continue to combine constraint Lie algebras with constraint 
complexes.
We state the definition directly as pairs of differential graded Lie 
algebras (DGLAs).

\begin{definition}[Constraint differential graded Lie algebra]\
	\label{def:ConDGLA}
	\index{constraint!differential graded Lie algebra}
\begin{definitionlist}
	\item \label{item:CoisoDGLA}
	A \emph{constraint DGLA} $\liealg{g}$ over $\field{k}$ is a 
	pair 
	of DGLAs
	$(\liealg{g}_\Total^\bullet, [\argument,\argument]_\Total, 
	\D_\Total)$ and
	$(\liealg{g}^\bullet_\Wobs,[\argument,\argument]_\Wobs,\D_\Wobs)$
	over $\field{k}$ together with a degree $0$ morphism
	$\iota_\liealg{g} \colon \liealg{g}_\Wobs^\bullet
	\to \liealg{g}_\Total^\bullet$ of DGLAs and a
	graded Lie ideal
	$\liealg{g}_\Null^\bullet \subset \liealg{g}_\Wobs^\bullet$
	such that
	$\D_\Wobs (\liealg{g}_\Null^\bullet) \subseteq
	\liealg{g}_\Null^{\bullet +1}$.
	\item \label{item:CoisoDGLAMorph}
	For two constraint DGLAs
	$\liealg{g}$ and $\liealg{h}$, a
	\emph{morphism
		$\Phi \colon \liealg{g}^\bullet \to	\liealg{h}^\bullet$
		of constraint DGLAs} is a pair of DGLA morphisms
	$\Phi_\Total \colon \liealg{g}^\bullet_\Total \to
	\liealg{h}^\bullet_\Total$ and
	$\Phi_\Wobs \colon \liealg{g}^\bullet_\Wobs \to
	\liealg{h}^\bullet_\Wobs$ such that
	$\Phi_\Total \circ \iota_\liealg{g} 
	= \iota_\liealg{h} \circ \Phi_\Wobs$
	and
	$\Phi_\Wobs(\liealg{g}^\bullet_\Null) 
	\subseteq \liealg{h}^\bullet_\Null$.
	\item \label{item:CoisoDGLACat}
	The category of constraint DGLAs will be denoted by
	\glsadd{ConDGLieAlg}$\ConDGLieAlg$.
\end{definitionlist}
\end{definition}

Note that a morphism of constraint DGLAs can equivalently be
understood as a morphism of constraint modules such that its
components are DGLA morphisms.
A constraint Lie algebra is a constraint DGLA with trivial 
differential and concentrated in degree $0$.
Similarly, a constraint graded Lie algebra can be defined as a 
constraint DGLA with trivial differential.

\index{constraint!cohomology}
Since every constraint DGLA $\liealg{g}$ is, in particular, a 
constraint complex we can always construct its corresponding 
cohomology $\functor{H}(\liealg{g})$.
Moreover, every morphism
$\Phi \colon \liealg{g}^\bullet \to \liealg{h}^\bullet$
of constraint DGLAs is a morphism of constraint complexes and 
therefore it induces a morphism
$\functor{H}(\Phi) \colon \functor{H}^\bullet(\liealg{g})
\to \functor{H}^\bullet(\liealg{h})$
on cohomology.
Clearly,
$\functor{H}(\liealg{g})$
is a constraint graded Lie algebra and every induced morphism
$\functor{H}(\Phi)$
is a morphism of constraint graded Lie algebras.
If $\functor{H}(\Phi)$ is an isomorphism we call $\Phi$ a 
\emph{quasi-isomorphism}.

A special case of a constraint Lie algebra which will be important in 
reformulating coisotropic reduction in constraint terms is that of a 
constraint Poisson algebra.

\begin{definition}[Constraint Poisson algebra]
	\label{def:ConPoissonAlg}
	\index{constraint!Poisson algebra}
A \emph{constraint Poisson algebra} is a constraint algebra
$\algebra{A}$ together with a constraint Lie bracket
\glsadd{PoissonBracket}$\{\argument, \argument\}$ such that
$\{a, \argument\} \in \ConDer(\algebra{A})_{\Total/\Wobs/\Null}$
for every $a \in \algebra{A}_{\Total/\Wobs/\Null}$.
\end{definition}

In other words, a constraint Poisson algebra consists of a morphism
$\iota_\algebra{A} \colon \algebra{A}_\Wobs \to \algebra{A}_\Total$
of Poisson algebras together with a Poisson ideal
$\algebra{A}_\Null \subseteq \algebra{A}_\Wobs$.

\begin{example}
	In the study of singular Riemannian foliation in \cite{nahari.strobl:2022a}
	so-called $\mathcal{I}$-Poisson manifolds are introduced.
	An $\mathcal{I}$-Poisson manifold is a Poisson manifold together with
	a locally finitely generated subsheaf $\mathcal{I}$ of Poisson subalgebras of $\Cinfty(M)$.
	Every such $\mathcal{I}$-Poisson manifold
	induces a constraint Poisson algebra
	$(\Cinfty(M), \normalizer(\mathcal{I}(M)), \mathcal{I}(M))$.
\end{example}

\subsubsection{Reduction}

For a constraint DGLA $(\liealg{g}, \D)$ the reduction
$\liealg{g}_\red \coloneqq \liealg{g}_\Wobs / \liealg{g}_\Null$
gives a well-defined functor
\index{reduction!differential graded Lie algebras}
\glsadd{DGLieAlg}
\begin{equation}
	\red \colon \ConDGLieAlg \to \DGLieAlg
\end{equation}
since by definition $\liealg{g}_\Null$ is a differential graded Lie ideal in
$\liealg{g}_\Wobs$.
It is then clear that reduction of constraint 
DGLAs preserves quasi-isomorphisms.
This functor clearly restricts to a reduction functor
\index{reduction!Lie algebras}
\glsadd{LieAlg}
\begin{equation}
	\red \colon \ConLieAlg \to \LieAlg.
\end{equation}
for constraint Lie algebras.
Together with the reduction of constraint derivations, see \autoref{ex:ReductionConDer}, this also
shows that a constraint Lie-Rinehart algebra $(\algebra{A},\liealg{g})$
can be reduced to a classical Lie-Rinehart algebra
\begin{equation}
	(\algebra{A},\liealg{g})_\red \coloneqq (\algebra{A}_\red, \liealg{g}_\red).
\end{equation}
Similarly, we obtain for a constraint Poisson algebra $(\algebra{A}, \{\argument, \argument\})$
a reduced Poisson algebra
\begin{equation}
	(\algebra{A},\{\argument, \argument\})_\red \coloneqq (\algebra{A}_\red, \{\argument,\argument\}_\red).
\end{equation}

%% file: geometric-structures.tex
Recall the situation of coisotropic reduction in Poisson geometry:
There we consider a coisotropic submanifold $C$ of a Poisson manifold $M$.
Then the reduced manifold $M_\red$ is given by the quotient of $C$ by its 
characteristic distribution $D \subseteq TC$, which is spanned by the Hamiltonian 
vector fields $X_f$ of functions $f$ vanishing on $C$,
and $M_\red$ carries a canonical Poisson structure.
If we forget about the Poisson structures, but keep the underlying geometric information
needed to construct the reduced manifold, we end up with a smooth version of 
constraint sets: A manifold $M$, together with a submanifold $C$ and an equivalence
relation on $C$ defined by the distribution $D$.
These so called constraint manifolds are the main object of study in this chapter.

In principle, more general notions of constraint manifolds would be possible:
The equivalence relation on $C$ need not be induced by a distribution $D$, we could as well study equivalence relations coming from discrete group actions, or even more general equivalence relations which may or may not allow for a smooth quotient space.
Nevertheless, since we are mainly interested in the coisotropic setting and there already non-trivial effects appear, we will stick with distributions in this thesis.

In \autoref{sec:ConMfld} we give a precise definition of constraint manifolds and study some first properties.
In particular, we will see that smooth functions $\ConCinfty(\mathcal{M})$ on a constraint manifold $\mathcal{M}$
carry the structure of an embedded strong constraint algebra, giving a first link
to the constraint algebraic objects from \autoref{chap:ConstraintAlgebraicStuctures}.
After introducing vector bundles over constraint manifolds in \autoref{sec:ConVectorBundles}
we will see in \autoref{sec:SectionsOfConVect} that sections of constraint vector bundles
form embedded strong constraint $\ConCinfty(\mathcal{M})$-modules.
Moreover, in \autoref{thm:strConSerreSwan} we will give a constraint version of the Serre-Swan Theorem,
showing that the category $\ConVect(\mathcal{M})$ of constraint vector bundles is 
equivalent to the category $\Proj(\ConCinfty(\mathcal{M}))$
of projective strong constraint modules.
Having established the strong relationship of constraint geometric structures with constraint algebras and modules we can proceed
to study differential forms and multivector fields on constraint manifolds in
\autoref{sec:ConCartanCalculus}, which will again carry rich algebraic structures.
Finally, in \autoref{sec:ConSymbolCalculus} we will consider differential operators on constraint manifolds
and use constraint covariant derivatives to establish a symbol calculus
on constraint manifolds, allowing to identify constraint (multi-)differential operators
with certain sections of constraint vector fields.

\section{Constraint Manifolds}
	\label{sec:ConMfld}

\subfile{constraint-manifolds}

\section{Constraint Vector Bundles}
	\label{sec:ConVectorBundles}

\subfile{constraint-vectorbundles}

\section{Sections of Constraint Vector Bundles}
	\label{sec:SectionsOfConVect}

\subfile{constraint-sections}

\section{Constraint Cartan Calculus}
	\label{sec:ConCartanCalculus}

\subfile{constraint-fieldsandforms}

\section{Constraint Symbol Calculus}
	\label{sec:ConSymbolCalculus}

\subfile{constraint-symbolcalculus}

%% file: constraint-manifolds.tex
Following our philosophy from 
\autoref{chap:ConstraintAlgebraicStuctures} we would like to define 
constraint manifolds
as some kind of manifold object internal to a category of constraint objects replacing
a classical category which is suitable for defining manifolds.
In the geometric situation it is not so clear how this can be achieved.
Looking at the classical situation there are various possibilities to generalize the definition of a smooth manifold
to the constraint setting:
We could use the classical definition by charts to define constraint manifolds.
For this we first need to introduce constraint topological spaces (which can be done, since we have a good notion of constraint subsets),
and establish an extensive theory to make sense of notions like constraint Hausdorff space, second countability etc.
Another approach could be to consider manifolds as sheaves, or more precisely, locally ringed spaces, which locally look like smooth functions on $\Reals^n$.
Then constraint manifolds should be understood as sheaves taking values in $\injstrConAlg$ which locally look like
constraint functions on $\Reals^{(n_\Total,n_\Wobs,n_\Null)}$.
All these strategies would need a considerable amount of theory building before we could even state the definition of a constraint manifold.
At one point it might be useful to develop such a theory in detail, but for our purposes it will be enough to simply define
constraint manifolds as constraint objects internal to the category $\Manifolds$ of smooth manifolds, i.e. as
a manifold $M$ together with a smooth embedded submanifold $C$ and a distribution $D \subseteq TC$ allowing
for a smooth quotient.
Such distributions will in particular be regular and integrable, and will be called \emph{simple}.
\index{simple distribution}

\begin{definition}[Constraint manifold]\
	\label{def:ConstraintManifold}
	\index{constraint!manifold}
\begin{definitionlist}
	\item A \emph{constraint manifold} $\mathcal{M} = (M_\Total,M_\Wobs,D_\mathcal{M})$
	consists of a smooth manifold $M_\Total$, a closed embedded submanifold
	$\iota_\mathcal{M} \colon M_\Wobs \to M_\Total$ and a simple distribution 
	$D_\mathcal{M} \subseteq TM_\Wobs$ on $M_\Wobs$.
	\item A smooth map
	$\phi \colon \mathcal{M} \to \mathcal{N}$
	(or \emph{constraint map}) between constraint 
	manifolds is given by a smooth map
	$\phi \colon M_\Total \to N_\Total$
	such that $\phi(M_\Wobs) \subseteq N_\Wobs$
	and $T\phi(D_\mathcal{M}) \subseteq D_\mathcal{N}$.
	\item The category of constraint manifolds and smooth maps is 
	denoted by
	\glsadd{ConMfld}$\ConMfld$.
\end{definitionlist}
\end{definition}

If we consider only a single constraint manifold we will often write
$\mathcal{M} = (M,C,D)$, with $D \subseteq TC$ the distribution on 
the closed submanifold $C \subseteq M$, instead of using subscripts.
Additionally, we will sometimes denote the inclusion of $C$ in $M$
by $\iota \colon C \to M$.

\begin{remark}\
	\label{rem:MoreGeneralConManifolds}
\begin{remarklist}
	\item So far constraint objects were also allowed to have non-injective maps from $\WOBS$- to $\TOTAL$-compo\-nents.
	Thus it would be natural to replace the submanifold $C \subseteq M$ by a smooth map $\iota \colon C \to M$.
	Nevertheless, we will stick to the simpler notion with $C$ being an embedded submanifold.
	\item There exist more equivalence relations on $C$ allowing for a smooth quotient than are given by simple distributions.
	For example, actions of discrete groups are not included in this setting.
	See \cite{serre:1964a} for Godement's theorem, which shows that the quotient by an equivalence relation $R \subseteq C \times C$ is smooth if and only if $R$ is a closed embedded submanifold and 
	$\pr_1 \colon R \to C$ is a surjective submersion.
	We chose to stick to our more special definition, since this is the situation dictated by coisotropic reduction in Poisson geometry.
	Implementing these more general features would, on one hand, lead to a more involved theory.
	On the other hand, it should be clear for most of the following results how these can be transferred to the general situation.
	\item From a geometric point of view it would be desirable to allow also for non-smooth quotients.
	In particular, one might be interested in integrable and regular distributions which are not simple.
	Many of the following results hold in this more general situation, and we will indicate whenever a result actually uses the simplicity of the distribution.
	
	In Vinogradov's secondary calculus \cite{vinogradov:1998a}, see also \cite{vitagliano:2014}, one treats the geometry of a possibly singular quotient by cohomological methods.
	This can be another way to enlarge the class of constraint manifolds.
	
	\item Another way to include non-smooth quotients would be to enlarge the category of smooth manifolds, e.g. to the category of diffeological spaces or differentiable stacks, such that quotients exists in more general situations.
	
	Then constraint objects in these categories can be studied instead.
\end{remarklist}
\end{remark}

We will call the finite constraint index set
\index{dimension!constraint manifold}
\glsadd{dim}
\begin{equation}
	\dim(\mathcal{M}) \coloneqq \big(\dim(M), \dim(C), \rank(D)\big)
\end{equation}
the \emph{constraint dimension} of the constraint manifold
$\mathcal{M} = (M,C,D)$.

\begin{example}\
	\label{ex:ConManifolds}
\begin{examplelist}
	\item \label{ex:ConManifolds_LieGroupAction}
	Let $\group{G}$ be a Lie group acting via
	$\Phi \colon \group{G} \times C \to C$ in a free and proper way 
	on a closed submanifold $C \subseteq M$.
	Then the images of the infinitesimal action
	$T_e\Phi_p \colon \liealg{g} \to TC$, for all $p \in C$, define 
	a simple distribution on $C$, inducing the structure of a 
	constraint manifold.
	\item \index{coisotropic submanifold}
	Let $C \subseteq M$ be a coisotropic submanifold of a Poisson manifold $(M,\pi)$.
	Then if the characteristic distribution
	$D$ is simple $\mathcal{M} = (M,C,D)$
	defines a constraint manifold.
	\item Every $b$-manifold \cite{guillemin.miranda.pires:2014a}, i.e. an oriented manifold $M$ together with an oriented codimension $1$ submanifold $Z$, is a constraint manifold
	$(M,Z,0)$.
	Morphisms between $b$-manifolds, so-called $b$-maps, are constraint maps with an additional transversality condition.
	\item \label{ex:ConManifolds_Euclidean}
	Let $n = (n_\Total, n_\Wobs, n_\Null)$ be a finite 
	constraint index set.
	Then $\mathbb{R}^{n_\Wobs} \subseteq \mathbb{R}^{n_\Total}$
	together with the distribution 
	$T\mathbb{R}^{n_\Null} \subseteq T\mathbb{R}^{n_\Wobs}$
	defines a constraint manifold.
	Note that by identifying $T\mathbb{R}^{n_\Null}$ with 
	$\mathbb{R}^{n_\Null}$ and 
	$T\mathbb{R}^{n_\Wobs}$ with $\mathbb{R}^{n_\Wobs}$
	this is simply a constraint vector space, see \autoref{sec:ConVectorSpaces}.
\end{examplelist}
\end{example}

As classical manifolds locally look like a patch of euclidean space, so 
do constraint manifolds locally look like a patch of
``constraint euclidean space''
$\Reals^{n} = (\Reals^{n_\Total}, \Reals^{n_\Wobs}, \Reals^{n_\Null})$,
as in \autoref{ex:ConManifolds} \ref{ex:ConManifolds_Euclidean}.
While for $p \in M \setminus C$ there is locally no additional information
to that of the manifold $M$, so we can find a neighbourhood
homeomorphic to $(\Reals^{\dim(M)}, \Reals^{\dim(M)},0)$,
this changes for $p \in C$.
In this case we can identify a neighbourhood isomorphic to
$\Reals^{\dim(\mathcal{M})}
= (\Reals^{\dim(M)}, \Reals^{\dim(C)}, \Reals^{\rank(D)})$,
as the next result shows.

\begin{lemma}[Local structure of constraint manifolds]
	\label{lem:LocalStructureConManfifold}
Let $\mathcal{M} = (M,C,D)$ be a constraint manifold.
\begin{lemmalist}
	\item If $U \subseteq M$ is open, then
	$\mathcal{M}\at{U} \coloneqq (U, U \cap C, D\at{U})$
	is a constraint manifold.
	\item For every $p \in C$ there exists a coordinate chart $(U,x)$ around $p$
	such that
	\begin{align}
		x(U \cap L_p) &= \big(\mathbb{R}^{n_\Null} \times \{0\}\big) \cap x(U)\\
		\shortintertext{and}
		x(U \cap C) &= \big(\mathbb{R}^{n_\Wobs} \times \{0\}\big) \cap x(U),
	\end{align}
	where $L_p$ denotes the leaf of the distribution $D$ through $p$ 
	and $n=(n_\Total,n_\Wobs,n_\Null) = \dim(\mathcal{M})$ is the dimension of $\mathcal{M}$.
\end{lemmalist}
\end{lemma}

\begin{proof}
	The first part is clear.
	For the second part choose a foliation chart on $C \cap U$ and extend it as a submanifold chart to $U$.
\end{proof}

We will call charts of the above form
\index{adapted chart}
\emph{adapted charts} for a given constraint manifold. 

\renewcommand{\thesubsubsection}{\thesection.\arabic{subsubsection}}
\subsubsection{Functions on Constraint Manifold}

Let $\mathcal{M} = (M,C,D)$ be a constraint manifold with 
distribution $D \subseteq TC$.
Forgetting the smooth structure on $M$ and $C$ and equipping $C$ with 
the equivalence relation induced by the foliation of $D$ gives a 
constraint set.
Obviously, this construction is functorial, giving the forgetful 
functor
\begin{equation*}
	\functor{U} \colon \ConMfld \to \injConSet.
\end{equation*}
By \autoref{ex:injstrConAlg} \ref{ex:injstrConAlg_1} the
$\mathbb{R}$-valued constraint functions on $\functor{U}(\mathcal{M})$ 
constitute an embedded strong constraint algebra.
When we equip $(\mathbb{R}, \mathbb{R},0)$ with its canonical smooth 
structure we can consider the constraint subalgebra of smooth 
functions on $\mathcal{M}$.

\begin{proposition}[Functions on constraint manifolds]\
	\label{prop:FunctionsOnConManifolds}
	\index{constraint!functions on constraint manifold}
Mapping every constraint manifold $\mathcal{M} = (M,C,D)$ to
\glsadd{ConCinfty}
\begin{equation}
\begin{split}
	\ConCinfty(\mathcal{M})_\Total
	&\coloneqq \Cinfty(M,\mathbb{R}),\\
	\ConCinfty(\mathcal{M})_\Wobs
	&\coloneqq \left\{ f \in \Cinfty(M,\mathbb{R})
	\mid \Lie_X f\at{C} = 0 \text{ for all } X \in \Secinfty(D) 
	\right\},\\
	\ConCinfty(\mathcal{M})_\Null
	&\coloneqq \left\{ f \in \Cinfty(M,\mathbb{R})
	\mid f\at{C} = 0 \right\},
\end{split}
\end{equation}
and every constraint morphism $\phi \colon \mathcal{M} \to \mathcal {N}$ between constraint manifolds to
\begin{equation}
	\phi^* \colon \ConCinfty(\mathcal{N}) \to \ConCinfty(\mathcal{M}),
	\qquad
	\phi^*(f) \coloneqq f \circ \phi
\end{equation}
defines a functor $\ConCinfty \colon \ConMfld \to 
\injstrConAlg^\opp$.
\end{proposition}

\begin{proof}
Note that $\ConCinfty(\mathcal{M})_\Wobs$
is a subalgebra of $\Cinfty(M,\mathbb{R})$ by the fact that
$\Lie_X$ is $\mathbb{R}$-linear and satisfies a Leibniz rule.
The $\NULL$-component is obviously contained in the $\WOBS$-component 
and, since it is just the vanishing ideal of $C$, it is a two-sided 
ideal in $\Cinfty(M,\mathbb{R})$.
This shows that $\ConCinfty(\mathcal{M})$ is indeed an embedded 
strong constraint algebra.
Now given a smooth constraint map
$\phi \colon \mathcal{M} \to \mathcal{N}$
we have
$(\phi^*f)(p) = f(\phi(p)) = 0$
for $f \in \ConCinfty(\mathcal{N})_\Null$
and all $p\in M_\Wobs$.
Thus $\phi^*(\ConCinfty(\mathcal{N})_\Null) \subseteq 
\ConCinfty(\mathcal{N})_\Null$.
To show that $\phi^*$ also preserves the $\WOBS$-component 
let $f \in \ConCinfty(\mathcal{N})_\Wobs$ be given.
Then for $X_p \in D_\mathcal{M}\at{p}$, $p \in \mathcal{M}_\Wobs$
we have
$X_p(\phi^*f) = T_p\phi(X_p)f = 0$
since $T_p\phi(X_p) \in D_\mathcal{N}\at{\phi(p)}$
by assumption.
This shows $\phi^*f \in \ConCinfty(\mathcal{M})_\Wobs$.
\end{proof}

\begin{example}\
	\label{ex:ConFunctionsOnConMfld}
\begin{examplelist}
	\item \index{adapted chart} 
	\label{ex:ConFunctionsOnConMfld_1}
	Let $\mathcal{M} = (M,C,D)$ be a constraint manifold
	of dimension $n = (n_\Total,n_\Wobs,n_\Null)$,
	$p \in C$ and $(U,x)$ an adapted  chart around $p$ as in
	\autoref{lem:LocalStructureConManfifold}.
	Then
	\begin{equation}
	\begin{split}
		x^i &\in \Cinfty(\mathcal{M}\at{U})_\Null \text{ if } i \in \{n_\Wobs+1, \dotsc, n_\Total\} = (n^*)_\Null, \\
		x^i &\in \Cinfty(\mathcal{M}\at{U})_\Wobs \text{ if } i \in \{n_\Null+1, \dotsc, n_\Total\} = (n^*)_\Wobs, \\
		x^i &\in \Cinfty(\mathcal{M}\at{U})_\Total \text{ if } i \in 
		\{1, \dotsc, n_\Total\} = (n^*)_\Total,
	\end{split}
	\end{equation}
	\item \label{ex:ConFunctionsOnConMfld_2}
	\index{coisotropic submanifold}
	Let $C \subseteq M$ be a coisotropic submanifold of a Poisson manifold $(M, \pi)$ and
	denote by $\mathcal{M} = (M,C,D)$ the corresponding constraint manifold.
	Then, as for any constraint manifold, 
	$\ConCinfty(\mathcal{M})_\Null = \vanishing_C$ is the vanishing ideal of $C$, and
	additionally
	\begin{equation}
		\ConCinfty(\mathcal{M})_\Wobs 
		= \Pnormalizer_C
		= \left\{ f \in \Cinfty(M) \mid \{f,g\} \in \vanishing_C \text{ for all } g \in \vanishing_C \right\}
	\end{equation}
	is the Poisson normalizer of $\vanishing_C$.
\end{examplelist}
\end{example}

\autoref{ex:ConFunctionsOnConMfld} \ref{ex:ConFunctionsOnConMfld_1}
hints at the fact that $\ConCinfty(\mathcal{M})$
can also be understood as a sheaf of embedded strong constraint 
algebras on the topological space $M$.
Let us denote the stalk of the sheaf of smooth functions on $M$ at 
the point $p$ by $\Cinfty_p(M) = \ConCinfty_p(\mathcal{M})_\Total$.
The subsets of $\ConCinfty_p(\mathcal{M})_\Total$
given by germs of functions in $\ConCinfty(\mathcal{M})_\Wobs$
and $\ConCinfty(\mathcal{M})_\Null$ will be denoted by
$\ConCinfty_p(\mathcal{M})_\Wobs$ and
$\ConCinfty_p(\mathcal{M})_\Null$, respectively.
Then it is easy to see that
\begin{equation}
	\ConCinfty_p(\mathcal{M})
	\coloneqq \big(\ConCinfty_p(\mathcal{M})_\Total,
	\ConCinfty_p(\mathcal{M})_\Wobs,
	\ConCinfty_p(\mathcal{M})_\Null\big)
\end{equation}
is the stalk of the sheaf $\ConCinfty$ of constraint functions on 
$\mathcal{M}$, and thus in particular an embedded strong constraint 
algebra.

\begin{remark} \label{rem:ConPartitionOfUnity}
	\index{partition of unity}
For any open cover $\{U_\alpha\}_{\alpha \in I}$ of a classical manifold $M$ there
exists a subordinate partition of unity given by compactly supported functions $\chi_\alpha \in \Cinfty_0(M)$.
This often allows to glue locally defined objects together by first extending every locally defined objects to a global one by multiplying with some $\chi_\alpha$.
For a constraint manifold not every open cover admits a partition of unity consisting of functions
$\chi_\alpha \in \ConCinfty(\mathcal{M})_\Wobs$.
In particular, every $U_\alpha \subseteq M$ with $U \cap C \neq \emptyset$ 
needs to be saturated.
\index{saturated subset}
\end{remark}

\begin{remark}
	\label{rem:strConManifolds}
	\index{strong constraint!manifold}
Recall that for algebraic objects we always considered the strong constraint notions alongside the
constraint ones.
The same can be done for manifolds by defining strong constraint manifolds as constraint manifolds
with a globally defined equivalence relation, i.e. with $D \subseteq TM$.
Functions on such a strong constraint manifold $\mathcal{M}$ would then be given by
\glsadd{strConCinfty}$\strConCinfty(\mathcal{M}) \in \injConAlg$ with $\strConCinfty(\mathcal{M})_\Wobs$ given by functions
globally constant along the leaves of $D$ and $\strConCinfty(\mathcal{M})_\Null$ given by those globally invariant functions vanishing on $C$.
Note that in general $\strConCinfty(\mathcal{M})$ will be a non-strong constraint algebra.
Such strong constraint manifolds appear for example in the Marsden-Weinstein reduction with a
Lie group $\group{G}$ acting on the manifold $M$.

Obviously, any strong constraint manifold $\mathcal{M}$ can be turned into a constraint manifold by forgetting the equivalence relation, i.e. the distribution $D$, outside of $C$.
This yields a forgetful functor
\glsadd{strConMfld}$\functor{U} \colon \strConMfld \to \ConMfld$.
On the algebraic side this corresponds to the strong hull of
$\strConCinfty(\mathcal{M})$, making the diagram
\begin{equation}
	\begin{tikzcd}
		\strConMfld
		\arrow[r,"\strConCinfty"]
		\arrow[d,"\functor{U}"]
		&\injConAlg^\opp
		\arrow[d,"\argument^\str"] \\
		\ConMfld
		\arrow[r,"\ConCinfty"]
		&\injstrConAlg^\opp
	\end{tikzcd}
\end{equation}
commute,
see also \autoref{prop:ForgettingVSstrongHull}.
\end{remark}

%
%
\subsubsection{Reduction}

On every constraint manifold $\mathcal{M} = (M,C,D)$
we have an equivalence relation $\sim_\mathcal{M}$ on $C$ for which 
equivalence classes coincide with the leaves of $D$.
Requiring a simple distribution simply means that
$C / D = C / \mathord{\sim}_\mathcal{M}$ is a smooth 
manifold and
$\pr \colon C \to \mathcal{M}_\red$ is a surjective submersion.
Hence by the definition of constraint manifolds the quotient
$M_\Wobs / D_\mathcal{M}$ is a smooth manifold, and smooth maps of constraint manifolds 
drop to smooth maps on the quotients:

\begin{definition}[Reduced manifold]
	\label{def:ReductionConManifolds}
	\index{reduction!manifold}
	\glsadd{Manifolds}
The functor $\red \colon \ConMfld \to \Manifolds$ given by mapping
a constraint manifold $\mathcal{M} = (M,C,D)$ to
$\mathcal{M}_\red \coloneqq C / D$ and a constraint morphism $\phi \colon \mathcal{M} \to \mathcal{N}$ to
	\begin{equation}
		\phi_\red \colon \mathcal{M}_\red \to \mathcal{N}_\red,
		\qquad
		\phi_\red([p]) \coloneqq [\phi(p)]
	\end{equation}
is called \emph{reduction functor}.
\end{definition}

Constructing the embedded strong constraint algebra of smooth functions on a constraint manifold
then commutes with reduction:

\begin{proposition}[Constraint functions vs. reduction]
	\label{prop:ConFunctionsVSReduction}
	\index{constraint!functions on constraint manifold}
There exists a natural isomorphism such that the following 
diagram commutes:
\begin{equation}
	\label{prop:ConFunctionsVSReduction_diag}
	\begin{tikzcd}
		\ConMfld
		\arrow[r,"\ConCinfty"]
		\arrow[d,"\red"{swap}]
		&\injstrConAlg^\opp
		\arrow[d,"\red"]\\
		\Manifolds
		\arrow[r,"\Cinfty"]
		&\Algebras^\opp
	\end{tikzcd}
\end{equation}
\end{proposition}

\begin{proof}
Observe that every
$f \in \ConCinfty(\mathcal{M})_\Wobs$
drops to a function $f_\red \in \Cinfty(\mathcal{M}_\red)$, 
and the kernel of this map is exactly given by the vanishing ideal
$\ConCinfty(\mathcal{M})_\Null$.
Hence we obtain an inclusion
$\ConCinfty(\mathcal{M})_\red \subseteq \Cinfty(\mathcal{M}_\red)$.
To show surjectivity of this map, choose a tubular neighbourhood $V$ of
$C$ with projection $\pr_V \colon V \to C$ and a bump function
$\chi \in \Cinfty(M,\mathbb{R})$ with $\chi\at{C} = 1$ and 
$\chi\at{M\setminus V} = 0$.
Note that the closedness of $C$ is needed for the existence of such a 
$\chi$.
Then every function $f \in \Cinfty(\mathcal{M}_\red)$
can first be pulled back to a function $\pi_\red^*f$ on $C$ and
afterwards pulled back to $V$ via $\pr_V^*(\pi_\red^*f)$,
where $\pi_\red \colon C \to \mathcal{M}_\red$ denotes the projection 
to the quotient.
Finally, we can extend it to all of $M$ using $\chi$ obtaining
$\hat{f} \coloneqq \chi\cdot (\pr_V^*(\pi_\red^*f))$.
Since $\hat{f}\at{C} = \pi_\red^*f$ we clearly get
$\hat{f} \in \ConCinfty(\mathcal{M})_\Wobs$
and $(\hat{f})_\red = f$.
Hence we get $\ConCinfty(\mathcal{M})_\red = 
\Cinfty(\mathcal{M}_\red)$.
For the naturality consider a smooth constraint map
$\phi \colon \mathcal{M} \to \mathcal{N}$.
Then for every $f \in \Cinfty(\mathcal{N})_\red$ we have
\begin{equation*}
	(\phi^*)_\red (f_\red)
	= (\phi^*(f))_\red
	= (f \circ \phi)_\red
	= f_\red \circ \phi_\red
	= (\phi_\red)^*(f_\red).
\end{equation*}
This shows that \eqref{prop:ConFunctionsVSReduction_diag} commutes up to a natural isomorphism.
\end{proof}

\begin{proposition}
	Let $\mathcal{M} = (M,C,D)$ be a constraint manifold.
	For every $p \in C$ there is a canonical isomorphism
	$\ConCinfty_p(\mathcal{M})_\red \simeq 
	\Cinfty_{[p]}(\mathcal{M}_\red)$.
\end{proposition}

\begin{proof}
Define $\eta \colon \ConCinfty_p(\mathcal{M})_\red \to 
\Cinfty_{[p]}(\mathcal{M}_\red)$ by
\begin{equation*}
	\eta([\germ_p f]) \coloneqq \germ_{[p]} f_\red.
\end{equation*}
It is obviously an algebra morphism.
To show that it is an isomorphism, we first assume that
$\eta([\germ_p f]) = 0$.
Thus there exists an open neighbourhood $U \subseteq 
\mathcal{M}_\red$ of $[p]$ such that
$f_\red\at{U} = 0$.
Then $\pi_\mathcal{M}^{-1}(U) \subseteq C$ is an open neighbourhood
of $p$ such that
$\pi_\mathcal{M}^*f_\red\at{\pi_\mathcal{M}^{-1}(U)} 
= 0$.
Since $\pi_\mathcal{M}^{-1}(U)$ is open in $C$ and $C$ is an embedded submanifold,
there exists an open 
neighbourhood $V$ of $p$ in $M$ such that
$V \cap C = \pi_\mathcal{M}^{-1}(U)$ and
$f\at{V \cap U} = 
\pi_\mathcal{M}^*f_\red\at{\pi_\mathcal{M}^{-1}(U)} = 0$.
Therefore, we have $\germ_p f \in 
\ConCinfty_p(\mathcal{M})_\Null$,
leading to $[\germ_p f] = 0$.
This shows that $\eta$ is injective.
For the surjectivity of $\eta$ recall that by
\autoref{prop:ConFunctionsVSReduction}
we have $\ConCinfty(\mathcal{M})_\red \simeq 
\Cinfty(\mathcal{M}_\red)$
and thus every
$\germ_{[p]} g \in \Cinfty_{[p]}(\mathcal{M}_\red)$
is of the form $\germ_{[p]} f_\red$
for some $f \in \ConCinfty(\mathcal{M})_\Wobs$.
\end{proof}

%% file: constraint-vectorbundles.tex
Fix a constraint manifold $\mathcal{M} = (M,C,D)$.
A vector bundle $E$ over $\mathcal{M}$ should now consist of a vector
bundle $E_\Total \to M_\Total$ which is compatible with reduction.
By our general philosophy for constructing constraint objects
we expect a subbundle $E_\Wobs \to C$ of $\iota^\#E_\Total \to C$ 
together with 
an equivalence relation on $E_\Wobs$ such that the quotient space 
defines a vector bundle over $\mathcal{M}_\red$.
This equivalence relation should be compatible with the geometry in 
two ways: First, it should identify points in a common fibre in a 
linear way, so that we obtain a linear fibre in the quotient.
And, second, it should identify different fibres over the same leaf, 
to give a well-defined vector bundle over the leaf space 
$\mathcal{M}_\red$ at all.
The first part can be implemented by requiring a subbundle 
$E_\Null \to C$ of $E_\Wobs$.
For the second part we need the notion of a partial connection 
(or partial covariant derivative), cf. \cite{bott:1972a}.

\begin{definition}[Partial connection]
	\label{def:PartialConnection}
	\index{partial connection}
	\index{covariant derivative!partial|see{partial connection}}
Let $E \to C$ be a vector bundle over a manifold $C$, and let
$D \subseteq TC$ be regular involutive distribution on $C$.
A \emph{partial $D$-connection} on $E$ is given by a bilinear map
\glsadd{CovDerivative}
\begin{equation}
	\nabla \colon \Secinfty(D) \tensor \Secinfty(E) \to \Secinfty(E)
\end{equation}
such that
\begin{equation}
	\nabla_{fX} s = f \nabla_X s
\end{equation}
and
\begin{equation}
	\nabla_X(fs) = (\Lie_Xf) s + f\nabla_Xs
\end{equation}
for all $f \in \Cinfty(C)$, $X \in \Secinfty(D)$
and $s \in \Secinfty(E)$.
\end{definition}

Note that partial $D$-connections always exist by restricting a 
connection on $E$ to $D$.
Moreover, every partial $D$-connection can be extended to a 
connection on $E$ by choosing a complement $D^\perp$ of $D$ inside
$TC$ and a partial $D^\perp$ connection, then taking the sum of those.
Given a curve $\gamma \colon I \to C$ inside a fixed leaf of $D$
connecting $p,q \in C$ we obtain corresponding parallel transport
$\Parallel_\gamma \colon E_p \to E_q$.
Let us show that this parallel transport is actually independent of the chosen extension of $\nabla$.

\begin{lemma}
	Let $E \to C$ be a vector bundle over a manifold $C$ and let
	$D,D^\perp \subseteq TC$ be subbundles such that
	$TC = D \oplus D^\perp$.
	Moreover, let $\nabla$ be a partial $D$-connection and $\nabla^\perp$ be a partial $D^\perp$-connection
	on $E$.
	For every smooth path $\gamma \colon I \to C$ such that $\dot{\gamma}(t) \in D_{\gamma(t)}$ for all $t \in I$
	the parallel transport along $\gamma$ of $\nabla + \nabla^\perp$ is independent of $\nabla^\perp$.
\end{lemma}

\begin{proof}
Let $\gamma \colon I \to C$ be a smooth curve with $\gamma(0) = p$, $\gamma(1) = q$
and $\dot{\gamma}(t) \in D_{\gamma(p)}$ for all $t \in I$.
For every $s_p \in E_p$ the parallel transport along $\gamma$ is given by
the unique $s \in \Secinfty(\gamma^\#E)$ with $s(p) = s_p$
and
$(\gamma^\#\nabla')_{\frac{\del}{\del t}}s = 0$
with $\nabla' \coloneqq \nabla + \nabla^\perp$.
The pullback covariant derivative is the unique covariant derivative on $\gamma^\#E$ such that
\begin{equation*} \label{eq:PullbackCovDer}
	\gamma^\#\big((\gamma^\#\nabla')_{\frac{\del}{\del t}} \gamma^\#u\big) = \nabla'_{\dot{\gamma}(t)} u
	\tag{$*$}
\end{equation*}
for all $u \in \Secinfty(E)$.
Since $\dot{\gamma}(t) \in D_{\gamma(p)}$ we have $\nabla'_{\dot{\gamma}(t)} u = \nabla_{\dot{\gamma}(t)} u$,
thus the right hand side of \eqref{eq:PullbackCovDer} and therefore the parallel transport does not depend on $\nabla^\perp$. 
\end{proof}

Thus every partial $D$-connection has a well-defined notion of parallel transport.
If this parallel transport is independent of the chosen (leafwise) 
path, we will call the $D$-connection $\nabla$ \emph{holonomy-free}.\index{holonomy free}
Note that every holonomy free partial connection is flat, but the converse does not hold in general.
With this we are now ready to define constraint vector bundles.

\begin{definition}[Constraint vector bundle]
	\label{def:ConVectorBundle}
Let constraint manifolds 
$\mathcal{M} = (M_\Total,M_\Wobs,D_\mathcal{M})$ and
$\mathcal{N} = (N_\Total,N_\Wobs,D_\mathcal{N})$
be given.
\begin{definitionlist}
	\item \index{constraint!vector bundle}
	A \emph{constraint vector bundle}
	$E = (E_\Total,E_\Wobs,E_\Null,\nabla)$ over $\mathcal{M}$
	consists of a vector bundle $E_\Total \to M_\Total$,
	a subbundle $E_\Wobs \to M_\Wobs$ of $\iota^\#E_\Total$,
	a subbundle $E_\Null \to M_\Wobs$ of $E_\Wobs$
	and a holonomy-free partial $D_\mathcal{M}$-connection on
	$E_\Wobs / E_\Null$.	
	\item Let $E = (E_\Total,E_\Wobs,E_\Null,\nabla^E)$
	and $F = (F_\Total, F_\Wobs, F_\Null, \nabla^F)$ be constraint 
	vector bundles over constraint manifolds 
	$\mathcal{M}$ and $\mathcal{N}$, respectively.
	A morphism $\Phi \colon E \to F$ of constraint vector 
	bundles over a smooth map 
	$\phi \colon \mathcal{M} \to \mathcal{N}$
	is given by a vector bundle morphism
	$\Phi_\Total \colon E_\Total \to F_\Total$ such that	
	\begin{propertieslist}
		\item $\iota^\#\Phi_\Total$ restricts to a vector bundle morphism
		$\Phi_\Wobs \colon E_\Wobs \to F_\Wobs$,
		 \item $\Phi_\Wobs(E_\Null) \subseteq F_\Null$ and 
		 \item $\Phi_\Wobs$ is compatible with the connections, i.e.
		 \begin{equation} \label{eq:ConstraintVBMorphism}
		 	\Phi_\Wobs^* \left( \nabla^{F^*}_{T_p\phi(v_p)} \alpha \right)
		 	= \nabla^{E^*}_{v_p}  (\Phi_\Wobs^*\alpha)
		 \end{equation}
		 for all $p \in M_\Wobs$, $v_p \in D_\mathcal{M}\at{p}$ and
		 $\alpha \in \Secinfty(F_\Wobs / F_\Null)^*$,
		 with the pullback of forms 
		 $\Phi_\Wobs^* \colon \Secinfty(E_\Wobs / E_\Null)^*
		 \to \Secinfty(F_\Wobs / F_\Null)^*$ induced by $\Phi_\Wobs$.
	\end{propertieslist}
	\item The category of constraint vector bundles is denoted by
	\glsadd{ConVect}$\ConVect.$
	For a fixed constraint manifold $\mathcal{M}$ we denote by
	\glsadd{ConVectM}$\ConVect(\mathcal{M})$ the category of constraint vector bundles
	over $\mathcal{M}$ with vector bundle morphisms over
	$\id_\mathcal{M}$.
\end{definitionlist}
\end{definition}

\begin{remark}
If we refrain from requiring simplicity of the distribution in the definition of constraint manifolds, it would be more natural to drop the holonomy-freeness in the definition of constraint vector bundles.
Instead, it seems reasonable to require a flat partial covariant derivative.
This would also bring us closer to the situation of infinitesimal ideal systems considered in \cite{jotzlean.ortiz:2014a}.
\end{remark}

For every constraint vector bundle $E$ over a constraint manifold
$\mathcal{M} = (M,C,D)$
and $p \in M$ we can consider the fibre
$E_\Total\at{p}$.
If $p \in C$ is a point in the submanifold, we have subspaces defined by the subbundles
$E_\Wobs$ and $E_\Null$, leading to a constraint vector space
\index{fiber of a constraint vector bundle}
\begin{equation}
	E\at{p} \coloneqq \big(E_\Total\at{p}, E_\Wobs\at{p}, E_\Null\at{p}\big).
\end{equation}
For $p \in M\setminus C$ we define $E\at{p} \coloneqq (E_\Total\at{p},0,0)$.
Since $M$ and $C$ are supposed to be connected the dimension of this constraint vector space is independent of the base point $p \in C$.
Thus we call the constraint index set
\index{rank of a constraint vector bundle}
\glsadd{rank}
\begin{equation}
	\rank(E) \coloneqq \big(\rank(E_\Total), \rank(E_\Wobs), \rank(E_\Null)\big)
\end{equation}
the \emph{rank of $E$}.

Note that for a morphism $\Phi \colon E \to F$ of constraint vector bundles over the identity
the requirement \eqref{eq:ConstraintVBMorphism} simplifies to
\begin{equation} \label{eq:CompWithCovDer}
	\nabla_{v_p} \Phi(s) = \Phi(\nabla_{v_p} s)
\end{equation}
for all $s \in \Secinfty(E_\Wobs / E_\Null)$ and $v_p \in D\at{p}$.
The following simple observation will be useful later on.

\begin{lemma}
	\label{lem:FiberwiseIso}
Let $\Phi \colon E  \to F$ be a morphism of constraint vector bundles over a constraint
manifold $\mathcal{M} = (M,C,D)$ covering the identity.
Then $\Phi$ is an isomorphism of constraint vector bundles if and only if
it is a fiberwise isomorphism, i.e.
$\Phi\at{p} \colon E\at{p} \to F\at{p}$
is an isomorphism of constraint vector spaces for all $p \in M$.
\end{lemma}

\begin{proof}
Since $\Phi_\Total$ is a vector bundle morphism over the identity we know that it is an isomorphism if and only if it is a fiberwise isomorphism by classical differential geometry.
The same holds for the restrictions to the subbundles $E_\Wobs$ and $E_\Null$.
The compatibility of $\Phi^{-1}$ with the covariant derivative is automatic, since using
\eqref{eq:CompWithCovDer} we have
$\Phi(\nabla_{v_p} \Phi^{-1}(t)) = \nabla_{v_p}t$,
from which
$\nabla_{v_p} \Phi^{-1}(t) = \Phi^{-1}(\nabla_{v_p} t)$
follows.
\end{proof}

\begin{example}
Instances of constraint vector bundles have, under different names, appeared in the literature before.
\begin{examplelist}
	\item In \cite[Def 2.2]{cabrera.ortiz:2022a} the notion of quotient 
	data \index{quotient data}
	$(q_M,K,\Delta)$ for a vector bundle
	$E \to M$ is introduced, see also \cite[§2.1]{mackenzie:2005a}.
	Here $q_M \colon M \to \tilde{M}$ denotes a surjective submersion with connected fibres,
	$K \subseteq E$ is a subbundle and
	$\delta$ is a smooth assignment taking a pair of points $x,y \in M$ on the same
	$q_M$-fibre to a linear isomorphism
	$\nabla_{x,y} \colon E_y / K_y \to E_x / K_x$.
	This directly gives a constraint vector bundle
	$(E,E,K)$ over $(M,M,\ker(Ty_M))$ with $\bar\nabla$ the 
	partial 
	connection induced by $\Delta$.
	\item By Batchelor's Theorem \cite{batchelor:1980a, bonavolonta.poncin:2013a}
	graded manifolds
	\index{graded manifold}
	of degree one correspond to vector 
	bundles 
	over manifolds.
	Under this identification a graded submanifold of a degree 
	one graded manifold corresponds
	to a constraint vector bundle $(E,\iota^\#E,F)$ over 
	$(M,C,0)$, see \cite{cueca:2019a}.
\end{examplelist}
\end{example}

\begin{example}[Trivial constraint vector bundle]\
	\label{ex:ConVectBundles}
	\label{ex:TrivialConVectorBundle}
Let $\mathcal{M} = (M,C,D)$ be a constraint manifold and 
$k = (k_\Total, k_\Wobs, k_\Null)$ a finite constraint index set.
Then
\begin{equation}
	\mathcal{M} \times \Reals^k \coloneqq 
	\big(M \times \mathbb{R}^{k_\Total}, C \times \mathbb{R}^{k_\Wobs}, C \times \mathbb{R}^{k_\Null}, \Lie\big)
\end{equation}
defines a constraint vector bundle.
Here $\Lie$ denotes the component-wise Lie derivative.
\end{example}

We will call a constraint vector bundle of that form \emph{trivial}.
\index{trivial constraint vector bundle}
Constraint vector bundles isomorphic to trivial vector bundles will be called
\emph{trivializable}.

As an important tool we need the existence of local frames adapted to the structure of
a constraint vector bundle.
For this observe that every constraint vector bundle $E = (E_\Total,E_\Wobs,E_\Null)$
over a constraint manifold $\mathcal{M} = (M,C,D)$
can be restricted to an open subset $U \subset M$ giving 
a constraint vector bundle
$E\at{U} = (E_\Total\at{U}, E_\Wobs\at{U \cap C}, E_\Null\at{U \cap C})$
over $\mathcal{M}\at{U}$.

\begin{lemma}
	\label{lem:AdaptedLocalFrames}
	\index{constraint!local frame}
Let $E = (E_\Total, E_\Wobs, E_\Null, \nabla)$ be a constraint vector bundle of rank
$k = (k_\Total, k_\Wobs, k_\Null)$ 
over a constraint manifold
$\mathcal{M} = (M,C,D)$.
Let furthermore $E_\Null^\perp \to C$ and $E_\Wobs^\perp \to C$ be subbundles of
$\iota^\#E_\Total$ such that
$E_\Wobs = E_\Null \oplus E_\Null^\perp$
and
$\iota^\# E_\Total = E_\Wobs \oplus E_\Wobs^\perp$.
Then for every $p \in C$ there exists a local frame $e_1, \dotsc, e_{k_\Total} \in \Secinfty(E_\Total\at{U})$
on an open neighbourhood $U \subseteq M$ around $p$ such that
\begin{lemmalist}
	\item $\iota^\#e_i \in \Secinfty(E_\Null\at{U \cap C})$ for all $i=1,\dotsc,k_\Null$,
	\item $\iota^\#e_i \in \Secinfty(E_\Null^\perp\at{U\cap C})$ and
	$\nabla_X \iota^\#e_i = 0$ for all  $X \in \Secinfty(D)$ and $i = k_\Null + 1, \dotsc, k_\Wobs$,
	\item $\iota^\#e_i \in \Secinfty(E_\Wobs^\perp\at{U \cap C})$ for all $i = k_\Wobs +1, \dotsc, k_\Total$.
\end{lemmalist}
\end{lemma}

\begin{proof}
	Take a local frame $g_1, \dotsc, g_{k_\Wobs - k_\Null}$ of $E_\red$ on an open neighbourhood
	$\check{V} \subseteq \mathcal{M}_\red$ of $\pi_\mathcal{M}(p)$.
	Using \autoref{prop:ReductionOfConVect} \ref{prop:ReductionOfConVect_2}
	we obtain a local frame $g_1, \dotsc, g_{k_\Wobs - k_\Null}$ of $E_\Null^\perp \simeq E_\Wobs / E_\Null$
	on the open neighbourhood $ \pi_\mathcal{M}^{-1}(\check{V})$
	 with $\nabla_X g_i = 0$ for all $X \in \Secinfty(D)$
	 and $i = 1, \dotsc, k_\Wobs - k_\Null$.
	Choose additionally local frames $f_1, \dotsc, f_{k_\Null}$ of $E_\Null$ and
	$h_1, \dotsc,  h_{k_\Total - k_\Wobs}$ of $E_\Wobs^\perp$ on a possibly smaller
	open neighbourhood $V$.
	Using a tubular neighbourhood $\pr_U \colon U \to C \cap V$ of $C \cap V$ inside $V$ we can pull back those local
	frames to a local frame 
	\begin{equation*}
		e_i \coloneqq \begin{cases}
			\pr_U^\#f_i &\text{ if } i = 1, \dotsc, k_\Null \\
			\pr_U^\#g_{i-n_\Null} &\text{ if } i = k_\Null + 1, \dotsc, k_\Wobs \\
			\pr_U^\#h_{i-n_\Wobs} &\text{ if } i = k_\Wobs  +1, \dotsc, k_\Total
		\end{cases}
	\end{equation*}
	of $E_\Total$ fulfilling the required properties.	
\end{proof}

\begin{remark}
	Even though the existence of a smooth reduced vector bundle was used in the proof of \autoref{lem:AdaptedLocalFrames} this result actually only depends on local considerations, and therefore could also be obtained for regular and integrable distributions $D$ on $C$ and 
	flat partial connections on $E_\Wobs/E_\Null$ by using foliation charts.
\end{remark}

Recall from \cite{bott:1972a} that for a given manifold $M$ with a regular involutive distribution
$D \subseteq TM$ there exists a canonical partial $D$-connection on the normal bundle
$TM / D$, the so-called \emph{Bott connection}, given by
\begin{equation} \label{eq:BottConnection}
	\nabla^\Bott_X \cc{Y} = \cc{[X,Y]}
\end{equation}
for $X \in \Secinfty(D)$ and $\cc{Y} \in \Secinfty(TM / D)$.
Here $\cc{Y}$ denotes the equivalence class of
$Y \in \Secinfty(TM)$.
With this we can now construct a constraint tangent bundle out of a constraint
manifold.

\begin{proposition}[Constraint tangent bundle]
	\index{constraint!tangent bundle}
Let $\mathcal{M} = (M, C, D)$ be a constraint manifold.
Then
\glsadd{TangentBundle}\glsadd{BottConnection}$T\mathcal{M} \coloneqq (TM,TC,D,\nabla^\Bott)$
is a constraint vector bundle over $\mathcal{M}$.
\end{proposition}

\begin{proof}
We clearly have $TC \subseteq TM$, and since $D$ is regular it is a subbundle
of $TC$.
It only remains to show that $\nabla^\Bott$ is holonomy-free.
For this let $p, q \in C$ in the same leaf be given.
Moreover, let $\gamma, \tilde{\gamma} \colon I \to C$ be paths in the leaf of $p$ and $q$
such that
$\gamma(0) = p = \tilde{\gamma}(0)$,
$\gamma(1) = q = \tilde{\gamma}(1)$.
In particular, we have $\pi_\mathcal{M} \circ \gamma = \pi_\mathcal{M} \circ \tilde{\gamma}$,
with $\pi_\mathcal{M} \colon C \to \mathcal{M}_\red$
the projection onto the leaf space.
We need to show that the parallel transport of $\cc{v_p} \in T_pC/D_p$
along $\gamma$ agrees with the parallel transport along $\gamma'$.
We have $\Parallel_{\gamma,p\to q}(\cc{v_p}) = \gamma^\#(s(1))$,
where $\gamma^\# \colon \gamma^\#(TC/D) \to TC/D$
is the canonical vector bundle morphism given by 
$\gamma^\#(t, \cc{v_{\gamma(t)}}) = \cc{v_\gamma(t)}$
and $s \in \Secinfty(\gamma^\#(TC/D))$
is the unique section with $\nabla^\#_{\frac{\del}{\del t}}s = 0$
and $\gamma^\#(s(0)) = \cc{v_p}$.
Similarly, we have $\Parallel_{\tilde{\gamma},p\to q}(\cc{v_p}) = \tilde{\gamma}^\#(\tilde{s}(1))$.
Since $D = \ker T\pi_\mathcal{M}$ we know that
$T\pi_\mathcal{M} \colon TC/D \to T\mathcal{M}_\red$ is well-defined
and induces an isomorphism
$TC/D \simeq \pi_\mathcal{M}^\#T\mathcal{M}_\red$.
With this we get a canonical isomorphism
\begin{align*}
	\gamma^\#(TC/D)
	&\simeq \gamma^\#\pi_\mathcal{M}^\#T\mathcal{M}_\red
	\simeq (\pi_\mathcal{M} \circ \gamma)^\#T\mathcal{M}_\red\\
	&\simeq (\pi_\mathcal{M} \circ \tilde{\gamma})^\#T\mathcal{M}_\red 
	\simeq \tilde{\gamma}^\#\pi_\mathcal{M}^\#T\mathcal{M}_\red
	\simeq \tilde{\gamma}^\#(TC/D)
\end{align*}
which is compatible with the pullback covariant derivatives on
$\gamma^\#(TC/D)$ and $\tilde{\gamma}^\#(TC/D)$,
respectively.
Hence $s$ and $s'$ solve the same initial value problem and therefore have to agree.
Then $\Parallel_{\gamma,p\to q}(\cc{v_p}) = \Parallel_{\tilde{\gamma},p\to q}(\cc{v_p})$,
showing that $\nabla^\Bott$ is holonomy-free. 
\end{proof}

We will call $T\mathcal{M} = (TM,TC,D,\nabla^\Bott)$ the
\emph{(constraint) tangent bundle} of $\mathcal{M}$
and write $T_p\mathcal{M}$ for the constraint tangent space $T\mathcal{M}\at{p}$ as usual.

\begin{proposition}[Constraint tangent bundle functor]
	Mapping constraint manifolds to
	their constraint tangent bundles and smooth maps
	$\phi \colon \mathcal{M} \to \mathcal{N}$ between constraint manifolds
	$\mathcal{M}$ and $\mathcal{N}$ to the tangent map
	$T\phi \colon T\mathcal{M} \to T\mathcal{N}$
	defines a functor
	\begin{equation}
		T \colon \ConMfld \to \ConVect.
	\end{equation}
\end{proposition}

\begin{proof}
For the $\TOTAL$-components the statement is clear, and since $T\phi$ is completely
determined by $T\phi \colon TM_\Total \to TN_\Total$ the only thing left to show is that $T\phi$ is actually a constraint morphism.
Since $\phi$ maps $M_\Wobs$ to $N_\Wobs$
we immediately get that $\iota^\#T\phi$ restricts to
$T\phi \colon TM_\Wobs \to TN_\Wobs$.
Moreover, by \autoref{def:ConstraintManifold} we have
$T\phi(D_\mathcal{M}) \subseteq D_\mathcal{N}$.
It remains to show that $T\phi$ is compatible with the Bott connections.
We check \eqref{eq:ConstraintVBMorphism} locally.
For this let $(\tilde{U},x)$ and $(\tilde{V},y)$ be adapted coordinates around $p$ and $\phi(p)$, respectively.
Since $\phi$ restricts to a smooth map $\phi \colon \mathcal{M}_\Wobs \to \mathcal{N}_\Wobs$
it is enough to consider $U \coloneqq \tilde{U} \cap \mathcal{M}_\Wobs$
and $V \coloneqq V \cap \mathcal{N}_\Wobs$.
Then $D_\mathcal{M}\at{U}$ is spanned by $\frac{\del}{\del x^1}, \dotsc, \frac{\del}{\del x^{n_\Null}}$
and 
$D_\mathcal{N}\at{V}$ is spanned by $\frac{\del}{\del y^j}, \dotsc \frac{\del}{\del y^{m_\Null}}$.
Thus we can identify $T\mathcal{M}_\Wobs / D_\mathcal{M}$ with the subbundle spanned by
$\frac{\del}{\del x^{n_\Null +1}}, \dotsc, \frac{\del}{\del x^{n_\Wobs}}$,
and similarly $T\mathcal{N}_\Wobs / D_\mathcal{N}$ with the subbundle spanned by
$\frac{\del}{\del y^{m_\Null +1}}, \dotsc, \frac{\del}{\del y^{m_\Wobs}}$.
Note that the projection $\cc{\argument}$ to $T\mathcal{M}_\Wobs / D_\mathcal{M}$
and to $T\mathcal{N}_\Wobs/D_\mathcal{N}$
is then given by projection on the corresponding subbundles.
We will denote these by $\pr_\Null$.
For $\phi^j_i \coloneqq \frac{\del (y^j \circ \phi \circ x^{-1})}{\del x^i}$ we have
\begin{equation*}
	T_p \phi\Big( \frac{\del}{\del x^i}\at[\Big]{p}\Big) 
	= \sum_{j=1}^{m_\Wobs} \phi^j_i(p) \frac{\del}{\del y^j}\at[\Big]{\phi(p)}
\end{equation*}
with $\phi^j_i$ constant along $\Reals^{n_\Null}$ for all $j > n_0$.
Then for any 
\begin{equation*}
	v_p = \sum_{k=1}^{n_\Null} v_p^k \frac{\del}{\del x^k} \in D_\mathcal{M}\at{p}
\end{equation*}
we have
\begin{align*}
	\left( \nabla^*_{T_p\phi(v_p)} \D y^j \right)\Big( \frac{\del}{\del y^i}\at[\Big]{\phi(p)} \Big)
	&= T_p\phi(v_p)\left( \D y^j\Big(\frac{\del}{\del y^i} \Big) \at[\Big]{\phi(p)}\right)
	- \D y^j\at{\phi(p)}\left( \nabla_{T_p\phi(v_p)} \frac{\del}{\del y^i} \right)\\
	&= - \sum_{k=1}^{n_\Null} \sum_{\ell=1}^{m_\Wobs} v_p^k \phi^\ell_k(p) \D y^j\at{\phi(p)}\left(  
	\left[ \frac{\del}{\del y^\ell}, \frac{\del}{\del y^i}  \right]\at[\Big]{\phi(p)} \right)
	= 0
\end{align*}
for all $n_\Null < i \leq n_\Wobs$ and $m_\Null < j \leq m_\Wobs$.
Thus the left hand side of \eqref{eq:ConstraintVBMorphism} vanishes.
For the right hand side we compute
\begin{align*}
	\left(\nabla_{v_p} (T\phi)^*\D y^j \right)\Big(\frac{\del}{\del x^i}\Big)\at[\Big]{p}
	&= v_p\bigg( \big((T\phi)^*\D y^j\big)\Big(\frac{\del}{\del x^i}\Big) \bigg)
	- \big((T\phi)^* \D y^j\big)\left(\nabla_{v_p} \frac{\del}{\del x^i} \right)\at[\Big]{p}
	\shortintertext{with}
	\big((T\phi)^* \D y^j\big)\left(\nabla_{v_p} \frac{\del}{\del x^i} \right)\at[\Big]{p}
	&= \sum_{k=1}^{n_\Null} v_p^k \big((T\phi)^* \D y^j\big)\left( \left[ \frac{\del}{\del x^k}, \frac{\del}{\del x^i} \right] \right)\at[\Big]{p}
	= 0\\
	\shortintertext{and}
	\big((T\phi)^*\D y^j\big)\Big(\frac{\del}{\del x^i}\Big)\at[\Big]{p}
	&= \D y^j\at{\phi(p)}\Big(T_p\phi \big(\frac{\del}{\del x^i}\at[\Big]{p}\big)\Big)\\
	&= \sum_{\ell=1}^{m_\Wobs} \phi^\ell_i(p) \D y^j\at{\phi(p)}\Big( \frac{\del}{\del y^\ell}\at[\Big]{\phi(p)}\Big)
	=  \phi^j_i(p).
\end{align*}
Since $v_p(\phi^j_i) = 0$ we see that also the right hand side of \eqref{eq:ConstraintVBMorphism}
vanishes.
Thus, $T \phi$ is indeed a morphism of constraint manifolds.
\end{proof}

As in classical differential geometry, we can lift the usual constructions known for constraint vector spaces, see
\autoref{sec:ConVectorSpaces}, to constraint vector bundles.
Even though we did not introduce constraint vector bundles using vector bundle charts, the following constructions correspond at least morally to what we expect from a fiberwise definition.
For the construction of the constraint homomorphism bundle we need the following lemma:

\begin{lemma}
	\label{lem:WobsNullQuotientHom}
Let $E,F$ be constraint vector bundles over a constraint manifold
$\mathcal{M} = (M,C,D)$.
Define vector bundles
\glsadd{ConHomVect}
\begin{align}
		\ConHom(E,F)_\Wobs 
		& \coloneqq \left\{ \Phi_p \in \Hom(\iota^\#E_\Total, \iota^\#F_\Total)
		\bigm| \Phi_p(E_\Wobs\at{p}) \subseteq F_\Wobs\at{p}
		\text{ and } \Phi_p(E_\Null\at{p}) \subseteq F_\Null\at{p} \right\} \\
		\shortintertext{and}
		\ConHom(E,F)_\Null 
		& \coloneqq \left\{ \Phi_p \in \Hom(\iota^\#E_\Total, \iota^\#F_\Total)
		\bigm| \Phi_p(E_\Wobs\at{p}) \subseteq F_\Null\at{p} \right\}
\end{align}
over $C$.
Then
\begin{align}
	\Theta \colon \ConHom(E,F)_\Wobs / \ConHom(E,F)_\Null
	&\to \Hom(E_\Wobs/ E_\Null, F_\Wobs/F_\Null)
	\shortintertext{defined by}
	\Theta(\cc{\Phi_p})(\cc{v_p}) &\coloneqq \cc{\Phi_p(v_p)}
\end{align}
is an isomorphism of vector bundles.
Here $\cc{\argument}$ denotes the projection to the quotient.
\end{lemma}

\begin{proof}
It is clear that $\Theta$ is a well-defined map.
Moreover, it is a vector bundle morphism 
since it is essentially given by evaluation.
The fiberwise injectivity is again clear by definition, while for the 
fiberwise surjectivity we need to choose complements
$E_\Null^\perp$ and $E_\Wobs^\perp$
 of
 $E_\Null$ inside $E_\Wobs$ and
$E_\Wobs$ inside $\iota^\#E_\Total$.
Thus $\iota^\#E_\Total = E_\Null \oplus E_\Null^\perp \oplus E_\Wobs^\perp$
and $E_\Wobs / E_\Null \simeq E_\Null^\perp$.
Then for every $\Psi_p \in \Hom(E_\Wobs/E_\Null, F_\Wobs/F_\Null)$
set $\Phi(v_p) = \Psi(v_p)$ for all $v_p \in E_\Null^\perp$
and $\Phi(v_p) = 0$ for all $v_p \in E_\Null$ or $v_p \in E_\Wobs^\perp$.
With this we have $\Theta(\Phi_p) = \Psi_p$.
Thus we have an isomorphism of vector bundles as claimed.
\end{proof}

\begin{proposition}
	\label{prop:ConstructionsCVect}
Let $\mathcal{M} = (M,C,D)$ be a constraint manifold and
$E = (E_\Total,E_\Wobs,E_\Null,\nabla^E)$ and
$F = (F_\Total,F_\Wobs,F_\Null,\nabla^F)$
constraint vector bundles over $\mathcal{M}$ 
with
$\rank(E) = (n_\Total,n_\Wobs, n_\Null)$
and
$\rank(F) = (m_\Total,m_\Wobs,m_\Null)$.
\begin{propositionlist}
	\item \label{prop:ConstructionsCVect_DirectSum}
	\index{direct sum!constraint vector bundles}
	Defining $E \oplus F$ by
	\glsadd{directSum}
	\begin{equation}
	\begin{split}
		(E \oplus F)_\Total
		&\coloneqq E_\Total \oplus F_\Total,\\
		(E \oplus F)_\Wobs
		&\coloneqq E_\Wobs \oplus F_\Wobs,\\
		(E \oplus F)_\Null
		&\coloneqq E_\Null \oplus F_\Null,\\
		\nabla^{E \oplus F}
		&\coloneqq \nabla^{E} \oplus \nabla^F
	\end{split}
	\end{equation}
	yields a constraint vector bundle over $\mathcal{M}$, called the
	\emph{direct sum}.
	For $p \in C$ it holds
	\begin{equation}
		(E \oplus F)\at{p} = E\at{p} \oplus F\at{p},
	\end{equation}
	and it follows 
	\begin{equation}
		\label{eq:ConstructionsCVect_SumRank}
		\rank (E \oplus F) = \rank(E) + \rank(F).
	\end{equation}
	\item \label{prop:ConstructionsCVect_Tensor}
	\index{tensor product!constraint vector bundles}
	Defining $E \tensor F$ by
	\glsadd{tensor}
	\begin{equation}
	\begin{split}
		(E \tensor F)_\Total
		&\coloneqq E_\Total \tensor F_\Total,\\
		(E \tensor F)_\Wobs
		&\coloneqq E_\Wobs \tensor F_\Wobs,\\
		(E \tensor F)_\Null
		&\coloneqq E_\Null \tensor F_\Wobs + E_\Wobs \tensor F_\Null,\\
		\nabla^{E \tensor F}
		&\coloneqq \nabla^{E} \tensor \id + \id \tensor \nabla^F
	\end{split}
	\end{equation}
	yields a constraint vector bundle over $\mathcal{M}$, called the \emph{tensor product}.
	For $p \in C$ it holds
	\begin{equation}
		(E \tensor F)\at{p} = E\at{p} \tensor F\at{p},
	\end{equation}
	and therefore 
	\begin{equation}
		\label{eq:ConstructionsCVect_TensorRank}
		\rank(E \tensor F) = \rank(E) \tensor \rank(F).
	\end{equation}
	\item \label{prop:ConstructionsCVect_StrTensor}
	\index{strong tensor product!constraint vector bundles}
	Defining $E \strtensor F$ by
	\glsadd{strtensor}
	\begin{equation}
	\begin{split}
		(E \strtensor F)_\Total
		&\coloneqq E_\Total \tensor F_\Total,\\
		(E \strtensor F)_\Wobs
		&\coloneqq E_\Wobs \tensor F_\Wobs
		+ E_\Null \tensor \iota^\#F_\Total 
		+ \iota^\#E_\Total \tensor F_\Null,\\
		(E \strtensor F)_\Null
		&\coloneqq E_\Null \tensor \iota^\#F_\Total 
		+ \iota^\#E_\Total \tensor F_\Null,\\
		\nabla^{E \strtensor F}
		&\coloneqq  \nabla^{E} \tensor \id + \id \tensor \nabla^F
	\end{split}
	\end{equation}
	yields a constraint vector bundle over $\mathcal{M}$, called 
	the \emph{strong tensor product}.
	For $p \in C$ it holds
	\begin{equation}
		(E \strtensor F)\at{p} = E\at{p} \strtensor F\at{p},
	\end{equation}
	and thus		
	\begin{equation}
		\label{eq:ConstructionsCVect_StrTensorRank}
		\rank(E \strtensor F) = \rank(E) \strtensor \rank(F).
	\end{equation}
	\item \label{prop:ConstructionsCVect_Dual}
	\index{dual!constraint vector bundle}
	Defining $E^*$ by
	\glsadd{DualVect}
	\begin{equation} \label{eq:ConDualBundle}
		\begin{split}
			(E^*)_\Total
			&= (E_\Total)^*,\\
			(E^*)_\Wobs
			&= \Ann_{\iota^\#E_\Total}(E_\Null),\\
			(E^*)_\Null
			&= \Ann_{\iota^\#E_\Total}(E_\Wobs),
		\end{split}
	\end{equation}
	with \glsadd{annihilator}$\Ann_{\iota^\#E_\Total}(E_\Null)$ and $\Ann_{\iota^\#E_\Total}(E_\Wobs)$ the annihilator subbundles
	of $E_\Null$ and $E_\Wobs$ with respect to $\iota^\#E_\Total$
	and $\nabla^{E^*}$ the dual covariant derivative,
	yields a constraint vector bundle over $\mathcal{M}$, called the
	\emph{dual vector bundle}.
	For $p \in C$ it holds
	\begin{equation}
		E^*\at{p} = (E\at{p})^*,
	\end{equation}
	and it follows
	\begin{equation}
		\label{eq:ConstructionsCVect_DualRank}
		\rank (E^*) 
		= \rank(E)^*.
	\end{equation}		
	\item \label{prop:ConstructionsCVect_Hom}
	\index{constraint!homomorphism bundle}
	Defining $\ConHom(E,F)$ by
	\begin{equation}
	\begin{split}
		\ConHom(E,F)_\Total 
		& \coloneqq \Hom(E_\Total, F_\Total),\\
		\ConHom(E,F)_\Wobs 
		& \coloneqq \left\{ \Phi_p \in \Hom(\iota^\#E_\Total, \iota^\#F_\Total)
		\bigm| \Phi_p(E_\Wobs\at{p}) \subseteq F_\Wobs\at{p}\right.\\
		&\qquad\qquad\qquad\quad\qquad\qquad\qquad\left.\text{ and } \Phi_p(E_\Null\at{p}) \subseteq F_\Null\at{p} \right\},\\
		\ConHom(E,F)_\Null 
		& \coloneqq \left\{ \Phi_p \in \Hom(\iota^\#E_\Total, \iota^\#F_\Total)
		\bigm| \Phi_p(E_\Wobs\at{p}) \subseteq F_\Null\at{p} \right\}, \\
		\nabla^{\Hom}_X A 
		&\coloneqq \nabla_X^F \circ A - A \circ \nabla_X^E,
	\end{split}
	\end{equation}
	where $A \in \Secinfty(\ConHom(E,F)_\Wobs / \ConHom(E,F)_\Null)$
	is identified with the module morphism\linebreak
	$A \colon \Secinfty(E_\Wobs/E_\Null) \to \Secinfty(F_\Wobs/F_\Null)$
	using \autoref{lem:WobsNullQuotientHom}
	and $X \in \Secinfty(D)$,
	yields a constraint vector bundle, called the \emph{homomorphism bundle}.
	For $p \in C$ it holds
	\begin{equation}
		\ConHom(E,F)\at{p} = \ConHom(E\at{p},F\at{p}),
	\end{equation}
	and thus
	\begin{equation}
		\label{eq:ConstructionsCVect_HomRank}
		\rank(\ConHom(E,F)) = \rank(E^*) \strtensor \rank(F).
	\end{equation}
\end{propositionlist}
\end{proposition}

\begin{proof}
	\ref{prop:ConstructionsCVect_DirectSum}: Note that the direct sum of subbundles is a subbundle of the direct sum 
	and
	$(E \oplus F)_\Wobs / (E \oplus F)_\Null \simeq (E_\Wobs / E_\Null) \oplus (F_\Wobs / F_\Null)$.
	Moreover, the parallel transport of $\nabla^{E \oplus F}$ is given by the direct sum of the parallel transports of $\nabla^E$ and $\nabla^F$, and thus it is holonomy-free.
	By definition we have $(E \oplus F)\at{p} = ((E_\Total \oplus F_\Total)\at{p}, (E_\Wobs \oplus F_\Wobs)\at{p}, (E_\Null \oplus F_\Null)\at{p})
	= E\at{p} \oplus F\at{p}$.
	And from this $\rank(E \oplus F) = \rank(E) + \rank(F)$ directly follows.
	
	\ref{prop:ConstructionsCVect_Tensor}: We need to show that $(E \tensor F)_\Null$ actually forms a subbundle
	of $E_\Wobs \tensor F_\Wobs$.
	Let $p \in C$ be given, then the dimension of
	$(E_\Null \tensor F_\Wobs)\at{p} \cap (E_\Wobs \tensor F_\Null)\at{p}
	= (E_\Null \tensor F_\Null )\at{p}$ is independent of $p$, and thus
	$(E \tensor F)_\Null$ has constant rank and therefore defines a subbundle
	of $E_\Wobs \tensor F_\Wobs$.
	The parallel transport of $\nabla^{E \tensor F}$ on 
	$(E \tensor F)_\Wobs / (E \tensor F)_\Null
	\simeq (E_\Wobs / E_\Null) \tensor (F_\Wobs / F_\Null)$
	is given by the tensor product of the parallel transports, and hence is holonomy-free.
	
	\ref{prop:ConstructionsCVect_StrTensor}: 
	With an analogous argument we see that $(E \strtensor F)_\Wobs$ and $(E \strtensor F)_\Null$
	are well-defined subbundles with
	$(E \strtensor F)_\Wobs / (E \strtensor F)_\Null
	\simeq (E_\Wobs / E_\Null) \tensor (F_\Wobs / F_\Null)$
	and holonomy-free covariant derivative.
	
	\ref{prop:ConstructionsCVect_Dual}:
	For the dual bundle we have by definition subbundles
	\begin{equation*}
		\Ann_{\iota^\#E_\Total}(E_\Wobs) \subseteq \Ann_{\iota^\#E_\Total}(E_\Null)\subseteq \iota^\#(E_\Total)^*
	\end{equation*}
	holds.
	Moreover, $\Ann_{\iota^\#E_\Total}(E_\Null) /\Ann_{\iota^\#E_\Total}(E_\Wobs) \simeq (E_\Wobs / E_\Null)^*$
	holds and since $\nabla^{E}$ is holonomy-free so is the dual covariant derivative $\nabla^{E^*}$.
	
	\ref{prop:ConstructionsCVect_Hom}: Finally, for the homomorphism bundle note that $\iota^\#\Hom(E_\Total, F_\Total) \simeq \Hom(\iota^\#E_\Total, \iota^\#F_\Total)$.
	By using adapted local frames as in \autoref{lem:AdaptedLocalFrames} it is then easy to see that
	$\ConHom(E,F)_\Wobs$ and $\ConHom(E,F)_\Null$ form subbundles of
	$\iota^\#\ConHom(E,F)_\Total$.
	Moreover, since $\nabla^{\Hom}$ is the covariant derivative obtained from the isomorphism
	$\Hom(E_\Wobs / E_\Null, F_\Wobs / F_\Null) \simeq (E_\Wobs / E_\Null)^* \tensor (F_\Wobs / F_\Null)$
	and duals as well as tensor products of holonomy-free covariant derivatives are again holonomy-free,
	so is $\nabla^{\Hom}$.
	All statements about the rank of the involved constructions follow from
	\autoref{prop:ConstuctionsCVectSpaces}.
\end{proof}

Recall from \eqref{eq:RebracketingConTensors} that the order of $\tensor$ and $\strtensor$ can in general not be changed arbitrarily.
We always have a constraint vector bundle morphism
\begin{equation}
	E \tensor (F \strtensor G) \to (E \tensor F) \strtensor G,
\end{equation}
for constraint vector bundles $E$, $F$ and $G$ over $\mathcal{M}$,
but it will in general not be an isomorphism.

\begin{example}
	\index{constraint!cotangent bundle}
Let $\mathcal{M} = (M,C,D)$ be a constraint manifold.
Then the constraint cotangent bundle is given by
\begin{equation}
\begin{split}
	(T^*\mathcal{M})_\Total 
	&= T^*M\\
	(T^*\mathcal{M})_\Wobs
	&= \Ann_{\iota^\#T^*M}(D)\\
	(T^*\mathcal{M})_\Null 
	&= \Ann_{\iota^\#T^*M}(TC).
\end{split}
\end{equation}
Note that we can canonically identify
$\Ann_{\iota^\#T^*M}(D) / \Ann_{\iota^\#T^*M}(TC) \simeq \Ann_{TC}(D)$.
Under this identification $\cc{\iota^\#\alpha}$ becomes just
the pullback (or restriction) $\iota^*\alpha \in \Ann_{TC}(D)$ of the form $\alpha$ to $C$.
Then the dual Bott connection is given by
\begin{equation}
	\nabla^{\Bott}_X \iota^*\alpha = \Lie_X \iota^*\alpha.
\end{equation}
\end{example}

Having two different notions of tensor products also leads to two 
separate notions of symmetric and 
antisymmetric powers.
We denote by
$\Sym_{\tensor}^k E$ and
$\Anti_{\tensor}^k E$ 
the symmetric and antisymmetric tensor powers with respect to 
$\tensor$ and by
$\Sym_{\strtensor}^k E$ and
$\Anti_{\strtensor}^k E$ 
the respective tensor powers with respect to $\strtensor$.
By definition of the tensor products we have
\glsadd{SymTensor}
\glsadd{AntiTensor}
\begin{equation}
\begin{split}
		(\Sym_{\tensor}^k E)_\Total &= \Sym^k E_\Total,\\
		(\Sym_{\tensor}^k E)_\Wobs &= \Sym^k E_\Wobs,\\
		(\Sym_{\tensor}^k E)_\Null &= \Sym^{k-1} E_\Wobs \vee E_\Null,
\end{split}
\qquad\qquad\qquad
\begin{split}
	(\Anti_{\tensor}^k E)_\Total &= \Anti^k E_\Total,\\
	(\Anti_{\tensor}^k E)_\Wobs &= \Anti^k E_\Wobs,\\
	(\Anti_{\tensor}^k E)_\Null &= \Anti^{k-1} E_\Wobs \wedge E_\Null,
\end{split}
\end{equation}
with $\vee$ denoting the symmetric tensor product.
Similarly, we have
\glsadd{SymStrTensor}
\glsadd{AntiStrTensor}
\begin{equation}
\begin{split}
	(\Sym_{\strtensor}^k E)_\Total &= \Sym^k E_\Total,\\
	(\Sym_{\strtensor}^k E)_\Wobs &= \Sym^k E_\Wobs +  \Sym^{k-1} E_\Total \vee E_\Null,\\
	(\Sym_{\strtensor}^k E)_\Null &= \Sym^{k-1} E_\Total \vee E_\Null,
\end{split}
\qquad\qquad\qquad
\begin{split}
	(\Anti_{\strtensor}^k E)_\Total &= \Anti^k E_\Total,\\
	(\Anti_{\strtensor}^k E)_\Total &= \Anti^k E_\Wobs + \Anti^{k-1} E_\Total \wedge E_\Null,\\
	(\Anti_{\strtensor}^k E)_\Total &= \Anti^{k-1} E_\Total \wedge E_\Null
\end{split}
\end{equation}
for the strong constraint tensor product.
Here we suppressed the pullback $\iota^\#$ for the $\TOTAL$-bundles, since
from the context it is clear that we only can take tensor products of vector bundles
over the submanifold.

We can now determine how these constructions interact.
To show these, we essentially apply the results from
\autoref{prop:ConstuctionsCVectSpaces} 
fiberwise.

\begin{proposition}
		\label{prop:PropsOfStrongDuals}
Let $\mathcal{M} = (M,C,D)$ be a constraint manifold and let 
$E$ and $F$ be constraint vector bundles over $\mathcal{M}$.
\begin{propositionlist}
	\item \label{prop:PropsOfStrongDuals_1}
	We have $(E \oplus F)^* \simeq E^* \oplus F^*$.
	\item \label{prop:PropsOfStrongDuals_2}
	We have $(E \tensor F)^* \simeq E^* \strtensor F^*$.
	\item \label{prop:PropsOfStrongDuals_3}
	We have $(E \strtensor F)^* \simeq E^* \tensor F^*$.
	\item \label{prop:PropsOfStrongDuals_4}
	We have $\ConHom(E,F) \simeq E^* \strtensor F$.
\end{propositionlist}
\end{proposition}

\begin{proof}
	\ref{prop:PropsOfStrongDuals_1}: We know from classical differential geometry, that
	$\Phi(\alpha,\beta)(x,y) \coloneqq \alpha(x) + \beta(y)$
	defines an isomorphism
	$\Phi \colon E_\Total^* \oplus F_\Total^* \to (E_\Total \oplus F_\Total)^*$.
	It preserves the $\WOBS$-component, since for
	$p \in C$ and $\alpha_p \in (E^*\at{p})_\Wobs = \Ann(E_\Null\at{p})$,
	$\beta_p \in (F^*\at{p})_\Wobs = \Ann(F_\Null\at{p})$
	we have $\Phi(\alpha_p, \beta_p)(v_p,w_p) = 0$
	for all $v_p \in E_\Null\at{p}$ and $w_p \in F_\Null\at{p}$.
	This shows $\Phi((E^* \oplus F^*)_\Wobs) \subseteq (E \oplus F)^*_\Wobs$.
	Similarly, we see that $\Phi((E^* \oplus F^*)_\Null) \subseteq (E \oplus F)^*_\Null$.
	Moreover, this clearly gives isomorphisms
	$\Phi\at{p} \colon (E^* \oplus F^*)\at{p} \to (E \oplus F)^*\at{p}$
	for all $p \in C$.
	To show that $\Phi$ is compatible with the partial derivatives, we compute
	\begin{align*}
		\left( \nabla^{(E \oplus F)^*}_X \Phi(\alpha,\beta) \right) (v,w)
		&= \Lie_X (\Phi(\alpha,\beta)(v,w))
		- \Phi(\alpha,\beta)\left( \nabla^{E \oplus F}_X (v,w) \right) \\
		&= \Lie_X (\alpha(v)) + \Lie_X(\beta(w))
		- \alpha(\nabla^E_X v) - \beta (\nabla^F_X w)\\
		&= ( \nabla^{E^*}_X \alpha)(v) + (\nabla^{F^*}_X \beta)(w) \\
		&= \Phi\left( \nabla^{E^* \oplus F^*}_X (\alpha,\beta) \right)(v,w).
	\end{align*}
	Thus $\Phi$ is a morphism of constraint vector bundles.
	Since $\Phi$ is injective and we know by
	\autoref{prop:ConstructionsCVect} that
	$\rank(E^* \oplus F^*) = \rank(E)^* + \rank(F)^* = \rank((E \oplus F)^*)$,
	showing that $\Phi$ is a fiberwise isomorphism, and therefore by \autoref{lem:FiberwiseIso} an isomorphism
	of constraint vector bundles.
	
	\ref{prop:PropsOfStrongDuals_2}:
	The map $\Phi \colon E^*_\Total \tensor F^*_\Total \to (E_\Total \tensor F_\Total)^*$ defined by
	$\Phi(\alpha_p \tensor \beta_p)(v_p \tensor w_p) \coloneqq \alpha_p(v_p) \cdot \beta_p(w_p)$
	is an isomorphism of vector bundles.
	Let $\alpha_p \tensor \beta_p \in (E^* \strtensor F^*)_\Null = \Ann_{\iota^\#E_\Total}(E_\Wobs) \tensor \iota^\#F_\Total + \iota^\#E_\Total \tensor \Ann_{\iota^\#F_\Total}(F_\Wobs)$.
	Then for all $v_p \tensor w_p \in E_\Wobs \tensor F_\Wobs$ it holds
	\begin{equation*}
		\Phi(\alpha_p \tensor \beta_p)(v_p \tensor w_p) \coloneqq \alpha_p(v_p) \cdot \beta_p(w_p) = 0.
	\end{equation*}
	Thus $\Phi$ preserves the $\NULL$-subbundle.
	For $\alpha_p \tensor \beta_p \in E^*_\Wobs \strtensor F^*_\Wobs
	= \Ann_{\iota^\#E_\Total}(E_\Null) \tensor \Ann_{\iota^\#F_\Total}(F_\Null)$ we have
	\begin{equation*}
		\Phi(\alpha_p \tensor \beta_p)(v_p \tensor w_p) \coloneqq \alpha_p(v_p) \cdot \beta_p(w_p) = 0 
	\end{equation*}
	for all $v_p \tensor w_p \in E_\Null \tensor \iota^\#F_\Total + \iota^\#E_\Total \tensor F_\Null$.
	Thus $\Phi$ also preserves the $\WOBS$-component.
	It remains to show that $\Phi$ is compatible with the partial derivatives:
	\begin{align*}
		\left( \nabla^{(E \tensor F)^*}_X \Phi(\alpha \tensor \beta) \right) (v \tensor w)
		&= \Lie_X (\Phi(\alpha \tensor \beta)(v \tensor w))
		- \Phi(\alpha \tensor \beta)\left( \nabla^{E \tensor F}_X (v \tensor w) \right) \\
		&= \Lie_X (\alpha(v))\beta(w) + \alpha(v)\Lie_X(\beta(w))\\
		&\qquad- \alpha(\nabla^E_X v) \tensor \beta(w) - \alpha(v) \tensor \beta (\nabla^F_X w)\\
		&= (\nabla^{E^*}_X \alpha)(v) \tensor \beta(w) + \alpha(v) \tensor (\nabla^{F^*}_X \beta)(w) \\
		&= \Phi\left( \nabla^{E^* \tensor F^*}_X (\alpha \tensor \beta) \right)(v \tensor w).
	\end{align*}
	It is now straightforward to show that $\Phi$ is a morphism of constraint vector bundles and an
	isomorphism in every fibre, and therefore an isomorphism of constraint vector bundles.
	
	\ref{prop:PropsOfStrongDuals_3}:
	Here we can use the same map $\Phi \colon E^*_\Total \tensor F^*_\Total \to (E_\Total \tensor F_\Total)^*$
	as before.
	Which is an isomorphism by the same arguments.
	
	\ref{prop:PropsOfStrongDuals_4}:
	Consider the isomorphism $\Phi \colon E^*_\Total \tensor F_\Total \to \Hom(E_\Total,F_\Total)$ given by
	$\Phi(\alpha_p \tensor w_p)(v_p) \coloneqq \alpha_p(v_p) \cdot w_p$.
	Again, $\Phi$ becomes a constraint morphism and a fiberwise isomorphism for reasons of rank, and therefore an isomorphism.
\end{proof}

\renewcommand{\thesubsubsection}{\thesection.\arabic{subsubsection}}
\subsubsection{Reduction}

On every constraint vector bundle $E$ the subbundle $E_\Null$ 
together with the partial $D$-connection $\nabla^E$ defines an 
equivalence relation on $E_\Wobs$ by 
$v_p \sim_E w_p$ if and only if
$p \sim_\mathcal{M} q$
and there exists a path $\gamma \colon I \to C$ in the leaf of $p$ 
such that $\cc{w_q} = P_{\gamma,p\to q}(\cc{v_p})$ is the parallel 
transport of $\cc{v_p}$ along $\gamma$.
Here $\cc{\argument}$ denotes the equivalence class in $E_\Wobs / 
E_\Null$.
Since $\nabla^E$ is holonomy-free this is independent of the chosen 
leafwise path, and thus indeed gives a well-defined equivalence 
relation.

\begin{proposition}[Reduction of constraint vector bundles]
	\label{prop:ReductionOfConVect}
	\index{reduction!vector bundles}
Let $E = (E_\Total,E_\Wobs,E_\Null,\nabla^E)$ be a constraint 
vector bundle over a constraint manifold
$\mathcal{M} = (M,C,D)$.
\begin{propositionlist}
	\item \label{prop:ReductionOfConVect_1}
	There exists a unique vector bundle structure on 
	\begin{equation}
		\pr_{E_\red} \colon E_\Wobs / \mathord{\sim}_E \to \mathcal{M}_\red,
		\qquad \pr_E([v_p]) = \pi_\mathcal{M}(p),
	\end{equation}
	with $\pi_\mathcal{M} \colon C \to \mathcal{M}_\red$,
	such that the quotient map
	\begin{equation}
		\pi_E \colon E_\Wobs \to E_\Wobs / \mathord{\sim}_E,
		\qquad \pi_E(v_p) = [v_p]
	\end{equation}
	is a submersion and a vector bundle morphism over 
	$\pi_\mathcal{M}$.
	\item \label{prop:ReductionOfConVect_2}
	There exists an isomorphism
	\begin{equation}
		\Theta \colon (E_\Wobs / E_\Null) \to 
		\pi_\mathcal{M}^\#(E_\Wobs / \sim_E),
		\qquad \Theta(\cc{v_p}) = (p,[v_p])
	\end{equation}
	of vector bundles fulfilling
	\begin{equation}
		\label{eq:ReductionOfConVect_Parallel}
		\Theta^{-1}(q,[v_p]) = \Parallel_{\gamma,p\to q}(\cc{v_p})
	\end{equation}
	for $v_p \in E_\Wobs\at{p}$ and $p \sim_\mathcal{M} q$.
\end{propositionlist}	
\end{proposition}

\begin{proof}
	We can split the quotient procedure into two steps.
	First we consider the quotient vector bundle
	$E_\Wobs / E_\Null \to C$ with quotient map
	$\pi_{E_\Null} \colon E_\Wobs \to E_\Wobs / E_\Null$
	being a submersion and vector bundle morphism.
	Now the partial $D$-connection $\nabla^E$ induces an equivalence 
	relation on $E_\Wobs / E_\Null$ by
	$\cc{v_p} \sim_{\nabla^E} \cc{w_p}$
	if and only if $p \sim_\mathcal{M} q$ and
	$\cc{w_p} = \Parallel_{\gamma,p\to q}(\cc{v_p})$.
	In the language of Lie groupoids it is easy to see that
	the parallel transport of $\nabla^E$ defines a linear action
	of the Lie groupoid
	$R(\pi_\mathcal{M}) = 
	C \deco{}{\pi_\mathcal{M}}{\times}{}{\pi_\mathcal{M}} C$
	on $(E_\Wobs / E_\Null)$.
	Then \cite[Lemma 4.1]{higgins.mackenzie:1990a}
	gives the existence of a unique vector bundle structure on 
	$\pr_\nabla \colon (E_\Wobs / E_\Null) / \sim_{\nabla^E} \to 
	\mathcal{M}_\red$ such that the quotient map
	$\pi_\nabla \colon (E_\Wobs / E_\Null) \to (E_\Wobs / E_\Null) / 
	\sim_{\nabla^E}$ is a submersion and a vector bundle morphism over
	$\pi_\mathcal{M}$.
	Combining these we obtain a unique vector bundle structure on
	$E_\Wobs /\mathbin{\sim_E} \simeq (E_\Wobs / E_\Null) / \sim_{\nabla^E}$
	such that $\pi_E = \pi_{\nabla^E} \circ \pi_{E_\Null}$
	is a submersion and vector bundle morphism over $\pi_\mathcal{M}$.
	The second part is again directly given by \cite[Lemma 
	4.1]{higgins.mackenzie:1990a}.
\end{proof}

We will mostly write $(E_\Wobs / E_\Null) / \nabla^E$ instead of 
$(E_\Wobs / E_\Null) / \mathbin{\sim_{\nabla^E}} = E_\Wobs / \sim_E$.

\begin{definition}[Reduced vector bundle]
	\label{def:ReducedVectorBundle}
	Let $E = (E_\Total, E_\Wobs,E_\Null,\nabla^E)$ be a constraint vector bundle
	over a constraint manifold
	$\mathcal{M} = (M,C,D)$.
	Then the vector bundle $E_\red \coloneqq (E_\Wobs / E_\Null) / \nabla^E$ 
	over $\mathcal{M}_\red$ is called the \emph{reduced vector bundle} of $E$.
\end{definition}

Morphisms of constraint vector bundles are designed to yield well-defined morphisms between the reduced vector bundles, allowing for a reduction functor, as expected.

\begin{proposition}[Reduction functor]
	Mapping constraint vector bundles to their reduced bundles defines a functor
	$\red \colon \ConVect \to \Vect$.
\end{proposition}

\begin{proof}
	We need to show that morphisms of constraint vector bundles induce morphisms between the respective reduced bundles.
	For this let $\Phi \colon E \to F$ be a morphism of constraint vector bundles
	$E \to \mathcal{M}$ and $F \to \mathcal{N}$
	over a smooth map $\phi \colon \mathcal{M} \to \mathcal{N}$.
	Since $\Phi$ restricts to a vector bundle morphism
	$\Phi_\Wobs \colon E_\Wobs \to F_\Wobs$
	which maps the subbundle $E_\Null$ to $F_\Null$
	we obtain a well-defined vector bundle morphism
	$\Phi_\Wobs \colon E_\Wobs / E_\Null \to F_\Wobs / F_\Null$,
	which is compatible with the covariant derivates in the sense of
	\eqref{eq:ConstraintVBMorphism}.
	Now suppose that $v_p \sim_E w_q$.
	Then, by definition of the equivalence relation, we have
	$\cc{w_q} = \Parallel_{\gamma,p\to q}(\cc{v_p})$,
	which means there exists a leafwise path
	$\gamma \colon I \to M$ 
	with $\gamma(a) = p$, $\gamma(b) = q$ for some $a,b \in I$ 
	and $s \in \Secinfty(\gamma^\#(E_\Wobs / E_\Null))$
	with $s(a) = \cc{v_p}$, $s(b) = \cc{w_q}$ such that
	$\gamma^\#\nabla_{\frac{\del}{\del t}} s = 0$.
	Define now 
	\begin{equation*}
		\hat{\gamma} \coloneqq \phi \circ \gamma 
		\colon I \to N
		\qquad\text{ and }\qquad
		\hat{s} \coloneqq \Phi \circ s \colon I \to \gamma^\#\phi^\# F = \hat{\gamma}^\#F.
	\end{equation*}
	Then it holds
	\begin{align*}
		\hat\gamma^\# \nabla_{\frac{\del}{\del t}} \hat{s}
		= \gamma^\#\phi^\# \nabla{\frac{\del}{\del t}}( \Phi \circ s)
		= \Phi\big( \gamma^\# \nabla_{\frac{\del}{\del t}} s \big)
		= 0,
	\end{align*}
	where we used the fact that the pullback covariant derivative satisfies the universal property of the pullback in the fibred category of vector bundles with covariant derivatives.
	Thus we get $\Phi(\cc{w_q})	 = \Parallel_{\hat\gamma,\phi(p)\to \phi(q)} \Phi(\cc{v_p})$,
	showing that $\Phi$ preserves the equivalence relation and thus drops to a map
	$\Phi_\red \colon E_\red \to F_\red$.
	It is smooth since locally there exist sections of the projection map 
	$\pi_\red \colon E_\Wobs \to E_\red$.
	Moreover, it is clearly fiberwise linear, hence defining a vector bundle morphism
	$\Phi_\red \colon E_\red \to F_\red$.
\end{proof}

\begin{example}
	\index{trivial constraint vector bundle}
Consider a trivial bundle $\mathcal{M} \times \Reals^k$ as in \autoref{ex:ConVectBundles}.
Then
\begin{equation}
	(C \times \Reals^{k_\Wobs}) / (C \times \Reals^{k_\Null}) \simeq C \times \Reals^{k_\Wobs - k_\Null}
\end{equation}
and since the $D$-connection is just given by the Lie derivative we get
$(\mathcal{M} \times \Reals^k)_\red \simeq \mathcal{M}_\red \times \Reals^{k_\red}$.
\end{example}

\begin{proposition}
	\label{prop:RedConstructionsConVect}
Let $\mathcal{M} = (M,C,D)$ be a constraint manifold, and let
$E, F \in \Vect(\mathcal{M})$ be constraint vector bundles over 
$\mathcal{M}$.
\begin{propositionlist}
	\item \label{prop:RedConstructionsConVect_1}
	There exists a canonical isomorphism
	$(E \oplus F)_\red \simeq E_\red \oplus F_\red$.
	\item \label{prop:RedConstructionsConVect_2}
	There exists a canonical isomorphism
	$(E \tensor F)_\red \simeq E_\red \tensor F_\red$.
	\item \label{prop:RedConstructionsConVect_3}
	There exists a canonical isomorphism
	$(E \strtensor F)_\red \simeq E_\red \tensor F_\red$.
	\item \label{prop:RedConstructionsConVect_4}
	There exists a canonical isomorphism
	$\ConHom(E,F)_\red \simeq \Hom(E_\red, F_\red)$.
	\item \label{prop:RedConstructionsConVect_5}
	There exists a canonical isomorphism
	$(E^*)_\red \simeq (E_\red)^*$.
\end{propositionlist}
\end{proposition}

\begin{proof}
	The idea is the same for all parts:
	We pull all involved vector bundles back to $C$ along $\pi_\mathcal{M} \colon C \to \mathcal{M}_\red$
	and then use \autoref{prop:ReductionOfConVect} \ref{prop:ReductionOfConVect_2} to compare them.
	Since $\pi_\mathcal{M}$ is a surjective submersion this will be enough to infer
	isomorphy on $\mathcal{M}_\red$.
	For the first part we use the following sequence of isomorphisms:
	\begin{equation*}
		\pi_\mathcal{M}^\#(E\oplus F)_\red
		\simeq \frac{(E \oplus F)_\Wobs}{(E \oplus F)_\Null}
		= \frac{E_\Wobs \oplus F_\Wobs}{E_\Null \oplus F_\Null}
		\simeq \frac{E_\Wobs}{E_\Null} \oplus \frac{F_\Wobs}{F_\Null}
		\simeq \pi_\mathcal{M}^\#(E_\red \oplus F_\red).
	\end{equation*}
	Similarly, we have
	\begin{equation*}
		\pi_\mathcal{M}^\#(E \tensor F)_\red
		\simeq \frac{(E \tensor F)_\Wobs}{(E \tensor F)_\Null}
		= \frac{E_\Wobs \tensor F_\Wobs}{E_\Wobs \tensor F_\Null + E_\Null \tensor F_\Wobs}
		\simeq \frac{E_\Wobs}{E_\Null} \tensor \frac{F_\Wobs}{F_\Null}
		\simeq \pi_\mathcal{M}^\#(E_\red \tensor F_\red)
	\end{equation*}
	and
	\begin{align*}
		\pi_\mathcal{M}^\#(E \strtensor F)_\red
		&\simeq \frac{(E \strtensor F)_\Wobs}{(E \strtensor F)_\Null}
		= \frac{E_\Wobs \tensor F_\Wobs + E_\Total \tensor F_\Null + E_\Null \tensor F_\Total}
		{E_\Total \tensor F_\Null + E_\Null \tensor F_\Total}
		\simeq \frac{E_\Wobs \tensor F_\Wobs}
		{E_\Total \tensor F_\Null + E_\Null \tensor F_\Total} \\
		&\simeq \frac{E_\Wobs}{E_\Null} \tensor \frac{F_\Wobs}{F_\Null}\\
		&\simeq \pi_\mathcal{M}^\#(E_\red \tensor F_\red),
	\end{align*}
	as well as
	\begin{align*}
		\pi_\mathcal{M}^\#\ConHom(E,F)_\red
		&\simeq \frac{\ConHom(E,F)_\Wobs}{\ConHom(E,F)_\Null}
		\simeq \Hom(\frac{E_\Wobs}{E_\Null},\frac{F_\Wobs}{F_\Null})\\
		&\simeq \Hom(\pi_\mathcal{M}^\#E_\red,\pi_\mathcal{M}^\#F_\red)\\
		&\simeq \pi_\mathcal{M}^\#\Hom(E_\red,F_\red).
	\end{align*}
	The last part follows by choosing for $F$ the trivial constraint line bundle in
	\ref{prop:RedConstructionsConVect_4}.
\end{proof}

The above isomorphisms can be shown to be part of natural isomorphisms,
turning the functor $\red \colon \ConVect(\mathcal{M}) \to \Vect(\mathcal{M}_\red)$
into an additive, closed and monoidal functor with respect to both tensor products.

\begin{proposition}
	\label{prop:ReductionOfTangentBundle}
There exists a natural isomorphism making the following diagram commute:
\begin{equation}
	\begin{tikzcd}
		\ConMfld
		\arrow[r,"T"]
		\arrow[d,"\red"{swap}]
		&\ConVect
		\arrow[d,"\red"]\\
		\Manifolds
		\arrow[r,"T"]
		& \Vect
	\end{tikzcd}
\end{equation}
\end{proposition}

\begin{proof}
We construct an isomorphism $\Psi \colon \pi_\mathcal{M}^\#(T\mathcal{M})_\red \to \pi_\mathcal{M}^\#T(\mathcal{M}_\red)$.
From \autoref{prop:ReductionOfConVect} we know that $\pi_\mathcal{M}^\#(T\mathcal{M})_\red \simeq TC/D$.
Moreover, recall that we can pull back every germ $f_{[p]} \in \Cinfty_{[p]}(\mathcal{M}_\red)$
to $\pi_\mathcal{M}^*f_{[p]} \in \Cinfty_p(C)$.
Thus we can define
\begin{equation*}
	\Psi \colon TC/D
	\ni [v_p]
	\mapsto (p, v_p \circ \pi_\mathcal{M}^*)
	\in \pi_\mathcal{M}^\#T(\mathcal{M}_\red),
\end{equation*}
giving a fiberwise injective vector bundle morphism since $D$ is obviously the kernel.
To show surjectivity let $(p,w_{[p]}) \in \pi_\mathcal{M}^\#T(\mathcal{M}_\red)$ be given.
Since $\pi_\mathcal{M}$ is a surjective submersion, there exists a local section
$\sigma \colon V \to C$ on an open neighbourhood $V \subseteq \mathcal{M}_\red$ around
$[p]$.
With this we can set $v_p(f_p) \coloneqq w_{[p]}((\sigma^*f)_{[p]})$ for any
 $f \in \Cinfty(\pi_\mathcal{M}^{-1}(V))$,
 and thus $\Psi([v_p]) = (p,w_{[p]})$.
 This shows that $\Psi$ is a fiberwise isomorphism and hence an isomorphism of vector bundles.
 Then $\Psi$ induces the isomorphism $(T\mathcal{M})_\red \simeq T(\mathcal{M}_\red)$
 as required.
\end{proof}

%% file: constraint-sections.tex
In order to motivate the definition of sections of constraint vector bundles consider the total space of a constraint vector bundle $E = (E_\Total,E_\Wobs,E_\Null)$ over a constraint manifold
$\mathcal{M} = (M,C,D)$ in the following way:
The vector bundle $E_\Total$ is clearly a smooth manifold, and since $C \subseteq M$ is a closed submanifold so is
$E_\Wobs \subseteq E_\Total$.
Additionally, by identifying $E_\Wobs / E_\Null$ with $E_\Null^\perp$ such that
$E_\Wobs \simeq E_\Null \oplus E_\Null^\perp$, there is a distribution $D_E$ on $E_\Wobs$ which is given by
$D_E \coloneqq TE_\Null \oplus \Hor(E_\Null^\perp)$,
with $\Hor(E_\Null^\perp) \subseteq TE_\Wobs$ denoting the horizontal bundle constructed out of $\nabla^E$.
Thus we can understand the total space of a constraint vector bundle as a constraint manifold.
The vector bundle projection $\pr \colon E \to \mathcal{M}$ turns out to be a smooth map of constraint manifolds.
Thus a constraint section of $E$ should be a constraint map $s \colon \mathcal{M} \to E$ such that 
$\pr \circ s = \id_\mathcal{M}$.
This means in particular that $s$ restricted to $C$ yields a section $\iota^\# s$ of $E_\Wobs$.
Moreover, $\iota^\# s$ should map equivalent points in $C$ to equivalent vectors in $E_\Wobs$.
In other words, $\iota^\# s$ should either map to $E_\Null$ or be covariantly constant along the leaves of $D$.
These considerations motivate the following definition of the constraint module of sections.

\begin{proposition}[Functor of constraint sections]
	\label{prop:FunctorConSections}
	\index{sections|see {constraint sections}}
	\index{constraint!sections}
Let $\mathcal{M} = (M,C,D)$ be a constraint manifold.
Mapping a constraint vector bundle 
	$E = (E_\Total,E_\Wobs,E_\Null,\nabla)$ to
	\glsadd{ConSecinfty}
	\begin{equation}
	\begin{split}
		\ConSecinfty(E)_\Total
		&= \Secinfty(E_\Total) \\
		\ConSecinfty(E)_\Wobs
		&= \left\{ s \in \Secinfty(E_\Total) \bigm|  \iota^\#s \in \Secinfty(E_\Wobs),
		\nabla_X\cc{\iota^\#s} = 0 \text{ for all } X \in \Secinfty(D) \right\} \\
		\ConSecinfty(E)_\Null
		&= \left\{ s \in \Secinfty(E_\Total) \bigm| \iota^\#s \in 
		\Secinfty(E_\Null) \right\},
	\end{split}
	\end{equation}
	and a constraint vector bundle morphism $\Phi \colon E \to F$ over the identity to
	\begin{equation}
		\Phi \colon \ConSecinfty(E) \to \ConSecinfty(F),
		\qquad
		\Phi(s)(p) \coloneqq \Phi\big(s(p)\big)
	\end{equation}
defines a functor 
$\ConSecinfty \colon \ConVect(\mathcal{M}) \to \injstrConRMod{\ConCinfty(\mathcal{M})}$.
\end{proposition}

\begin{proof}
	First note that $\ConSecinfty(E)_\Total$ is clearly a $\ConCinfty(\mathcal{M})_\Total$-module.
	In general we have $\iota^\#(f\cdot s) = \iota^*f \cdot  \iota^\#s$ for
	$s \in \Secinfty(E_\Total)$ and $f \in \Cinfty(M)$.
	Thus for $f \in \ConCinfty(\mathcal{M})_\Wobs$ and
	$s \in \ConSecinfty(E)_\Wobs$ we have
	$\iota^\#(f\cdot s) \in \Secinfty(E_\Wobs)$
	and
	\begin{align*}
		\nabla^E_X(\cc{\iota^\#(f \cdot s)})
		= \nabla^E_X(\iota^*f \cdot \cc{\iota^\#s})
		= \Lie_X\iota^*f \cdot \cc{\iota^\#s} + \iota^*f \cdot \nabla^E_X \cc{\iota^\#s}
		= \iota^*f \cdot \nabla^E_X \cc{\iota^\#s}
		= 0
	\end{align*}
	for all $X \in \Secinfty(D)$, where we used $\Lie_X \iota^*f = 0$.
	Now let $s \in \Secinfty(E_\Total)$ and $f \in \ConCinfty(\mathcal{M})_\Null$ be given, then
	$\iota^\#(f \cdot s) = \iota^*f \cdot \iota^\#s = 0 \in \Secinfty(D)$.
	If $s \in \Secinfty(E)_\Null$ and $f \in \Cinfty(M)$, we get again
	$\iota^\# (f \cdot s) \in \Secinfty(E_\Null)$.
	Hence we see that $\ConSecinfty(E)$ is indeed a strong $\ConCinfty(\mathcal{M})$-module.
	Let now $\Phi \colon E \to F$ be a constraint morphism of constraint vector bundles.
	Then $\Phi$ can be restricted to a morphism between the $\WOBS$- or $\NULL$-components, meaning that
	$\Phi$ commutes with $\iota^\#$.
	Moreover, since $\Phi$ is by definition compatible with the partial connections, it maps flat sections
	to flat sections.
	Hence $\Phi$ induces a constraint module morphism between the modules of sections.
\end{proof}

It should be stressed that $\ConSecinfty(E)_\Wobs$ and $\ConSecinfty(E)_\Null$ consist
of globally defined sections, with additional properties on $C$.
In particular, $\ConSecinfty(E)_\Null$ consists of those sections of $E_\Total$ which on $C$
are sections of the subbundle $E_\Null$,
while $\ConSecinfty(E)_\Wobs$ consists of sections of $E_\Total$ such that on $C$
it is a section of the subbundle $E_\Wobs$ whose $E_\Null$ component can be arbitrary, but everything complementary to $E_\Null$ needs to be covariantly constant along the leaves.

\begin{example}[(Co-)Tangent bundle]
Let $\mathcal{M} = (M,C,D)$ be a constraint manifold.
\begin{examplelist}
	\item For the constraint tangent bundle $T\mathcal{M}$ we get
	\begin{equation}
	\begin{split}
		\ConSecinfty(T\mathcal{M})_\Total 
		&= \Secinfty(TM),\\
		\ConSecinfty(T\mathcal{M})_\Wobs
		&= \left\{ X \in \Secinfty(TM) \bigm| X\at{C} \in \Secinfty(TC) \text{ and}\right.\\
		&\hspace{8.5em}\left.[X,Y] \in \Secinfty(D) \text{ for all } Y \in \Secinfty(D) \right\},\\
		\ConSecinfty(T\mathcal{M})_\Null
		&= \left\{ X \in \Secinfty(TM) \bigm| X\at{C} \in \Secinfty(D) \right\},
	\end{split}
	\end{equation}
	by the definition of the Bott connection, see \eqref{eq:BottConnection}.
	\item For the constraint cotangent bundle $T^*\mathcal{M}$ we get
	\begin{equation}
		\begin{split}
			\ConSecinfty(T^*\mathcal{M})_\Total 
			&= \Secinfty(T^*M),\\
			\ConSecinfty(T^*\mathcal{M})_\Wobs
			&= \left\{ \alpha \in \Secinfty(T^*M) \bigm| \ins_X \iota^*\alpha = 0 \text{ and } \right. \\
			&\hphantom{= \left\{\alpha \in \Secinfty(T^*M) \bigm|\right.}
			\left. \Lie_X \iota^*\alpha = 0 \text{ for all } X \in \Secinfty(D) \right\},\\
			\ConSecinfty(T^*\mathcal{M})_\Null
			&= \left\{ \alpha \in \Secinfty(T^*M) \bigm| \iota^*\alpha = 0 \right\},
		\end{split}
	\end{equation}
	by the definition of the dual vector bundle in \eqref{eq:ConDualBundle}.
	In other words $\ConSecinfty(T^*\mathcal{M})_\Wobs$ are exactly those one-forms on $M$ which are basic
	when restricted to $C$, and $\ConSecinfty(T^*\mathcal{M})_\Null$ are those which vanish on $C$.
	Here we have to carefully distinguish between the pullback $\iota^\#\alpha$ as a section of the pullback bundle $\iota^\#T^*M$
	and the pullback (or restriction) $\iota^*\alpha \in \Secinfty(T^*C)$ of the form $\alpha$ along $\iota$.
\end{examplelist}
\end{example}

\begin{example}
Given a $b$-manifold $M$ with codimension $1$ submanifold $Z \subseteq M$
the constraint vector fields $\ConSecinfty(T\mathcal{M})$ are given by those vector fields on 
$M$ which are tangent to $Z$, hence they agree with the $b$-vector fields,
see \cite{guillemin.miranda.pires:2014a}.
Note that the $b$-vector fields are always sections of the so called $b$-tangent bundle.
In contrast, we will later see that $\ConSecinfty(\mathcal{M})_\Wobs$ is in general not given
by all sections of a vector bundle on $M$, since it will in general not be projective.
Thus we can also interpret constraint manifolds as generalization of $b$-manifolds
to higher codimensions.
\end{example}

\begin{example}[Constraint Lie algebroid]
	\label{ex:ConLieAlgebroid}
	\index{constraint!Lie algebroid}
We can now define a constraint Lie algebroid
as a morphism $\rho \colon E \to T\mathcal{M}$
of constraint vector bundles together with a constraint Lie bracket
$[\argument, \argument ]$ on $\ConSecinfty(E)$ such that
$\ConCinfty(\mathcal{M})$ together with
$\ConSecinfty(E)$ becomes a constraint Lie-Rinehart algebra,
see \autoref{def:ConLieRinehartAlg}.
Particular instances of constraint Lie algebroids have been introduced in
\cite{jotzlean.ortiz:2014a} as
\index{infinitesimal ideal system}\emph{infinitesimal ideal systems}.
These are equivalent to constraint Lie algebroids of the form
$\mathcal{M} = (M,M,D)$
and
$E = (A,A,K,\nabla)$.
Note that even though the $\TOTAL$- and $\WOBS$-components of $\mathcal{M}$
and $E$ agree, this is not the case for $\ConSecinfty(T\mathcal{M})$
and $\ConSecinfty(E)$.
It is then clear that the reduction of constraint Lie algebroids, and hence also
infinitesimal ideal systems, yields classical Lie algebroids over $\mathcal{M}_\red$.
Such constraint Lie algebroids will be studied in \cite{dippell.kern:2022a}.

Another example of constraint Lie algebroids is given by so-called \emph{Lie pairs}, i.e.
pairs $(A,L)$ of Lie algebroids with $L \subseteq A$ a Lie subalgebroid over a common manifold $M$.
Multivector fields and differential operators on Lie pairs have been studied in \cite{bandiera.stienon.xu:2021a,stienon.vitagliano.xu:2022a}
using methods from the theory of $L_\infty$- and $A_\infty$-algebras.
\end{example}

\begin{example}\
	\label{ex:ConSections}
Let $\mathcal{M} = (M,C,D)$ be a constraint manifold
of dimension $n = (n_\Total,n_\Wobs,n_\Null)$,
$p \in C$ and $(U,x)$ an adapted  chart around $p$ as in
\autoref{lem:LocalStructureConManfifold}.
Then
\begin{equation}
\begin{split}
	\frac{\del}{\del x^i} &\in 
	\ConSecinfty(\mathcal{M}\at{U})_\Total
	\quad\text{ if }\quad i \in \{1, \dotsc, n_\Total\},\\
	\frac{\del}{\del x^i} &\in 
	\ConSecinfty(\mathcal{M}\at{U})_\Wobs
	\quad\text{ if }\quad i \in \{1, \dotsc, n_\Wobs\}, \\
	\frac{\del}{\del x^i} &\in 
	\ConSecinfty(\mathcal{M}\at{U})_\Null 
	\quad\text{ if }\quad i \in \{1, \dotsc, n_\Null \}.
\end{split}
\end{equation}
\end{example}

This example motivates the definition of a constraint local frame.

\begin{definition}[Constraint local frame]
	\label{def:ConLocalFrame}
	\index{constraint!local frame}
Let $E = (E_\Total, E_\Wobs, E_\Null)$ be a constraint vector bundle of
rank $k = \rank(E)$ over
a constraint manifold $\mathcal{M} = (M,C,D)$.
A \emph{local frame} of $E$ on an open $U \subseteq M$, is
a local frame $e_1, \dotsc, e_{k_\Total}$ of $E_\Total$ on $U$,
such that
\begin{definitionlist}
	\item $e_1, \dotsc, e_{k_\Wobs} \in \ConSecinfty(E\at{U})_\Wobs$
	and $\iota^\#e_1,\dotsc, \iota^\#e_{k_\Wobs}$ is a local 
	frame for $E_\Wobs$ on $U \cap C$, and
	\item $e_1, \dotsc, e_{k_\Null} \in \ConSecinfty(E\at{U})_\Null$
	and $\iota^\#e_1,\dotsc, \iota^\#e_{k_\Null}$ is a local 
	frame for $E_\Null$ on $U \cap C$.
\end{definitionlist}
\end{definition}

The existence of local frames for constraint vector bundles is 
guaranteed by
\autoref{lem:AdaptedLocalFrames}.
To show that every $v_p \in E_\Total\at{p}$ is the value of some section $s \in \Secinfty(E_\Total)$
one can simply extend a local frame for $E_\Total$ to all of $M$ by means of a cut-off function.
Now for $v_p \in E_\Wobs\at{p}$ this is not so easy any more, since a cut-off function would now need to
be an element of $\ConCinfty(\mathcal{M})_\Wobs$ itself, to end up with a section in $\ConSecinfty(E)_\Wobs$.
Recall from \autoref{rem:ConPartitionOfUnity} that the existence of such cut-off functions for arbitrary open subsets can not be guaranteed in general.
Nevertheless, we can use the reduced manifold to construct such constraint sections as follows:

\begin{corollary}
	\label{lem:ExtendingToConSections}
	Let $E = (E_\Total,E_\Wobs,E_\Null,\nabla^E)$ be a constraint
	vector bundle over a constraint manifold $\mathcal{M} = (M,C,D)$.
	\begin{lemmalist}
		\item \label{lem:ExtendingToConSections_1}
		For each $p \in C$ and $v_p \in E_\Null\at{p}$ there exists an
		$s \in \ConSecinfty(E)_\Null$ such that
		$s(p) = v_p$.
		\item  \label{lem:ExtendingToConSections_2}
		For each $p \in C$ and $v_p \in E_\Wobs\at{p}$ there exists an
		$s \in \ConSecinfty(E)_\Wobs$ such that
		$s(p) = v_p$.
	\end{lemmalist}
\end{corollary}

\begin{proof}
	For the first part choose a local frame $e_1, \dotsc, e_{n_0}$
	of $E_\Null$ around $p$ with $n_0 = \rank(E_\Null)$.
	Then using $v_p = \sum_{k=1}^{n_0} v_p^k e_k(p)$
	we can define a local section $\sum_{k=1}^{n_0} v_p^k e_k$
	which we extend to a section $\tilde{s} \in \Secinfty(E_\Null) 
	\subseteq \Secinfty(\iota^\#E_\Total)$
	by means of a bump function.
	In order to extend $\tilde{s}$ to a section of $E_\Total$ choose 
	a tubular neighbourhood $V \subseteq M$ of $C$ with bundle 
	projection $\pi_V \colon V \to C$.
	Then pulling back $\tilde{s}$ to $V$ via $\pi_V$ and afterwards 
	extending to all of $M$ using a suitable bump function gives
	a globally defined section $s \in \Secinfty(E_\Total)$
	with $\iota^\#s = \tilde{s} \in \Secinfty(E_\Null)$
	and $s(p) = \tilde{s}(p) = v_p$.
	Note that the existence of such a bump function requires the closedness of $C$.
	For \ref{lem:ExtendingToConSections_2} choose a complementary 
	vector bundle $E_\Null^\perp \to C$ to $E_\Null$ inside of 
	$E_\Wobs$, i.e. $E_\Wobs = E_\Null \oplus E_\Null^\perp$ and hence
	$E_\Null^\perp \simeq E_\Wobs / E_\Null$.
	Then $v_p = v_p^0 + v_p^\perp$ with $v_p^0 \in E_\Null\at{p}$
	and $v_p^\perp \in E_\Null^\perp\at{p}$.
	By \ref{lem:ExtendingToConSections_1} we find a section
	$s_0 \in \Secinfty(E)_\Null$ such that $s_0(p) = v_p^0$.
	Now choose $\check{s} \in \Secinfty(E_\red)$ such that
	$\check{s}(\pi_\mathcal{M}(p)) = [v_p^\perp]$.
	Then by \autoref{prop:ReductionOfConVect} \ref{prop:ReductionOfConVect_2}
	we can identify $\pi_\mathcal{M}^\#\check{s}$ with a section
	$s^\perp \in \Secinfty(E_\Wobs / E_\Null)$ such that
	$\nabla_X s^\perp = 0$ for all $X \in \Secinfty(D)$.
	Then using a tubular neighbourhood as before to extend $s_0 + s^\perp$ to all of
	$M$ we obtain the desired section.
\end{proof}

\begin{remark}
	\label{rem:NonSimpleGlobalSections}
Note that the proof of \autoref{lem:ExtendingToConSections}
\ref{lem:ExtendingToConSections_1} still works if we refrain from $D$ being simple.
In the proof of the second part, however, we crucially used the smooth structure on $\mathcal{M}_\red$.
In particular, the holonomy-freeness of $\nabla$ is needed in order to extend $s$ to a section in $\ConSecinfty(E)_\Wobs$.
Thus it is not clear if this statement still holds for non-simple distributions.
\end{remark}

As a first important property of the sections functor we show that it is compatible with direct sums.

\begin{proposition}
	\label{prop:ConSectionsDirectSum}
Let $E = (E_\Total, E_\Wobs, E_\Null, \nabla^E)$ and
$F = (F_\Total, F_\Wobs, F_\Null, \nabla^F)$ be constraint vector bundles over a constraint manifold
$\mathcal{M} = (M,C,D)$.
Then
\begin{equation}
	\ConSecinfty(E \oplus F) \simeq \ConSecinfty(E) \oplus \ConSecinfty(F)
\end{equation}
as strong constraint $\ConCinfty(\mathcal{M})$-modules.
\end{proposition}

\begin{proof}
From classical differential geometry we know that
$\Phi \colon \Secinfty(E_\Total) \oplus \Secinfty(F)_\Total \to \Secinfty(E \oplus F)$
given by $\Phi(s,s')(p) \coloneqq s(p) \oplus s'(p)$
is an isomorphism of $\Cinfty(M)$-modules.
Now let $s \in \ConSecinfty(E)_\Wobs$ and $s' \in \ConSecinfty(F)_\Wobs$ be given.
Then clearly $\Phi(s,s')(p) = s(p) \oplus s'(p) \in E_\Wobs\at{p} \oplus F_\Wobs\at{p}$ for all $p \in C$.
Moreover, it holds
\begin{align}
	\nabla^\oplus_X \cc{\Phi(s,s')}
	= \nabla^E_X \cc{s} \oplus \nabla^F_X \cc{s'}
	= 0
\end{align}
by the definition of $\nabla^\oplus$ in \autoref{prop:ConstructionsCVect}.
Thus $\Phi$ preserves the $\WOBS$-component.
Next, let $s \in \ConSecinfty(E)_\Null$ and $s' \in \ConSecinfty(F)_\Null$ be given.
Then $\Phi(s,s')(p) = s(p) \oplus s'(p) \in E_\Null\at{p} \oplus F_\Null\at{p}$
for all $p \in C$
shows that $\Phi$ also preserves the $\NULL$-components.
For $\Phi$ to be a constraint isomorphism it remains to show that
$\Phi^{-1}(\ConSecinfty(E \oplus F)_\Null) = \ConSecinfty(E)_\Null \oplus \ConSecinfty(F)_\Null$, cf. \autoref{lem:indConIsos}.
For this let $t \in \ConSecinfty(E \oplus F)_\Null$ be given.
Then we know that $t = s \oplus s'$ for some
$s \in \Secinfty(E_\Total)$ and $t \in \Secinfty(F_\Total)$.
For all $p \in C$ we have
\begin{equation*}
	s(p) \oplus s'(p) = (s \oplus s')(p) = t(p)
	\in (E \oplus F)_\Null\at{p} = E_\Null\at{p} \oplus F_\Null\at{p},
\end{equation*} 
and thus $s \in \ConSecinfty(E)_\Null$ and $s' \in \ConSecinfty(F)_\Null$.
Therefore, $\Phi$ is a constraint isomorphism.
\end{proof}

Similarly, sections of constraint vector bundles are compatible with internal homs:

\begin{proposition}
	\label{prop:ConSectionsHom}
	\index{constraint!homomorphism bundle}
Let $E = (E_\Total, E_\Wobs, E_\Null, \nabla^E)$ and
$F = (F_\Total, F_\Wobs, F_\Null, \nabla^F)$ be constraint vector bundles over a constraint manifold
$\mathcal{M} = (M,C,D)$.
Then
\begin{equation}
	\ConSecinfty(\ConHom(E,F)) \simeq \ConHom_{\ConCinfty(\mathcal{M})}(\ConSecinfty(E), \ConSecinfty(F))
\end{equation}
as strong constraint $\ConCinfty(\mathcal{M})$-modules.
\end{proposition}

\begin{proof}
On the $\TOTAL$-component we have the isomorphism
\begin{equation*}
	\eta \colon \Secinfty(\Hom(E_\Total,F_\Total)) \to \Hom_{\Cinfty(M)}(\Secinfty(E_\Total),\Secinfty(F_\Total))
\end{equation*}
given by
\begin{equation*}
	\eta(A)(s)\at{p} \coloneqq A\at{p}(s\at{p})
\end{equation*}
for all $p \in M$ and $s \in \Secinfty(E_\Total)$.
We first show that $\eta$ is indeed a constraint morphism:
If $A \in \ConSecinfty(\ConHom(E,F))_\Null$, then for every $p \in C$ and $s \in \ConSecinfty(E)_\Null$
we have $\eta(A)(s)\at{p} = A\at{p}(s\at{p}) \in F_\Null\at{p}$
since $s\at{p} \in E_\Null\at{p}$.
Thus $\eta$ preserves the $\NULL$-component.
Consider now $A \in \ConSecinfty(\ConHom(E,F))_\Wobs$.
For all $p \in C$ and $s \in \ConSecinfty(E)_\Null$ we have
$\eta(A)(s)\at{p} = A\at{p}(s\at{p}) \in F_\Null\at{p}$
since $s\at{p} \in E_\Null\at{p}$.
Moreover, if $s \in \ConSecinfty(E)_\Wobs$, then
$\eta(A)(s)\at{p} = A\at{p}(s\at{p}) \in F_\Wobs\at{p}$
and
\begin{equation*}
	\nabla^F_X \cc{\eta(A)(s)\at{C}} = \cc{\eta}\big( \underbrace{\nabla^\ConHom_X \cc{A\at{C}}}_{=0}\big)(\cc{s\at{C}})
	+ \cc{\eta(A)}(\underbrace{\nabla^E_X \cc{s\at{C}}}_{=0})
	= 0.
\end{equation*}
Thus $\eta(A)(s) \in \ConSecinfty(F)_\Wobs$.
Summarizing, this shows that $\eta$ is a constraint morphism.

It remains to show that $\eta$ is regular surjective.
For this recall from classical differential geometry that for every
$A \in \Hom_{\Cinfty(M)}(\Secinfty(E_\Total),\Secinfty(F_\Total))$
the corresponding preimage is given by
$A(p)(s_p) \coloneqq A(s)(p)$
for all $s_p \in E_\Total\at{p}$ and $s \in \Secinfty(E_\Total)$
such that $s(p) = s_p$.
Note that this does not depend on the choice of the section $s$.
Here we use the usual abuse of notation.

Now let $A \in \ConHom_{\ConCinfty(\mathcal{M})}(\ConSecinfty(E),\ConSecinfty(F))_\Wobs$ be given.
Then for every $p \in C$ and $s_p \in E_\Wobs\at{p}$ there exists a section
$s \in \ConSecinfty(E)_\Wobs$ with $s(p) = s_p$ by \autoref{lem:ExtendingToConSections}.
Then we have $A(p)(s_p) = A(s)(p) \in F_\Wobs\at{p}$
since $A(s) \in \ConSecinfty(F)_\Wobs$.
Similarly, if $s_p \in E_\Null\at{p}$, then
there exists $s \in \ConSecinfty(E)_\Null$ with $s(p) = s_p$
and thus
$A(p)(s_p) = A(s)(p) \in F_\Null\at{p}$.
We also need to show that $\nabla^{\ConHom}_X\cc{A\at{C}} = 0$
for all $X \in \Secinfty(D)$.
For this let $p \in C$ and $\cc{s_p} \in (E_\Wobs / E_\Null)\at{p}$
be given.
Again by \autoref{lem:ExtendingToConSections} we find a section $s \in \ConSecinfty(E)_\Wobs$
such that $\cc{s}(p) = \cc{s_p}$.
Then
\begin{align*}
	\big(\nabla^{\ConHom}_X\cc{A\at{C}}\big)(p)(\cc{s_p})
	&= \big(\nabla^\ConHom_X \cc{A\at{C}}\big)(\cc{s\at{C}})(p) \\
	&= \nabla^F_X\big(\cc{A\at{C}}(\cc{s\at{C}})\big)(p) - \cc{A\at{C}}\big(\nabla^E_X \cc{s\at{C}}\big)\\
	&= \nabla^F_X\big(\cc{A(s)\at{C}}\big)(p)
	= 0,
\end{align*}
since $\nabla^E_X\cc{s\at{C}} = 0$
and $A(s)\at{C} \in \ConSecinfty(F)_\Wobs$.
This shows $A \in \ConSecinfty(\ConHom(E,F))_\Wobs$, and hence
$\eta$ is surjective on the $\WOBS$-component.

Recall from \autoref{lem:indConIsos} that we additionally have to check that $\eta$ is also surjective on the $\NULL$-component.
Thus let $A \in \ConHom_{\ConCinfty(\mathcal{M})}(\ConSecinfty(E),\ConSecinfty(F))_\Null$ be given.
Then for $s \in C$ and $s_p \in E_\Wobs\at{p}$ we find again by \autoref{lem:ExtendingToConSections}
a section $s \in \ConSecinfty(E)_\Wobs$ with $s(p) = s_p$.
Then
$A(p)(s_p) = A(s)(p) \in F_\Null\at{p}$
since $A(s) \in \ConSecinfty(F)_\Null$.
This finally shows that $\eta$ is a regular epimorphism and hence a constraint isomorphism.
\end{proof}

\begin{remark}
Note that we used \autoref{lem:ExtendingToConSections} to prove \autoref{prop:ConSectionsHom}.
Hence by \autoref{rem:NonSimpleGlobalSections} it is not clear if the \autoref{prop:ConSectionsHom}
remains valid for non-simple distributions.
\end{remark}

\begin{corollary}
	\index{dual!constraint vector bundle}
	Let $E = (E_\Total,E_\Wobs,E_\Null,\nabla^E)$ be a constraint vector bundle over a constraint
	manifold $\mathcal{M} = (M,C,D)$.
	Then 
	\begin{equation}
		\ConSecinfty(E^*) \simeq \ConSecinfty(E)^*
	\end{equation}
	as strong constraint $\ConCinfty(\mathcal{M})$-modules.
\end{corollary}

\begin{proof}
	Choose $F = \mathcal{M} \times \Reals$ in \autoref{prop:ConSectionsHom}.
\end{proof}

In classical differential geometry the famous Serre-Swan Theorem states that the category
\glsadd{VectM}$\Vect(M)$ of vector bundles over a fixed manifold $M$ is equivalent to the category
$\Proj(\Cinfty(M))$ of finitely generated projective $\Cinfty(M)$-modules.
By \autoref{prop:FunctorConSections} we know that sections of constraint vector bundles form
strong constraint modules over the strong constraint algebra $\ConCinfty(\mathcal{M})$
of functions on the constraint manifold $\mathcal{M}$.
Thus for a constraint analogue of the Serre-Swan Theorem we expect projective strong constraint modules to be the correct algebraic notion.

Before tackling the full Serre-Swan Theorem, let us take a look at the case of free strong constraint $\ConCinfty(\mathcal{M})$-modules.
As in classical differential geometry these relate to trivial vector bundles, now in the sense of
\autoref{ex:ConVectBundles}.
Recall from \autoref{lem:AdaptedLocalFrames} that every constraint vector bundle admits local frames
adapted to the constraint structure.

\begin{proposition}
	\label{prop:SectionsTrivStrConVectorBundle}
	\index{trivial constraint vector bundle}
	\index{free strong constraint modules}
Let $\mathcal{M} = (M,C,D)$ be a constraint manifold and let
$E = (E_\Total,E_\Wobs,E_\Null,\nabla)$ be a constraint vector bundle over $\mathcal{M}$
of rank $k = (k_\Total,k_\Wobs,k_\Null)$.
Then the following statements are equivalent:
\begin{propositionlist}
	\item \label{prop:SectionsTrivStrConVectorBundle_1}
	The constraint vector bundle $E$ is trivializable.
	\item \label{prop:SectionsTrivStrConVectorBundle_2} There exists 
	a global frame of $E$.
	\item \label{prop:SectionsTrivStrConVectorBundle_3}
	The strong constraint module
	$\ConSecinfty(E)$ is free
	and
	$\ConSecinfty(E) \simeq \ConCinfty(\mathcal{M})^k$.
\end{propositionlist}	
\end{proposition}

\begin{proof}
	\ref{prop:SectionsTrivStrConVectorBundle_1}\,$\Rightarrow$\,\ref{prop:SectionsTrivStrConVectorBundle_2}:
	If $E$ is trivializable there exists a constraint vector bundle isomorphism
	$\Phi \colon E \to \mathcal{M} \times \mathbb{R}^k$
	inducing an isomorphism
	\begin{equation*}
		\Phi \colon \ConSecinfty(E) \to \ConSecinfty(\mathcal{M} \times \mathbb{R}^k)
	\end{equation*} on sections.
	Let $f_1, \dotsc, f_{k_\Total} \in \Secinfty(M \times 
	\mathbb{R}^{k_\Total})$ be the canonical global frame.
	Then $e_i \coloneqq \Phi^{-1}(f_i)$ is a global frame for
	$E_\Total$, such that $e_1, \dotsc, e_{k_\Wobs} \in \ConSecinfty(E)_\Wobs$
	and $e_1,\dotsc, e_{k_\Null} \in \ConSecinfty(E)_\Null$.
	Moreover, since $\iota^\#f_1, \dotsc, \iota^\#f_{k_\Wobs}$
	and $\iota^\#f_1, \dotsc,\iota^\#f_{k_\Null}$ form global frames 
	for the trivial vector bundles $C \times \mathbb{R}^{k_\Wobs}$
	and $C \times \mathbb{R}^{k_\Null}$, respectively, and since
	$\Phi$ induces isomorphisms on $E_\Wobs$ and $E_\Null$, we see 
	that
	$\iota^\#e_1, \dotsc, \iota^\#e_{k_\Wobs}$
	and
	$\iota^\#e_1,\dotsc, \iota^\#e_{k_\Null}$ form global frames for
	$E_\Wobs$ and $E_\Null$, respectively.
	
	\ref{prop:SectionsTrivStrConVectorBundle_2}\,$\Rightarrow$\,\ref{prop:SectionsTrivStrConVectorBundle_3}:
	Every $s \in \ConSecinfty(E)_\Total$ 
	can be written as $s = \sum_{i=1}^{k_\Total} s_i e_i$ with
	$s_i \in \Cinfty(M)$.
	If $s \in \ConSecinfty(E)_\Wobs$, then
	from
	\begin{equation*}
		\iota^\#s = \sum_{i=1}^{k_\Total} \iota^*s_i \cdot 
		\iota^\#e_i \in \Secinfty(E_\Wobs)
	\end{equation*}
	it follows that $s_{k_\Wobs +1}, \dotsc, s_{k_\Total} \in \ConCinfty(\mathcal{M})_\Null$.
	Moreover, since
	\begin{equation*}
		0 = \Lie_X \cc{s} = \sum_{i=k_\Null +1}^{k_\Wobs} (\Lie_X \iota^*s_i) \cdot e_i
	\end{equation*}
	for all $X \in \Secinfty(D)$ we get
	$s_{k_\Null +1},\dotsc, s_{k_\Wobs} \in \ConCinfty(\mathcal{M})_\Wobs$.
	And thus $\ConSecinfty(E)_\Wobs \simeq 
	(\ConCinfty(\mathcal{M})^k)_\Wobs$.
	If $s \in \ConSecinfty(E)_\Null$ we have
	\begin{equation*}
		\iota^\#s = \sum_{i=1}^{k_\Total} \iota^*s_i \cdot 
		\iota^\#e_i \in \Secinfty(E_\Null)
	\end{equation*}
	and thus
	$s_{k_\Null + 1}, \dotsc, s_{k_\Total} \in \ConCinfty(\mathcal{M})_\Null$,
	giving $\ConSecinfty(E)_\Null \simeq 
	(\ConCinfty(\mathcal{M})^k)_\Null$.
	Together this yields
	$\ConSecinfty(E) \simeq \ConCinfty(\mathcal{M})^k$.
	
	\ref{prop:SectionsTrivStrConVectorBundle_3}\,$\Rightarrow$\,\ref{prop:SectionsTrivStrConVectorBundle_1}:
	Suppose we have an isomorphism
	$\Phi \colon \ConSecinfty(E) \to \ConCinfty(\mathcal{M})^k$.
	From classical differential geometry we know that $E_\Total \simeq M \times \mathbb{R}^{k_\Total}$
	by mapping $v_p \in E\at{p}$ to $\Psi(v_p) \coloneqq \Phi(s)(p)$ for any $s  \in \Secinfty(E_\Total)$
	with $s(p) = v_p$,
	and $\Psi$ does not depend on the choice of $s$.
	We need to check that $\Psi$ is an isomorphism of constraint vector bundles.
	For this let $p \in C$ and $v_p \in E_\Wobs\at{p}$ be given.
	By \autoref{lem:ExtendingToConSections} \ref{lem:ExtendingToConSections_2}
	there exists $s \in \ConSecinfty(E)_\Wobs$ with $s(p) = v_p$.
	Hence $\Psi(v_p) \in \Phi(s)(p) \in (C \times \mathbb{R}^{k_\Wobs})\at{p}$
	since $\Phi(s) \in \ConSecinfty(M \times \mathbb{R}^k)_\Wobs$.
	Similarly, if $v_p \in E_\Null\at{p}$ then by
	\autoref{lem:ExtendingToConSections} \ref{lem:ExtendingToConSections_1}
	there exists $s \in \ConSecinfty(E)_\Null$ such that $s(p) = v_p$.
	Then $\Psi(v_p) = \Phi(s)(p) \in (C \times \mathbb{R}^{n_\Wobs})\at{p}$,
	since $\Phi(s) \in \ConSecinfty(M\times \mathbb{R}^{n_\Total})_\Null$.
	The same arguments show that $\Psi^{-1} \colon (M \times \mathbb{R}^{k_\Total}) \to E_\Total$
	preserves the $\WOBS$- and $\NULL$-components, hence inducing isomorphisms
	$\Psi_\Wobs \colon E_\Wobs \to (C \times \mathbb{R}^{k_\Wobs})$
	and
	$\Psi_\Null \colon E_\Null \to (C \times \mathbb{R}^{k_\Null})$.
	To show that $\Psi$ is compatible with the covariant derivatives note that it induces also an isomorphism
	$\Psi_{\Wobs / \Null} \colon (E_\Wobs / E_\Null) \to (C \times \mathbb{R}^{k_\Wobs - k_\Null})$.
	Then for $s \in \Secinfty(E_\Wobs / E_\Null)$ we have
	\begin{align*}
		\Psi \big( \nabla^E_{v_p} s\big)
		= \Psi \Big(\sum_{i=n_0+1}^{k_\Wobs - k_\Null}  \nabla^E_{v_p}  (s^i e_i) \Big)
		= \Psi \Big( \sum_{i=n_0+1}^{k_\Wobs - k_\Null} (\Lie_Xs^i) e_i \Big)
		= \sum_{i=n_0+1}^{k_\Wobs - k_\Null} (\Lie_Xs^i) \Psi(e_i)
	\end{align*}
	for all $X \in \Secinfty(D)$, showing that $\Psi$ is indeed an isomorphism of constraint vector bundles.
\end{proof}

\begin{remark}
We again used \autoref{lem:ExtendingToConSections} in the above proof.
Hence by \autoref{rem:NonSimpleGlobalSections} it is not clear if the the above equivalences still hold
for non-simple distributions.
\end{remark}

The existence of local frames for constraint vector bundles can therefore be understood
as local freeness of $\ConSecinfty(E)$.

As a first step towards the constraint Serre-Swan Theorem we show that every finitely generated projective strong constraint module over the constraint algebra of functions can be realized as sections of a constraint vector bundle.

\begin{proposition}
	\label{prop:CSecIsEssSurj}
Let $\mathcal{M} = (M,C,D)$ be a constraint manifold
and $\module{P} \in \strConProj(\ConCinfty(\mathcal{M}))$
a finitely generated projective strong constraint 
$\ConCinfty(\mathcal{M})$-module.
Then there exists a constraint vector bundle
$E = (E_\Total,E_\Wobs,E_\Null,\nabla)$ over $\mathcal{M}$
such that $\ConSecinfty(E) \simeq \module{P}$.
\end{proposition}

\begin{proof}
	Since $\module{P}$ is finitely generated projective there exists a finite constraint index set
	$n \in \injConIndSet$ and a projection
	$e \in \ConEnd_\algebra{A}(\algebra{A}^{(M)})$ with $e^2 = e$
	such that $\module{P} \simeq e\ConCinfty(\mathcal{M})^n$.
	By \autoref{prop:ConSectionsHom} the projection
	$e$ can be viewed as a constraint section of $\ConEnd(\mathcal{M} \times \mathbb{R}^n)$.
	Moreover, since $e$ is completely determined by its $\TOTAL$-component
	we can identify it with a matrix
	$e \in \Mat_{n_\Total}(\Cinfty(M))$.
	This leads to a vector bundle morphism
	\begin{align*}
		\Phi_\Total \colon
		M \times \field{R}^{n_\Total} \to M \times \field{R}^{n_\Total},
		\quad
		(p,v)
		&\mapsto
		(p, e(p) v)
	\end{align*}
	of constant rank.
	And therefore we can define $E_\Total \coloneqq \image(e_\Total)$
	as a subbundle of $M \times \mathbb{R}^{n_\Total}$.
	Since $e \in \ConSecinfty(\End(\mathcal{M} \times \mathbb{R}^n))_\Wobs$
	we know
	$\iota^\#e \in \Secinfty(\End(C \times \mathbb{R}^{n_\Wobs}))$
	leading to a constant rank vector bundle morphism
		\begin{align*}
				\Phi_\Wobs \colon
				C \times \field{R}^{n_\Wobs} \to C \times \field{R}^{n_\Wobs},
				\quad
				(p,v)
				&\mapsto
				(p, \iota^\#e(p) v).
	\end{align*}
	This allows us to define $E_\Wobs \coloneqq \image(e\at{C})$
	as a subbundle of $\iota^\#E_\Total$.
	Moreover, since $\iota^\#e_p$ preserves also the $\NULL$-component of the fibre
	we can restrict $\Phi_\Wobs$ to $C \times \mathbb{R}^{n_0}$,
	giving a subbundle $E_\Null \coloneqq \image(\Phi_\Wobs\at{C\times \mathbb{R}^{n_\Null}})$ of $E_\Wobs$.
	Finally, we can define a partial $D$-connection on $E_\Wobs/E_\Null$
	by
	\begin{equation*}
		\nabla_X \cc{s} \coloneqq \sum_{i=n_0+1}^{n_\Wobs} \Phi_\Wobs(b_i) \cdot \Lie_X s^i 
	\end{equation*}
	for all $s \in \Secinfty(E)_\Wobs$ with
	$s = \sum_{i=1}^{n_\Wobs} \Phi_\Wobs(b_i) s^i$ and
	$X \in \Secinfty(D)$.
	Here the $b_i$ denote the canonical basis sections of $C \times \mathbb{R}^{n_\Wobs}$.
	This clearly gives a well-defined covariant derivative.
	To show that $\nabla$ is path-independent consider
	$s_p = \sum_{i=n_\Null+1}^{n_\Wobs} \Phi_\Wobs(b_i)(p) s_p^i \in (E_\Wobs / E_\Null)\at{p}$.
	Then the section $s = \sum_{i=n_\Null+1}^{n_\Wobs} \Phi_\Wobs(b_i) s_p^i$ is clearly covariantly constant and
	thus induces the parallel transport along any leafwise curve $\gamma \colon I \to C$.
	It remains to show that $\ConSecinfty(E)$ is isomorphic to $\module{P}$
	as a strong constraint $\ConCinfty(\mathcal{M})$-module.
	It is straightforward to check that $\Psi \colon \image e \to \ConSecinfty(E)$
	defined by
	\begin{align*}
		\Psi_\Total (s) \coloneqq (p \mapsto (p,s(p)))
	\end{align*}
	is an isomorphism of constraint modules.
	And hence $\ConSecinfty(E) \simeq \image e \simeq \module{P}$ follows.
\end{proof}

To show that sections of constraint vector bundles are always finitely generated projective we actually need the requirement of a simple distribution:

\begin{proposition}
	\label{prop:ConSectionsAreFGP}
Let $E = (E_\Total,E_\Wobs,E_\Null,\nabla)$ be a constraint vector bundle over a constraint manifold $\mathcal{M} =(M,C,D)$.
Then $\ConSecinfty(E)$ is a finitely generated projective strong constraint $\ConCinfty(\mathcal{M})$-module.
\end{proposition}

\begin{proof}
We construct a dual basis in the sense of \autoref{prop:StrDualBasis}.
For this we first choose a complement $E_\Null^\perp$ of
$E_\Null$ inside $E_\Wobs$, hence we get
$E_\Wobs = E_\Null \oplus E_\Null^\perp$ with $E_\Null^\perp \simeq E_\Wobs /E_\Null$,
and additionally a complement $E_\Wobs^\perp$ of $E_\Wobs$ inside
$\iota^\#E_\Total$.
This yields $\iota^\#E_\Total = E_\Null \oplus E_\Null^\perp \oplus E_\Wobs^\perp$.
Now choose a finite dual basis of $\Secinfty(E_\red)$ given by
$g_j \in \Secinfty(E_\red)$ and $g^j \in \Secinfty(E_\red^*)$, for $j \in J_\Null^\perp$.
By \autoref{prop:ReductionOfConVect} we can pull back the dual basis to a dual basis
of $E_\Null^\perp$, which we still denote by $g_j \in \Secinfty(E_\Null^\perp)$ and
$g^j\in \Secinfty((E_\Null^\perp)^*)$.
Note that these sections fulfil $\nabla_X g_j = 0$ and $\nabla_X^* g^j = 0$ for
$X \in \Secinfty(D)$.
Additionally, choose a dual basis $\{f_j,f^j\}_{j \in J_\Null}$ of $E_\Null$
and a dual basis $\{h_j,h^j\}_{j \in J_\Wobs^\perp}$ of $E_\Wobs^\perp$.
This way we obtain a dual basis
$\{c_j, c^j\}_{i \in J_C}$ of $\iota^\#E_\Total$ with
$J_C = J_\Null \sqcup J_\Null^\perp \sqcup J_\Wobs^\perp$
\begin{equation*}
	c_j = \begin{cases}
		f_j & \text{ if } j \in J_\Null \\
		g_j & \text{ if } j \in J_\Null^\perp \\
		h_j & \text{ if } j \in J_\Wobs^\perp
	\end{cases}\quad \text{ and }\quad
c^j = \begin{cases}
	f^j & \text{ if } j \in J_\Null \\
	g^j & \text{ if } j \in J_\Null^\perp \\
	h^j & \text{ if } j \in J_\Wobs^\perp
\end{cases}.
\end{equation*}
To extend the dual basis to all of $M$ we choose a tubular neighbourhood 
$\pr_V \colon V \to C$, with $\iota_V \colon V \hookrightarrow M$
an open neighbourhood of $C$.
Then we can pull back the $c_j$ and $c^j$ to obtain a dual basis of 
$\iota_V^\#E_\Total$, which we again denote by $\{c_j,c^j\}_{j\in J_C}$.
On the open subset $\iota_{M\setminus C} \colon M\setminus C \hookrightarrow M$
choose another dual basis $\{d_k,d^k\}_{k \in K}$ of $\iota_{M\setminus C}^\#E_\Total$.
We now need to patch these dual bases together.
For this choose a quadratic partition of unity $\chi_1, \chi_2 \in \Cinfty(M)$
with $\chi_1^2 + \chi_2^2 = 1$ and $\supp \chi_1 \subseteq V$ and
$\supp \chi_2 \subseteq M \setminus C$.
Then $\{e_i, e^i\}_{i \in I_\Total}$ with $I_\Total = J_C \sqcup K$ defined by
\begin{equation*}
	e_i = \begin{cases}
		\chi_1\cdot c_i & \text{ if } i \in J_C \\
		\chi_2 \cdot d_i & \text{ if } i \in K
	\end{cases}
	\quad \text{ and }\quad
	e^i = \begin{cases}
		\chi_1 \cdot c^i & \text{ if } i \in J_C \\
		\chi_2 \cdot d^i & \text{ if } i \in K
	\end{cases}
\end{equation*}
forms a dual basis for $E_\Total$.
It remains to show that this dual basis fulfils the properties of
\autoref{prop:StrDualBasis}.
For this consider the constraint set $I$ with $I_\Total$ as above,
$I_\Wobs = J_\Null \sqcup J_\Null^\perp \sqcup K$ and $I_\Null = J_\Null \sqcup K$.
By construction we have $e_i \in \ConSecinfty(E)_\Wobs$ for $i \in I_\Wobs$
and $e_i \in \ConSecinfty(E)_\Null$ for $i \in I_\Null$.
From the fact that $g^j \in \Secinfty(\Ann E_\Null)$ and $h^j \in \Secinfty(\Ann E_\Wobs)$
it follows that $e^i \in \ConSecinfty(E^*)_\Wobs$ for $i \in I_\Total \setminus I_\Null$
and $e^i \in \ConSecinfty(E^*)_\Null$ for $i \in I_\Total \setminus I_\Wobs$.	
\end{proof}

\begin{remark}
In the above proof we heavily used the existence of a reduced vector bundle on a smooth reduced manifold.
Thus it is not clear if the above statement still holds for non-simple distributions.
Nevertheless, this situation is of great interest for its geometric applications,
hence the question if all modules of sections are projective, even for non-simply distributions, deserves further attention.
\end{remark}

The above results lead us now to a constraint version of the Serre-Swan Theorem:

\begin{theorem}[Constraint Serre-Swan]
	\label{thm:strConSerreSwan}
	\index{Serre-Swan Theorem}
Let $\mathcal{M} = (M,C,D)$ be a constraint manifold.
The functor $\ConSecinfty \colon \ConVect(\mathcal{M}) \to \strConProj({\ConCinfty(\mathcal{M})})$
is an equivalence of categories.
\end{theorem}

\begin{proof}
\autoref{prop:ConSectionsAreFGP} shows that $\ConCinfty$ actually 
maps to $\strConProj(\ConCinfty(\mathcal{M}))$, while
\autoref{prop:ConSectionsHom} proves that $\ConCinfty$
is fully faithful.
Finally, by \autoref{prop:CSecIsEssSurj}
it is essentially surjective, and therefore an equivalence of 
categories.
\end{proof}

\begin{remark}
In \cite{dippell.menke.waldmann:2022a, menke:2020a} a similar result for
non-strong projective constraint modules over
$\ConCinfty(\mathcal{M})$ as a non-strong constraint algebra was found.
The geometric objects used there are similar but not identical to the notion of constraint vector bundles, in particular the vector bundle $E_\Wobs$ is a subbundle of $E_\Total$ defined on all of $M$, and $\nabla$ is a partial connection on $\iota^\#E_\Wobs$ instead of $E_\Wobs / E_\Null$.
\end{remark}

The constraint Serre-Swan Theorem finally justifies the study of projective strong constraint modules, and their predecessors in \autoref{sec:RegularProjectiveModules}.
This important result allows us now to examine the compatibility of the sections functor with the different notions of tensor products.
Consider vector bundles $E$ and $F$ over a constraint manifold $\mathcal{M}$.
By the Serre-Swan Theorem we know that
$\ConSecinfty(E)$ and $\ConSecinfty(F)$ are finitely generated projective strong constraint modules,
moreover, \autoref{prop:TensorStrConProj} 
and \autoref{prop:StrTensorStrConProj}
tell us that also 
$\ConSecinfty(E) \injstrtensor[\ConCinfty(\mathcal{M})] \ConSecinfty(F)$
and
$\ConSecinfty(E) \strcirctensor[\ConCinfty(\mathcal{M})] \ConSecinfty(F)$
are finitely generated projective
and hence embedded.
This is something we cannot expect for finitely generated projective modules over arbitrary embedded strong constraint algebras,
since their (strong) tensor products need in general not be embedded.
From now on we will write 
$\ConSecinfty(E) \strtensor[\ConCinfty(\mathcal{M})] \ConSecinfty(F)$
and
$\ConSecinfty(E) \tensor[\ConCinfty(\mathcal{M})] \ConSecinfty(F)$
instead, since on $\strConProj(\ConCinfty(\mathcal{M}))$ there will be no risk of confusion.

\begin{lemma}
	\label{lem:ConSecinftyLaxMonoidal}
	Let $\mathcal{M} = (M,C,D)$ be a constraint manifold and
	let $E, F \in \ConVect(\mathcal{M})$ be constraint vector bundles over $\mathcal{M}$.
	\begin{lemmalist}
		\item 
		\label{lem:ConSecinftyLaxMonoidal_1}
		Setting 
		\begin{align}
			\label{eq:ConSecinftyLaxMonoidal_1}
			I_{E,F}(s \tensor t)(p) &\coloneqq s(p) \tensor t(p)
		\end{align}
		defines a constraint morphism
		$I_{E,F} \colon \ConSecinfty(E) \tensor[\ConCinfty(\mathcal{M})] \ConSecinfty(F)
		\to \ConSecinfty(E \tensor F)$.
		These morphisms constitute a natural transformation
		$I \colon \tensor[\ConCinfty(\mathcal{M})] \circ\, (\ConSecinfty \times \ConSecinfty)
		\Rightarrow \ConSecinfty \circ \tensor$.
		\item \label{lem:ConSecinftyLaxMonoidal_2}
		Setting
		\begin{align}
			\label{eq:ConSecinftyLaxMonoidal_2}
			J_{E,F}(s \tensor t)(p) &\coloneqq s(p) \tensor t(p)
		\end{align}
		defines a constraint morphism
		$J_{E,F} \colon \ConSecinfty(E) \strtensor[\ConCinfty(\mathcal{M})] \ConSecinfty(F)
		\to \ConSecinfty(E \strtensor F)$.
		These morphisms constitute a natural transformation
		$J \colon \strtensor[\ConCinfty(\mathcal{M})] \circ\, (\ConSecinfty \times \ConSecinfty)
		\Rightarrow \ConSecinfty \circ \strtensor$.
	\end{lemmalist}
\end{lemma}

\begin{proof}
In both cases we know from classical differential geometry that the $I_{E,F}$ 
are a morphism of $\algebra{A}_\Total$-modules on the $\TOTAL$-components,
forming natural transformations.
It remains to show that $I$ and $J$ are constraint morphisms, meaning that
they preserve the $\WOBS$- and $\NULL$-components.

For the first part let
$s \tensor t \in \ConSecinfty(E) \ConGrid[2][2][0][2] \ConSecinfty(F)
= \big(\ConSecinfty(E) \tensor[\ConCinfty(\mathcal{M})] \ConSecinfty(F)\big)_\Null$.
Then
\begin{equation*}
	I(s \tensor t)(p)
	= s(p) \tensor t(p) 
	\in E_\Null\at{p} \tensor F_\Wobs\at{p}
	+ E_\Wobs\at{p} \tensor F_\Null\at{p}
\end{equation*}
for all $p \in C$, and therefore
$I(s \tensor t) \in \ConSecinfty(E \tensor F)_\Null$.
Now consider 
$s \tensor t \in \ConSecinfty(E) \ConGrid[0][0][0][0][2] \ConSecinfty(F)
\subseteq (\ConSecinfty(E) \tensor[\ConCinfty(\mathcal{M})] \ConSecinfty(F))_\Wobs$, then
\begin{equation*}
	I(s \tensor t)(p)
	= s(p) \tensor t(p) 
	\in E_\Wobs\at{p} \tensor F_\Wobs\at{p}
\end{equation*}
for all $p \in C$, hence $I(s \tensor t)\at{C} \in \Secinfty(\big(E \tensor F)_\Wobs\big)$.
Moreover, for $X \in \Secinfty(D)$ we have
\begin{equation*}
	\nabla_X \cc{I(s \tensor t)\at{C}}
	= \nabla_X \cc{s\at{C} \tensor t\at{C}}
	= \nabla_X \cc{s\at{C}} \tensor \cc{t\at{C}}
	+ \cc{s\at{C}} \tensor \nabla_X \cc{t\at{C}}
	= 0,
\end{equation*}
showing that $I(s \tensor t) \in \ConSecinfty(E \tensor F)_\Wobs$.

For the second part we start with
$s \tensor t \in \ConSecinfty(E) \ConGrid[2][2][2][2][0][0][2] \ConSecinfty(F)
= \big(\ConSecinfty(E) \strtensor[\ConCinfty(\mathcal{M})] \ConSecinfty(F)\big)_\Null$.
Then
\begin{equation*}
	J(s \tensor t)(p)
	= s(p) \tensor t(p) 
	\in E_\Null\at{p} \tensor \iota^\#F_\Total\at{p}
	+ \iota^\#E_\Total\at{p} \tensor F_\Null\at{p}
\end{equation*}
for all $p \in C$, and therefore
$J(s \tensor t) \in \ConSecinfty(E \strtensor F)_\Null$.
Now consider 
$s \tensor t \in \ConSecinfty(E) \ConGrid[0][0][0][0][2] \ConSecinfty(F)
\subseteq \big(\ConSecinfty(E) \strtensor[\ConCinfty(\mathcal{M})] \ConSecinfty(F)\big)_\Wobs$, then
\begin{equation*}
	J(s \tensor t)(p)
	= s(p) \tensor t(p) 
	\in E_\Wobs\at{p} \tensor F_\Wobs\at{p}
\end{equation*}
for all $p \in C$, hence $J(s \tensor t)\at{C} \in \Secinfty\big((E \strtensor F)_\Wobs\big)$.
Moreover, for $X \in \Secinfty(D)$ we have
\begin{equation*}
	\nabla_X \cc{J(s \tensor t)\at{C}}
	= \nabla_X \cc{s\at{C} \tensor t\at{C}}
	= \nabla_X \cc{s\at{C}} \tensor \cc{t\at{C}}
	+ \cc{s\at{C}} \tensor \nabla_X \cc{t\at{C}}
	= 0,
\end{equation*}
showing that $J(s \tensor t) \in \ConSecinfty(E \strtensor F)_\Wobs$.
\end{proof}

The canonical morphisms \ref{eq:ConSecinftyLaxMonoidal_1}
and \ref{eq:ConSecinftyLaxMonoidal_2} can be constructed without using
the Serre-Swan Theorem.
But to see that these are in fact isomorphisms we need constraint dual bases.

\begin{proposition}
	\index{tensor product!constraint vector bundles}
	\index{strong tensor product!constraint vector bundles}
Let $\mathcal{M} = (M,C,D)$ be a constraint manifold.
\begin{propositionlist}
	\item The sections functor
	$\ConSecinfty \colon (\ConVect(\mathcal{M}), \tensor)
	\to \big(\strConProj\big(\ConCinfty(\mathcal{M})\big), \tensor[\ConCinfty(\mathcal{M})]\big)$
	is monoidal.
	\item The sections functor
	$\ConSecinfty \colon (\ConVect(\mathcal{M}), \strtensor)
	\to \big(\strConProj\big(\ConCinfty(\mathcal{M})\big), \strtensor[\ConCinfty(\mathcal{M})]\big)$
	is monoidal.
\end{propositionlist}
\end{proposition}

\begin{proof}
We first show that the natural transformations 
$I$ and $J$ from \autoref{lem:ConSecinftyLaxMonoidal}
are in fact natural isomorphisms.
For this we construct inverses.
Let $E, F \in \ConVect(\mathcal{M})$ and let
$(\{e_i\}_{i \in M}, \{e^i\}_{i \in M^*})$ as well as
$(\{f_j\}_{j \in N}, \{f^j\}_{j \in N^*})$ be finite dual bases of $E$ and $F$, respectively.
From classical differential geometry we know that
$(\{e_i \tensor f_j\}_{(i,j) \in M_\Total \times N_\Total}, \{e^i \tensor f^j\}_{(i,j) \in M_\Total \times N_\Total})$
is a dual basis of $\Secinfty(E_\Total \tensor F_\Total)$
and that
\begin{equation*}
	K(X) = \sum_{i \in M_\Total} \sum_{j \in N_\Total}
	(e^i \tensor f^j)(X) \cdot e_i \tensor[\Cinfty(M)] f_j,
\end{equation*}
for $X \in \Secinfty(E_\Total \tensor F_\Total)$,
defines an inverse 
$K \colon \Secinfty(E \tensor F) \to \Secinfty(E) \tensor[\Cinfty(M)] \Secinfty(F)$
to $I$.
Here we use $\tensor[\Cinfty(M)]$ as the algebraic tensor product to separate it from the
geometric tensor product $\tensor$ of sections.
To show that $K$ is a constraint morphism we prove that 
the families
$(\{e_i \tensor f_j\}_{(i,j) \in M \tensor N}, \{e^i \tensor f^j\}_{(i,j) \in (M \tensor N)^*})$
form a dual basis for $\ConSecinfty(E \tensor F)$:
\begin{cptitem}
	\item $(i,j) \in (M \tensor N)_\Wobs = M \ConGrid[2][2][0][2][2] N$:
	Then by \autoref{lem:ConSecinftyLaxMonoidal} we know that
	\begin{equation*}
		e_i \tensor f_j = I_{E,F}(e_i \tensor[\Cinfty(M)] f_j) \in \ConSecinfty(E \tensor F)_\Wobs.
	\end{equation*}
	\item $(i,j) \in (M \tensor N)_\Null = M \ConGrid[2][2][0][2] N$:
	Then we know that
	$e_i \tensor f_j = I_{E,F}(e_i \tensor[\Cinfty(M)] f_j) \in \ConSecinfty(E \tensor F)_\Null$.
	\item $(i,j) \in (M \tensor N)^*_\Wobs = M \ConGrid[0][0][2][0][2][2][2][2][2] N = (M^* \strtensor N^*)_\Wobs$:
	Then we know that
	\begin{equation*}
		e^i \tensor f^j = J_{E^*,F^*}(e^i \tensor[\Cinfty(M)] f^j) \in \ConSecinfty(E^* \strtensor F^*)_\Wobs \simeq \ConSecinfty(E \tensor F)^*_\Wobs.
	\end{equation*}
	\item $(i,j) \in (M \tensor N)^*_\Null = M \ConGrid[0][0][2][0][0][2][2][2][2] N$:
	Then we know that
	\begin{equation*}
		e^i \tensor f^j = J_{E^*,F^*}(e^i \tensor[\Cinfty(M)] f^j) \in \ConSecinfty(E^* \strtensor F^*)_\Null
		\simeq \ConSecinfty(E \tensor F)^*_\Null.
	\end{equation*}
\end{cptitem}
This shows that $K$ is a constraint morphism, and therefore $I$ is an isomorphism.
With completely analogous arguments, on can show that $J$ is an isomorphism as well.
The unit object in $\ConVect(\mathcal{M})$ is for both products given by 
$\mathcal{M} \times \Reals$.
Since $\ConSecinfty(\mathcal{M} \times \Reals) \simeq \ConCinfty(\mathcal{M})$
the section functor preserves the monoidal units, and hence gives a monoidal functor in both cases.
\end{proof}

Since $\ConSecinfty$ is monoidal and compatible with direct sums, we also get
\begin{align}
	\Sym^\bullet_{\tensor} \ConSecinfty(E)
	\simeq \ConSecinfty(\Sym^\bullet_{\tensor} E),
	\qquad
	\Anti^\bullet_{\tensor} \ConSecinfty(E)
	\simeq \ConSecinfty(\Anti^\bullet_{\tensor} E),
\shortintertext{as well as}
	\Sym^\bullet_{\strtensor} \ConSecinfty(E)
	\simeq \ConSecinfty(\Sym^\bullet_{\strtensor} E),
	\qquad
	\Anti^\bullet_{\strtensor} \ConSecinfty(E)
	\simeq \ConSecinfty(\Anti^\bullet_{\strtensor} E)
\end{align}
for any constraint vector bundle $E$.
For sections of constraint vector bundles we can make the relation of the strong and non-strong tensor products precise.
In particular, their difference will be located on the submanifold $C$ only.

\begin{proposition}
	\index{tensor product!constraint vector bundles}
	\index{strong tensor product!constraint vector bundles}
Let $\mathcal{M} = (M,C,D)$ be a constraint manifold and let
$E,F \in \ConVect(\mathcal{M})$ be constraint vector bundles
over $\mathcal{M}$.
Then there exists an isomorphism of constraint 
$\ConCinfty(\mathcal{M})$-modules such that
\begin{equation} \label{eq:StrTensorForSections}
\begin{split}
	\left(\ConSecinfty(E) \strtensor[\ConCinfty(\mathcal{M})] \ConSecinfty(F) \right)_\Total
	&\simeq \left( \ConSecinfty(E)  \tensor[\ConCinfty(\mathcal{M})] \ConSecinfty(F) \right)_\Total,\\
	\left( \ConSecinfty(E) \strtensor[\ConCinfty(\mathcal{M})] \ConSecinfty(F) \right)_\Wobs
	&\simeq \left( \ConSecinfty(E) \tensor[\ConCinfty(\mathcal{M})] \ConSecinfty(F) \right)_\Wobs\\
	&\quad \oplus \Secinfty(E_\Null) \tensor[\Cinfty(C)] \Secinfty(\iota^\#F_\Total / F_\Wobs)\\
	&\quad \oplus \Secinfty(\iota^\#E_\Total / E_\Wobs) \tensor[\Cinfty(C)] \Secinfty(F_\Null), \\
	\left( \ConSecinfty(E) \strtensor[\ConCinfty(\mathcal{M})] \ConSecinfty(F) \right)_\Null 
	& \simeq \left( \ConSecinfty(E) \tensor[\ConCinfty(\mathcal{M})] \ConSecinfty(F) \right)_\Null\\
	&\quad \oplus \Secinfty(E_\Null) \tensor[\Cinfty(C)] \Secinfty(\iota^\#F_\Total / F_\Wobs)\\
	&\quad \oplus \Secinfty(\iota^\#E_\Total / E_\Wobs) \tensor[\Cinfty(C)] \Secinfty(F_\Null).
\end{split}
\end{equation}
\end{proposition}

\begin{proof}
Choose complementary vector bundles $E_\Null^\perp$ and $E_\Wobs^\perp$ over $C$ such that
$\iota^\#E_\Total = E_\Wobs \oplus E_\Wobs^\perp = E_\Null \oplus E_\Null^\perp \oplus E_\Wobs^\perp$.
In particular we have $E_\Wobs^\perp \simeq E_\Wobs / E_\Null$ and $E_\Total^\perp \simeq \iota^\#E_\Total / E_\Wobs$.
Similarly, choose complementary vector bundles $F_\Null^\perp$ and $F_\Wobs^\perp$.
Additionally, we need a tubular neighbourhood $\pr_V \colon V \to C$, with
$\iota_V \colon V \hookrightarrow M$ an open neighbourhood of $C$.
Using this we can extend the vector bundles $E_\Null$, $E_\Total^\perp$, $F_\Null$ and $F_\Total^\perp$
to $V$ by pulling them back along $\pr_V$.
Finally, we need a bump function $\chi$ such that $\chi\at{C} = 1$ and $\chi\at{M\setminus V} = 0$.
With this we can turn the right hand side of \eqref{eq:StrTensorForSections} 
into a constraint module by defining $\iota$
on $\Secinfty(E_\Null) \tensor[\Cinfty(C)] \Secinfty(\iota^\#F_\Total / F_\Wobs) \oplus \Secinfty(\iota^\#E_\Total / E_\Wobs) \tensor[\Cinfty(C)] \Secinfty(F_\Null)$
as 
\begin{equation*}
	\iota(t_1 \tensor t_2) \coloneqq \chi \cdot \pr_V^\#(t_1 \tensor t_2).
\end{equation*}	
Now let
$s = \sum_{i \in I} s_1^i \tensor s_2^i \in (\ConSecinfty(E)
\strtensor[\ConCinfty(\mathcal{M})] \ConSecinfty(F) )_\Wobs
= \ConSecinfty(E \strtensor F)_\Wobs$ be given.
Since
$(E \strtensor F)_\Wobs \simeq (E_\Wobs \tensor F_\Wobs)
\oplus (E_\Null \tensor F_\Total^\perp)
\oplus (E_\Total^\perp \tensor F_\Null)$
we can write $I$ as $I = I_{\Wobs\Wobs} \sqcup I_{\Null\Total} \sqcup I_{\Total\Null}$
with
\begin{align*}
	i \in I_{\Wobs\Wobs} &\iff s_1^i\at{C} \tensor s_2^{i}\at{C} \in \Secinfty(E_\Wobs) \tensor[\Cinfty(C)] \Secinfty(F_\Wobs), \\
	i \in I_{\Null\Total} &\iff s_1^i\at{C} \tensor s_2^{i}\at{C} \in \Secinfty(E_\Null) \tensor[\Cinfty(C)] \Secinfty(F_\Total^\perp), \\
	i \in I_{\Total\Null} &\iff s_1^i\at{C} \tensor s_2^{i}\at{C} \in \Secinfty(E_\Total^\perp) \tensor[\Cinfty(C)] \Secinfty(F_\Null).
\end{align*}
Extending the sections first to $V$ by pullback and then to $M$ by use of $\chi$ yields
\begin{equation*}
	s - \sum_{i \in I_{\Null\Total} \sqcup I_{\Total\Null}}\iota\big(s_1^i\at{C} \tensor s_2^i\at{C}\big)
	\in \big(\ConSecinfty(E) \tensor[\ConCinfty(\mathcal{M})] \ConSecinfty(F)\big)_\Wobs.
\end{equation*}
Thus we can define
\begin{equation*}
	\Psi_\Wobs(s)
	\coloneqq \Big(s - \sum_{i \in I_{\Null\Total} \sqcup I_{\Total\Null}}\iota(s_1^i\at{C} \tensor s_2^i\at{C}),\,\,
	\sum_{i \in I_{\Null\Total}} s_1^i\at{C} \tensor s_2^i\at{C},\,\,
	\sum_{i \in I_{\Total\Null}} s_1^i\at{C} \tensor s_2^i\at{C} \Big).
\end{equation*}
It is then easy to see that $\Psi_\Wobs$ preserves the $\NULL$-component and together with 
the canonical isomorphism $\Psi_\Total$ on the $\TOTAL$-component defines a constraint
module morphism.
Moreover, the inverse of $\Psi_\Wobs$ is given by
\begin{equation*}
	\Psi^{-1}_\Wobs(s,t,u) = s + \iota(t) + \iota(u).
\end{equation*}
Thus we get indeed an isomorphism as required.
\end{proof}

With \eqref{eq:StrTensorForSections} it becomes obvious that the tensor products $\tensor$ and $\strtensor$ indeed differ the moment that
$E_\Null$ is non-trivial and $F_\Wobs \subsetneq \iota^\#F_\Wobs$ is an honest subbundle, since then
\begin{equation}
	\Secinfty(E_\Null) \tensor[\Cinfty(C)] \Secinfty(\iota^\#F_\Total / F_\Wobs) \simeq \Secinfty(E_\Null \tensor \iota^\#F_\Total/F_\Wobs)
\end{equation}
does not vanish.

\renewcommand{\thesubsubsection}{\thesection.\arabic{subsubsection}}
\subsubsection{Reduction}

As closure of this section we can show that the constraint Serre-Swan Theorem reduces to the classical Serre-Swan Theorem.
More precisely, taking sections commutes with reduction as shown in the following:

\begin{proposition}[Constraint sections vs. reduction]
	\label{prop:ConSecVSReduction}
	\index{reduction!sections}
	\index{constraint!sections}
Let $\mathcal{M} = (M,C,D)$ be a constraint manifold.
There exists a natural isomorphism making the following diagram commute:
\begin{equation}
\begin{tikzcd}
	\ConVect(\mathcal{M})
		\arrow[r,"\ConSecinfty"]
		\arrow[d,"\red"{swap}]
	& \strConProj(\ConCinfty(\mathcal{M}))
		\arrow[d,"\red"] \\
	\Vect(\mathcal{M}_\red)
		\arrow[r,"\Secinfty"]
	& \Proj(\Cinfty(\mathcal{M}_\red))
\end{tikzcd}
\end{equation}
\end{proposition}

\begin{proof}
	Our goal is to construct an isomorphism
	$\eta_E \colon \ConSecinfty(E)_\red \to \Secinfty(E_\red)$
	for every constraint vector bundle $E$ over $\mathcal{M}$.
	Thus let $s \in \ConSecinfty(E)_\Wobs$ be given.
	For any (possibly non-smooth) section
	$\sigma \colon \mathcal{M}_\red \to C$
	of the quotient map $\pi_\mathcal{M}$
	we can define a map
	$\eta_E(s) \colon \mathcal{M}_\red \to E_\red$
	by
	$\eta_E(s)(p) \coloneqq [s(\sigma(p))]$,
	which is a section of
	the vector bundle projection $\pr_{E_\red}$.
	Note that this map is independent of the choice
	of the section $\sigma$, since $s \in \ConSecinfty(E)_\Wobs$.
	Thus $\eta_E(s)$ is also smooth, since locally we can choose
	$\sigma$ to be smooth.
	So we end up with $\eta_E(s) \in \Secinfty(E_\red)$.
	Note also that $\eta_E$ is clearly $\ConCinfty(\mathcal{M})_\Wobs$-linear
	along the projection
	$\pi_{\ConCinfty(\mathcal{M})} \colon \ConCinfty(\mathcal{M})_\Wobs
	\to \Cinfty(M_\red)$.
	Now suppose $\eta_E(s) = 0$.
	Then $[s(\sigma(p))] = 0$ for all $p \in \mathcal{M}_\red$ and every section $\sigma$.
	Thus $\iota^\#s \in \Secinfty(E_\Null)$.
	This means that $\ker \eta_E = \ConSecinfty(E)_\Null$
	and therefore it induces an injective morphism
	$\eta_E \colon \ConSecinfty(E)_\red \to \Secinfty(E_\red)$
	of $\ConCinfty(\mathcal{M})_\red \simeq \Cinfty(\mathcal{M}_\red)$-modules.
	It remains to show that $\eta_E$ is also surjective.
	For this let $t \in \Secinfty(E_\red)$ be given.
	Now choose a splitting $E_\Wobs \simeq E_\Null \oplus \iota^\#E_\red$
	using \autoref{prop:ReductionOfConVect} \ref{prop:ReductionOfConVect_2}
	and define
	$s(q) \coloneqq \Theta^{-1}(t(\pi_\mathcal{M}(q)))$ for all $q \in C$
	and extend it to a section of $E_\Total$ by use of a tubular neighbourhood.
	By \eqref{eq:ReductionOfConVect_Parallel} $s$ is covariantly constant
	as a section of $E_\Wobs / E_\Null$, and therefore we have
	$s \in \ConSecinfty(E)_\Wobs$.
	Finally, we have $\eta_E(s) = t$, showing that $\eta_E$ is surjective,
	and thus an isomorphism.	
\end{proof}

%% file: constraint-fieldsandforms.tex
The close relationship between constraint vector bundles and 
constraint modules as established by the constraint Serre-Swan 
Theorem allows us to introduce further analogues of classical 
geometric structures on constraint manifolds, such as differential 
forms and multivector fields.
For both differential forms and multivector fields we can choose 
between the strong and non-strong tensor product, leading to two 
different graded constraint modules each.
In \autoref{sec:ConDifferentialForms} we will see that the classical 
de Rham differential is only well-defined on 
$\ConForms_{\strtensor}(\mathcal{M})$, but not on 
$\ConForms_{\tensor}(\mathcal{M})$.
In fact, the classical Cartan calculus, including the insertion of 
vector fields and Lie derivative of forms, is canonically given on
$\ConForms_{\strtensor}(\mathcal{M})$, thus singeling out 
$\ConForms_{\strtensor}$ as the correct constraint analogue of 
differential forms.
When we study constraint multivector fields in 
\autoref{sec:ConMultiVect} we find 
that here the situation is quite different,
since both $\ConVecFields_{\strtensor}(\mathcal{M})$
and $\ConVecFields_{\tensor}(\mathcal{M})$ carry the structure
of a constraint Gerstenhaber algebra.
Moreover, while $\ConVecFields_{\tensor}(\mathcal{M})$
seems to be the reasonable choice for constraint multivector fields,
since these are dual to the constraint forms 
$\ConForms_{\strtensor}(\mathcal{M})$, in the study of coisotropic 
reduction we will need $\ConVecFields_{\strtensor}(\mathcal{M})$.
Thus for constraint multivector fields there does not seem to be a 
preferred choice.

To ease notation we will, when considering constraint modules given 
by sections of constraint vector bundles, drop the subscript
$\ConCinfty(\mathcal{M})$ from the tensor products and simply write
$\tensor$ and $\strtensor$ instead of 
$\tensor[\ConCinfty(\mathcal{M})]$
or $\strtensor[\ConCinfty(\mathcal{M})]$.
Only when taking tensor products over other algebras or the base 
ring, we will use the usual subscripts.

%
%
\subsection{Differential Forms}
	\label{sec:ConDifferentialForms}

Before studying constraint differential forms we need to better 
understand constraint vector fields.
The following lemma shows how constraint vector fields can locally be 
characterized by their coefficient functions.

\begin{lemma}
	\label{lem:LocalConVect}
	\index{constraint!vector fields}
Let $\mathcal{M} = (M,C,D)$ be a constraint manifold of dimension 
$n = (n_\Total,n_\Wobs,n_\Null)$ and consider
$X \in \Secinfty(T M)$.
\begin{lemmalist}
	\item \label{lem:LocalConVect_1}
	We have
	$X \in \ConSecinfty(T\mathcal{M})_\Wobs$
	if and only if for every adapted chart $(U,x)$ around $p \in C$ it holds
	\begin{align}
		X^i &\in \ConCinfty(\mathcal{M}\at{U})_\Wobs \text{ if } i \in (n^*)_\Wobs,\\
		X^i &\in \ConCinfty(\mathcal{M}\at{U})_\Null \text{ if } i \in (n^*)_\Null,
	\end{align}
	where $X\at{U} = \sum_{i=1}^{n_\Total} X^i \frac{\del}{\del x^i}$.
	\item \label{lem:LocalConVect_2}
	We have $X \in \ConSecinfty(T\mathcal{M})_\Null$
	if and only if for every adapted chart around $p \in C$ it holds
	\begin{align}
		X^i &\in \ConCinfty(\mathcal{M}\at{U})_\Null \text{ if } 
		i \in (n^*)_\Wobs,
	\end{align}
	where $X\at{U} = \sum_{i=1}^{n_\Total} X^i \frac{\del}{\del x^i}$.
\end{lemmalist}
\end{lemma}

\begin{proof}
	By \autoref{ex:ConSections}
	locally we always find adapted coordinates such that
	\begin{align*}
		\iota^\#(\frac{\del}{\del x^i}) &\in \Secinfty(D\at{U}), \text{ if } i \in \{1, \dotsc, n_\Null\}, \\
		\shortintertext{and}
		\iota^\#(\frac{\del}{\del x^i}) &\in \Secinfty(TC\at{U}), \text{ if } i \in \{n_\Null +1, \dotsc, n_\Wobs\}.
	\end{align*}
	We have $X \in \ConSecinfty(T \mathcal{M})_\Wobs$
	if and only if $\iota^\# X \in \Secinfty(TC)$
	and $[Y, \iota^\#X] \in \Secinfty(D)$ hold for all $Y \in \Secinfty(D)$.
	The first condition exactly means that locally we have
	$X^i \in \ConCinfty(\mathcal{M}\at{U})_\Null$ for all $i \in \{ n_\Wobs +1, \dotsc, n_\Total \} = (n^*)_\Null$.
	Moreover, since $D$ is locally spanned by $\frac{\del}{\del x^1}, \dotsc, \frac{\del}{\del x^{n_\Null}}$
	the second condition shows $X^i \in \ConCinfty(\mathcal{M})_\Wobs$ for
	$i \in \{ n_\Null +1, \dotsc, n_\Total \} = (n^*)_\Wobs$.
	This shows the first part.
	The second part follows since $X \in \ConSecinfty(T\mathcal{M})_\Null$
	if and only if $\iota^\#X \in \Secinfty(D)$.
\end{proof}

With the help of this local characterization we can now identify 
constraint vector fields with constraint derivations, see 
\autoref{prop:ConDer},
using the Lie derivative:

\begin{proposition}
	\label{prop:ConVectorFieldsAreConDer}
	\index{constraint!Lie derivative}
Let $\mathcal{M} = (M,C,D)$ be a constraint manifold.
Then
\glsadd{LieDer}
\begin{equation}
	\Lie \colon \ConSecinfty(T\mathcal{M})
	\to \ConDer(\ConCinfty(\mathcal{M}))
\end{equation}
given by the Lie derivative is an isomorphism of constraint
$\ConCinfty(\mathcal{M})$-modules.
\end{proposition}

\begin{proof}
	From classical differential geometry we know that $\Lie$ is an isomorphism on the
	$\TOTAL$-com\-po\-nents.
	To show that $\Lie$ is a constraint morphism consider
	$X \in \ConSecinfty(T \mathcal{M})_\Null$
	and $f \in \ConCinfty(\mathcal{M})_\Wobs$.
	Then
	\begin{equation*}
		(\Lie_X f)\at{C}
		= \Lie_{\iota^\#X} f\at{C}
		= 0,
	\end{equation*}
	since $X\at{C} \in \Secinfty(D)$.
	Thus $\Lie$ maps the $\NULL$-component to the $\NULL$-component.
	Now let $X \in \ConSecinfty(T \mathcal{M})_\Wobs$ be given.
	Then for $f \in \ConCinfty(\mathcal{M})_\Wobs$
	we get
	\begin{equation*}
		\Lie_Y (\Lie_X f)\at{C}
		= \Lie_{[Y,\iota^\#X]}f\at{C} 
		+ \Lie_{\iota^\#X} \underbrace{\Lie_Y f\at{C}}_{=0}
		= \Lie_{[Y,\iota^\#X]} f\at{C} = 0,
	\end{equation*}
	for all $Y \in \Secinfty(D)$, since
	$[Y,\iota^\#X] \in \Secinfty(TC)$.
	Finally, for $f \in \ConCinfty(\mathcal{M})_\Null$
	we have $f\at{C} = 0$ and therefore
	\begin{equation*}
		(\Lie_X f)\at{C} = \Lie_{\iota^\#X}f\at{C} = 0,
	\end{equation*}
	which shows that $\Lie$ is a constraint morphism.
	Since the $\TOTAL$-component of $\Lie$ is just the classical Lie 
	derivative, which is an isomorphism, $\Lie$ is a constraint 
	monomorphism.
	To show that $\Lie$ is also a regular epimorphism let
	$D \in \ConDer(\ConCinfty(\mathcal{M}))_\Wobs$ be given.
	Since $D$ is in particular a derivation of $\Cinfty(M)$
	we know that there exists $X \in \Secinfty(TM)$ such that
	$\Lie_X = D$.
	Choose an adapted chart $(U,x)$ around $p \in C$, then
	$X\at{U} = \sum_{i=1}^{n_\Total} X^i \frac{\del}{\del x^i}$.
	Since $\Lie_X$ is a constraint derivation we get
	$X^i = \Lie_X(x^i) \in \ConCinfty(\mathcal{M}\at{U})_\Null$
	for all
	$i \in \{n_\Wobs +1, \dotsc, n_\Total\} = (n^*)_\Null$
	and 
	$X^i = \Lie_X(x^i) \in \ConCinfty(\mathcal{M}\at{U})_\Wobs$
	for all $i \in \{n_\Null +1, \dotsc, n_\Total\} = (n^*)_\Wobs$,
	by \autoref{ex:ConFunctionsOnConMfld} 
	\ref{ex:ConFunctionsOnConMfld_1}.
	And thus $X \in \ConSecinfty(T\mathcal{M})_\Wobs$
	using \autoref{lem:LocalConVect} \ref{lem:LocalConVect_1}.
	With the same line of reasoning we obtain $X \in 
	\ConSecinfty(T\mathcal{M})_\Null$ if
	$D \in \ConDer(\ConCinfty(\mathcal{M}))_\Null$,
	showing that $\Lie$ is a regular epimorphism, and
	therefore an isomorphism.
\end{proof}

With this we can transport the constraint Lie algebra 
structure from
$\ConDer(\ConCinfty(\mathcal{M}))$
to $\ConSecinfty(T \mathcal{M})$.
This is just the usual Lie bracket of vector fields, but now we see 
that it is actually compatible with the constraint structure.
Alternatively, one could directly check that the classical Lie bracket of vector fields
yields a constraint Lie algebra structure.

With this at hand let us introduce constraint differential forms.
Since there are two tensor products available we can define 
constraint differential
forms in two ways.

\begin{definition}[Constraint Differential Forms]
	\index{constraint!differential forms}
Let $\mathcal{M} = (M,C,D)$ be a constraint manifold.
We denote by
\glsadd{ConFormsTensor}
\glsadd{ConFormsStrTensor}
\begin{align}
	\ConForms_{\tensor}^\bullet(\mathcal{M}) 
	\coloneqq \Anti^\bullet_{\tensor} \ConSecinfty(T^* \mathcal{M})
	= \bigoplus_{k=0}^{\infty} \Anti^k_{\tensor} 
	\ConSecinfty(T^*\mathcal{M})
	\shortintertext{and}
	\ConForms_{\strtensor}^\bullet(\mathcal{M}) 
	\coloneqq \Anti^\bullet_{\strtensor} \ConSecinfty(T^* \mathcal{M})
	= \bigoplus_{k=0}^{\infty} \Anti^k_{\strtensor} 
	\ConSecinfty(T^*\mathcal{M})
\end{align}
the graded strong constraint modules of
\emph{constraint differential forms} on $\mathcal{M}$.
\end{definition}

Note that $\ConForms_{\tensor}^\bullet(\mathcal{M}) \simeq (\Anti_{\strtensor}^\bullet \ConSecinfty(T\mathcal{M}))^*$
and $\ConForms_{\strtensor}^\bullet(\mathcal{M}) \simeq (\Anti_{\tensor}^\bullet \ConSecinfty(T\mathcal{M}))^*$.
Thus $\alpha \in \ConForms_{\tensor}^k(\mathcal{M})$ can be evaluated at
$X_1 \tensor \dots \tensor X_k \in \Anti_{\strtensor}^\bullet \ConSecinfty(T\mathcal{M})$,
while $\alpha \in \ConForms_{\strtensor}^k(\mathcal{M})$ can be evaluated at
$X_1 \tensor \dots \tensor X_k \in \Anti_{\tensor}^\bullet \ConSecinfty(T\mathcal{M})$.
For $\ConForms_{\strtensor}^\bullet(\mathcal{M})$
there is a good constraint Cartan calculus as we see in the following.


\begin{proposition}[Cartan calculus]
	\label{prop:ConCartanCalculus}
	\index{constraint!Cartan calculus}
Let $\mathcal{M} = (M,C,D)$ be a constraint manifold.
\begin{propositionlist}
	\item $\ConForms_{\strtensor}^\bullet(\mathcal{M})$ is an embedded 
	graded commutative strong constraint algebra with respect to the 
	wedge product $\wedge$.
	\item \index{insertion}
	The insertion of vector fields into forms defines a 
	constraint $\ConCinfty(\mathcal{M})$-module morphism
	\glsadd{insertion}
	\begin{equation}
		\ins \colon \ConSecinfty(T\mathcal{M}) \to 
		\ConDer^{-1}(\ConForms_{\strtensor}^\bullet(\mathcal{M})),
	\end{equation}
	with $\ConDer^{-1}(\ConForms_{\strtensor}^\bullet(\mathcal{M}))$
	denoting the graded constraint derivations of degree $-1$.
	\item \index{constraint!Lie derivative}
	The Lie derivative defines a $\Reals$-linear constraint morphism
	\begin{equation}
		\Lie \colon \ConSecinfty(T\mathcal{M}) \to 
		\ConDer^{0}(\ConForms_{\strtensor}^\bullet(\mathcal{M}))
	\end{equation}
	into the graded constraint derivations of degree $0$ of 
	$\ConForms_{\strtensor}^\bullet(\mathcal{M})$.
	\item \index{constraint!de Rham differential}
	The de Rham differential defines a graded constraint 
	derivation
	\glsadd{deRham}
	\begin{equation}
		\D \colon \ConForms_{\strtensor}^\bullet(\mathcal{M}) \to 
		\ConForms_{\strtensor}^{\bullet+1}(\mathcal{M})
	\end{equation}
	of degree $+1$.
\end{propositionlist}
\end{proposition}

\begin{proof}
In all cases we only need to show that the involved maps are 
actually constraint maps.
For the first part this is clear by the definition of
$\ConForms^\bullet(\mathcal{M})$.
For the insertion consider
$X \in \ConSecinfty(T\mathcal{M})_\Null$.
Then $\ins_X \alpha \in \ConSecinfty(T\mathcal{M})_\Null$
for all
$\alpha \in \ConSecinfty(T^*\mathcal{M})_\Wobs$.
Since $\ins_X$ is a derivation of the wedge product
it maps $\ConForms^\bullet(\mathcal{M})_\Wobs$
to $\ConForms^\bullet(\mathcal{M})_\Null$.
Now consider $X \in \ConSecinfty(T\mathcal{M})_\Wobs$.
Then again by the derivation property it is easy to see that
$\ins_X(\ConForms^\bullet(\mathcal{M})_\Wobs) \subseteq 
\ConForms^\bullet(\mathcal{M})_\Wobs$
and
$\ins_X(\ConForms^\bullet(\mathcal{M})_\Null) \subseteq 
\ConForms^\bullet(\mathcal{M})_\Null$.
Thus $\ins$ is a constraint morphism.
Since the Lie derivative is again a derivation and we know, by
\autoref{prop:ConVectorFieldsAreConDer}
and from the fact that
$\Lie_X Y = [X,Y]$,
that $\Lie_X$ is a constraint endomorphism of 
$\ConSecinfty(T\mathcal{M})$,
it follows that $\Lie$ is a constraint morphism.
For the de Rham differential we can argue with the 
formula
\begin{align*}
	(\D\alpha)(X_0 \tensor \cdots \tensor X_k)
	&= \sum_{i=0}^{k} (-1)^k \Lie_{X_i}(\alpha(X_0,\dotsc, 
	\overset{i}{\wedge}, \dotsc, X_k))\\
	&\qquad + \sum_{i < j} (-1)^{i+j}
	\alpha([X_i,X_j],X_0,\dotsc, 
	\overset{i}{\wedge}, \dotsc, \overset{j}{\wedge}, \dotsc, 
	X_k),
\end{align*}
for $X_0 \tensor \dotsc \tensor X_k \in (\Anti_{\tensor}^\bullet 
\ConSecinfty(T\mathcal{M}))_\Total$
to see that $\D$ is a constraint morphism.
For example, if $\alpha \in \ConForms_{\strtensor}^k(\mathcal{M})_\Null$ is given, we have
\begin{equation*}
	(\D \alpha)(X_0 \tensor \cdots \tensor X_k) \in \ConCinfty(\mathcal{M})_\Null
\end{equation*}
for all
$X_0, \dotsc, X_k \in \ConSecinfty(T\mathcal{M})_\Wobs$,
since from
\begin{equation*}
	\alpha(X_0, \dotsc, \overset{i}{\wedge},\dotsc, X_k) \in \ConCinfty(\mathcal{M})_\Null
\end{equation*}
it follows that
\begin{equation*}
	\Lie_{X_i} \alpha(X_0, \dotsc, \overset{i}{\wedge},\dotsc, X_k) \in \ConCinfty(\mathcal{M})_\Null
\end{equation*}
and from
$[X_i,X_j] \in \ConSecinfty(T\mathcal{M})_\Null$
it follows
\begin{equation*}
	\alpha([X_i,X_j],X_0,\dotsc,\overset{i}{\wedge},\dotsc,\overset{j}{\wedge},\dotsc,X_k) \in \ConCinfty(\mathcal{M})_\Null.
\end{equation*}
Thus we have $\D\alpha \in \ConForms_{\strtensor}^{k+1}(\mathcal{M})_\Null$.
In a similar way we can argue for $\alpha \in \ConForms_{\strtensor}^k(\mathcal{M})_\Wobs$.
\end{proof}

Since $\ins$,  $\Lie$ and $\D$ are completely determined by their $\TOTAL$-components, we immediately get all the usual formulas from the classical Cartan calculus, such as e.g. Cartan's magic formula
\begin{equation}
	\Lie_X = [\ins_X,\D].
\end{equation}
We cannot expect a similarly well behaved Cartan calculus on
$\ConForms_{\tensor}^\bullet(\mathcal{M})$, since in this case
the de Rham differential is not well-defined, as the next example shows.

\begin{example}
Consider $\mathcal{M} = (\Reals^{n_\Total}, 
\Reals^{n_\Wobs},\Reals^{n_\Null})$
with $n_\Null \geq 1$
and let $\alpha = x^1 \D x^{n_\Total} \in 
\ConSecinfty(T^*\mathcal{M})_\Null$.
Then we have
\begin{equation}
	\D \alpha = \D x^1 \wedge \D x^{n_\Total}
	\in \ConSecinfty(T^*\mathcal{M})_\Total \wedge 
	\ConSecinfty(T^*\mathcal{M})_\Null
	\nsubseteq \ConForms_{\tensor}^2(\mathcal{M})_\Null.
\end{equation}
\end{example}

\begin{remark}
	\index{constraint!Lie algebroids}
The constraint de Rham differential can also be understood 
as the constraint Lie algebroid differential \cite[Chapter 7]{mackenzie:2005a} for 
the constraint Lie algebroid $T\mathcal{M}$, see
\autoref{ex:ConLieAlgebroid}.
\end{remark}

Even though $\ConForms_{\tensor}^\bullet(\mathcal{M})$ does not carry as rich an algebraic structure
as $\ConForms_{\strtensor}^\bullet(\mathcal{M})$ it will nevertheless play an important
role as those constraint forms which are dual to constraint multivector fields
of the form $\Anti_{\strtensor}^\bullet \ConSecinfty(T\mathcal{M})$.

\begin{definition}[Constraint de Rham cohomology]
	\label{def:ConDeRhamCohomology}
	\index{constraint!de Rham cohomology}
Let $\mathcal{M} = (M,C,D)$ be a constraint\linebreak
manifold.
We call the constraint complex 
$(\ConForms_{\strtensor}^\bullet(\mathcal{M}), \D)$
the \emph{(constraint) de Rham complex} of $\mathcal{M}$.
Its cohomology is called the \emph{constraint de Rham cohomology} of $\mathcal{M}$
and will be denoted by \glsadd{deRhamCohomology}$\functor{H}_\deRham^\bullet(\mathcal{M})$.
\end{definition}

Recall from \autoref{prop:ConCohomology}
and \autoref{def:ConQuotientkModule} that
$\functor{H}^k_\deRham(\mathcal{M})_\Wobs
= \frac{\ker\D^k_\Wobs}{\image \D^{k-1}_\Wobs}$ with
\begin{equation}
	\image \D^{k-1}_\Wobs = \left\{ \omega \in \Forms^k(M) \bigm| \exists \eta \in \ConForms_{\strtensor}^{k-1}(\mathcal{M})_\Wobs : \D \eta = \omega \right\}.
\end{equation}
Thus there might exist $\omega \in \ConForms_{\strtensor}^k(\mathcal{M})_\Wobs$
which is exact in the classical sense, but not exact as a constraint form.
This means that even if $\functor{H}_\deRham(M)$ is trivial,
$\functor{H}_\deRham(\mathcal{M})$ might be non-trivial.
Here the fact that the category $\injConMod_{\field{k}}$ is not closed under colimits enters crucially, since this allows for a non-embedded cohomology, cf. \autoref{ex:QuotientEmbeddedConModk}.
Nevertheless, there is a constraint Poincaré Lemma:

\begin{proposition}[Constraint Poincar\'e Lemma]
	\index{Poincar\'e Lemma}
For $\Reals^n = \Reals^{(n_\Total, n_\Wobs, n_\Null)}$ the constraint de Rham cohomology is given by
\begin{equation}
	\functor{H}_{\deRham}^k(\Reals^n) = 
	\begin{cases}
		(\Reals,\Reals,0) &\text{ for k = 0} \\
		(0,0,0) &\text{ for } k \geq 1
	\end{cases}.
\end{equation}
\end{proposition}

\begin{proof}
The $\TOTAL$-component is exactly the classical Poincar\'e Lemma.
For $k=0$ forms are just functions and hence there do not exist exact ones.
A function is closed if and only if it is constant.
Thus $\functor{H}_{\deRham}^0(\Reals^n)_\Wobs$ consists of the constant functions, which are on
$\Reals^{n_\Wobs}$ constant along $\Reals^{n_\Null}$.
But this is fulfilled by every constant function, hence 
$\functor{H}_{\deRham}^0(\Reals^n)_\Wobs = \Reals$.
The only constant function that vanishes on $\Reals^{n_\Null}$
is the zero function, hence
$\functor{H}_{\deRham}^0(\Reals^n)_\Null = 0$.
Now let $k \geq 1$ be given and consider a closed
$\omega \in \ConForms_{\strtensor}^k(\Reals^n)_\Total$.
From the classical Poincar\'e Lemma we know that $\omega = \D \eta$ is exact with
\begin{equation*} \label{eq:PoincareHomotopy}
	\eta\at{x}(v_1,\dotsc,v_{k-1}) = \int_{0}^{1} t^{k-1}\omega\at{tx}(x,v_1, \dotsc, v_{k-1}) \D t.
	 \tag{$*$}
\end{equation*}
Now if $\omega \in \ConForms_{\strtensor}^k(\Reals^n)_\Null$, then
$\iota^*\eta(v_1,\dotsc,v_{k-1})$ vanishes for $v_1 \tensor \dots \tensor v_{k-1} \in ((\Reals^n)^{\tensor k})_\Wobs$, since so does $\omega$.
If $\omega \in \ConForms_{\strtensor}^k(\Reals^n)_\Wobs$
we know $\iota^*\omega(w_1, \dotsc, w_k) = 0$ for all
$w_1 \tensor \dots \tensor w_k \in ((\Reals^n)^{\tensor k})_\Null$.
Then clearly $\iota^*\eta(v_1, \dots v_{k-1}) = 0$
for all $v_1 \tensor \dots \tensor v_{k-1} \in ((\Reals^n)^{\tensor k-1})_\Null$.
Moreover,
$\iota^*\omega$ is constant along $\Reals^{n_0}$, thus
\begin{equation*}
	\omega\at{t(x+y)}(x+y,v_1, \dotsc, v_{k-1}) = \omega\at{tx}(x,v_1, \dotsc, v_{k-1})
\end{equation*}
for all $x,v_1, \dotsc, v_{k-1} \in \Reals^{n_\Wobs}$
and $y \in \Reals^{n_\Null}$.
Then $\iota^*\eta$ is constant along $\Reals^{n_\Null}$
by \eqref{eq:PoincareHomotopy}.
\end{proof}

\subsubsection{Reduction}

Both types of constraint forms reduce to the classical forms on the reduced manifolds:

\begin{proposition}[Constraint forms vs. reduction]
	\label{prop:ConFormsVSReduction}
	\index{reduction!differential forms}
Let $\mathcal{M} = (M,C,D)$ be a constraint\linebreak
manifold.
\begin{propositionlist}
	\item There exists a canonical isomorphism
	$\ConForms_{\tensor}^\bullet(\mathcal{M})_\red
	\simeq \Forms^\bullet(\mathcal{M}_\red)$
	of graded $\ConCinfty(\mathcal{M})_\red$-modules.
	\item There exists a canonical isomorphism
	$\ConForms_{\strtensor}^\bullet(\mathcal{M})_\red
	\simeq \Forms^\bullet(\mathcal{M}_\red)$
	of complexes.
\end{propositionlist}
\end{proposition}

\begin{proof}
We combine established results to the following chain of canonical isomorphisms:
\begin{align*}
	\ConForms_{\tensor}^\bullet(\mathcal{M})_\red
	&= \big( \bigoplus_{k=0}^\infty \Anti_{\tensor}^k \ConSecinfty(T^*\mathcal{M}) \big)_\red
	\simeq \bigoplus_{k=0}^\infty (\Anti_{\tensor}^k \ConSecinfty(T^*\mathcal{M}))_\red \\
	&\simeq  \bigoplus_{k=0}^\infty \Anti^k \ConSecinfty(T^*\mathcal{M})_\red
	\simeq \bigoplus_{k=0}^\infty \Anti^k \Secinfty(T^*\mathcal{M}_\red)
	= \Forms^k(\mathcal{M}_\red).
\end{align*}
Since we know that the reduction of $\strtensor$ and $\tensor$ agree and the reduced 
de Rham differential $\D_\red$ fulfils the same local characterization as the de Rham differential on 
$\mathcal{M}_\red$ the second part follows.
\end{proof}

Since by \autoref{prop:CohomologyCommutesWithReduction} cohomology commutes with reduction in general, we obtain
as a special case
\begin{equation}
	\functor{H}_\deRham(\mathcal{M})_\red \simeq \functor{H}_\deRham(\mathcal{M}_\red).
\end{equation}

%
%
\subsection{Multivector Fields and Poisson Manifolds}
\label{sec:ConMultiVect}

Let us now turn our attention to constraint multivector fields.
As for constraint differential forms we can define multivector fields using both tensor products available.

\begin{definition}[Constraint multivector fields]
	\index{constraint!multivector fields}
Let $\mathcal{M} = (M,C,D)$ be a constraint manifold.
Then we denote by
\glsadd{ConMultVectTensor}
\glsadd{ConMultVectStrTensor}
\begin{align}
	\ConVecFields_{\tensor}^\bullet(\mathcal{M}) 
	\coloneqq \Anti^\bullet_{\tensor} \ConSecinfty(T\mathcal{M})
	\simeq \ConSecinfty(\Anti^\bullet_{\tensor} T\mathcal{M})
	\shortintertext{and}
	\ConVecFields_{\strtensor}^\bullet(\mathcal{M}) 
	\coloneqq \Anti^\bullet_{\strtensor} \ConSecinfty(T\mathcal{M})
	\simeq \ConSecinfty(\Anti^\bullet_{\strtensor} T\mathcal{M})
\end{align}
the graded strong constraint modules of
\emph{constraint multivector fields}
on $\mathcal{M}$.
\end{definition}

In low degrees we can easily characterize constraint multivector fields in local charts.

\begin{lemma}
	\label{lem:LocalConBiVect}
	\index{constraint!bivector fields}
Let $\mathcal{M} = (M,C,D)$ be a constraint manifold of dimension 
$n = (n_\Total,n_\Wobs,n_\Null)$ and consider
$\pi \in \VecFields^2(M)$.
\begin{lemmalist}
	\item \label{lem:LocalConBiVect_3}
	We have $\pi \in \ConVecFields_{\tensor}^2(\mathcal{M})_\Wobs$
	if and only if for every $p \in C$ there exists a local chart 
	$(U,x)$ around $p$
	such that
	$\pi\at{U} = \sum_{i,j=1}^{n_\Total} \pi^{ij} \frac{\del}{\del 
		x^i} \wedge \frac{\del}{\del x^j}$
	with
	\begin{equation}
	\begin{split}
		\pi^{ij} &\in \ConCinfty(\mathcal{M}\at{U})_\Null
		\text{ if } (i,j) \in  (n^* \strtensor n^*)_\Null
		= n \ConGrid[0][0][2][0][0][2][2][2][2] n,\\
		\pi^{ij} &\in \ConCinfty(\mathcal{M}\at{U})_\Wobs
		\text{ if } (i,j) \in (n^* \strtensor n^*)_\Wobs
		= n \ConGrid[0][0][2][0][2][2][2][2][2] n.
	\end{split}
	\end{equation}
	\item \label{lem:LocalConBiVect_4}
	We have $\pi \in \ConVecFields_{\tensor}^2(\mathcal{M})_\Null$
	if and only if for every $p \in C$ there exists a local chart $(U,x)$ 
	around $p$
	such that
	$\pi\at{U} = \sum_{i,j=1}^{n_\Total} \pi^{ij} \frac{\del}{\del x^i} 
	\wedge \frac{\del}{\del x^j}$
	with
	\begin{align}
		\pi^{ij} &\in \ConCinfty(\mathcal{M}\at{U})_\Null \text{ if } 
		(i,j) \in  (n^* \strtensor n^*)_\Wobs = n 
		\ConGrid[0][0][2][0][2][2][2][2][2] n.
	\end{align}
	\item \label{lem:LocalConBiVect_1}
	We have $\pi \in \ConVecFields_{\strtensor}^2(\mathcal{M})_\Wobs$
	if and only if for every $p \in C$ there exists a local chart 
	$(U,x)$ around $p$
	such that
	$\pi\at{U} = \sum_{i,j=1}^{n_\Total} \pi^{ij} \frac{\del}{\del 
	x^i} \wedge \frac{\del}{\del x^j}$
	with
	\begin{equation} \label{eq:ConditionsConBivectField_N}
	\begin{split}
		\pi^{ij} &\in \ConCinfty(\mathcal{M}\at{U})_\Null
		\text{ if } (i,j) \in  (n^* \tensor n^*)_\Null
		= n \ConGrid[0][0][0][0][0][2][0][2][2] n,\\
		\pi^{ij} &\in \ConCinfty(\mathcal{M}\at{U})_\Wobs
		\text{ if } (i,j) \in (n^* \tensor n^*)_\Wobs
		= n \ConGrid[0][0][0][0][2][2][0][2][2] n.
	\end{split}
	\end{equation}
	\item \label{lem:LocalConBiVect_2}
	We have $\pi \in \ConVecFields_{\strtensor}^2(\mathcal{M})_\Null$
	if and only if for every $p \in C$ there exists a local chart $(U,x)$ 
	around $p$
	such that
	$\pi\at{U} = \sum_{i,j=1}^{n_\Total} \pi^{ij} \frac{\del}{\del x^i} 
	\wedge \frac{\del}{\del x^j}$
	with
	\begin{align}
		\pi^{ij} &\in \ConCinfty(\mathcal{M}\at{U})_\Null \text{ if } 
		(i,j) \in  (n^* \tensor n^*)_\Wobs = n 
		\ConGrid[0][0][0][0][2][2][0][2][2] n.
	\end{align}
\end{lemmalist}
\end{lemma}

\begin{proof}
Let us proof \ref{lem:LocalConBiVect_3}, the other statements follow analogously:
For every $p \in C$ we find by \autoref{ex:ConSections} an adapted chart $(U,x)$
such that
\begin{equation*}
	\frac{\del}{\del x^i} \in \ConSecinfty(TU)_\Null \text{ if } i \in n_\Null \\
	\quad\text{and}\quad
	\frac{\del}{\del x^i} \in \ConSecinfty(TU)_\Wobs \text{ if } i \in n_\Wobs.
\end{equation*}
By definition we have $\pi \in \ConVecFields_{\tensor}^2(\mathcal{M})$ if and only
if for every $p \in C$ it holds that $\pi(p) \in T_pC \wedge T_pC$
and $\Lie_X \cc{\pi\at{C}} = 0$ for all $X \in \Secinfty(D)$.
Thus, using $\pi\at{U} = \sum_{i,j=1}^{n_\Total} \pi^{ij} \frac{\del}{\del x^i} \wedge \frac{\del}{\del x^j}$
we see that $\pi^{ij} = 0$ whenever $i <n_\Wobs$ or $j > n_\Wobs$, i.e. whenver
$(i,j) \in n \ConGrid[0][0][2][0][0][2][2][2][2] n$.
Moreover, we have
\begin{align*}
	\Lie_{\frac{\del}{\del x^k}} \cc{\pi\at{U\cap C}}
	= \sum_{i,j=n_\Null +1}^{n_\Wobs} \frac{\del}{\del x^k}\pi^{ij} \frac{\del}{\del x^i} \wedge \frac{\del}{\del x^j}
\end{align*}
for all $k = 1,\dotsc, n_\Null$,
which vanishes if and only if $\pi^{ij} \in \ConCinfty(\mathcal{M}\at{U})_\Wobs$
for all $i,j = n_\Null + 1, \dotsc, n_\Wobs$, i.e. if
$(i,j) \in n \ConGrid[0][0][0][0][0][2] n$.
Together this yields \ref{lem:LocalConBiVect_3}.
\end{proof}

The following example shows that every constraint manifold constructed from a coisotropic submanifold of a Poisson manifold
carries a constraint bivector field in $\ConVecFields_{\strtensor}^2(\mathcal{M})$,
while a Poisson submanifold yields a constraint bivector field in
$\ConVecFields_{\tensor}^2(\mathcal{M})$.

\begin{example}
	\label{ex:ConBivectorsFromPoissonManifolds}
	Let $(M,\pi)$ be a Poisson manifold.
\begin{examplelist}
	\item \label{ex:ConPoissonMfld}
	\index{coisotropic submanifold}
	If $C \subseteq M$ is a closed coisotropic submanifold allowing for a smooth reduction
	we denote by $\mathcal{M} = (M,C,D)$ the constraint manifold with $D$ the characteristic distribution
	of the coisotropic submanifold $C$.
	Let $n = (n_\Total,n_\Wobs,n_\Null)$ be its constraint dimension.
	Then $\pi \in \Anti^2 \Secinfty(TM)$ is a bivector field, fulfilling 
	$\iota^\#\pi \in \Secinfty(TC \wedge TC + \iota^\#TM \wedge D)
	= \Secinfty((\Anti_{\strtensor} T \mathcal{M})_\Wobs)$.
	In an adapted coordinate chart $(U,x)$ around $p \in C$, cf. \autoref{lem:AdaptedLocalFrames},
	we have
	\begin{equation}
		\cc{\iota^\#\pi\at{U \cup C}} = \sum_{i,j = n_\Null +1}^{n_\Wobs} \pi^{ij} \frac{\del}{\del x^i} \wedge \frac{\del}{\del x^j},
	\end{equation}
	and thus for all $\ell = 1, \dotsc n_\Null$
	\begin{align*}
		\nabla_{\frac{\del}{\del x^\ell}} \cc{\iota^\#\pi\at{U \cup C}}
		= \sum_{i,j = n_\Null +1}^{n_\Wobs} \pi^{ij} \left[\frac{\del}{\del x^\ell}, \frac{\del}{\del x^i}\right] \wedge \frac{\del}{\del x^j}
		+ \sum_{i,j = n_\Null +1}^{n_\Wobs} \pi^{ij} \frac{\del}{\del x^i} \wedge \left[\frac{\del}{\del x^\ell},\frac{\del}{\del x^j}\right]
		= 0
	\end{align*}
	holds.
	Here we crucially use that $\pi^{ij} \in \ConCinfty(U)_\Wobs$ for all $i,j = n_\Null, \dotsc, n_\Wobs$.
	Since $D$ is locally spanned by $\frac{\del}{\del x^1},\dotsc, \frac{\del}{\del x^{n_\Null}}$
	we have $\pi \in \ConVecFields_{\strtensor}^2(\mathcal{M})_\Wobs$.
	\item \index{Poisson submanifold}	
	Since every Poisson submanifold is in particular coisotropic, every closed Poisson submanifold
	gives a constraint manifold $\mathcal{M} = (M,C,0)$ the constraint manifold with trivial distribution.
	Let $n = (n_\Total,n_\Wobs,0)$ be its constraint dimension.
	Then $\pi \in \Anti^2\Secinfty(TM)$ restricts to a bivector field
	$\pi\at{C} \in \Anti^2\Secinfty(TC) = \Secinfty((\Anti^2 T\mathcal{M})_\Wobs)$.
	Since $D$ is trivial we thus get $\pi \in 
	\ConVecFields_{\tensor}^2(\mathcal{M})_\Wobs$.
	\item Every closed Poisson submanifold $C$ of a Poisson manifold $M$ can also be equipped with
	another distribution $D$ given by the symplectic leaves of $C$.
	In general, the leaf space will not be smooth, but e.g. for certain types of Poisson manifolds of compact type at least an orbifold structure on the leaf space can be achieved, see \cite{crainic.fernandes.martinez:2019b,crainic.fernandes.martinez:2019a}.
	Note that in the case of a smooth leaf space we obtain a constraint manifold
	$\mathcal{M} = (M,C,D)$ with a constraint Poisson structure $\pi \in \Anti^2\ConVecFields_{\tensor}^2(\mathcal{M})_\Wobs$.
	The reduced space then describes the transversal structure.
\end{examplelist}
\end{example}

This suggests that a constraint manifold equipped with a constraint bivector field
$\pi \in \ConVecFields_{\strtensor}^2(\mathcal{M})_\Wobs$
fulfilling the Jacobi identity
induces a coisotropic structure on its submanifold.
On the other hand $\pi \in \ConVecFields_{\tensor}^2(\mathcal{M})_\Wobs$
fulfilling the Jacobi identity seems to induce a Poisson structure on 
$C$, which drops to $\mathcal{M}_\red$.
To make this precise we first introduce the Schouten bracket for constraint multivector fields.

\begin{proposition}[Constraint Schouten bracket]
	\index{constraint!Schouten bracket}
Let $\mathcal{M} = (M,C,D)$ be a constraint manifold.
The classical Schouten bracket defines constraint graded Lie algebra structures
\glsadd{SchoutenBracket}
\begin{align}
	\Schouten{\argument, \argument} 
	\colon \ConVecFields_{\tensor}^{k+1}(\mathcal{M}) \tensor[\field{k}] \ConVecFields_{\tensor}^{\ell+1}(\mathcal{M})
	\to \ConVecFields_{\tensor}^{k+\ell + 1}(\mathcal{M})
	\shortintertext{and}
	\Schouten{\argument, \argument} 
	\colon \ConVecFields_{\strtensor}^{k+1}(\mathcal{M}) \tensor[\field{k}] \ConVecFields_{\strtensor}^{\ell+1}(\mathcal{M})
	\to \ConVecFields_{\strtensor}^{k+\ell + 1}(\mathcal{M})
\end{align}
with respect to the degree shifted by $1$.	
\end{proposition}

\begin{proof}
	This follows directly from the formula
	\begin{equation*}
		\Schouten{X_0 \wedge \cdots \wedge X_k, Y_0 \wedge \cdots \wedge Y_\ell}
		= \sum_{i=0}^k\sum_{j=0}^\ell (-1)^{i+j} [X_i,Y_j] \wedge X_1 \cdots \overset{i}{\wedge} \cdots X_k \wedge Y_0 \wedge \cdots \overset{j}{\wedge} \cdots \wedge Y_\ell
	\end{equation*}
	and the fact that $[\argument, \argument]$ is a constraint Lie bracket on
	$\ConVecFields^1(\mathcal{M})$.
\end{proof}

It is important to note that even for 
$\ConVecFields_{\strtensor}^\bullet(\mathcal{M})$
we do \emph{not} obtain a strong constraint Lie algebra structure.
One way to see this is to note that $\ConDer(\ConCinfty(\mathcal{M}))$
is only a constraint Lie algebra, even though $\ConCinfty(\mathcal{M})$
is a strong constraint algebra.
Ultimately, this comes from the fact that $\ConHom$ is adjoint
to $\tensor$ and not $\strtensor$.

\begin{corollary}
	\label{cor:ConDGLAOfConMultVectFields}
Let $\mathcal{M} = (M,C,D)$ be a constraint manifold.
Then
\begin{corollarylist}
	\item $\big(\ConVecFields_{\tensor}^{\bullet+1}(\mathcal{M}), \D = 0, \Schouten{\argument, \argument}\big)$
	is a constraint DGLA.
	\item $\big(\ConVecFields_{\strtensor}^{\bullet+1}(\mathcal{M}), \D = 0, \Schouten{\argument, \argument}\big)$
	is a constraint DGLA.
\end{corollarylist}
\end{corollary}

In contrast to constraint differential forms there is no preferred choice of the tensor products, at least from the point of view of available algebraic structure.
Nevertheless, \autoref{ex:ConBivectorsFromPoissonManifolds} \ref{ex:ConPoissonMfld} shows that if we are interested
in coisotropic submanifolds we are forced to consider
$\ConVecFields_{\strtensor}^2(\mathcal{M})$
instead of $\ConVecFields_{\tensor}^2(\mathcal{M})$.
Thus we define the constraint analogue of a Poisson manifold as follows.

\begin{definition}[Constraint Poisson manifold]
	\label{def:ConPoissonManifold}
	\index{constraint!Poisson manifold}
A \emph{constraint Poisson manifold} consists of a constraint manifold
$\mathcal{M} = (M,C,D)$ together with a constraint bivector field
$\pi \in \ConVecFields_{\strtensor}^2(\mathcal{M})_\Wobs$
such that $\Schouten{\pi,\pi} = 0$.
\end{definition}

We can characterize constraint Poisson manifolds as Poisson manifolds with coisotropic submanifolds
and a compatible distribution containing the characteristic distribution.

\begin{proposition}
	\label{prop:CharacterizationConPoissonMfld}
	\index{coisotropic submanifold}
Let $\mathcal{M} = (M,C,D)$ be a constraint manifold
and $\pi \in \Secinfty(\Anti^2TM)$.
Then the following statements are equivalent:
\begin{propositionlist}
	\item \label{prop:CharacterizationConPoissonMfld_1}
	$(\mathcal{M}, \pi)$ is a constraint Poisson manifold.
	\item \label{prop:CharacterizationConPoissonMfld_2}
	$\{f,g\} \coloneqq \pi(\D f, \D g)$
	defines a constraint Poisson bracket on $\ConCinfty(\mathcal{M})$.
	\item \label{prop:CharacterizationConPoissonMfld_3}
	$(M,\pi)$ is a Poisson manifold,
	$C \subseteq M$ is a coisotropic submanifold with characteristic distribution $D_C \subseteq D$
	and
	\begin{equation}
		\Cinfty_D(M) = \{ f \in \Cinfty(M) \mid \Lie_X f\at{C} = 0 \text{ for all } X \in \Secinfty(D) \}
	\end{equation}
	is closed under the Poisson bracket.
\end{propositionlist}
\end{proposition}

\begin{proof}
We first show the equivalence of 
\ref{prop:CharacterizationConPoissonMfld_1}
and \ref{prop:CharacterizationConPoissonMfld_2}.
Thus assume $(\mathcal{M},\pi)$ is a constraint Poisson manifold.
Then $\{\argument,\argument\}$ is a Poisson bracket on 
$\ConCinfty(\mathcal{M})_\Total$ by classical results.
It remains to show that it is a constraint map.
For this recall that
$\pi \in \ConVecFields_{\strtensor}^2(\mathcal{M})_\Wobs$
and
\begin{equation*}
	\D \tensor \D \colon \ConCinfty(\mathcal{M}) \tensor 
	\ConCinfty(\mathcal{M}) \to \ConSecinfty(T^*\mathcal{M}) \tensor 
	\ConSecinfty(T^*\mathcal{M})
	= \ConSecinfty(T\mathcal{M} \strtensor T\mathcal{M})^*.
\end{equation*}
Since the evaluation is a constraint map we see that
$\{\argument,\argument\}$ is constraint.
On the other hand, if $\pi$ induces a constraint Poisson bracket,
then $\pi$ is a classical Poisson structure on $M$.
It remains to show that
$\pi \in \ConVecFields^2(\mathcal{M})_\Wobs$.
For this consider local adapted coordinates, such that
\begin{equation*}
	\pi\at{U} = \sum_{i,j=1}^{n_\Total} \pi^{ij} \frac{\del}{\del 
	x^i} \wedge \frac{\del}{\del x^j}.
\end{equation*}
Then we have $\{x^i,x^j\} = \pi^{ij}$
showing that \eqref{eq:ConditionsConBivectField_N}
holds, and therefore $\pi$ is a constraint Poisson bivector.

Next we show the equivalence of 
\ref{prop:CharacterizationConPoissonMfld_2}
and \ref{prop:CharacterizationConPoissonMfld_3}.
Assume $(\mathcal{M},\pi)$ is a constraint Poisson manifold.
Since $\pi \in \Anti^2 \Secinfty(TM)$ is a bivector field on $M$ with
$\Schouten{\pi,\pi} = 0$ it is a Poisson structure on $M$.
Moreover, for $f \in \vanishing_C$ we have 
\begin{equation*}
	X_f = \pi(\argument, \D f) = -\ins_{\D f} \pi \in \ConSecinfty(T 
	\mathcal{M})_\Null
\end{equation*}
by \autoref{prop:ConCartanCalculus}.
This shows $X_f(p) \in T_pC$ for $p \in C$, and thus $C \subseteq M$ 
is a coisotropic submanifold.
To show that $D$ is the corresponding characteristic distribution 
consider $p \in C$ and let
$(U,x)$ be an adapted chart around $p$ as in \autoref{lem:LocalStructureConManfifold}.
Since $\ConCinfty(\mathcal{M}\at{U})_\Null$ is generated by 
$x^{n_\Wobs+1}, \dotsc, x^{n_\Total}$
the characteristic distribution is spanned by
\begin{equation*}
	X_{x^i} = -\ins_{\D x^i} \pi\at{U} = \sum_{j=1}^{n_\Total} 
	\pi^{ij}\at{U} \frac{\del}{\del x^j} \in 
	\ConSecinfty(T\mathcal{M}\at{U})_\Null.
\end{equation*}
From \autoref{lem:LocalConBiVect} \ref{lem:LocalConBiVect_2} it 
follows that 
$X_{x^i}(p) = \sum_{j=1}^{n_\Null} \pi^{ij}(p) \frac{\del}{\del 
x^j}\at{p}$
and therefore the characteristic distribution is included in $D$.
The reverse implication follows from
\autoref{ex:ConBivectorsFromPoissonManifolds} \ref{ex:ConPoissonMfld}
and the fact that $\Cinfty_D(M) \subseteq \Cinfty_{D_C}(M)$ is a Poisson subalgebra.
\end{proof}

\begin{remark}
This result will have far reaching consequences for the deformation 
quantization of Poisson manifolds equipped with coisotropic 
submanifolds as considered in \autoref{sec:ConStarProducts}.
This will be discussed in more detail later on, but let us mention 
here that requiring
$\pi \in \ConVecFields_{\tensor}^2(\mathcal{M})_\Wobs$
instead of
$\pi \in \ConVecFields_{\strtensor}^2(\mathcal{M})_\Wobs$
would correspond to $C \subseteq M$ being a Poisson submanifold,
or equivalently $\ConCinfty(\mathcal{M})$ being
a strong constraint Poisson algebra.
Thus the choice of the tensor product, $\tensor$ or $\strtensor$ amounts to the choice
between a Poisson and a coisotropic submanifold.
\end{remark}

\subsubsection{Reduction}

Both types of constraint multivector fields are well behaved under reduction:

\begin{proposition}[Multivector fields vs. reduction]
	\label{prop:ConVecFieldsVSRedcution}
	\index{reduction!multivector fields}
Let $\mathcal{M} = (M,C,D)$ be a constraint manifold.
\begin{propositionlist}
	\item There exists a canonical isomorphism
	$\ConVecFields_{\tensor}^\bullet(\mathcal{M})_\red \simeq \VecFields^\bullet(\mathcal{M}_\red)$
	of DGLAs.
	\item There exists a canonical isomorphism
	$\ConVecFields_{\strtensor}^\bullet(\mathcal{M})_\red \simeq \VecFields^\bullet(\mathcal{M}_\red)$
	of DGLAs.
\end{propositionlist}

\end{proposition}

\begin{proof}
Similar to the proof of \autoref{prop:ConFormsVSReduction} this is a chain of canonical isomorphisms introduced before:
\begin{align*}
	\ConVecFields_{\tensor}^\bullet(\mathcal{M})_\red
	&= \big( \bigoplus_{k=0}^\infty \Anti_{\tensor}^k \ConSecinfty(T\mathcal{M}) \big)_\red
	\simeq \bigoplus_{k=0}^\infty (\Anti_{\tensor}^k \ConSecinfty(T\mathcal{M}))_\red \\
	&\simeq  \bigoplus_{k=0}^\infty \Anti^k \ConSecinfty(T\mathcal{M})_\red
	\simeq \bigoplus_{k=0}^\infty \Anti^k \Secinfty(T\mathcal{M}_\red)
	= \VecFields^k(\mathcal{M}_\red).
\end{align*}
Since the defining equation of the Schouten bracket holds for the reduced Schouten bracket we get an isomorphism of (differential) graded Lie algebras.
The second part follows since $\strtensor$ and $\tensor$ agree after reduction, and the Schouten bracket is given by the same formula.
\end{proof}

Since the reduced Schouten bracket is defined on representatives
we can infer that constraint Poisson manifolds reduce to classical Poisson manifolds.

\begin{corollary}
	\index{reduction!Poisson manifold}
Let $(\mathcal{M}, \pi)$ be a constraint Poisson manifold.
Then $(\mathcal{M}_\red,  \pi_\red)$ is a Poisson manifold.
\end{corollary}

\begin{example}
	\label{ex:ConBivectorsFromPoissonManifoldsReduction}
Let us revisit the examples of \autoref{ex:ConBivectorsFromPoissonManifolds}.
For this let $(M,\pi)$ be a Poisson manifold.
\begin{examplelist}
	\item For every closed Poisson submanifold $C \subseteq M$ the reduction
	of the constraint Poisson manifold $((M,C,0),\pi)$ is given by
	$(C, \pi\at{C})$.
	\item For every closed coisotropic submanifold $C \subseteq M$ the reduction
	of the constraint Poisson manifold 
	$((M,C,D),\pi)$ agrees with the classical coisotropic reduction.
	\item 
	Since every Poisson submanifold is in particular coisotropic we also get
	for every closed Poisson submanifold $C \subseteq M$ a constraint Poisson manifold
	$\mathcal{M} = (M,C,D)$ with $\pi \in \ConVecFields_{\strtensor}^2(\mathcal{M})$.
\end{examplelist}
\end{example}

Even though we can always reduce Poisson structures, it is not clear
that, in general, all Poisson structures on $\mathcal{M}_\red$
come from a constraint Poisson structure on $\mathcal{M}$,
since, even though we can always lift a bivector field to $\mathcal{M}$,
it is not obvious how it can be extended from $C$ to $M$ such that it
still fulfils $\Schouten{\pi,\pi}=0$,
see also \autoref{rem:ReductionConSets} \ref{rem:RedcutionOfEquations}.

%% file: constraint-symbolcalculus.tex
In this last section about constraint geometry we want to study (multi-)differential operators  on a manifold which are compatible with reduction, i.e. constraint (multi-)differential operators on constraint manifolds.
We start in \autoref{sec:ConDiffOps} by introducing algebraic constraint differential operators and study the particular case of constraint differential operators on sections of constraint vector bundles.
This will lead to a constraint leading symbol.
In order to find a full constraint symbol we define constraint covariant derivatives in \autoref{sec:ConCovDerivatives}, which we use in \autoref{sec:ConSymbCalculus} to establish a constraint symbol calculus.
Finally, \autoref{sec:ConMultiDiffops} is concerned with the generalization of the constraint symbol calculus to constraint multidifferential operators.
	
\subsection{Differential Operators}
\label{sec:ConDiffOps}

By an approach of Grothendieck, first introduced in \cite{grothendieck:1967a},
for a classical commutative algebra $\algebra{A}$ differential operators can be defined recursively as
\glsadd{DiffOp}
\begin{equation} \label{eq:RecursiveDefinitionDiffOp}
	\Diffop^k(\algebra{A})
	\coloneqq \left\{ D \in \End_\field{k}(\algebra{A}) \bigm|
	[L_a,D] \in \Diffop^{k-1}(\algebra{A}) \text{ for all } a \in \algebra{A} \right\},
\end{equation}
for $k \geq 0$ and
\begin{equation}
	\Diffop^{-1}(\algebra{A}) \coloneqq \{0\},
\end{equation}
where $L_a$ denotes the left multiplication with the fixed element
$a \in \algebra{A}$.

Instead of repeating the classical definitions internal to our categories
of constraint algebras and modules, let us directly give the following definition.

\begin{definition}[Constraint differential operators]
	\label{def:ConDiffop}
	\index{constraint!differential operators}
Let $\algebra{A} \in \injConAlg$ be a commutative embedded constraint algebra, and let
$\module{E}, \module{F}$ be embedded constraint
$\algebra{A}$-modules.
For $k \in \Integers$ we define
the \emph{constraint differential operators} as
\glsadd{ConDiffOp}
\begin{equation}
	\begin{split}
		\ConDiffop^k(\module{E}; \module{F})_\Total
		&\coloneqq \Diffop^k(\module{E}_\Total; \module{F}_\Total) \\
		\ConDiffop^k(\module{E}; \module{F})_\Wobs
		&\coloneqq \left\{ D \in \Diffop^k(\module{E}_\Total; \module{F}_\Total)
		\bigm| D \in \ConHom_\field{k}(\module{E}, \module{F})_\Wobs \right\}, \\
		\ConDiffop^k(\module{E}; \module{F})_\Null
		&\coloneqq \left\{ D \in \Diffop^k(\module{E}_\Total; \module{F}_\Total)
		\bigm| D \in \ConHom_\field{k}(\module{E}, 
		\module{F})_\Null \right\},
	\end{split}
\end{equation}
and write
\begin{equation}
	\ConDiffop^\bullet(\module{E};\module{F}) \coloneqq 
	\bigoplus_{k \in \Integers} 
	\ConDiffop^k(\module{E};\module{F}).
\end{equation}
\end{definition}

Note that $\ConDiffop^k(\module{E};\module{F})$,
and also $\ConDiffop^\bullet(\module{E};\module{F})$,
become strong constraint $\algebra{A}$-bimodules with 
respect to the classical $\algebra{A}_\Total$-bimodule structure given by
$(a \cdot D)(b) = a\cdot D(b)$
and
$(D \cdot a)(b) = D(a \cdot b)$.

Let us now focus on the case of differential operators on the sections of 
constraint vector bundles.
We will write $\ConDiffop^\bullet(E;F)$ instead of
$\ConDiffop^\bullet\big(\ConSecinfty(E); \ConSecinfty(F)\big)$
and $\ConDiffop^\bullet(\mathcal{M})$
for 
$\ConDiffop^\bullet\big(\ConCinfty(\mathcal{M});\ConCinfty(\mathcal{M})\big)$.

\begin{example}
	\label{ex:ParitalsAsConDiffOps}
Consider $\mathcal{M} = \Reals^n = 
(\Reals^{n_\Total},\Reals^{n_\Wobs}, \Reals^{n_\Null})$.
Then for any multi index $I = (i_1, \dotsc, i_r) \in \Naturals_0^r$
we write
\begin{equation}
	\del_I = \frac{\del^r}{\del x^{i_1} \cdots \del x^{i_r}}
	\in \ConDiffop^r(\Reals^n)_\Total.
\end{equation}
We have $\del_I \in \ConDiffop^r(\Reals^n)_\Wobs$
if and only if it only differentiates in direction of
the subspace $\Reals^{n_\Wobs}$, since then it preserves
$\ConCinfty(\Reals^{n})_\Wobs$
and $\ConCinfty(\Reals^n)_\Null$.
Similarly, we have $\del_I \in \ConDiffop^r(\Reals^n)_\Null$
if and only if it only differentiates in direction of 
$\Reals^{n_\Wobs}$ and at least once in direction of the distribution
$\Reals^{n_\Null}$.
In other words
\begin{equation}
	n^{\tensor r} \ni I \mapsto \del_I \in \ConDiffop^r(\Reals^n),
\end{equation}
with $n^{\tensor r}$ as defined in \autoref{def:TensorProdDualConIndSet},
is a constraint map.
\end{example}

This example leads to the following useful observation.

\begin{lemma}
	\label{lem:ConPartialDerivatives}
Let $E$ be a constraint vector bundle over a constraint manifold $\mathcal{M} = (M,C,D)$
of dimension $n = (n_\Total,n_\Wobs,n_\Null)$
and let $e_1, \dotsc, e_{\rank E} \in \ConSecinfty(E)_\Total$
be a constraint local frame.
For all $r \in \Naturals$ the following statements hold:
\begin{lemmalist}
	\item If $s \in \ConSecinfty(E)_\Wobs$, then the map
	\begin{equation}
		\varphi \colon n^{\tensor r} \tensor (\rank E)^*
		\ni (I,\alpha)
		\mapsto \del_I s^\alpha
		\in \ConCinfty(\mathcal{M}),
	\end{equation}
	with $s^\alpha = e^\alpha(s)$,
	is constraint, i.e. $\varphi \in \ConMap\big(n^{\tensor r} \tensor (\rank E)^*, \ConCinfty(\mathcal{M})\big)_\Wobs$.
	\item If $s \in \ConSecinfty(E)_\Null$, then it holds
	$\varphi \in \ConMap(n^{\tensor r} \tensor (\rank E)^*, \ConCinfty(\mathcal{M}))_\Null$.
\end{lemmalist}	
\end{lemma}

In this case we can locally characterize differential operators as follows.

\begin{proposition}[Local form of constraint differential operators]
	\label{prop:LocalConDiffop}
Let $E$ and $F$ be constraint vector bundles over a constraint manifold
$\mathcal{M} =(M,C,D)$ of dimension
$n = (n_\Total,n_\Wobs,n_\Null)$
and let
$D \in \ConDiffop^k(E; F)_\Total$, for $k \in \Naturals_0$.
Consider local adapted coordinates $(U,x)$ on $\mathcal{M}$ and let
$e_1, \dotsc, e_{\rank(E_\Total)} \in \Secinfty(E_\Total)$
be a constraint local frame.
Then
\begin{equation}
	\label{eq:LocalConDiffop}
	D\at{U}(s) = \sum_{r=0}^{k} 
	\sum_{\alpha=1}^{\rank(E_\Total)}\frac{1}{r!} D_{U,\alpha}^I 
	\cdot\del_I s^\alpha
\end{equation}
with $D_U^I \in \Secinfty(F_\Total\at{U})$
and $s^\alpha = e^\alpha(s)$.
\begin{propositionlist}
	\item It holds $D \in \ConDiffop^k(E; F)_\Wobs$ if and only if
	\begin{equation}
		\label{prop:LocalConDiffop_Wobs}
	\begin{split}
		D_{U,\alpha}^{I} \in \ConSecinfty(F)_\Wobs
		&\qquad\text{ if } (I,\alpha) \in \big((n^*)^{\strtensor r} \strtensor \rank E\big)_\Wobs\\
		D_{U,\alpha}^{I} \in \ConSecinfty(F)_\Null
		&\qquad\text{ if } (I,\alpha) \in \big((n^*)^{\strtensor r} \strtensor \rank E\big)_\Null,
	\end{split}
	\end{equation}
	\item It holds $D \in \ConDiffop^k(E;F)_\Null$ if and only if
	\begin{equation}
		\label{prop:LocalConDiffop_Null}
		D_{U,\alpha}^{I} \in \ConSecinfty(F)_\Null
		\qquad\text{ if } (I,\alpha) \in \big((n^*)^{\strtensor r} \strtensor \rank E\big)_\Wobs.
	\end{equation}
\end{propositionlist}
\end{proposition}

\begin{proof}
Evaluating $D\at{U}$ on
$x^{i_1}\cdots x^{i_r}\cdot e_\alpha$ yields
$D_{U,\alpha}^I$.
Now from \autoref{ex:ConFunctionsOnConMfld} 
\ref{ex:ConFunctionsOnConMfld_1}
it follows that
$x^{i_1}\cdots x^{i_r} \in \ConCinfty(\mathcal{M})_\Wobs$ if
$I \in ((n^*)^{\strtensor r})_\Wobs$
and that 
$x^{i_1}\cdots x^{i_r} \in \ConCinfty(\mathcal{M})_\Null$
if $I \in ((n^*)^{\strtensor r})_\Null$.
Moreover, since we use a constraint local frame we know
$e_\alpha \in \ConSecinfty(E)_\Wobs$
if and only if $\alpha \in \rank(E)_\Wobs$
and
$e_\alpha\in \ConSecinfty(E)_\Null$
if and only if $\alpha \in \rank(E)_\Null$.
Then for $D \in \ConDiffop^k(E;F)_\Wobs$ we immediately
get \eqref{prop:LocalConDiffop_Wobs}.
And similarly we obtain for $D \in \ConDiffop^k(E;F)_\Null$
directly \eqref{prop:LocalConDiffop_Null}.
For the other implication assume \eqref{prop:LocalConDiffop_Wobs} 
holds.
Let $s \in \ConSecinfty(E)_\Null$.
Then all terms of
\eqref{eq:LocalConDiffop}
end up in $\ConSecinfty(F)_\Null$:
By \autoref{lem:ConPartialDerivatives}
we have either $(I,\alpha) \in (n^{\tensor r} \tensor \rank E)_\Null$
and thus $\del_I s^\alpha \in \ConCinfty(\mathcal{M})_\Null$,
or $(I,\alpha) \in (n^{\tensor r} \tensor (\rank E)^*)^*_\Wobs
= ((n^*)^{\strtensor r} \strtensor \rank E)_\Wobs$
and thus $D^I_{U,\alpha} \in \ConSecinfty(F)_\Null$.
For
$s \in \ConSecinfty(E)_\Wobs$
we have $\del_I s^\alpha \in \ConCinfty(\mathcal{M})_\Null$
if 
$(I,\alpha) \in (n^{\tensor r} \tensor (\rank E)^*)_\Null$
and 
$\del_I s^\alpha \in \ConCinfty(\mathcal{M})_\Wobs$
if 
$(I,\alpha) \in (n^{\tensor r} \tensor (\rank E)^*)_\Wobs$.
Thus if 
$(I,\alpha) \in (n^{\tensor r} \tensor (\rank E)^*)_\Total 
\setminus (n^{\tensor r} \tensor (\rank E)^*)_\Wobs
= ((n^*)^{\strtensor r} \strtensor \rank E)_\Null$,
then $D^I_{U,\alpha} \in \ConSecinfty(F)_\Null$,
and if
$(I,\alpha) \in (n^{\tensor r} \tensor (\rank E)^*)_\Wobs
\setminus (n^{\tensor r} \tensor (\rank E)^*)_\Null
\subseteq ((n^*)^{\strtensor r} \strtensor \rank E)_\Wobs$,
then $D^I_{U,\alpha} \in \ConSecinfty(F)_\Wobs$.
This gives the first part.
The second part follows by completely analogous considerations.
\end{proof}

If $E = F = \mathcal{M} \times \Reals$ the local formula simplifies as follows.

\begin{corollary}
	\label{prop:LocalConDiffopCinfty}
Let $\mathcal{M} =(M,C,D)$ be a constraint manifold of dimension
$n = (n_\Total,n_\Wobs,n_\Null)$ and let
$D \in \ConDiffop^k(\mathcal{M})_\Total$, for $k \in \Naturals_0$.
Locally we can write
\begin{equation}
	D\at{U} = \sum_{r=0}^{k} 
	\frac{1}{r!} D_U^I \del_I
\end{equation}
with $D_U^I \in \ConCinfty(\mathcal{M})_\Total$.
\begin{corollarylist}
	\item It holds $D \in \ConDiffop^k(\mathcal{M})_\Wobs$ if and only if
	\begin{equation}
		\label{prop:LocalConDiffopCinfty_Wobs}
		\begin{split}
			D_U^I \in \ConCinfty(\mathcal{M})_\Wobs
			&\qquad\text{ if } I \in \left((n^*)^{\strtensor 
			r}\right)_\Wobs \\
			D_U^I \in \ConCinfty(\mathcal{M})_\Null
			&\qquad\text{ if } I \in \left((n^*)^{\strtensor 
			r}\right)_\Null.
		\end{split}
	\end{equation}
	\item It holds $D \in \ConDiffop^k(\mathcal{M})_\Null$ if and only if
	\begin{equation}
		\label{prop:LocalConDiffopCinfty_Null}
		\begin{split}
			D_U^I \in \ConCinfty(\mathcal{M})_\Null
			&\qquad\text{ if } I \in \left((n^*)^{\strtensor 
			r}\right)_\Wobs
		\end{split}
	\end{equation}
\end{corollarylist}
\end{corollary}

For every differential operator $D \in \Diffop^k(E;F)$ the classical leading symbol
\begin{equation}
	\sigma_k(D) \in \Secinfty(\Sym^k TM) \tensor E^* \tensor F
\end{equation}
is locally given by
\begin{equation}
	\sigma_k(D)\at{U} = \sum_{i_1,\dotsc,i_k=1}^{n_\Total} 
	\sum_{\alpha=1}^{\rank E} \frac{1}{k!} \frac{\del}{\del x^{i_1}} 
	\vee \dots \vee  \frac{\del}{\del x^{i_k}} \tensor e^\alpha 
	\tensor D_{U,\alpha}^{(i_1, \dotsc, i_k)}.
\end{equation}
For constraint differential operators this becomes a constraint section.

\begin{proposition}[Constraint leading symbol]
	\index{constraint!leading symbol}
Let $E$ and $F$ be constraint vector bundles over a constraint manifold
$\mathcal{M} = (M,C,D)$.
\begin{propositionlist}
	\item The leading symbol defines a constraint morphism
	\glsadd{LeadingSymbol}
	\begin{equation}
		\sigma_k \colon \ConDiffop^k(E;F) \to \ConSecinfty\big((\Sym^k_{\tensor} T\mathcal{M} \tensor E^*) \strtensor F\big)
	\end{equation}
	of strong constraint $\ConCinfty(\mathcal{M})$-modules. 
	\item If $E = F = \mathcal{M} \times \Reals$ the leading symbol becomes a constraint morphism
	\begin{equation}
		\sigma_k \colon \ConDiffop^k(\mathcal{M}) \to \ConSecinfty(\Sym^k_{\tensor} T\mathcal{M} )
	\end{equation}
	of strong constraint $\ConCinfty(\mathcal{M})$-modules.
\end{propositionlist}
\end{proposition}

\begin{proof}
	The $\TOTAL$-component of $\sigma_k$ is just the classical leading symbol.
	So it only remains to show that $\sigma_k$ is a constraint morphism.
	For this let $D \in \ConDiffop^k(E;F)_\Wobs$ be given.
	If $(I,\alpha) \in ((n^*)^{\strtensor k} \strtensor \rank E )_\Null$ we have
	$D^{(i_1, \dotsc, i_k)}_{U,\alpha} \in \ConSecinfty(F)_\Null$ by
	\autoref{prop:LocalConDiffop}.
	If 
	\begin{equation*}
		(I,\alpha) \in
		\big((n^*)^{\strtensor k} \strtensor \rank E \big)_\Wobs
		\setminus \big((n^*)^{\strtensor k} \strtensor \rank E 
		\big)_\Null
		\subseteq \big(n^{\tensor k} \tensor (\rank E)^*\big)_\Wobs
	\end{equation*}
	we have $\frac{\del}{\del x^{i_1}} 
	\vee \dots \vee  \frac{\del}{\del x^{i_k}} \tensor e^\alpha \in 
	\ConSecinfty(\Sym_{\tensor}^k T\mathcal{M} \tensor E^*)_\Wobs$
	and
	$D^{(i_1,\dotsc,i_k)}_{U,\alpha} \in \ConSecinfty(F)_\Wobs$.
	Moreover, for
	\begin{equation*}
		(I,\alpha) \in 
		\big((n^*)^{\strtensor k} \strtensor \rank E \big)_\Total 
		\setminus \big((n^*)^{\strtensor k} \strtensor \rank E 
		\big)_\Wobs
		= \big(n^{\tensor k} \tensor (\rank E)^*\big)_\Null
	\end{equation*}
	we obtain
	$\frac{\del}{\del x^{i_1}} \vee \dots \vee  \frac{\del}{\del x^{i_k}} \tensor e^\alpha
	\in \ConSecinfty(\Sym_{\tensor}^k T\mathcal{M} \tensor E^*)_\Null$.
	Thus $\sigma_k$ preserves the $\WOBS$-component.
	Let now $D \in \ConDiffop^k(E;F)_\Null$ be given.
	Then for $(I,\alpha) \in ((n^*)^{\strtensor k} \strtensor \rank E )_\Wobs$ we have
	$D^{(i_1, \dotsc, i_k)}_{U,\alpha} \in \ConSecinfty(F)_\Null$ by
	\autoref{prop:LocalConDiffop}, and 
	for
	\begin{equation*}
		(I,\alpha) \in 
		\big((n^*)^{\strtensor k} \strtensor \rank E \big)_\Total
		\setminus \big((n^*)^{\strtensor k} \strtensor \rank E 
		\big)_\Wobs
		= \big(n^{\tensor k} \tensor (\rank E)^*\big)_\Null
	\end{equation*}
	we get
	$\frac{\del}{\del x^{i_1}} \vee \dots \vee
	\frac{\del}{\del x^{i_k}} \tensor e^\alpha
	\in \ConSecinfty(\Sym_{\tensor}^k T\mathcal{M}
	\tensor E^*)_\Null$
	as before.
	The second part is just a special case of the first.
\end{proof}

Restricting the leading symbol to
$\ConDer(\ConCinfty(\mathcal{M}))$
gives the inverse
\begin{equation}
	\sigma_1\at{\ConDer\big(\ConCinfty(\mathcal{M})\big)}
	\colon \ConDer(\ConCinfty(\mathcal{M}))
	\to \ConSecinfty(T\mathcal{M})
\end{equation}
of the Lie derivative, see \autoref{prop:ConVectorFieldsAreConDer}.
Observe, that the local formula for constraint differential operators recovers the local formulas 
of constraint vector fields from \autoref{lem:LocalConVect}.

\subsubsection{Reduction}

Let us first consider the reduction of constraint differential operators on
general constraint modules.

\begin{proposition}[Constraint differential operators vs. reduction]
	\label{prop:ConDiffopsVSReduction}
	\index{reduction!differential operators}
Let $\algebra{A} \in \injConAlg$ be a commutative embedded constraint 
algebra, and let
$\module{E}, \module{F}$ be embedded constraint
$\algebra{A}$-modules.
For each $k \in \Naturals_0$ there is a natural injective morphism
\begin{equation} \label{eq:RedConDiffOp}
	\ConDiffop^k\left(\module{E};\module{F}\right)_\red 
	\to \Diffop^k\left(\module{E}_\red; \module{F}_\red\right)
\end{equation}
of $\algebra{A}_\red$-modules.
\end{proposition}

\begin{proof}
Since by definition we have
$\ConDiffop^k(\module{E};\module{F}) \subseteq 
\ConHom_{\field{k}}(\module{E},\module{F})$
and
$\ConHom_\field{k}(\module{E},\module{F})_\red
\simeq \Hom_\field{k}(\module{E}_\red,\module{F}_\red)$
by \autoref{prop:ReductionOnCModk}
we obtain an injective morphism
$\ConDiffop^k(\module{E};\module{F})_\red \to 
\Hom_\field{k}(\module{E}_\red,\module{F}_\red)$.
The recursive condition in
\eqref{eq:RecursiveDefinitionDiffOp} still holds after reduction,
cf. \autoref{rem:ReductionConSets} \ref{rem:RedcutionOfEquations},
which shows that we obtain the required morphism.
\end{proof}

Note again that we can in general not expect the morphism
\eqref{eq:RedConDiffOp} to be an isomorphism, cf. \autoref{rem:ReductionConSets} \ref{rem:RedcutionOfEquations}.
Let us now take a look at reduction of constraint differential operators of sections:

\begin{proposition}[Constraint differential operators of sections vs. reduction]{\ \\} 
	\label{prop:ConDiffopsSecVSReduction}
Let $E$ and $F$ be constraint vector bundles over a constraint manifold 
$\mathcal{M} = (M,C,D)$.
\begin{propositionlist}
	\item Let $D \in \ConDiffop^k(E; F)_\Wobs$, then locally
	\begin{equation}
		\label{eq:ReducedLocConDiffop}
		\big(D_\red\at{U}\big)(s) = \sum_{r=0}^{k} \sum_{n_0 < I \leq 
		n_\Wobs} \frac{1}{r!} (D_{U,\alpha}^I)_\red
		\cdot \del_I s^\alpha,
	\end{equation}
	for $s \in \Secinfty(E_\red)$.
	\item The constraint leading symbol $\sigma_k$ induces the classical leading symbol
	\begin{equation}
		(\sigma_k)_\red \colon \Diffop^k(E_\red; F_\red) \to \Secinfty(\Sym^kT\mathcal{M}_\red \tensor E_\red ^* \tensor F_\red)
	\end{equation}
	on the reduced manifold $\mathcal{M}_\red$.
\end{propositionlist}
\end{proposition}

\begin{proof}
Let $\check{e}_1,\dotsc,\check{e}_{\rank E_\red} \in 
\ConSecinfty(E_\red\at{U})$
be a local frame.
Then by the same construction as in the proof of
\autoref{lem:AdaptedLocalFrames}
we obtain a constraint local frame
$e_1,\dotsc, e_{\rank E_\Total} \in 
\ConSecinfty(E_\Total\at{\pi_\red^{-1}(U)})$.
Moreover, from \autoref{prop:ConSecVSReduction} we know that there 
exists
$\hat{s} \in \ConSecinfty(E)_\Wobs$ such that
$s = [\hat{s}]$.
Then it follows from \autoref{prop:LocalConDiffop} 
that 
\begin{equation*}
	(D_\red\at{U})(s)
	= \big(D\at{\pi_\red^{-1}(U)}(\hat{s})\big)_\red
	= \sum_{r=0}^{k} 
	\sum_{\alpha=1}^{\rank(E_\Total)}\frac{1}{r!} 
	(D_{U,\alpha}^I)_\red
	\cdot[\del_I \hat{s}^\alpha]
\end{equation*}
holds.
From \autoref{prop:LocalConDiffop}
we also know that $(D^I_{U,\alpha})_\red = 0$
if $(I,\alpha) \in ((n^*)^{\strtensor r} \strtensor \rank E)_\Null$.
If $(I,\alpha) \in (n^{\tensor r} \tensor (\rank E)^*)_\Null$
then $[\del_I \hat{s}^\alpha] = 0$
by \autoref{lem:ConPartialDerivatives}.
The remaining summands give \eqref{eq:ReducedLocConDiffop}.
Since the leading symbol is characterized by the highest order terms 
of the local expression it follows from \eqref{eq:ReducedLocConDiffop}
that $(\sigma_k)_\red$ is indeed the leading symbol for the reduced 
differential operators.
\end{proof}

\subsection{Covariant Derivatives}
\label{sec:ConCovDerivatives}

We introduce covariant derivatives by copying the classical definition.

\begin{definition}[Constraint covariant derivative]
	\label{def:ConCovariantDerivative}
	\index{constraint!covariant derivative}
	\glsadd{CovDerivative}
Let $E = (E_\Total,E_\Wobs,E_\Null)$ be a constraint vector bundle over a constraint manifold
$\mathcal{M} = (M,C,D)$.
A \emph{constraint covariant derivative} is a morphism 
\begin{equation}
	\nabla \colon \ConSecinfty(T \mathcal{M}) \tensor[\Reals] \ConSecinfty(E)
	\to \ConSecinfty(E)
\end{equation}
of constraint $\Reals$-modules such that
\begin{align}
		\label{def:ConCovariantDerivative_1}
	\nabla_{fX}s = f \nabla_X s
	\shortintertext{and}
		\label{def:ConCovariantDerivative_2}
	\nabla_X{fs} = (\Lie_X f) s + f \nabla_Xs
\end{align}
for all $X \in \ConSecinfty(T\mathcal{M})_\Total$,
$s \in \ConSecinfty(E)_\Total$
and $f \in \ConCinfty(\mathcal{M})_\Total$.
\end{definition}

Condition \eqref{def:ConCovariantDerivative_1}
could be rephrased as saying that $\nabla$ is a left
$\ConCinfty(\mathcal{M})$-module morphism.

\begin{remark}
	\label{rem:DefinitionConCovDer}
The question arises why we use $\tensor[\Reals]$ and not $\strtensor[\Reals]$
in the definition of constraint covariant derivatives.
One way to answer this is by observing that the Lie derivative
$\Lie \colon \ConSecinfty(T\mathcal{M}) \tensor[\field{k}] \ConCinfty(\mathcal{M}) \to \ConCinfty(\mathcal{M})$
is not well-defined if we would use $\strtensor[\field{k}]$ instead of $\tensor[\field{k}]$,
and hence \eqref{def:ConCovariantDerivative_2} would cause problems.
Another justification comes from the fact that we can rephrase a classical covariant derivative as a map
$\nabla \colon \Secinfty(E) \to \Hom_\field{k}(\Secinfty(TM), \Secinfty(E))$.
Using this as a starting point, we could define a constraint covariant derivative by a constraint map
$\nabla \colon \ConSecinfty(E) \to \ConHom_{\field{k}}(\ConSecinfty(T\mathcal{M}),\ConSecinfty(E))$.
Using \autoref{prop:DualTensorHomIsos} 
this translates to
\begin{equation}
	\nabla \colon \ConSecinfty(E) \to \ConSecinfty(T^*\mathcal{M}) \strtensor[\field{k}] \ConSecinfty(E)
\end{equation}
and applying
\autoref{cor:StrHomOfTensor}
yields our definition of constraint covariant derivative using
$\tensor[\field{k}]$.
\end{remark}

\begin{corollary}
	\label{cor:CharacterizationConCovDer}
Let $E = (E_\Total,E_\Wobs,E_\Null)$ be a constraint vector bundle over a constraint
manifold $\mathcal{M} = (M,C,D)$.
Let $\nabla$ be a covariant derivative on $E_\Total$.
Then $\nabla$ is a constraint covariant derivative on $E$ if and only if the following properties hold:
\begin{corollarylist}
	\item $\nabla_X s \in \ConSecinfty(E)_\Wobs$ for all $X \in \ConSecinfty(T \mathcal{M})_\Wobs$
	and $s \in \ConSecinfty(E)_\Wobs$.
	\item $\nabla_X s \in \ConSecinfty(E)_\Null$ for all $X \in \ConSecinfty(T \mathcal{M})_\Wobs$
	and $s \in \ConSecinfty(E)_\Null$.
	\item $\nabla_X s \in \ConSecinfty(E)_\Null$ for all $X \in \ConSecinfty(T \mathcal{M})_\Null$
	and $s \in \ConSecinfty(E)_\Wobs$.
\end{corollarylist}
\end{corollary}

\begin{example}
	\label{ex:CanonicalConCovDerivative}
	\index{trivial constraint vector bundle}
Let $E = \mathcal{M} \times \Reals^k$ be a trivial constraint vector bundle
over a constraint manifold $\mathcal{M} = (M,C,D)$ as
in \autoref{ex:ConVectBundles}.
By \autoref{prop:SectionsTrivStrConVectorBundle} we know that
$\ConSecinfty(E) \simeq \ConCinfty(\mathcal{M})^k$
is a free strong constraint module.
The componentwise Lie derivative then defines a constraint covariant
derivative on $E$.
\end{example}

This example also shows the local existence of constraint covariant derivates.
Global existence can be shown using the constraint \hyperref[thm:strConSerreSwan]{Serre-Swan Theorem}:

\begin{proposition}[Existence of constraint covariant derivatives]
On any constraint
vector bundle 
$E = (E_\Total,E_\Wobs,E_\Null,\nabla^E)$
over a constraint manifold $\mathcal{M} = (M,C,D)$
exists a constraint covariant derivative.
\end{proposition}

\begin{proof}
By \autoref{thm:strConSerreSwan} we know that
$\ConSecinfty(E)$ is finitely generated projective.
Let $\{e_i,e^i\}_{i \in n}$ be a constraint dual basis of $\ConSecinfty(E)$
as in \autoref{prop:StrDualBasis}, then every $s \in \ConSecinfty(E)_\Total$ can be written as
$s = \sum_{i=1}^{n_\Total} e^i(s) e_i$.
We define
\begin{equation*} 
	\nabla_X s \coloneqq \sum_{i=1}^{n_\Total} \Lie_X(e^i(s)) e_i
\end{equation*} 
for every $X \in \Secinfty(TM)$ and $s \in \Secinfty(E_\Total)$.
An easy computation shows that $\nabla$ defines indeed a covariant derivative on $E_\Total$.
We still need to show that $\nabla$ is compatible with the constraint structure.
By \autoref{prop:StrDualBasis} we know that
\begin{equation*}\label{eq:CovDerProjective}
	\nabla_X s 
	= \sum_{i=1}^{n_\Null} \Lie_X(e^i(s)) \cdot\mkern-18mu\underbrace{e_i}_{\in \ConSecinfty(E)_\Null}
	+\sum_{i=n_\Null +1}^{n_\Wobs} \Lie_X(e^i(s)) \cdot\mkern-18mu\underbrace{e_i}_{\in \ConSecinfty(E)_\Wobs}
	+\sum_{i=n_\Wobs + 1}^{n_\Total} \Lie_X(e^i(s)) \cdot\mkern-18mu\underbrace{e_i}_{\in \ConSecinfty(E)_\Total}\mkern-30mu.
	\tag{$*$}
\end{equation*}
Thus the first term in \eqref{eq:CovDerProjective} is always in $\ConSecinfty(E)_\Null$.
Now let $X \in \ConSecinfty(T\mathcal{M})_\Wobs$
be given.
Again by \autoref{prop:StrDualBasis} we get the following case by case study:
\begin{cptitem}
	\item For $s \in \ConSecinfty(E)_\Wobs$ we have $e^i(s) \in \ConCinfty(\mathcal{M})_\Wobs$
	and hence $\Lie_X(e^i(s)) \in \ConCinfty(\mathcal{M})_\Wobs$ for all $i = n_\Null +1,\dotsc, n_\Wobs$.
	Hence the second term of \eqref{eq:CovDerProjective} is in $\ConSecinfty(E)_\Wobs$.
	Moreover, for $i = n_\Wobs +1, \dotsc, n_\Total$ we have $e^i(s) \in \ConCinfty(\Null)$
	and hence $\Lie_X(e^i(s)) \in \ConCinfty(\Null)$.
	Therefore also the third term is in $\ConSecinfty(E)_\Wobs$.
	\item If $s \in \ConSecinfty(E)_\Null$, then $e^i(s) \in \ConCinfty(\mathcal{M})_\Null$ for all
	$i = n_\Null + 1, \dotsc, n_\Total$.
	Thus both the second and third term of \eqref{eq:CovDerProjective} are elements in
	$\ConSecinfty(E)_\Wobs$.
\end{cptitem}
Suppose $X \in \ConSecinfty(T\mathcal{M})_\Null$.
\begin{cptitem}
	\item For all $s \in \ConSecinfty(E)_\Wobs$ we have $e^i(s) \in \ConCinfty(\mathcal{M})_\Null$
	and hence $\Lie_X(e^i(s)) \in \ConCinfty(\mathcal{M})_\Null$, showing that
	\eqref{eq:CovDerProjective} ends up in $\ConSecinfty(E)_\Null$.
\end{cptitem}
Thus $\nabla$ is a constraint covariant derivative.
\end{proof}

\begin{proposition}
Let $E = (E_\Total,E_\Wobs,E_\Null, \nabla^E)$
be a constraint vector bundle over a constraint manifold
$\mathcal{M} = (M,C,D)$.
\begin{propositionlist}
	\item If $\nabla$ and $\tilde{\nabla}$ are constraint covariant
	derivatives for $E$ then,
	$\nabla - \tilde{\nabla}$ is $\ConCinfty(\mathcal{M})$-bilinear, 
	hence
	$\nabla - \tilde{\nabla} \in 
	\ConSecinfty(T^*\mathcal{M} \strtensor \ConEnd(E))_\Wobs$.
	\item If $\nabla$ is a constraint covariant derivative
	on $E$ and $A \in \ConSecinfty(T^*\mathcal{M} \strtensor 
	\ConEnd(E))_\Wobs$, then
	\begin{equation}
		\tilde{\nabla}_Xs\coloneqq \nabla_X s + A(X \tensor s),
	\end{equation}
	with $X \in \ConSecinfty(T\mathcal{M})_\Total$,
	$s \in \ConSecinfty(E)_\Total$,
	defines another constraint covariant derivative on $E$.
\end{propositionlist}
\end{proposition}

\begin{proof}
	For the first part, a quick check or the well-known classical 
	statement shows that $\nabla - \tilde{\nabla}$ is bilinear,
	i.e. $\nabla - \tilde{\nabla} \in \ConSecinfty(T\mathcal{M}) 
	\tensor \ConSecinfty(E) \to \ConSecinfty(E)$.
	By \autoref{cor:StrHomOfTensor} this is equivalently given by an 
	element in
	$\ConSecinfty(T^*\mathcal{M} \strtensor \ConEnd(E))$.
	For the second part note that
	\begin{equation*}
		\ConSecinfty\big(T^*\mathcal{M} \strtensor \ConEnd(E)\big)
		\simeq \ConSecinfty(T\mathcal{M})^* \strtensor 
		\ConEnd_{\ConCinfty(\mathcal{M})}\big(\ConSecinfty(E)\big),
	\end{equation*}
	thus the evaluation of $A$ on $X \tensor s$ is indeed a 
	constraint morphism by \eqref{eq:EvaluationMap}, showing that $\tilde{\nabla}$
	is a constraint covariant derivative.
\end{proof}

The above shows that the set of constraint covariant derivatives 
forms an affine space over
$\ConSecinfty\big(T^*\mathcal{M} \strtensor \ConEnd(E)\big)_\Wobs$.
\begin{remark}
Even though we have not formally introduced constraint 
affine spaces, it becomes clear that the constraint set of 
covariant derivatives on $E$ is a constraint affine space
over $\ConSecinfty(T^*\mathcal{M} \strtensor \ConEnd(E))$.
In particular $\nabla$ and $\tilde{\nabla}$ are equivalent, if and 
only if $\nabla - \tilde{\nabla} \in \ConSecinfty(T^*\mathcal{M} 
\strtensor \ConEnd(E))_\Null$.
\end{remark}

\begin{proposition}
	\label{prop:ConstructionsConCovDer}
Let $E = (E_\Total,E_\Wobs,E_\Null)$ and
$F = (F_\Total,F_\Wobs,F_\Null)$
be constraint vector bundles over a constraint manifold
$\mathcal{M} = (M,C,D)$.
\begin{propositionlist}
	\item \label{prop:ConstructionsConCovDer_1}
	\index{dual!constraint covariant derivative}
	Suppose $\nabla$ is a constraint covariant derivative on $E$, then
	$\nabla^*$ defined by
	\begin{equation}
		(\nabla^*_X \alpha) (s) \coloneqq \Lie_X \big(\alpha(s)\big) - \alpha(\nabla_X s),
	\end{equation}
	for $X \in \ConSecinfty(T\mathcal{M})_\Total$,
	$\alpha \in \ConSecinfty(E^*)$ and 
	$s \in \ConSecinfty(E)_\Total$,
	defines a constraint covariant derivative on $E^*$.
	\item \label{prop:ConstructionsConCovDer_2}
	\index{tensor product!constraint covariant derivatives}
	Suppose $\nabla^E$ and $\nabla^F$ are constraint covariant derivative on $E$ and $F$, respectively.
	Then
	$\nabla^{E \tensor F}$ defined by
	\begin{equation}
		\nabla^{E \tensor F}_X (s \tensor t)
		\coloneqq (\nabla^E_X s) \tensor t + s \tensor (\nabla^F_X t)
	\end{equation}
	for $X \in \ConSecinfty(T\mathcal{M})_\Total$, 
	$s \in \ConSecinfty(E)_\Total$ 
	and
	$t \in \ConSecinfty(F)_\Total$
	defines a constraint covariant derivative on $E \tensor F$.
	\item \label{prop:ConstructionsConCovDer_3}
	\index{strong tensor product!constraint covariant derivatives}
	Suppose $\nabla^E$ and $\nabla^F$ are constraint covariant derivative on $E$ and $F$, respectively.
	Then
	$\nabla^{E \strtensor F}$ defined by
	\begin{equation}
		\nabla^{E \strtensor F}_X (s \tensor t)
		\coloneqq (\nabla^E_X s) \tensor t + s \tensor (\nabla^F_X t)
	\end{equation}
	for $X \in \ConSecinfty(T\mathcal{M})_\Total$, 
	$s \in \ConSecinfty(E)_\Total$ 
	and
	$t \in \ConSecinfty(F)_\Total$
	defines a constraint covariant derivative on $E \strtensor F$.
\end{propositionlist}
\end{proposition}

\begin{proof}
	On the $\TOTAL$-components these constructions are just
	given by the usual canonical constructions for covariant
	derivatives on $E^*_\Total$ and
	$E_\Total \tensor F_\Total$.
	Thus a straightforward check of the three properties
	from \autoref{cor:CharacterizationConCovDer} shows that these are indeed constraint
	covariant derivatives.
\end{proof}

By \autoref{rem:DefinitionConCovDer} a constraint covariant derivative
$\nabla^E$ on a constraint vector bundle $E$ can be understood
as a constraint map
$\nabla \colon \ConSecinfty(E) \to \ConSecinfty(T^* \mathcal{M}) \strtensor[\field{k}] \ConSecinfty(E)$.
If we additionally choose a constraint covariant derivative on $T\mathcal{M}$,
then by \autoref{prop:ConstructionsConCovDer} \ref{prop:ConstructionsConCovDer_3}
we obtain a constraint covariant derivative on
$T^*\mathcal{M} \strtensor E$.
Thus we obtain an iterated covariant derivative
\begin{equation}
	\underbrace{\nabla \circ \dots \circ \nabla}_{k\text{-times}} \colon \ConSecinfty(E) \to \ConSecinfty(T^*\mathcal{M})^{\strtensor k} \strtensor[\field{k}] \ConSecinfty(E).
\end{equation}
Symmetrizing yields the following notion of symmetrized covariant derivative.

\begin{definition}[Symmetrized constraint covariant derivative]
	\index{symmetrized constraint covariant derivative}
Let $E = (E_\Total,E_\Wobs,E_\Null)$ be a constraint vector bundle 
over a constraint 
manifold $\mathcal{M} = (M,C,D)$.
Moreover, let $\nabla^E$ and $\nabla$ be constraint covariant 
derivatives
on $E$ and $T\mathcal{M}$, respectively.
The constraint morphism
\begin{equation}
	\SymCoDer^E \colon \ConSecinfty(\Sym^\bullet_{\strtensor} 
	T^*\mathcal{M} \strtensor E)
	\to \ConSecinfty(\Sym^{\bullet+1}_{\strtensor} T^*\mathcal{M} 
	\strtensor E)
\end{equation}
defined by
\begin{equation}
\begin{split}
	(\SymCoDer^E \alpha)(X_0,\dotsc, X_k)
	\coloneqq &\sum_{i=0}^{k} \nabla_{X_i}^E
	\big( \alpha(X_0, \dotsc, \overset{i}{\wedge}, \dotsc, X_k) \big)\\
	&- \sum_{i\neq j} 
	\alpha \big( \nabla_{X_i}X_j, X_0,\dotsc, 
	\overset{i}{\wedge},\dotsc, \overset{j}{\wedge},\dotsc, X_k 
	\big),
\end{split}
\end{equation}
for $\alpha \in \ConSecinfty(\Sym^k_{\strtensor} T^*\mathcal{M} 
\strtensor E)_\Total$
and $X_0, \dotsc, X_k \in \ConSecinfty(T\mathcal{M})_\Total$,
is called \emph{symmetrized constraint covariant derivative}.
\end{definition}

Since $\SymCoDer^E$ is defined as a composition of constraint maps, it is itself 
a constraint morphism.

If we consider the trivial bundle $E = \mathcal{M} \times \mathbb{R}$ with its canonical constraint covariant derivative from
\autoref{ex:CanonicalConCovDerivative}, then we denote the symmetrized covariant derivative corresponding
to a constraint covariant derivative $\nabla$ for $T\mathcal{M}$ simply by
\begin{equation}
	\SymCoDer \colon \ConSecinfty(\Sym^\bullet_{\strtensor} T^*\mathcal{M}) \to \ConSecinfty(\Sym^{\bullet+1}_{\strtensor}T^*\mathcal{M}).
\end{equation}

\subsubsection{Reduction}

As expected, the reduction of a constraint covariant derivative yields a covariant derivative
on the reduced bundle.

\begin{proposition}[Covariant derivative vs. Reduction]
	\label{prop:CovDerVSReduction}
	\index{covariant derivative}
Let $E$ be a constraint vector bundle
over
$\mathcal{M} = (M,C,D)$.
Moreover, let $\nabla$ be a constraint covariant derivative on $E$.
\begin{propositionlist}
	\item \label{prop:CovDerVSReduction_CovDer}
	The reduction $\nabla_\red \colon 
	\Secinfty(T\mathcal{M}_\red) 
	\tensor[\field{k}] \Secinfty(E_\red) \to \Secinfty(E_\red)$
	defines a covariant derivative on $E_\red$.
	\item \label{prop:CovDerVSReduction_Dual}
	For the dual covariant derivative it holds
	$(\nabla^*)_\red = (\nabla_\red)^*$.
	\item \label{prop:CovDerVSReduction_Tensor}
	If $F$ is another constraint vector bundle over 
	$\mathcal{M}$ with covariant derivative $\tilde{\nabla}$, we get
	\begin{equation}
		(\nabla^{E\tensor F})_\red = \nabla^{E_\red \tensor F_\red}
		= (\nabla^{E \strtensor F})_\red.
	\end{equation}
	\item \label{prop:CovDerVSReduction_SymCoDer}
	For the symmetrized constraint covariant derivative
	it holds $(\SymCoDer^E)_\red = \SymCoDer^{E_\red}$.
\end{propositionlist}
\end{proposition}

\begin{proof}
	Since taking sections commutes with reduction by
	\autoref{prop:ConSecVSReduction},
	all reduced maps are defined for the correct domains and 
	codomains.
	The defining equations for all involved morphism carry over to 
	the reduction by \autoref{rem:ReductionConSets} \ref{rem:RedcutionOfEquations}.
\end{proof}


\subsection{Symbol Calculus}
\label{sec:ConSymbCalculus}

We define
\index{insertion}
\begin{align}
	\ins_s \colon \ConSecinfty((\Sym_{\tensor}^k T\mathcal{M} \tensor E^*) \strtensor F)
	\tensor
	\ConSecinfty(\Sym_{\strtensor}^\ell T^*\mathcal{M} \strtensor E)
	\to \ConSecinfty(\Sym_{\tensor}^{k-\ell} T\mathcal{M} \strtensor F)
\end{align}
by using
$\ConSecinfty(\Sym_{\strtensor}^\ell T^*\mathcal{M} \strtensor E)
\simeq \left( \ConSecinfty(\Sym_{\tensor}^k T\mathcal{M} \tensor E^*) \right)^*$.
More precisely, on factorizing tensors we have
\begin{equation}
	\ins_s(X \tensor \alpha \tensor t)(\omega \tensor s) = \alpha(s) \cdot \ins_s(X)(\omega) \tensor t,
\end{equation}
where $X \in \ConSecinfty(\Sym^k_{\tensor} T\mathcal{M})_\Total$,
$\alpha \in \ConSecinfty(E^*)_\Total$,
$t \in \ConSecinfty(F)_\Total$
$\omega \in \ConSecinfty(\Sym^\ell_{\strtensor} T^*\mathcal{M})_\Total$
and
$s \in \ConSecinfty(E)_\Total$.
This can then be extended to a constraint morphism
\begin{equation}
	\ins_s \colon \ConSecinfty((\Sym_{\tensor}^\bullet T\mathcal{M} \tensor E^*) \strtensor F)
	\tensor
	\ConSecinfty(\Sym_{\strtensor}^\bullet T^*\mathcal{M} \strtensor E)
	\to \ConSecinfty(\Sym_{\tensor}^\bullet T\mathcal{M} \strtensor F),
\end{equation}
see also \eqref{eq:EvaluationMap} for the evaluation morphism for constraint modules.
With this, and with the help of the symmetrized constraint covariant derivative we can now introduce the full constraint symbol:

\begin{theorem}[Constraint symbol calculus]
	\label{thm:ConSymbolCalculus}
	\index{constraint!symbol calculus}
Let $E$ and $F$ be constraint vector bundles over a constraint 
manifold $\mathcal{M} = (M,C,D)$.
Moreover, let $\nabla^E$ be a constraint covariant derivative on $E$ and 
$\nabla$ a constraint covariant derivative on $T\mathcal{M}$.
\begin{theoremlist}
	\item \index{constraint!symbol map}
	Then
	\begin{equation}
		\Op \colon \ConSecinfty\big((\Sym_{\tensor}^\bullet T\mathcal{M} \tensor E^*) \strtensor F)\big)
		\to \ConDiffop^\bullet(E; F)
	\end{equation}
	defined by
	\begin{equation}
		\Op(X)s \coloneqq \frac{1}{k!} \ins_s(X) (\SymCoDer^E)^ks,
	\end{equation}
	for $X \in \ConSecinfty((\Sym^k_{\tensor}T\mathcal{M} \tensor E^*) \strtensor F))$
	and $s \in \ConSecinfty(E)$,
	is a morphism of strong constraint $\ConCinfty(\mathcal{M})$-modules.
	\item For $X \in \ConSecinfty((\Sym^k_{\tensor} T\mathcal{M} \tensor E^*) \strtensor F))_\Total$ we have
	\begin{equation}
		\sigma_k\big(\Op(X)\big) = X,
	\end{equation}
	where $\sigma_k$ denotes the leading symbol as usual.
	\item $\Op$ is an isomorphism of strong constraint $\ConCinfty(\mathcal{M})$-modules.
\end{theoremlist}

\end{theorem}

\begin{proof}
By the classical theory we know that $\Op$ fulfils all the above 
properties on the $\TOTAL$-compo\-nent.
For the first part note that $\Op$ is defined as a composition of 
constraint
morphism, and thus defines itself a constraint morphism to
$\ConHom_\field{k}(\ConSecinfty(E), \ConSecinfty(F))$.
Since we know that $\Op(X)$ is a differential operator
it follows that $\Op$ actually maps to the constraint submodule
$\ConDiffop^\bullet(E;F)$
of $\ConHom_\field{k}(\ConSecinfty(E), \ConSecinfty(F))$.
The second part is just the classical statement.
Nevertheless, from this follows directly that $\Op$ is a monomorphism.
To show that it is also a regular epimorphism we repeat the classical
argument for constructing preimages.
Let $D \in \ConDiffop^k(\ConSecinfty(E); \ConSecinfty(F))_\Total$ be 
given.
Then $D - \Op(\sigma_k(D))$ has order $k-1$.
We write $X_k \coloneqq \sigma_k(D)$, then by induction we obtain
$D = \Op(X_k + \dots + X_0)$,
and thus $\Op$ is surjective.
If $D \in \ConDiffop^k(E; F)_\Wobs$,
then we have
$X_k \in \ConSecinfty((\Sym^k_{\tensor}T\mathcal{M} \tensor E^*) \strtensor F))_\Wobs$
and therefore $\Op$ is also surjective on the $\WOBS$-component.
Similarly, using $D \in \ConDiffop^k(E; F)_\Null$
we get that $\Op$ is indeed a regular epimorphism.
Hence $\Op$ is a regular epimorphism and monomorphism, and therefore
a constraint isomorphism.
\end{proof}

For $E = F = \mathcal{M} \times \Reals$ we immediately get
the following isomorphism for differential operators of
$\ConCinfty(\mathcal{M})$.

\begin{corollary}
	Let $\mathcal{M} = (M,C,D)$ be a constraint manifold
	and let $\nabla$ be a constraint covariant derivative for
	$T\mathcal{M}$.
	Then 
	\begin{align}
		\Op &\colon \ConSecinfty(\Sym^\bullet_{\tensor}T\mathcal{M})
		\to \ConDiffop^\bullet(\mathcal{M}),\\
		\shortintertext{with}
		\Op(X)(f) &= \frac{1}{k!}\ins_s(X) \SymCoDer^kf
		\quad\text{for}\quad
		X \in \ConSecinfty(S^k_{\tensor}T\mathcal{M})
	\end{align}
	is an isomorphism of constraint $\ConCinfty(\mathcal{M})$-modules.
\end{corollary}

\subsubsection{Reduction}

It turns out that the reduction of the constraint full symbol map yields
the full symbol on the reduced vector bundles:

\begin{proposition}
	\index{reduction!symbol map}
	\index{reduction!symbol calculus}
Let $E$ and $F$ be constraint vector bundles over a constraint 
manifold
$\mathcal{M} = (M,C,D)$.
Moreover, let $\nabla^E$ be a constraint covariant derivative on
$E$ and $\nabla$ a constraint covariant derivative on 
$T\mathcal{M}$.
Then
\begin{equation}
	\label{eq:ReductionOfSymbolCalculus}
	\Op_\red \colon \Secinfty\big(\Sym^\bullet T\mathcal{M}_\red 
	\tensor E^*_\red \tensor F_\red)\big)
	\to \Diffop^\bullet(E_\red;F_\red)
\end{equation}
is the symbol calculus associated to the vector bundles
$E_\red$ and $F_\red$ equipped with the covariant derivatives
$(\nabla^E)_\red$ and $\nabla_\red$.
\end{proposition}

\begin{proof}
Since reduction commutes with tensor products and taking sections we 
have
\begin{equation*}
	\ConSecinfty\big((\Sym_{\tensor}^\bullet T\mathcal{M} \tensor E^*) \strtensor F)\big)_\red 
	\simeq \Secinfty\big(\Sym^\bullet T\mathcal{M}_\red 
	\tensor E^*_\red \tensor F_\red)\big).
\end{equation*}
Together with \autoref{prop:ConDiffopsVSReduction} this shows that 
$\Op_\red$ is indeed of the form 
\eqref{eq:ReductionOfSymbolCalculus}.
The reduced map is given by
\begin{equation*}
	\Op_\red([X]) = \frac{1}{k!} \ins_s([X]) (\SymCoDer^E)_\red^k
	= \frac{1}{k!} \ins_s([X]) (\SymCoDer^{E_\red})^k
\end{equation*}
due to \autoref{prop:CovDerVSReduction} 
\ref{prop:CovDerVSReduction_SymCoDer}.
This shows that $\Op_\red$ is the associated symbol calculus on the 
reduced manifold.	
\end{proof}

The full symbol map allows us to improve the canonical morphism
\eqref{eq:ReductionOfSymbolCalculus} for the reduced differential operators to an isomorphism:

\begin{corollary}
Let $\mathcal{M} = (M,C,D)$ be a constraint manifold.
\begin{corollarylist}
	\item If $E$ and $F$ are constraint vector bundles, then
	\begin{equation}
		\ConDiffop^k(E;F)_\red \simeq \Diffop^k(E_\red;F_\red)
	\end{equation}
	for all $k \in \Naturals_0$.
	\item It holds
	\begin{equation}
		\ConDiffop^k(\mathcal{M})_\red \simeq \Diffop^k(\mathcal{M}_\red)
	\end{equation}
	for all $k \in \Naturals_0$.
\end{corollarylist}
\end{corollary}

\subsection{Multidifferential Operators}
\label{sec:ConMultiDiffops}

\index{multidifferential operator}
Grothendieck's definition of differential operators
can be extended to define multidifferential operators
$\Diffop^K(\module{E}_1, \dotsc, \module{E}_\ell; \module{F})$
of order $K = (k_1, \dotsc, k_\ell)$
on $\algebra{A}$-modules $\module{E}_1, \dotsc, \module{E}_\ell$
with values in an $\algebra{A}$-module
$\module{F}$.
We write $I \leq K$ for $I=(i_1, \dotsc, i_r)$ and $K = (k_1,\dotsc, k_r)$
if $i_\ell \leq k_\ell$ for all $\ell \in \{1, \dotsc r\}$.
Moreover, we write $\len(I) = r$ for the length of the multi index.

Note that
$\Diffop^K(\module{E}_1,\dotsc, \module{E}_\ell; \module{F})
\subseteq \Hom_\field{k}(\module{E}_1 \tensor[\field{k}] \dots \tensor[\field{k}] \module{E}_\ell; \module{F})$.
With this we can define constraint multidifferential operators
as those classical multidifferential operators compatible with the constraint structure.

\begin{definition}[Constraint multidifferential operators]
	\label{def:ConMultiDiffop}
	\index{constraint!multidifferential operator}
Let $\algebra{A} \in \injConAlg$ be a commutative embedded constraint algebra, and let
$\module{E}_1, \dotsc, \module{E}_\ell, \module{F}$ be embedded constraint
$\algebra{A}$-modules.
For a multi index
$K = (k_1, \dotsc, k_\ell) \geq 0$ we define
the \emph{constraint multidifferential operators} as
\begin{equation}
	\begin{split}
		\ConDiffop^K(\module{E}_1, \dotsc, \module{E}_\ell; \module{F})_\Total
		&\coloneqq \Diffop^K\big((\module{E}_1)_\Total,\dotsc, (\module{E}_\ell)_\Total; \module{F}_\Total\big) \\
		\ConDiffop^K(\module{E}_1, \dotsc, \module{E}_\ell; \module{F})_\Wobs
		&\coloneqq \big\{ D \in \Diffop^K\big((\module{E}_1)_\Total,\dotsc, (\module{E}_\ell)_\Total; \module{F}_\Total\big) \mid\\
		&\qquad\qquad\qquad D \in \ConHom_\field{k}(\module{E}_1 \tensor[\field{k}] \cdots \tensor[\field{k}] \module{E}_\ell, \module{F})_\Wobs \big\}, \\
		\ConDiffop^K(\module{E}_1, \dotsc, \module{E}_\ell; \module{F})_\Null
		&\coloneqq \big\{ D \in \Diffop^K\big((\module{E}_1)_\Total,\dotsc, (\module{E}_\ell)_\Total; \module{F}_\Total\big) \mid\\
		&\qquad\qquad\qquad D \in \ConHom_\field{k}(\module{E}_1 \tensor[\field{k}] \cdots \tensor[\field{k}] \module{E}_\ell, \module{F})_\Null \big\}.
	\end{split}
\end{equation}
\end{definition}

Note that $\ConDiffop^K(\module{E}_1,\dotsc,\module{E}_\ell;\module{F})$
becomes a strong constraint $\algebra{A}$-bimodule with 
respect to the classical left $\algebra{A}_\Total$-module structure given by
$(a \cdot D)(b) = a\cdot D(b)$.
Moreover, the constraint module of all multidifferential operators
$\ConDiffop^\bullet(\module{E}_1,\dotsc, \module{E}_\ell;\module{F})$
is filtered, in the sense that for multi indices $0 \leq L \leq K$ we have
\begin{equation}
	\ConDiffop^L(\module{E}_1,\dotsc, \module{E}_\ell;\module{F})
	\subseteq \ConDiffop^K(\module{E}_1,\dotsc, \module{E}_\ell;\module{F}).
\end{equation}

We now want to find a symbol calculus for constraint multidifferential operators
taking as arguments sections of constraint vector bundles.
For this we first need to study the local form of constraint multidifferential operators.

\begin{proposition}[Local form of constraint multidifferential operators]
	\label{prop:LocalConMultDiffop}
	{\ \\}
Let $E_1,\dotsc, E_\ell$ and $F$ be constraint vector bundles over a constraint manifold
$\mathcal{M} =(M,C,D)$ of dimension
$n = (n_\Total,n_\Wobs,n_\Null)$
and let
$D \in \ConDiffop^K(E_1,\dotsc,E_\ell; F)_\Total$ with
$K = (k_1, \dotsc, k_\ell)$.
Consider local adapted coordinates $(U,x)$ on $\mathcal{M}$ and let
$e_1^{(i)}, \dotsc, e_{n^i_\Total}^{(i)} \in \ConSecinfty(E_i)_\Total$
be constraint local frames with
$n^{i} = \rank E_i$.
Then
\begin{equation}
	\label{eq:LocalConMultDiffop}
	D\at{U}(s_1,\dotsc,s_\ell)
	= \sum_{0\leq R \leq K} \frac{1}{R!} 
	D_{U,\alpha_1,\dotsc,\alpha_\ell}^{R,I_1, \dotsc, I_\ell} 
	\del_{I_1} s_1^{\alpha_1} 
	\cdots \del_{I_\ell} s_\ell^{\alpha_\ell}
\end{equation}
for all $s_i \in \ConSecinfty(E_i)_\Total$,
with $D_{U,\alpha_1,\dotsc,\alpha_\ell}^{R,I_1,\dotsc,I_\ell} \in \ConSecinfty(F_\Total\at{U})$,
$I_j = (i_1^{(j)},\dotsc, i_{r_j}^{(j)})$,
and $s_j^{\alpha_j} = e^{\alpha_j}_{(j)}(s_j)$.
\begin{propositionlist}
	\item \label{prop:LocalConMultDiffop_1}
	It holds
	$D \in \ConDiffop^K(E_1,\dotsc, E_\ell;F)_\Wobs$
	if and only if
	\begin{equation}
		\label{prop:LocalConMultDiffop_Wobs}
		D_{U,\alpha_1,\dotsc,\alpha_\ell}^{R,I_1,\dotsc,I_\ell}
		\in \ConSecinfty(F)_\Wobs
	\end{equation}
	for $(I_1, \alpha_1,\dotsc,I_\ell,\alpha_\ell)
		\in \left( \left((n^*)^{\strtensor r_1} \strtensor \rank E_1\right) \tensor \cdots \tensor 
		\left((n^*)^{\strtensor r_\ell} \strtensor  \rank E_\ell \right) \right)_\Wobs$,
	and
	\begin{equation}
		D_{U,\alpha_1,\dotsc,\alpha_\ell}^{R,I_1,\dotsc,I_\ell}
		\in \ConSecinfty(F)_\Null
	\end{equation}
	for $(I_1, \alpha_1,\dotsc,I_\ell,\alpha_\ell) \in \left( \left( (n^*)^{\strtensor r_1} \strtensor \rank E_1\right) \tensor \cdots \tensor 
		\left((n^*)^{\strtensor r_\ell} \strtensor  \rank E_\ell \right) \right)_\Null$.
	\item \label{prop:LocalConMultDiffop_2}
	It holds
	$D \in \ConDiffop^K(E_1,\dotsc,E_\ell;F)_\Null$
	if and only if
	\begin{equation}
		\label{prop:LocalConMultDiffop_Null}
			D_{U,\alpha_1,\dotsc,\alpha_\ell}^{R,I_1,\dotsc,I_\ell} 
			\in \ConSecinfty(F)_\Null
	\end{equation}
	for $(I_1, \alpha_1,\dotsc,I_\ell,\alpha_\ell)
	\in \left( \left((n^*)^{\strtensor r_1} \strtensor \rank E_1\right) \tensor \cdots \tensor 
	\left((n^*)^{\strtensor r_\ell} \strtensor  \rank E_\ell\right) \right)_\Wobs$
\end{propositionlist}
\end{proposition}

\begin{proof}
We have
\begin{align*}
	D_{U,\alpha_1, \dotsc, \alpha_\ell}^{R,I_1,\dotsc,I_\ell}
 = D\at{U}(x^{i_1^{(1)}}\cdots x^{i_{r_1}^{(1)}} \cdot 
 e^{(1)}_{\alpha_1},\dotsc,
 x^{i_1^{(\ell)}}\cdots x^{i_{r_\ell}^{(\ell)}} \cdot 
 e^{(\ell)}_{\alpha_\ell}).
\end{align*}
Now from \autoref{ex:ConFunctionsOnConMfld} 
\ref{ex:ConFunctionsOnConMfld_1}
it follows that
\begin{equation*}
	x^{i_1^{(j)}}\cdots x^{i_{r_j}^{(j)}} \cdot e^{(j)}_{\alpha_j} \in 
	\ConSecinfty(E_j)_\Wobs
	\quad\text{ if }\quad
	(I_j \strtensor \alpha_j) \in ((n^*)^{\strtensor r_j} \strtensor 
	\rank E_j)_\Wobs
\end{equation*}
and that
\begin{equation*}
	x^{i_1^{(j)}}\cdots x^{i_{r_j}^{(j)}} \cdot e^{(j)}_{\alpha_j} \in 
	\ConSecinfty(E_j)_\Null
	\quad\text{ if }\quad
	(I_j \strtensor \alpha_j) \in ((n^*)^{\strtensor r_j} \strtensor 
	\rank E_j)_\Null.
\end{equation*}
Then for $D \in \ConDiffop^K(E_1, \dotsc, E_\ell; F)_\Wobs$ we 
immediately
get \eqref{prop:LocalConMultDiffop_Wobs}.
And similarly we obtain for $D \in \ConDiffop^K(E_1,\dotsc, E_\ell, 
F)_\Wobs$
directly \eqref{prop:LocalConMultDiffop_Null}.
For the other implication assume \eqref{prop:LocalConMultDiffop_Wobs} 
holds.
Let $s_1 \tensor \cdots \tensor s_\ell \in (\ConSecinfty(E_1) \tensor 
\cdots \tensor \ConSecinfty(E_\ell))_\Null$.
Then all terms of \eqref{eq:LocalConMultDiffop}
end up in $\ConSecinfty(F)_\Null$:	
We write 
\begin{equation*}
	S \coloneqq  \left( \left(n^{\tensor r_1} \tensor (\rank 
	E_1)^*\right) \strtensor \cdots\strtensor \left(n^{\tensor 
	r_\ell} \tensor (\rank E_\ell)^*\right) \right).
\end{equation*}
Recall that 
$(I_1, \alpha_1,\dotsc,I_\ell,\alpha_\ell) \in S_\Wobs$
if at least one of the pairs $(I_j,\alpha_j)$ has
$I_j = (i_1^{(j)}, \dotsc, i_{r_j}^{(j)}) \in n_\Wobs^{r_j}$
and $\alpha_j \in \rank(E_j)_\Total \setminus \rank(E_j)_\Null$
such that for at least one $m \in \{1, \dotsc, r_j\}$ 
it holds
$i_m^{(j)} \in n_\Null$
or $\alpha_j \in \rank(E_j)_\Total \setminus \rank(E_j)_\Wobs$.
Thus either
$(I_1, \alpha_1,\dotsc,I_\ell,\alpha_\ell) \in S_\Wobs$,
and hence
$\del^{I_j}s^{\alpha_j}_j \in \ConCinfty(\mathcal{M})_\Null$
for at least one $j \in \{1, \dotsc, \ell\}$, or
\begin{equation*}
	(I_1, \alpha_1,\dotsc,I_\ell,\alpha_\ell) \in 
	S_\Total \setminus S_\Wobs
	= \left( \left((n^*)^{\strtensor r_1} \strtensor \rank E_1\right) 
	\tensor \cdots\tensor
	\left((n^*)^{\strtensor r_\ell} \strtensor \rank E_\ell \right) 
	\right)_\Null,
\end{equation*}
and hence $D^{R,I_1,\dotsc,I_\ell}_{U,\alpha_1, \dotsc, \alpha_\ell}$.
For $s_1 \tensor \dots \tensor s_\ell \in (\ConSecinfty(E_1) \tensor 
\cdots \tensor \ConSecinfty(E_\ell))_\Wobs$
all terms of \eqref{eq:LocalConMultDiffop} end up in 
$\ConCinfty(F)_\Wobs$:
If $(I_1, \alpha_1,\dotsc,I_\ell,\alpha_\ell) \in S_\Null$,
then $\del^{I_j}s^{\alpha_j}_j \in \ConCinfty(\mathcal{M})_\Null$ for 
at least one $j \in \{1, \dotsc, \ell\}$.
If
\begin{equation*}
	(I_1, \alpha_1,\dotsc,I_\ell,\alpha_\ell)
	\in S_\Wobs \setminus S_\Null \subseteq 
	\left( \left((n^*)^{\strtensor r_1} \strtensor \rank E_1\right) 
	\tensor \cdots\tensor
	\left((n^*)^{\strtensor r_\ell} \strtensor \rank E_\ell \right) 
	\right)_\Wobs,
\end{equation*}
we have $\del^{I_j}s^{\alpha_j}_j \in \ConCinfty(\mathcal{M})_\Wobs$ 
for all $j$ and
$D^{R,I_1,\dotsc,I_\ell}_{U,\alpha_1,\dotsc,\alpha_\ell}
\in \ConSecinfty(F)_\Wobs$.
Finally if
\begin{equation*}
	(I_1, \alpha_1,\dotsc,I_\ell,\alpha_\ell)
	\in S_\Total \setminus S_\Wobs
	= \left( \left((n^*)^{\strtensor r_1} \strtensor \rank E_1\right)
	\tensor \cdots\tensor 
	\left((n^*)^{\strtensor r_\ell} \strtensor \rank E_\ell \right) 
	\right)_\Null,
\end{equation*}
then
$D^{R,I_1,\dotsc,I_\ell}_{U,\alpha_1,\dotsc,\alpha_\ell}
\in \ConSecinfty(F)_\Null$.
This gives the first part.
The second part follows by completely analogous considerations.
\end{proof}

The classical leading symbol
$\sigma_k(D) \in \Secinfty\left( (\Sym^{k_1} TM \tensor E_1^*) \tensor \cdots \tensor (\Sym^{k_\ell} TM \tensor E_\ell^*) \tensor F\right)$
for a multidifferential operator
$D \in \Diffop^K(E_1,\dotsc,E_\ell;F)$ is locally given by
\begin{equation}
	\sigma_K(D)\at{U} = \frac{1}{K!} 
	(\del^{\tensor}_{I_1} \tensor  e_{(1)}^{\alpha_1}) \tensor \cdots \tensor 
	(\del^{\tensor}_{I_\ell} \tensor e_{(\ell)}^{\alpha_\ell}) \tensor D_{U,\alpha_1,\dotsc,\alpha_\ell}^{K,I_1, \dotsc, I_\ell},
\end{equation}
with
\begin{equation}
	\del^{\tensor}_{I_j} \coloneqq \frac{\del}{\del x^{i^{(j)}_1}} 
	\vee \cdots \vee \frac{\del}{\del x^{i^{(j)}_{k_j}}}
	\in \Secinfty(\Sym^{k_j} TM).
\end{equation}
For constraint differential operators this will become a constraint 
section:

\begin{proposition}[Constraint leading symbol]
	\index{constraint!leading symbol}
Let $E_1, \dotsc, E_\ell$ and $F$ be constraint vector bundles over a 
constraint manifold
$\mathcal{M} = (M,C,D)$.
\begin{propositionlist}
	\item The leading symbol defines a constraint morphism
	\begin{equation} \label{eq:ConLeadingSymbol}
		\sigma_K \colon \ConDiffop^K(E_1, \dotsc,E_\ell;F)
		\to	\ConSecinfty\big( (\Sym_{\tensor}^{k_1} T\mathcal{M}
		\tensor E_1^*) \strtensor \cdots \strtensor
		(\Sym_{\tensor}^{k_\ell} T\mathcal{M} \tensor E_\ell^*)
		\strtensor F\big)
	\end{equation}
	of strong constraint 
	$\ConCinfty(\mathcal{M})$-modules. 
	\item If $E_1= \dots = E_\ell = F = \mathcal{M} \times \Reals$ 
	the leading symbol becomes a constraint morphism
	\begin{equation}
		\sigma_K \colon \ConDiffop^K(\mathcal{M}) \to 
		\ConSecinfty\big(\Sym^{k_1}_{\tensor} T\mathcal{M} \strtensor 
		\cdots \strtensor \Sym_{\tensor}^{k_\ell} T\mathcal{M}\big)
	\end{equation}
	of strong constraint $\ConCinfty(\mathcal{M})$-modules.
\end{propositionlist}
\end{proposition}

\begin{proof}
It only remains to shows that $\sigma_K$ is actually a constraint 
morphism.
We use again the shorthand
\begin{equation*}
	S \coloneqq  \left( \left(n^{\tensor r_1} \tensor
	(\rank E_1)^*\right) \strtensor \cdots\strtensor
	\left(n^{\tensor r_\ell} \tensor (\rank E_\ell)^*\right)
	\right).
\end{equation*}
First suppose
$D \in \ConDiffop^K(E_1, \dotsc, E_\ell;F)_\Wobs$.
For $(I_1,\alpha_1,\dotsc,I_\ell,\alpha_\ell) \in S_\Null$
it holds
\begin{equation*}
	\del_{I_j}^{\tensor} \tensor e_{(j)}^{\alpha_j}
	\in (\Sym_{\tensor}^{k_j} T\mathcal{M} \tensor E_j^*)_\Null
\end{equation*}
for one $j \in \{1, \dotsc, \ell\}$.
If
\begin{equation*}
	(I_1, \alpha_1,\dotsc,I_\ell,\alpha_\ell)
	\in S_\Wobs \setminus S_\Null \subseteq 
	\left( \left((n^*)^{\strtensor r_1} \strtensor \rank E_1\right) 
	\tensor \cdots\tensor
	\left((n^*)^{\strtensor r_\ell} \strtensor \rank E_\ell \right) 
	\right)_\Wobs,
\end{equation*}
then 
$D_{U,\alpha_1,\dotsc,\alpha_\ell}^{R,I_1,\dotsc,I_\ell}
\in \ConSecinfty(F)_\Wobs$
by \autoref{prop:LocalConMultDiffop}
and
$\del_{I_j}^{\tensor} \tensor e_{(j)}^{\alpha_j}
\in (\Sym_{\tensor}^{k_j} T\mathcal{M} \tensor E_j^*)_\Wobs$
for all $j=1,\dotsc,\ell$ by the definition of $S$.
Next, consider
\begin{equation*}
	(I_1, \alpha_1,\dotsc,I_\ell,\alpha_\ell)
	\in S_\Total \setminus S_\Wobs
	= \left( \left((n^*)^{\strtensor r_1} \strtensor \rank E_1\right) 
	\tensor \cdots\tensor
	\left((n^*)^{\strtensor r_\ell} \strtensor \rank E_\ell \right) 
	\right)_\Null.
\end{equation*}
Then
$D_{U,\alpha_1,\dotsc,\alpha_\ell}^{R,I_1,\dotsc,I_\ell}
\in \ConSecinfty(F)_\Null$
holds again by \autoref{prop:LocalConMultDiffop}.
Therefore, $\sigma_K$ preserves the $\WOBS$-component.
Finally, let
$D \in \ConDiffop^K(E_1, \dotsc, E_\ell;F)_\Null$
be given.
Then for 
$(I_1,\alpha_1,\dotsc,I_\ell,\alpha_\ell) \in S_\Null$
it holds again
$\del_{I_j}^{\tensor} \tensor e_{(j)}^{\alpha_j}
\in (\Sym_{\tensor}^{k_j} T\mathcal{M} \tensor E_j^*)_\Null$
for one $j \in \{1, \dotsc, \ell\}$.
And if
\begin{equation*}
	(I_1, \alpha_1,\dotsc,I_\ell,\alpha_\ell)
	\in S_\Total \setminus S_\Null
	= \left( \left((n^*)^{\strtensor r_1} \strtensor \rank E_1\right) 
	\tensor \cdots\tensor
	\left((n^*)^{\strtensor r_\ell} \strtensor \rank E_\ell \right) 
	\right)_\Wobs,
\end{equation*}
we obtain 
$D_{U,\alpha_1,\dotsc,\alpha_\ell}^{R,I_1,\dotsc,I_\ell}
\in \ConSecinfty(F)_\Null$
by \autoref{prop:LocalConMultDiffop}.
Thus, $\sigma_K$ is a constraint morphism.
The second part is a direct consequence.
\end{proof}

It is important to note that in the constraint leading symbol \eqref{eq:ConLeadingSymbol}
both kinds of tensor products appear.
In particular, we cannot rearrange the factors on the right hand side of \eqref{eq:ConLeadingSymbol}
in an arbitrary way, see also \eqref{eq:RebracketingConTensors}.

If constraint covariant derivatives $\nabla^{E_i}$ on constraint vector bundles $E_i$, $i = 1, \dotsc, \ell$,
are given, we define
\begin{equation}
\begin{split}
	D^K \coloneqq (D^{E_1})^{k_1} \tensor \dots \tensor (D^{E_\ell})^{k_\ell} 
	&\colon	\ConSecinfty(E_1) \tensor \cdots \tensor 
	\ConSecinfty(E_\ell) \\
	&\to \ConSecinfty(\Sym_{\strtensor}^{k_1}T^*\mathcal{M} 
	\strtensor E_1) \tensor	\cdots \tensor 
	\ConSecinfty(\Sym_{\strtensor}^{k_\ell}T^*\mathcal{M} \strtensor 
	E_\ell)
\end{split}	
\end{equation}
Note that $\ConSecinfty\left((\Sym_{\tensor}^\bullet T\mathcal{M} \tensor E_1^*) \strtensor \cdots \strtensor 
(\Sym_{\tensor}^\bullet T\mathcal{M} \tensor E_\ell^*) \strtensor F\right)$,
which is the dual of target space of $D^K$, is filtered by multi indices $K = (k_1, \dotsc, k_\ell)$.
With this we can now give the full symbol calculus for constraint multidifferential operators.

\begin{theorem}[Constraint multisymbol calculus]
	\label{thm:ConMultSymbolCalculus}
	\index{constraint!symbol calculus}
Let $E_1, \dotsc, E_\ell$ and $F$ be constraint vector bundles over a 
constraint manifold
$\mathcal{M} = (M,C,D)$.
Moreover, let $\nabla^{E_1}, \dotsc, \nabla^{E_\ell}$
be constraint covariant derivatives for $E_1, \dotsc, E_\ell$
and let $\nabla$ be a constraint covariant derivative for
$T\mathcal{M}$.
\begin{theoremlist}
	\item Then
	\begin{equation}
	\begin{split}
		\Op &\colon \ConSecinfty\left((\Sym_{\tensor}^\bullet 
		T\mathcal{M} \tensor E_1^*) \strtensor \cdots \strtensor 
		(\Sym_{\tensor}^\bullet T\mathcal{M} \tensor E_\ell^*) 
		\strtensor F\right) \\
		&\qquad\qquad\qquad\qquad\quad
		\longrightarrow \ConDiffop^\bullet(E_1, \dotsc, E_\ell;F),
	\end{split}
	\end{equation}
	defined by
	\begin{equation}
		\Op(X_1 \tensor \cdots \tensor X_\ell)(s_1,\dotsc,s_\ell)
		\coloneqq \frac{1}{k_1!\dots k_\ell!}\ins_s(X_1 \tensor \cdots 
		\tensor X_\ell) D^K(s_1\tensor \cdots \tensor s_\ell)
	\end{equation}
	for 
	$X_j \in \ConSecinfty(\Sym_{\tensor}^{k_j} T\mathcal{M} \tensor 
	E_j^*)_\Total$
	and 
	$s_j \in \ConSecinfty(E_j)_\Total$, with $K = (k_1, \dotsc, 
	k_\ell)$ and $j=1,\dotsc,\ell$,
	is a filtration preserving morphism of strong constraint 
	$\ConCinfty(\mathcal{M})$-modules.
	\item For $X_j \in \ConSecinfty(\Sym_{\tensor}^{k_j} T\mathcal{M} 
	\tensor 
	E_j^*)_\Total$, $j=1,\dotsc,\ell$ we have
	\begin{equation}
		\sigma_K\big(\Op(X_1 \tensor \cdots \tensor X_\ell)\big)
		= X_1 \tensor\cdots \tensor X_\ell,
	\end{equation}
where $\sigma_K$ denotes the leading symbol.
\item $\Op$ is a filtration preserving isomorphism of strong constraint 
$\ConCinfty(\mathcal{M})$-modules.
\end{theoremlist}
\end{theorem}

\begin{proof}
By the classical theory we know that $\Op$ fulfils all the above properties on the $\TOTAL$-compo\-nent,
see \cite[Chap. IV, §9]{palais:1965a} for a version of the symbol calculus for multidifferential operators.
For the first part note that $\Op$ is defined as a composition of constraint
morphism, and thus defines itself a constraint morphism to
$\ConHom_\field{k}(\ConSecinfty(E_1) \tensor \cdots \tensor
\ConSecinfty(E_\ell),\ConSecinfty(F))$.
Since we know that
$\Op(X_1 \tensor \cdots \tensor X_\ell)$
is a multidifferential operator it follows that
$\Op$ actually maps to the constraint submodule
$\ConDiffop^\bullet(E_1,\dotsc,E_\ell;F)$.
The second part is just the classical statement.
Nevertheless, from this follows directly that $\Op$ is a monomorphism.
To show that it is also a regular epimorphism we repeat the classical
argument for constructing preimages:
Let
$D \in \ConDiffop^K(E_1,\dotsc,E_\ell;F)_\Total$
be given.
Then $D - \Op(\sigma_K(D))$ has total order $\abs{K}-1$.
We write $X_K \coloneqq \sigma_K(D)$, then by induction we obtain
$D = \Op(\sum_{0 \leq R \leq K} X_R )$,
and thus $\Op$ is surjective.
If $D \in \ConDiffop^K(E_1,\dotsc,E_\ell; F)_\Wobs$,
then we have
\begin{equation*}
	X_K \in \ConSecinfty\big((\Sym_{\tensor}^{k_1} 
	T\mathcal{M} \tensor E_1^*) \strtensor \cdots \strtensor 
	(\Sym_{\tensor}^{k_\ell} T\mathcal{M} \tensor E_\ell^*) 
	\strtensor F\big)
\end{equation*}
and therefore $\Op$ is also surjective on the $\WOBS$-component.
Similarly, for a differential operators $D \in \ConDiffop^K(E_1,\dotsc,E_\ell; F)_\Null$
we get that $\Op$ is indeed a regular epimorphism.
Hence $\Op$ is a regular epimorphism and monomorphism, and therefore
a constraint isomorphism.
\end{proof}

For $E_1 = \dots = E_\ell = F = \mathcal{M} \times \Reals$
we immediately get the following isomorphism for multidifferential 
operators of
$\ConCinfty(\mathcal{M})$.

\begin{corollary}
	\label{cor:ConMultDiffOpSymbolCalculusFunctions}
Let $\mathcal{M} = (M,C,D)$ be a constraint manifold
and let $\nabla$ be a constraint covariant derivative for
$T\mathcal{M}$.
Then
\begin{align}
	\Op &\colon \ConSecinfty\left(\Sym_{\tensor}^\bullet
	T\mathcal{M} \strtensor \cdots \strtensor 
	\Sym_{\tensor}^\bullet T\mathcal{M}\right)
	\to \ConDiffop^\bullet(\mathcal{M}), \\
	\shortintertext{given by}
	\Op&(X_1 \tensor \cdots \tensor X_\ell)
	\coloneqq \frac{1}{k_1!\cdots k_\ell!}\ins_s(X_1 \tensor \cdots 
	\tensor X_\ell) D^K
\end{align}
for $X_1 \tensor \cdots \tensor X_\ell \in \ConSecinfty\big(\Sym_{\tensor}^{k_1}
T\mathcal{M} \strtensor \cdots \strtensor 
\Sym_{\tensor}^{k_\ell} T\mathcal{M}\big)$,
is a filtration preserving isomorphism of constraint $\ConCinfty(\mathcal{M})$-modules.
\end{corollary}

\subsubsection{Reduction}

The various compatibilities of constraint multidifferential operators
with reduction are by now quite obvious and 
can be proven in a completely analogous fashion to 
those of constraint differential operators in 
\autoref{sec:ConDiffOps} and
\autoref{sec:ConSymbCalculus}.
We will therefore just give the statements without repeating the 
proofs.

\begin{proposition}[Constraint multidifferential operators vs. reduction]{\ \\}
	\label{prop:ConMultDiffopsVSReduction}
	\index{reduction!symbol calculus}Let $\algebra{A}$ be a commutative embedded constraint 
algebra, and let
$\module{E}^1, \dotsc, \module{E}^\ell, \module{F}$ be embedded 
constraint
$\algebra{A}$-modules.
For any multi index
$K = (k_1, \dotsc, k_\ell) \in \Naturals_0^\ell$
there is a natural injective morphism
\begin{equation}
	\ConDiffop^K\left(\module{E}_1,\dotsc,\module{E}_\ell;
	\module{F}\right)_\red 
	\to 
	\Diffop^K\left((\module{E}_1)_\red,\dotsc,(\module{E}_\ell)_\red; 
	\module{F}_\red\right)
\end{equation}
of $\algebra{A}_\red$-modules.
\end{proposition}

For multidifferential operators this becomes an isomorphism:

\begin{proposition}[Constraint multidifferential operators of 
sections vs. reduction]
	\label{prop:ConSecMultDiffopsVSReduction}
Let $E_1,\dotsc, E_\ell$ and $F$ be constraint vector bundles 
over a constraint manifold
$\mathcal{M} =(M,C,D)$ of dimension
$n = (n_\Total,n_\Wobs,n_\Null)$.
\begin{propositionlist}
	\item Let
	$D \in \ConDiffop^K(E_1,\dotsc,E_\ell; F)_\Wobs$
	of order
	$K = (k_1, \dotsc, k_\ell)$ be given.
	Then locally
	\begin{equation}
		\label{eq:LocalConMultDiffop_red}
		D_\red\at{U}(s_1,\dotsc,s_\ell)
		= \sum_{0\leq R \leq K} 
		\sum_{n_0<I_1,\dotsc,I_\ell \leq n_\Wobs}\frac{1}{R!} 
		(D_{U,\alpha_1,\dotsc,\alpha_\ell}^{R,I_1, \dotsc, 
			I_\ell})_\red 
		\del_{I_1} s_1^{\alpha_1} 
		\cdots \del_{I_\ell} s_\ell^{\alpha_\ell}
	\end{equation}
	for all $s_i \in \ConSecinfty((E_i)_\red\at{U})$.
	\item The constraint leading symbol $\sigma_K$ for constraint 
	multidifferential operators induces the classical leading symbol
	\begin{equation}
	\begin{split}
		\sigma_K &\colon \Diffop^K\big((E_1)_\red, 
		\dotsc,(E_\ell)_\red;F_\red\big) \\
		&\qquad\qquad \longrightarrow	\Secinfty\big( \Sym^{k_1} T\mathcal{M}_\red
		\tensor (E_1)_\red^* \tensor \cdots \tensor
		\Sym^{k_\ell} T\mathcal{M}_\red \tensor (E_\ell)_\red^*
		\tensor F_\red\big)
	\end{split}
	\end{equation}
	on the reduced manifold $\mathcal{M}_\red$.
	\item Let $\nabla^{E_1}, \dotsc, \nabla^{E_\ell}$
	be constraint covariant derivatives for $E_1, \dotsc, E_\ell$
	and let $\nabla$ be a constraint covariant derivative for
	$T\mathcal{M}$.
	Then
	\begin{equation}
		\label{eq:ReductionOfMultSymbolCalculus}
	\begin{split}
		\Op_\red \colon& \Secinfty\big(\Sym^\bullet 
		T\mathcal{M}_\red \tensor (E_1)_\red^* \tensor \cdots \tensor
		\Sym^\bullet T\mathcal{M}_\red \tensor (E_\ell)_\red^*
		\tensor F_\red\big) \\
		&\longrightarrow \Diffop^\bullet\big((E_1)_\red, \dotsc, 
		(E_\ell)_\red;F_\red\big)
	\end{split}
	\end{equation}
	is the symbol calculus associated to the vector bundles
	$(E_1)_\red,\dotsc, (E_\ell)_\red$ and $F_\red$ equipped with the 
	covariant derivatives
	$(\nabla^{E_1})_\red, \dotsc, (\nabla^{E_\ell})_\red$ and 
	$\nabla_\red$.
	\item It holds
	\begin{equation}
		\ConDiffop^\bullet(E_1, \dotsc, E_\ell;F)_\red
		\simeq \Diffop^\bullet\big((E_1)_\red, \dotsc, (E_\ell)_\red;F_\red\big)
	\end{equation}
	as filtered strong constraint $\ConCinfty(\mathcal{M})$-modules.
	\item It holds
	\begin{equation}
		\ConDiffop^\bullet(\mathcal{M})_\red
		\simeq \Diffop^\bullet(\mathcal{M}_\red )
	\end{equation}
	as filtered strong constraint $\ConCinfty(\mathcal{M})$-modules.
\end{propositionlist}
\end{proposition}

%% file: deformation-theory.tex
Formal deformation quantization aims to construct a quantum analogue 
of a given classical mechanical system by deforming the 
multiplication of the algebra $\Cinfty(M)$ of smooth functions on a 
Poisson manifold $(M,\pi)$
into a non-commutative multiplication $\star$ on the algebra
$\Cinfty(M)\formal{\lambda}$ of formal power series.
An important observation from classical deformation quantization is 
that given such a star product $\star$ we can always reconstruct
a Poisson bracket on $\Cinfty(M)$ from the $\star$-commutator.
It is now reasonable to only consider those star products which 
recover the Poisson structure on $M$.
If we start with a constraint Poisson manifold $(\mathcal{M},\pi)$
we obtain a commutative strong constraint algebra
$\ConCinfty(\mathcal{M})$, 
and there are two reasonable ways to define a constraint star product:
either as a deformation into a strong constraint algebra 
multiplication or, more generally, as a deformation into a constraint 
multiplication.
If we would consider deformations as strong constraint algebras, 
it can be shown that the induced Poisson 
bracket on $\ConCinfty(\mathcal{M})$ would be a strong constraint 
Poisson bracket.
And thus, from our discussion after 
\autoref{prop:CharacterizationConPoissonMfld},
the submanifold
$C \subseteq M$ would need to be a Poisson submanifold.
Thus if we want to consider star products which are compatible with 
honest coisotropic submanifolds we are forced to consider 
deformations of $\ConCinfty(\mathcal{M})$ as, in general, non-strong 
constraint algebras.

Following these ideas we will introduce constraint star products and deformations of constraint
algebras in
\autoref{sec:ConStarProducts}.
Then, following general ideas from deformation theory, we will study Maurer-Cartan elements and
their equivalence for constraint DGLAs in \autoref{sec:ConDeformationFunctor}
before we introduce constraint Hochschild cohomology in \autoref{sec:ConHochschildCohomology}.
Then in \autoref{sec:FormalDeformations} we will identify constraint Hochschild cohomology as the cohomology theory governing
the deformation problem of constraint algebras.
Finally, in \autoref{sec:SecondConHochschildCohomology} we will take some first steps into the direction of a constraint Hochschild-Kostant-Rosenberg theorem.
In particular, we will compute the second constraint Hochschild cohomology of the constraint functions on $\Reals^n$ in \autoref{sec:SecondConHochschildCohomology}, which already exhibits unexpected contributions.
\section{Constraint Star Products}
\label{sec:ConStarProducts}

Recall from \autoref{prop:FunctionsOnConManifolds} the definition of the constraint algebra of functions $\ConCinfty(\mathcal{M})$ on a constraint manifold $\mathcal{M} = (M,C,D)$.
Let us define the strong constraint algebra 
$\ConCinfty(\mathcal{M})\formal{\lambda}$ of formal power series
by
\index{formal power series}
\begin{equation}
	\ConCinfty(\mathcal{M})\formal{\lambda}
	\coloneqq \big(\ConCinfty(\mathcal{M})_\Total\formal{\lambda},\,\,
	\ConCinfty(\mathcal{M})_\Wobs\formal{\lambda},\,\,
	\ConCinfty(\mathcal{M})_\Null\formal{\lambda}\big).
\end{equation}
With this we can state our definition of constraint star product:

\begin{definition}[Constraint star product]
	\label{def:ConStarProduct}
	\index{constraint!star product}
Let $(\mathcal{M},\pi)$ be a constraint Poisson manifold.
A \emph{(formal) constraint star product} $\star$ on
$(\mathcal{M},\pi)$ is a $\ComplexNum\formal{\lambda}$-linear 
constraint map
\begin{equation}
	\star \colon \ConCinfty(\mathcal{M})\formal{\lambda}
	\tensor[\ComplexNum\formal{\lambda}] \ConCinfty(\mathcal{M})\formal{\lambda}
	\to 
	\ConCinfty(\mathcal{M})\formal{\lambda}
\end{equation}
of the form
\begin{equation}
	\star = \sum_{r=0}^{\infty} \lambda^r C_r,
\end{equation}
with $\ComplexNum$-bilinear constraint maps
$C_r \colon \ConCinfty(\mathcal{M}) \tensor[\ComplexNum] 
\ConCinfty(\mathcal{M}) 
\to \ConCinfty(\mathcal{M})$,
fulfilling:
\begin{definitionlist}
	\item \label{def:ConStarProduct_1}
	$\star$ is associative.
	\item \label{def:ConStarProduct_2}
	$1 \star f = f = f \star 1$.
	\item \label{def:ConStarProduct_3}
	$\star = \mu_0 + \sum_{r=1}^\infty \lambda^r C_r$,
	with $\mu_0$ the pointwise multiplication on 
	$\ConCinfty(\mathcal{M})$.
	\item \label{def:ConStarProduct_4}
	$\frac{1}{\I \hbar}[f,g]_\star = \{f,g\} + 
	\lambda(\dots)$.
	\item \label{def:ConStarProduct_5}
	$C_r$ is a constraint bidifferential operator, for all
	$r \in \Naturals_0$.
\end{definitionlist}
\end{definition}

\begin{example}[Standard ordered star product]
Consider the classical standard-ordered\linebreak
star product
\begin{equation} \label{eq:StdOrderedStarProduct}
	f \star_{\standard} g = \sum_{r=0}^{\infty} \frac{1}{r!} \left(\frac{\hbar}{\I}\right)^r
	\sum_{i_1, \dotsc, i_r} \frac{\del^r f}{\del p_{i_1} \dotsc \del p_{i_r}}
	 \frac{\del^r g}{\del q^{i_1} \dotsc \del q^{i_r}}
\end{equation}
on $T^*\Reals^n \simeq \Reals^{2n}$
with coordinates $(q^1,\dotsc,q^{n}, p_1,  \dotsc, p_{n})$.
By a change of coordinates any coisotropic subspace $C$ can be identified with
$\Reals^{n + k}$
with coordinates $(q^1, \dotsc q^{n}, p_1, \dotsc p_k)$
and its characteristic distribution is then given by 
$\Reals^{n_\Total-k}$ with coordinates
$(q^{n-k+1},\dotsc,q^{n})$.
Thus we can consider the constraint vector space
$\mathcal{M} = (\Reals^{2n}, \Reals^{n+k} \oplus \{0\}^{n - k},\{0\}^{k} \oplus \Reals^{n-k} \oplus \{0\}^{n})$.
From classical deformation quantization we know that $\star_{\standard}$ defines a 
star product on $\Reals^{2n}$ and it is straightforward to check that
$\star_{\standard}$ indeed defines a constraint multiplication.
For this it is important to note that for $i_1,\dotsc,i_r \leq k$ we have
\begin{equation}
	\frac{\del^r}{\del p_{i_1} \dotsc \del p_{i_r}} \in \ConDiffop^r(\ConCinfty(\mathcal{M}))_\Wobs
	\quad\text{ and }\quad
	\frac{\del^r}{\del q^{i_1} \dotsc \del q^{i_r}} \in \ConDiffop^r(\ConCinfty(\mathcal{M}))_\Wobs.
\end{equation}
But if there is one $\ell \in \{1, \dotsc, r\}$ such that $i_\ell > k$, then
$\frac{\del^r}{\del p_{i_1} \dotsc \del p_{i_r}}$ is \emph{not} constraint any more.
Nevertheless, in this case we have
\begin{equation}
	\frac{\del^r}{\del q^{i_1} \dotsc \del q^{i_r}} \in \ConDiffop^r(\ConCinfty(\mathcal{M}))_\Null
\end{equation}
by \autoref{ex:ParitalsAsConDiffOps},
making \eqref{eq:StdOrderedStarProduct} a constraint star product.
\end{example}

We want to study constraint star products using a constraint version
of Gerstenhaber's theory of deformation of associative algebras.
Thus in the rest of this section we will consider possibly non-unital constraint algebras.
Then we want to consider deformations of a constraint algebra 
$\algebra{A}$ with respect to the (constraint) ring
$\field{k}\formal{\lambda} 
= (\field{k}\formal{\lambda},\field{k}\formal{\lambda},0)$.
In general, the constraint module of formal power series of a given
constraint module $\module{E}$ is defined as
\begin{equation}
	\module{E}\formal{\lambda} \coloneqq 
	(\module{E}_\Total\formal{\lambda},
	\module{E}_\Wobs\formal{\lambda},
	\module{E}_\Null\formal{\lambda})
\end{equation}
with $\iota_{\module{E}\formal{\lambda}}$ given by the $\lambda$-linear 
extension of $\iota_\module{E}$.
\index{classical limit}
The \emph{classical limit} of a given constraint 
$\field{k}\formal{\lambda}$-module $\algebra{E}$
is defined by
\begin{equation}
	\cl(\module{E}) \coloneqq \module{E} / \lambda \module{E}.
\end{equation}
This defines a functor
$\cl \colon \ConMod_{\field{k}\formal{\lambda}} \to \ConMod_{\field{k}}$,
and it can be shown that taking the classical limit commutes with 
reduction, i.e. there is a natural isomorphism making the diagram
\begin{equation}
\begin{tikzcd}
	\ConMod_{\field{k}\formal{\lambda}}
		\arrow[r,"\cl"]
		\arrow[d,"\red"{swap}]
	& \ConMod_{\field{k}}
		\arrow[d,"\red"] \\
	\Modules_{\field{k}\formal{\lambda}}
		\arrow[r,"\cl"]
	&\Modules_{\field{k}}
\end{tikzcd}
\end{equation}
commute, see our work \cite[Thm. 7.13]{dippell.esposito.waldmann:2019a} for details.

Now we can define a formal associative deformation of a constraint algebra 
$\algebra{A}$ to be a constraint
$\field{k}\formal{\lambda}$-algebra $\algebra{B}$
together with an isomorphism
$\alpha \colon \cl(\algebra{B}) \to \algebra{A}$.
It is easy to see that this definition agrees with the one from deformation 
via Artin rings, see e.g. \cite{manetti:2009a}.
Usually, one is interested in more specific deformations, namely those that 
are e.g. free $\field{k}$-modules.
This leads us to the following definition:

\begin{definition}[Deformation of constraint algebra]
	\label{def:DeformationConAlg}%
	\index{deformation of constraint algebra}
Let $\algebra{A} \in \ConAlg_\field{k}$ be a (possibly non-unital) constraint algebra.
A \emph{(associative formal) deformation} of $\algebra{A}$ is given by an 
associative multiplication
$\mu \colon \algebra{A}\formal{\lambda} \tensor[\field{k}\formal{\lambda}] 
\algebra{A}\formal{\lambda} \to \algebra{A}\formal{\lambda}$
on $\algebra{A}\formal{\lambda}$
turning it into a constraint
$\field{k}\formal{\lambda}$-algebra, such that
$\cl(\algebra{A}\formal{\lambda},\mu) \simeq \algebra{A}$.
\end{definition}

Note that we have two formal associative deformations
$\mu_\Total$ and $\mu_\Wobs$ for
$\algebra{A}_\Total\formal{\lambda}$ and
$\algebra{A}_\Wobs\formal{\lambda}$ of the form
$\mu_\Total 
= (\mu_\Total)_0 + \lambda (\mu_\Total)_1 + \lambda^2(\dotsc)$
and
$\mu_\Wobs 
= (\mu_\Wobs)_0 + \lambda (\mu_\Wobs)_1 + \lambda^2(\dotsc)$,
respectively, such that the \emph{undeformed} map
$\iota_{\algebra{A}}$ is an algebra homomorphism and such that
$\algebra{A}_\Null\formal{\lambda}$ is a two-sided ideal in
$\algebra{A}_\Wobs\formal{\lambda}$ with respect to $\mu_\Wobs$.  
We insist on the $\algebra{A}_\Wobs$ and $\algebra{A}_\Null$
being the \emph{same} up to taking formal series.
Also the algebra morphism $\iota_{\algebra{A}}$ is \emph{not} deformed.

\begin{remark}
There are different approaches to study the deformations of diagrams
of associative algebras, e.g. via derived bracket as in 
\cite{fregier.zambon:2015} or with an operadic approach as in 
\cite{fregier.markl.yau:2009}.
See also \cite{gerstenhaber.schack:1983a}.
Nevertheless, our goal is to deform the multiplicative structure of a 
constraint algebra, but not the morphism it contains.
A thorough comparison to these deformations of diagrams needs to be done.
\end{remark}

\index{equivalence!deformations}
We say that two formal associative deformations $\mu$ and
$\mu'$ of $(\algebra{A},\mu_0)$
are \emph{equivalent} if there exists
$T = \id + \lambda(\ldots) \in
\ConAut_{\field{k}\formal{\lambda}}(\algebra{A}\formal{\lambda})_\Wobs$
such that $T \circ \mu = \mu' \circ (T \tensor T)$, i.e. we have
\begin{equation} \label{eq:EquivalenceInComponents}
	T_\Total(\mu_\Total (a,b))
	= \mu'_\Total (T_\Total(a), T_\Total(b))
	\quad\textrm{and}\quad
	T_\Wobs(\mu_\Wobs (a,b))
	= \mu'_\Wobs (T_\Wobs(a), T_\Wobs(b))
\end{equation}
for $a, b \in \algebra{A}_{\Total/\Wobs}$.
Thus, as in the case of associative algebras, there exists a unique
$D = \sum_{k=0}^\infty \lambda^k D_k \in
\ConHom_{\field{k}\formal{\lambda}}(\algebra{A}\formal{\lambda},
\algebra{A}\formal{\lambda})_\Wobs$
such that
$T = \exp(\lambda D)$.

\begin{remark}
	\label{rem:UnitalStarProduct}
Suppose that $(\algebra{A},\mu_0)$ is a unital constraint algebra
with unit $1$.
Then we know from classical deformation theory, see
\cite[Sec. 20]{gerstenhaber:1968a},
that any deformation of $\algebra{A}_\Wobs$ is again unital.
This unit also serves as a unit for the deformation of the constraint algebra
$\algebra{A}$.
Then the multiplication with this unit yields an equivalence to
a deformation of $\algebra{A}$
with unit given by $1$.
Thus in the following we can always assume to deform
unital constraint algebras into unital constraint algebras.
\end{remark}

Let us show that a constraint deformation of a commutative 
constraint algebra always induces a constraint Poisson structure on it:

\begin{proposition}
	\label{prop:ConPoissonStructureFromDeformation}
	\index{constraint!Poisson algebra}
Let $(\algebra{A},\mu_0)$ be a commutative 
constraint algebra, and let
$\mu = \mu_0 + \lambda \mu_1 + \cdots$
be an associative formal deformation of $\algebra{A}$.
Then
\begin{equation}
	\{ \argument, \argument \} \coloneqq \mu_1 - \mu_1 \circ \tau,
\end{equation}
with $\tau$ denoting the flip,
defines a constraint Poisson structure on $\algebra{A}$.
\end{proposition}

\begin{proof}
From classical deformation theory we know that
$\{\argument, \argument\}_\Total$
and $\{\argument, \argument\}_\Wobs$
define Poisson structures on $\algebra{A}_\Total$ 
and $\algebra{A}_\Wobs$.
Since $\{\argument, \argument\}$ is defined by a composition
of constraint maps it gives a constraint Poisson structure on
$\algebra{A}$.
\end{proof}

\index{constraint!Poisson manifold}
When considering $\algebra{A} = \ConCinfty(\mathcal{M})$
the above result shows that every deformation
induces the structure of a constraint Poisson manifold on
$\mathcal{M}$.
In this sense, property \ref{def:ConStarProduct_4} of 
\autoref{def:ConStarProduct}
is always fulfilled for \emph{some} constraint Poisson structure.
Together with \autoref{rem:UnitalStarProduct} we see that a constraint star product is nothing but a formal associative deformation
of the strong constraint algebra $\ConCinfty(\mathcal{M})$ by
bidifferential operators.

One particular scenario we will be interested in the context of
deformation quantization of phase space reduction is the 
following.
This, and the following two examples are taken from our work \cite[Sec. 4]{dippell.esposito.waldmann:2022a}.

\begin{example}
	\label{exmp:DeformationStrConAlg}%
We will work over a field $\field{K}$ instead of a general ring.
Let
$\algebra{A} = (\algebra{A}_\Total, \algebra{A}_\Wobs, \algebra{A}_\Null)$
be a unital embedded constraint algebra such that additionally
$\algebra{A}_\Null \subseteq \algebra{A}_\Total$ is a left 
ideal, then
$\algebra{A}_\Wobs \subseteq \normalizer(\algebra{A}_\Null)$
is a unital subalgebra of the normalizer of this left ideal.
Consider now a formal associative deformation $\mu_\Total$ of
$\algebra{A}_\Total$ with the additional property that the 
formal series
$\algebra{A}_\Null\formal{\lambda}$
are still a left ideal inside
$\algebra{A}_\Total\formal{\lambda}$
with respect to $\mu_\Total$.
Then we know that the normalizer
$\qalgebra{A}_\Wobs \coloneqq
\normalizer_{\mu_\Total}(\algebra{A}_\Null\formal{\lambda}) 
\subseteq \algebra{A}_\Total\formal{\lambda}$
with respect to $\mu_\Total$ satisfies
$\cl(\qalgebra{A}_\Wobs) \subseteq \normalizer(\algebra{A}_\Null)$.
We assume additionally 
$\cl(\qalgebra{A}_\Wobs) \subseteq \algebra{A}_\Wobs$.
This would be automatically true if $\algebra{A}_\Wobs$ coincides with the 
undeformed normalizer but poses an additional condition otherwise.

It is now easy to check that
$\qalgebra{A}_\Wobs \subseteq \algebra{A}_\Total\formal{\lambda}$
is a \emph{closed} subspace with respect to the $\lambda$-adic
topology.
Moreover, if $\lambda a \in \qalgebra{A}_\Wobs$ for some
$a \in \algebra{A}_\Total\formal{\lambda}$
we can conclude
$a \in \qalgebra{A}_\Wobs$.
Hence
$\qalgebra{A}_\Wobs \subseteq \algebra{A}_\Total\formal{\lambda}$
is a deformation of a subspace in the sense of 
\cite[Def.~30]{bordemann.herbig.waldmann:2000a},
i.e. we have a subspace
$\algebra{D} \subseteq \algebra{A}_\Total$
and linear maps
$q_r\colon \algebra{D} \longrightarrow \algebra{A}_\Total$,
for $r \in \mathbb{N}$, such that
$\qalgebra{A}_\Wobs = q (\algebra{D}\formal{\lambda})$,
where
$q = \iota_{\algebra{D}} + \sum_{r=1}^\infty \lambda^r q_r$
with $\iota_{\algebra{D}}$ being the canonical inclusion of 
the subspace.
By our assumption 
$\algebra{D} \subseteq \algebra{A}_\Wobs$,
but the inclusion could be proper.
Moreover, since by our assumption
$\algebra{A}_\Null\formal{\lambda} \subseteq 
\normalizer(\algebra{A}_\Null\formal{\lambda})
= \qalgebra{A}_\Wobs$,
we have
$\algebra{A}_\Null \subseteq \algebra{D}$.

Since we work over a field, we can find a complement
$\algebra{C} \subseteq \algebra{D}$
such that
$\algebra{A}_\Null \oplus \algebra{C} = \algebra{D}$.
This allows to redefine the maps $q_r$ to
\begin{equation} \label{eq:qprimer}
	q'_r\at{\algebra{C}} = q_r\at{\algebra{C}}
	\quad \textrm{and} \quad
	q'_r\at{\algebra{A}_\Null} = 0.
\end{equation}
The resulting map $q'$ then satisfies
$q'(\algebra{D}\formal{\lambda}) = \qalgebra{A}_\Wobs$
and
$q'\at{\algebra{A}_\Null} = \id_{\algebra{A}_\Null}$.
We can then use $q'$ to pass to a new deformation $\mu'_\Total$ of 
$\algebra{A}_\Total$ with the property that
$\algebra{A}_\Null\formal{\lambda}$
is still a left ideal in
$\algebra{A}_\Total\formal{\lambda}$
with respect to $\mu'_\Total$ and the normalizer of this left ideal is now 
given by
$\algebra{D}\formal{\lambda} \subseteq \algebra{A}_\Total\formal{\lambda}$.
It follows that $\mu'_\Total$ provides a deformation of the constraint 
algebra
$(\algebra{A}_\Total, \algebra{D}, \algebra{A}_\Null)$
in the sense of \autoref{def:DeformationConAlg}.

Of course, it might happen that
$\algebra{D} \ne \algebra{A}_\Wobs$
and hence this construction will not provide a deformation of the original 
constraint algebra in general.
It turns out that this can be controlled as follows:
We assume in addition that the deformed normalizer
$\qalgebra{A}_\Wobs$ is \emph{large enough} in the sense that the classical 
limit
\begin{equation} \label{eq:LargeEnoughNormalizer}
	\cl\colon \qalgebra{A}_\red
	= \qalgebra{A}_\Wobs \big/ (\algebra{A}_\Null\formal{\lambda})
	\to \algebra{A}_\red 
	= \algebra{A}_\Wobs \big/ \algebra{A}_\Null
\end{equation}
between the reduced algebras is \emph{surjective}.
As $\field{K}$ is a field, this gives us a split
$Q\colon \algebra{A}_\red \longrightarrow \qalgebra{A}_\red$ 
which we can extend $\lambda$-linearly to
$Q\colon \algebra{A}_\red\formal{\lambda} \to \qalgebra{A}_\red$.
It is then easy to see that this is in fact a 
$\field{K}\formal{\lambda}$-linear isomorphism.
It follows, that in this case we necessarily have
$\algebra{D} = \algebra{A}_\Wobs$.
Thus the previous construction gives indeed a deformation
$\mu'_\Total$ of the original constraint algebra.
This seemingly very special situation will turn out to be responsible for 
one of the main examples from deformation quantization.
\end{example}

In the following we present two examples from deformation quantization
which can be interpreted as deformations of constraint algebras in the above sense.
These show that even though we will mainly be interested in the abstract deformation theory of constraint algebras,
these actually appear in well-studied situations.
Note, however, that both examples should illustrate the
concept of a deformation of a constraint algebra \emph{without}
actually computing the corresponding Hochschild cohomology.
Even in these examples it seems to be a rather difficult task to compute 
the constraint Hochschild cohomology of a constraint manifold $\mathcal{M} = (M,C,D)$.

\renewcommand{\thesubsubsection}{\thesection.\arabic{subsubsection}}
\subsubsection{Example I: BRST Reduction}
\label{subsec:ExampleIBRST}

The first example comes from BRST reduction of star products.
We recall the situation of
\cite{bordemann.herbig.waldmann:2000a,gutt.waldmann:2010a}.
Consider a Poisson manifold $M$ with a strongly Hamiltonian action
of a connected Lie group $G$ and momentum map
$J\colon M \longrightarrow \liealg{g}^*$,
where $\liealg{g}$ is the Lie algebra of $G$.
One assumes that the classical level surface
$C = J^{-1}(\{0\}) \subseteq M$
is a non-empty (necessarily coisotropic) submanifold by requiring $0$
to be a regular value of $J$.
Moreover, we assume that the action on $C$ is free and proper.
Then we have a constraint manifold
$\mathcal{M} = (M,C,D)$
with $D$ the characteristic distribution on $C$.
This leads to the strong constraint algebra
\begin{equation}
	\algebra{A} \coloneqq \ConCinfty(\mathcal{M}) = \big(\Cinfty(M),\, \Pnormalizer_C,\, \vanishing_C\big),
\end{equation}
where
$\vanishing_C = \ker \iota^* \subseteq \Cinfty(M)$
is the vanishing ideal of the constraint surface
$C \subseteq M$ and $\Pnormalizer_C$ its Poisson normalizer,
cf. \autoref{ex:ConFunctionsOnConMfld} \ref{ex:ConFunctionsOnConMfld_2}.
Next, we assume to have a star product $\star$ strongly invariant under the action
of $G$ which admits a deformation $\deform{J}$ of $J$ into a quantum momentum map.
In the symplectic case such star products always exist since we assume the action 
of $G$	to be proper, see \cite{reichert.waldmann:2016a} for a complete
classification and further references.
In the general Poisson case the	situation is less clear.

Out of this a constraint $\ComplexNum\formal{\lambda}$-algebra
$\qalgebra{A} = (\Cinfty(M)\formal{\lambda}, \qPnormalizer_C,\qvanishing_C)$
is then constructed, where
$\qvanishing_C \coloneqq \ker \deform{\iota}^*\subseteq \Cinfty(M)\formal{\lambda}$
is the quantum vanishing ideal given by the kernel of the deformed restriction
$\deform{\iota}^* \coloneqq \iota^* \circ S$.
Here
$S = \id + \sum_{k=1}^{\infty} \lambda^k S_k$
is a formal power series of differential operators guaranteeing that $\qvanishing_C$ is 
indeed a left ideal with respect to $\star$.
In fact, $S$ can be chosen to be $G$-invariant.

We now want to construct a constraint algebra structure on
$\algebra{A}\formal{\lambda} = 
(\Cinfty(M)\formal{\lambda}, \algebra{B}_C\formal{\lambda},\algebra{J}_C\formal{\lambda})$
which is isomorphic to
$\qalgebra{A}$.  
For
this, note that
$S\colon \Cinfty(M)\formal{\lambda} \longrightarrow
\Cinfty(M)\formal{\lambda}$ is invertible, hence we get a star 
product
\begin{equation}
	f \star^S g \coloneqq S(S^{-1}f \star S^{-1}g)
\end{equation}
on $\Cinfty(M)\formal{\lambda}$.
From
$\deform{\iota}^* = \iota^* \circ S$
it directly follows that $S$ maps $\qvanishing_C$ to
$\vanishing_C\formal{\lambda}$.
It is slightly less evident, but follows from the characterization of the normalizer
$\qPnormalizer_C$ as those functions whose restriction to $C$ are
$G$-invariant, that $S$ maps the normalizer $\qPnormalizer_C$ to 
the normalizer $\qPnormalizer_C^S$ of $\qvanishing_C$ 
with respect to $\star^S$.
Finally, we know that $f \in \qPnormalizer_C$ if
and only if for all $\xi \in \liealg{g}$ it holds that
$0 = \Lie_{\xi_C} \deform{\iota}^* f = \Lie_{\xi_C} \iota^*Sf$,
where $\Lie_{\xi_C}$ denotes the Lie derivative in the direction 
of the fundamental vector field $\xi_C$.
Hence $f \in \qPnormalizer_C$ if and only if
$Sf \in \Pnormalizer_C\formal{\lambda}$.
Thus $S$ is an isomorphism of constraint algebras
\begin{equation}
	S \colon (\qalgebra{A},\star)
	\to
	(\ConCinfty(\mathcal{M})\formal{\lambda},\star^S).
\end{equation}
In particular, we have a deformation
of the classical constraint algebra in this case, and the 
constraint algebra
$\qalgebra{A}$ is isomorphic to it.

%
%
\subsubsection{Example II: Coisotropic Reduction in the Symplectic Case}
\label{subsec:ExampleIIBRST}

While the previous example makes use of a Lie group symmetry, the
following relies on a coisotropic submanifold only.
However, at the present state, we have to restrict ourselves to a symplectic manifold
$(M, \omega)$.
Thus let $\iota\colon C \longrightarrow M$ be a coisotropic submanifold.
We assume that the classical reduced phase space
$M_\red = C \big/ \mathord{\sim}$ is smooth with the projection map
$\pi\colon C \longrightarrow M_\red$ being a surjective submersion.
In other words, we consider a constraint symplectic manifold
$\mathcal{M} = (M,C,D)$,
with $D$ the characteristic distribution as before,
and $\mathcal{M}_\red = M_\red$.
It follows that there is a unique symplectic form
$\omega_\red$ on $\mathcal{M}_\red$ with
$\pi^*\omega_\red = \iota^*\omega$.
We follow closely the construction of Bordemann in 
\cite{bordemann:2005a,bordemann:2004a:pre}
to construct a deformation of the classical constraint algebra
$\ConCinfty(\mathcal{M}) = 
(\Cinfty(M), \algebra{B}_C,\algebra{J}_C )$
as before.

To this end, one considers the product
$M \times \mathcal{M}_\red^-$
with the symplectic structure
$\pr^*_M \omega - \pr_{\mathcal{M}_\red}^*\omega_\red$.
Then
\begin{equation} \label{eq:CintoMMred}
	I \colon C 
	\ni p \mapsto
	(\iota(p), \pi(p))
	\in M \times \mathcal{M}_\red
\end{equation}
is an embedding of $C$ as a Lagrangian submanifold.
By Weinstein's Lagrangian neighbourhood theorem \cite{weinstein:1971a}
one has a tubular neighbourhood
$U \subseteq M \times \mathcal{M}_\red$
and an open neighbourhood $V \subseteq T^*C$
of the zero section $\iota_C\colon C \longrightarrow T^*C$ in the
cotangent bundle $\pi_C\colon T^*C \longrightarrow C$ with a
symplectomorphism $\Psi\colon U \longrightarrow V$, where $T^*C$ 
is equipped with its canonical symplectic structure, such that
$\Psi \circ I = \iota_C$.

In the symplectic case, star products $\star$ are classified by 
their characteristic or Fedosov class $c(\star)$ in
$\HdR^2(M, \field{C})\formal{\lambda}$.
The assumption of having a smooth reduced phase space
allows us now to choose star products $\star$ on
$M$ and $\star_\red$ on $\mathcal{M}_\red$ in such a way that
$\iota^*c(\star\at{U}) = \pi^* c(\star_\red)$.
Note that this is a non-trivial condition on the relation between
$\star$ and $\star_\red$ which, nevertheless, always has solutions.
Given such a matching pair we have a star product $\star \tensor \star_\red^\opp$
on $M \times \mathcal{M}_\red^-$
by taking the tensor product of the individual ones.
Note that we need to take the opposite star product on the
second factor as we also took the negative of $\omega_\red$ 
needed to have a Lagrangian embedding in \eqref{eq:CintoMMred}.
It follows that the characteristic class
$c\big((\star \tensor \star_\red^\opp)\at{U}\big) = 0$ is trivial.

On the cotangent bundle $T^*C$ the choice of a covariant 
derivative induces a standard-ordered star product $\star_\std$
together with a left module structure on $\Cinfty(C)\formal{\lambda}$
via the corresponding symbol calculus, see
\cite{bordemann.neumaier.waldmann:1998a}.
The characteristic class of
$\star_\std$ is known to be trivial, $c(\star_\std) = 0$, see
\cite{bordemann.neumaier.pflaum.waldmann:2003a}.
Hence the pullback star product $\Psi^*(\star_\std\at{V})$
is equivalent to
$(\star \tensor \star_\red^\opp)\at{U}$.
Thus we find an equivalence transformation between
$\Psi^*(\star_\std)$ and $\star \tensor \star_\red$
on the tubular neighbourhood $U$.
Using this, we can also pullback the left module structure to obtain a 
left module structure on $\Cinfty(C)\formal{\lambda}$ for the algebra
$\Cinfty(M \times M_\red)\formal{\lambda}$. Note that here we 
even get
an extension to all functions since the left module structure with
respect to $\star_\std$ coming from the symbol calculus is by
differential operators and $\Psi \circ I = \iota_C$. Hence the 
module structure with respect to $\star \tensor \star_\red^\opp$ is by
differential operators as well. This ultimately induces a left 
module structure $\acts$ on $\Cinfty(C)\formal{\lambda}$ with respect to
$\star$ and a right module structure $\racts$ with respect to
$\star_\red$ such that the two module structures commute: We have 
a bimodule structure. Moreover, it is easy to see that the module
endomorphisms of the left $\star$-module are given by the right
multiplications with functions from 
$\Cinfty(\mathcal{M}_\red)\formal{\lambda}$,
i.e.
\begin{equation} \label{eq:EndosCinftyC}
	\End_{(\Cinfty(M)\formal{\lambda}, \star)}
	(\Cinfty(C)\formal{\lambda})^\opp
	\cong
	\Cinfty(\mathcal{M}_\red)\formal{\lambda}.
\end{equation}
Moreover, one can construct from the above equivalences a formal
series $S = \id + \sum_{r=1}^\infty \lambda^r S_r$ of differential
operators $S_r$ on $M$ such that the left module structure is 
given by
\begin{equation}
	\label{eq:LeftModuleOnC}
	f \acts \psi
	= \iota^*(S(f) \star \prol(\psi)),
\end{equation}
for $f \in \Cinfty(M)\formal{\lambda}$ and
$\psi \in \Cinfty(C)\formal{\lambda}$, where
$\prol\colon \Cinfty(C)\formal{\lambda} \longrightarrow
\Cinfty(M)\formal{\lambda}$ is the prolongation coming from the
tubular neighbourhood $U$
and the choice of a bump function.

The left module structure is cyclic with cyclic vector
$1 \in \Cinfty(C)\formal{\lambda}$.
This means that
\begin{equation}
	\label{eq:KerModule}
	\qvanishing_C
	= 	\left\{	f \in \Cinfty(M)\formal{\lambda} \;\big|\; f \acts 1 = 0 \right\}
\end{equation}
is a left $\star$-ideal and
$\Cinfty(C)\formal{\lambda} \cong \Cinfty(M)\formal{\lambda} \big/ \qvanishing_C$
as left $\star$-modules.
Moreover, the normalizer
\begin{equation} \label{eq:qBCDef}
	\qalgebra{B}_C = \normalizer_\star(\qvanishing_C)
\end{equation}
with respect to $\star$ gives first
$\qalgebra{B}_C \big/ \qvanishing_C \cong \End_{(\Cinfty(M)\formal{\lambda}, \star)}(\Cinfty(C)\formal{\lambda})^\opp$
for general reasons.
Then this yields the algebra isomorphism
$\qalgebra{B}_C \big/ \qvanishing_C \cong \Cinfty(M_\red)\formal{\lambda}$.

Thanks to the explicit formula for $\acts$ we can use the series 
$S$ to pass to a new equivalent star product $\star'$ such that
$\qvanishing'_C = \vanishing_C\formal{\lambda}$.
We see that this brings us precisely in the situation of
\autoref{exmp:DeformationStrConAlg}:
The constraint algebra
$\qalgebra{A} = (\Cinfty(M)\formal{\lambda}, \qalgebra{B}_C, \qvanishing_C)$
is isomorphic to a deformation of the classical constraint algebra
$\ConCinfty(\mathcal{M})$ we started with.
Note that it might not be directly a deformation of
$\ConCinfty(\mathcal{M})$ as we still might have to untwist first 
$\qvanishing_C$
using $S$ and then $\qalgebra{B}_C$ as in
\autoref{exmp:DeformationStrConAlg}.
This way we can give a re-interpretation of Bordemann's construction in the language of
deformations of constraint algebras.

%
%

\section{Constraint Deformation Functor}
\label{sec:ConDeformationFunctor}

By a well-known principle of classical deformation theory, a
deformation problem is controlled by a certain differential 
graded Lie algebra, see e.g. \cite{manetti:2009a}.
Thus, the first step to discuss the deformation theory of constraint 
algebras consists in introducing a constraint deformation functor for
a constraint DGLA.
For this we will need constraint Maurer-Cartan elements and a
notion of gauge equivalence. 

\index{Maurer-Cartan element}
Recall that a Maurer-Cartan element in a DGLA
$\liealg{g}^\bullet$ is an element $\xi \in \liealg{g}^1$ 
satisfying the Maurer-Cartan equation
\begin{equation}
	\D \xi + \frac{1}{2}[\xi, \xi] = 0.
\end{equation}
While up to here we did not have to make any further assumption 
about the ring $\field{k}$ of scalars, from now on we assume
$\field{Q} \subseteq \field{k}$ in order to have a well-defined
Maurer-Cartan equation and gauge action later on.
We denote by $\MCset(\liealg{g})$ the set of all Maurer-Cartan 
elements of a DGLA.

\begin{definition}[Constraint set of Maurer-Cartan elements]
	\label{def:ConMCElements}%
	\index{constraint!Maurer-Cartan element}
	\index{constraint!differential graded Lie algebra}
Let $\liealg{g}$ be a constraint\linebreak
DGLA.
The \emph{constraint set $\MCset(\liealg{g})$ of Maurer-Cartan 
	elements} of $\liealg{g}$ is given by
\glsadd{MaurerCartan}
\begin{equation}
	\label{eq:MCset}
	\MCset(\liealg{g})
	= \big(	\MCset(\liealg{g}_\Total), \; \MCset(\liealg{g}_\Wobs), \;
	\sim_\MC
	\big),
\end{equation}
together with
$\iota_\MC \colon \MCset(\liealg{g}_\Wobs) \to 
\MCset(\liealg{g}_\Total)$
given by the map
$\iota_\liealg{g} \colon \liealg{g}^\bullet_\Wobs \to 
\liealg{g}^\bullet_\Total$
of $\liealg{g}$ and where the relation $\sim_\MC$ is defined by
\begin{equation}
	\label{eq:MCEquivalence}
	\xi_1 \sim_\MC \xi_2 \iff \xi_1 - \xi_2 \in \liealg{g}_\Null^1
\end{equation}
for $\xi_1, \xi_2 \in \MCset(\liealg{g}_\Wobs)$.
\end{definition}

\begin{example}[Constraint multivector fields]
	\label{ex:ConPoissonStructureAsMCElements}
	\index{constraint!multivector fields}
Let $\mathcal{M} = (M,C,D)$ be a constraint manifold.
By \autoref{cor:ConDGLAOfConMultVectFields} we know that
$(\ConVecFields_{\strtensor}^{\bullet+1}(\mathcal{M}), \D = 0, \Schouten{\argument, \argument})$
is a constraint DGLA.
Then $\MCset(\ConVecFields_{\strtensor}^\bullet(\mathcal{M}))_\Total$
is the set of Poisson structures on $M$, 
and, by \autoref{def:ConPoissonManifold}, 
$\MCset(\ConVecFields_{\strtensor}^\bullet(\mathcal{M}))_\Wobs$
is exactly the set of constraint Poisson structures on $\mathcal{M}$.
Two such constraint Poisson structures
$\pi_1$ and $\pi_2$ are equivalent as Maurer-Cartan elements if
and only if $\pi_1 - \pi_2 \in \ConVecFields_{\strtensor}^\bullet(\mathcal{M})_\Null$,
i.e. if at least one leg of the bivector $\pi_1 - \pi_2$ points into the direction
of the distribution, and therefore the bivector vanishes after reduction,
c.f. \autoref{lem:LocalConBiVect} \ref{lem:LocalConBiVect_2}.
\end{example}

\begin{lemma}[Maurer-Cartan functor]
	\label{lem:MaurerCartanFunctor}%
	\index{constraint!Maurer-Cartan functor}
Mapping constraint DGLAs to their constraint sets of Maurer-Cartan 
elements defines a functor
\begin{equation}
	\MCset\colon \ConDGLieAlg \to \ConSet.
\end{equation}
\end{lemma}

\begin{proof}
Every morphism
$\Phi \colon \liealg{g} \to \liealg{h}$
of constraint DGLAs induces maps
$\Phi_\Total \colon \MCset(\liealg{g}_\Total) \to 
\MCset(\liealg{h}_\Total)$
and
$\Phi_\Wobs \colon \MCset(\liealg{g}_\Wobs) \to 
\MCset(\liealg{h}_\Wobs)$.
Moreover, since
$\Phi_\Wobs \colon \liealg{g}_\Wobs \to \liealg{h}_\Wobs$
preserves the $\NULL$-component its induced map on
$\MCset(\liealg{g}_\Wobs)$ maps equivalent elements to equivalent 
elements.
\end{proof}

As in the setting of classical DGLAs, for a given constraint DGLA
$(\liealg{g}, [\argument, \argument], \D)$
and a given Maurer-Cartan element
$\xi_0 \in \MCset(\liealg{g})_\Wobs$
we can always obtain a twisted constraint DGLA by
$\liealg{g}_{\xi_0} = (\liealg{g}, [\argument,\argument], \D_{\xi_0})$
with
\begin{equation}
	\label{eq:TwistedDifferentialMC}
	\D_{\xi_0} := \D + [\xi_0, \argument].
\end{equation}
Here we are using the tensor-hom adjunction in the sense of
\eqref{eq:InternalTensorHomCModk}.

Note that for any constraint DGLA $\liealg{g}$ and constraint algebra 
$\algebra{A}$ the tensor product 
$\liealg{g} \tensor \algebra{A}$
is again a constraint DGLA by the usual construction.
For this observe that
$\liealg{g}_\Null \tensor \algebra{A}_\Wobs + \liealg{g}_\Wobs 
\tensor \algebra{A}_\Null$
is indeed a Lie ideal in
$\liealg{g}_\Wobs \tensor \algebra{A}_\Wobs$.

Reformulating the equivalence of deformations of a given Maurer-Cartan
element in terms of its twisted constraint DGLA requires a notion of
a constraint gauge group.
To define this we either need to assume that the DGLA we are starting 
with has additional properties, e.g. being nilpotent, or we can use 
formal power series instead.
Since later on we are interested in formal deformation theory, we 
will choose the latter option.
It is now easy to see that
$\liealg{g}\formal{\lambda}$
is a constraint DGLA for any constraint DGLA $\liealg{g}$ by 
$\lambda$-linear extension of all structure maps.

Note that the gauge action will require to have
$\field{Q} \subseteq \field{k}$ since we need the (formal) exponential
series and the (formal) Baker-Campbell-Hausdorff (BCH) series.

\begin{proposition}[Gauge group]
	\label{prop:ConGaugeGroup}%
	\index{gauge group}
	\glsadd{GaugeGroup}
Let $\liealg{g}$ be a constraint Lie algebra.
Then
\begin{equation}
	\functor{G}(\liealg{g})= 
	\big(\lambda\liealg{g}_\Total\formal{\lambda},\,\,
	\lambda\liealg{g}_\Wobs\formal{\lambda},\,\,
	\lambda\liealg{g}_\Null\formal{\lambda} \big)
\end{equation}
with multiplication $\bullet$ given by the Baker-Campbell-Hausdorff
formula \cite[Eq. 2.4.8.]{esposito:2015a}
\begin{equation}
	\lambda \xi \bullet \lambda \eta 
	= \lambda \xi + \lambda \eta + \frac 12 [\lambda \xi ,\lambda \eta] + \cdots
\end{equation}
is a constraint group.
\end{proposition}

\begin{proof}
	The additional prefactor $\lambda$ makes all the BCH series
	$\lambda$-adically convergent.
	The well-known group structures on
	$\liealg{g}_\Total\formal{\lambda}$
	and
	$\liealg{g}_\Wobs\formal{\lambda}$
	are given by the BCH formula and we clearly have a group morphism
	$\liealg{g}_\Wobs\formal{\lambda} \to 
	\liealg{g}_\Total\formal{\lambda}$.
	Finally, we need to show that
	$\lambda\liealg{g}_\Null\formal{\lambda}$
	is a normal subgroup of
	$\lambda\liealg{g}_\Wobs\formal{\lambda}$.
	For this let
	$\lambda g \in \lambda\liealg{g}_\Wobs\formal{\lambda}$
	and
	$\lambda h \in \lambda\liealg{g}_\Null\formal{\lambda}$
	be given.
	Since by the BCH formula
	$\lambda g \bullet \lambda h \bullet (\lambda g)^{-1} 
	= \lambda g_0 + \lambda h_0 - \lambda g_0 + \lambda^2(\cdots)$,
	where all higher order terms are given by Lie brackets and
	$\liealg{g}_\Null$
	is a Lie ideal in $\liealg{g}_\Wobs$, we see that
	$\lambda g \bullet \lambda h \bullet (\lambda g)^{-1}
	\in \lambda\liealg{g}_\Null\formal{\lambda}$.
\end{proof}

By abuse of notation we will write
$\functor{G}(\liealg{g}) = \functor{G}(\liealg{g}^0)$
for every constraint DGLA $\liealg{g}$.
With the composition $\bullet$ on
$\functor{G}(\liealg{g})$
defined by the Baker-Campbell-Hausdorff formula it is immediately 
clear that every morphism
$\Phi \colon \liealg{g} \to \liealg{h}$
of constraint DGLAs induces a morphism
$\functor{G}(\Phi) \colon \functor{G}(\liealg{g})
\to \functor{G}(\liealg{h})$
of the corresponding gauge groups, given by the $\lambda$-linear 
extension of $\Phi$.
In other words, we obtain a functor
$\functor{G} \colon \ConDGLieAlg \to \ConGroups$.

The usual gauge action of the formal group on the (formal) 
Maurer-Cartan elements can now be extended to a constraint DGLA as 
follows:

\begin{proposition}[Gauge action]
	\label{prop:GaugeAction}%
	\index{gauge action}
Let
$(\liealg{g}, [\argument,\argument], \D)$
be a constraint DGLA.
Then the constraint group
$\functor{G}(\liealg{g})$
acts on the constraint set
$\MCset(\lambda\liealg{g}\formal{\lambda})$
by
\begin{equation}
	\label{eq:ActsTotal}
	\lambda g \acts_\Total \xi
	\coloneqq \E^{\lambda \ad_\Total(g)}(\xi) -
	\lambda \sum_{k=0}^{\infty}
	\frac{(\lambda \ad_\Total(g))^k}{(1+k)!}(\D_\Total g)
\end{equation}
for
$\lambda g \in \functor{G}(\liealg{g})_\Total$
and
$\xi \in \MCset(\lambda\liealg{g}\formal{\lambda})_\Total$
as well as
\begin{equation}
	\label{eq:ActsWobs}
	\lambda g \acts_\Wobs \xi
	\coloneqq
	\E^{\lambda \ad_\Wobs(g)}(\xi) - 
	\lambda \sum_{k=0}^{\infty}
	\frac{(\lambda \ad_\Wobs(g))^k}{(1+k)!}(\D_\Wobs g)
\end{equation}
for
$\lambda g \in \functor{G}(\liealg{g})_\Wobs$
and
$\xi \in \MCset(\lambda\liealg{g}\formal{\lambda})_\Wobs$.
\end{proposition}

\begin{proof}
	Clearly, $\acts_\Total$ and $\acts_\Wobs$ define actions of
	$\functor{G}(\liealg{g})_\Total$
	and 
	$\functor{G}(\liealg{g})_\Wobs$
	on
	$\MCset(\lambda\liealg{g}\formal{\lambda})_\Total$
	and
	$\MCset(\lambda\liealg{g}\formal{\lambda})_\Wobs$, 
	respectively,
	by classical results, see \cite{esposito:2015a}.
	Moreover, writing out the exponential series and using the fact that
	$\ad(g) = [g, \argument]$
	and $\D$ commute with $\iota_\liealg{g}$
	directly yields
	\begin{align*}
		\iota_\liealg{g}
		\left(\lambda g \acts_\Wobs \xi\right)
		&= \E^{\lambda 
			\ad_\Total(\iota_\liealg{g}(g))}(\iota_\liealg{g}(\xi))
		- \lambda \sum_{k=0}^{\infty}
		\frac{\big(\lambda \ad_\Total(\iota_\liealg{g}(g))\big)^k}{(1+k)!}
		(\D_\Total \iota_\liealg{g}(g))	\\
		&= \lambda \iota_\liealg{g}(g) \acts_\Total \iota_\liealg{g}(\xi).
	\end{align*}
	Finally, we have for any
	$\lambda g \in \functor{G}(\liealg{g})_\Null$
	and
	$\xi \in \MCset(\lambda\liealg{g}\formal{\lambda})_\Wobs$
	\begin{align*}
		\E^{\lambda \ad_\Wobs(g)}(\xi) - \xi
		&= \sum_{k=0}^{\infty} \frac{\lambda^k}{k!} (\ad_\Wobs(g))^k(\xi)
		- \lambda \sum_{k=0}^{\infty}
		\frac{(\lambda \ad_\Wobs(g))^k}{(1+k)!}(\D_\Wobs g)
		- \xi \\
		&= \sum_{k=1}^{\infty} \frac{\lambda^k}{k!} (\ad_\Wobs(g))^k (\xi)
		- \lambda \sum_{k=0}^{\infty}
		\frac{(\lambda \ad_\Wobs(g))^k}{(1+k)!}(\D_\Wobs g)
		\in \lambda\liealg{g}_\Null\formal{\lambda},
	\end{align*}
	since
	$\D_\Wobs g \in \liealg{g}_\Null\formal{\lambda}$
	and
	$\ad_\Wobs(g)(\xi) \in \liealg{g}_\Null\formal{\lambda}$.
\end{proof}

This shows that the constraint sets of Maurer-Cartan elements admit 
more structure, namely that of an action of the associated gauge 
group.
This suggests that the functor $\MCset$ of 
\autoref{lem:MaurerCartanFunctor}
factors through $\ConGroupAct$, cf. \autoref{def:ConGroupAction}.

\begin{corollary}
	Mapping constraint DGLAs
	$(\liealg{g}, [\argument, \argument], \D)$
	to their corresponding gauge action of $\functor{G}(\liealg{g})$
	on $\MCset(\lambda \liealg{g}\formal{\lambda})$
	defines a functor
	$\MCset \colon \ConDGLieAlg \to \ConGroupAct$.	
\end{corollary}

\begin{proof}
	Let $\Phi \colon \liealg{g} \to \liealg{h}$
	be a morphism of constraint DGLAs.
	Its $\lambda$-linear extension gives morphisms
	$\Phi \colon \MCset(\lambda\liealg{g}\formal{\lambda})
	\to \MCset(\lambda\liealg{h}\formal{\lambda})$
	and
	$\Phi \colon \functor{G}(\liealg{g}) \to \functor{G}(\liealg{h})$.
	With this we get
	\begin{align*}
		\Phi_\Total (\lambda g \acts_\Total \xi)
		&= \Phi_\Total\Big( \E^{\lambda \ad_\Total(g)}(\xi)
		- \lambda \sum_{k=0}^\infty \frac{(\lambda 
			\ad_\Total(g))^k}{(1+k)!} (\D_\Total g) \Big) \\
		&= \E^{\lambda \ad_\Total(\Phi_\Total(g))}(\Phi_\Total(\xi))
		- \lambda \sum_{k=0}^\infty
		\frac{\big(\lambda \ad_\Total(\Phi_\Total(g))\big)^k}{(1+k)!} 
		(\D_\Total \Phi_\Total(g)) \\
		&= \lambda \Phi_\Total(g) \acts_\Total \Phi_\Total(\xi),
	\end{align*}
	for all $\lambda g \in \functor{G}(\liealg{g})_\Total$
	and $\xi \in \MCset(\lambda \liealg{g}\formal{\lambda})_\Total$.
	With an analogous computation we find that also
	$\Phi_\Wobs (\lambda g \acts_\Wobs \xi)
	= \lambda \Phi_\Wobs(g) \acts_\Wobs \Phi_\Total(\xi)$
	holds for all
	$\lambda g \in \functor{G}(\liealg{g})_\Wobs$
	and $\xi \in \MCset(\lambda \liealg{g}\formal{\lambda})_\Total$,
	showing that $\Phi$ is an equivariant map.
\end{proof}

Maurer-Cartan elements are said to be \emph{equivalent} if they 
lie in the same orbit of the gauge action.
Hence the object of interest for deformation theory is not the set of 
Maurer-Cartan elements itself but its set of equivalence classes.
More precisely, let us denote by
\index{constraint!deformation functor}
\glsadd{DeformationFunctor}
\begin{equation}
	\Def(\liealg{g}) \coloneqq 
	\MCset(\lambda\liealg{g}\formal{\lambda})
	/ \functor{G}(\liealg{g})
\end{equation}
the orbit space of the gauge action of the gauge group 
$\functor{G}(\liealg{g})$ 
on the constraint set 
$\MCset(\lambda\liealg{g}\formal{\lambda})$
of Maurer-Cartan elements.
The corresponding functor 
$\Def \colon \ConDGLieAlg \to \ConSet$
is called \emph{deformation functor}.

\subsubsection{Reduction}

The question arises if the above constructions of the constraint 
set of Maurer-Cartan elements, the constraint gauge group and the
deformation functor commute with reduction. 
The next theorem shows that this is partially true, in the sense that 
at least an injective natural transformation exists,
see \cite[Thm. 3.14]{dippell.esposito.waldmann:2022a}.

\begin{theorem}[Deformation functor vs. reduction]\
	\label{thm:defvsred}%
	\index{reduction!deformation functor}
\begin{theoremlist}
	\item \label{item:MCredNatural}
	There exists an injective natural transformation
	$\eta \colon \red \circ \MCset
	\Longrightarrow \MCset \circ \red$, i.e.
	\begin{equation}
		\begin{tikzcd}
			\ConDGLieAlg
			\arrow{r}{\MCset}
			\arrow{d}[swap]{\red}
			& \ConSet
			\arrow{d}{\red}
			\arrow[Rightarrow]{dl}[swap]{\eta}\\
			\DGLieAlg
			\arrow{r}{\MCset}
			&\Sets
		\end{tikzcd}
	\end{equation}
	commutes with $\eta$ injective.
	\item \label{item:GredNatural} 
	There exists a natural isomorphism
	such that the diagram
	\begin{equation}
		\begin{tikzcd}
			\ConDGLieAlg
			\arrow{r}{\functor{G}}
			\arrow{d}[swap]{\red}
			& \ConGroups
			\arrow{d}{\red}
			\arrow[Rightarrow]{dl}[swap]{\eta}\\
			\DGLieAlg
			\arrow{r}{\functor{G}}
			&\Groups
		\end{tikzcd}
	\end{equation}
	commutes with $\eta$ bijective.
	\item \label{item:DefredNatural}
	There exists an injective natural transformation
	$\eta \colon \red \circ \Def 
	\Longrightarrow \Def \circ \red$, i.e.
	\begin{equation}
		\begin{tikzcd}
			\ConDGLieAlg
			\arrow{r}{\Def}
			\arrow{d}[swap]{\red}
			& \ConSet
			\arrow{d}{\red}
			\arrow[Rightarrow]{dl}[swap]{\eta}\\
			\DGLieAlg
			\arrow{r}{\Def}
			&\Sets
		\end{tikzcd}
	\end{equation}
	commutes with $\eta$ injective.
\end{theoremlist}
\end{theorem}

\begin{proof}
	For this proof we need to construct natural transformations
	$\eta$, consisting of $\TOTAL$- and $\WOBS$-components.
	Since the computations are identical in both cases, we omit the subscripts.
	\begin{theoremlist}
		\item In the following we denote by $[\argument]_\MC$ the
		equivalence classes of elements in 
		$\MCset(\liealg{g}_\Wobs)$
		and by $[\argument]_\liealg{g}$ the equivalence classes of
		elements in $\liealg{g}_\Wobs$.
		For any constraint DGLA
		$\liealg{g}$ define
		\begin{equation*}
			\eta_\liealg{g} \colon \MCset(\liealg{g})_\red
			\to	\MCset(\liealg{g}_\red)
			\qquad\text{ by }\qquad
			\eta_\liealg{g}([\xi]_\MC) = [\xi]_\liealg{g}.
		\end{equation*}
		This map is well-defined since
		$[\xi]_\MC \subseteq [\xi]_\liealg{g}$ 
		and
		\begin{equation*}
			\D_\red [\xi]_\liealg{g}
			+ \big[ [\xi]_\liealg{g}, [\xi]_\liealg{g} \big]_\red
			= \big[ \D_\Wobs \xi + [\xi,\xi]_\Wobs \big]_\liealg{g}
			= [0]_\liealg{g}
		\end{equation*}
		for every
		$\xi \in \MCset(\liealg{g}_\Wobs)$.
		To show that $\eta_\liealg{g}$ is injective let
		$[\xi_1]_\MC, [\xi_2]_\MC \in \MCset(\liealg{g})_\red$ be
		given such that
		$[\xi_1]_\liealg{g} = [\xi_2]_\liealg{g}$.
		Then $\xi_2 \in [\xi_1]_\liealg{g}$ and hence
		$\xi_1-\xi_2 \in \liealg{g}_\Null^1$. 
		Thus by definition
		$\xi_1 \sim_\MC \xi_2$
		and therefore
		$[\xi_1]_\MC = [\xi_2]_\MC$.
		To show naturality of $\eta$ let a morphism
		$\Phi \colon \liealg{g} \to \liealg{h}$
		of constraint DGLAs be given.
		This induces morphisms
		$\Phi \colon \MCset(\liealg{g})_\red 
		\to	\MCset(\liealg{h})_\red$
		and
		$\Phi \colon \MCset(\liealg{g}_\red)
		\to	\MCset(\liealg{h}_\red)$
		by applying $\Phi_\Wobs$ to	representatives.
		Then we have
		\begin{equation*}
			(\eta_\liealg{h} \circ \Phi)([\xi]_\MC)
			= \eta_\liealg{h}([\Phi_\Wobs(\xi)]_\MC)
			= [\Phi_\Wobs(\xi)]_\liealg{h}
			= \Phi([\xi]_\liealg{g})
			= \Phi(\eta_\liealg{g}([\xi]_\MC)),
		\end{equation*}
		showing that $\eta$ is natural.
		\item Then
		$\eta_\liealg{g} \colon \functor{G}(\liealg{g})_\red
		\to	\functor{G}(\liealg{g}_\red)$
		given by
		$[\lambda g]_{\functor{G}}\mapsto \lambda [g]_\liealg{g}$,
		where $[g]_\liealg{g}$ denotes the equivalence class of $g$ in
		$\liealg{g}_\red$, is well-defined.
		Indeed, $\eta_\liealg{g}$ is just the $\lambda$-linear extension 
		of the obvious identity
		$\liealg{g}_\Wobs / \liealg{g}_\Null = \liealg{g}_\red$.
		Moreover, $\eta_\liealg{g}$ is a group morphism, since
		$[\argument]_\liealg{g} \colon \liealg{g}_\Wobs
		\to	\liealg{g}_\red$
		is a morphism of DGLAs and $\bullet$ is given by sums of iterated 
		brackets.
		Naturality follows directly.
		\item By definition $\Def$ factors as
		$\Def = \ConOrb \circ \MCset$, 
		with functors $\MCset \colon \ConDGLieAlg \to \ConGroupAct$
		and $\ConOrb \colon \ConGroupAct \to \ConSet$
		as in \autoref{prop:OrbitSpaceFunctor}.
		By \autoref{prop:OrbitSpaceVSReduction} $\ConOrb$ commutes with 
		reduction,
		so we only need to consider $\MCset$.
		For this we show that $\eta$ from \ref{item:MCredNatural}
		is equivariant:
		\begin{equation*}
			\eta_\liealg{g}\left( [\lambda g]_\functor{G} \acts [\xi]_\MC 
			\right)
			= \eta_\liealg{g}([\lambda g \acts \xi]_\MC)
			= [\lambda g \acts \xi]_\liealg{g}
			= [\lambda g]_\functor{G} \acts [\xi]_\liealg{g}.
		\end{equation*}
		Here we implicitly used \ref{item:GredNatural}.
		Now composing $\eta$ with the natural isomorphism from
		\autoref{prop:OrbitSpaceVSReduction}
		yields the wanted injective natural transformation.
	\end{theoremlist}
\end{proof}

The missing surjectivity in \autoref{thm:defvsred} \ref{item:MCredNatural}
comes again from the fact that the reduction functor does not 
reflect limits, cf. \autoref{rem:ReductionConSets}.

%
%

\section{Constraint Hochschild Cohomology}
\label{sec:ConHochschildCohomology}

We now want to introduce a constraint version of Hochschild cohomology for associative algebras.
This constraint Hochschild complex will turn out to be the constraint DGLA which controls the deformation problem of constraint algebras.

In this section we assume that $\field{Q} \subseteq \field{k}$.
\index{constraint!Hochschild complex}
Let $\module{M}, \module{N} \in \ConMod_\field{k}$
be constraint $\field{k}$-modules.
We define for any $n \in \Naturals$
\glsadd{HochschildComplex}
\begin{equation} \label{eq:CnDef}
	\mathrm{C}^n(\module{M}, \module{N})
	\coloneqq \ConHom_\field{k}(\module{M}^{\tensor n},\module{N})
\end{equation}
with $\ConHom_\field{k}$ denoting the internal hom as usual.
Recall that
\begin{align*}
	\mathrm{C}^n(\module{M},\module{N})_\Total
	&= \Hom_\field{k}(\module{M}_\Total^{\tensor n}, \module{N}_\Total),\\
	\mathrm{C}^n(\module{M},\module{N})_\Wobs
	&= \Hom_\field{k}(\module{M}^{\tensor n},\module{N}), \\
	\mathrm{C}^n(\module{M},\module{N})_\Null
	&= \big\{ (f_\Total,f_\Wobs) 
	\in \Hom_{\field{k}}(\module{M}^{\tensor n},\module{N})
	\mathbin{\big|}
	f_\Wobs(\module{M}^{\tensor n}_\Wobs)
	\subseteq \module{N}_\Null \big\},
\end{align*}
with
$\iota_n \colon \mathrm{C}^n(\module{M},\module{N})_\Wobs
\ni (f_\Total,f_\Wobs) \mapsto f_\Total
\in \mathrm{C}^n(\module{M},\module{N})_\Total$.
Note that a morphism
$f = (f_\Total, f_\Wobs) \in \mathrm{C}^n(\module{M},\module{N})_\Wobs$ 
fulfils
$f_\Wobs((\module{M}^{\tensor n})_\Null) \subseteq \module{N}_\Null$
where, by definition of the tensor product, we have
\begin{equation} \label{eq:NullOfTensorPower}
	(\module{M}^{\tensor n})_\Null
	= \sum_{i=1}^{n} \module{M}^{\tensor i-1}_\Wobs \tensor	\module{M}_\Null
	\tensor \module{M}_\Wobs^{\tensor n-i}.
\end{equation}
In other words, $f_\Wobs$ maps to $\module{N}_\Null$ if at least one tensor 
factor comes from $\module{M}_\Null$.
This clearly defines a graded constraint $\field{k}$-module
$\mathrm{C}^\bullet(\module{M},\module{N})$.
Since $\module{M}^{\tensor 0} \simeq (\field{k}, \field{k},0)$
it holds $C^0(\module{M},\module{N}) = \module{N}$.

Let us now consider the case $\module{N} = \module{M}$.
Then we write
$\mathrm{C}^\bullet(\module{M}) 
= \mathrm{C}^\bullet(\module{M},\module{M})$.
We now want to transfer the Gerstenhaber algebra structure of the classical 
Hochschild complex to $\mathrm{C}^\bullet(\module{M})$.
For this denote by
$[\argument,\argument]^{\module{M}_\Total}$
and
$[\argument,\argument]^{\module{M}_\Wobs}$
the Gerstenhaber brackets for the modules
$\module{M}_\Total$ and $\module{M}_\Wobs$, respectively.
Then we need to show that
$[\argument,\argument]^{\module{M}_\Wobs}$
preserves the $\NULL$-components.
This follows directly from the usual formula for the Gerstenhaber bracket,
see \cite{gerstenhaber:1963a}.

\begin{definition}[Gerstenhaber bracket]
	\label{definition:CoisoGerstenhaberBracket}%
	\index{constraint!Gerstenhaber bracket}
	\index{LieBracket}
Let $\module{M} \in \ConMod_\field{k}$.
Then the morphism
\begin{equation}
	[\argument, \argument]
	\colon \mathrm{C}^\bullet(\module{M}) \tensor \mathrm{C}^\bullet(\module{M})
	\to \mathrm{C}^\bullet(\module{M})
\end{equation}
of constraint $\field{k}$-modules defined by
\begin{equation}
	[\argument, \argument]_\Total =
	[\argument,\argument]^{\module{M}_\Total}
	\quad\text{ and }\quad
	[\argument, \argument]_\Wobs
	=\big([\argument,\argument]^{\module{M}_\Total},
	[\argument,\argument]^{\module{M}_\Wobs}\big)
\end{equation}
is called the
\emph{constraint Gerstenhaber bracket}.
\end{definition}

Since
$[\argument, \argument]^{\module{M}_\Total}$
and
$[\argument, \argument]^{\module{M}_\Wobs}$
induce graded Lie algebra structures on the classical Hochschild complexes 
of
$\module{M}_\Total$ and $\module{M}_\Wobs$
it is easy to see that
$\mathrm{C}^\bullet(\module{M})$
together with the constraint Gerstenhaber bracket $[\argument, \argument]$ 
forms a graded constraint Lie algebra.

\begin{remark}
	\label{remark:GerstenhaberFromPreLie}%
The constraint Gerstenhaber bracket can also be derived from a constraint 
pre-Lie algebra structure on $\mathrm{C}^\bullet(\module{M})$,
which in turn results from a sort of partial composition.
These partial compositions can be interpreted as the usual endomorphism 
operad structure of $\module{M}$ in $\ConMod_\field{k}$.
\end{remark}
	
As in the classical theory of deformations of associative algebras, we can 
characterize associative multiplications by using the Gerstenhaber bracket.
\begin{lemma}
	\label{lem:AssoConMultiplication}%
Let $\module{M} \in \ConMod_\field{k}$
be a constraint module.
Then a morphism
$\mu \colon \module{M} \tensor \module{M} \to \module{M}$
of constraint $\field{k}$-modules is an associative constraint algebra 
structure on $\module{M}$ if and only if
\begin{equation} \label{eq:AssociativityViaMC}
	[\mu, \mu]_\Wobs = 0.
\end{equation}
\end{lemma}

\begin{proof}
First, note that a constraint morphism
$\mu\colon \module{M} \tensor \module{M} \to \module{M}$
is an element in $\mathrm{C}^2(\module{M})_\Wobs$
and hence consists of a pair $(\mu_\Total, \mu_\Wobs)$
and
$[\argument, \argument]_\Wobs
= ([\argument,\argument]^{\module{M}_\Total},
[\argument,\argument]^{\module{M}_\Wobs})$.
From the classical theory for associative algebras we know that
$\mu_\Total$ and $\mu_\Wobs$
are associative multiplications if and only if
$[\mu_\Total,\mu_\Total]^{\module{M}_\Total} = 0$
and
$[\mu_\Wobs, \mu_\Wobs]^{\module{M}_\Wobs} = 0$
hold.
\end{proof}

Note that \eqref{eq:AssociativityViaMC} only involves the
$\WOBS$-component of the constraint Gerstenhaber bracket
$[\argument, \argument]$.
Using the constraint structure of
$\mathrm{C}^2(\module{M})$ we get
$\iota_2(\mu) = \mu_\Total \in \mathrm{C}^2(\module{M})_\Total$, 
from which directly $[\mu_\Total,\mu_\Total]_\Total = 0$ follows.

Let us now move from a module $\module{M}$ to an algebra
$(\algebra{A},\mu)$.
Then we can use the multiplication to construct a differential on 
$\mathrm{C}^\bullet(\algebra{A})$.

\begin{proposition}[Constraint Hochschild differential]
	\label{prop:ConHochschildDifferential}%
	\index{constraint!Hochschild differential}
	\glsadd{HochschildDifferential}
Let $(\algebra{A}, \mu) \in \ConAlg_\field{k}$ 
be a constraint algebra.
Then the morphism
$\delta\colon \mathrm{C}^\bullet(\algebra{A}) \to
\mathrm{C}^{\bullet +1}(\algebra{A})$
of constraint $\field{k}$-modules, defined by its components
\begin{equation} \label{eq:HochschildDef}
	\delta_\Total = -[\argument, \mu_\Total]_\Total
	\quad\textrm{ and }\quad
	\delta_\Wobs = -[\argument, \mu]_\Wobs,
\end{equation}
is a constraint chain map of degree $1$ with $\delta^2 = 0$.
\end{proposition}

\begin{proof}
Since $\mu_\Total$ is an associative multiplication on
$\algebra{A}_\Total$ we know that
$\delta_\Total \colon \mathrm{C}^\bullet(\algebra{A}_\Total)
\to \mathrm{C}^{\bullet+1}(\algebra{A}_\Total)$ 
is a differential.
Moreover, it is clear that
$\delta_\Wobs \colon \mathrm{C}^\bullet(\algebra{A})_\Wobs
\to \mathrm{C}^\bullet(\algebra{A})_\Wobs$
is also a differential and it preserves the $\NULL$-component by the 
definition of $[\argument, \argument]_\Wobs$.
Finally, we have for
$(\Phi_\Total,\Phi_\Wobs) \in \mathrm{C}^n(\algebra{A})_\Wobs$
that
$(\delta_\Total \circ \iota_n)(\Phi_\Total,\Phi_\Wobs)
= \delta_\Total(\Phi_\Total)
= \iota_{n+1}(\delta_\Wobs((\Phi_\Total,\Phi_\Wobs)))$
holds, and hence
$(\delta_\Total,\delta_\Wobs)$
is a constraint morphism.
\end{proof}

Note that $\delta$ can be understood as 
$\delta = -[\argument, \mu]$
using the tensor-hom adjunction \eqref{eq:InternalTensorHomCModk}.
The constraint Hochschild differential can be interpreted as twisting the 
constraint DGLA
$(\mathrm{C}^\bullet(\algebra{A}), [\argument, \argument],0)$ 
with the Maurer-Cartan element
$\mu \in \mathrm{C}^2(\algebra{A})_\Wobs$, 
but with signs chosen in such a way that it corresponds to the usual
Hochschild differential.
More explicitly we have the following result.

\begin{corollary}
	\label{corollary:HochschildExplicit}%
Let $(\algebra{A},\mu) \in \ConAlg_\field{k}$
be a constraint algebra.
Then the constraint Hochschild differential
$\delta \colon \mathrm{C}^\bullet(\algebra{A})
\to \mathrm{C}^{\bullet +1}(\algebra{A})$
is given by
$\delta = (\delta^{\algebra{A}_\Total},
(\delta^{\algebra{A}_\Total},\delta^{\algebra{A}_\Wobs}))$, 
where
$\delta^{\algebra{A}_\Total}$ and $\delta^{\algebra{A}_\Wobs}$
denote the Hochschild differentials of the algebras
$(\algebra{A}_\Total, \mu_\Total)$and
$(\algebra{A}_\Wobs, \mu_\Wobs)$, respectively.
In particular, for $\phi \in \mathrm{C}^{n+1}$
and $a_0,\dotsc, a_n \in \algebra{A}_\Total$
we have
\begin{equation} \label{eq:HochschildDifferentialFormula}
	\begin{split}
		(\delta \phi)(a_0,\dotsc, a_n)
		&= a_0\phi(a_1,\dotsc, a_n)
		+ (-1)^n \phi(a_0,\dotsc,a_{n-1}) a_n\\
		&\quad+ \sum_{i=0}^{n} (-1)^{i+1} \phi(a_0,\dotsc, a_ia_{i+1},\dotsc, a_n).
	\end{split}	
\end{equation}
\end{corollary}

From this explicit characterization of the constraint Hochschild 
differential in terms of the classical Hochschild differentials it becomes 
clear that
$(\mathrm{C}^\bullet(\algebra{A}),[\argument,\argument],\delta)$ 
is a constraint DGLA.

\begin{definition}[Constraint Hochschild complex]
	\label{definition:CoisoHochschildComplex}%
	\index{constraint!Hochschild complex}
Let $(\algebra{A},\mu) \in \ConAlg_\field{k}$
be a constraint algebra.
The constraint DGLA
$\big(\mathrm{C}^\bullet(\algebra{A}),[\argument,\argument],\delta\big)$
is called the \emph{constraint Hochschild complex} of $\algebra{A}$.
\end{definition}

As we would expect, the constraint Hochschild complex also carries an additional
multiplication, the so-called cup product.

\begin{definition}[Constraint cup product]
	\index{constraint!cup product}
	\glsadd{union}
Let $(\algebra{A},\mu) \in \ConAlg_{\field{k}}$ be a constraint algebra.
The constraint morphism
$\cup \colon \mathrm{C}^\bullet(\algebra{A}) \tensor \mathrm{C}^\bullet(\algebra{A})
\to \mathrm{C}^\bullet(\algebra{A})$
defined by
\begin{equation}
	\phi \cup_\Total \psi \coloneqq
	\mu_\Total \circ (\phi \tensor \psi)
	\qquad\text{and}\qquad
	\phi' \cup_\Wobs \psi' \coloneqq
	(\mu_\Total,\mu_\Wobs) \circ (\phi' \tensor \psi'),
\end{equation}
for $\phi,\psi \in \mathrm{C}^\bullet(\algebra{A})_\Total$
and $\phi',\psi' \in \mathrm{C}^\bullet(\algebra{A})_\Wobs$,
is called the \emph{constraint cup product}.
\end{definition}

Let us quickly summarize the properties for the constraint cup product.

\begin{proposition}
Let $(\algebra{A},\mu) \in \ConAlg_{\field{k}}$ be a constraint algebra.
\begin{propositionlist}
	\item The cup product $\cup$ turns $\mathrm{C}^\bullet(\algebra{A})$
	into a graded constraint algebra.
	\item If $\algebra{A}$ is a strong constraint algebra, then 
	$\mathrm{C}^\bullet(\algebra{A})$ is a strong constraint algebra
	with respect to $\cup$.
	\item The Hochschild differential $\delta$ is a graded derivation of degree $1$
	with respect to the cup product $\cup$.
\end{propositionlist}
\end{proposition}

\begin{proof}
	The first and the last part follow directly from the fact that these properties
	hold on $\TOTAL$- and $\WOBS$-component separately by the classical theory.
	The second part follows, since in this case 
	$\mu$ is well-defined on $\strtensor$ by definition of a strong constraint algebra.
\end{proof}

Now let us turn to the cohomology of the constraint Hochschild complex.

\begin{definition}[Constraint Hochschild cohomology]
	\label{definition:CoisoHCohom}%
	\index{constraint!Hochschild cohomology}
	\glsadd{ConstraintHochschildCohomology}
Let $(\algebra{A},\mu) \in \ConAlg_\field{k}$
be a constraint algebra.
The cohomology
$\Hochschild^\bullet(\algebra{A})
= \ker \delta / \image \delta$
of the Hochschild complex $\mathrm{C}^\bullet(A)$ is called the
\emph{constraint Hochschild cohomology} of $\algebra{A}$.
\end{definition}

Using the definition of kernel, image and quotient in
$\ConMod_\field{k}$, as given in
\autoref{sec:ConKModules}, we can express the constraint Hochschild 
cohomology more explicitly as follows.

\begin{lemma}
	\label{lem:ConHochschildCohomology}%
The constraint Hochschild cohomology of
$\algebra{A} \in \ConAlg_\field{k}$ is given by
\begin{equation}
\begin{split}
	\Hochschild^\bullet(\algebra{A})_\Total
	&= \Hochschild^\bullet(\algebra{A}_\Total),\\
	\Hochschild^\bullet(\algebra{A})_\Wobs
	&= \ker \delta_\Wobs / \image \delta_\Wobs,\\
	\Hochschild^\bullet(\algebra{A})_\Null
	&= \ker (\delta_\Wobs\at{\Null}) / \image \delta_\Wobs,
\end{split}
\end{equation}
with
\begin{align}
	\ker \delta_\Wobs^{n+1}
	&= \big\{ (f_\Total, f_\Wobs) \in \mathrm{C}^{n+1}(\algebra{A})_\Wobs
	\bigm|\delta^{\algebra{A}_\Total} f_\Total = 0
	\textrm{ and }
	\delta^{\algebra{A}_\Wobs} f_\Wobs = 0 \big\}\\
	&\subseteq \ker \delta_{\algebra{A}_\Total}^{n+1}
	\times \ker \delta_{\algebra{A}_\Wobs}^{n+1},\nonumber\\
	\image \delta_\Wobs^n
	&= \big\{ (f_\Total, f_\Wobs) \in \mathrm{C}^{n+1}(\algebra{A})_\Wobs
	\bigm| \exists (g_\Total, g_\Wobs) \in 
	\mathrm{C}^n(\algebra{A})_\Wobs: 
	\delta^{\algebra{A}_\Total}g_\Total = f_\Total \\
	&\hphantom{= \big\{ (f_\Total, f_\Wobs) \in \mathrm{C}^{n+1}(\algebra{A})_\Wobs \bigm|}
	\textrm{ and }
	\delta^{\algebra{A}_\Wobs}g_\Wobs = f_\Wobs	\big\},\nonumber\\
	\shortintertext{and}
	\ker (\delta^n_\Wobs\at{\Null} )
	&= \big\{ (f_\Total, f_\Wobs) \in \mathrm{C}^{n+1}(\algebra{A})_\Null
	\bigm| \delta^{\algebra{A}_\Total} f_\Total = 0
	\textrm{ and }
	\delta^{\algebra{A}_\Wobs} f_\Wobs = 0 \big\}\\
	&\subseteq \ker \delta_{\algebra{A}_\Total}^n 
	\times \ker \delta_{\algebra{A}_\Wobs}^n.\nonumber
\end{align}
\end{lemma}

With this we can compute the zeroth and first constraint Hochschild 
cohomology of a given constraint algebra.
For this recall the characterization of centre of a constraint algebra in \autoref{prop:ConCenter}
and of constraint derivations from
\autoref{prop:ConDer}, and define the constraint inner derivations of a 
given constraint algebra $\algebra{A}$ by
\index{constraint!inner derivations}
\glsadd{ConInnerDer}
\begin{equation}
\begin{split}
	\ConInnDer(\algebra{A})_\Total
	&\coloneqq \InnDer(\algebra{A}_\Total), \\
	\ConInnDer(\algebra{A})_\Wobs
	&\coloneqq \big\{ (D_\Total,D_\Wobs) \in \ConDer(\algebra{A})_\Wobs
	\bigm| \exists a \in \algebra{A}_\Wobs : D_\Wobs = [\argument, a]_\Wobs\\
	&\qquad\qquad\qquad\qquad\qquad\qquad\quad\textrm{ and }
	D_\Total = [\argument, \iota_\algebra{A}(a)]_\Total \big\}, \\
	\ConInnDer(\algebra{A})_\Null
	&\coloneqq \big\{ (D_\Total,D_\Wobs) \in \ConDer(\algebra{A})_\Null
	\bigm| \exists a \in \algebra{A}_\Null : D_\Wobs = [\argument, a]_\Wobs\\
	&\qquad\qquad\qquad\qquad\qquad\qquad\quad\textrm{ and }
	D_\Total = [\argument, \iota_\algebra{A}(a)]_\Total \big\}.
\end{split}
\end{equation}
The following also shows that in low degrees the interpretation of the 
constraint Hochschild cohomology is analogous to that for usual algebras.

\begin{proposition}
	\label{prop:ZerothFirstHochschild}
Let $\algebra{A} \in \ConAlg_\field{k}$
be a constraint algebra.
\begin{propositionlist}
	\item \label{prop:ZerothFirstHochschild_1}
	We have
	\begin{equation}
	\begin{split}
		\Hochschild^0(\algebra{A})_\Total
		&= \Center(\algebra{A}_\Total), \\
		\Hochschild^0(\algebra{A})_\Wobs
		&= \big\{ a \in \algebra{A}_\Wobs 
		\bigm|
		a \in \Center(\algebra{A}_\Wobs)
		\textrm{ and }
		\iota_\algebra{A}(a) \in \Center(\algebra{A}_\Total)\big\},	\\
		\Hochschild^0(\algebra{A})_\Null
		&= \big\{a_0 \in \algebra{A}_\Null
		\bigm|
		a_0 \in \Center(\algebra{A}_\Wobs)
		\textrm{ and }
		\iota_\algebra{A}(a_0) \in \Center(\algebra{A}_\Total) \big\}.
	\end{split}
	\end{equation}
	Hence $\Hochschild^0(\algebra{A}) = \Center(\algebra{A})$.
	\item \label{item:HHoneIsDer}
	We have
	\begin{equation}
	\begin{split}
		\Hochschild^1(\algebra{A})_\Total
		&= \Der(\algebra{A}_\Total) / \InnDer(\algebra{A}_\Total), \\
		\Hochschild^1(\algebra{A})_\Wobs
		&= \Der(\algebra{A})_\Wobs / \big\{	(D_\Total, D_\Wobs)
		\in \Der(\algebra{A})_\Wobs
		\bigm|\\
		&\qquad\qquad\qquad\exists a \in \algebra{A}_\Wobs:
		D_\Total = [\argument, \iota_\algebra{A}(a)],
		D_\Wobs = [\argument, a] \big\}, \\
		\Hochschild^1(\algebra{A})_\Null
		&= \Der(\algebra{A})_\Null / \big\{	(D_\Total, D_\Wobs)
		\in \Der(\algebra{A})_\Wobs
		\bigm|\\
		&\qquad\qquad\qquad\exists a \in \algebra{A}_\Wobs:
		D_\Total = [\argument, \iota_\algebra{A}(a)],
		D_\Wobs = [\argument, a]\big\}.
	\end{split}
	\end{equation}
	Hence
	$\Hochschild^1(\algebra{A}) 
	= \ConDer(\algebra{A}) / \ConInnDer(\algebra{A})$.
\end{propositionlist}
\end{proposition}

\begin{proof}
The first claim is clear by
\autoref{lem:ConHochschildCohomology} and $\delta_{-1} = 0$.
The $\TOTAL$-component of the second part is clear by the 
classical result for the first Hochschild cohomology of the classical 
algebra $\algebra{A}_\Total$.
For the $\WOBS$-component consider
$D = (D_\Total,D_\Wobs) \in \ker \delta_\Wobs^1$.
Then $\delta^{\algebra{A}_\Total}D_\Total = 0$
and 
$\delta^{\algebra{A}_\Wobs}D_\Wobs = 0$,
hence $D_\Total$ and $D_\Wobs$ are derivations and it follows 
$D \in \Der(\algebra{A})_\Wobs$.
Similarly, we get
$D \in \Der(\algebra{A})_\Null$ for
$D \in \ker (\delta_\Wobs^1\at{\Null})$.
Now let
$D \in \image \delta_\Wobs^0$,
then there exists $a \colon \field{k} \to \algebra{A}$
with
$D_\Total 
= \delta^{\algebra{A}_\Total}a_\Total 
= [\argument, a_\Total]$
and
$D_\Wobs = \delta^{\algebra{A}_\Wobs}a_\Wobs 
= [\argument, a_\Wobs]$.
Since $a_\Total = \iota(a_\Wobs)$ the second part holds.
\end{proof}
%



\subsubsection{Reduction}
Assigning the (constraint) Hochschild complex to a given (constraint)
algebra is not functorial on all of $\ConAlg_\field{k}$.
But if we restrict ourselves to the subcategory
$\ConAlg_\field{k}^\times$
of constraint algebras with invertible morphisms we get a functor
$\mathrm{C}^\bullet \colon \ConAlg_\field{k}^\times \to \ConDGLieAlg$
by mapping each constraint algebra to its constraint Hochschild 
complex and 
every algebra isomorphism
$\phi \colon \algebra{A} \to \algebra{B}$ to
$\mathrm{C}^\bullet(\phi) \colon \mathrm{C}^\bullet(\algebra{A}) 
\to \mathrm{C}^\bullet(\algebra{B})$
given by
$\mathrm{C}^\bullet(\phi)(f) = \phi \circ f \circ 
(\phi^{-1})^{\tensor n}$
for $f \in \mathrm{C}^n(\algebra{A})_{\Total/\Wobs}$.
A similar construction clearly also works for usual algebras.
We can now show that this functor commutes with reduction up to an 
injective natural transformation.

\begin{proposition}[Hochschild complex vs. reduction]
	\label{prop:HochschildVSReduction}%
	\index{reduction!Hochschild cohomology}
There exists an injective natural\linebreak
transformation
$\eta \colon \red \circ \mathrm{C}^\bullet \Longrightarrow
\mathrm{C}^\bullet \circ \red$, i.e.
\begin{equation}
	\begin{tikzcd}
		\ConAlg_\field{k}^\times
		\arrow{r}{\mathrm{C}^\bullet}
		\arrow{d}[swap]{\red}
		& \ConDGLieAlg
		\arrow{d}{\red}
		\arrow[Rightarrow]{dl}[swap]{\eta}\\
		\Algebras_\field{k}^\times
		\arrow{r}{\mathrm{C}^\bullet}
		&\DGLieAlg
	\end{tikzcd}
\end{equation}
commutes with $\eta$ injective.
\end{proposition}

\begin{proof}
	For every constraint algebra $\algebra{A}$ define
	$\eta_\algebra{A} \colon \mathrm{C}^\bullet(\algebra{A})_\red 
	\to \mathrm{C}^\bullet(\algebra{A}_\red)$ by
	\begin{equation*}
		\eta_\algebra{A}([f])([a_1], \dotsc, [a_n])
		= \left[ f_\Wobs(a_1, \dotsc, a_n) \right].
	\end{equation*}
	for $[f]=[(f_\Total,f_\Wobs)] \in \mathrm{C}^n(\algebra{A})_\red$.
	First note that
	$\eta_\algebra{A}([f]) \colon \algebra{A}_\red^{\tensor n}
	\to \algebra{A}_\red$ is well-defined since if
	$a_i \in \algebra{A}_0$
	for any $i = 1,\dotsc,n$ we have
	$f_\Wobs(a_1, \dotsc, a_n) \in \algebra{A}_\Null$
	and hence
	$[f_\Wobs(a_1, \dotsc, a_n)] = 0$.
	Moreover, 
	$\eta_\algebra{A}$
	is well-defined since for
	$f \in \mathrm{C}^n(\algebra{A})_\Null$
	we have
	$f_\Wobs(a_1, \dotsc, a_n) \in \algebra{A}_\Null$
	and thus $\eta([f]) = 0$.
	To see that $\eta$ is indeed a natural transformation we need to 
	show that 
	for every isomorphism
	$\phi \colon \algebra{A} \to \algebra{B}$
	we have
	$\eta_\algebra{B} \circ \mathrm{C}^\bullet(\phi)_\red 
	= \mathrm{C}^\bullet([\phi]) \circ \eta_\algebra{A}$.
	But it is clear after inserting the definitions.
	Finally, suppose
	$\eta_\algebra{A}([f]) = \eta_\algebra{A}([g])$.
	This means that
	$(f_\Wobs-g_\Wobs)(a_1, \dotsc, a_n) \in \algebra{A}_\Null$ 
	and therefore $[f] = [g]$.
	Thus $\eta_\algebra{A}$ is injective.
\end{proof}

Combining \autoref{prop:CohomologyCommutesWithReduction} with
\autoref{prop:HochschildVSReduction} immediately
yields the following compatibility of Hochschild cohomology with 
reduction:

\begin{corollary}
	\label{cor:HHandReduction}%
There exists an injective natural transformation
$\eta \colon \red \circ \Hochschild^\bullet \Longrightarrow
\Hochschild^\bullet \circ \red$.
In particular, for any constraint algebra $\algebra{A}$ we have
\begin{equation} \label{eq:HHReduction}
	\Hochschild^\bullet(\algebra{A})_\red
	\subseteq \Hochschild^\bullet(\algebra{A}_\red).
\end{equation}
\end{corollary}

\section{Formal Deformations via Hochschild Cohomology}
	\label{sec:FormalDeformations}

Throughout this section we will again assume that the scalars satisfy
$\field{Q} \subseteq \field{k}$ in order to make use of the
description of deformations by Maurer-Cartan elements.

\index{constraint!deformation}
Let $(\algebra{A},\mu_0) \in \ConAlg_\field{k}$
be a constraint $\field{k}$-algebra.
By \autoref{def:DeformationConAlg} a formal associative
deformation
$(\algebra{A}\formal{\lambda},\mu)$
is given by an associative multiplication
$\mu \colon \algebra{A}\formal{\lambda} \tensor
\algebra{A}\formal{\lambda} \to
\algebra{A}\formal{\lambda}$
making $\algebra{A}\formal{\lambda}$ 
a constraint $\field{k}\formal{\lambda}$-algebra such that
$\cl(\algebra{A},\mu)$ is given by $(\algebra{A},\mu_0)$,
or in other words
\begin{equation}
	\label{eq:MuIsFormalSeries}
	\mu = \mu_0 + \sum_{k=1}^{\infty} \lambda^k \mu_k
\end{equation}
with
$\mu_k \colon \algebra{A} \tensor \algebra{A} \to \algebra{A}$.

Such deformations can now be understood as Maurer-Cartan elements 
in the constraint DGLA
$\lambda \mathrm{C}^\bullet(\algebra{A})\formal{\lambda}$
corresponding to $(\algebra{A}\formal{\lambda},\mu_0)$.

\begin{lemma}
	\label{lem:AssoDefIsMC}%
Let $(\algebra{A},\mu) \in \ConAlg_\field{k}$
be a constraint $\field{k}$-algebra.
A multiplication
$\mu = \mu_0 + M$,
with
$M = \sum_{k=1}^{\infty} \lambda^k \mu_k$
is a formal associative deformation of $\mu_0$ if and only if
\begin{equation}
	\label{eq:MisMC}
	\delta M + \frac{1}{2}[M, M] = 0.
\end{equation}
\end{lemma}

\begin{proof}
By \autoref{lem:AssoConMultiplication} we know that 
we have to check that
$[\mu_\Total, \mu_\Total]_{\algebra{A}_\Total} = 0$
and $[\mu_\Wobs, \mu_\Wobs]_{\algebra{A}_\Wobs} = 0$. 
Thus, consider the total component of $\mu$ as
$\mu_\Total = (\mu_0)_\Total + M_\Total$.
We have
\begin{equation*}
	[\mu_\Total, \mu_\Total]
	= [(\mu_0)_\Total + M_\Total, (\mu_0)_\Total + M_\Total]
	= 2 \delta M_\Total + [M_\Total, M_\Total],
\end{equation*}
where we used the associativity of $(\mu_0)_\Total$ and the 
graded skew-symmetry of Gerstenhaber bracket.
The very same holds for the $\WOBS$-component.
\end{proof}

Equivalence of formal deformations can be phrased using the constraint Gerstenhaber bracket as follows:

\begin{proposition}
	\label{prop:EquivalenceOfDeformations}
	Let $(\algebra{A}, \mu_\Null)$ be a constraint algebra and let
	$\mu$ and $\mu^\prime$ be deformations of $\algebra{A}$.
	Then $\mu$ and $\mu^\prime$ are equivalent via
	$T = \exp(\lambda D)$ if and only if
	\begin{equation}
		\E^{\lambda[D,\argument]}(\mu) = \mu^\prime
	\end{equation}
	holds, where $[\argument,\argument]$ denotes the constraint Gerstenhaber bracket.
\end{proposition}

\begin{proof}
	This follows directly since the statement holds in the $\TOTAL$- and $\WOBS$-components separately by classical deformation theory.
\end{proof}

This allows us to conclude that the equivalence of deformations coincides with
gauge equivalence of Maurer-Cartan elements:

%

\begin{theorem}[Equivalence classes of deformations]
	\label{theorem:EquivalenceClassesOfDefs}%
	\index{equivalence!deformations}
Let $\field{k}$ be a commutative ring with\linebreak
$\field{Q} \subseteq \field{k}$.
Let $(\algebra{A},\mu_0)$ be a constraint $\field{k}$-algebra.
Then the constraint set of equivalence classes of formal associative 
deformations of $\algebra{A}$ coincides with 
$\Def(\mathrm{C}^\bullet(\algebra{A}))$,
where
$\mathrm{C}^\bullet(\algebra{A})$
is the constraint Hochschild DGLA of $\algebra{A}$.
\end{theorem}

\begin{proof}
Let $\mu$ and $\mu^\prime$
be two deformations of $\mu_0$.
By \autoref{prop:EquivalenceOfDeformations} we know that
$\mu$ and $\mu^\prime$ are equivalent deformations of $\mu_0$
if and only if there exists $D \in \lambda C^1(\algebra{A}\formal{\lambda})$
such that
\begin{equation*} \label{eq:EquivalenceOfDeformations}
	\E^{\lambda[D,\argument]}(\mu) = \mu'. \tag{$*$}
\end{equation*}
Using  $\mu = \mu_0 + M$ and $\mu^\prime = \mu_0 + M^\prime$
with
$M = \sum_{k=1}^{\infty} \lambda^k \mu_k$
as well as $\delta_\Wobs = [\mu_0, \argument]_\Wobs$
for the Hochschild differential it is easy to see that 
\eqref{eq:EquivalenceOfDeformations}
is equivalent to
\begin{equation*}
	\lambda D \acts_\Wobs M = M^\prime,
\end{equation*}
meaning that the Maurer-Cartan elements
$M$ and $M^\prime$ are gauge equivalent.
\end{proof}

Finally, we can reformulate the classical theorem about the 
extension of a deformation up to a given order for constraint 
algebras.

\begin{theorem}[Obstructions]
	\label{thm:Obstructions}%
Let $\field{k}$ be a commutative ring with
$\field{Q} \subseteq \field{k}$.
Let
$(\algebra{A},\mu_0) \in \ConAlg_\field{k}$
be a constraint $\field{k}$-algebra.
\begin{theoremlist}
	\item Furthermore, let
	$\mu^{(k)} = \mu_0 + \cdots + \lambda^k \mu_k
	\in \mathrm{C}^2(\algebra{A})_\Wobs$
	be an associative deformation of $\mu_0$ up to order $k$.
	Then
	\begin{equation} \label{eq:thm:ExtendingDeformation}
		R_{k+1}
		= \Big( \frac{1}{2} \sum_{\ell=1}^k 
		\big[(\mu_\ell)_\Total,(\mu_{k+1-\ell})_\Total	
		\big]^{\algebra{A}_\Total},\;
		\frac{1}{2} \sum_{\ell=1}^{k}
		\big[(\mu_\ell)_\Wobs,(\mu_{k+1-\ell})_\Wobs		
		\big]^{\algebra{A}_\Wobs}\Big)
		\in \mathrm{C}^3(\algebra{A})_\Wobs
	\end{equation}
	is a constraint Hochschild cocycle, i.e.
	$\delta_\Wobs R_{k+1} = 0$. 
	The deformation $\mu^{(k)}$ can be extended to order $k+1$ if and 
	only if
	$R_{k+1}= \delta_\Wobs \mu_{k+1}$.
	In this case every such	$\mu_{k+1}$ yields an extension
	$\mu^{(k+1)} = \mu^{(k)} + \lambda^{k+1}\mu_{k+1}$.
	\item Let
	$\mu_1 \in \mathrm{C}^2(\algebra{A})_\Wobs$.
	Then
	$\mu = \mu_0 + \lambda \mu_1$
	is an associative deformation of $\mu_0$ up to order $1$ if and 
	only if
	$\delta_\Wobs \mu_1 = 0$.
	Moreover, if $\mu'_1$ is another deformation up to order $1$
	of $\mu_0$ then these two deformations are equivalent up to
	order $1$ if and only if $\mu_1 - \mu'_1$ is exact.
\end{theoremlist}
\end{theorem}

\begin{proof}
By classical deformation theory of associative algebras 
it is clear that
\eqref{eq:thm:ExtendingDeformation}
is closed since
$\delta_\Wobs =
(\delta^{\algebra{A}_\Total},\delta^{\algebra{A}_\Wobs})$.
If $R_{k+1}$ is exact, we know that $\mu^{(k)}_\Total$ and
$\mu^{(k)}_\Wobs$ can be extended via $(\mu_{k+1})_\Total$ and
$(\mu_{k+1})_\Wobs$, respectively.
Thus $\mu_{k+1}$ yields an extension of $\mu^{(k)}$.
On the other hand, if $\mu^{(k)}$ can be extended, we know that
$(R_{k+1})_\Total = \delta^{\algebra{A}_\Total}(\mu_{k+1})_\Total$
and
$(R_{k+1})_\Wobs = \delta^{\algebra{A}_\Wobs}(\mu_{k+1})_\Wobs$.
Hence, $R_{k+1} = \delta_\Wobs \mu_{k+1}$.
For the second part, consider the first part for $k=0$, then
$\delta_\Wobs \mu_1 = R_1 = 0$
follows directly.
By \autoref{prop:EquivalenceOfDeformations} two deformations
$\mu = \mu_0 +\mu_1$ and $\mu' = \mu_0 + \mu'_1$
are equivalent if and only if there exists
$D \in 
\ConHom_\field{k}(\algebra{A}\formal{\lambda},\algebra{A}\formal{\lambda})_\Wobs$
such that $\E^{\ad(D)}(\mu) = \mu'$.
If we only want to consider deformations up to order $1$ we 
can restrict to the case
$D = D_0 \in \ConHom_\field{k}(\algebra{A},\algebra{A})_\Wobs$.
Then we get equivalently $\mu + \lambda[D_0,\mu] = \mu'$.
The first order term then directly yields
$\mu'_1 - \mu_1 = -\delta_\Wobs D_0$.
\end{proof}

Thus $\Hochschild^2(\algebra{A})_\Wobs$ classifies infinitesimal
constraint deformations while
$\Hochschild^3(\algebra{A})_\Wobs$
gives the obstructions to extending
such deformations in a constraint way.
The constraint module $\Hochschild^3(\algebra{A})$ carries more
information than just the obstructions to deformations of the
constraint algebra $\algebra{A}$.
Since
$\Hochschild^3(\algebra{A})_\Total =
\Hochschild^3(\algebra{A}_\Total)$
it also encodes the obstructions of deformations of the classical 
algebra $\algebra{A}_\Total$.
Moreover,
$\Hochschild^3(\algebra{A})_\Null$
is important for the reduction of
$\Hochschild^3(\algebra{A})$
and hence controls which obstructions on
$\algebra{A}$ descend to obstructions on $\algebra{A}_\red$.
In particular, we have seen in
\autoref{cor:HHandReduction} that
$\Hochschild^3(\algebra{A})_\red \subseteq
\Hochschild^3(\algebra{A}_\red)$.
The components of $\Hochschild^2(\algebra{A})$ can be interpreted 
in a similar fashion.

\section{Second Constraint Hochschild Cohomology on $\Reals^n$}
\label{sec:SecondConHochschildCohomology}

Let us now turn again to constraint star products on a constraint manifold $\mathcal{M}$.
We have seen in \autoref{sec:ConStarProducts} that such a constraint star product is nothing but a differentiable 
formal deformation of the constraint algebra $\ConCinfty(\mathcal{M})$.
Thus following \autoref{sec:FormalDeformations} we are interested in the subcomplex
$\mathrm{C}^\bullet_\diff(\ConCinfty(\mathcal{M})) \subseteq \mathrm{C}^\bullet(\ConCinfty(\mathcal{M}))$
of differential constraint Hochschild cochains.
Thus we want to be able to compute
$\Hochschild_\diff^\bullet(\ConCinfty(\mathcal{M}))$.
In classical deformation theory the Hochschild-Kostant-Rosenberg Theorem
computes the differential Hochschild cohomology for a given smooth manifold $M$,
see \cite{hochschild.kostant.rosenberg:1962a} for the original result.

\begin{theorem}[HKR Theorem]
	\label{thm:ClassicalHKR}
	\index{Hochschild-Kostant-Rosenberg Theorem}
	\glsadd{HKRMap}
Let $M$ be a smooth manifold. Then
\begin{equation} \label{eq:ClassicalHKRMap}
	\mathcal{U} \colon \VecFields^\bullet(M) \to \Hochschild^\bullet_\diff(\Cinfty(M)),
	\qquad
	\mathcal{U}(X)(f_1,\dotsc,f_k) \coloneqq \frac{1}{k!} \ins_{\D f_1,\dotsc, \D f_k} X
\end{equation}
is an isomorphism of Gerstenhaber algebras.
\end{theorem}

Here the Gerstenhaber algebra structure on $\VecFields^\bullet(M)$ 
is given by $\wedge$ and the Schouten bracket $\Schouten{\argument,\argument}$,
while on $\Hochschild_\diff^\bullet(\Cinfty(M))$
it is given by the cup product $\cup$ and the Gerstenhaber bracket $[\argument,\argument]$.

For a constraint manifold $\mathcal{M} = (M,C,D)$ we know from
\autoref{prop:ZerothFirstHochschild}
that
\begin{align}
	\Hochschild_\diff^0\big(\ConCinfty(\mathcal{M})\big) 
	&= \ConCinfty(\mathcal{M})\\
	\shortintertext{and}
	\Hochschild_\diff^1\big(\ConCinfty(\mathcal{M})\big) 
	&= \ConDer\big(\ConCinfty(\mathcal{M})\big)
	\simeq \ConSecinfty(T \mathcal{M}).
\end{align}
This suggests that constraint multivector fields might compute constraint Hochschild cohomology.
But for higher constraint multivector fields we have to choose between
$\ConVecFields_{\tensor}^\bullet(\mathcal{M})$
and
$\ConVecFields_{\strtensor}^\bullet(\mathcal{M})$.
Now \autoref{ex:ConPoissonStructureAsMCElements}
shows, that if we are interested in deforming not merely Poisson but coisotropic
submanifolds, we need to go with 
$\ConVecFields_{\strtensor}^\bullet(\mathcal{M})$.
The next result shows that the constraint version of
\eqref{eq:ClassicalHKRMap}
is well-defined at the level of cochains and also
yields an injection in cohomology.

\begin{proposition}
	\label{prop:ConHKRmap}
	\index{constraint!multivector fields}
Let $\mathcal{M} = (M,C,D)$ be a constraint manifold.
\begin{propositionlist}
	\item The map
	\begin{equation} \label{eq:ConHKRmap}
		\mathcal{U} \colon \ConVecFields_{\strtensor}^\bullet(\mathcal{M})
		\to \mathrm{C}^\bullet_\diff(\Cinfty(\mathcal{M})),
		\qquad
		\mathcal{U}(X)(f_1,\dotsc,f_k) \coloneqq \frac{1}{k!} \ins_{\D f_1,\dotsc, \D f_k} X
	\end{equation}
	is a morphism between the constraint complexes
	$\big(\VecFields_{\strtensor}^\bullet(M),\D=0\big)$
	and
	$\big(\mathrm{C}^\bullet_\diff(\ConCinfty(\mathcal{M})),\delta\big)$.
	\item The induced morphism 
	\begin{equation} \label{eq:ConHKRmap_Cohom}
		\mathcal{U} \colon \ConVecFields_{\strtensor}^\bullet(\mathcal{M}) \to \Hochschild_\diff^\bullet\big(\ConCinfty(\mathcal{M})\big)
	\end{equation}
	is a regular monomorphism.
\end{propositionlist}
\end{proposition}

\begin{proof}
The map $\mathcal{U}$ can be seen as the lowest order
of $\Op$ from \autoref{cor:ConMultDiffOpSymbolCalculusFunctions}.
Note that in this case $\SymCoDer f = \D f$, and hence
this restriction of
$\Op$ is indeed independent of the chosen constraint
covariant derivative.
Thus $\mathcal{U}$ is a constraint regular monomorphism.
Moreover, from classical theory we know that
$\delta \circ \mathcal{U} = 0$, and hence
$\mathcal{U}$ is a morphism of constraint
complexes.
For the second part note that
$\mathcal{U}_\Total \colon \VecFields^\bullet(M) \to \mathrm{C}^\bullet_\diff(\ConCinfty(M))$
is an isomorphism by the classical HKR theorem.
Since quotients of embedded constraint modules need not necessarily be embedded,
$\mathcal{U}_\Wobs \colon \ConVecFields_{\strtensor}^\bullet(\mathcal{M})_\Wobs
\to \Hochschild^\bullet_\diff(\ConCinfty(M))_\Wobs$
is not given by the restriction of $\mathcal{U}_\Total$, but
it fulfils $\iota_{\Hochschild} \circ \mathcal{U}_\Wobs = \mathcal{U}_\Total \circ \iota_{\ConVecFields_{\strtensor}}$.
Since the right hand side is injective so is $\mathcal{U}_\Wobs$.
Thus $\mathcal{U}$ is a monomorphism.
To show that it is regular, consider
$\mathcal{U}_\Wobs(X) = [\mathcal{U}(X)] \in \Hochschild^\bullet_\diff(\ConCinfty(\mathcal{M}))_\Null$.
Then by definition
$\mathcal{U}(X) \in \mathrm{C}^\bullet_\diff(\ConCinfty(\mathcal{M}))_\Null$
and thus since $\mathcal{U}$ is a constraint monomorphism
on cochain level we get $X \in \ConVecFields_{\strtensor}^\bullet(\mathcal{M})_\Null$.
\end{proof}

Even though $\ConVecFields_{\strtensor}^\bullet(\mathcal{M})$
already yields interesting, and perhaps even unexpected, contributions
to $\Hochschild_\diff^\bullet(\ConCinfty(\mathcal{M}))_\Wobs$,
we cannot hope for \eqref{eq:ConHKRmap_Cohom} to be an isomorphism:

\begin{example}
Let
$\mathcal{M} = \Reals^n = (\Reals^{n_\Total}, \Reals^{n_\Wobs},\Reals^{n_\Null})$
with $0<n_\Null<n_\Wobs<n_\Total$.
and consider
\begin{equation}
	\del_{(1,n_\Total)} = \frac{\del^2}{\del x^1 \del x^{n_\Total}}
	\in \mathrm{C}^1_\diff(\ConCinfty(\Reals^n))_\Total.
\end{equation}
This differential operator is clearly not constraint:
We have $x^1x^{n_\Total} \in \ConCinfty(\Reals^n)_\Null$
but
\begin{equation}
	\del_{(1,n_\Total)}(x^1 x^{n_\Total}) = 1 \notin \ConCinfty(\Reals^n)_\Null.
\end{equation}
Nevertheless, applying the Hochschild differential yields
\begin{equation}
	\delta(\del_{(1,n_\Total)}) = - \del_1 \cup \del_{n_\Total} - \del_{n_\Total} \cup \del_1,
\end{equation}
which \emph{is} constraint. 
In fact, $\delta(\del_{(1,n_\Total)}) \in \mathrm{C}^2_\diff(\ConCinfty(\Reals^n))_\Null$,
since for $f,g \in \ConCinfty(\Reals^n)_\Wobs$ we have
\begin{equation}
	\delta(\del_{(1,n_\Total)})(f,g)
	= - \underbrace{\frac{\del f}{\del x^1}}_{=0} \cdot \frac{\del g}{\del x^{n_\Total}}
	- \frac{\del f}{\del x^{n_\Total}}\cdot\underbrace{\frac{\del g}{\del x^{1}}}_{=0}
	= 0.
\end{equation}
To show that $\delta(\del_{(1,n_\Total)})$ defines a non-trivial class in cohomology assume that
$\delta(\del_{(1,n_\Total)}) = \delta(D)$ for some
$D = \sum_{r=0}^{\infty}\sum_{I \in n^{\tensor r}}D^I \del_I \in \mathrm{C}^1_\diff(\ConCinfty(\Reals^n))_\Total$.
Then
$D - \delta(\del_{(1,n_\Total)})$ is closed, hence a derivation,
and it follows
$D = \del_{(1,n_\Total)} + \sum_{i=1}^{n_\Total} D^i \del_i$.
Evaluating on $x^1x^{n_\Total}$ shows that
$D$ is not constraint.
Finally, since $\delta(\del_{(1,n_\Total)})$ is symmetric it cannot be in the image
of $\mathcal{U}$.
Thus we have found a non-trivial cohomology class, not coming
from constraint multivector fields.
\end{example}

This example can easily be generalized to construct non-vanishing
symmetric cohomology classes with arbitrary order of differentiation:
For this consider $\del_I \in \mathrm{C}_\diff^1(\ConCinfty(\Reals^n))_\Total$
with $I = (i_1, \dotsc, i_r)$ such that for one $\ell \in \{1,\dotsc, r\}$
it holds $n_\Wobs < i_\ell$ and $i_k \leq n_0$ for all $k \neq \ell$.
Then $\del_I$ is not constraint, but in $\delta(\del_I)$
there appears in every term at least one $\cup$-factor from
$\mathrm{C}^1_\diff(\ConCinfty(\Reals^n))_\Null$, showing
that $\delta(\del_I)$ is constraint.
It is then straightforward to see that it also yields a non-vanishing
class in cohomology.

In the following we will concentrate on the local case 
$\mathcal{M} = \Reals^n = (\Reals^{n_\Total}, \Reals^{n_\Wobs},\Reals^{n_\Null})$
and its second constraint Hochschild cohomology.
By the product rule from classical calculus we have for 
$i_1, \dotsc, i_r \in \{1, \dotsc n_\Total\}$ 
and $f,g \in \Cinfty(M)$ that
\begin{equation} \label{eq:PartialProductRule}
\begin{split}
	\frac{\del^r(f\cdot g)}{\del x^{i_1} \dots \del x^{i_r}}
	&= \sum_{s = 0}^{r} \sum_{\sigma \in S_r} \frac{1}{s!(r-s!)}
	\frac{\del^s f}{\del x^{i_{\sigma(1)}} \dots \del x^{i_{\sigma(s)}}}
	\cdot \frac{\del^{(r-s)} g}{\del x^{i_{\sigma(s+1)}} \dots \del x^{i_{\sigma(r)}}} \\
	&= \sum_{s=0}^{r} \sum_{\sigma \in \Shuffle(s,r-s)}
	\frac{\del^s f}{\del x^{i_{\sigma(1)}} \dots \del x^{i_{\sigma(s)}}}
	\cdot \frac{\del^{(r-s)} g}{\del x^{i_{\sigma(s+1)}} \dots \del x^{i_{\sigma(r)}}}.
\end{split}	
\end{equation}
Here $\Shuffle(s,r-s)$ denotes the set of $(s,r-s)$-shuffle permutations, i.e.
$\sigma \in S_{r}$ such that $\sigma(1) < \dots < \sigma(s)$ and
$\sigma(s+1) < \dots < \sigma(r)$.
In order to write \eqref{eq:PartialProductRule} in a more concise fashion we use the following notation:
For a multi index $I = (i_1, \dotsc, i_r)$ and $s \in \{1,\dotsc, r\}$
we define
\begin{equation}
	\deco{}{}{I}{}{s} \coloneqq (i_1, \dotsc, i_s)\\
	\qquad\text{and}\qquad
	\deco{}{s}{I}{}{} \coloneqq (i_{s+1},\dotsc, i_r).
\end{equation}
Moreover, for a permutation $\sigma \in S_r$ we set
$\sigma(I) \coloneqq (i_{\sigma(1)}, \dotsc, i_{\sigma(r)})$.
With this \eqref{eq:PartialProductRule} reads
\begin{equation}
	\del_I (f \cdot g) = \sum_{s=0}^{r} \sum_{\sigma \in \Shuffle(s,r-s)}
	\del_{\deco{}{}{\sigma(I)}{}{s}} f \cdot \del_{\deco{}{{s}}{\sigma(I)}{}{}} g.
\end{equation}
We can now use \eqref{eq:HochschildDifferentialFormula} to express the Hochschild differential applied to some
$\del_I$ as
\begin{equation}
	\delta(\del_I) = \sum_{s=1}^{r-1} \sum_{\sigma \in \Shuffle(s,r-s)}
	\del_{\deco{}{}{\sigma(I)}{}{s}} \cup \del_{\deco{}{{s}}{\sigma(I)}{}{}}.
\end{equation}

\begin{lemma}
Let $I = (i_1,\dotsc, i_r) \in (n^{\tensor r})_\Wobs$.
\begin{lemmalist}
	\item It holds $\delta(\del_I) \in \mathrm{C}^2_\diff(\ConCinfty(\Reals^n))_\Wobs$
	if and only if
	$I \in (n^{\tensor r})_\Wobs$ or
	\begin{equation}
	\begin{split}
		\exists \ell \in \{1,\dotsc, r\} : i_\ell \in n_\Total\setminus n_\Wobs
		\qquad\text{and}\qquad\forall k \neq \ell : i_k \in n_\Null.
	\end{split}
	\end{equation}
	\item It holds $\delta(\del_I) \in \mathrm{C}^2_\diff(\ConCinfty(\Reals^n))_\Null$
	if and only if
	$I \in (n^{\tensor r})_\Null$ or
	\begin{equation}
		\begin{split}
			\exists \ell \in \{1,\dotsc, r\} : i_\ell \in n_\Total\setminus n_\Wobs
			\qquad\text{and}\qquad\forall k \neq \ell : i_k \in n_\Null.
		\end{split}
	\end{equation}
\end{lemmalist}
\end{lemma}

\begin{proof}
Let us first show the second part:
By \autoref{prop:LocalConMultDiffop} \ref{prop:LocalConMultDiffop_2} the terms with
\begin{equation*}
	(\sigma(I)_s,\deco{}{s}{\sigma(I)}{}{}) 
	\in ((n^*)^{\strtensor} \tensor (n^*)^{\strtensor})_\Wobs
	= (n^{\tensor s} \strtensor n^{\tensor r-s})^*_\Wobs
\end{equation*}
need to vanish.
Thus we have 
$\delta(\del_I) \in \mathrm{C}^2_\diff(\ConCinfty(\Reals^n))_\Null$
if and only if
\begin{equation*}
	(\deco{}{}{\sigma(I)}{}{s}, \deco{}{s}{\sigma(I)}{}{})
	\in (n^{\tensor s} \strtensor n^{\tensor r-s})_\Null
	\qquad\text{for all }
	s = 1, \dotsc, r-1
	\text{ and all }
	\sigma \in \Shuffle(s,r-s).
\end{equation*}
By \autoref{lem:DecomposeStrTensorConIndSets} we can write
\begin{align*}
	(n^{\tensor s} \strtensor n^{\tensor r-s})_\Null
	&= (n^{\tensor s} \tensor n^{\tensor r-s})_\Null
	\sqcup \left((n^{\tensor s})^*_\Null \times (n^{\tensor r-s})_\Null\right)
	\sqcup \left((n^{\tensor s})_\Null \times (n^{\tensor r-s})^*_\Null\right)\\
	&= (n^{\tensor r})_\Null
	\sqcup \left((n^{\tensor s})^*_\Null \times (n^{\tensor r-s})_\Null\right)
	\sqcup \left((n^{\tensor s})_\Null \times (n^{\tensor r-s})^*_\Null\right).
\end{align*}
Now for $(\deco{}{}{\sigma(I)}{}{s}, \deco{}{s}{\sigma(I)}{}{})$
to end up in
$\left((n^{\tensor s})^*_\Null \times (n^{\tensor r-s})_\Null\right)
\sqcup \left((n^{\tensor s})_\Null \times (n^{\tensor r-s})^*_\Null\right)$
we clearly need at least one $\ell \in \{1, \dotsc, r\}$ with
$i_\ell \in n^*_0 = n_\Total \setminus n_\Wobs$.
If there is one other $k \in \{1,\dotsc r\}$
with $i_k \in n_\Total\setminus n_\Null$ then
the permutation $\tau \in \Shuffle(1,r-1)$ which moves
$i_k$ to the first or last position gives a contradiction.
This shows the second part.
For the first part it follows from \autoref{prop:LocalConMultDiffop} \ref{prop:LocalConMultDiffop_1}
that only the terms in 
\begin{align*}
	(n^{\tensor s} \strtensor n^{\tensor r-s})_\Wobs
	&= (n^{\tensor s} \tensor n^{\tensor r-s})_\Wobs
	\sqcup \left((n^{\tensor s})^*_\Null \times (n^{\tensor r-s})_\Null\right)
	\sqcup \left((n^{\tensor s})_\Null \times (n^{\tensor r-s})^*_\Null\right)\\
	&= (n^{\tensor r})_\Wobs
	\sqcup \left((n^{\tensor s})^*_\Null \times (n^{\tensor r-s})_\Null\right)
	\sqcup \left((n^{\tensor s})_\Null \times (n^{\tensor r-s})^*_\Null\right).
\end{align*}
need to vanish.
Then the same arguments as before apply.
\end{proof}

\begin{proposition}
	\label{prop:LocalConClosedOneCochains}
Let
$D = \sum_{\len(I)\leq r} D^I \del_I \in \mathrm{C}^1_\diff(\ConCinfty(\mathcal{\Reals}^n))_\Total$
be given.
\begin{propositionlist}
	\item It holds $\delta(D) \in \mathrm{C}^2_\diff(\ConCinfty(\Reals^n))_\Wobs$
	if and only if
	\begin{equation}
	\begin{split}
		D^I  \in \ConCinfty(\Reals^n)_\Wobs
		\quad\text{if}\quad &\forall \ell \in \{1,\dotsc, r\} : i_\ell \in n_\Wobs \setminus n_\Null
	\end{split}	
	\end{equation}
	and
	\begin{equation}
	\begin{split}
		D^I  \in \ConCinfty(\Reals^n)_\Null
		\quad\text{if}\quad &\exists \ell \in \{1,\dotsc, r\} : i_\ell \in n_\Total \setminus n_\Wobs\\
		\quad\text{and}\quad&\exists k \neq \ell :  i_k \in n_\Total \setminus n_\Null.
	\end{split}
	\end{equation}
	\item It holds $\delta(D) \in \mathrm{C}^2_\diff(\ConCinfty(\Reals^n))_\Null$
	if and only if
	\begin{equation}
		\begin{split}
			D^I  \in \ConCinfty(\Reals^n)_\Null
			\quad\text{if}\quad &\forall \ell \in \{1,\dotsc, r\} : i_\ell \in n_\Wobs\setminus n_\Null
		\end{split}	
	\end{equation}
	and
	\begin{equation}
		\begin{split}
			D^I  \in \ConCinfty(\Reals^n)_\Null
			\quad\text{if}\quad &\exists \ell \in \{1,\dotsc, r\} : i_\ell \in n_\Total\setminus n_\Wobs\\
			\quad\text{and}\quad&\exists k \neq \ell :  i_k \in n_\Total \setminus n_\Null.
		\end{split}
	\end{equation}
\end{propositionlist}
\end{proposition}

\begin{proof}
	We have
	\begin{equation*}
		\delta(D) = - \sum_{\len(I)\leq r} D^I \delta(\del_I)
		= - \sum_{\len(I) \leq r} \sum_{s=1}^{\len(I)-1} \sum_{\sigma \in \Shuffle(s, \len(I)-s)}
		D^I \cdot \del_{\deco{}{}{\sigma(I)}{}{s}} \cup \del_{\deco{}{s}{\sigma(I)}{}{}}.
	\end{equation*}
	Assume $\delta(D) \in \mathrm{C}^2_\diff(\ConCinfty(\Reals^n))_\Wobs$.
	By \autoref{prop:LocalConMultDiffop} this holds if and only if 
	$D^I \in \ConCinfty(\Reals^n)_\Null$
	for
	\begin{align*}
		\big(\deco{}{}{\sigma(I)}{}{s},\deco{}{s}{\sigma(I)}{}{}\big) \in
		\left((n^*)^{\strtensor s} \tensor (n^*)^{\strtensor r-s}\right)_\Null
		= (n^{\tensor s})^*_\Null \times (n^{\tensor r-s})^*_\Wobs
		\sqcup (n^{\tensor s})^*_\Wobs \times (n^{\tensor r-s})^*_\Null,
	\end{align*}
	and $D^I \in \ConCinfty(\Reals^n)_\Wobs$
	for
	\begin{align*}
		\big(\deco{}{}{\sigma(I)}{}{s},\deco{}{s}{\sigma(I)}{}{}\big) &\in
		\left((n^*)^{\strtensor s} \tensor (n^*)^{\strtensor r-s}\right)_\Wobs \setminus \left((n^*)^{\strtensor s} \tensor (n^*)^{\strtensor r-s}\right)_\Null\\
		&= (n^{\tensor s})^*_\red \times (n^{\tensor r-s})^*_\red.
	\end{align*}
	Here we used \autoref{lem:DecomposeStrTensorConIndSets}.
	This shows the first part.
	The second part follows then directly from
	\autoref{prop:LocalConMultDiffop}.
\end{proof}

Suppose $D \in \mathrm{C}^1_\diff(\ConCinfty(\mathcal{\Reals}^n))_\Total$
such that $\delta(D)$ is constraint.
Then we are interested in those parts of $D$ which are not constraint.
To separate the non-constraint part,
denote by
\begin{equation}
	\prol \colon \Cinfty(\Reals^{n_\Wobs}) \to \Cinfty(\Reals^{n_\Total})
\end{equation}
the constant extension of functions on $\Reals^{n_\Wobs}$ to functions
on $\Reals^{n_\Total}$.
With this we can always write
\begin{equation}
	f = \big(f - \prol(f\at{\Reals^{n_\Wobs}})\big) + \prol(f\at{\Reals^{n_\Wobs}}),
\end{equation}
splitting $f \in \Cinfty(\Reals^{n_\Total})$ into a part vanishing on the submanifold and the rest, thus we obtain a direct sum decomposition
\begin{equation}
	\Cinfty(\Reals^{n_\Total}) \simeq \vanishing_{\Reals^{n_\Wobs}} \oplus \Cinfty(\Reals^{n_\Wobs}).
\end{equation}
Since we can view $\Reals^{n_\red} \simeq \{0\}^{n_\Null} \oplus \Reals^{n_\Wobs \setminus n_\Null}$ as a subspace
of $\Reals^{n_\Wobs}$, we can similarly decompose $\Cinfty(\Reals^{n_\Wobs})$ to obtain
\begin{equation}
	\Cinfty(\Reals^{n_\Total}) \simeq \vanishing_{\Reals^{n_\Wobs}} \oplus \Cinfty(\Reals^{n_\red}) \oplus \vanishing_{\Reals^{n_\red}}(\Reals^{n_\Wobs}),
\end{equation}
with $\vanishing_{\Reals^{n_\red}}(\Reals^{n_\Wobs})$ denoting those functions on $\Reals^{n_\Wobs}$ vanishing
on the subspace $\Reals^{n_\red}$.
Note that $\ConCinfty(\Reals^n)_\Null = \vanishing_{\Reals^{n_\Wobs}}$
and $\ConCinfty(\Reals^n)_\Wobs = \vanishing_{\Reals^{n_\Wobs}} \oplus \Cinfty(\Reals^{n_\red})$,
thus $\vanishing_{\Reals^{n_\red}}(\Reals^{n_\Wobs})$ should be understood as a complement to
$\ConCinfty(\Reals^n)_\Wobs$ in $\ConCinfty(\Reals^n)_\Total$.
We will denote the projections to these summands by
\begin{align}
	\pr_\Null &\colon \Cinfty(\Reals^{n_\Total}) \to \vanishing_{\Reals^{n_\Wobs}}, \\
	\pr_\Null^\perp &\colon \Cinfty(\Reals^{n_\Total}) \to  \Cinfty(\Reals^{n_\red}), \\
	\pr_\Wobs \coloneqq \pr_\Null + \pr_\Null^\perp&\colon \Cinfty(\Reals^{n_\Total}) \to \vanishing_{\Reals^{n_\Wobs}} \oplus \Cinfty(\Reals^{n_\red}), \\
	\pr_\Wobs^\perp &\colon \Cinfty(\Reals^{n_\Total}) \to \vanishing_{\Reals^{n_\red}}(\Reals^{n_\Wobs}).	
\end{align}

We can find a similar decomposition of $\Diffop^r(\Reals^n)$:

\begin{proposition}
	\label{prop:DecompositionLocalDiffOps}
The $\Reals$-module maps
$\pr_\Null, \,\pr_\Null^\perp, \,\pr_\Wobs^\perp \colon C^1_\diff(\Cinfty(\Reals^{n_\Total})) \to C^1_\diff(\Cinfty(\Reals^{n_\Total}))$
defined by
\begin{align}
	\pr_\Null(D) &\coloneqq \sum_{I \in (n^{\tensor r})_\Null} D^I \del_I\\
	\pr_\Null^\perp(D) &\coloneqq \sum_{I \in (n_\red)^r} \pr_\Wobs(D^I) \del_I
	+ \sum_{{I \in (n^{\tensor r})^*_\Null}} \pr_\Null(D^I) \del_I \\
	\pr_\Wobs^\perp(D) &\coloneqq \sum_{I \in (n_\red)^r} \pr_\Wobs^\perp(D^I) \del_I
	+ \sum_{{I \in (n^{\tensor r})^*_\Null}} \pr_\Null^\perp(D^I) \del_I
\end{align}
are projections with
\begin{equation}
	\pr_\Null + \pr_\Null^\perp + \pr_\Wobs^\perp = \id,
\end{equation}
as well as
\begin{align}
	\image(\pr_\Null) = C^1_\diff(\ConCinfty(\Reals^n))_\Null \\
	\shortintertext{and}
	\image(\pr_\Wobs) = C^1_\diff(\ConCinfty(\Reals^n))_\Wobs
\end{align}
for $\pr_\Wobs \coloneqq \pr_\Null + \pr_\Null^\perp$.
\end{proposition}

\begin{proof}
Note that
\begin{equation*}
\begin{split}
	n_\Total^r = (n^{\tensor r})_\Null \sqcup \Big((n^{\tensor r})_\Wobs \setminus (n^{\tensor r})_\Null \Big)
	\sqcup \big((n^{\tensor r})_\Total \setminus (n^{\tensor r})_\Wobs \Big)
	= (n^{\tensor r})_\Null \sqcup (n_\red)^r	\sqcup (n^{\tensor r})^*_\Null,
\end{split}
\end{equation*}
with $(n^{\tensor r})_\Wobs = (n^{\tensor r})_\Null \sqcup (n_\red)^r$.
Then, with the help of \autoref{ex:ParitalsAsConDiffOps},
every $D \in C^1_\diff(\ConCinfty(\Reals^n))_\Total$ of order $r$ can uniquely be written as
\begin{equation*}
\begin{split}
	D = &\sum_{I \in (n^{\tensor r})_\Null} D^I \del_I
	+ \sum_{I \in (n_\red)^r} \pr_\Wobs(D^I) \del_I
	+ \sum_{{I \in (n^{\tensor r})^*_\Null}} \pr_\Null(D^I) \del_I \\
	&+ \sum_{I \in (n_\red)^r} \pr_\Wobs^\perp(D^I) \del_I
	+ \sum_{{I \in (n^{\tensor r})^*_\Null}} \pr_\Null^\perp(D^I) \del_I
\end{split}
\end{equation*}
with
\begin{align*}
	\sum_{I \in (n^{\tensor r})_\Null} D^I \del_I &\in C^1_\diff(\ConCinfty(\Reals^n))_\Null, \\
	\sum_{I \in (n_\red)^r} \pr_\Wobs(D^I) \del_I &\in C^1_\diff(\ConCinfty(\Reals^n))_\Wobs \\
	\shortintertext{and}
	\sum_{{I \in (n^{\tensor r})^*_\Null}} \pr_\Null(D^I) \del_I &\in C^1_\diff(\ConCinfty(\Reals^n))_\Wobs.
\end{align*}
Thus $\pr_\Null$, $\pr_\Null^\perp$ and $\pr_\Wobs^\perp$ are indeed projections with
$\image(\pr_\Null) \subseteq C^1_\diff(\ConCinfty(\Reals^n))_\Null$
and $\image(\pr_\Wobs) \subseteq C^1_\diff(\ConCinfty(\Reals^n))_\Wobs$.
The surjectivity of these maps follows from evaluating at
$x^{i_1} \cdots x^{i_r}$.
\end{proof}

This shows that $D \in C^1_\diff(\ConCinfty(\Reals^n))_\Wobs$ is constraint if and only if
$\pr_\Wobs^\perp(D) = 0$.
Suppose again that $\delta(D)$ is constraint, then by \autoref{prop:LocalConClosedOneCochains} we know that
$\pr_\Wobs^\perp(D^I) = 0$ for all $I \in n_\red^r$
and
$\pr_\Null(D^I) = 0$ whenever there exist $k,\ell \in \{1, \dotsc, r\}$ with $k \neq \ell$
such that $i_\ell \in n_\Total \setminus n_\Wobs$ and $i_k \in n_\Total \setminus n_\Null$.
Hence \autoref{prop:DecompositionLocalDiffOps}
shows that that the constraint Hochschild $2$-cochains which are exact but not constraint exact, are
those differential operators of order $r$ with
\begin{equation}
	\pr_\Wobs^\perp(D) = \sum_{I \in S} \pr_\Null^\perp(D^I) \del_I \neq 0
\end{equation}
with
\begin{equation}
	S_r \coloneqq \left\{ I \in n_\Total^r \mid \exists \ell \in \{1,\dotsc, r \}: i_\ell \in n_\Total \setminus n_\Wobs
		\text{ and } \forall k\neq \ell: i_k \in n_\Null \right\},
\end{equation}
i.e. which differentiate once in a direction perpendicular to the subspace $\Reals^{n_\Wobs}$
and $(r-1)$-times in direction of the distribution $\Reals^{n_\Null}$.
Using the constraint symbol calculus from \autoref{sec:ConSymbCalculus} leads us to the following definition.

\begin{definition}[Extended constraint bivector fields]
	\label{def:ExtendedConBivector}
	\index{constraint!HKR map}
	\glsadd{ExtConBivect}
For the constraint manifold $\Reals^n=(\Reals^{n_\Total},\Reals^{n_\Wobs},\Reals^{n_\Null})$
we define the strong constraint $\ConCinfty(\Reals^n)$-module
$\ConVecFields_\ext^k(\Reals^n)$ of \emph{extended constraint bivector fields} by
\begin{equation} \label{eq:ExtendedConBivector}
\begin{split}
	\ConVecFields_\ext^2(\Reals^n)_{\Total}
	&\coloneqq \ConVecFields_{\strtensor}^2(\Reals^n)_\Total, \\
	\ConVecFields_\ext^2(\Reals^n)_{\Wobs}
	&\coloneqq \ConVecFields_{\strtensor}^2(\Reals^n)_{\Wobs}
	\oplus \Big(\bigoplus_{k=1}^\infty \Sym^k \Secinfty(T\Reals^{n_\Null}\at{\Reals^{n_\Wobs}}) \vee \Secinfty(T\Reals^{n_\Total - n_\Wobs}\at{\Reals^{n_\Wobs}}) \Big),\\
	\ConVecFields_\ext^2(\Reals^n)_{\Null}
	&\coloneqq \ConVecFields_{\strtensor}^2(\Reals^n)_{\Null}
	\oplus \Big(\bigoplus_{k=1}^\infty \Sym^k \Secinfty(T\Reals^{n_\Null}\at{\Reals^{n_\Wobs}}) \vee \Secinfty(T\Reals^{n_\Total - n_\Wobs}\at{\Reals^{n_\Wobs}}) \Big),\\
\end{split}
\end{equation}
with $\iota_{\ext} \colon \ConVecFields_\ext^2(\Reals^n)_\Wobs \ni (X,D) \mapsto X \in \ConVecFields_\ext^2(\Reals^n)_\Total$.
\end{definition}

It is important to remark that $\ConVecFields_\ext^2(\Reals^n)$
is not embedded.
The additional terms in \eqref{eq:ExtendedConBivector}
should be interpreted as certain higher order differential operators 
living only on the submanifold $\Reals^{n_\Wobs}$.
To make this identification precise, define for every 
$D = D^{i_1,\dotsc,i_r} \frac{\del}{\del x^{i_1}} \vee \dots \vee \frac{\del}{\del x^{i_r}}\in \Secinfty(\Sym^rT\Reals^{n_\Total}\at{\Reals^{n_\Wobs}})$ on $\Reals^{n_\Wobs}$
its \emph{prolongation}
\begin{equation}
	\prol(D) \coloneqq \prol(D^{i_1,\dotsc,i_r}) \frac{\del}{\del x^{i_1}} \vee \dots \vee \frac{\del}{\del x^{i_r}}
	\in \Secinfty(\Sym^r T\Reals^{n_\Total})
\end{equation}
by extending the coefficient functions to $\Reals^{n_\Total}$ in a constant fashion.
Since the constraint manifold $\Reals^n$ carries a
canonical constraint covariant derivative, see \autoref{ex:CanonicalConCovDerivative},
we can then identify $\prol(D)$ with a
differential operator.

We now want to extend the morphism 
$\mathcal{U}$ from \autoref{prop:ConHKRmap} to include these new terms:

\begin{proposition}[Extended constraint HKR map]
	\label{prop:ExtendedConHKRMap}
	\glsadd{ExtConHKRMap}
Consider the constraint manifold
$\Reals^n = (\Reals^{n_\Total}, \Reals^{n_\Wobs}, \Reals^{n_\Null})$.
	The map $\mathcal{U}_\ext \colon \ConVecFields_{\ext}^2(\Reals^n) \to \Hochschild_\diff^2(\ConCinfty(\Reals^n))$
	defined by 
	\begin{equation} \label{eq:ExtendedConHKRmap_Cohom}
	\begin{split}
		(\mathcal{U}_\ext)_\Total(X) &\coloneqq \mathcal{U}(X) \\
		(\mathcal{U}_\ext)_\Wobs(X,D) &\coloneqq \mathcal{U}(X) + \delta\big(\Op(\prol(D))\big)
	\end{split}
	\end{equation}
	is a morphism between constraint $\field{k}$-modules
	and a regular monomorphism.
\end{proposition}

\begin{proof}

Recall from \autoref{thm:ClassicalHKR}
that $\mathcal{U}_\ext$ is an isomorphism on the $\TOTAL$-components.
Moreover, $(\mathcal{U}_\ext)_\Wobs$
clearly preserves the $\NULL$-component
and $\delta\big(\Op(\prol(D))\big)$
vanishes in $\Hochschild_\diff^2(\Cinfty(\Reals^{n_\Total}))$.
Thus $\mathcal{U}_\ext$ is a constraint morphism.
Now assume that $[(\mathcal{U}_\ext)_\Wobs(X,D)] = 0$.
Since $\mathcal{U}(X)$ is an antisymmetric bidifferential operator
and $\delta(\Op(\prol(D)))$ is a symmetric bidifferential operator
these two parts have to vanish separately in cohomology.
Then from \autoref{prop:ConHKRmap} it follows
$X = 0$.
To show that also $[\delta(\Op(\prol(D)))] = 0$ assume that
there exists
$\tilde{D} = \sum_{r=0}^{k} \sum_{I \in n^{\tensor r}} \frac{1}{r!} \tilde{D}^I \del_I \in \mathrm{C}^1_\diff(\ConCinfty(\Reals^n))_\Wobs$
such that $\delta(\Op(\prol(D))) = \delta(\tilde{D})$.
Then $\Op(\prol(D)) - \tilde{D}$ is closed and hence
a derivation.
Since $\Op(\prol(D))$ is a differential operator of order at least $2$,
we obtain
\begin{equation*}
	\Op(\prol(D)) = \sum_{r=2}^{k} \sum_{I \in n^{\tensor r}} \frac{1}{r!} \tilde{D}^I \del_I.
\end{equation*}
From \autoref{prop:LocalConDiffopCinfty} it follows that $\Op(\prol(D))$, and thus also $\tilde{D}$, is not constraint, 
giving a contradiction to $\tilde{D} \in \mathrm{C}^1_\diff(\ConCinfty(\Reals^n))_\Wobs$.
This shows that \eqref{eq:ExtendedConHKRmap_Cohom} is a monomorphism.
For its regularity suppose that $[(\mathcal{U}_\ext)_\Wobs(X,D)] \in \Hochschild_\diff^2(\ConCinfty(\Reals^n))_\Null$.
By \autoref{def:ExtendedConBivector} we have
$D \in \ConVecFields_\ext^2(\Reals^n)_\Null$, and thus
$[\delta(\Op(\prol(D)))] \in \Hochschild_\diff^2(\ConCinfty(\Reals^n))_\Null$.
Then from
\begin{equation*}
	[\mathcal{U}(X)] = [(\mathcal{U}_\ext)_\Wobs(X,D)] - [\delta(\Op(\prol(D)))] \in \Hochschild_\diff^2(\ConCinfty(\Reals^n))_\Null
\end{equation*}
it follows from the fact that $\mathcal{U}$ is a regular monomorphism,
see \autoref{prop:ConHKRmap},
that $X \in \ConVecFields_\ext^2(\Reals^n)_\Null$.
\end{proof}

With this we have found contributions to the second constraint Hochschild cohomology which go beyond 
the classical Hochschild cohomology as computed by the HKR theorem.
The next and final theorem shows that no other contributions appear.

\begin{theorem}[Second constraint Hochschild cohomology on $\Reals^n$]
	\label{thm:LocalSecondConHochschild}
	\index{constraint!Hochschild cohomology}
The morphism
\begin{equation}
	\mathcal{U}_\ext \colon \ConVecFields_{\ext}^2(\Reals^n) \to \Hochschild_\diff^2\big(\ConCinfty(\Reals^n)\big)
\end{equation}
as defined in \autoref{prop:ExtendedConHKRMap} is an isomorphism of constraint $\Reals$-modules.
\end{theorem}

\begin{proof}
It remains to show that $\mathcal{U}_\ext$ is an epimorphism.
On the $\TOTAL$-component it is an epimorphism by \autoref{thm:ClassicalHKR}.
To show the surjectivity on the $\WOBS$-component let 
$B \in \mathrm{C}_\diff^2(\ConCinfty(\Reals^n))_\Wobs$ be given with
$\delta(B) = 0$.
Then the classical HKR theorem tells us that we can write
$B = \delta(D) + \AntiSymmetrizer(B)$
with $D \in \mathrm{C}_\diff^1(\ConCinfty(\Reals^n))_\Total$
and $\AntiSymmetrizer(B) \in \ConDiffop^{(1,1)}(\Reals^n)_\Wobs$
the antisymmetric part of $B$.
From this it follows $\delta(D) \in \mathrm{C}_\diff^1(\ConCinfty(\Reals^n))_\Wobs$.
By \autoref{prop:DecompositionLocalDiffOps} $D$ splits as
$D = \pr_\Wobs(D) \pr_\Wobs^\perp(D)$, with
$\pr_\Wobs(D) \in \mathrm{C}_\diff^1(\ConCinfty(\Reals^n))_\Wobs$
and 
\begin{equation*}
	\pr_\Wobs^\perp(D) = \Op\bigg(\prol\Big( \sum_{I \in S} D^I\at[\Big]{\Reals^{n_\Wobs}} \del_{i_1} \vee \dots \vee \del_{i_r}\Big)\bigg),
\end{equation*}
where $S = \bigcup_{r=0}^\infty \left\{ I \in n_\Total^r \mid \exists \ell \in \{1,\dotsc, r\}: i_\ell \in n_\Total \setminus n_\Wobs
\text{ and } \forall k\neq \ell: i_k \in n_\Null \right\}$.
Thus
\begin{equation*}
	B = \delta(\pr_\Wobs(D)) + (\mathcal{U}_\ext)_\Wobs\Big(\sigma(X), \sum_{I \in S} D^I\at[\Big]{\Reals^{n_\Wobs}} \del_{i_1} \vee \dots \vee \del_{i_r}\Big),
\end{equation*}
showing that $\mathcal{U}_\ext \colon \ConVecFields_{\ext}^2(\Reals^n) \to \Hochschild_\diff^2(\ConCinfty(\Reals^n))$
is surjective on the $\WOBS$-components, and therefore an isomorphism.
\end{proof}

\autoref{thm:Obstructions} shows that the second constraint Hochschild cohomology can be interpreted as the constraint set of equivalence classes of infinitesimal deformations.
More precisely,\linebreak
$\Hochschild_\diff^2(\ConCinfty(\Reals^n))_\Total$ is the set of equivalence classes of classical infinitesimal deformations of $\Cinfty(M)$,
while $\Hochschild_\diff^2(\ConCinfty(\Reals^n))_\Wobs$ are equivalence classes
of constraint infinitesimal deformations, i.e. deformations which respect the reduction information.
In the local case of $\mathcal{M} = \Reals^n$ we see that
$\Hochschild_\diff^2(\ConCinfty(\Reals^n))$ is not embedded, which means there are
non-equivalent constraint deformations which are equivalent when we forget about the reduction data.
And these equivalence classes are exactly characterized by
the additional symmetric parts in $\ConVecFields_\ext^2(\Reals^n)_\Wobs$, see \eqref{eq:ExtendedConBivector}.

\subsubsection{Reduction}

Observe that these symmetric contributions also appear in $\ConVecFields_\ext^2(\Reals^n)_\Null$,
and hence should vanish after reduction.
More precisely, we have the following statement:

\begin{proposition}
	\index{reduction!HKR map}
Consider the constraint manifold $\Reals^n =  (\Reals^{n_\Total}, \Reals^{n_\Wobs}, \Reals^{n_\Null})$.
\begin{propositionlist}
	\item The morphism $\mathcal{U}_\ext \colon \ConVecFields_\ext^2(\Reals^n) \to \mathrm{C}_\diff^2(\ConCinfty(\Reals^n))$
	reduces to the classical HKR map
	\begin{equation}
		(\mathcal{U}_\ext)_\red \colon \VecFields^2(\Reals^{n_\red}) \to \mathrm{C}_\diff^2(\Cinfty(\Reals^{n_\red}))
	\end{equation}
	on $\Reals^{n_\red}$.
	\item The isomorphism $\mathcal{U}_\ext \colon \ConVecFields_{\ext}^2(\Reals^n) \to \Hochschild_\diff^2(\ConCinfty(\Reals^n))$
	reduces to the classical HKR isomorphism
	\begin{equation}
		(\mathcal{U}_\ext)_\red \colon \VecFields^2(\Reals^{n_\red}) \to \Hochschild_\diff^2(\Cinfty(\Reals^{n_\red}))
	\end{equation}
	on $\Reals^n$.
\end{propositionlist}
\end{proposition}

\begin{proof}
	For the first part note that 
	$(\ConVecFields_\ext^2(\Reals^n))_\red \simeq \VecFields^2(\Reals^{n_\red})$ holds since
	$(\ConVecFields_{\strtensor}^2(\Reals^n))_\red \simeq \VecFields^2(\Reals^{n_\red})$
	by \autoref{prop:ConVecFieldsVSRedcution} and the additional symmetric terms vanish after reduction.
	Moreover, $(\mathrm{C}_\diff^2(\ConCinfty(\Reals^n)))_\red \simeq \mathrm{C}_\diff^2(\Cinfty(\Reals^{n_\red}))$
	holds by \autoref{prop:HochschildVSReduction} and the fact that every multidifferential operator on 
	$\Reals^{n_\Wobs -n_\Null}$ can be extended to a constraint multidifferential operator on
	$\Reals^n$.
	Then $(\mathcal{U}_\ext)_\red$ becomes the classical HKR map, by its explicit definition in
	\eqref{eq:ConHKRmap}.
	
	The second part follows since taking cohomology commutes with reduction as we know from
	by \autoref{prop:CohomologyCommutesWithReduction}.
\end{proof}

%% file: outlook.tex
We have established in this thesis a general framework which allows to treat geometric and algebraic features of coisotropic reduction on equal footing.
This allowed us to introduce constraint star products, which are essentially star products compatible with reduction.
These induce automatically star products on the reduced spaces, and therefore quantization commutes with reduction in this setting.
Nevertheless, the existence of such constraint star products is not obvious, and we adapted classical techniques from deformation theory to establish constraint Hochschild cohomology, which governs the deformation problem of constraint algebras.
As a first step towards a constraint HKR Theorem we were able to compute the zeroth and first constraint Hochschild cohomologies in the general situation and the second constraint Hochschild cohomology in the flat case.
This second constraint Hochschild cohomology turned out to contain symmetric terms of arbitrary differentiation order, which are unexpected from the point of view of the classical HKR Theorem.
This leads to the following open questions, that should be studied in future projects:
\begin{cptitem}
	\item \index{constraint!Hochschild-Kostant-Rosenberg Theorem}
	The explicit characterization in \autoref{thm:LocalSecondConHochschild} of the second constraint Hochschild cohomology 
	$\Hochschild_\diff^2(\ConCinfty(\Reals^n))$
	gives strong hints on how the higher constraint Hochschild cohomologies may be
	described.
	Besides the constraint multivector fields $ \ConVecFields_{\strtensor}^\bullet(\mathcal{M})$
	we expect contributions given by constraint Hochschild cochains which are exact with non-constraint potentials.
	Such a potential $\phi$ should differentiate only $k$ times in the direction of
	$\Reals^{n_\Total - n_\Wobs}$, where $k$ is the number of slots, and at least once in the direction of the distribution $\Reals^{n_\Null}$, since then $\delta(\phi)$ will have at least one factor in the $\NULL$-component of the constraint differential operators, making $\delta(\phi)$ itself constraint.
	It then needs to be shown that all additional contributions appearing in higher orders of constraint Hochschild cohomology are of this special form.
	\item \index{constraint!Hochschild-Kostant-Rosenberg Theorem}
	Globalizing a constraint HKR Theorem for $\Reals^n$ to an arbitrary constraint manifold $\mathcal{M}$ will not always be possible, since there need not exist partitions of unity compatible with the constraint structure.
	Thus classical proofs for the HKR Theorem that use such a glueing procedure, as can be found e.g. in \cite{gutt.rawnsley:1999a}, cannot directly be applied in the constraint situation.
	Instead it seems reasonable to take a classical proof of the HKR Theorem which is inherently global \cite{dewilde.lecomte:1995a}, and reformulate this in the constraint framework.
	The case of $\Reals^n$ already suggests that a constraint HKR map depends on the choice of a constraint covariant derivative.
	Whether the resulting isomorphism in cohomology really depends on that choice remains to be seen.
	\item A constraint algebra $\algebra{A}$ can equivalently be understood as a span
	$\algebra{A}_\red \twoheadleftarrow \algebra{A}_\Wobs \to \algebra{A}_\Total$ of associative algebras.
	Deformations of such diagrams of algebras have been studied e.g. in \cite{fregier.markl.yau:2009, fregier.zambon:2015,gerstenhaber.schack:1983a}.
	This deformation theory of diagrams deforms the algebras as well as the morphisms of the diagram, while
	for a deformation of constraint algebras we only want to deform the algebras.
	Moreover, the category of modules over such diagrams is abelian, while the category of constraint modules is not.
	Thus, even though the deformation theory of constraint algebras is obviously linked to the deformation theory of diagrams,
	we have to expect differences in the details.
	The exact relationship between these deformation theories is yet to be uncovered.
	\item Constraint manifolds were introduced using simple distributions, but as already discussed in \autoref{rem:MoreGeneralConManifolds}
	it would be useful to allow for more general quotient procedures.
	On one hand we could allow for general equivalence relations which still provide a smooth quotient space.
	In this situation most of the results of constraint differential geometry as presented in \autoref{chap:ConstraintGeometricStructures} should still hold.
	On the other hand, we might want to allow for more singular reduction.
	In this case, one might abandon the geometry completely and instead focus on its algebraic description using
	constraint algebras, or one could enlarge the categories of geometric objects we allow.
	For example we could study constraint versions of orbifolds, diffeological spaces etc.
	The properties of these constraint objects will then greatly rely on the categories of objects they depend on.
	\item Based on the differential geometry of constraint manifolds, as introduced in \autoref{chap:ConstraintGeometricStructures}, the reduction of more sophisticated
	geometric objects, such as Lie (bi-)algebroids, can be investigated, see \cite{dippell.kern:2022a}.
	\item Strong constraint manifolds, i.e. constraint manifolds with globally defined equivalence relations, are natural objects to study.
	These can be understood as generalizations of Marsden-Weinstein reduction, instead of coisotropic reduction, where
	the global distribution comes from a well-behaved global group action of a Lie group $\group{G}$ on a manifold $M$.
	Functions on such strong constraint manifolds coming from Marsden-Weinstein reduction would form non-strong constraint
	algebras, consisting of globally invariant functions $\Cinfty(M)^\group{G}$ in the $\WOBS$-component and 
	globally invariant functions vanishing on the submanifold $\vanishing_C \cap \Cinfty(M)^\group{G}$
	in the $\NULL$-component.
	See \cite{sniatycki.weinstein:1983a} for a formulation of Marsden-Weinstein reduction in terms of these classes of functions.
	\item The reduction of differential operators and multivector fields in the setting of Hamiltonian Lie group actions was studied in \cite{esposito.kraft.schnitzer:2022a,esposito.kraft.schnitzer:2022b} using $L_\infty$-algebras.
	There, reduction of differential operators and multivector fields is encoded in an $L_\infty$-morphism to the reduced objects.
	To bring this in contact with our constraint reduction scheme it should be useful to introduce constraint
	$L_\infty$-algebras and morphisms, based on our notion of constraint DGLAs.
	\item \index{Morita equivalence}
	The bicategories $\ConBimod$ and $\strConBimod$ suggest to study the representation theory of (strong) constraint algebras from a Morita theoretic perspective, see \autoref{rem:MoritaEquivalence}.
	This has been done for a special class of constraint algebras in \cite{dippell.esposito.waldmann:2019a}.
	Besides the purely algebraic insights this will entail, representation theories of constraint algebras is also interesting from the point of view of deformation quantization.
	To bring a formal deformation of a (constraint) algebra of functions into contact with physics we need to choose a suitable representation, hence it would be desirable to compare the representation theories via Morita theory.
	\item \index{K-Theory}
	The introduction of projective constraint modules in \autoref{sec:RegularProjectiveModules} suggests to define a constraint version of algebraic $K$-theory, which might be the first step towards a constraint algebraic index theorem, i.e. an algebraic index theorem compatible with reduction.
\end{cptitem}

%% file: appendix-monoids.tex
We will give the basic definitions of category theory here, not least to fix our notation.
See \cite{maclane:1998a} for the standard textbook on category theory or for example \cite{kashiwara.schapira:2006a, 
	brandenburg:2016a} for more modern introductions.
\autoref{sec: categories and morphisms} to \autoref{sec:MonoidsModules} are mainly taken from
\cite{dippell:2018a}.

Category theory is a branch of mathematics that tries to reveal the underlying mechanics of constructions done in different 
branches of mathematics, in order to uncover the common features and to allow to transfer techniques from one field of 
mathematics to another.
As such category theory takes a bird's eye perspective of mathematics, leading us to consider such things as the collection of all 
vector spaces or of all sets, etc.
Here one might get suspicious, since this sounds a lot like we immediately run into Russel's paradox.
To avoid this we do not consider the set of all sets, but the collection of all sets.
What we mean by collection is now depending on the foundations of category theory we choose.
For our purposes it will be enough to be aware that a collection can be bigger than a set, and does not need to share all of the 
properties we are used to from axiomatic systems like ZFC.
For an overview over possible foundations of category theory see \cite{shulman:2008a}.

\section{Categories and Morphisms}
\label{sec: categories and morphisms}
In this section we will give the basic definitions of categories and examine some important properties of morphisms.

\begin{definition}[Category]
	\label{def: category}
	\index{category}
A \emph{category} $\category{C}$ consists of the following data:
\begin{definitionlist}
	\item A collection $\category{C}_0$ of \emph{objects}.
	\item For any two objects $A,B \in \category{C}_0$ a set
	$\category{C}(B,A) = \Hom(A,B)$ of \emph{morphisms} from $A$ to $B$, called \emph{hom-set},
	where $f \in \category{C}(B,A)$ will be written $f \colon A \longrightarrow B$.
	\item For any three objects $A,B,C \in \category{C}_0$ a map
	$\circ \colon \category{C}(C,B) \times  \category{C}(B,A) \longrightarrow \category{C}(C,A)$,
	which assigns to any appropriate pair of morphisms $f,g$ their \emph{composition} $f \circ g$.
	\item For each object $A \in \category{C}_0$ a morphism $\id_A \in \category{C}(A,A)$, called the
	\emph{identity morphism} at $A$.
\end{definitionlist}
These data are required to fulfil the following properties:
\begin{definitionlist}
	\item \emph{Associativity}: For any four objects $A,B,C,D \in \category{C}_0$ and any
	$f \in \category{C}(D,C)$,
	$g \in \category{C}(C,B)$ and
	$h \in \category{C}(B,A)$ it holds
	\begin{equation}
		(f \circ g) \circ h = f \circ (g \circ h).
	\end{equation}
	\item \emph{Left and right identity laws}: For any $A,B \in \category{C}_0$ and any $f \in \category{C}(B,A)$
	it holds
	\begin{equation}
		\id_B \circ f = f = f \circ \id_A.
	\end{equation}
\end{definitionlist}
\end{definition}

If it is clear that we are talking about objects of a given category we will often drop the subscript and simply write 
$\category{C}$ instead of $\category{C}_0$.
Thus by $C \in \category{C}$ we mean an object of the category $\category{C}$.
Note also that the order of objects in our notation of hom-sets is different from the standard notation.
What we call category is sometimes called a locally-small category in the literature, but since we will not need categories with 
hom-sets being mere classes instead of sets we will stick to this convention.
A category $\category{C}$ where the collection $\category{C}_0$ of objects is a set is called a \emph{small category}, whereas a 
category with $\category{C}_0$ not being a set is called a \emph{large category}.

\begin{example}[Categories]\
	\begin{examplelist}
		\item The \emph{trivial category} $\mathsf{1}$ consists of one object $* \in \mathsf{1}_0 = \{*\}$ and one morphism 
		$\id_* 	\in \mathsf{1}_1(*,*)$.
		\item The \emph{interval category} $\mathsf{2}$ consists of two objects $0$ and $1$ and three morphisms;
		the identities on $0$ and $1$ and exactly one morphism $0 \longrightarrow 1$.
		\item Any class of objects can be turned into a category by adding only the identity morphisms for every object.
		Categories of this kind are called \emph{discrete}.
		\item Given any category $\category{C}$ with composition $\circ$ we can build the \emph{opposite category}
		$\category{C}^\opp$ by keeping the objects $\category{C}^\opp_0 = \category{C}_0$ but using the inverted hom-sets
		$\category{C}^\opp(B,A) = \category{C}(A,B)$ with composition $f \circ_\opp g = g \circ f$.
	\end{examplelist}
\end{example}

One important way to construct a category out of two given categories is by taking their product.

\begin{definition}[Product category]
	\index{product!category}
Let $\category{C}$ and $\category{D}$ be two categories.
The \emph{product category} $\category{C} \times \category{D}$ is the category with
\begin{definitionlist}
	\item objects being ordered pairs $(C,D)$ of objects $C \in \category{C}$ and $D \in \category{D}$,
	\item morphisms being pairs $(f,g) \colon (C,D) \longrightarrow (C^\prime,D^\prime)$ of morphisms
	$f \in \category{C}(C^\prime,C)$ and $g \in \category{D}(D^\prime, D)$,
	\item composition of morphisms $(f,g) \colon (C,D) \longrightarrow (C^\prime, D^\prime)$ and
	$(f^\prime,g^\prime ) \colon (C^\prime, D^\prime) \longrightarrow (C^{\prime\prime},D^{\prime\prime})$
	given by the componentwise composition
	$(f^\prime \circ f , g^\prime \circ g) \colon  (C,D) \longrightarrow (C^{\prime\prime},D^{\prime\prime})$ and
	\item identity morphisms given by pairs $(\id_C,\id_D)$ of the identity morphisms in $\category{C}$ and $\category{D}$.
\end{definitionlist}
\end{definition}

We think of morphisms between objects as different ways to relate these objects, and the hom-sets consist only of those 
morphisms that respect the inner structure of the objects.
This is how we construct most categories, but following the philosophy of category theory we should actually think of this the other 
way around: The inner structure of an object is determined by all possible ways of relating it to other objects.
Following this idea we cannot distinguish two objects that behave in the same way in relation to all other objects.
Thus, we should not think of objects being equal, but only isomorphic in the following sense.

\begin{definition}[Isomorphism]
	Let $\category{C}$ be a category. Two objects $a,b \in \category{C}_0$ are \emph{isomorphic} if there exist morphisms
	$f \colon A \longrightarrow B$ and 
	$g \colon B \longrightarrow A$ such that
	$g \circ f = \id_A$ and $f \circ g = \id_B$ hold.
	The morphisms $f$ and $g$ are called \emph{isomorphisms}.
\end{definition}

The idea that the only notion of sameness in a category is that of being isomorphic is sometimes called the principle of equivalence.

In the category $\Sets$ of sets it is easy to show that a function between sets is injective if and only if it is left cancellable and a 
function is surjective if and only if it is right cancellable.
This allows us to transfer these notions to arbitrary categories, where we cannot talk about elements of an object, but only about 
morphisms between objects.

\begin{definition}[Monomorphism]
	\label{def:Monomorphism}
	\index{monomorphism}
Let $\category{C}$ be a category. A morphism $f \colon B \longrightarrow C$ is called \emph{monomorphism} (or 
\emph{mono} for 
short) if for all morphisms $g_1, g_2 \colon A \longrightarrow B$ it holds
\begin{equation}
	f \circ g_1 = f \circ g_2 \implies g_1 = g_2.
\end{equation}
\end{definition}

If we want to highlight the fact that a morphism is a monomorphism we will depict it as
\begin{equation*}
	\begin{tikzcd}
		A \arrow[hookrightarrow]{r}{f} &B
	\end{tikzcd}
\end{equation*}
in diagrams.

\begin{proposition}
	In a category $\category{C}$ the following statements hold:
	\begin{propositionlist}
		\item Every isomorphism is a monomorphism.
		\item The composition of two monomorphisms is a monomorphism.
		\item If the composition $f \circ g$ is a monomorphism, then $g$ is a monomorphism.
	\end{propositionlist}
\end{proposition}

\begin{definition}[Epimorphism]
	\label{def:Epimorphism}
	\index{epimorphism}
Let $\category{C}$ be a category. A morphism $f \colon A \longrightarrow B$ is called \emph{epimorphism} (or \emph{epi} 
for short) if 
for all morphisms $g_1,g_2 \colon B \longrightarrow C$ it holds
\begin{equation}
	g_1 \circ f = g_2 \circ f \implies g_1 = g_2.
\end{equation}
\end{definition}

If we want to stress that a morphism is an epimorphism we will depict it as
\begin{equation*}
	\begin{tikzcd}
		A \arrow[twoheadrightarrow]{r}{f} &B
	\end{tikzcd}
\end{equation*}
in diagrams.

\begin{proposition}
	In a category $\category{C}$ the following statements hold:
	\begin{propositionlist}
		\item Every isomorphism is an epimorphism.
		\item The composition of two epimorphisms is an epimorphism.
		\item If the composition $f \circ g$ is an epimorphism, then $f$ is an epimorphism.
	\end{propositionlist}
\end{proposition}

Note that in general categories every isomorphism is a mono and epi, but not every morphism that is mono and epi has to be an 
isomorphism.

The following more special classes of monos and epis occur quite often.

\begin{definition}[Section \& retraction]
	Let $\category{C}$ be a category and let morphisms $f \colon A \longrightarrow B$ and $g \colon B \longrightarrow A$
	be given, such that
	$g \circ f = \id_A$.
	Then $f$ is called a \emph{section of} $g$ and $g$ is called a \emph{retraction of} $f$.
	Furthermore, $A$ is called a \emph{retract} of $B$.
\end{definition}

\begin{lemma}
	Every section is a monomorphism and every retraction is an epimorphism.
\end{lemma}

\begin{proof}
	Let $f \colon A \longrightarrow B$ and $g \colon B \longrightarrow A$ such that $g \circ f = \id_A$.
	Let furthermore $h_1,h_2 \colon X \longrightarrow A$ such that
	$f \circ h_1 = f \circ h_2$, then
	$g \circ f \circ h_1 = g \circ f \circ h_2$ and thus $h_1 = h_2$.
	Hence $f$ is a monomorphism.
	Now let $k_1, k_2 \colon A \longrightarrow Y$ such that $k_1 \circ g = k_2 \circ g$.
	Then from $k_1 \circ g \circ f = k_2 \circ g \circ f$ follows $k_1 = k_2$ and hence $g$ is an epimorphism.
\end{proof}

We will also use the terms \emph{split monomorphism} for sections and \emph{split epimorphism} for retractions.
In many categories we encounter objects with very small hom-sets for every other object.
\filbreak
\begin{definition}[Initial, terminal \& zero object]
Let $\category{C}$ be a category, $C \in \category{C}$.
\begin{definitionlist}
	\item \index{initial object}
	$C$ is called \emph{initial object} if for every object $D \in \category{C}$ there exists a unique morphism
	$f \colon C \longrightarrow D$.
	\item \index{terminal object}
	$C$ is called \emph{terminal object} if for every object $B \in \category{C}$ there exists a unique morphism
	$g \colon B \longrightarrow C$.
	\item $C$ is called \emph{zero object} if it is initial and terminal.
\end{definitionlist}
\end{definition}

We will mostly use $0$ for zero objects.
Note that initial and terminal objects, and hence zero objects as well, are unique up to isomorphisms, so we often speak of
\emph{the} initial, terminal or zero object. 

The existence of a zero object also allows to speak of zero morphisms.

\begin{definition}[Zero morphism]
Let $\category{C}$ be a category with zero object $0$.
The \emph{zero morphism}
$0_{A,B} \colon A \longrightarrow B$
between two objects $A$ and $B$ is the unique morphism that factors through $0$.
\end{definition}

In general we say that a morphism $f \colon A \longrightarrow C$ \emph{factors through} $B$ if there exist morphisms
$g \colon A \longrightarrow B$ and $h \colon B \longrightarrow C$ such that $f = h \circ g$.

Using the existence of a zero object we can generalize the concept of kernel of a linear map between vector spaces.

\begin{definition}[Kernel]
	\index{kernel}
Let $\category{C}$ be a category with zero object $0$ and let $f \colon A \longrightarrow B$. An object $K$ together with a 
morphism 
$k \colon K \longrightarrow A$ is called \emph{kernel of} $f$ if it satisfies the following universal property: it holds $f \circ k = 
0_{K,B}$ 
and for any morphism $k^\prime \colon K^\prime \longrightarrow A$ such that $f \circ k^\prime = 0_{K^\prime ,B}$ there is a 
unique 
morphism $u \colon K^\prime \longrightarrow K$ such that $k \circ u = k^\prime$. Expressed as a diagram:
\begin{equation}
	\begin{tikzcd}
		{ }
		& A
		\arrow{dr}{f}
		&{ } \\
		{ }
		& K
		\arrow{u}{k}
		\arrow{r}{0_{K,B}}
		& B \\
		K^\prime
		\arrow{uur}{k^\prime}
		\arrow[dashed]{ur}{u}
		\arrow{urr}[swap]{0_{K^\prime, B}}
		& { }
		& { }				
	\end{tikzcd}.
\end{equation}
\end{definition}

It is clear that the kernel is unique up to isomorphism if it exists at all.
We will also write $\ker(f)$ for the kernel morphism of $f$ and $\Ker(f)$ for the kernel object of $f$.

\begin{corollary} \label{cor: kernel is mono}
	Every kernel is a monomorphism.
\end{corollary}

A useful observation is that if $f$ is a monomorphism its kernel is the zero object together with the zero morphism.
Another important case is that if $0 \colon A \longrightarrow B$ is the zero morphism, then the kernel is clearly $A$ together with 
the identity morphism.

\begin{definition}[Cokernel]
	\index{cokernel}
Let $\category{C}$ be a category with zero object $0$ and let $f \colon A \longrightarrow B$. An object $C$ together with a 
morphism 
$c \colon B \longrightarrow C$ is called \emph{cokernel of} $f$ if it satisfies the following universal property: it holds $c \circ f 
= 0_{A,C}$ 
and for any morphism $c^\prime \colon B \longrightarrow C^\prime$ such that $c^\prime \circ f = 0_{A,C^\prime}$ there is a 
unique morphism 
$u \colon C \longrightarrow C^\prime$ such that $u \circ c = c^\prime$. Expressed as a diagram:
\begin{equation}
	\begin{tikzcd}
		{ }
		& B
		\arrow{d}{c}
		\arrow{ddr}{c^\prime}
		& { } \\
		A
		\arrow{ur}{f}
		\arrow{r}{0_{A,C}}
		\arrow{drr}[swap]{0_{A,C^\prime}}
		& C
		\arrow[dashed]{dr}{u}\\
		{ }
		& { }
		& C^\prime
	\end{tikzcd}.
\end{equation}
\end{definition}

The cokernel is unique up to isomorphism if it exists.
We will also write $\coker(f)$ for the cokernel morphism of $f$ and $\Coker(f)$ for the cokernel object of $f$.

\begin{corollary} \label{cor: cokernel is epi}
	Every cokernel is an epimorphism.
\end{corollary}

It can easily be seen that the cokernel of an epimorphism is the zero object together with the zero morphism and the cokernel of 
the zero morphism $0 \colon A \longrightarrow B$ is $B$ together with the identity morphism.

\section{Functors and Natural Transformations}
\label{sec: functors and natural transformations}

In this section we will take a step back and instead of investigating the relation of objects in a given category using morphisms, we 
want to study how we can relate categories using so-called functors.
As it turns out, in contrast to objects and morphisms in categories, there is even a way to relate morphisms between categories by  
natural transformations.

\begin{definition}[Functor]
	\label{def: functor}
	\index{functor}
A \emph{(covariant) functor} $\functor{F}$ from a category $\category{C}$ to a category $\category{D}$, written
$\functor{F} \colon \category{C} \longrightarrow  \category{D}$,
is a map sending each object 
$A \in \category{C}_0$ to an object
$\functor{F}A \in \category{D}$ and each morphism
$f \in \category{C}(B,A)$ to a morphism
$\functor{F}f \in \category{D}(\functor{F}B, \functor{F}A)$
such that
\begin{definitionlist}
	\item \label{def: functor_comp}
	$\functor{F}$ preserves composition, i.e.
	$\functor{F}(f \circ g) = \functor{F}f \circ \functor{F}g$,
	for any $f \in \category{C}(C,B)$ and $g \in \category{C}(B,A)$,
	\item \label{def: functor_id}
	$\functor{F}$ preserves identity morphisms, i.e.
	$\functor{F}\id_A = \id_{\functor{F}A}$, for each object
	$A \in \category{C}_0$.
\end{definitionlist}
A \emph{contravariant functor} is a functor $\functor{F}$, where instead of
$\functor{F}f \in \category{D}(\functor{F}B, \functor{F}A)$ we have
$\functor{F}f \in \category{D}(\functor{F}A, \functor{F}B)$ and instead of \ref{def: functor_comp} it holds
$\functor{F}(f \circ g) = \functor{F}g \circ \functor{F}f$.
\end{definition}

We will usually only use the term functor for covariant functors.

\begin{example}[Functors]
	\label{ex: functors}
Let $\category{C}$, $\category{D}$ be a categories.
\begin{examplelist}
	\item The map
	$\id_\category{C} \colon \category{C} \longrightarrow \category{C}$
	sending each object and each morphism to itself is the so-called \emph{identity functor on $\category{C}$}.
	\item For every object $A \in \category{C}_0$ there is a functor
	$\Id_A \colon \mathsf{1} \longrightarrow \category{C}$ by $\Id_A(*) = A$
	and $\Id_A(\id_*) = \id_A$.
	\item Fix an object $B \in \category{C}_0$.
	Then mapping each object $C \in \category{C}_0$ to the set $\Hom(B,C)$ and each morphism
	$f \colon X \longrightarrow Y$ to the map
	\begin{equation}
		\Hom(B,f) \colon \Hom(B,X) \longrightarrow \Hom(B,Y); \quad g \longmapsto f \circ g
	\end{equation}
	is a covariant functor $\Hom(B, \argument) \colon \category{C} \longrightarrow \Sets$.
	Similarly, mapping each object $A \in \category{C}_0$ to the set $\Hom(A,B)$ and each morphism
	$f \colon X \longrightarrow Y$ to the map
	\begin{equation}
		\Hom(f,B) \colon \Hom(Y,B) \longrightarrow \Hom(X,B); \quad g \longmapsto g \circ f
	\end{equation}
	is a contravariant functor $\Hom(\argument, B) \colon \category{C} \longrightarrow \Sets$.
	These functors are called \emph{Hom-functors}.
	\item The functor $\tau \colon \category{C} \times \category{D} \longrightarrow \category{D} \times \category{C}$
	given by $\tau(A,B) = (B,A)$ and $\tau(f,g) = (g,f)$ is a functor called \emph{flip}.
\end{examplelist}
\end{example}

We will call a functor
$\functor{F} \colon \category{C} \longrightarrow\category{D}$
\emph{faithful} if for any pair of objects $A,B \in \category{C}$ the map
$\functor{F} \colon \category{C}(B,A) \longrightarrow \category{D}(\functor{F}B, \functor{F}A)$
is injective, and $\functor{F}$ is called \emph{full} if it is surjective on hom-sets.
A functor is called \emph{fully faithful} if it is full and faithful.

Composing two functors by composing the maps on objects and morphisms yields again a functor.

\begin{proposition}[Composition of functors]
	Let
	$\functor{F} \colon \category{A} \longrightarrow\category{B}$
	and
	$\functor{G} \colon \category{B} \longrightarrow\category{C}$
	be functors between categories $\category{A}$, $\category{B}$, $\category{C}$.
	Mapping each object $A \in \category{A}$ to
	$\functor{G}\functor{F}A \in \category{C}$
	and each morphism $f \in \category{A}(B,A)$ to 
	$\functor{G}\functor{F}f \in \category{C}(\functor{G}\functor{F}B, \functor{G}\functor{F}A)$
	defines a functor $\functor{G} \circ \functor{F} \colon \category{A} \longrightarrow \category{C}$ called \emph{composition}.
\end{proposition}

Ignoring all issues that arise by taking categories of large categories,
this enables us to view functors as morphisms between categories.
We will denote the category of categories with functors as morphisms by $\Cat$.

Next we want to define the notion of subcategory.
This is not as straightforward as it first seems, and in general there does not seem to exist a universally accepted definition.
We will follow the idea that a subcategory $\category{B}$ of a category $\category{C}$ should be a subcollection
$\category{B}_0$ of the objects $\category{C}_0$, and for every pair
$B,B^\prime \in \category{B}$ a subset $\Hom_{\category{B}}(B,B^\prime)$ 
of the hom-set $\Hom_\category{C}(B,B^\prime)$.
This leads us to the following definition of embedding.

\begin{definition}[Embedding of Categories]
	\label{def: embedding of categories}
	A functor $F \colon \category{C} \longrightarrow \category{D}$ between categories $\category{C}$ and $\category{D}$ is 
	called an 
	\emph{embedding of categories}, if it is injective on objects and faithful.
\end{definition}

Obviously, a functor is an embedding in this sense if and only if it is injective on all morphisms.
Moreover, for every embedding
$\mathsf{F} \colon \category{C} \longrightarrow \category{D}$
the image $\mathsf{F}(\category{C})$ of $\category{C}$ in $\category{D}$ is isomorphic to $\category{C}$,
see \cite[p. 14]{maclane:1998a}.
This notion of embedding allows us to define subcategories.

\begin{definition}[Subcategory]
	\index{subcategory}
A category $\category{C}$ together with an embedding $\functor{I} \colon \category{C} \longrightarrow \category{D}$ is 
called a \emph{subcategory of} $\category{D}$.
\end{definition}

A subcategory $\category{C}$ of a category $\category{D}$ is called \emph{full} if the embedding
$\functor{I} \colon \category{C} \longrightarrow \category{D}$ is full.
It should be noted that this definition of subcategory does violate the principle of equivalence, since being injective on objects 
requires to identify objects.

On top of comparing categories by functors there is also a way to compare functors between the same categories.

\begin{definition}[Natural transformation]
	\label{def: natural transformation}
	\index{natural transformation}
Let $\category{C}$ and $\category{D}$ be categories
and $\functor{F}, \functor{G} \colon \category{C} \longrightarrow \category{D}$ be functors.
A \emph{natural transformation} $\eta$ from $\functor{F}$ to $\functor{G}$, written
$\eta \colon \functor{F} \Longrightarrow \functor{G}$,
is an assignment of a morphism
$\eta(A) \colon \functor{F}(A) \longrightarrow \mathsf{G}(A)$ in $\category{D}$  to every object $A \in \category{C}_0$,
such that for each morphism
$f \in \category{C}(B,A)$ 
the following diagram commutes
\begin{equation}
	\label{eq: def: natural transformation}
	\begin{tikzcd}[row sep=large, column sep = large]
		\mathsf{F}A
		\arrow{r}{\functor{F}f}
		\arrow{d}[swap]{\eta(A)}
		& \mathsf{F}B
		\arrow{d}{\eta(B)} \\
		\mathsf{G}A
		\arrow{r}{\functor{G}f}
		& \mathsf{G}B
	\end{tikzcd}.
\end{equation}
The morphisms $\eta(A)$ are called \emph{components} of $\eta$.
If all components $\eta(A) \colon \functor{F}A \longrightarrow \functor{F}B$ of a natural transformation are isomorphisms it is 
called a \emph{natural isomorphism}.
\end{definition}

Natural transformations can roughly be seen as a consistent choice of turning images under
$\mathsf{F}$ into images under $\mathsf{G}$.
We will depict a natural transformation $\eta \colon \functor{F} \Longrightarrow \functor{G}$ between functors
$\functor{F}, \functor{G} \colon \category{C} \longrightarrow \category{D}$ as
\begin{equation*}
	\begin{tikzcd}
		\category{C}
		\arrow[r, bend left=50,"\functor{F}", " "{name=F, below}]
		\arrow[r, bend right=50, "\functor{G}"{below}, " "{name=G}]
		& \category{D}
		\arrow[Rightarrow, "\eta", from=F, to=G	]
	\end{tikzcd}.
\end{equation*}

Using natural transformations we can define the notion of equivalence of categories.

\begin{definition}[Equivalence of categories]
	\index{equivalence!categories}
Let $\category{C}$ and $\category{D}$ be categories.
An \emph{equivalence of the categories $\category{C}$ and $\category{D}$} is a pair of functors
$\functor{F} \colon \category{C} \longrightarrow \category{D}$ and
$\functor{G} \colon \category{D} \longrightarrow \category{C}$ 
together with natural isomorphisms
$\eta \colon \functor{G} \circ \functor{F} \Longrightarrow \id_\category{C}$
and
$\varepsilon \colon \functor{F} \circ \functor{G} \Longrightarrow \id_\category{D}$.
We call $\category{C}$ and $\category{D}$ \emph{equivalent} if there exists an equivalence between them.
\end{definition}

Sometimes, one just states that $\functor{F} \colon \category{C} \longrightarrow \category{D}$ is an equivalence of
$\category{C}$ and $\category{D}$, implying the existence of a suitable functor
$\functor{G} \colon \category{D} \longrightarrow \category{C}$.
Equivalent categories share the same categorical properties.
Some first results are gathered in the next proposition.
But we will see later that equivalent categories share a lot more properties.

\begin{proposition}
	\label{prop: equivalence of categories}
	Let $\category{C}$ and $\category{D}$ be categories.
	Let furthermore $\functor{F} \colon \category{C} \longrightarrow \category{D}$,
	$\functor{G} \colon \category{D} \longrightarrow \category{C}$ be an equivalence of categories with natural isomorphisms
	$\eta \colon \functor{G} \circ \functor{F} \Longrightarrow \id_\category{C}$ and
	$\varepsilon \colon \functor{F} \circ \functor{G} \Longrightarrow \id_\category{D}$.
	\begin{propositionlist}
		\item \label{prop: equivalence of categories_1}
		The functors $\functor{F}$ and $\functor{G}$ are faithful and full.
		\item \label{prop: equivalence of categories_2}
		A morphism $f \colon B \longrightarrow C$ in $\category{C}$ is a monomorphism if and only if
		$\functor{F}f \colon \functor{F}B \longrightarrow \functor{F}C$ is a monomorphism in $\category{D}$.
		\item \label{prop: equivalence of categories_3}
		A morphism $f \colon A \longrightarrow B$ in $\category{C}$ is an epimorphism if and only if
		$\functor{F}f \colon \functor{F}A \longrightarrow \functor{F}B$ is an epimorphism in $\category{D}$.
	\end{propositionlist}
\end{proposition}

%
%

Again we can compose natural transformations.
But this time there are actually two different versions of composition.

\begin{proposition}[Vertical composition of natural transformations]
	Let $\eta \colon \functor{F} \Longrightarrow \functor{G}$ and
	$\mu \colon \functor{G} \Longrightarrow \functor{H}$ be natural transformations between functors
	$\functor{F}, \functor{G}, \functor{H} \colon \category{C} \longrightarrow \category{D}$.
	Their \emph{vertical composition}
	$\mu \circ \eta \colon \functor{F} \Longrightarrow \functor{H}$
	is a natural transformation given by morphisms
	\begin{equation}
		(\mu \circ \eta)(A) = \mu(A) \circ \eta(A) 
		\colon
		\functor{F}A
		\longrightarrow
		\functor{H}A
	\end{equation}
	for any $A \in \category{C}$.
\end{proposition}

\begin{proof}
	In the diagram
	\begin{equation*}
		\begin{tikzcd}
			\functor{F}A
			\arrow{r}{\functor{F}f}
			\arrow{d}[swap]{\eta(A)}
			& \functor{F}B
			\arrow{d}{\eta(B)} \\
			\functor{G}A
			\arrow{r}{\functor{G}f}
			\arrow{d}[swap]{\mu(A)}
			& \functor{G}B
			\arrow{d}{\mu(B)} \\
			\functor{H}A
			\arrow{r}{\functor{H}f}
			& \functor{H}B
		\end{tikzcd}
	\end{equation*}
	the upper and lower squares commute since $\eta$ and $\mu$ are natural transformations.
	Hence the big rectangle commutes, showing that $\mu \circ \eta$ is a natural transformation.
\end{proof}

The reason this is called vertical composition is that we can illustrate it as
\begin{equation*}
	\begin{tikzcd}
		\mathstrut \category{C}
		\arrow[r, bend left=50,"\functor{F}", " "{name=F, yshift=-.5ex}]
		\arrow[r, "\functor{G}"{near start,xshift=-.8ex}, " "{name=Gb,yshift=-.5ex}, " "{name=Gu,yshift=-.5ex}]
		\arrow[r, bend right=50, "\functor{H}"{below}, " "{name=H,yshift=-.5ex}]
		& \mathstrut \category{D}
		\arrow[Rightarrow, "\eta", from=F, to=Gu]
		\arrow[Rightarrow, "\mu", from=Gb, to=H]
	\end{tikzcd}
	\rightsquigarrow
	\begin{tikzcd}
		\category{C}
		\arrow[r, bend left=50,"\functor{F}", " "{name=F, below}]
		\arrow[r, bend right=50, "\functor{H}"{below}, " "{name=H}]
		& \category{D}
		\arrow[Rightarrow, "\mu \circ\eta", from=F, to=H]
	\end{tikzcd}.
\end{equation*}

\begin{proposition}[Horizontal composition of natural transformations]{\ \\}
	\label{prop: horizontal composition of natural transformations}
	Let $\eta \colon \functor{F}_1 \Longrightarrow \functor{G}_1$ and
	$\mu \colon \functor{F}_2 \Longrightarrow \functor{G}_2$
	be natural transformations between functors
	$\functor{F}_1, \functor{G}_1 \colon \category{A} \longrightarrow \category{B}$ and
	$\functor{F}_2, \functor{G}_2 \colon \category{B} \longrightarrow \category{C}$.
	Their \emph{horizontal composition}
	$\mu * \eta \colon (\functor{F}_2 \circ \functor{F}_1) \Longrightarrow (\functor{G}_2 \circ \functor{G}_1)$
	is a natural transformation given by the morphisms
	\begin{equation}
		\label{eq: prop: horizontal composition of natural transformations_1}
		(\mu * \eta) (A) = \mu(\functor{G}_1A) \circ \functor{F}_2\eta(A)
	\end{equation}
	for each $A \in \category{A}$.
\end{proposition}

\begin{proof}
	We need to show that the diagram
	\begin{equation*}
		\begin{tikzcd}
			(\functor{F}_2 \circ \functor{F}_1)A
			\arrow{d}[swap]{\mu(\functor{G}_1 A) \circ \functor{F}_2 \eta(A)}
			\arrow{r}{(\functor{F}_2 \circ \functor{F}_1)f}
			&(\functor{F}_2 \circ \functor{F}_1) B
			\arrow{d}{\mu(\functor{G}_1 B) \circ \functor{F}_2 \eta(B)} \\
			(\functor{G}_2 \circ \functor{G}_1)A
			\arrow{r}{(\functor{G}_2 \circ \functor{G}_1)f}
			& (\functor{G}_2 \circ \functor{G}_1)B
		\end{tikzcd}
	\end{equation*}
	commutes for all $f \colon A \longrightarrow B$.
	We get
	\begin{align*}
		( \functor{G}_2 \circ \functor{G}_1)f \circ \mu(\functor{G}_1A) \circ \functor{F}_2 \eta(A)
		&\overset{(a)}{=} \mu(\functor{G}_1 B) \circ \functor{F}_2\functor{G}_1f \circ \functor{F}_2 \eta(A) 
		= \mu(\functor{G}_1 B)  \circ \functor{F}_1(\functor{G}_1 f \circ \eta(A)) \\
		&\overset{(b)}{=} \mu(\functor{G}_1 B) \circ \functor{F}_2 (\eta(B) \circ \functor{F}_1 f) 
		= \mu(\functor{G}_1 B) \circ \functor{F}_2 \eta(B) \circ \functor{F}_2 \functor{F}_1 f,
	\end{align*}
	where we used in $(a)$ the diagram \eqref{eq: def: natural transformation} for the natural transformation $\mu$ and 
	in $(b)$ for the natural transformation $\eta$.
\end{proof}

The horizontal composition can be visualized as
\begin{equation*}
	\begin{tikzcd}
		\category{A}
		\arrow[r, bend left=50,"\functor{F}_1", " "{name=F1, below}]
		\arrow[r, bend right=50, "\functor{G}_1"{below}, " "{name=G1}]
		& \category{B}
		\arrow[r, bend left=50,"\functor{F}_2", " "{name=F2, below}]
		\arrow[r, bend right=50, "\functor{G}_2"{below}, " "{name=G2}]
		\arrow[Rightarrow, "\eta", from=F1, to=G1]
		& \category{C}
		\arrow[Rightarrow, "\mu", from=F2, to=G2]
	\end{tikzcd}
	\rightsquigarrow
	\begin{tikzcd}
		\category{A}
		\arrow[r, bend left=50,"\functor{F}_2 \circ \functor{F}_1", " "{name=F, below}]
		\arrow[r, bend right=50, "\functor{G}_2 \circ \functor{G}_1"{below}, " "{name=G}]
		& \category{C}
		\arrow[Rightarrow, "\mu * \eta", from=F, to=G]
	\end{tikzcd}.
\end{equation*}

\begin{remark}
Since $\mu$ is a natural transformation the definition of $\mu * \eta$ by
\eqref{eq: prop: horizontal composition of natural transformations_1}
is equivalent to $(\mu * \eta) (A) = \functor{G}_2 \eta(A) \circ \mu(\functor{F}_1A)$.
\end{remark}

\begin{definition}[Adjoint functors]
	\label{def:AdjointFunctors}
	\index{adjoint functors}
Let $\category{C}$ and $\category{D}$ be categories.
An \emph{adjunction} between these categories consists of functors
$\functor{F} \colon \category{C} \to \category{D}$
and $\functor{G} \colon \category{D} \to \category{C}$
as well as natural transformations
$\varepsilon \colon \functor{F} \circ \functor{G} \Longrightarrow \id_\category{D}$
and
$\eta \colon \id_\category{D} \Longrightarrow \functor{G} \circ \functor{F}$
such that for each $A \in \category{C}$ and each $B \in \category{C}$
\begin{align}
	\id_{\functor{F}A} &= \varepsilon(\functor{F}A) \circ \functor{F}(\eta(A))
	\shortintertext{and}
	\id_{\functor{G}Y} &= \functor{G}(\varepsilon(Y)) \circ \eta(\functor{G}(Y))
\end{align} 
holds.
We call $\functor{F}$ \emph{left adjoint} to $\functor{G}$, and reversely 
$\functor{G}$ \emph{right adjoint} to $\functor{F}$.
\end{definition}

We sometimes write $\functor{F} \dashv \functor{G}$ if $\functor{F}$ is left adjoint to $\functor{G}$.
In some contexts $\ev$ is used instead of $\varepsilon$
and $\coev$ is used instead of $\eta$.

\section{Limits and Colimits}
\label{sec: limits and colimits}
Many of the standard notions in category theory we have seen so far -- like initial/terminal object, kernel/cokernel, but also 
equalizer, pullbacks etc. -- are special cases of a more general notion, so-called limits and colimits.
For this we first have to give a precise definition of a diagram in a category.

\begin{definition}[Diagram]
	\index{diagram}
Let $\category{C}$ be a category and $\category{I}$ a small category.
A functor $\functor{D} \colon \category{I} \longrightarrow \category{C}$ is called a \emph{diagram of shape $\category{I}$}.
\end{definition}

For a diagram $\functor{D} \colon \category{I} \longrightarrow \category{C}$ we will often write $D_I$ instead of $D(I)$ for $I \in 
\category{I}$ to indicate that one should think of a diagram as an indexed class of objects and morphisms.
Before defining limits and colimits for diagrams we introduce the general notions of sources and sinks.

\begin{definition}[Sources and sinks]
Let $\category{C}$ be a category.
\begin{definitionlist}
	\item A \emph{source} is an object $C \in \category{C}_0$ together with a family of morphisms
	$(f_i \colon C \longrightarrow C_i)_{i \in I}$ indexed by some class $I$.
	\item A \emph{sink} is an object $C \in \category{C}_0$ together with a family of morphisms
	$(f_i \colon C_i \longrightarrow C)_{i \in I}$ indexed by some class $I$.
\end{definitionlist}
\end{definition}

We will also denote a source simply by its family of morphisms, since then the corresponding object is clear.
Given a diagram $\functor{D} \colon \category{I} \longrightarrow \category{C}$ and a source
$\left(f_I \colon C \longrightarrow D_I\right)_{I \in \category{I}}$,
we say that this source is a \emph{source of the diagram $\functor{D}$} if for all
$u \colon I \longrightarrow J$ in $\category{I}$
the triangle
\begin{equation}
	\begin{tikzcd}
		\functor{D}_I
		\arrow{rr}{\functor{D}u}
		& { }
		& \functor{D}_J  \\
		{ }
		& C 
		\arrow{ul}{f_I}
		\arrow{ur}[swap]{f_J}
		& { }
	\end{tikzcd}
\end{equation}
commutes.
We will denote such a source of a diagram by the pair $(C, f_I)_{I \in \category{I}}$ since the domain and codomain of each $f_I$ 
are clear from the diagram $\functor{D}$.
Note that sources of a diagram are also often called \emph{cones} in the literature.
Dually we can define \emph{sinks of a diagram}.
With this we can define for any diagram in a category the limit and colimit of it.

\begin{definition}[Limit]
	\label{def: limit}
	\index{limit}
Let $\functor{D} \colon \category{I} \longrightarrow \category{C}$ be a diagram in $\category{C}$ of shape $\category{I}$.
A \emph{limit} of $\functor{D}$ is a source
$(L , \ell_I)_{I \in \category{I}}$
of the diagram $\functor{D}$ with the following universal property:
For any other source $(C,f_I )_{I \in \category{I}}$ of the diagram $\functor{D}$ there exists a unique morphism
$f \colon C \longrightarrow L$ with $\ell_I \circ f = f_I$ for all $I \in \category{I}$.
Put diagrammatically:
\begin{equation}
	\begin{tikzcd}
		\functor{D}_I
		\arrow{rr}{\functor{D}u}
		& { }
		& \functor{D}_J  \\
		{ }
		& L
		\arrow{ul}[swap]{\ell_I}
		\arrow{ur}{\ell_J}
		& { } \\
		{ }
		& C 
		\arrow[dashrightarrow]{u}{f}
		\arrow{uul}{f_I}
		\arrow{uur}[swap]{f_J}
		& { } 
	\end{tikzcd}
\end{equation}
commutes.
\end{definition}

\index{colimit}
The dual notion of a limit uses sinks of a diagram instead and is called colimit.
Limits and colimits of diagrams need not exist, but if they do they are unique up to unique isomorphisms as usual.

Limits appear everywhere in category theory.
To some degree one could even say that category theory is the study of limits.
To see why this is the case we list some limits we encounter during this thesis.

\begin{example}[Limits and colimits]
	\label{ex:limitsandcolimits}
Let $\category{C}$ be a category.
\begin{definitionlist}
	\item \index{initial object}
	An initial object is the limit of the empty diagram.
	\index{terminal object}
	Dually, a terminal object is the colimit of the empty diagram.
	\item \index{product}
	The limit of a diagram $\functor{D}$ given by a discrete category $\category{I}$ is the product of the objects
	$\functor{D}_I \in \category{C}$.
	\index{coproduct}
	Dually, the colimit of such a discrete diagram is the coproduct.	
	\item \index{kernel}
	If $\category{C}$ has a zero object $0$, the kernel of a morphism $f \colon A \longrightarrow B$ is the limit of the 
	diagram
	\begin{equation}
		\label{diag: kernel as limit}
		\begin{tikzcd}
			A
			\arrow[r,shift left, "f"]
			\arrow[r, shift right,  "0"{below}]
			&B
		\end{tikzcd}.
	\end{equation}
	\index{cokernel}
	Dually, the cokernel of $f$ is the colimit of \eqref{diag: kernel as limit}.		
	\item \index{equalizer}
	Generalizing the last example we call the limit of a diagram
	\begin{equation}
		\label{diag: equalizer}
		\begin{tikzcd}
			A
			\arrow[r,shift left, "f"]
			\arrow[r, shift right,  "g"{below}]
			&B
		\end{tikzcd}
	\end{equation}
	the \emph{equalizer} of $f$ and $g$.
	\index{coequalizer}
	Dually, the colimit of \eqref{diag: equalizer} is called \emph{coequalizer} of $f$ and $g$.
	\item \index{pullback}
	The limit of a diagram of the form
	\begin{equation}
		\label{diag: pullback}
		\begin{tikzcd}
			{ }
			&B
			\arrow{d}{g}\\
			A
			\arrow{r}{f}
			&C
		\end{tikzcd}
	\end{equation}
	is called \emph{pullback} of $f$ and $g$ and denoted by $f \times_C g$.
	\index{pushout}
	Dually, the colimit of \eqref{diag: pullback} is called \emph{pushout}.
\end{definitionlist}
\end{example}

\begin{proposition}
Let $\category{C}$ be a category.
\begin{propositionlist}
	\item Every equalizer is a monomorphism.
	\item Every coequalizer is an epimorphism.
\end{propositionlist}
\end{proposition}

The reverse implication need not hold in general.

\begin{definition}[Regular epi- and monomorphisms]
	\label{def:RegularEpiMono}
Let $\category{C}$ be a category.
\begin{definitionlist}
	\item \index{monomorphism!regular}
	A morphism $f \colon A \to B$ is called \emph{regular monomorphism}
	if it is the equalizer of some pair of morphisms.
	\item \index{epimorphism!regular}
	A morphism $f \colon A \to B$ is called \emph{regular epimorphism}
	if it is the coequalizer of some pair of morphisms.
\end{definitionlist}	
\end{definition}

We have already seen that initial objects, kernels and in general limits need not exist.
Thus it is of interest if we can transport existing limits to other categories via functors.
\filbreak
\begin{definition}[Preservation of limits]
	Let $\functor{F} \colon \category{C} \longrightarrow \category{D}$ be a functor and
	$\functor{D} \colon \category{I} \longrightarrow \category{C}$
	a diagram.
	\begin{definitionlist}
		\item The functor $\functor{F}$ is said to \emph{preserve a limit}
		$\left(f_I \colon L \longrightarrow \functor{D}_I\right)_{I \in \category{I}}$
		if
		$(\functor{F}f_I \colon \functor{F}L \longrightarrow \functor{F}D_I)_{I \in \category{I}}$
		is a limit of the diagram
		$\functor{F} \circ \functor{D} \colon \category{I} \longrightarrow \category{D}$.
		\item The functor $\functor{F}$ is said to \emph{preserve limits} if it preserves limits of all shapes.
	\end{definitionlist}
	
\end{definition}

An important example of limit preserving functors is an equivalence of categories.
Although this is not surprising, since equivalent categories should have the same categorical properties, the proof actually 
needs knowledge about adjunctions.
Thus we simply state the theorem here and refer to \cite[Chap. 9]{hott:2013} for a proof.

\begin{proposition}[Equivalences preserve limits]
	\label{prop: equivalences preserve limits}
	\index{equivalence!categories}
Every equivalence $\functor{F} \colon \category{C} \longrightarrow \category{D}$ of catego\-ries preserves limits.
\end{proposition}
	
\section{Monoids and Modules}
\label{sec:MonoidsModules}

We collect basic definitions and constructions of monoids and their modules
internal to a given monoidal category, see e.g. \cite{kashiwara.schapira:2006a}
for more about monoidal categories.

\begin{definition}[Monoidal category] 
	\label{def: MonoidalCategory}
	\index{monoidal category}
A \emph{monoidal category} is a category $\category{C}$ equipped with the following data:
\begin{definitionlist}
	\item A functor 
	\begin{equation}
		\tensor \colon \category{C} \times \category{C} \to \category{C}
	\end{equation} called \emph{tensor product}.
	\item An object $\Unit \in \category{C}_0$ called \emph{unit}.
	\item A natural isomorphism 
	\begin{equation}
		\label{eq: def: Monoidal Category_associativity}
		\asso \colon \tensor \circ \operatorname{(\tensor \times \id)} \Longrightarrow \tensor \circ (\id \times \tensor)
	\end{equation} 
	called \emph{associativity}.
	Diagrammatically:
	\begin{equation}
		\begin{tikzcd}
			\category{C} \times \category{C} \times \category{C}
			\arrow{d}[swap]{\tensor \times \id}
			\arrow{r}{\id \times \tensor}
			& \category{C} \times \category{C}
			\arrow{d}{\tensor} \\
			\category{C} \times \category{C}
			\arrow{r}{\tensor}
			\arrow[Rightarrow]{ur}{\asso}		
			& \category{C}					
		\end{tikzcd}
	\end{equation}
	\item Two natural isomorphisms
	\begin{equation}
		\label{eq: def: Monoidal Category_left identity}
		\lidentity \colon \tensor \circ \operatorname{(\Id_\Unit \times \id )} \Longrightarrow \id
	\end{equation}
	called \emph{left identity} and
	\begin{equation}
		\label{eq: def: Monoidal Category_right identity}
		\ridentity \colon \tensor \circ \operatorname{( \id \times \Id_\Unit)} \Longrightarrow \id
	\end{equation}
	called \emph{right identity}.
	Diagrammatically:
	\begin{equation}
		\begin{tikzcd}
			\mathsf{1} \times \category{C}
			\arrow{d}[swap]{\Id_\Unit \times \id}
			\arrow[Rightarrow]{dr}{\lidentity}
			& \category{C}
			\arrow{l}[swap]{\simeq}
			\arrow{r}{\simeq}
			\arrow{d}{\id}
			& \category{C} \times \mathsf{1} 
			\arrow{d}{\id \times \Id_\Unit}
			\arrow[Rightarrow]{dl}[swap]{\ridentity} \\
			\category{C} \times \category{C}
			\arrow{r}{\tensor}
			& \category{C}
			& \category{C}
			\arrow{l}[swap]{\tensor}
		\end{tikzcd}
	\end{equation}
\end{definitionlist}
These data are required to fulfil the following coherence conditions:
\begin{definitionlist}[resume]
	\item \emph{Associativity coherence:}
	the diagram
	\begin{equation}
		\begin{tikzcd}
			((A \tensor B) \tensor C) \tensor D
			\arrow{rr}{\asso(A,B,C) \tensor \id_D}
			\arrow{d}[swap]{\asso(A \tensor B, C ,D)}
			& { }
			& (A \tensor (B \tensor C)) \tensor D 	
			\arrow{d}{\asso(A, B \tensor C , D)} \\
			(A \tensor B) \tensor (C \tensor D)
			\arrow{dr}[swap]{\asso(A,B,C \tensor D)}
			& { }
			& A \tensor ((B \tensor C) \tensor D)
			\arrow{dl}{\id_A \tensor \asso(B,C,D)} \\
			{ }
			& (A \tensor (B \tensor (C \tensor D)))
			& { }
		\end{tikzcd}
	\end{equation}
	commutes for all objects $A,B,C,D \in \category{C}_0$.
	\item \emph{Identity coherence:}
	the diagram
	\begin{equation}
		\begin{tikzcd}
			(A \tensor \Unit) \tensor B	
			\arrow{rr}{\asso(A,\Unit, B)}
			\arrow{dr}[swap]{\ridentity(A) \tensor \id_B}
			& { }
			& A \tensor (\Unit \tensor B) 
			\arrow{dl}{\id_A \tensor \lidentity(B)} \\
			{ }
			& A \tensor B
			& { }							
		\end{tikzcd}
	\end{equation}
	commutes for all objects $A,B \in \category{C}_0$.
\end{definitionlist}
\end{definition}

A monoidal category is called \emph{strict} if the associativity, as well as the left and right identity are not mere isomorphisms but strict equalities.

\begin{example}\
	\begin{examplelist}
		\item The category $\Sets$ of sets together with the cartesian product and any one-point set as unit is a monoidal 
		category.
		\item The category $\Ab$ of abelian groups together with the tensor product of groups and the group of integers 
		$\mathbb{Z}$ as unit is a monoidal category.
		\item The category $\Bimodules(\ring{R},\ring{R})$ of bimodules over some ring $\ring{R}$ together with the
		tensor product of bimodules and	$\ring{R}$ as unit is a monoidal category.
	\end{examplelist}			
\end{example}

The commutativity of a monoidal category is captured by the following notion.
Here $\tau$ denotes the flip functor $\tau \colon \category{C} \times \category{C} \to \category{C} \times \category{C}$,
$\tau(A,b) = (B,A)$.

\begin{definition}[Symmetric monoidal category]
	\label{def: symmetric monoidal category}
	\index{monoidal category!symmetric}
	A monoidal category $\category{C}$ together with a\linebreak
	natural isomorphism
	\begin{equation}
		\mathsf{B} \colon \tensor \Longrightarrow \tensor \circ \tau
	\end{equation}
	such that
	\begin{equation}
		\mathsf{B} \circ \mathsf{B} = \id
	\end{equation}
	holds,
	is called \emph{symmetric} if the diagram
	\begin{equation}
		\begin{tikzcd}[column sep = small]
			{ }
			& (A \tensor B ) \tensor C
			\arrow{rr}{\asso(A,B,C)}
			\arrow{dl}[swap]{\mathsf{B}(A,B) \tensor \id}
			& {\qquad}
			& A \tensor (B \tensor C)
			\arrow{dr}{\mathsf{B}(A,B \tensor C)}
			& { } \\
			(B \tensor A) \tensor C
			\arrow{dr}[swap]{\asso(B,A,C)}
			& { }
			& { }
			& { }
			& (B \tensor C) \tensor A
			\arrow{dl}{\asso(B,C,A)} \\
			{ }
			& B \tensor (A \tensor C)
			\arrow{rr}{\id \tensor \mathsf{B}(A,C)}
			& { }
			& B \tensor (C \tensor A)
			& { }
		\end{tikzcd}
	\end{equation}
	commutes for all $A,B,C \in \category{C}_0$.
	The natural isomorphism $\mathsf{B}$ is called \emph{symmetric braiding}.
	\index{symmetric braiding}
\end{definition}

\begin{example}\
	\begin{examplelist}
		\item The monoidal category $\Sets$ with the cartesian product is symmetric.
		\item The monoidal category $\Ab$ with the tensor product of groups is symmetric.
		\item For every commutative ring $\ring{R}$ the monoidal category $\Bimodules_\sym(\ring{R},\ring{R})$ of symmetric bimodules is symmetric.
	\end{examplelist}
\end{example}

\begin{definition}[Lax Monoidal Functor]
	\label{def:LaxMonoidalFunctor}
Let $\category{C}$ and $\category{D}$ be monoidal categories.
A \emph{lax monoid\-al functor} is a functor
$\functor{F} \colon \category{C} \to \category{D}$
together with the following data:
\begin{definitionlist}
	\item A morphism
	$\varepsilon \colon \Unit_\category{D} \to	\functor{F}(\Unit_\category{C})$.
	\item A natural transformation 
	$\eta \colon \tensor[\category{D}] \circ \operatorname{\functor{F} \times \functor{F}}
	\Longrightarrow \functor{F} \circ \tensor[\category{C}]$.
\end{definitionlist}
These data are required to make the following diagrams commute for all
$A,B,C \in \category{C}$:
\begin{definitionlist}
	\item \emph{Associativity}:
	\begin{equation}
	\begin{tikzcd}[column sep = 3.5cm]
		\big(\functor{F}(A) \tensor[\category{D}] \functor{F}(B) \big) \tensor[\category{D}] \functor{F}(C)
			\arrow[d,"{\eta(A,B) \tensor \id}"{swap}] 
			\arrow[r,"{\asso_{\category{D}}(\functor{F}(A),\functor{F}(B),\functor{F}(C))}"]
		& \functor{F}(A) \tensor[\category{D}] \big(\functor{F}(B) \tensor[\category{D}] \functor{F}(C)\big)
			\arrow[d,"{\id \tensor \eta(B,C)}"] \\
		\functor{F}(A \tensor[\category{C}] B) \tensor[\category{D}] \functor{F}(C)
			\arrow[d,"{\eta(A \tensor[\category{C}] B,C)}"{swap}]
		& \functor{F}(A) \tensor[\category{D}] \functor{F}(B \tensor[\category{C}] C)
			\arrow[d,"{\eta(A,B \tensor[\category{C}] C)}"] \\
		\functor{F}\big((A \tensor[\category{C}] B) \tensor[\category{C}] C\big)
			\arrow[r,"{\functor{F}(\asso_\category{C}(A,B,C))}"] 
		&\functor{F}\big(A \tensor[\category{C}] (B \tensor[\category{C}] C)\big)
	\end{tikzcd}
	\end{equation}
	\item \emph{Unitality}:
	\begin{equation}
	\begin{tikzcd}
		\Unit_\category{D} \tensor[\category{D}] \functor{F}(A)
			\arrow[r,"\varepsilon \tensor \id"]
			\arrow[d,"\lidentity_\category{D}(\functor{F}(A))"{swap}]
		& \functor{F}(\Unit_\category{C}) \tensor[\category{D}] \functor{F}(A)
			\arrow[d,"{\eta(\Unit_\category{C},A)}"] \\
		\functor{F}(A)
		&\functor{F}(\Unit_\category{C} \tensor[\category{C}] A)
			\arrow[l, "\functor{F}(\lidentity_\category{C}(A))"]
	\end{tikzcd}
	\end{equation}
	and
	\begin{equation}
	\begin{tikzcd}
		\functor{F}(A)  \tensor[\category{D}]  \Unit_\category{D}
			\arrow[r,"\id \tensor \varepsilon"]
			\arrow[d,"\ridentity_\category{D}(\functor{F}(A))"{swap}]
		& \functor{F}(A) \tensor[\category{D}] \functor{F}(\Unit_\category{C}) 
			\arrow[d,"{\eta(A,\Unit_\category{C})}"] \\
		\functor{F}(A)
		&\functor{F}(A \tensor[\category{C}] \Unit_\category{C})
			\arrow[l, "\functor{F}(\ridentity_\category{C}(A))"]
	\end{tikzcd}
	\end{equation}	
\end{definitionlist}
\end{definition}

A lax monoidal functor with invertible $\varepsilon$ and $\eta$ is called monoidal,
while an oplax monoidal functor is a lax monoidal functor between the opposite categories.

By the microcosm principle \cite{baez.dolan:1998a}
monoidal categories are the correct categorical setting to define monoids.

\begin{definition}[Monoid]
	\label{def:Monoid}
	\index{monoid object}
Let $\category{C}$ be a monoidal category.
A \emph{monoid object} (or simply \emph{monoid}) is an object $A \in \category{C}_0$ equipped with a morphism
\begin{equation}
	\mu \colon A \tensor A \longrightarrow A
\end{equation}
called \emph{multiplication} and a morphism
\begin{equation}
	\eta \colon \Unit \longrightarrow A
\end{equation}
called \emph{unit} such that
\begin{equation}
	\begin{tikzcd}
		(A \tensor A) \tensor A 
		\arrow{rr}{\asso}
		\arrow{d}[swap]{\mu \tensor \id_A}
		& {}
		& A \tensor (A \tensor A)
		\arrow{d}{\id_A \tensor \mu} \\
		A \tensor A
		\arrow[swap]{dr}{\mu}
		& {}
		& A \tensor A
		\arrow{dl}{\mu} \\
		{ } & A &
	\end{tikzcd}
\end{equation}
and
\begin{equation}
	\begin{tikzcd}
		\Unit \tensor A
		\arrow{r}{\eta \tensor \id_A}
		\arrow{dr}[swap]{\lidentity}
		& A \tensor A
		\arrow{d}{\mu}
		& A \tensor \Unit
		\arrow{l}[swap]{\id_A \tensor \eta}
		\arrow{dl}{\ridentity}\\
		{ }
		& A
		& { }
	\end{tikzcd}
\end{equation}
commute.
\end{definition}

This definition is just the usual definition of monoids written in terms of the maps that are involved, instead of in terms of elements.
In a symmetric monoidal category we can also define commutative monoids.

\begin{definition}[Commutative monoid]
	A monoid $A$ in a symmetric monoidal category $\category{C}$ is called \emph{commutative} if the diagram
	\begin{equation}
		\begin{tikzcd}
			A \tensor A
			\arrow{rr}{B}
			\arrow{dr}[swap]{\mu}
			&{ }
			& A \tensor A
			\arrow{dl}{\mu} \\
			{ }
			& A
			&{ }
		\end{tikzcd}
	\end{equation}
	commutes.
	Here $B$ denotes the symmetric braiding of $\category{C}$ and $\mu$ denotes the multiplication of $A$.
\end{definition}

A morphism of monoids can then be phrased as follows.

\begin{definition}[Morphism of monoids]
	\label{def:MorphismOfMonoids}
	Let $\category{C}$ be a monoidal category and let $A$ and $A^\prime$ be monoids with multiplications $\mu$, 
	$\mu^\prime$ and units $\eta$, $\eta^\prime$, respectively.
	A morphism
	$f \colon A \longrightarrow A^\prime$ is a \emph{morphism of monoids} if
	\begin{equation}
		\label{eq: def: morphism of monoids_multiplicativity}
		\begin{tikzcd}
			A \tensor A
				\arrow[r, "f \tensor f"]
				\arrow[d, "\mu"{swap}]
			& A' \tensor A'
				\arrow[d,"\mu'"] \\
			A
				\arrow[r,"f"]
			&A'
		\end{tikzcd}
	\end{equation}
	and
	\begin{equation}
		\label{eq: def: morphism of monoids_unit}
		\begin{tikzcd}
			\Unit
				\arrow[d,"\eta"{swap}]
				\arrow[dr,"\eta'"]
			& {} \\			
			A
				\arrow[r,"f"]
			&A'
		\end{tikzcd}
	\end{equation}
	commute.
\end{definition}

Written in elements \eqref{eq: def: morphism of monoids_multiplicativity} is just the compatibility with multiplication and
\eqref{eq: def: morphism of monoids_unit} is the preservation of the unit.

\begin{corollary}[Category of monoids]
	\label{cor: category of monoids}
	The monoids of a monoidal category $\category{C}$ together with morphisms of monoids as morphisms form a category, 
	called \emph{category of monoids of} $\category{C}$ and denoted by $\Monoids(\category{C})$.
\end{corollary}

\begin{proof}
	For any monoid $A$ the identity morphism $\id_A$ is obviously a morphism of monoids.
	Let $f \colon A \longrightarrow A^\prime$ and $g \colon A^\prime \longrightarrow A^{\prime\prime}$
	be two morphisms of monoids.
	Then
	\begin{equation*}
		(g \circ f) \circ \mu 
		=g \circ \mu^\prime \circ  (f \tensor f)
		= \mu^{\prime\prime} \circ (g \tensor g) \circ (f \tensor f)
		= \mu^{\prime\prime} \circ \left( (g \circ f) \tensor (g \circ f) \right)
	\end{equation*}
	and
	\begin{equation*}
		(g \circ f) \circ \eta = g \circ \eta^\prime = \eta^{\prime\prime}
	\end{equation*}
	show that $g \circ f \colon A \longrightarrow A^{\prime\prime}$ is a morphism of monoids.
\end{proof}

Any lax monoidal functor $\functor{F} \colon \category{C} \to \category{D}$
induces a functor $\functor{F} \colon \Monoids(\category{C}) \to \Monoids(\category{D})$.
For a symmetric monoidal category $\category{C}$ the full subcategory of $\Monoids(\category{C})$ consisting of commutative 
monoids is denoted by
$\Monoids_\commutative(\category{C})$.

\begin{example}\
	\begin{examplelist}
		\item In $\Sets$ (commutative) monoid objects are usual (commutative) monoids and morphisms of such monoid
		objects are the	usual monoid homomorphisms.
		\item In $\Ab$ (commutative) monoid objects are unital (commutative) rings and morphisms of monoids are
		ring morphisms.
		Hence $\Monoids(\Ab)$ is $\Rings$.
		\item In $\Bimodules(\ring{R}, \ring{R})$ (commutative) monoid objects are unital (commutative) associative algebras over
		the (commutative) ring $\ring{R}$ and morphisms of monoids are unital algebra homomorphisms.
		Hence $\Monoids(\Bimodules(\ring{R},\ring{R}))$ is $\Algebras_{\ring{R}}$.
	\end{examplelist}
\end{example}

Thinking of monoids in a monoidal category as rings or algebras suggests how to proceed from here.
We can define now modules over monoids by a categorical version of the usual definition.

\begin{definition}[Right module over a monoid]
	\label{def:RightModuleOverMonoid}
	\index{right module}
Let $\category{C}$ be a monoidal category and let $A \in \Monoids(\category{C})$ be a monoid.
A \emph{right module over} $A$ is an object $M \in \category{C}$ equipped with a morphism
\begin{equation}
	\rho \colon M \tensor A \longrightarrow M,
\end{equation}
such that the diagrams
\begin{equation}
	\label{diag:ModuleMorphism_1}
	\begin{tikzcd}
		M \tensor A \tensor A
		\arrow{r}{\rho \tensor \id_A}
		\arrow{d}[swap]{\id_M \tensor \mu}
		& M \tensor A
		\arrow{d}{\rho} \\
		M \tensor A
		\arrow{r}{\rho}
		& M
	\end{tikzcd}
\end{equation}
and
\begin{equation}
	\begin{tikzcd}
		M \tensor \Unit
		\arrow{rr}{\id_M \tensor \eta}
		\arrow{dr}[swap]{\ridentity}
		& { }
		& M \tensor A
		\arrow{dl}{\rho}\\
		{ }
		& M
		&{ }
	\end{tikzcd}
\end{equation}
commute.
\end{definition}

A \emph{left module} over a monoid $A$ is then defined analogously.
Note that we implicitly used the associativity isomorphism in \eqref{diag:ModuleMorphism_1}.
Requiring an additional compatibility between given left and right module structures gives a bimodule.

\begin{definition}[Bimodule over monoids]
	\label{def: bimodule over monoids}
	\index{bimonoid}
Let $\category{C}$ be a monoidal category and let
$A,B \in \Monoids(\category{C})$ be monoids.
A \emph{$(B,A)$-bimodule} is an object $M \in \category{C}$ together with morphisms
\begin{equation}
	\lambda_B \colon B \tensor M \longrightarrow M
\end{equation}
and 
\begin{equation}
	\rho_A \colon M \tensor A \longrightarrow M,
\end{equation}
such that $M$ is a left $B$-module with respect to $\lambda_B$ and a right $A$-module with respect to
$\rho_A$ and
\begin{equation}
	\begin{tikzcd}
		B \tensor M \tensor A
		\arrow{r}{\lambda_B \tensor \id_A}
		\arrow{d}[swap]{\id_B \tensor \rho_A}
		& M \tensor A
		\arrow{d}{\rho_A}\\
		B \tensor M 
		\arrow{r}{\lambda_B}
		& M
	\end{tikzcd}
\end{equation}
commutes.
\end{definition}

Note that any monoid $A \in \Monoids(\category{C})$ can be seen as a $(A,A)$-bimodule by taking as left and right actions
the multiplication of the monoid.
As before also the notion of morphisms can be transferred to the categorical case without any problems.

\begin{definition}[Morphism of right modules]
	\label{def: morphism of modules}
Let $\category{C}$ be a monoidal category.
Moreover, let $A \in \Monoids(\category{C})$ be a monoid and let
$(M, \rho_M)$, $(M^\prime, \rho_{M^\prime})$
be right $A$-modules.
A morphism $f \colon M \longrightarrow M^\prime$ is called \emph{morphism of right modules} if 
\begin{equation}
	\begin{tikzcd}
		M \tensor A
		\arrow{r}{f \tensor \id_A}
		\arrow{d}[swap]{\rho_M}
		& M^\prime \tensor A
		\arrow{d}{\rho_{M^\prime}} \\
		M
		\arrow{r}{f}
		& M^\prime
	\end{tikzcd}
\end{equation}
commutes.
\end{definition}

A morphism of left modules is defined analogously and a morphism of bimodules is simply a morphism that respects both 
the left and right module structure.

\begin{corollary}[Categories of modules]
	\label{cor: categories of modules}
	Let $\category{C}$ be a monoidal category and consider mo\-noids $A,B \in \Monoids(\category{C})$.
	Left $B$-modules, right $A$-modules and $(B,A)$-bimodules together with the appropriate notion of morphism
	form categories, denoted by $\LeftModules{A}_{\category{C}}$, $\RightModules[\category{C}]{B}$ and 
	$\Bimodules_\category{C}(B,A)$, respectively.
\end{corollary}

\begin{proof}
	We only prove the statement for right modules, the other cases can be done similarly.
	First, it is clear that for any right $A$-module $(M,\rho_M)$ the identity morphism $\id_M$ 
	in $\category{C}$ is a module morphism.
	Moreover, for right $A$-module morphisms
	$f \colon (M,\rho_M) \longrightarrow (M',\rho_{M'})$
	and
	$g \colon (M', \rho_{M'}) \longrightarrow (M'', \rho_{M''})$
	we have
	\begin{align*}
		\rho_{M''} \circ (g \circ f ) \tensor \id_A
		= \rho_{M''} \circ (g \tensor \id_A) \circ (f \tensor \id_A)
		= g \circ \rho_{M'} \circ (f \tensor \id_A)
		= (g \circ f) \circ \rho_M,
	\end{align*}
	and thus $g \circ f$ is a morphism of right $A$-modules.
	Finally, from the associativity of the composition in $\category{C}$ follows directly the associativity of 
	the composition of module morphisms.
\end{proof}

These categories of modules indeed reproduce the various notions of modules we know from algebra.

\begin{example}\
	\begin{examplelist}
		\item In $\Sets$ a module over a monoid is just an action of a monoid on a set.
		\item In $\Ab$ a module over a monoid is a module over a ring in the usual sense of algebra.
		\item In $\Bimodules(\ring{R})$ a module over a monoid is a module over an algebra over $\ring{R}$.
	\end{examplelist}
\end{example}

The main feature of the tensor product of modules over a given ring (or algebra) is that we are able to either let a ring element 
act from the right on the left component of the tensor product or from the left on the right component of the tensor product.
This suggests to define the tensor product as the coequalizer of these two actions.
But to define the tensor product this way we need an additional requirement for the monoidal category we start with.

\begin{definition}[Tensor product of modules]
	\label{def: tensor product of modules}
	\index{tensor product!module objects}
Let $\category{C}$ be a monoidal category with coequalizers and let $B \in \Monoids(\category{C})$ be a monoid.
For any left $B$-module
$M \in \LeftModules{B}_{\category{C}}$
and right $B$-module
$N \in \RightModules[\category{C}]{B}$
the coequalizer of the left and right actions
\begin{equation}
	\begin{tikzcd}
		N \tensor B \tensor M 
		\arrow[r, shift left,"\rho_N \tensor \id_M"]
		\arrow[r, shift right, "\id_N \tensor \lambda_M", swap]
		& N \tensor M
		\arrow{r}{}
		& N \tensor[B] M
	\end{tikzcd}
\end{equation}
defines	$N \tensor[B] M \in \category{C}_0$, which is called \emph{tensor product of $N$ and $M$ over $B$}.
\end{definition}

This construction can actually be seen as a functor as follows.

\begin{proposition}
	\label{prop: tensor product as functor}
	Let $\category{C}$ be a monoidal category with coequalizers and let
	$B \in\Monoids(\category{C})$ be a monoid.
	The tensor product of modules over $B$ defines a functor
	\begin{equation}
		\tensor[B] \colon \RightModules[\category{C}]{B} \times \LeftModules{B}_\category{C} \longrightarrow \category{C}.
	\end{equation}
\end{proposition}

\begin{proof}
	Let
	$f \colon M \longrightarrow M^\prime$
	and
	$g \colon N \longrightarrow N^\prime$
	be morphisms between left $B$-modules
	$M,M^\prime \in \LeftModules{B}_\category{C}$
	and right $B$-modules
	$N, N^\prime \in \RightModules[\category{C}]{B}$,
	respectively.
	Together with the coequalizer property of
	$N \tensor[B] M$
	and
	$N^\prime \tensor[B] M^\prime$ 
	we get
	\begin{equation*}
		\begin{tikzcd}
			N \tensor B \tensor M
			\arrow[r, shift left,"\rho_N \tensor \id_M"]
			\arrow[r, shift right, "\id_N \tensor \lambda_M", swap]
			\arrow{d}[swap]{g \tensor \id_B \tensor f}
			&N \tensor M
			\arrow{r}{p}
			\arrow{d}{g \tensor f}
			& N \tensor[B] M \\
			N^\prime \tensor B \tensor M^\prime
			\arrow[r, shift left,"\rho_{N^\prime} \tensor \id_{M^\prime}"]
			\arrow[r, shift right, "\id_{N^\prime} \tensor \lambda_{M^\prime}", swap]
			& N^\prime \tensor M^\prime
			\arrow{r}{p^\prime}
			& N^\prime \tensor[B] M^\prime
		\end{tikzcd},
	\end{equation*}
	where $p$ and $p^\prime$ denote the coequalizers of the given actions.
	The left square of this diagram commutes for both actions.
	Indeed, 
	\begin{align*}
		(\id_{N^\prime} \tensor \lambda_{M^\prime}) \circ (g \tensor \id_B \tensor f)
		&= (\id_{N^\prime} \circ g) \tensor (\lambda_{M^\prime} \circ \id_B \tensor f) \\
		&=  (g \circ \id_N) \tensor (f \circ \lambda_M) \\
		&= (g \tensor f) \circ (\id_N \tensor \lambda_M)
	\end{align*}
	and
	\begin{align*}
		(\rho_N^\prime \tensor \id_{M^\prime}) \circ (g \tensor \id_B \tensor f)	
		&= (\rho_{N^\prime} \circ g \tensor \id_B) \tensor (\id_{M^\prime} \circ f) \\
		&= (g \circ \rho_N) \tensor (f \circ \id_M) \\
		&= (g \tensor f) \circ ( \rho_N \tensor \id_M)
	\end{align*}
	holds due to $f$ and $g$ being module morphisms and the functoriality of $\tensor$.
	Thus it follows that
	\begin{equation*}
		(p^\prime \circ g \tensor f) \circ (\rho_N \tensor \id_M)
		= (p^\prime \circ g \tensor f) \circ (\id_N \tensor \lambda_M),
	\end{equation*}
	and hence by the universal property of $N \tensor[B] M$ there exists a unique morphism
	$g \tensor[B] f \colon N \tensor[B] M \to N^\prime \tensor[B] M^\prime$.
	Then $\id_N \tensor[B] \id_M = \id(N \tensor[B] M)$ is clear, and the fact that $\tensor[B]$ respects composition can be 
	deduced by adding another coequalizer diagram at the bottom of the diagram above.
\end{proof}

In case that $N$ is not only a right $B$-module but a $(C,B)$-bimodule we would like to transfer the left $C$-bimodule 
structure onto the tensor product $N \tensor[B] M$.
The same would be desirable for $M$ being a $(B,A)$-bimodule.
To achieve this we need that the tensor product of the monoidal category preserves coequalizers in the following sense:
we say that $\tensor$ \emph{preserves coequalizers} in the first component if for every coequalizer
\begin{equation}
	\begin{tikzcd}
		A
		\arrow[r, shift left,"f"]
		\arrow[r, shift right, "g", swap]
		& B
		\arrow{r}{p}
		& Q
	\end{tikzcd}
\end{equation}
and every $C \in \category{C}_0$
\begin{equation}
	\begin{tikzcd}
		A \tensor C
		\arrow[r, shift left,"f \tensor \id_C"]
		\arrow[r, shift right, "g \tensor \id_C", swap]
		& B \tensor C
		\arrow{r}{p \tensor \id_C}
		& Q \tensor C
	\end{tikzcd}
\end{equation}
is a coequalizer.
Analogously, we define preserving coequalizers in the second component.

\begin{proposition}[Tensor product of bimodules]
	\label{prop: tensor product of bimodules}
	\index{tensor product!module objects}
Let $\category{C}$ be a monoidal category with coequalizers such that $\tensor$ preserves coequalizers in both 
components.
For
$M \in \Bimodules_\category{C}(B,A)$ and
$N \in \Bimodules_\category{C}(C,B)$
the tensor product
$N \tensor[B] M$
over $B$ is a $(C,A)$-bimodule with actions given by
$\id_N \tensor[B] \rho_M$
and
$\lambda_N \tensor[B] \id_M$.
\end{proposition}

\begin{proof}
	We only construct the right $A$-module structure on $N \tensor[B] M$.
	The left $C$-module structure can then be defined in a completely analogous fashion.
	Since $\tensor$ preserves coequalizers in the first component we get two coequalizers
	\begin{equation*}
		\begin{tikzcd}[row sep = small]
			{ }
			& { }
			& { }
			& N \tensor[B] (M \tensor C) \\
			N \tensor B \tensor M \tensor C
			\arrow[rr, shift left,"\rho_N \tensor \id_M \tensor \id_C"]
			\arrow[rr, shift right, "\id_N \tensor \lambda_M \tensor \id_C", swap]
			& { }
			& N \tensor M \tensor C
			\arrow[ru, start anchor = north east, end anchor = south west, "p^\prime"]
			\arrow[rd, start anchor = south east, end anchor = north west, "p \tensor \id_C", swap]
			&{ } \\
			{ }
			&{ }
			& { }
			& (N \tensor[B] M) \tensor C
		\end{tikzcd}.
	\end{equation*}
	Then if follows from the universal property of coequalizers that
	$N \tensor[B] ( M \tensor C) \simeq (N \tensor[B] M) \tensor C$.
	Thus we can define the right $C$ action of $N \tensor[B] M$ by
	$\id_N \tensor[B] \rho_M
	\colon  (N \tensor[B] M) \tensor C \simeq N \tensor[B] ( M \tensor C)
	\longrightarrow N \tensor[B] M$.
	Using the associativity isomorphisms of the tensor product one can easily verify that the such defined
	left and right actions indeed commute.
\end{proof}

With this the functor $\tensor[B]$ from \autoref{prop: tensor product as functor} gives directly a functor
\begin{equation}
	\tensor[B] \colon \Bimodules_\category{C}(C,B) \times \Bimodules_\category{C}(B,A)
	\longrightarrow \Bimodules_\category{C}(C,A).
\end{equation}

Putting all of these constructions together we are finally able to construct a bicategory with objects given by monoids in a monoidal
category. The 1-morphisms of this bicategory are the bimodules over the given monoids and 2-morphisms are bimodule 
morphisms.

\begin{theorem}[The bicategory $\Bimodules_\category{C}$]
	\label{thm: the bicategory Bimod_C}
	Let $\category{C}$ be a monoidal category with coequalizers such that $\tensor$ preserves coequalizers in both components.
	Then the following data defines a bicategory $\Bimodules_\category{C}$:
	\begin{theoremlist}
		\item The objects $(\Bimodules_\category{C})_0 = \Monoids(\category{C})$ are given by the monoids in $\category{C}$.
		\item For any two monoids $A,B \in \Monoids(\category{C})$ the category of 1- and 2-morphisms from $B$ to $A$ is given by
		$\Bimodules_\category{C}(B,A)$, the category of $(B,A)$-bimodules together with bimodule homomorphisms, see
		\autoref{cor: categories of modules}.
		\item For any three monoids $A,B,C \in \Monoids(\category{C})$ the tensor product functor
		\begin{equation}
			\tensor[B] \colon \Bimodules_\category{C}(C,B) \times \Bimodules_\category{C}(B,A)
			\longrightarrow \Bimodules_\category{C}(C,A)
		\end{equation}
		given as in \autoref{prop: tensor product of bimodules}.
		\item For each monoid $A \in \Monoids(\category{C})$ the identity over $A$ is $A$ itself considered as an
		$(A,A)$-bimodule.
	\end{theoremlist}
\end{theorem}

\begin{proof}
	The statements listed hold true due to the referred results.
	In order to define a bicategory the existence of natural transformations for associativity and left/right units are required.
	These can be constructed by transferring the natural isomorphisms of the monoidal category $\category{C}$,
	but this involves a lot of small statements to be checked, as well as some not so trivial arguments for functor categories.
	Surprisingly, this result seems to be common knowledge among category theorists, but there is no publication doing exactly 
	these computations (as far as the author knows), so we cannot give a reference for this proof.
	Nevertheless, there are some generalizations of this construction available in the literature.
	See \cite[Sec. 5.3]{leinster:2004a} for a discussion using generalized multicategories or
	\cite[Sec. 2]{haugseng:2017a}
	using $\infty$-categories.
\end{proof}

%
\begin{example}\
	\begin{examplelist}
		\item The bicategory $\Bimodules_\Ab$ is the bicategory of rings, modules and module homomorphisms.
		\item The bicategory $\Bimodules_{\Bimodules(\ring{R})}$ is the bicategory of $\ring{R}$-algebras, bimodules and module
		homomorphisms.
	\end{examplelist}
\end{example}
\section{Reflection Theorems}
\label{sec:Reflection}

Consider a functor $\functor{U} \colon \category{X} \to \category{C}$
from an arbitrary category to a monoidal category $(\category{C}, \tensor)$.
Under which conditions can we ``pull back'' the monoidal structure from $\category{C}$
to a monoidal structure on $\category{X}$?
This is clearly possible if $\functor{U}$ is an equivalence of categories, but it is also possible in a more general setup.

\begin{definition}[Reflective subcategory]
	\index{reflective subcategory}
Let $\category{C}$ be a category.
\begin{definitionlist}
	\item A \emph{full subcategory} of $\category{C}$
	is a category $\category{X}$
	together with a fully faithful functor
	$\functor{U} \colon \category{X} \to \category{C}$.
	\item A \emph{reflective subcategory} of $\category{C}$
	is a full subcategory $\functor{U} \colon \category{X} \to \category{C}$
	such that $\functor{U}$ admits a left adjoint.
\end{definitionlist}
\end{definition}

Having a reflective subcategory is close to having an equivalence of categories
as the following lemma shows.

\begin{lemma}
	Let $\functor{U} \colon \category{X} \to \category{C}$
	be a functor.
	Then the following statements are equivalent:
	\begin{lemmalist}
		\item $\category{X}$ is a reflective subcategory of $\category{C}$ via $\functor{U}$.
		\item The functor $\functor{U}$ has a left adjoint $\functor{F} \colon \category{C} \to \category{X}$ such that the counit
		$\varepsilon \colon \functor{F}\functor{U} \Rightarrow \id_{\category{X}}$
		is a natural isomorphism.
	\end{lemmalist}
\end{lemma}

In \cite{day:1970a,day:1972a} Day gave a list of equivalent conditions under which a closed symmetric monoidal structure on $\category{C}$ induces a closed symmetric monoidal structure on $\category{X}$ such that $\functor{F}$ becomes monoidal.
We give a simplified version of Day's reflection theorem for monoidal categories without any additional closedness or symmetry requirement.
The proof follows the idea of \cite{day:1970a}.

\begin{theorem}[Reflection Theorem]
	\label{thm:ReflectionTheorem}
	\index{reflection theorem}
Let $\category{C}$ be a monoidal category and let
$\functor{U} \colon \category{X} \to \category{C}$
be a reflective subcategory with unit denoted by
$\eta \colon \id_\category{C} \Rightarrow \functor{U}\functor{F}$.
If
\begin{equation}
	\functor{F}(\eta \tensor \eta)_{A,B} \colon \functor{F}(A \tensor B) \to \functor{F}(\functor{U}\functor{F}(A) \tensor \functor{U}\functor{F}(B))
\end{equation}
is an isomorphism for all $A,B \in \category{C}$,
then there exists a monoidal structure on $\category{X}$
such that $\functor{F}$ becomes a monoidal functor,
which is unique up to monoidal equivalence.	
\end{theorem}

\begin{proof}
	Define $\hat\tensor \colon \category{X} \times \category{X} \to \category{X}$
	by 
	$\hat\tensor \coloneqq \functor{F} \circ \tensor \circ (\functor{U} \times \functor{U})$
	and $\hat\Unit \coloneqq \functor{F}\Unit$.
	We define the natural isomorphisms
	$\widehat\lidentity \colon \hat\tensor \circ (\Id_\Unit \times \id) \rightarrow \id$
	and
	$\widehat\ridentity \colon \hat\tensor \circ (\id \times \Id_\Unit) \rightarrow \id$
	by setting $\widehat\lidentity_X$ and $\widehat\ridentity_X$
	for every $X \in \category{X}$ as the unique isomorphisms 
	making
	\begin{equation*}
	\begin{tikzcd}
		\functor{F}(\Unit \tensor \functor{U}X)
			\arrow[r,"\functor{F}\lidentity_{\functor{U}X}"]
		& \functor{F}\functor{U}X
			\arrow[dd,"\varepsilon_X"]\\
		\functor{F}(\functor{U}\functor{F}\Unit \tensor \functor{U}\functor{F}\functor{U}X)
			\arrow[u,"\functor{F}(\eta \tensor \eta)_{\Unit,\functor{U}X}^{-1}"]
		& \\
		\functor{F}(\functor{U}\functor{F}\Unit \tensor \functor{U}X)
			\arrow[u,"\functor{F}(\id \tensor \functor{U}\varepsilon^{-1}_{X})"]
			\arrow[r, "\widehat{\lidentity}_X"]
		& X
	\end{tikzcd}
	\quad\text{ and }\quad
	\begin{tikzcd}
		\functor{F}( \functor{U}X \tensor \Unit)
			\arrow[r,"\functor{F}\ridentity_{\functor{U}X}"]
		& \functor{F}\functor{U}X
			\arrow[dd,"\varepsilon_X"]\\
		\functor{F}(\functor{U}\functor{F}\functor{U}X \tensor \functor{U}\functor{F}\Unit )
			\arrow[u,"\functor{F}(\eta \tensor \eta)_{\functor{U}X,\Unit}^{-1}"]
		& \\
		\functor{F}(\functor{U}X \tensor \functor{U}\functor{F}\Unit)
			\arrow[u,"\functor{F}(\functor{U}\varepsilon^{-1}_{X} \tensor \id)"]
			\arrow[r, "\widehat{\ridentity}_X"]
		& X
	\end{tikzcd}
	\end{equation*}
commute.
Similarly, we define
$\widehat\asso \colon \hat\tensor \circ ( \hat\tensor \times \id ) \Rightarrow \hat\tensor \circ (\id \times \hat\tensor)$
by setting $\widehat\asso_{X,Y,Z}$ as the unique isomorphism making
\begin{equation*}
\begin{tikzcd}
	\functor{F}((\functor{U}X \tensor \functor{U}Y) \tensor \functor{U}Z)
		\arrow[r,"\functor{F}\asso_{\functor{U}X,\functor{U}Y,\functor{U}Z}"]
	&\functor{F}(\functor{U}X \tensor (\functor{U}Y \tensor \functor{U}Z))
		\arrow[d,"\functor{F}(\eta \tensor\eta)_{\functor{U}X,\functor{U}Y\tensor\functor{U}Z}"]\\
	\functor{F}\left( \functor{U}\functor{F}(\functor{U}X \tensor \functor{U}Y) \tensor \functor{U}\functor{F}\functor{U}Z \right)
		\arrow[u,"\functor{F}(\eta \tensor \eta)^{-1}_{\functor{U}X \tensor \functor{U}Y,\functor{U}Z}"]
	&\functor{F}\left( \functor{U}\functor{F}\functor{U}X \tensor \functor{U}\functor{F}(\functor{U}Y \tensor \functor{U}Z )\right)
		\arrow[d,"\functor{F}(\functor{U}\varepsilon_{X} \tensor \id)"]\\
	\functor{F}(\functor{U}\functor{F}(\functor{U}X \tensor \functor{U}Y) \tensor \functor{U}Z)
		\arrow[u,"\functor{F}(\id \tensor \functor{U}\varepsilon^{-1}_{Z})"]
		\arrow[r,"\widehat{\asso}_{X,Y,Z}"]	
	&\functor{F}(\functor{U}X \tensor \functor{U}\functor{F}(\functor{U}Y \tensor \functor{U}Z))
\end{tikzcd}
\end{equation*}
commute.
We need to check the coherences.
For the identity coherence take the above defining diagrams for $\widehat{\ridentity}$
and $\widehat{\lidentity}$ and take the $\hat\tensor$-product with $X$ and $Y$, respectively.
Gluing the resulting diagrams together yields the following:
\begin{equation*}
\begin{tikzcd}[column sep=tiny, row sep=huge]
	\functor{F}((\functor{U}X \tensor \Unit) \tensor \functor{U}Y)
		\arrow[rddd, start anchor={west}, end anchor={west},bend right=60pt,"\functor{F}(\ridentity_{\functor{U}X} \tensor \id)"{swap}] 
		\arrow[rr,"\functor{F}\asso_{\functor{U}X,\Unit,\functor{U}Y}"] 
	& {}                                                      
	& \functor{F}(\functor{U}X \tensor (\Unit \tensor \functor{U}Y))
		\arrow[lddd, start anchor={east}, end anchor={east}, 
		bend left=60pt,"\functor{F}(\id \tensor \lidentity_{\functor{U}Y})"] \\
	\functor{F}(\functor{U}\functor{F}(\functor{U}X \tensor \functor{U}\hat\Unit) \tensor \functor{U}Y)
		\arrow[rd,"\functor{F}(\functor{U}\widehat{\ridentity}_X \tensor \id)"{swap}] 
		\arrow[rr,"\widehat{\asso}_{X,\hat\Unit,Y}"]
		\arrow[u] 
	& {}
	& \functor{F}(\functor{U}X \tensor \functor{U}\functor{F}(\functor{U}\hat\Unit \tensor \functor{U}Y))  
		\arrow[ld,"\functor{F}(\id \tensor \functor{U}\widehat{\lidentity}_{Y})"]
		\arrow[u] \\
	{}                               
	& \functor{F}(\functor{U}X \tensor \functor{U}Y)
		\arrow[d] 
	&   {} \\
	{}                                    
	& \functor{F}(\functor{U}X \tensor \functor{U}Y)           
	& {}          
\end{tikzcd}
\end{equation*}
The inner triangle is the identity coherence we need to verify,
while the outer triangle commutes by the identity coherence for $\tensor$.
The unlabelled morphisms from the inner triangle to the outer triangle are given by the sides of
defining diagrams for $\widehat{\asso}$, $\widehat{\lidentity}$ and $\widehat{\ridentity}$.
Since these commute we obtain the identity coherence for $\hat{\tensor}$.
With the same strategy we can glue the defining diagram for $\widehat{\asso}$
to every edge of the associativity coherence for $\hat\tensor$.
Then the outer pentagon is the associativity coherence for $\tensor$, showing that the inner pentagon also commutes.
Since this diagram becomes too large, we refer instead to
\cite{day:1970a}.

To show the uniqueness suppose $\tilde{\tensor}$ is another monoidal structure on $\category{X}$ such that $\functor{F}$ becomes monoidal.
Then the identity functor on $\category{X}$ yields a monoidal equivalence, since
we have natural isomorphism implementing
\begin{equation*}
	X \tilde\tensor Y
	\simeq \functor{F}\functor{U}X \tilde\tensor \functor{F}\functor{U}Y
	\simeq \functor{F}(\functor{U}X \tensor \functor{U}Y)
	\simeq X \hat\tensor Y
\end{equation*}
for all $X,Y \in \category{X}$.
\end{proof}

%% file: poisson-geometry.tex
For the convenience of the reader we give some basic definitions and results from Poisson geometry 
and coisotropic reduction.
All of this can be found in a similar fashion in standard text books like
\cite{marsden.ratiu:1999a} and
\cite{crainic.fernandes.marcut:2021a}.
See also \cite{ortega.ratiu:2004a} for an in-depth treatment of various reductions schemes.

\begin{definition}[Poisson algebra]
	\label{def: poisson algebra}
	Let $\mathbb{k}$ be a commutative unital ring with $1 \neq -1$.
	A \emph{Poisson algebra} is a pair $(\algebra{A}, \{ \argument , \argument \})$ of an associative algebra 
	$\algebra{A}$ over $\mathbb{k}$ and a $\mathbb{k}$-linear map
	\begin{equation}
		\{ \argument, \argument \} \colon \algebra{A} \tensor[\field{k}] \algebra{A} \longrightarrow \algebra{A}
	\end{equation}
	fulfilling the following properties for all $a,b,c \in \algebra{A}$:
	\begin{definitionlist}
		\item Antisymmetry:
		\begin{equation}
			\{ a , b \} = - \{b , a \}
		\end{equation}
		\item Jacobi identity:
		\begin{equation}
			\{ a, \{b , c\} \} = \{ \{ a, b \} , c\} + \{b , \{a,c \} \}
		\end{equation}
		\item Leibniz rule:
		\begin{equation}
			\{a , bc \} = \{a ,b\}c + b \{a , c\}
		\end{equation}
	\end{definitionlist}
	The map $\{ \argument , \argument \}$ is called \emph{Poisson bracket}.
\end{definition}

A \emph{morphism of Poisson} algebras
$\Phi \colon (\algebra{A}, \{ \argument, \argument \}_\algebra{A}) 
\longrightarrow (\algebra{B} , \{ \argument, \argument \}_\algebra{B})$
between two Poisson algebras is an algebra homomorphism
$\Phi \colon \algebra{A} \longrightarrow \algebra{B}$
with
$\Phi(\{ a, a^\prime \}_\algebra{A}) = \{ \Phi(a) , \Phi(a^\prime) \}_\algebra{B}$
for all $a, a^\prime \in \algebra{A}$.

\begin{definition}[Poisson manifold]
	\label{def: poisson manifold}
	A \emph{Poisson manifold} is a pair $(M, \{ \argument , \argument \})$ of a smooth manifold $M$ together with a map
	$\{ \argument , \argument \} \colon \Cinfty(M) \tensor[\field{k}] \Cinfty(M) \longrightarrow \Cinfty(M)$ turning
	$\Cinfty(M)$ into a Poisson algebra.
\end{definition}

A Poisson map
$\Phi \colon (M_1, \{ \argument , \argument \}_1) \longrightarrow (M_2, \{ \argument , \argument \}_2)$
between two Poisson manifolds is a smooth map $\Phi \colon M_1 \to M_2$
such that
$\Phi^* \colon \Cinfty(M_2) \longrightarrow \Cinfty(M_1)$
is a morphism of Poisson algebras.
By antisymmetry and Leibniz rule every Poisson bracket $\{ \argument , \argument \}$ is a biderivation, hence for every 
Poisson manifold $(M, \{ \argument , \argument \})$ there exists a bivector field $\pi \in \Secinfty(\Anti^2 TM)$, called 
\emph{Poisson tensor}, such that
\begin{equation}
	\{f , g \} = \pi(\D f \tensor \D g)
\end{equation}
for all $f,g \in \Cinfty(M)$.
Hence we will also denote a Poisson manifold by $(M, \pi)$ if we want to stress the Poisson tensor.
Every such Poisson tensor induces a musical homomorphism
\begin{equation}
	\argument^\# \colon T^*M \ni \alpha_p \mapsto \alpha^\#_p \coloneqq \pi_p(\argument,\alpha_p) \in TM,
\end{equation}
which allows us to define the \emph{Hamiltonian vector field} 
\begin{equation}
	X_f \coloneqq (\D f)^\#
\end{equation}
for every $f \in \Cinfty(M)$.

\begin{definition}[Coisotropic submanifold]
	\label{def: coisotropic submanifold_poisson}
	\index{coisotropic submanifold}
Let $(M, \pi)$ be a Poisson manifold. A submanifold $C \subseteq M$ is called \emph{coisotropic} if
$T_p C \subseteq T_p M$ is a coisotropic subspace for all $p \in C$, i.e. if
\begin{equation}
	\label{eq: def: coisotropic submanifold_poisson}
	(T_p C^\ann )^\sharp \subseteq T_p C
\end{equation}
holds for all $p \in C$.
\end{definition}

%

We can always view a Poisson manifold in two different ways:
we can either focus on a geometric description as a pair $(M, \pi)$ of a manifold with additional structure, or on an algebraic 
description by considering instead the Poisson algebra $(\Cinfty(M), \{\argument, \argument \})$.
In a similar way we can assign to a submanifold $C \subseteq M$ an algebraic object, the
\emph{vanishing ideal}
of $C$
\begin{equation}
	\label{eq: vanishing ideal}
	\vanishing_C := \left\{ f \in \Cinfty(M) \mid f\at{C} = 0 \right\}.
\end{equation}

\begin{proposition}[Vanishing ideal]
	\label{prop: vanishing ideal}
Let $M$ be a manifold with closed submanifold $\iota \colon C \longrightarrow M$.
\begin{propositionlist}
	\item The vanishing ideal $\vanishing_C$ of $C$ is an ideal inside the algebra $\Cinfty(M)$.
	\item The algebras $\Cinfty(M) / \vanishing_C$ and $\Cinfty(C)$ are isomorphic via the map
	\begin{equation}
		\label{eq: prop: vanishing ideal}
		\Cinfty(M) / \vanishing_C 
		\ni
		[f]
		\longmapsto
		\iota^*f
		\in
		\Cinfty(C).
	\end{equation}
\end{propositionlist}
\end{proposition}
The inverse of \eqref{eq: prop: vanishing ideal} can be constructed using a tubular neighbourhood.
The following proposition gives a geometric and an algebraic characterization of coisotropic submanifolds using the vanishing 
ideal.

\begin{proposition}[Coisotropic submanifolds]
	\label{prop: coisotropic subamanifolds}
	Let $(M, \pi)$ be a Poisson manifold and let $C \subseteq M$ be a submanifold.
	Then the following statements are equivalent:
	\begin{propositionlist}
		\item \label{prop: coisotropic subamanifolds_item 1}
		The submanifold $C$ is coisotropic.
		\item \label{prop: coisotropic subamanifolds_item 2}
		For all $f \in \vanishing_C$ it holds $X_f(p) \in T_p C$ for all $p \in C$.
		\item \label{prop: coisotropic subamanifolds_item 3}
		The vanishing ideal $\vanishing_C$ is a Poisson subalgebra of $\Cinfty(M)$.
	\end{propositionlist}
\end{proposition}

The distribution generated by the Hamiltonian vector fields of functions vanishing on the coisotropic submanifold
will play an important role in coisotropic reduction.

\begin{definition}[Characteristic distribution]
	\label{def: characteristic distribution}
	Let $(M, \pi)$ be a Poisson manifold with coisotropic submanifold $C \subseteq M$.
	The distribution on $C$ spanned by the Hamiltonian vector fields $X_f$ of a function
	$f \in \vanishing_C$ is called the \emph{characteristic distribution} of $C$.
\end{definition}

It can then be shown that this is in fact an integrable distribution on $C$, given by the subspace
$(T_pC^\ann)^\sharp \subseteq T_pC$ at every point $p \in C$.
Again, we would like to have an equivalent algebraic description of this quite geometric notion, similar to the ones we 
presented for Poisson manifolds and coisotropic submanifolds.
For this we first need the following construction.

\begin{definition}[Poisson normalizer]
	\label{def: poisson normalizer}
Let $\algebra{A}$ be a Poisson algebra and let $\vanishing \subseteq \algebra{A}$ be an ideal for the associative and 
commutative product as well as a Poisson subalgebra.
The Poisson subalgebra given by
\begin{equation}	
	\label{eq: def: poisson normalizer}
	\Pnormalizer_\vanishing = \{ a \in \algebra{A} \mid \{a ,\vanishing \} \subseteq \vanishing \}
\end{equation}
is called the \emph{Poisson normalizer} of $\vanishing$.
\end{definition}

It is clear that $\Pnormalizer_\vanishing$ is the largest Poisson subalgebra containing $\vanishing$ as a Poisson ideal.
If $\vanishing = \vanishing_C$ is the vanishing ideal of a submanifold $C \subseteq M$ we will simply write $\Pnormalizer_C$.
Thus we can always assign a Poisson normalizer $\Pnormalizer_C \subseteq \Cinfty(M)$ to a coisotropic submanifold 
$C \subseteq M$ of a Poisson manifold $M$.
In general one even calls an ideal $\vanishing \subseteq \algebra{A}$ in a Poisson algebra \emph{coisotropic} if $\vanishing$ is 
in addition a Poisson subalgebra.
This Poisson normalizer now encodes the same information as the characteristic distribution of a coisotropic submanifold, thus 
giving us the algebraic formulation we were searching for.

\begin{proposition}
	\label{prop: functions in the poisson normalizer}
Let $(M,\pi)$ be a Poisson manifold with coisotropic submanifold $C \subseteq M$.
For a function $f \in \Cinfty(M)$ the following statements are equivalent:
\begin{propositionlist}
	\item \label{item: prop: functions in the poisson normalizer_1}
	One has $f \in \Pnormalizer_C$.
	\item \label{item: prop: functions in the poisson normalizer_2}
	The Hamiltonian vector field $X_f$ is tangent to $C$.
	\item \label{item: prop: functions in the poisson normalizer_3}
	The function $\iota^*f \in \Cinfty(C)$ is constant along the leaves of the characteristic foliation of $C$.
	\item \label{item: prop: functions in the poisson normalizer_4}
	One has $\Lie_{\iota^*X_g} \iota^*f = 0$ for all $g \in \vanishing_C$, where $\iota^*X_g$ denotes the restriction of 
	$X_g$ to $C$, which is possible since $C$ is coisotropic.
\end{propositionlist}
\end{proposition}

Finally, we want to identify all points along the leaves of the characteristic distribution, thus obtaining a quotient
$M_\red = C / \operatorname{\sim}$.
If this quotient is indeed a manifold we can equip the algebra of functions $\Cinfty(M_\red)$ with a Poisson structure 
$\pi_\red$.
The corresponding quotient on the algebraic side is given by $\Pnormalizer_C / \vanishing_C$, which is a Poisson 
algebra consisting of the functions on $C$ that are constant along the leaves of the characteristic distribution.
The observation that the geometric and the algebraic description lead to essentially the same reduction is made precise by the 
following theorem.

\begin{theorem}[Coisotropic reduction]
	\label{thm: coisotropic reduction}
Let $(M,\pi)$ be a Poisson manifold with a closed coisotropic submanifold $C \subseteq M$ such that the projection
\begin{equation}
	\pi \colon M \longrightarrow C / \operatorname{\sim} = M_\red
\end{equation}
is a surjective submersion for a smooth structure on the leaf space $M_\red$ of the characteristic distribution of $C$.
Then there exists a unique Poisson structure $\pi_\red$ on $M_\red$ such that
\begin{equation}
	\Pnormalizer_C / \vanishing_C 
	\ni
	[f]
	\longmapsto
	\iota^* f
	\in
	\pi^* \Cinfty(M_\red)
	\subseteq \Cinfty(C)
\end{equation}
is an isomorphism of Poisson algebras.
\end{theorem}